\documentclass[reqno]{amsart}

\usepackage{fullpage,graphicx}
\usepackage[margin=1in]{geometry}
\usepackage{amsmath,amsthm}
\usepackage{amsfonts,amssymb}
\usepackage{newlfont,mathrsfs,mathtools}
\usepackage{txfonts}
\usepackage{etoolbox}
\usepackage{comment}
\usepackage{tikz-cd}
\usepackage{float}

\makeatletter
\let\old@tocline\@tocline
\let\section@tocline\@tocline
\newcommand{\subsection@dotsep}{4.5}
\newcommand{\subsubsection@dotsep}{4.5}
\patchcmd{\@tocline}
{\hfil}
{\nobreak
	\leaders\hbox{$\m@th
		\mkern \subsection@dotsep mu\hbox{.}\mkern \subsection@dotsep mu$}\hfill
	\nobreak}{}{}
\let\subsection@tocline\@tocline
\let\@tocline\old@tocline

\patchcmd{\@tocline}
{\hfil}
{\nobreak
	\leaders\hbox{$\m@th
		\mkern \subsubsection@dotsep mu\hbox{.}\mkern \subsubsection@dotsep mu$}\hfill
	\nobreak}{}{}
\let\subsubsection@tocline\@tocline
\let\@tocline\old@tocline

\let\old@l@subsection\l@subsection
\let\old@l@subsubsection\l@subsubsection

\def\@tocwriteb#1#2#3{%
	\begingroup
	\@xp\def\csname #2@tocline\endcsname##1##2##3##4##5##6{%
		\ifnum##1>\c@tocdepth
		\else \sbox\z@{##5\let\indentlabel\@tochangmeasure##6}\fi}%
	\csname l@#2\endcsname{#1{\csname#2name\endcsname}{\@secnumber}{}}%
	\endgroup
	\addcontentsline{toc}{#2}%
	{\protect#1{\csname#2name\endcsname}{\@secnumber}{#3}}}%

\newlength{\@tocsectionindent}
\newlength{\@tocsubsectionindent}
\newlength{\@tocsubsubsectionindent}
\newlength{\@tocsectionnumwidth}
\newlength{\@tocsubsectionnumwidth}
\newlength{\@tocsubsubsectionnumwidth}
\newcommand{\settocsectionnumwidth}[1]{\setlength{\@tocsectionnumwidth}{#1}}
\newcommand{\settocsubsectionnumwidth}[1]{\setlength{\@tocsubsectionnumwidth}{#1}}
\newcommand{\settocsubsubsectionnumwidth}[1]{\setlength{\@tocsubsubsectionnumwidth}{#1}}
\newcommand{\settocsectionindent}[1]{\setlength{\@tocsectionindent}{#1}}
\newcommand{\settocsubsectionindent}[1]{\setlength{\@tocsubsectionindent}{#1}}
\newcommand{\settocsubsubsectionindent}[1]{\setlength{\@tocsubsubsectionindent}{#1}}

\renewcommand{\l@section}{\section@tocline{1}{\@tocsectionvskip}{\@tocsectionindent}{}{\@tocsectionformat}}%
\renewcommand{\l@subsection}{\subsection@tocline{2}{\@tocsubsectionvskip}{\@tocsubsectionindent}{}{\@tocsubsectionformat}}%
\renewcommand{\l@subsubsection}{\subsubsection@tocline{3}{\@tocsubsubsectionvskip}{\@tocsubsubsectionindent}{}{\@tocsubsubsectionformat}}%
\newcommand{\@tocsectionformat}{}
\newcommand{\@tocsubsectionformat}{}
\newcommand{\@tocsubsubsectionformat}{}
\expandafter\def\csname toc@1format\endcsname{\@tocsectionformat}
\expandafter\def\csname toc@2format\endcsname{\@tocsubsectionformat}
\expandafter\def\csname toc@3format\endcsname{\@tocsubsubsectionformat}
\newcommand{\settocsectionformat}[1]{\renewcommand{\@tocsectionformat}{#1}}
\newcommand{\settocsubsectionformat}[1]{\renewcommand{\@tocsubsectionformat}{#1}}
\newcommand{\settocsubsubsectionformat}[1]{\renewcommand{\@tocsubsubsectionformat}{#1}}
\newlength{\@tocsectionvskip}
\newcommand{\settocsectionvskip}[1]{\setlength{\@tocsectionvskip}{#1}}
\newlength{\@tocsubsectionvskip}
\newcommand{\settocsubsectionvskip}[1]{\setlength{\@tocsubsectionvskip}{#1}}
\newlength{\@tocsubsubsectionvskip}
\newcommand{\settocsubsubsectionvskip}[1]{\setlength{\@tocsubsubsectionvskip}{#1}}

\patchcmd{\tocsection}{\indentlabel}{\makebox[\@tocsectionnumwidth][l]}{}{}
\patchcmd{\tocsubsection}{\indentlabel}{\makebox[\@tocsubsectionnumwidth][l]}{}{}
\patchcmd{\tocsubsubsection}{\indentlabel}{\makebox[\@tocsubsubsectionnumwidth][l]}{}{}

\newcommand{\@sectypepnumformat}{}
\renewcommand{\contentsline}[1]{%
	\expandafter\let\expandafter\@sectypepnumformat\csname @toc#1pnumformat\endcsname%
	\csname l@#1\endcsname}
\newcommand{\@tocsectionpnumformat}{}
\newcommand{\@tocsubsectionpnumformat}{}
\newcommand{\@tocsubsubsectionpnumformat}{}
\newcommand{\setsectionpnumformat}[1]{\renewcommand{\@tocsectionpnumformat}{#1}}
\newcommand{\setsubsectionpnumformat}[1]{\renewcommand{\@tocsubsectionpnumformat}{#1}}
\newcommand{\setsubsubsectionpnumformat}[1]{\renewcommand{\@tocsubsubsectionpnumformat}{#1}}
\renewcommand{\@tocpagenum}[1]{%
	\hfill {\mdseries\@sectypepnumformat #1}}

\let\oldappendix\appendix
\renewcommand{\appendix}{%
	\leavevmode\oldappendix%
	\addtocontents{toc}{%
		\protect\settowidth{\protect\@tocsectionnumwidth}{\protect\@tocsectionformat\sectionname\space}%
		\protect\addtolength{\protect\@tocsectionnumwidth}{2em}}%
}
\makeatother

\makeatletter
\settocsectionnumwidth{2em}
\settocsubsectionnumwidth{2.5em}
\settocsubsubsectionnumwidth{3em}
\settocsectionindent{1pc}%
\settocsubsectionindent{\dimexpr\@tocsectionindent+\@tocsectionnumwidth}%
\settocsubsubsectionindent{\dimexpr\@tocsubsectionindent+\@tocsubsectionnumwidth}%
\makeatother

\settocsectionvskip{0pt}
\settocsubsectionvskip{0pt}
\settocsubsubsectionvskip{0pt}

\settocsectionformat{\bfseries}
\settocsubsectionformat{\mdseries}
\settocsubsubsectionformat{\mdseries}
\setsectionpnumformat{\bfseries}
\setsubsectionpnumformat{\mdseries}
\setsubsubsectionpnumformat{\mdseries}

\let\oldtableofcontents\tableofcontents
\renewcommand{\tableofcontents}{%
	\vspace*{-\linespacing}
	\oldtableofcontents}

\allowdisplaybreaks

\usepackage[colorlinks,
linkcolor=blue,
anchorcolor=blue,
citecolor=blue
]{hyperref}

\numberwithin{equation}{section}

\newtheorem{thm}{Theorem}[section]
\newtheorem{prop}[thm]{Proposition}
\newtheorem{lem}[thm]{Lemma}
\newtheorem{df}[thm]{Definition}
\newtheorem{co}[thm]{Corollary}
\newtheorem{rem}[thm]{Remark}

\newcommand{\R}{\mathbb{R}}

\newcommand{\T}{\mathbb{T}}
\newcommand{\N}{\mathbb{N}}
\newcommand{\Z}{\mathbb{Z}}

\newcommand{\cC}{\mathcal{C}}
\newcommand{\cF}{\mathcal{F}}

\newcommand{\mcI}{\mathcal{I}}

\newcommand{\mC}{\mathcal{C}}
\newcommand{\mN}{\mathcal{N}}

\newcommand{\rd}{\mathrm{d}}
\newcommand{\Id}{\operatorname{Id}}
\newcommand{\Div}{\operatorname{div}}
\newcommand{\curl}{\operatorname{curl}}
\newcommand{\tr}{\operatorname{tr}}
\newcommand{\supp}{\operatorname{supp}}
\newcommand{\cof}{\operatorname{cof}}

\newcommand{\dlt}{\delta}

\newcommand{\sgm}{\sigma}

\newcommand{\VV}{\upsilon}

\newcommand{\pa}{\partial}

\newcommand{\tm}{\tilde{m}}
\newcommand{\trho}{\tilde{\rho}}

\newcommand{\MM}[1]{\ensuremath{\mathcal{M}}\left(#1\right)}
\newcommand{\La}[1]{\ensuremath{\Xi_\alpha}\left(#1\right)}
\newcommand{\Lat}[1]{\ensuremath{\Xi_\alpha^2}\left(#1\right)}

\newcommand{\om}{\overline{m}}

\newcommand{\orho}{\overline{\rho}}

\newcommand{\PLq}{P_{\leq\ell_q^{-1}}}

\newcommand{\DTL}{\overline{D}_t}
\newcommand{\Dtq}{D_{t,q}}

\newcommand{\inc}{^{(\mathrm{inc})}}
\newcommand{\comp}{^{(\mathrm{comp})}}

\newcommand{\JX}{J_{X_\upsilon}}
\newcommand{\JPhi}{J_{\Phi_\upsilon}}

\newcommand{\nrm}[1]{\Vert#1\Vert}
\newcommand{\Nrm}[1]{\left\Vert#1\right\Vert}

\newcommand{\Lab}[2]{\ensuremath{\Xi_{#2}}\left(#1\right)}

\newcommand{\Lk}[2]{\ensuremath{\Xi_\alpha^{#2}}\left(#1\right)}
\newcommand{\Lkb}[3]{\ensuremath{\Xi_{#3}^{#2}}\left(#1\right)}



\begin{document}

	\title[] {An Onsager-type theorem for the compressible Euler equations}
	
	\author{Shunkai Mao}
	\address{School of Mathematics and Statistics, Donghua University, China.}
	\email[Shunkai Mao]{skmao@dhu.edu.cn}
	
	\author{Peng Qu}
	\address{School of Mathematical Sciences $\&$ Shanghai Key Laboratory for Contemporary Applied Mathematics, Fudan University, China.}
	\email[Peng Qu]{pqu@fudan.edu.cn}
	\thanks{}
	\begin{abstract}
		For any $0\leq\gamma_m<\frac13$ and $0\leq\gamma_\rho<\frac23$, we use convex integration to construct a nontrivial weak solution $(\rho,m)$ to the three-dimensional isentropic compressible Euler equations such that $\rho-\rho^*\in C^{\gamma_\rho}(\R\times\T^3)$ and $m\in C^{\gamma_m}(\R\times\T^3)$. Here $\rho^*>0$ is a constant, and the density stays bounded away from vacuum. Both $\rho-\rho^*$ and $m$ have compact support in time. The construction applies to smooth pressure laws with positive derivative throughout the density range. This yields nonuniqueness for the constant initial data $(\rho^*,0)$.
		
		Each stage of the convex integration scheme consists of finitely many Newton steps followed by a Nash step. Kirchhoff's formula in Lagrangian coordinates and the Piola transform yield coupled density and momentum corrections that treat the acoustic part of the linearized system while preserving the continuity equation exactly. In the Newton step, we combine these corrections with an incompressible correction adapted from Giri--Radu~\cite{GR24} to reduce the error. In the Nash step, a divergence-free perturbation built from spatially separated Mikado flows treats the gluing stresses on disjoint time intervals without changing the density.
	\end{abstract}
	
	\maketitle
	{
		\tableofcontents
	}

	\section{Introduction}\label{Introduction}

	The isentropic compressible Euler equations describe the motion of an inviscid fluid \cite{Dafermos16,Majda84}.
	We consider this system on $\R_t\times\T_x^3$, where
	$
	\T^3=[-\pi,\pi]^3
	$ is the three-dimensional torus with opposite faces identified:
	\begin{equation}
		\left\{
		\begin{aligned}
			&\partial_t\rho+\Div(\rho u)=0,\\
			&\partial_t(\rho u)+\Div(\rho u\otimes u)+\nabla p(\rho)=0.
		\end{aligned}
		\right.
		\label{Euler System}
	\end{equation}
	Here $\rho$ denotes the mass density,
    $
    \rho:\R\times\T^3\to(0,\infty)
    $,
	and $u$ denotes the velocity field,
    $
    u:\R\times\T^3\to\R^3.
    $
	We assume that the density is bounded away from vacuum.
	The pressure is prescribed by a fixed smooth law
    $
    p=p(\rho).
    $
	A standard example is the polytropic law $p(\rho)=\kappa\rho^{\gamma}$,
	where $\kappa>0$ and $\gamma>1$.

	We work with the density and momentum $(\rho,m)$, where
	$
	m:=\rho u.
	$
	Thus
	$
	u=m/\rho.
	$
	In these variables, the system takes the form
	\begin{equation}
		\left\{
		\begin{aligned}
			&\partial_t\rho+\Div m=0,\\
			&\partial_t m+\Div\left(\frac{m\otimes m}{\rho}\right)+\nabla p(\rho)=0.
		\end{aligned}
		\right.
		\label{compressible Euler conservative form}
	\end{equation}

	The Cauchy problem is locally well posed for smooth initial data with density bounded away from vacuum,
	provided $p'>0$ on the initial density range \cite{Majda84}.
	Classical solutions may develop shocks in finite time \cite{Dafermos16},
	which motivates the study of weak solutions.
	We use the following definition of a weak solution on $\R\times\T^3$.
	\begin{df}[Weak solution]
		A pair $(\rho,m)\in L^\infty(\R\times\T^3;\R\times\R^3)$
		with $\rho\geq c_0>0$ almost everywhere is a weak solution of
		\eqref{compressible Euler conservative form} if, for every
		$\varphi\in C_c^\infty(\R\times\T^3)$ and every
		$\psi\in C_c^\infty(\R\times\T^3;\R^3)$, one has
		\begin{align}
		&\int_{\R}\int_{\T^3}\left(\rho\partial_t\varphi+m\cdot\nabla\varphi\right)\,\rd x\,\rd t=0,
		\label{weak continuity identity}\\
		&\int_{\R}\int_{\T^3}\left(m\cdot\partial_t\psi+\frac{m\otimes m}{\rho}:\nabla\psi+p(\rho)\Div\psi\right)\,\rd x\,\rd t=0.
		\label{weak momentum identity}
		\end{align}
	\end{df}

	In this paper, we construct nontrivial H{\"o}lder continuous weak solutions
	of \eqref{compressible Euler conservative form} on $\R\times\T^3$
	for smooth pressure laws with positive derivative.
	These solutions agree with a constant state outside a bounded time interval.

	Our construction builds on convex integration methods for the incompressible Euler equations.
	In his work on turbulent flows \cite{Onsager49}, Onsager conjectured that weak solutions to the incompressible Euler
	equations with H{\"o}lder regularity greater than $\frac{1}{3}$ must conserve kinetic energy,
	while solutions with lower regularity
	need not do so.
	Following earlier work of Eyink \cite{Eyink94}, Constantin, E, and Titi \cite{CET94} proved energy conservation in
	$L^3(0,T;B_{3,\infty}^{\sigma})$ for $\sigma>\frac{1}{3}$. In \cite{CCFS08},
	Cheskidov, Constantin, Friedlander, and Shvydkoy extended this result to
	 $L^3(0,T;B_{3,c(\mathbb N)}^{\frac{1}{3}})$.

	On the flexible side, De Lellis and Sz\'ekelyhidi \cite{DS09} formulated the equations as a differential inclusion
	and used convex integration to prove nonuniqueness
	of bounded weak solutions in $L^\infty(\R\times\R^d)$
	for $d\geq2$.
	See also earlier results by Scheffer \cite{Scheffer93}
	and Shnirelman \cite{Shnirelman97,Shnirelman00}.
	De Lellis and Sz\'ekelyhidi \cite{DS13} then constructed continuous dissipative solutions on $[0,1]\times\T^3$.
	The H{\"o}lder regularity was subsequently improved to $C^{\frac{1}{10}-}$ by De Lellis and Sz\'ekelyhidi \cite{DS14},
	and to $C^{\frac{1}{5}-}$ by Buckmaster, De Lellis, Isett, and Sz\'ekelyhidi \cite{BDIS15}.
	After a series of developments \cite{Buc15,BDS16,DS17,Isett16},
	Isett \cite{Isett18} resolved the flexible part of Onsager's conjecture by constructing nontrivial
	solutions in $C(\R;C^\alpha(\T^3))$ with compact temporal support for every $0<\alpha<\frac{1}{3}$,
	see \cite{BDSV19} for the admissible case. The two-dimensional case was resolved by Giri and Radu \cite{GR24}. They introduced a Newton--Nash iteration to separate perturbations in different directions in time and constructed nontrivial $C^\gamma(\R\times\T^2)$ weak solutions with compact temporal support for $0\leq\gamma<\frac{1}{3}$,
	see \cite{DuLiYe25} for the admissible case.

	Buckmaster and Vicol \cite{BV19} introduced intermittent convex integration to prove nonuniqueness of finite-energy weak solutions to the three-dimensional Navier--Stokes equations.
	For the Euler equations, Buckmaster, Masmoudi, Novack, and Vicol developed a spatially intermittent scheme and constructed nonconservative weak solutions in
	$C([0,T];H^s(\T^3))$ for $0<s<\frac{1}{2}$ \cite{BMNV21}.
	Novack and Vicol \cite{NV23} obtained such solutions in $C([0,T];B_{3,\infty}^{\beta}(\T^3))$ for $0<\beta<\frac{1}{3}$.
	In \cite{GKN26Helicity}, Giri, Kwon, and Novack constructed weak solutions with prescribed kinetic energy and generalized helicity.
	Bru\`e, Colombo, and Kumar \cite{BCK26} proved nonuniqueness of weak solutions to the Euler equations on $\T^2$ with velocity in $C([0,1];L^2(\T^2))$ and vorticity in $C([0,1];L^p(\T^2))$ for some $p>1$.

	Convex integration can also produce solutions satisfying the local energy inequality, a stronger requirement than nonincreasing total kinetic energy.
	In \cite{DK22}, De Lellis and Kwon constructed weak solutions in
	$C^\alpha([0,T]\times\T^3)$, $0\leq\alpha<\frac{1}{7}$, which satisfy the local energy inequality and strictly dissipate total kinetic energy.
	Giri, Kwon, and Novack \cite{GKN26Strong} subsequently constructed solutions in $C([0,T];B_{3,\infty}^{\beta}(\T^3))$, $0<\beta<\frac{1}{3}$, satisfying the local energy inequality and dissipating total kinetic energy.
	This proves the $L^3$-based strong Onsager conjecture using the wavelet-inspired framework developed in \cite{GKN24}.

	Convex integration has also been adapted to other fluid and elastic systems.
	For the surface quasi-geostrophic (SQG) equation on $\T^2$, Buckmaster, Shkoller, and Vicol \cite{BSV19}
	proved nonuniqueness of weak solutions by convex integration.
	Dai, Giri, and Radu \cite{DGR24} and Isett and Looi \cite{IL24SQG}
	independently resolved the flexible part of Onsager's conjecture by constructing weak solutions $\theta$
	that fail to conserve the Hamiltonian and satisfy $(-\Delta)^{-\frac{1}{2}}\theta\in C(\R;C^\gamma(\T^2))$ for $\frac{1}{2}\leq\gamma<1$.
	Further flexibility results for the SQG equation on $\T^2$ were obtained in \cite{BJN26}.
	Further applications of Newton--Nash iterations concern general active scalar equations \cite{Zhao25}.
	For the inviscid Euler--Boussinesq system, Onsager-type results were obtained on $\T^3$ \cite{MNY24}
	and subsequently on $\T^2$ \cite{Koley26}.
	
	For ideal magnetohydrodynamics, Beekie, Buckmaster, and Vicol \cite{BBV20} used intermittent convex integration to construct finite-energy weak solutions on $\T^3$ that do not conserve magnetic helicity.
	Miao, Nie, and Ye subsequently constructed weak solutions in
	$C^\alpha([0,T]\times\T^3)$ with total energy dissipation and nonconstant cross helicity for $0<\alpha<\frac{1}{3}$ \cite{MNY25MHD}.
	Using volume-preserving diffeomorphisms, Enciso, Pe\~nafiel-Tom\'as, and Peralta-Salas \cite{EPP26MHD} constructed H{\"o}lder continuous ideal MHD solutions with nonzero constant magnetic helicity and nonconstant total energy and cross helicity.
	Giardi and Sz\'ekelyhidi developed this approach in \cite{GS26MHD} and reached $C^{\frac{1}{3}-}(\R\times\T^3)$ in \cite{GS26Onsager}.
	
	Following earlier work on three-dimensional viscous and resistive MHD \cite{LZZ22MHD},
	Li, Zeng, and Zhang \cite{LZZ24MHD} proved sharp nonuniqueness in supercritical spaces beyond the Lions exponent $\frac{5}{4}$.
	After a series of developments \cite{LT20,BCV22}, Cheskidov and Luo \cite{CL22} proved nonuniqueness for the Navier--Stokes equations in
	$L^p(0,T;L^\infty(\T^d))$ for $d\geq2$ and $1\leq p<2$,
	showing the sharpness of the Ladyzhenskaya--Prodi--Serrin criterion at $(2,\infty)$.
	In two dimensions, they also proved nonuniqueness in $C([0,T];L^p(\T^2))$ for $1\leq p<2$ \cite{CL23NS}.
	Li, Qu, Zeng, and Zhang \cite{LQZZ24HyperNS} subsequently proved sharp nonuniqueness for
	the hyperdissipative Navier--Stokes equations on $\T^3$ in supercritical spaces beyond the Lions exponent.
	Temporal intermittency also yields nonuniqueness for hypo-viscous compressible Navier--Stokes \cite{LQZZ26}.
	Further applications concern geophysical models \cite{BMT26}, compressible Euler--Maxwell \cite{MQ23EM}, and nonlinear elastodynamics \cite{MQ25Lame,MQ25Null}.

	For the isentropic compressible Euler equations, bounded weak solutions away from vacuum
	satisfying the entropy inequality in the sense of distributions are called admissible or entropy solutions.
	For $p'>0$, the relative entropy method gives weak--strong uniqueness:
	an entropy solution coincides with a classical solution with the same initial data as long as the classical solution exists \cite{Dafermos16,Wie18}.
	For suitable initial data, however, De Lellis and Sz\'ekelyhidi \cite{DS10} and Chiodaroli \cite{Chiodaroli14}
	constructed infinitely many bounded entropy solutions by convex integration.
	For certain two-dimensional Riemann data, nonuniqueness holds when the corresponding one-dimensional self-similar solution contains shocks \cite{CDLK15,KM18}
	or consists of a single contact discontinuity \cite{KS25,Hor26}.
	For $p(\rho)=\rho^2$, Chiodaroli, Kreml, M\'acha, and Schwarzacher \cite{CKMS21} obtained nonunique admissible continuations from smooth data after shock formation.
	Nonuniqueness for dense sets of initial data was established by Chen, Vasseur, and Yu \cite{CVY21}
	and Chiodaroli and Feireisl \cite{CF24}.
	Boutros and Markfelder \cite{BM26} obtained global-in-time nonuniqueness under the local energy inequality,
	with wild data dense in $L^r(\T^3)$ for $1\leq r<\infty$.

	The preceding results do not address nonuniqueness in H{\"o}lder spaces.
	Giri and Kwon \cite{GK22} constructed infinitely many entropy solutions satisfying the local energy equality,
	with the same initial data, a common smooth density, and momentum in $C^\alpha([0,T]\times\T^3)$ for
	$0\leq\alpha<\frac{1}{7}$.
	Their results allow arbitrary smooth initial densities away from vacuum and include examples with strict local energy dissipation.
	Subsequently, Giri and Koley \cite{GKo26} proved nonuniqueness for the inhomogeneous incompressible Euler equations
	on $[0,T]\times\T^3$ with positive nonconstant density and both density and velocity in $C^\alpha$, $0\leq\alpha<\frac{1}{7}$.

	For the compressible Euler equations, we construct weak solutions with higher momentum regularity, without imposing the entropy inequality.
	Our Newton--Nash iteration perturbs both density and momentum,
	yielding momentum and velocity with any H{\"o}lder exponent below $\frac{1}{3}$
	and density with any exponent below $\frac{2}{3}$.
	More precisely, we prove the following result.
	\begin{thm}
		\label{thm:main}
		Let $0\leq\gamma_m<\frac13$ and $0\leq\gamma_\rho<\frac23$, and let
		$\rho^*>0$ be a constant reference density. Assume that
		$p\in C^\infty((0,\infty))$ is a fixed pressure law satisfying
		$
		p'(\rho)>0
		\text{ for all } \rho\in\left[\frac{\rho^*}{2},\frac{3\rho^*}{2}\right].
		$
		Then there exists a nontrivial weak solution $(\rho,m)$ of the isentropic compressible Euler system on $\R\times\T^3$
		\begin{equation*}
			\left\{
			\begin{aligned}
				&\partial_t\rho+\Div m=0,\\
				&\partial_t m+\Div\left(\frac{m\otimes m}{\rho}\right)+\nabla p(\rho)=0,
			\end{aligned}
			\right.
		\end{equation*}
		such that
		$$
		\rho-\rho^*\in C^{\gamma_\rho}(\R_t\times\T^3_x),
		\qquad
		m\in C^{\gamma_m}(\R_t\times\T^3_x),
		$$
		and
		$$
		\rho(t,x)\geq\frac{\rho^*}{2}
		\qquad
		\text{for all }(t,x)\in\R\times\T^3.
		$$
		Moreover, its deviation from the constant state $(\rho^*,0)$ is compactly
		supported in time, in the sense that
		$$
		(\rho,m)=(\rho^*,0)
		\qquad
		\text{for }t\notin[-2,2].
		$$
	\end{thm}

	In particular, for every $0\leq\gamma<1/3$, one may take
	$\gamma_m=\gamma$ and $\gamma_\rho=2\gamma$. Thus the density may be
	taken in $C^{2\gamma}$ while the momentum is in $C^\gamma$. Since the density
	stays bounded away from vacuum, the corresponding velocity $u=m/\rho$ also belongs
	to $C^\gamma$.
	By time translation, Theorem~\ref{thm:main} also gives nonuniqueness
	for the constant initial data $(\rho^*,0)$ in the stated H{\"o}lder classes.

	\begin{rem}
	For each fixed $\frac{1}{3}<\gamma_\rho<\frac{2}{3}$, the momentum range
	$\gamma_m<\frac{1}{3}$ in Theorem~\ref{thm:main} is sharp up to the endpoint
	for nontrivial solutions agreeing with $(\rho^*,0)$ outside a bounded time interval.
	Indeed, if $\gamma_m>\frac13$, then $\rho$, $m$, and $u=m/\rho$ have
	space--time H{\"o}lder regularity above $\frac13$, since $\rho$ stays away from vacuum.
	Thus the energy criterion in \cite{FGSW17,GMS18} gives local energy equality.
	Fix $t_0$ before the temporal support of $(\rho-\rho^*,m)$.
	Since $p'>0$ near $\rho^*$, weak--strong uniqueness gives agreement
	with $(\rho^*,0)$ on a short interval starting at $t_0$.
	At any finite right endpoint, continuity allows the same argument
	to be repeated.
	Hence $(\rho,m)\equiv(\rho^*,0)$ on $\R\times\T^3$,
	contrary to the assumption that $(\rho-\rho^*,m)\not\equiv(0,0)$.
	\end{rem}

	\subsection{Main difficulties and ideas of the proof}
	\label{subsec:proof-ideas}

	In a convex integration scheme, we construct a sequence of smooth approximate solutions converging to a weak solution.
	Here the approximate solutions satisfy the following compressible Euler--Reynolds system.

	\begin{df}
	A triple of smooth fields $(\rho_q,m_q,R_q)$ on $\R\times\T^3$
	with positive density is called a solution of the compressible
	Euler--Reynolds system if
	\begin{equation}
		\left\{
		\begin{aligned}
			&\partial_t\rho_q+\Div m_q=0,\\
			&\partial_t m_q+\Div\left(\frac{m_q\otimes m_q}{\rho_q}\right)+\nabla p(\rho_q)=\Div R_q,
		\end{aligned}
		\right.
		\label{eq:compressible-Euler-Reynolds}
	\end{equation}
	where the Reynolds stress $R_q$ is a possibly nonsymmetric $3\times3$ matrix field.
	\end{df}

	At stage $q$, the frequency and amplitude parameters are
	$\lambda_q=\lceil a^{b^q}\rceil$ and $\delta_q=\lambda_q^{-2\beta}$.
	We take $0<\beta<1/3$, $b>1$ close to $1$, and $a>1$ sufficiently large.
	We choose the perturbations so that $R_q\to0$ uniformly and
	$(\rho_q,m_q)$ converges in the stated H{\"o}lder spaces.

	\subsubsection{Gluing and the Newton--Nash method}

	In earlier Beltrami-based constructions for incompressible Euler, short time cutoffs
	led to transport errors that limited the attainable H{\"o}lder regularity.
	In three dimensions, the Mikado flows of Daneri--Sz\'ekelyhidi~\cite{DS17}
	in different directions can be supported in pairwise disjoint tubes,
	which remain disjoint under transport by a common flow map.
	Isett~\cite{Isett18} combined Mikado flows with a gluing method to allow longer transport intervals and a larger time cutoff scale,
	reducing the transport error and reaching the Onsager range.
	The gluing construction joins exact Euler solutions and localizes the Reynolds stress to separated time intervals.
	This step uses short-time existence and H{\"o}lder stability estimates for exact Euler solutions~\cite{BDSV19}.
	In Isett's construction, the Mikado perturbations associated with distinct intervals have disjoint temporal supports.

	In two dimensions, Mikado flows cannot achieve the required non-interaction between different directions,
	since tubes in nonparallel directions necessarily intersect.
	Temporal separation provides an alternative. Cheskidov and Luo introduced
	temporally oscillatory and intermittent perturbations for the Navier--Stokes equations~\cite{CL22}
	and used them to separate different directions in time~\cite{CL23NS}.
	Part of the stress is canceled by the time derivative of a temporal corrector,
	and the rest by the quadratic self-interaction of the Nash perturbation.

	For two-dimensional Euler, Giri and Radu~\cite{GR24} reach the Onsager range
	by replacing the temporal corrector with a solution of the linearized Euler equations.
	All terms of the linearized operator then contribute to stress cancellation.
	The principal Nash perturbations for different directions have disjoint temporal supports,
	so the corresponding quadratic cross terms vanish.
	The forcing terms in the Newton equation contain zero-mean periodic factors.
	Their frequency is chosen much larger than the material derivative scale of the coefficients.
	Integration by parts along the background flow then gives a factor of the inverse frequency in the Newton perturbation.

	For compressible Euler, gluing exact solutions would require H{\"o}lder stability estimates
	for both density and momentum, compatible with the frequency and time scales of the iteration.
	Smooth local existence for the hyperbolic system does not by itself provide these estimates.
	We instead construct localized Newton corrections around a mollified background.
	The coupled density--momentum correction is obtained from a scalar wave equation
	with locally frozen sound speed in Lagrangian coordinates.
	The wave source does not contain the rapidly oscillating time factors used in \cite{GR24}.
	The Newton iteration cancels the Reynolds stress up to smaller errors and a symmetric gluing stress
	generated by the time cutoffs. The gluing stress is canceled in the Nash step.

	\subsubsection{The coupled Newton correction}

	At stage $q$, let $(\rho_{q,j},m_{q,j})$ denote the state after $j$ Newton
	substeps, with $(\rho_{q,0},m_{q,0})=(\rho_q,m_q)$.
	Writing its Reynolds stress as $R_{q,j}+S_{q,j}+P_{q+1,j}$, we have
	$$
	\left\{
	\begin{aligned}
	&\partial_t\rho_{q,j}+\Div m_{q,j}=0,\\
	&\partial_t m_{q,j}+\Div\left(\frac{m_{q,j}\otimes m_{q,j}}{\rho_{q,j}}\right)+\nabla p(\rho_{q,j})
	=\Div R_{q,j}+\Div S_{q,j}+\Div P_{q+1,j}.
	\end{aligned}
	\right.
	$$
	Here $R_{q,j}$ is the interior stress.
	The tensor $S_{q,j}$ is the sum of the symmetric gluing stresses already produced, and $P_{q+1,j}$ contains the flux remainders and mollification errors.
	The number of Newton substeps $J_*$ is fixed independently of $q$ and $a$.
	The system holds for $0\leq j\leq J_*$; for $j<J_*$, we correct $R_{q,j}$.
	Initially, $R_{q,0}$ is the mollified stress, $S_{q,0}=0$, and $P_{q+1,0}=R_q-R_{q,0}$.
	
	Fix the background $(\orho_q,\om_q)$ mollified in space at length
	$\ell_q=\lambda_q^{-\frac34}\lambda_{q+1}^{-\frac14}$, and set
	$
	\overline u_q=\frac{\om_q}{\orho_q},\
	\DTL=\partial_t+\overline u_q\cdot\nabla.
	$
	At each Newton substep, let $Q[R_{q,j}]$ and $W[R_{q,j}]$ denote the zero-mean scalar and vector Hodge potentials.
	They give the decomposition
	$$
	\Div R_{q,j}=\nabla Q[R_{q,j}]+\curl W[R_{q,j}],\qquad
	\Div W[R_{q,j}]=0.
	$$
	The potentials are defined by
	$$
	Q[R_{q,j}]:=\Delta^{-1}\Div\Div R_{q,j},\qquad
	W[R_{q,j}]:=-\Delta^{-1}\curl\Div R_{q,j}.
	$$
	Since the pressure is prescribed by the density, the gradient part cannot be absorbed into an independent pressure.
	We call the corrections of the gradient and curl parts the compressible and incompressible perturbations, respectively.
	Both perturbations are constructed in Section~\ref{sec:Newton steps}.

	We first describe the compressible perturbation.
	We define the density and momentum corrections through a vector potential $V$ by setting
	$$
	\begin{aligned}
	\varrho[V]=\Div V,\quad
	w[V]=-\partial_tV
	+\Div(\overline u_q\otimes V-V\otimes\overline u_q).
	\end{aligned}
	$$
	The continuity equation then holds for every smooth $V$.
	For the gradient part, we choose $(\varrho,w)=(\varrho[V],w[V])$ so that
	$$
	\begin{aligned}
	\partial_tw
	+\Div\bigl(\overline u_q\otimes w+w\otimes\overline u_q
	-\varrho\overline u_q\otimes\overline u_q\bigr)
	+\nabla\bigl(p'(\orho_q)\varrho\bigr)
	\approx-\nabla Q[R_{q,j}].
	\end{aligned}
	$$
	We construct $V$ from a scalar wave equation along
	the background flow.

	Let $\alpha>0$ be a small exponent. The time scale has the form
	$\tau_q\sim\lambda_{q+1}^{-4\alpha}(\lambda_q\delta_q^{\frac12})^{-1}$.
	Let $X_u$ and $\Phi_u$ be the forward and backward flows of $\overline u_q$
	based at the center $t_u$ of a time interval of length comparable to $\tau_q$.
	Set
	$
	J_u=\det\nabla\Phi_u,\ \mathbf C_u=\cof\nabla\Phi_u.
	$
	For the local construction, we extend the background fields periodically to $\R^3$ and use the corresponding flow maps.
	On a cube of side length $O(\tau_q)$ centered at $x_\kappa$, we freeze the squared sound speed at
	$c_{q,u,\kappa}^2=p'(\orho_q)(t_u,x_\kappa)$.
	Choose time and spatial cutoffs $\theta_u$ and $\chi_\kappa$ adapted to these intervals and cubes, respectively, with $\sum_u\theta_u^2=1$ and $\sum_\kappa\chi_\kappa^2=1$.

	To supply the source derivatives required by the wave estimates, we use the along-flow mollification
	from Section~\ref{Mollification}, as in \cite{GK22,GR24}.
	The temporal mollification scale is
	$\ell_{t,q}=\delta_q^{-\frac12}\lambda_q^{-\frac12}\lambda_{q+1}^{-\frac12}$.
	Let $\eta_{\ell_{t,q}}^{\mathrm t}$ be the temporal mollifier at this scale.
	Write $X_{\overline u_q}(t+s,x;t)$ for the flow of $\overline u_q$ from $(t,x)$ to time $t+s$.
	For a smooth field $f$, set
	$$
	\bigl(\eta_{\ell_{t,q}}^{\mathrm t}\ast_{X_{\overline u_q}}f\bigr)(t,x)
	:=\int_{\R}f(t+s,X_{\overline u_q}(t+s,x;t))
	\eta_{\ell_{t,q}}^{\mathrm t}(s)\,\rd s.
	$$
	At each Newton substep, we localize and mollify the scalar Hodge potential of $R_{q,j}$ by setting
	$$
	G_{q,j,u,\kappa}^{\mathrm{mol}}
	=\theta_u\chi_\kappa^2(\Phi_u)
	\bigl(\eta_{\ell_{t,q}}^{\mathrm t}\ast_{X_{\overline u_q}}
	(J_u^{-1}Q[R_{q,j}])\bigr),
	$$
	and mollify the vector Hodge potential by
	$$
	W_{q,j}^{\mathrm{mol}}
	=\eta_{\ell_{t,q}}^{\mathrm t}\ast_{X_{\overline u_q}}W[R_{q,j}].
	$$

	For the scalar source, we solve the wave equation in Lagrangian coordinates:
	$$
	(\partial_t^2-c_{q,u,\kappa}^2\Delta_y)
	\widehat I_{q,u,\kappa}[G_{q,j,u,\kappa}^{\mathrm{mol}}](t,y)
	=G_{q,j,u,\kappa}^{\mathrm{mol}}(t,X_u(t,y))
	$$
	with zero initial data at the preceding time center $t=t_{u-1}$. The cutoff $\theta_u$ makes $G_{q,j,u,\kappa}^{\mathrm{mol}}$ vanish near this time.
	For a smooth scalar source $H$ vanishing near $t_{u-1}$, Kirchhoff's formula gives
	$$
	\begin{aligned}
	\widehat I_{q,u,\kappa}[H](t,y)
	=\int_{t_{u-1}}^t\frac{t-s}{4\pi}
	\int_{\mathbb S^2}
	H\left(s,X_u\left(s,y+c_{q,u,\kappa}(t-s)\xi\right)\right)
	\,\rd S_\xi\,\rd s.
	\end{aligned}
	$$

	Write $Y_y=X_u(s,y+c_{q,u,\kappa}(t-s)\xi)$.
	Integration by parts gives
	$$
	\begin{aligned}
	\nabla_y\widehat I_{q,u,\kappa}[H](t,y)
	=\frac{1}{4\pi c_{q,u,\kappa}}
	\int_{t_{u-1}}^t\int_{\mathbb S^2}\xi
	\left[H(s,Y_y)+(t-s)(\DTL H)(s,Y_y)\right]\,\rd S_\xi\,\rd s,
	\end{aligned}
	$$
	and since $\DTL\Phi_u=0$ and $H$ vanishes near $t_{u-1}$, for every $M\in\N_0$ we have
	$$
	\DTL^M\left[\widehat I_{q,u,\kappa}[H](t,\Phi_u)\right]
	=\widehat I_{q,u,\kappa}[\DTL^M H](t,\Phi_u).
	$$
	The gradient formula applied to $\DTL^M H$ therefore uses source material derivatives
	through order $M+1$, without increasing the spatial order.

	Applying the Piola transform to the gradient of the wave solution, we set
	$$
	\widetilde V_{q,j+1,u,\kappa}
	=\mathbf C_u^T\nabla_y
	\widehat I_{q,u,\kappa}[G_{q,j,u,\kappa}^{\mathrm{mol}}](t,\Phi_u).
	$$
	Using the wave equation, we obtain
	$$
	\begin{aligned}
	-\DTL^2\widetilde V_{q,j+1,u,\kappa}
	+\nabla\bigl(p'(\orho_q)\varrho[\widetilde V_{q,j+1,u,\kappa}]\bigr)
	+\nabla(J_uG_{q,j,u,\kappa}^{\mathrm{mol}})=\nabla\bigl((p'(\orho_q)-c_{q,u,\kappa}^2)
	\varrho[\widetilde V_{q,j+1,u,\kappa}]\bigr)
	+\mathcal E_{q,j,u,\kappa}^{\sharp}.
	\end{aligned}
	$$
	Here
	$$
	\begin{aligned}
	\mathcal E_{q,j,u,\kappa}^{\sharp}
	:=&(\nabla J_u)\widehat I_{q,u,\kappa}[\DTL^2G_{q,j,u,\kappa}^{\mathrm{mol}}](t,\Phi_u)
	+(J_u\nabla\Phi_u-\mathbf C_u)^T
	\nabla_y\widehat I_{q,u,\kappa}[\DTL^2G_{q,j,u,\kappa}^{\mathrm{mol}}](t,\Phi_u)\\
	&-2(\DTL\mathbf C_u)^T
	\nabla_y\widehat I_{q,u,\kappa}[\DTL G_{q,j,u,\kappa}^{\mathrm{mol}}](t,\Phi_u)
	-(\DTL^2\mathbf C_u)^T
	\nabla_y\widehat I_{q,u,\kappa}[G_{q,j,u,\kappa}^{\mathrm{mol}}](t,\Phi_u).
	\end{aligned}
	$$
	By finite propagation, for $t\in\supp\theta_u$, the spatial support of
	$\widetilde V_{q,j+1,u,\kappa}(t,\cdot)$ lies in the image under $X_u(t,\cdot)$
	of a ball centered at $x_\kappa$ with radius $O(\tau_q)$.

	To apply a local right inverse of divergence, each momentum residual must have zero spatial integral.
	Choose a profile $\phi_{u,\kappa}$ with unit spatial integral, supported in the same transported region and satisfying
	$$
	\partial_t\phi_{u,\kappa}+\Div(\overline u_q\phi_{u,\kappa})=0.
	$$
	We enforce the zero-integral condition by correcting the mean of $\widetilde V_{q,j+1,u,\kappa}$ before inserting the time cutoff:
	$$
	\begin{aligned}
	V_{q,j+1,u,\kappa}^{\circ}
	=\widetilde V_{q,j+1,u,\kappa}
	-\left(\int_{\R^3}\widetilde V_{q,j+1,u,\kappa}(t,x)\,\rd x\right)
	\phi_{u,\kappa},\quad
	\int_{\R^3}\phi_{u,\kappa}(t,x)\,\rd x=1,\qquad
	\int_{\R^3}V_{q,j+1,u,\kappa}^{\circ}(t,x)\,\rd x=0.
	\end{aligned}
	$$
	Write $\overline V_{q,j+1,u,\kappa}(t):=\int_{\R^3}\widetilde V_{q,j+1,u,\kappa}(t,x)\,\rd x$.
	Set $V_{q,j+1,u,\kappa}\comp=\theta_uV_{q,j+1,u,\kappa}^{\circ}$ and define the local density and momentum perturbations by
	$(\varrho,w):=(\varrho[V_{q,j+1,u,\kappa}\comp],w[V_{q,j+1,u,\kappa}\comp])$.

	Their linearized momentum equation takes the form
	$$
	\begin{aligned}
	\partial_tw+\Div\bigl(\overline u_q\otimes w+w\otimes\overline u_q
	-\varrho\overline u_q\otimes\overline u_q\bigr)
	+\nabla\bigl(p'(\orho_q)\varrho\bigr)=-\theta_u\nabla(J_uG_{q,j,u,\kappa}^{\mathrm{mol}})
	+F_{q,j+1,u,\kappa}\comp+F_{q,j+1,u,\kappa}^{\mathrm{cut}},
	\end{aligned}
	$$
	where the interior and cutoff residuals are
	$$
	\begin{aligned}
	F_{q,j+1,u,\kappa}\comp
	:=&\theta_u\nabla\Bigl[
	\bigl(p'(\orho_q)-c_{q,u,\kappa}^2\bigr)\varrho[\widetilde V_{q,j+1,u,\kappa}]
	-p'(\orho_q)\varrho[\overline V_{q,j+1,u,\kappa}\phi_{u,\kappa}]\Bigr]
	+\theta_u\Div\bigl((\DTL\overline u_q)\otimes V_{q,j+1,u,\kappa}^{\circ}\bigr)\\
	&+\theta_u\mathcal E_{q,j,u,\kappa}^{\sharp}
	-\theta_u\DTL\bigl((\Div\overline u_q)\widetilde V_{q,j+1,u,\kappa}\bigr)
	-\theta_u(\Div\overline u_q)\DTL\widetilde V_{q,j+1,u,\kappa}
	-\theta_u(\Div\overline u_q)^2\widetilde V_{q,j+1,u,\kappa}
	+\theta_u(\partial_t^2\overline V_{q,j+1,u,\kappa})\phi_{u,\kappa},\\
	F_{q,j+1,u,\kappa}^{\mathrm{cut}}
	=&-2(\partial_t\theta_u)
	\bigl(\DTL V_{q,j+1,u,\kappa}^{\circ}
	+(\Div\overline u_q)V_{q,j+1,u,\kappa}^{\circ}\bigr)
	-(\partial_t^2\theta_u)V_{q,j+1,u,\kappa}^{\circ}.
	\end{aligned}
	$$

	To incorporate these residuals into the Reynolds stress, we write them as divergences of matrix fields.
	After separating the explicit gradient and divergence terms, we use a local right inverse of divergence for the remaining vector terms.

	\subsubsection{A local right inverse of divergence}

	The local inverse constructed in Section~\ref{sec:local-inverse-divergence} preserves support in the image of a fixed ball under the flow.
	Its spatial and material derivative estimates gain a factor equal to the radius of that ball.
	We adapt the construction of Isett--Oh~\cite{Isett16}
	to the nonsymmetric divergence equation, see also \cite{IMOT25}.
	
	For a smooth velocity $\upsilon$, write $D_{t,\upsilon}:=\partial_t+\upsilon\cdot\nabla$.
	Let $X_\upsilon,\Phi_\upsilon$ be the forward and backward flow maps of
	$\upsilon$ based at $t_0$, and set
	$
	\JX=\det\nabla X_\upsilon,\ \JPhi=\det\nabla\Phi_\upsilon,\
	P=\JPhi(\nabla\Phi_\upsilon)^{-1}.
	$
	
	Let $U$ be a smooth vector field satisfying $\int U(t,x)\,\rd x=0$.
	Assume that its spatial support at each time $t$ under consideration is contained in $X_\upsilon(t,B(x_0,\bar\varepsilon))$, with $0<\bar\varepsilon\leq1$.
	Let $R^*$ be a linear local right inverse of divergence.
	It is constructed with a fixed smooth averaging kernel of integral one, supported in $B(x_0,\bar\varepsilon)$.
	Pull $U$ back to the fixed ball $B(x_0,\bar\varepsilon)$ and define
	$$
	V=\JX(U\circ X_\upsilon),\qquad
	R_{ij}^{(\upsilon)}[U]
	=P_{jl}\bigl(R^*_{il}[V]\circ\Phi_\upsilon\bigr).
	$$
	For every smooth zero-mean vector field $V$ supported in this ball, $R^*[V]$ is a smooth matrix field satisfying
	$$
	\partial_{y_l}R^*_{il}[V]=V_i,\qquad
	\supp R^*[V]\subset B(x_0,\bar\varepsilon).
	$$
	For $N\in\N_0$, we have
	$$
	\nrm{R^*[V]}_N
	\lesssim_N\bar\varepsilon
	\sum_{r=0}^N\bar\varepsilon^{-(N-r)}\nrm{V}_r.
	$$

	The Jacobian factor in the definition of $V$ preserves the zero-integral condition even when the flow does not preserve volume.
	The fields $\partial_t^kV$ therefore have zero spatial integral
	and remain supported in the fixed ball.
	The kernel is independent of time, so $\partial_t^kR^*[V]=R^*[\partial_t^kV]$ for every $k\in\N_0$.
	The Piola transform gives $\Div R^{(\upsilon)}[U]=U$.
	Since $R^*$ commutes with $\partial_t$ and $D_{t,\upsilon}\Phi_\upsilon=0$,
	$$
	D_{t,\upsilon}^MR_{ij}^{(\upsilon)}[U]
	=\sum_{k=0}^M\binom Mk
	(D_{t,\upsilon}^{M-k}P_{jl})
	\bigl(R^*_{il}[\partial_t^kV]\circ\Phi_\upsilon\bigr).
	$$
	The resulting spatial and material derivative estimates in $C^N$ include the factor $\bar\varepsilon$.
	For the Newton step, taking $\bar\varepsilon\sim\tau_q$
	gives the factor $\tau_q(\tau_q\ell_q^{-1})^\alpha$ in the H{\"o}lder estimates.

	We apply this operator to the interior and cutoff residuals of the compressible perturbation.
	The mean correction and finite propagation ensure the required zero-integral and support conditions.
	Writing this operator as $R_{u,\kappa}^{(\overline u_q)}$, we collect the explicit gradient and divergence terms in the matrix
	$$
	\begin{aligned}
	\mathcal S_{q,j+1,u,\kappa}
	:=&\theta_u\Bigl[
	\bigl(p'(\orho_q)-c_{q,u,\kappa}^2\bigr)\varrho[\widetilde V_{q,j+1,u,\kappa}]
	-p'(\orho_q)\varrho[\overline V_{q,j+1,u,\kappa}\phi_{u,\kappa}]\Bigr]\Id
	+\theta_u(\DTL\overline u_q)\otimes V_{q,j+1,u,\kappa}^{\circ}.
	\end{aligned}
	$$
	The remaining vector field is
	$$
	\begin{aligned}
	H_{q,j+1,u,\kappa}^{\mathrm{int}}
	:=\theta_u\Bigl[\mathcal E_{q,j,u,\kappa}^{\sharp}
	-\DTL\bigl((\Div\overline u_q)\widetilde V_{q,j+1,u,\kappa}\bigr)
	-(\Div\overline u_q)\DTL\widetilde V_{q,j+1,u,\kappa}
	-(\Div\overline u_q)^2\widetilde V_{q,j+1,u,\kappa}
	+(\partial_t^2\overline V_{q,j+1,u,\kappa})\phi_{u,\kappa}\Bigr].
	\end{aligned}
	$$
	We then write
	$$
	\begin{aligned}
	F_{q,j+1,u,\kappa}\comp
	&=\Div\mathcal S_{q,j+1,u,\kappa}
	+H_{q,j+1,u,\kappa}^{\mathrm{int}},&&
	R_{q,j+1}\comp
	=\sum_{u,\kappa}\left(
	\mathcal S_{q,j+1,u,\kappa}
	+R_{u,\kappa}^{(\overline u_q)}
	[H_{q,j+1,u,\kappa}^{\mathrm{int}}]\right).
	\end{aligned}
	$$

	For the cutoff residual, we use the same local inverse and set
	$$
	S_{q,j+1,u}^{\mathrm{cut}}
	=\sum_{\kappa\in\mathbb Z^3}R_{u,\kappa}^{(\overline u_q)}
	[F_{q,j+1,u,\kappa}^{\mathrm{cut}}],\qquad
	\Div S_{q,j+1,u}^{\mathrm{cut}}
	=\sum_{\kappa\in\mathbb Z^3}F_{q,j+1,u,\kappa}^{\mathrm{cut}}.
	$$
	The interior stress $R_{q,j+1}\comp$ and the cutoff stress $S_{q,j+1,u}^{\mathrm{cut}}$ need not be symmetric.
	The incompressible correction below cancels the divergence of the antisymmetric part of $S_{q,j+1,u}^{\mathrm{cut}}$, leaving a symmetric gluing stress for the Nash step.

	\subsubsection{The incompressible part of the Newton correction}

	The incompressible perturbation changes only the momentum and corrects the curl part of the stress divergence.
	We decompose the cutoff stress as $S_{q,j+1,u}^{\mathrm{cut}}
	=S_{q,j+1,u}^{\mathrm{cut,sym}}+S_{q,j+1,u}^{\mathrm{cut,skew}},$ where
	$$
	\begin{aligned}
	S_{q,j+1,u}^{\mathrm{cut,sym}}
	&:=\frac{1}{2}\bigl(S_{q,j+1,u}^{\mathrm{cut}}+(S_{q,j+1,u}^{\mathrm{cut}})^T\bigr),&&
	S_{q,j+1,u}^{\mathrm{cut,skew}}
	:=\frac{1}{2}\bigl(S_{q,j+1,u}^{\mathrm{cut}}-(S_{q,j+1,u}^{\mathrm{cut}})^T\bigr).
	\end{aligned}
	$$
	For an antisymmetric periodic matrix $A$, one has
	$$
	\Div\Div A=0,\qquad \Div A=\curl W[A].
	$$
	We therefore combine the divergence of the antisymmetric part of the cutoff stress with the curl forcing
	in a divergence-free momentum correction. The density is unchanged.

	For each time cell, we solve
	$$
	\left\{
	\begin{aligned}
	&\partial_t\omega_{q,j,u}
	+\Div(\overline u_q\otimes\omega_{q,j,u}
	+\omega_{q,j,u}\otimes\overline u_q)+\nabla\pi_{q,j,u}=-\theta_u^2\curl W_{q,j}^{\mathrm{mol}}
	-\curl W[S_{q,j+1,u}^{\mathrm{cut,skew}}],\\
	&\Div\omega_{q,j,u}=0,\qquad \omega_{q,j,u}|_{t=t_u}=0.
	\end{aligned}
	\right.
	$$
	We use the enlarged temporal cutoff $\tilde\theta_u$, which equals one on a neighborhood of
	$\supp\theta_u\cup\supp_t S_{q,j+1,u}^{\mathrm{cut}}$.
	Multiplication by this cutoff defines the divergence-free momentum correction and the corresponding pressure:
	$$
	\tm_{q,j+1}\inc=\sum_u\tilde\theta_u\omega_{q,j,u},\qquad
	\pi_{q,j+1}\inc=\sum_u\tilde\theta_u\pi_{q,j,u}.
	$$
	Using $\sum_u\theta_u^2=1$, we obtain
	$$
	\left\{
	\begin{aligned}
	&\partial_t\tm_{q,j+1}\inc
	+\Div\bigl(\overline u_q\otimes\tm_{q,j+1}\inc
	+\tm_{q,j+1}\inc\otimes\overline u_q\bigr)
	+\nabla\pi_{q,j+1}\inc=-\curl W_{q,j}^{\mathrm{mol}}
	-\sum_u\curl W[S_{q,j+1,u}^{\mathrm{cut,skew}}]
	+\sum_u(\partial_t\tilde\theta_u)\omega_{q,j,u},\\
	&\Div\tm_{q,j+1}\inc=0.
	\end{aligned}
	\right.
	$$
	We include the pressure $\pi_{q,j+1}\inc$ in the interior stress by setting
	$
	R_{q,j+1}=R_{q,j+1}\comp-\pi_{q,j+1}\inc\Id.
	$

	The curl forcing and zero initial data also give zero spatial mean for
	$\omega_{q,j,u}$.
	Set $\mathcal B=-\Delta^{-1}\curl$.
	With $\Psi_{q,j,u}=\mathcal B\omega_{q,j,u}$ and
	$\curl\Psi_{q,j,u}=\omega_{q,j,u}$, we have
	$
	(\partial_t\tilde\theta_u)\omega_{q,j,u}
	=\curl\bigl((\partial_t\tilde\theta_u)\Psi_{q,j,u}\bigr),
	$
	since $\tilde\theta_u$ depends only on time.
	Let $\mathcal R$ denote the non-local periodic symmetric inverse-divergence operator.
	The remaining gluing stress is therefore
	$$
	R_{q,j+1}^{\mathrm{gluing}}
	=\sum_uS_{q,j+1,u}^{\mathrm{cut,sym}}
	+\mathcal R\curl\left(\sum_u(\partial_t\tilde\theta_u)\Psi_{q,j,u}\right).
	$$
	Both terms are symmetric, as required for cancellation by the quadratic self-interaction of the Nash perturbation.
	Their temporal supports lie in the separated transition intervals.
	The two corrections update the state within the same Newton substep:
	$$
	\begin{aligned}
	\rho_{q,j+1}
	=\rho_{q,j}+\sum_{u,\kappa}\varrho[V_{q,j+1,u,\kappa}\comp],\quad
	m_{q,j+1}
	=m_{q,j}+\sum_{u,\kappa}w[V_{q,j+1,u,\kappa}\comp]
	+\tm_{q,j+1}\inc.
	\end{aligned}
	$$

	\subsubsection{Higher material derivative estimates}
	The interior stress from one Newton substep is the source for the next, so its material derivatives must also be controlled.
	At material order $M$, the coupled correction requires source derivatives through order $M+2$,
	and its interior residual requires order $M+3$. The Nash transport error similarly requires $M+1$ material derivatives of the perturbation coefficients.
	The along-flow mollification supplies these additional source derivatives.
	Each Newton substep preserves the range of material orders and lowers the maximum total derivative order by three.

	The spatial mollification supplies extra spatial derivatives but does not commute with material differentiation.
	For $P_{\leq\ell^{-1}}f=\eta_\ell*f$, we have
	$$
	[D_{t,\upsilon},P_{\leq\ell^{-1}}]f(x)
	=\int_{\R^3}\eta_\ell(h)
	\bigl(\upsilon(x)-\upsilon(x-h)\bigr)\cdot\nabla f(x-h)\,\rd h.
	$$
	The velocity difference gives a factor of $\ell$, as in the commutator estimates of Giri--Kwon~\cite{GK22}.
	For higher material derivatives, we estimate the iterated commutators using their kernel representations.
	Let $N_{\mathrm{ind}}$ denote the maximum total derivative order for density and momentum in the inductive assumptions of Section~\ref{sec:main-iteration}.
	We choose a mollifier satisfying
	$$
	\int_{\R^3}\eta(h)\,\rd h=1,\qquad
	\int_{\R^3}h^\sigma\eta(h)\,\rd h=0
	\quad(1\leq|\sigma|\leq N_{\mathrm{ind}}).
	$$
	These moments cancel low-order Taylor terms and give higher powers of $\ell$ in the mollification error,
	also after material differentiation.
	These mollification estimates are established in Appendix~\ref{app:mollification-estimates} and used in Section~\ref{Mollification}.

	For the Hodge potentials and the periodic operators $\mathcal R\curl$ and $\mathcal B\curl$,
	we use the div--curl argument of Glass, Sueur, and Takahashi~\cite{GST12}.
	The estimates are proved in Section~\ref{sec:periodic-Hodge-and-Newton-stresses}.
	We differentiate the Hodge decomposition and the divergence constraint,
	then apply the periodic div--curl estimate at each fixed time.
	For a smooth periodic vector field $Y$, this estimate reads
	$$
	\nrm{Y}_{N+1+\alpha}
	\lesssim_{N,\alpha}\nrm{\Div Y}_{N+\alpha}+\nrm{\curl Y}_{N+\alpha}
	+\left|\fint_{\T^3}Y\,\rd x\right|,
	$$
	for $N\in\mathbb N_0$ and $0<\alpha<1$.
	The identities for iterated material derivatives in \cite{GST12}
	express the commutators in terms of lower material derivatives.
	For a smooth field $h$, we also use the commutator and mean identities for a compressible background flow:
	$$
	\begin{aligned}
	[D_{t,\upsilon},\partial_i]h
	=-(\partial_i\upsilon_j)\partial_jh,\quad
	\fint_{\T^3}D_{t,\upsilon}h\,\rd x
	=\frac{\rd}{\rd t}\fint_{\T^3}h\,\rd x
	-\fint_{\T^3}h\Div\upsilon\,\rd x.
	\end{aligned}
	$$
	In particular, the zero mean of a Hodge potential does not imply
	zero mean for its material derivatives.
	We estimate the commutator products at a finite sequence of decreasing H{\"o}lder exponents
	and control the potentials and their means by induction on the material order.

	The same div--curl induction gives higher material derivative estimates for the incompressible perturbation.
	Short-time transport estimates first control $\omega_{q,j,u}$ and its vector potential $\Psi_{q,j,u}$
	at material order zero.
	The pressure $\pi_{q,j,u}$ is recovered from $\Psi_{q,j,u}$ by
	$$
	\pi_{q,j,u}
	=-2\partial_i\partial_k\Delta^{-1}\bigl(\Psi_{q,j,u}\times\nabla(\overline u_q)_i\bigr)_k,
	\qquad \fint_{\T^3}\pi_{q,j,u}\,\rd x=0.
	$$
	The pressure estimates follow by induction on the material order, using this identity,
	the estimates for $\Psi_{q,j,u}$, and the commutator and mean identities above.

	\subsubsection{Error estimates for the Newton iteration}

	All Newton substeps are linearized about the same mollified state $(\orho_q,\om_q)$ and use the same time scale $\tau_q$.
	Set $\Lambda_q=\delta_q^{\frac12}\lambda_q$, so that $\tau_q\Lambda_q\leq\lambda_{q+1}^{-4\alpha}$.
	Choose a sufficiently large constant $C_*>1$, independently of $q$ and $a$, and set
	$$
	\begin{aligned}
		A_{q,0}&:=C_*\lambda_q^{-8\alpha}\delta_{q+1},&&
		A_{q,j+1}:=C_*\tau_q\Lambda_q(\tau_q\ell_q^{-1})^\alpha\ell_q^{-\alpha}A_{q,j}
		\leq\lambda_q^{-\frac{8\alpha}{3}}A_{q,j},
	\end{aligned}
	$$
	where $0\leq j<J_*$ and $a$ is sufficiently large.
	The interior stress satisfies
	$$
	\begin{aligned}
	\nrm{R_{q,j+1}}_0&\leq A_{q,j+1}\ell_q^{-\alpha}.
	\end{aligned}
	$$
	For $0\leq j<J_*$, the density and momentum increments satisfy
	$$
	\nrm{\rho_{q,j+1}-\rho_{q,j}}_0+\nrm{m_{q,j+1}-m_{q,j}}_0
	\lesssim A_{q,j}\delta_q^{-\frac12}\lambda_{q+1}^{-3\alpha}.
	$$

	For each local correction, the coefficient estimates give
	$$
	\nrm{p'(\orho_q)-c_{q,u,\kappa}^2}_0+\nrm{J_u-1}_0
	+\nrm{\nabla\Phi_u-\Id}_0+\nrm{\mathbf C_u-\Id}_0
	\lesssim\tau_q\Lambda_q,
	$$
	for $t\in[t_u-\tau_q,t_u+\tau_q]$.
	The norm of $p'(\orho_q)-c_{q,u,\kappa}^2$ is restricted to the spatial support of the local correction.
	The interior residuals satisfy
	$$
	\nrm{\mathcal S_{q,j+1,u,\kappa}}_\alpha
	+\tau_q\nrm{H_{q,j+1,u,\kappa}^{\mathrm{int}}}_\alpha
	\lesssim\tau_q\Lambda_q A_{q,j}\ell_q^{-2\alpha}.
	$$

	For the localization radius $O(\tau_q)$, the local inverse applied to $H_{q,j+1,u,\kappa}^{\mathrm{int}}$ contributes the factor $\tau_q(\tau_q\ell_q^{-1})^\alpha$.
	For the incompressible correction, the pressure identity above gives the factor $\tau_q\Lambda_q$ through the estimates for $\Psi_{q,j,u}$ and $\nabla\overline u_q$.
	The compressible cutoff residual has a factor $\tau_q^{-1}$, which is compensated by the factor $\tau_q$ in the local inverse.
	For the incompressible cutoff term, the factor $\tau_q^{-1}$ from $\partial_t\tilde\theta_u$ is compensated by the factor $\tau_q$ in the estimate for $\Psi_{q,j,u}$.
	Summing the local estimates and using the bounded overlap of the cutoffs, we obtain
	$$
	\begin{aligned}
	\nrm{R_{q,j+1}}_\alpha
	&\lesssim\tau_q\Lambda_q(\tau_q\ell_q^{-1})^\alpha A_{q,j}\ell_q^{-2\alpha},&&
	\nrm{R_{q,j+1}^{\mathrm{gluing}}}_\alpha
	\lesssim(\tau_q\ell_q^{-1})^\alpha A_{q,j}\ell_q^{-2\alpha}.
	\end{aligned}
	$$
	The factor $(\tau_q\ell_q^{-1})^\alpha$ is the loss in the H{\"o}lder estimate for the local inverse.
	The spatial H{\"o}lder factors have been estimated by $\ell_q^{-2\alpha}$.
	The additional factor $\tau_q\Lambda_q$ provides the gain for reducing the interior stress.

	Compared with the source estimate for the next Newton substep, the full derivative estimates have an additional H{\"o}lder factor bounded by $\ell_q^{-\alpha}$.
	We include this factor in $A_{q,j+1}$, so $R_{q,j+1}$ satisfies that source estimate; see Proposition~\ref{prop:Newton-interior-stress-estimates}.
	The amplitudes decrease geometrically while $\tau_q$ remains fixed.
	
	The tensor $P_{q+1,J_*}$ contains the initial mollification error $R_q-R_{q,0}$ and the source mollification errors accumulated over the Newton substeps.
	It also contains the remainders in $m\otimes m/\rho$ and $p(\rho)$ after subtracting the terms linearized about $(\orho_q,\om_q)$.
	The stopping depth $J_*$ and the parameter choices then give
	$$
	\nrm{R_{q,J_*}}_0+\nrm{P_{q+1,J_*}}_0\lesssim A_{q,J_*}\ell_q^{-\alpha}+A_{q,0}^2\delta_q^{-1}\lambda_{q+1}^{-6\alpha}\lesssim\lambda_q^{-\alpha}\lambda_{q+1}^{-8\alpha}\delta_{q+2}.
	$$

	After $J_*$ substeps, the Newton iteration gives the stress
	$
	R_{q,J_*}+P_{q+1,J_*}+S_{q,J_*},
	$ with $S_{q,J_*}=\sum_{j=0}^{J_*-1}R_{q,j+1}^{\mathrm{gluing}}$.
	To estimate $S_{q,J_*}$ in $C^0$, we use the gluing estimates at a fixed H{\"o}lder exponent sufficiently small relative to $\alpha$ and $\beta$, keeping the amplitudes $A_{q,j}$ unchanged.
	The resulting H{\"o}lder and localization factors are bounded by $(\tau_q\ell_q^{-1})^\alpha$.
	Summing the gluing stresses then gives
	$$
	\nrm{S_{q,J_*}}_0\lesssim A_{q,0}(\tau_q\ell_q^{-1})^\alpha.
	$$

	All substeps use the same temporal cutoffs, so
	$$
	\supp_tS_{q,J_*}
	\subset\bigcup_u\bigl(\supp\partial_t\theta_u\cup\supp\partial_t\tilde\theta_u\bigr).
	$$
	These supports are contained in separated transition intervals of length comparable to $\tau_q$, with gaps comparable to $\tau_q$.
	We include $R_{q,J_*}+P_{q+1,J_*}$ in the new Reynolds stress $R_{q+1}$ and use the Nash perturbation to cancel the symmetric gluing stress $S_{q,J_*}$.
	The Newton iteration thus provides the temporal localization needed for the Nash step, as does the gluing construction in~\cite{BDSV19}.
	The time scale $\tau_q$ remains unchanged throughout the Newton iteration and can also be used for the Nash cutoffs.
	The Newton estimates are proved in Section~\ref{sec:Newton steps}.

	\subsubsection{The Nash step}

	Let $\bar u_{q,J_*}$ and $\bar\rho_{q,J_*}$ denote the spatially mollified terminal velocity and density.
	We write $\bar S_{q,J_*,u}$ for the part of $S_{q,J_*}$ on the transition interval indexed by $u$,
	mollified along the flow of $\bar u_{q,J_*}$.
	Choose time cutoffs $\theta_u^{\mathrm N}$ with disjoint supports and $\theta_u^{\mathrm N}=1$ on $\supp_t\bar S_{q,J_*,u}$.
	We use spatially separated Mikado flows \cite{DS17,GK22}, composed on each selected interval with a common backward flow map of $\bar u_{q,J_*}$.
	
	The momentum perturbation $w_{q+1}=w_{q+1}^{\mathrm p}+w_{q+1}^{\mathrm c}$ (principal part and divergence corrector) is a curl.
	With $\rho_{q+1}=\rho_{q,J_*}$ and $m_{q+1}=m_{q,J_*}+w_{q+1}$, the continuity equation is preserved.
	Spatial separation and the disjoint time cutoffs eliminate cross terms between different time intervals or different Mikado directions.
	The amplitudes are chosen so that averaging the periodic profiles in $w_{q+1}^{\mathrm p}\otimes w_{q+1}^{\mathrm p}/\bar\rho_{q,J_*}$ cancels the mollified gluing stress up to a multiple of $\Id$ depending only on time.

	We construct $R_{q+1}^{\mathrm{trans}}$ and $R_{q+1}^{\mathrm{osc}}$ to satisfy
	$$
	\begin{aligned}
	\Div R_{q+1}^{\mathrm{trans}}
	&=(\partial_t+\bar u_{q,J_*}\cdot\nabla)w_{q+1}
	+(w_{q+1}\cdot\nabla)\bar u_{q,J_*}
	+(\Div\bar u_{q,J_*})w_{q+1},\\
	\Div R_{q+1}^{\mathrm{osc}}
	&=\Div\left(
	\frac{w_{q+1}^{\mathrm p}\otimes w_{q+1}^{\mathrm p}}{\bar\rho_{q,J_*}}
	-\sum_u(\theta_u^{\mathrm N})^2(A_{q,0}(\tau_q\ell_q^{-1})^\alpha\Id-\bar S_{q,J_*,u})
	\right).
	\end{aligned}
	$$
	We localize a vector potential of $w_{q+1}$ and the tensor inside the divergence in the second identity. 
	In Lemma~\ref{lem:Nash-transported-oscillatory-inverse-divergence}, a periodic potential for each profile gives an explicit divergence with a factor $\lambda_{q+1}^{-1}$.
	The local right inverse is then applied to the remainder, giving the same factor in the stress estimates.
	For sufficiently large $a$, we obtain
	$$
	\begin{aligned}
	\nrm{R_{q+1}^{\mathrm{trans}}}_0
	&\lesssim A_{q,0}^{\frac12}(\tau_q\ell_q^{-1})^{\frac\alpha2}\tau_q^{-1}\lambda_{q+1}^{-1},&&
	\nrm{R_{q+1}^{\mathrm{osc}}}_0
	\lesssim A_{q,0}(\tau_q\ell_q^{-1})^\alpha\lambda_q\lambda_{q+1}^{-1}.
	\end{aligned}
	$$

	To close the inductive stress estimate, we compare the transport and oscillation errors with $\lambda_{q+1}^{-8\alpha}\delta_{q+2}$.
	Ignoring the losses involving $\alpha$, we obtain for the transport error
	$$
	\frac{(\delta_q\delta_{q+1})^{\frac12}\lambda_q}
	{\lambda_{q+1}\delta_{q+2}}
	\lesssim\lambda_q^{-(b-1)(1-\beta(2b+1))}.
	$$
	For the oscillation error, we similarly have
	$$
	\frac{\delta_{q+1}\lambda_q}{\lambda_{q+1}\delta_{q+2}}
	\lesssim\lambda_q^{-(b-1)(1-2\beta b)}.
	$$
	Both exponents are negative if $\beta<\min\{\frac1 {2b+1},\frac1 {2b}\}=\frac1 {2b+1}$.
	For every $0<\beta<\frac 13$, this holds for $b>1$ sufficiently close to $1$.
	Combining these estimates with those for $R_{q,J_*}$ and $P_{q+1,J_*}$, we obtain, for sufficiently large $a$,
	$$
	\begin{aligned}
	\nrm{R_{q+1}}_0
	\lesssim
	\bigl(\lambda_q^{-\alpha}+\lambda_q^{-15\alpha}\bigr)
	\lambda_{q+1}^{-8\alpha}\delta_{q+2}
	\leq\lambda_{q+1}^{-8\alpha}\delta_{q+2}.
	\end{aligned}
	$$

	\subsubsection*{Organization of the paper}

	In Section~\ref{sec:main-iteration}, we fix the parameters and inductive assumptions,
	state the main iterative proposition, and derive Theorem~\ref{thm:main}.
	The mollified background and cutoffs are introduced in Section~\ref{Mollification},
	where we also prove the flow and mollification estimates.
	Section~\ref{sec:local-inverse-divergence} constructs a local right inverse of divergence
	with prescribed support and proves its spatial and material derivative estimates.
	For the Newton step, Section~\ref{sec:periodic-Hodge-and-Newton-stresses}
	establishes higher material derivative estimates for periodic nonlocal operators and the linearized Euler equations.
	The Newton corrections are constructed and estimated in Section~\ref{sec:Newton steps}.
	Section~\ref{sec:Nash-step} constructs the Nash perturbation and estimates the resulting errors.
	We also obtain estimates for the updated density and momentum, their increments,
	and the Reynolds stress with the material derivative at stage $q+1$.
	Section~\ref{sec:main-iteration-proof} completes the proof of the main iterative proposition.
	The appendices collect the tools used in the construction.
	\section{The main iterative proposition}
	\label{sec:main-iteration}

	We reduce Theorem~\ref{thm:main} to an iteration for smooth solutions
	of the compressible Euler--Reynolds system~\eqref{eq:compressible-Euler-Reynolds}.
	In this part, we first fix the notation and the inductive assumptions used in the iteration.
	\subsection{Notation and conventions}
	\label{subsec:notation-conventions}

	\begin{itemize}
		\item Unless stated otherwise, implicit constants in $A\lesssim_{M,N,s}B$
		may depend on $M,N,s$ and the parameters specified with the estimate,
		but not on $q$ or $a$.
		On a fixed finite range of derivative orders, we use one constant that
		may depend on the highest orders. Dependence on the fixed pressure law $p$ is suppressed in constants and parameter thresholds.
		\item For a vector $w\in\R^3$ and a multi-index
		$\kappa=(\kappa_1,\kappa_2,\kappa_3)\in\N_0^3$, we define
		$
		w^\kappa:=w_1^{\kappa_1}w_2^{\kappa_2}w_3^{\kappa_3}.
		$
		We also write
		$
		|\kappa|=\kappa_1+\kappa_2+\kappa_3
		$
		and
		$
		\partial^\kappa=\partial_1^{\kappa_1}\partial_2^{\kappa_2}\partial_3^{\kappa_3},
		$
		where $\partial_i=\partial_{x_i}$ for $i=1,2,3$. For an integer $N\in\N_0$,
		the notation $\nabla^N u$ denotes the collection of spatial derivatives
		of order $N$, with $\nabla^0u=u$. We use the same notation in dimension $d$.

		\item For a vector field $v=(v_1,v_2,v_3)$ on $\R^3$, we set
		$(\nabla v)_{ij}:=\partial_jv_i$. For a matrix field $R=(R_{ij})$, we use
		$(\Div R)_i:=\sum_{j=1}^3\partial_jR_{ij}$. The three-dimensional curl is
		$$
		\curl v
		:=
		(\partial_2v_3-\partial_3v_2,\partial_3v_1-\partial_1v_3,\partial_1v_2-\partial_2v_1).
		$$

		\item For a time interval $\mcI\subset\R$, $|\mcI|$ denotes its length.
		For a set $E\subset\R$ and $\tau>0$, we write
		$$
		\mN_\tau(E)
		:=
		\{\,t\in\R:\operatorname{dist}(t,E)<\tau\,\}
		$$
		for its open $\tau$-neighborhood.
		For temporal mollification at scale $\ell_t$ on $\mcI$, the
		field is assumed to be defined on $\mN_{\ell_t}(\mcI)$.
		We use the spatial norms defined in Subsection~\ref{sec:holder-spaces}.
		
		\item For a function $f=f(t,x)$ on $\R\times\R^3$, or for the periodic
		extension of a function on $\R\times\T^3$, we denote its space-time support by
		$
		\supp_{t,x} f:=\supp f\subset\R\times\R^3.
		$
		We also define
		$$
		\supp_t f
		:=
		\{\,t\in\R:\exists x\in\R^3\text{ such that }(t,x)\in\supp_{t,x} f\,\},
		$$
		and, for each fixed $t\in\R$,
		$$
		\supp_x f(t,\cdot)
		:=
		\supp(f(t,\cdot))
		=
		\overline{\{\,x\in\R^3:f(t,x)\neq0\,\}}.
		$$
		For a finite tuple of fields, $\supp_t$ denotes the union of their temporal supports.
		
		\item $B(x_0,\varepsilon)$ denotes the spatial ball centered at
		$x_0$ with radius $\varepsilon$. On the torus, these balls are understood
		through a fixed local lift when the radius is sufficiently small.
	\end{itemize}

	\subsection{Parameters and inductive assumptions}
	\label{subsec:parameters-inductive}

	For $q\in\mathbb N_0$, we define the frequency and amplitude parameters by
	\begin{equation}
		\lambda_q=\lceil a^{b^q}\rceil,
		\qquad
		\delta_q=\lambda_q^{-2\beta}.
		\label{def:lambda-delta}
	\end{equation}
	We will choose $a>1$ sufficiently large and $b>1$ close to $1$.
	The exponent $0<\beta<\frac13$ determines the target H{\"o}lder
	regularity of the density and momentum, with exponents below
	$2\beta$ and $\beta$, respectively.

	Assume that a smooth solution $(\rho_q,m_q,R_q)$ of
	\eqref{eq:compressible-Euler-Reynolds} has been constructed at stage
	$q$. We seek density and momentum perturbations
	$(\rho_{q+1}-\rho_q,m_{q+1}-m_q)$ that cancel $R_q$ up to
	a new stress $R_{q+1}$ with a smaller prescribed $C^0$ upper bound.
	The corrected triple $(\rho_{q+1},m_{q+1},R_{q+1})$ must solve
	the same system and satisfy the inductive estimates below with
	$q$ replaced by $q+1$.

	Fix a reference density $\rho^*>0$. Let $0<\alpha<1$,
	$\mathfrak M=\mathfrak M(\beta,\alpha,N_{\mathrm{max}})>1$, and integers
	$4\leq N_R<N_{\mathrm{ind}}<N_{\mathrm{max}}$ be parameters
	to be chosen in Proposition~\ref{prop:main-iteration}.
	We assume the following inductive estimates for $M,N\in\mathbb N_0$:
	\begin{align}
		&\|\rho_q-\rho^*\|_0
		\leq
		\mathfrak M\sum_{i=1}^q\delta_i,
		\label{ind:rho-zero}
		\\
		&\|D_{t,q}^{M}\rho_q\|_N
		\leq
		\delta_q
		\lambda_q^N
		\left(\lambda_q\delta_q^{\frac12}\right)^M,
		&&
		1\leq N+M\leq N_{\mathrm{ind}},
		\label{ind:rho-material}
		\\
		&\|m_q\|_0
		\leq
		\mathfrak M(1-\delta_q^{\frac12}),
		\label{ind:m-zero}
		\\
		&\|D_{t,q}^{M}m_q\|_N
		\leq
		\mathfrak M\delta_q^{\frac12}
		\lambda_q^N
		\left(\lambda_q\delta_q^{\frac12}\right)^M,
		&&
		1\leq N+M\leq N_{\mathrm{ind}},
		\label{ind:m-material}
		\\
		&\|R_q\|_0
		\leq
		\lambda_q^{-8\alpha}\delta_{q+1},
		\label{ind:R-zero}
		\\
		&\|D_{t,q}^{M}R_q\|_N
		\leq
		\delta_{q+1}
		\lambda_q^{N-8\alpha}
		\left(\lambda_q\delta_q^{\frac12}\right)^M,
		&&
		1\leq N+M\leq N_R,
		\label{ind:R-material}
	\end{align}
	where the material derivative associated with $(\rho_q,m_q)$ is defined by
	$D_{t,q}:=\partial_t+\frac{m_q}{\rho_q}\cdot\nabla$.
	To keep the density away from vacuum, we also assume
	\begin{equation}
		\mathfrak M\sum_{i=1}^{\infty}\delta_i
		\leq
		\frac{\rho^*}{2},\qquad \frac{\rho^*}{2}
		\leq
		\rho_q
		\leq
		\frac{3\rho^*}{2}.
		\label{ind:rho-gap}
	\end{equation}
	Finally, as in the hypotheses of \cite[Proposition~2.1]{GR24}, we impose
	\begin{equation}
		\supp_t R_q
		\subset
		\bigl[-2+\bigl(\lambda_q\delta_q^{\frac12}\bigr)^{-1},
		-1-\bigl(\lambda_q\delta_q^{\frac12}\bigr)^{-1}\bigr]
		\cup
		\bigl[1+\bigl(\lambda_q\delta_q^{\frac12}\bigr)^{-1},
		2-\bigl(\lambda_q\delta_q^{\frac12}\bigr)^{-1}\bigr].
		\label{ind:R-support}
	\end{equation}
	We take $a$ sufficiently large so that
	$(\lambda_0\delta_0^{\frac12})^{-1}<\frac14$. 
	
	To propagate these estimates, we require
	$b<(1+3\beta)/(6\beta)$
	and also take $b<1+2\beta$. For $\alpha$, we assume
	\begin{equation}
		0<\alpha<1,
		\qquad
		8\alpha
		<\frac{b-1}{b}\left(\frac13-\beta(2b-1)\right),
		\qquad
		(7-b)\alpha>\beta(b-1)^2.
		\label{def:admissible-alpha-losses}
	\end{equation}
	Moreover, we define the number of Newton substeps by
	\begin{equation}\label{def:Newton-stopping-depth}
		J_*:=\left\lceil\frac{b-1}{4\alpha}\right\rceil+1.
	\end{equation}
	\subsection{Main iterative proposition and proof of the main theorem}
	\label{subsec:main-proposition}
	\label{subsec:proof-main-theorem}
	
	\begin{prop}
		\label{prop:main-iteration}
		Fix $\rho^*>0$ and set
		$
		I_{\rho^*}:=[\rho^*/2,3\rho^*/2].
		$
		Let $p\in C^\infty((0,\infty))$ be a pressure law satisfying
		$p'(\rho)>0$ for every $\rho\in I_{\rho^*}$.
		Let $0<\beta<\frac13$ and
		$1<b<\min\{1+2\beta,(1+3\beta)/(6\beta)\}$.
		Fix $\alpha$ satisfying \eqref{def:admissible-alpha-losses}.
		Let $J_*$ be defined by \eqref{def:Newton-stopping-depth}.
		Choose integers $4\leq N_R<N_{\mathrm{ind}}<N_{\mathrm{max}}$ satisfying
		\eqref{eq:Newton-top-material-reserve} and the inequalities
		$N_R\geq4J_*+12$,
		$N_{\mathrm{ind}}\geq N_R+4J_*+12$, and
		$N_{\mathrm{max}}\geq N_{\mathrm{ind}}+4J_*+12$.
		Then there exist $\mathfrak M=\mathfrak M(\beta,\alpha,N_{\mathrm{max}})>1$ and
		$
		a_0
		=
		a_0
		(
		\beta,\allowbreak b,\allowbreak N_R,\allowbreak N_{\mathrm{ind}},\allowbreak N_{\mathrm{max}},\allowbreak \rho^*,
		\allowbreak \mathfrak M,\allowbreak \alpha
		)
		>
		1
		$
		with the following property.
		For every $a>a_0$ and $q\in\mathbb N_0$, let
		$(\rho_q,m_q,R_q)$ be a smooth solution of
		\eqref{eq:compressible-Euler-Reynolds} on $\R\times\T^3$ satisfying
		the inductive estimates \eqref{ind:rho-zero}--\eqref{ind:R-material}
		in their stated ranges, the density condition \eqref{ind:rho-gap},
		and the support condition \eqref{ind:R-support}.
		Then there exists a smooth solution $(\rho_{q+1},m_{q+1},R_{q+1})$ of
		\eqref{eq:compressible-Euler-Reynolds} satisfying
		\eqref{ind:rho-zero}--\eqref{ind:R-material}, \eqref{ind:rho-gap}, and
		\eqref{ind:R-support} with $q$ replaced by $q+1$.
		Moreover, the perturbations satisfy
		\begin{align}
			\|\rho_{q+1}-\rho_q\|_0
			+
			\lambda_{q+1}^{-1}\|\rho_{q+1}-\rho_q\|_1
			+
			\lambda_{q+1}^{-1}\|\partial_t(\rho_{q+1}-\rho_q)\|_0
			&\leq
			\mathfrak M\delta_{q+1},
			\label{inc:rho}
			\\
			\|m_{q+1}-m_q\|_0
			+
			\lambda_{q+1}^{-1}\|m_{q+1}-m_q\|_1
			+
			\lambda_{q+1}^{-1}\|\partial_t(m_{q+1}-m_q)\|_0
			&\leq
			\mathfrak M\delta_{q+1}^{\frac12}.
			\label{inc:m}
		\end{align}
		The perturbations also satisfy
		\begin{equation}
			\supp_t(\rho_{q+1}-\rho_q,m_{q+1}-m_q)
			\subset
			(-2,-1)\cup(1,2).
			\label{inc:support}
		\end{equation}
	\end{prop}

	\begin{rem}
	Starting from an initial triple satisfying the hypotheses of
	Proposition~\ref{prop:main-iteration}, we obtain smooth Euler--Reynolds triples
	$(\rho_q,m_q,R_q)$ for every $q\geq0$, with the same $a$ and $\mathfrak M$.
	Estimates~\eqref{inc:rho}--\eqref{inc:m} give convergence of $\rho_q$ and $m_q$
	in $C^{\gamma_\rho}(\R\times\T^3)$ and $C^{\gamma_m}(\R\times\T^3)$
	for $0\leq\gamma_\rho<2\beta$ and $0\leq\gamma_m<\beta$, respectively.
	By \eqref{ind:R-zero} and \eqref{ind:rho-gap}, the limit solves the
	compressible Euler equations with positive density.
	Condition~\eqref{inc:support} preserves $(\rho_0,m_0)$ on $[-1,1]$
	and outside $[-2,2]$. Thus the proof of Theorem~\ref{thm:main} reduces
	to choosing admissible parameters and an initial triple with
	$m_0\not\equiv0$ on $[-1,1]\times\T^3$ and
	$(\rho_0,m_0)=(\rho^*,0)$ for $t\notin[-2,2]$.
	\end{rem}

	\begin{rem}[Parameter compatibility]
		\label{subsec:closure-parameter-hierarchy}
		
		We verify the compatibility of the parameter conditions in the range
		$$
		1<b<\min\left\{1+2\beta,\frac{1+3\beta}{6\beta}\right\}.
		$$
		Assume \eqref{def:admissible-alpha-losses}.
		These parameter conditions imply
		$$
		4b\alpha<\frac{b-1}{6}<\frac{\beta}{3}<\beta,
		\qquad
		\alpha<\frac{b-1+2\beta}{8b}.
		$$
		
		An explicit choice valid for every
		$0<\beta<1/3$ is
		\begin{equation}\label{eq:closure-explicit-parameter-family}
			b:=1+\beta\left(\frac13-\beta\right),\qquad \alpha:=\frac{\beta\left(\frac13-\beta\right)^2}{12}.
		\end{equation}
		Indeed, $0<b-1<2\beta$, $b<7/6$, and
		$1/3-\beta(2b-1)=\left(\frac13-\beta\right)(1-2\beta^2)>0$. The latter identity
		also gives $b<(1+3\beta)/(6\beta)$. Since $2b/3<7/9<1-2\beta^2$,
		$$
		8b\alpha=\frac{2b}{3}\beta\left(\frac13-\beta\right)^2
		<\beta\left(\frac13-\beta\right)^2(1-2\beta^2)
		=(b-1)\left(\frac13-\beta(2b-1)\right).
		$$
		Finally,
		$$
		(7-b)\alpha-\beta(b-1)^2
		=\alpha(7-b-12\beta^2)>0.
		$$
		Thus \eqref{eq:closure-explicit-parameter-family} satisfies
		\eqref{def:admissible-alpha-losses}.
		For the choice \eqref{eq:closure-explicit-parameter-family},
		\eqref{def:Newton-stopping-depth} gives
		$$
		J_*=\left\lceil\frac{b-1}{4\alpha}\right\rceil+1
		=\left\lceil\frac{3}{\frac13-\beta}\right\rceil+1.
		$$
		For the remaining derivative orders in Proposition~\ref{prop:main-iteration}, we may take
		$N_R:=4J_*+12$, $N_{\mathrm{ind}}:=8J_*+24$, and
		$N_{\mathrm{max}}:=12J_*+36$.
		Since $J_*>3/(\frac13-\beta)$ and $b<7/6$, we have
		$$
		\begin{aligned}
			(b-1)\left(\frac12-\beta\right)(N_{\mathrm{ind}}-1)=\beta\left(\frac13-\beta\right)\left(\frac12-\beta\right)(8J_*+23)
			>24\beta\left(\frac12-\beta\right)>4\beta>(2b-1)\beta.
		\end{aligned}
		$$
		This also verifies \eqref{eq:Newton-top-material-reserve}.
	\end{rem}

	\begin{proof}[Proof of Theorem~\ref{thm:main}]
		Let $0\leq\gamma_m<\frac13$ and $0\leq\gamma_\rho<\frac23$. Choose
		$\max\{\gamma_m,\gamma_\rho/2\}<\beta<\frac13$ and
		$1<b<\min\{1+2\beta,(1+3\beta)/(6\beta)\}$ sufficiently close to $1$.
		Choose $\alpha$ satisfying \eqref{def:admissible-alpha-losses}.
		Define $J_*$ by \eqref{def:Newton-stopping-depth}, and choose integers
		$N_R\geq4J_*+12$,
		$N_{\mathrm{ind}}\geq N_R+4J_*+12$, and
		$N_{\mathrm{max}}\geq N_{\mathrm{ind}}+4J_*+12$
		satisfying \eqref{eq:Newton-top-material-reserve}.
		Let $\mathfrak M=\mathfrak M(\beta,\alpha,N_{\mathrm{max}})>1$ and $a_0$
		be given by Proposition~\ref{prop:main-iteration}.
		We will fix $a>a_0$ after verifying the estimates for the initial triple $(\rho_0,m_0,R_0)$.

		We first construct the initial triple. Fix a smooth function $f:\R\to[0,1]$ with
		$\supp f\subset\left[-\frac74,\frac74\right]$ and $f=1$ on
		$\left[-\frac54,\frac54\right]$.
		Let $e_2=(0,1,0)$ and define
		\begin{equation}
			\rho_0(t,x):=\rho^*,
			\qquad
			m_0(t,x):=\delta_0^{\frac12}f(t)\cos(\lambda_0x_1)e_2.
			\label{start:rho-m}
		\end{equation}
		Since $m_0$ is parallel to $e_2$ and independent of $x_2$,
		$\Div m_0=0$ and
		$\Div(m_0\otimes m_0/\rho_0)=0$. Also, $\nabla p(\rho_0)=0$. Set
		\begin{equation}
			R_0(t,x)
			:=
			\frac{\delta_0^{\frac12}f'(t)}{\lambda_0}
			\begin{pmatrix}
				0 & \sin(\lambda_0x_1) & 0\\
				\sin(\lambda_0x_1) & 0 & 0\\
				0 & 0 & 0
			\end{pmatrix}.
			\label{start:R}
		\end{equation}
		From \eqref{start:rho-m}--\eqref{start:R}, we have
		$$
		\Div R_0=\delta_0^{\frac12}f'(t)\cos(\lambda_0x_1)e_2
		=\partial_t m_0.
		$$
		Thus $(\rho_0,m_0,R_0)$ solves \eqref{eq:compressible-Euler-Reynolds}
		at stage $q=0$.

		Since $\rho_0=\rho^*$, \eqref{ind:rho-zero} and the pointwise bounds in
		\eqref{ind:rho-gap} hold at $q=0$.
		For $M\geq1$, $D_{t,0}^M\rho_0=0$.
		For $M=0$ and $N\geq1$, $\|\rho_0\|_N=\rho^*$.
		Since $1-2\beta>0$, increasing $a$ gives
		$\rho^*\leq\lambda_0^{1-2\beta}\leq\delta_0\lambda_0^N$.
		This verifies \eqref{ind:rho-material} at $q=0$.

		To estimate the material derivatives of $m_0$ and $R_0$, write
		$$
		D_{t,0}
		=
		\partial_t+\frac{m_0}{\rho_0}\cdot\nabla
		=
		\partial_t+\frac{\delta_0^{\frac12}f(t)\cos(\lambda_0x_1)}{\rho^*}\partial_2.
		$$
		Since $m_0$ and $R_0$ are independent of $x_2$, $D_{t,0}$ acts as $\partial_t$
		on both fields and all their iterated material derivatives.
		For $M,N\in\mathbb N_0$ with $1\leq N+M\leq N_{\mathrm{ind}}$,
		\eqref{start:rho-m} gives
		\begin{equation}\label{start:m-est}
			\|D_{t,0}^{M}m_0\|_N
			=\delta_0^{\frac12}\|f^{(M)}\|_{L^\infty(\R)}
			\sum_{r=0}^N\lambda_0^r
			\leq\frac{\|f^{(M)}\|_{L^\infty(\R)}}{1-\lambda_0^{-1}}
			\delta_0^{\frac12}\lambda_0^N
			\leq\mathfrak M\delta_0^{\frac12}
			\lambda_0^N\left(\lambda_0\delta_0^{\frac12}\right)^M.
		\end{equation}
		The last inequality holds for sufficiently large $a$. For $M=0$, use
		$\|f\|_{L^\infty(\R)}=1<\mathfrak M$. For $M\geq1$, use
		$\lambda_0\delta_0^{\frac12}\to\infty$.
		Moreover, $\|m_0\|_0=\delta_0^{\frac12}\leq\mathfrak M(1-\delta_0^{\frac12})$
		provided $\delta_0^{\frac12}\leq1/2$.
		Thus \eqref{ind:m-zero}--\eqref{ind:m-material} hold at $q=0$.

		Similarly, for every $M\in\mathbb N_0$, \eqref{start:R} gives
		$$
		D_{t,0}^{M}R_0(t,x)
		=\partial_t^{M}R_0(t,x)
		=\frac{\delta_0^{\frac12}f^{(M+1)}(t)}{\lambda_0}
		\begin{pmatrix}
			0 & \sin(\lambda_0x_1) & 0\\
			\sin(\lambda_0x_1) & 0 & 0\\
			0 & 0 & 0
		\end{pmatrix}.
		$$
		Hence, for $M,N\in\mathbb N_0$ with $M+N\leq N_R$,
		\begin{equation}
			\|D_{t,0}^{M}R_0\|_N
			\lesssim_{N_R}
			\delta_0^{\frac12}\lambda_0^{N-1}.
			\label{start:R-est-pre}
		\end{equation}
		By \eqref{def:admissible-alpha-losses},
		$\eta_R:=1+\beta-2b\beta-8\alpha>0$.
		Using \eqref{def:lambda-delta}, we have
		$
		\delta_1\lambda_0^{N-8\alpha}
		\geq
		2^{-2\beta}\lambda_0^{N-1-\beta}\lambda_0^{\eta_R}.
		$
		For sufficiently large $a$, \eqref{start:R-est-pre} therefore gives
		\begin{equation}
			\|R_0\|_0
			\leq
			\lambda_0^{-8\alpha}\delta_1,
			\label{start:R-zero-est}
		\end{equation}
		and, for all $M,N\in\mathbb N_0$ with $1\leq N+M\leq N_R$,
		\begin{equation}
			\|D_{t,0}^{M}R_0\|_N
			\leq
			\delta_1
			\lambda_0^{N-8\alpha}
			\left(\lambda_0\delta_0^{\frac12}\right)^M.
			\label{start:R-material-est}
		\end{equation}

		Finally, since $\supp_tR_0\subset\bigl[-\frac74,-\frac54\bigr]\cup
		\bigl[\frac54,\frac74\bigr]$, condition \eqref{ind:R-support} holds at $q=0$ whenever
		$\left(\lambda_0\delta_0^{\frac12}\right)^{-1}<\frac14$.
		We now fix $a>a_0$ so that the preceding estimates and support condition hold and
		$\mathfrak M\sum_{i=1}^{\infty}\delta_i\leq\rho^*/2$.
		Then $(\rho_0,m_0,R_0)$ satisfies all the assumptions of
		Proposition~\ref{prop:main-iteration} at stage $q=0$.
		
		\begin{samepage}
		Applying Proposition~\ref{prop:main-iteration} inductively, we obtain a sequence
		$(\rho_q,m_q,R_q)$ of smooth solutions to \eqref{eq:compressible-Euler-Reynolds}.
		Estimates \eqref{inc:rho}--\eqref{inc:m} and
		$\sum_{q=0}^{\infty}\delta_{q+1}^{\frac12}<\infty$ give uniform convergence of $\rho_q$ and $m_q$.
		For the space-time regularity, the same estimates give
		$$
		\begin{aligned}
		\|\rho_{q+1}-\rho_q\|_{C^1(\R\times\T^3)}
		&\lesssim\mathfrak M\delta_{q+1}\lambda_{q+1},&&
		\|m_{q+1}-m_q\|_{C^1(\R\times\T^3)}
		\lesssim\mathfrak M\delta_{q+1}^{\frac12}\lambda_{q+1}.
		\end{aligned}
		$$
		For positive exponents, H{\"o}lder interpolation yields
		\begin{align*}
		\sum_{q=0}^{\infty}
		\|\rho_{q+1}-\rho_q\|_{C^{\gamma_\rho}(\R\times\T^3)}
		&\lesssim_{\gamma_\rho}
		\sum_{q=0}^{\infty}
		\|\rho_{q+1}-\rho_q\|_0^{1-\gamma_\rho}
		\|\rho_{q+1}-\rho_q\|_{C^1(\R\times\T^3)}^{\gamma_\rho}
		\lesssim_{\gamma_\rho}
		\mathfrak M\sum_{q=0}^{\infty}\lambda_{q+1}^{\gamma_\rho-2\beta}<\infty,\\
		\sum_{q=0}^{\infty}
		\|m_{q+1}-m_q\|_{C^{\gamma_m}(\R\times\T^3)}
		&\lesssim_{\gamma_m}
		\sum_{q=0}^{\infty}
		\|m_{q+1}-m_q\|_0^{1-\gamma_m}
		\|m_{q+1}-m_q\|_{C^1(\R\times\T^3)}^{\gamma_m}
		\lesssim_{\gamma_m}
		\mathfrak M\sum_{q=0}^{\infty}\lambda_{q+1}^{\gamma_m-\beta}<\infty.
		\end{align*}
		\end{samepage}
		Thus $\rho_q\to\rho$ in $C^{\gamma_\rho}(\R\times\T^3)$ and
		$m_q\to m$ in $C^{\gamma_m}(\R\times\T^3)$. In particular,
		$
		\rho-\rho^*\in C^{\gamma_\rho}(\R\times\T^3),\ m\in C^{\gamma_m}(\R\times\T^3). 
		$

		By \eqref{ind:rho-gap}, both $\rho_q$ and its uniform limit $\rho$
		take values in $[\rho^*/2,3\rho^*/2]$.
		Consequently, $(m_q\otimes m_q)/\rho_q\to(m\otimes m)/\rho$ and
		$p(\rho_q)\to p(\rho)$ uniformly. Moreover, \eqref{ind:R-zero} gives
		$\|R_q\|_0\to0$ as $q\rightarrow \infty$. Passing to the limit in
		\eqref{eq:compressible-Euler-Reynolds} in the sense of distributions, we obtain
		\begin{equation*}
			\left\{
			\begin{aligned}
				&\partial_t\rho+\Div m=0,\\
				&\partial_t m+\Div\left(\frac{m\otimes m}{\rho}\right)+\nabla p(\rho)=0.
			\end{aligned}
			\right.
		\end{equation*}
		By \eqref{start:rho-m} and \eqref{inc:support},
		$(\rho,m)=(\rho^*,0)$ for $t\notin[-2,2]$.
		For $t\in[-1,1]$, \eqref{inc:support} and $f\equiv1$ give
		$\rho(t,x)=\rho^*$ and
		$m(t,x)=\delta_0^{\frac12}\cos(\lambda_0x_1)e_2$.
		Since $\delta_0>0$, the limiting momentum is not identically zero.
		This proves Theorem~\ref{thm:main}.
	\end{proof}

	\section{Mollification and cutoffs}
\label{Mollification}

\subsection{Spatial mollification of the Euler--Reynolds triple}

At stage $q$, we spatially mollify $(\rho_q,m_q,R_q)$ before the Newton
iteration.

We use a compactly supported spatial mollifier with finitely many vanishing moments,
as in \cite[Lemma~3.1]{GKN24} and \cite{NV23}.
Let $\eta\in C_c^\infty(\R^3;\R)$ be the spatial mollifier
from Appendix~\ref{app:mollification-estimates}, with $d=3$ and
$N_{\mathrm{mom}}=N_{\mathrm{ind}}$. The function $\eta$ need not be nonnegative.
It satisfies $\supp\eta\subset B_1(0)$ and
$$
\int_{\R^3}\eta(x)\rd x=1,
\qquad
\int_{\R^3}\eta(x)x^\sgm\rd x=0,
\quad 1\leq|\sgm|\leq N_{\mathrm{ind}}.
$$
Set
$
\eta_\ell(x):=\ell^{-3}\eta(x/\ell).
$
For a spatially periodic function $f$ on $\mcI\times\T^3$, we define
\begin{equation}
	\label{def:spatial-mollifier}
	P_{\leq\ell^{-1}}f:=\eta_\ell*f,
	\qquad
	P_{>\ell^{-1}}f:=f-P_{\leq\ell^{-1}}f.
\end{equation}
Here the convolution is taken only in the spatial variables, after periodically
extending $f$ to $\R^3$.

We use the vanishing moments up to order three in the Nash step and
fix the spatial mollification scale by
\begin{equation}
	\label{def:parameter-mollification}
	\ell_q=\lambda_q^{-\frac34}\lambda_{q+1}^{-\frac14}.
\end{equation}
In particular,
$
(\ell_q\lambda_q)^4=\lambda_q\lambda_{q+1}^{-1},\
\ell_q^{-2}\lambda_q^{-1}\lambda_{q+1}^{-1}=\lambda_q^{\frac12}\lambda_{q+1}^{-\frac12}.
$

We also fix the temporal kernel used for mollification along the flow. Let
$\eta^{\mathrm t}\in C_c^\infty((-1,1);\mathbb R)$ be a one-dimensional
kernel obtained from the construction in
Appendix~\ref{app:mollification-estimates}, with $d=1$ and
$N_{\mathrm{mom}}=N_{\mathrm{ind}}$. Thus
$$
\int_{\mathbb R}\eta^{\mathrm t}(s)\,\rd s=1,
\qquad
\int_{\mathbb R}s^k\eta^{\mathrm t}(s)\,\rd s=0,
\quad 1\leq k\leq N_{\mathrm{ind}}.
$$
We set
\begin{equation}
	\ell_{t,q}:=\delta_q^{-\frac12}\lambda_q^{-\frac12}
	\lambda_{q+1}^{-\frac12},
	\qquad
	\eta_{\ell_{t,q}}^{\mathrm t}(s):=\ell_{t,q}^{-1}\eta^{\mathrm t}(s/\ell_{t,q}).
	\label{def:material-mollification-scale}
\end{equation}

We define the spatially mollified fields by
\begin{equation}
	\label{def:regularized-terms}
	\orho_q:=\PLq\rho_q,
	\qquad
	\om_q:=\PLq m_q,
	\qquad
	R_{q,0}^{\mathrm{sp}}:=\PLq R_q.
\end{equation}
The associated velocity fields are
$$
u_q:=\frac{m_q}{\rho_q},
\qquad
\overline{u}_q:=\frac{\om_q}{\orho_q}
$$
and we write
$
\DTL:=\partial_t+\overline{u}_q\cdot\nabla
$
for the material derivative generated by $\overline{u}_q$. 

We write $(\cdot)_+:=\max\{\cdot,0\}$.
We use the following two-scale notation. For $K\in\N_0$,
$\lambda,\widetilde\lambda\geq1$, and $N\geq0$, set
$$
\MM{N,K,\lambda,\widetilde\lambda}
:=
\lambda^{\min\{N,K\}}
\widetilde\lambda^{(N-K)_+}.
$$
We also set
$
\Lambda_q:=\lambda_q\delta_q^{\frac12}
$
and
$$
\La{N,K,\lambda,\widetilde\lambda}
:=
\begin{cases}
	\lambda^\alpha, & 0\leq N<K,\\
	\widetilde\lambda^\alpha, & N\geq K.
\end{cases}
$$
The spatial mollification satisfies the following estimates. 
\begin{prop}
	\label{prop:uniform-mollification-estimates}
	Suppose that $(\rho_q,m_q,R_q)$ is a smooth solution of the
	compressible Euler--Reynolds system~\eqref{eq:compressible-Euler-Reynolds} and satisfies the inductive estimates
	\eqref{ind:rho-zero}--\eqref{ind:R-material} in their stated ranges and
	the density condition \eqref{ind:rho-gap}.
	There exists $a_{\mathrm{mol}}=a_{\mathrm{mol}}(\beta,b,\alpha,\rho^*,\mathfrak M,N_R,N_{\mathrm{ind}},N_{\mathrm{max}})>1$.
	For $a>a_{\mathrm{mol}}$, the following conclusions hold.
	The mollified density satisfies
	$$
	\frac{\rho^*}{2}
	\leq
	\orho_q
	\leq
	\frac{3\rho^*}{2}.
	$$

	The spatially mollified density, momentum, and stress satisfy
	\begin{align}
		\|\DTL^M(\orho_q-\rho^*)\|_N
		&\lesssim
		\delta_q\Lambda_q^M
		\MM{N,N_{\mathrm{ind}}-M,\lambda_q,\ell_q^{-1}},
		&&0\leq M\leq N_{\mathrm{ind}},
		\quad
		1\leq M+N\leq N_{\mathrm{max}},
		\label{est:uniform-orho}\\
		\|\DTL^M\om_q\|_N
		&\lesssim
		\delta_q^{\frac12}
		\Lambda_q^M
		\MM{N,N_{\mathrm{ind}}-M,\lambda_q,\ell_q^{-1}},
		&&0\leq M\leq N_{\mathrm{ind}},
		\quad
		1\leq M+N\leq N_{\mathrm{max}},
		\label{est:uniform-om}\\
		\|\DTL^M R_{q,0}^{\mathrm{sp}}\|_N
		&\lesssim
		\lambda_q^{-8\alpha}\delta_{q+1}
		\Lambda_q^M
		\MM{N,N_R-M,\lambda_q,\ell_q^{-1}},
		&&0\leq M\leq N_R,
		\quad
		M+N\leq N_{\mathrm{max}}.
		\label{est:uniform-Rq0}
	\end{align}
	For the undifferentiated fields, we have
	\begin{align}
		&\|\orho_q-\rho^*\|_0
		\lesssim
		\mathfrak M\sum_{i=1}^q\delta_i,
		\label{est:uniform-orho-zero}\\
		&\|\om_q\|_0
		\leq
		\mathfrak M\left(1-\frac12\delta_q^{\frac12}\right),
		\qquad
		\|\overline{u}_q\|_0
		\leq
		\frac{2\mathfrak M}{\rho^*}.
		\label{est:uniform-zero-velocity}
	\end{align}
	The mollified velocity $\overline{u}_q$, its gradient, and its divergence satisfy
	\begin{align} 
		\|\DTL^M\overline{u}_q\|_N 
		&\lesssim 
		\delta_q^{\frac12} 
		\Lambda_q^M 
		\MM{N,N_{\mathrm{ind}}-M,\lambda_q,\ell_q^{-1}}, 
		&&0\leq M\leq N_{\mathrm{ind}}, 
		\quad 
		1\leq M+N\leq N_{\mathrm{max}}, 
		\label{est:uniform-om-over-orho}\\ 
		\|\DTL^M\nabla\overline{u}_q\|_N 
		&\lesssim 
		\delta_q^{\frac12}\Lambda_q^M 
		\MM{N+1,N_{\mathrm{ind}}-M,\lambda_q,\ell_q^{-1}}, 
		&&0\leq M\leq N_{\mathrm{ind}}, 
		\quad 
		M+N\leq N_{\mathrm{max}}-1, 
		\label{est:uniform-grad-om-over-orho}\\ 
		\|\DTL^M\Div\overline{u}_q\|_N 
		&\lesssim 
		\delta_q\Lambda_q^{M+1} 
		\MM{N,N_{\mathrm{ind}}-M-1,\lambda_q,\ell_q^{-1}}, 
		&&0\leq M\leq N_{\mathrm{ind}}-1, 
		\quad 
		M+N\leq N_{\mathrm{max}}-1. 
		\label{est:uniform-div-om-over-orho} 
	\end{align}
	The differences between the original and spatially mollified fields satisfy
	\begin{align}
		\|\DTL^M(\rho_q-\orho_q)\|_N
		&\lesssim
		\delta_q\Lambda_q^M
		\ell_q^{N_{\mathrm{ind}}-M-N}\lambda_q^{N_{\mathrm{ind}}-M},
		&&0\leq M\leq N_{\mathrm{ind}},
		\quad
		M+N\leq N_{\mathrm{ind}},
		\label{est:uniform-rho-minus-orho}\\
		\|\DTL^M(m_q-\om_q)\|_N
		&\lesssim
		\delta_q^{\frac12}\Lambda_q^M
		\ell_q^{N_{\mathrm{ind}}-M-N}\lambda_q^{N_{\mathrm{ind}}-M},
		&&0\leq M\leq N_{\mathrm{ind}},
		\quad
		M+N\leq N_{\mathrm{ind}},
		\label{est:uniform-m-minus-om}\\
		\|\DTL^M(R_q-R_{q,0}^{\mathrm{sp}})\|_N
		&\lesssim
		\lambda_q^{-8\alpha}\delta_{q+1}\Lambda_q^M
		\ell_q^{N_R-M-N}\lambda_q^{N_R-M},
		&&0\leq M\leq N_R,
		\quad
		M+N\leq N_R,
		\label{est:uniform-R-minus-Rq0}\\
		\|\DTL^M(\overline{u}_q-u_q)\|_N
		&\lesssim
		\delta_q^{\frac12}\Lambda_q^M
		\ell_q^{N_{\mathrm{ind}}-M-N}\lambda_q^{N_{\mathrm{ind}}-M},
		&&0\leq M\leq N_{\mathrm{ind}},
		\quad
		M+N\leq N_{\mathrm{ind}}.
		\label{est:uniform-velocity-difference}
	\end{align}
	
	Moreover, the following H{\"o}lder estimates hold.
	For the undifferentiated fields, one has
	$$
	\begin{aligned}
		\|\orho_q-\rho^*\|_\alpha
		\lesssim
		\delta_q^\alpha\lambda_q^\alpha,\quad
		\|\om_q\|_\alpha
		\lesssim
		\delta_q^{\frac\alpha2}\lambda_q^\alpha,\quad
		\|\overline{u}_q\|_\alpha
		\lesssim
		\delta_q^{\frac\alpha2}\lambda_q^\alpha.
	\end{aligned}
	$$
	The spatially mollified density, momentum, and stress satisfy
	\begin{align} 
		\|\DTL^M(\orho_q-\rho^*)\|_{N+\alpha} 
		&\lesssim 
		\delta_q\Lambda_q^M 
		\La{N,N_{\mathrm{ind}}-M,\lambda_q,\ell_q^{-1}}
		\notag\\
		&\quad\cdot
		\MM{N,N_{\mathrm{ind}}-M,\lambda_q,\ell_q^{-1}}, 
		\hspace{40pt}0\leq M\leq N_{\mathrm{ind}}, 
		\quad 
		1\leq M+N\leq N_{\mathrm{max}}-1, 
		\label{est:holder-orho}\\ 
		\|\DTL^M\om_q\|_{N+\alpha} 
		&\lesssim 
		\delta_q^{\frac12}\Lambda_q^M 
		\La{N,N_{\mathrm{ind}}-M,\lambda_q,\ell_q^{-1}}
		\notag\\
		&\quad\cdot
		\MM{N,N_{\mathrm{ind}}-M,\lambda_q,\ell_q^{-1}}, 
		\hspace{40pt}0\leq M\leq N_{\mathrm{ind}}, 
		\quad 
		1\leq M+N\leq N_{\mathrm{max}}-1, 
		\label{est:holder-om}\\ 
		\|\DTL^M R_{q,0}^{\mathrm{sp}}\|_{N+\alpha} 
		&\lesssim 
		\lambda_q^{-8\alpha}\delta_{q+1}\Lambda_q^M 
		\La{N,N_R-M,\lambda_q,\ell_q^{-1}}
		\notag\\
		&\quad\cdot
		\MM{N,N_R-M,\lambda_q,\ell_q^{-1}}, 
		\hspace{46pt}0\leq M\leq N_R, 
		\quad 
		M+N\leq N_{\mathrm{max}}-1. 
		\label{est:holder-Rq0} 
	\end{align} 
	The mollified velocity $\overline{u}_q$, its gradient, and its divergence satisfy
	\begin{align} 
		\|\DTL^M\overline{u}_q\|_{N+\alpha} 
		&\lesssim 
		\delta_q^{\frac12}\Lambda_q^M 
		\La{N,N_{\mathrm{ind}}-M,\lambda_q,\ell_q^{-1}}
		\notag\\
		&\quad\cdot
		\MM{N,N_{\mathrm{ind}}-M,\lambda_q,\ell_q^{-1}}, 
		&&0\leq M\leq N_{\mathrm{ind}}, 
		\quad 
		1\leq M+N\leq N_{\mathrm{max}}-1, 
		\label{est:holder-om-over-orho}\\ 
		\|\DTL^M\nabla\overline{u}_q\|_{N+\alpha} 
		&\lesssim 
		\delta_q^{\frac12}\Lambda_q^M 
		\La{N+1,N_{\mathrm{ind}}-M,\lambda_q,\ell_q^{-1}}
		\notag\\
		&\quad\cdot
		\MM{N+1,N_{\mathrm{ind}}-M,\lambda_q,\ell_q^{-1}}, 
		&&0\leq M\leq N_{\mathrm{ind}}, 
		\quad 
		M+N\leq N_{\mathrm{max}}-2, 
		\label{est:holder-grad-om-over-orho}\\ 
		\|\DTL^M\Div\overline{u}_q\|_{N+\alpha} 
		&\lesssim 
		\delta_q\Lambda_q^{M+1} 
		\La{N,N_{\mathrm{ind}}-M-1,\lambda_q,\ell_q^{-1}}
		\notag\\
		&\quad\cdot
		\MM{N,N_{\mathrm{ind}}-M-1,\lambda_q,\ell_q^{-1}}, 
		&&0\leq M\leq N_{\mathrm{ind}}-1, 
		\quad 
		M+N\leq N_{\mathrm{max}}-2. 
		\label{est:holder-div-om-over-orho} 
	\end{align}
	The differences between the original and spatially mollified fields satisfy
	\begin{align}
		\|\DTL^M(\rho_q-\orho_q)\|_{N+\alpha}
		&\lesssim
		\delta_q\Lambda_q^M
		\ell_q^{N_{\mathrm{ind}}-M-N-\alpha}\lambda_q^{N_{\mathrm{ind}}-M},
		&&0\leq M\leq N_{\mathrm{ind}}-1,
		&&
		M+N\leq N_{\mathrm{ind}}-1,
		\label{est:holder-rho-minus-orho}\\
		\|\DTL^M(m_q-\om_q)\|_{N+\alpha}
		&\lesssim
		\delta_q^{\frac12}\Lambda_q^M
		\ell_q^{N_{\mathrm{ind}}-M-N-\alpha}\lambda_q^{N_{\mathrm{ind}}-M},
		&&0\leq M\leq N_{\mathrm{ind}}-1,
		&&
		M+N\leq N_{\mathrm{ind}}-1,
		\label{est:holder-m-minus-om}\\
		\|\DTL^M(R_q-R_{q,0}^{\mathrm{sp}})\|_{N+\alpha}
		&\lesssim
		\lambda_q^{-8\alpha}\delta_{q+1}\Lambda_q^M
		\ell_q^{N_R-M-N-\alpha}\lambda_q^{N_R-M},
		&&0\leq M\leq N_R-1,
		&&
		M+N\leq N_R-1,
		\label{est:holder-R-minus-Rq0}\\
		\|\DTL^M(\overline{u}_q-u_q)\|_{N+\alpha}
		&\lesssim
		\delta_q^{\frac12}\Lambda_q^M
		\ell_q^{N_{\mathrm{ind}}-M-N-\alpha}\lambda_q^{N_{\mathrm{ind}}-M},
		&&0\leq M\leq N_{\mathrm{ind}}-1,
		&&
		M+N\leq N_{\mathrm{ind}}-1.
		\label{est:holder-velocity-difference}
	\end{align}
	The implicit constants may depend on $\rho^*,\mathfrak M,N_{\mathrm{max}},\alpha$.
	They are independent of $q$ and $a$.

\end{prop}

\begin{proof}
	By \eqref{ind:rho-zero} and
	\eqref{est:appendix-spatial-mollification} at $N=0$, we have
	\begin{equation}\label{eq:mollification-density-smallness}
	\|\orho_q-\rho^*\|_0
	\leq
	\|\eta\|_{L^1(\mathbb R^3)}\mathfrak M\sum_{i=1}^q\delta_i\leq\frac{\rho^*}{2}.
	\end{equation}
	For $a>a_{\mathrm{mol}}$, with $a_{\mathrm{mol}}$ chosen large enough, the last inequality holds uniformly in $q$.
	Hence
	$$
	\frac{\rho^*}{2}\leq\orho_q\leq\frac{3\rho^*}{2},
	$$
	and \eqref{est:uniform-orho-zero} follows. By
	\eqref{ind:m-zero}, \eqref{ind:rho-gap}, and the preceding inequalities,
	the $C^0$ assumptions in
	Lemma~\ref{lem:appendix-mollified-quotient-velocity} hold with
	$$
	c_{\rho^*}:=\min\left\{1,\frac{\rho^*}{2},\frac{2}{3\rho^*}\right\},
	\qquad
	\widetilde{\mC}_0:=1+\mathfrak M+\rho^*.
	$$
	
	Since $\mathfrak M>1$, \eqref{ind:rho-material}, \eqref{ind:m-material}, and
	\eqref{ind:rho-gap} give
	$$
	\begin{aligned}
		\|\Dtq^M(\rho_q-\rho^*)\|_N
		\leq
		\mathfrak M\delta_q\Lambda_q^M\lambda_q^N,\quad
		\|\Dtq^Mm_q\|_N
		\leq
		\mathfrak M\delta_q^{\frac12}\Lambda_q^M\lambda_q^N,
		\qquad
		1\leq M+N\leq N_{\mathrm{ind}}.
	\end{aligned}
	$$
	We apply
	Lemma~\ref{lem:appendix-mollified-quotient-velocity} with
	$$
	\begin{aligned}
		&\varrho=\rho_q,
		&&
		w=m_q,
		&&
		\overline\varrho=\orho_q,
		&&
		\overline w=\om_q,
		&&K^*=N_{\mathrm{ind}},
		&&
		K^\sharp=N_{\mathrm{max}},\\
		&
		\lambda=\lambda_q,
		&&
		\ell=\ell_q,
		&&
		\mu=\mathfrak M\Lambda_q,
		&&\mC_\varrho=\mathfrak M\delta_q,
		&&
		\mC_w=\mathfrak M\delta_q^{\frac12},
		&&
		\mC_\upsilon=\max\{\mC_\varrho,\mC_w\}.
	\end{aligned}
	$$
	Here $\varrho_*=\rho^*$ and $K^*=N_{\mathrm{ind}}=N_{\mathrm{mom}}$.
	Increasing $a$ once more, we ensure that
	$\mC_\varrho,\mC_w\leq1$ and
	$
	\mC_\upsilon\lambda_q\leq\mathfrak M\Lambda_q=\mu,
	\
	\ell_q\lambda_q\leq1,
	\
	\mu\geq1.
	$
	Thus the remaining hypotheses are satisfied.
	
	Estimate~\eqref{est:appendix-mollified-quotient-main} proves
	\eqref{est:uniform-orho}--\eqref{est:uniform-om} and
	\eqref{est:uniform-om-over-orho}. The gradient estimate in
	Lemma~\ref{lem:appendix-mollified-quotient-velocity} yields
	\eqref{est:uniform-grad-om-over-orho}. Its divergence estimate,
	using the continuity equation, proves \eqref{est:uniform-div-om-over-orho}.
	Estimates~\eqref{est:appendix-mollified-quotient-high-frequency} and
	\eqref{est:appendix-mollified-quotient-material-high-frequency} imply
	\eqref{est:uniform-rho-minus-orho}--\eqref{est:uniform-m-minus-om} and
	\eqref{est:uniform-velocity-difference} for $M=0$ and
	$1\leq M\leq N_{\mathrm{ind}}$, respectively.
	For \eqref{est:uniform-zero-velocity}, use \eqref{est:uniform-m-minus-om} to obtain
	$$
	\|m_q-\om_q\|_0
	\leq C
	(\ell_q\lambda_q)^{N_{\mathrm{ind}}}\delta_q^{\frac12}\leq\frac12\delta_q^{\frac12}.
	$$
	Here $C$ is at least as large as the implicit constant in \eqref{est:uniform-m-minus-om} at $M=N=0$.
	By \eqref{def:parameter-mollification}, the last inequality holds uniformly in $q$ after increasing $a_{\mathrm{mol}}$.
	Hence
	$$
	\|\om_q\|_0
	\leq
	\|m_q\|_0+\|m_q-\om_q\|_0
	\leq
	\mathfrak M\left(1-\frac12\delta_q^{\frac12}\right).
	$$
	Since $\orho_q\geq\rho^*/2$, the first estimate in \eqref{est:uniform-zero-velocity} yields
	$$
	\|\overline{u}_q\|_0
	\leq\frac{2}{\rho^*}\|\om_q\|_0
	\leq\frac{2\mathfrak M}{\rho^*}\left(1-\frac12\delta_q^{\frac12}\right)
	\leq\frac{2\mathfrak M}{\rho^*}.
	$$
	
	To apply Lemma~\ref{lem:appendix-change-material-derivative}, note that
	$
	\DTL=\Dtq+(\overline{u}_q-u_q)\cdot\nabla.
	$
	The required coefficient estimates are
	\begin{equation}\label{eq:mollification-material-change-coefficient-estimates}
	\begin{aligned}
		\|\DTL^j(\overline{u}_q-u_q)\|_N
		&\lesssim
		\delta_q^{\frac12}\Lambda_q^j\lambda_q^N,\\
		\|\DTL^j\nabla\overline{u}_q\|_N
		+
		\|\nabla\DTL^j\overline{u}_q\|_N
		&\lesssim
		\delta_q^{\frac12}\Lambda_q^j\lambda_q^{N+1},
		\qquad
		0\leq j<N_R,
		\quad
		j+N\leq N_R-1.
	\end{aligned}
	\end{equation}
	These follow from \eqref{est:uniform-velocity-difference},
	$\ell_q\lambda_q\leq1$, and the gradient estimate in
	Lemma~\ref{lem:appendix-mollified-quotient-velocity}.
	After increasing $a$, we apply
	Lemma~\ref{lem:appendix-change-material-derivative} with
	$$
	\begin{aligned}
		&K=J=N_R,
		&&
		u=u_q,
		&&
		\overline u=\overline{u}_q,
		&&
		H=R_q,
		&&\mC_u=\mathfrak M\delta_q^{\frac12},
		&&
		\mC_H=\mathfrak M\lambda_q^{-8\alpha}\delta_{q+1},
		&&
		\lambda=\lambda_q,
		&&
		\mu=\mathfrak M\Lambda_q.
	\end{aligned}
	$$
	Since $j+N+1\leq N_R<N_{\mathrm{ind}}$, \eqref{ind:R-material} and
	\eqref{eq:mollification-material-change-coefficient-estimates} verify the hypotheses. Moreover,
	$\mC_u\leq1$ and $\mC_u\lambda_q=\mu$. Hence
	\eqref{ind:R-zero} at zero total order and
	\eqref{est:appendix-change-material-derivative} at positive total orders, together with
	\eqref{est:appendix-change-material-derivative-gradient}, give
	\begin{equation}
		\begin{aligned}
			\|\DTL^{M_1}\Dtq^{M_2}R_q\|_N
			&\lesssim
			\lambda_q^{-8\alpha}\delta_{q+1}
			\Lambda_q^{M_1+M_2}\lambda_q^N,
			&&M_1+M_2+N\leq N_R,\\
			\|\DTL^{M_1}\nabla\Dtq^{M_2}R_q\|_N
			+
			\|\nabla\DTL^{M_1}\Dtq^{M_2}R_q\|_N
			&\lesssim
			\lambda_q^{-8\alpha}\delta_{q+1}
			\Lambda_q^{M_1+M_2}\lambda_q^{N+1},
			&&M_1+M_2+N\leq N_R-1.
		\end{aligned}
		\label{proof:DtL-unmollified-stress}
	\end{equation}
	Taking $M_2=0$ proves the required estimates for $R_q$.
	For $0\leq M\leq N_R-1$ and $M+N\leq N_{\mathrm{max}}-1$,
	the $C^N$ estimates for the gradient in
	Lemma~\ref{lem:appendix-mollified-quotient-velocity} give
	$$
	\begin{aligned}
		\|\DTL^M\nabla\overline u_q\|_N+\|\nabla\DTL^M\overline u_q\|_N
		\lesssim\delta_q^{\frac12}\Lambda_q^M\lambda_q\MM{N,N_{\mathrm{ind}}-M-1,\lambda_q,\ell_q^{-1}}
		\leq\delta_q^{\frac12}\Lambda_q^M\lambda_q\MM{N,N_R-M-1,\lambda_q,\ell_q^{-1}}.
	\end{aligned}
	$$
	Since $N_R\leq N_{\mathrm{ind}}=N_{\mathrm{mom}}$, we apply
	Proposition~\ref{prop:appendix-material-mollification} to $f=R_q$ and the
	velocity $\overline{u}_q$, with $\ell=\ell_q$ and
	$$
	\begin{aligned}
		&K^*=N_R,
		&&
		K^\sharp=N_{\mathrm{max}},
		&&
		K_f=0,&&
		\mC_f=\mathfrak M\lambda_q^{-8\alpha}\delta_{q+1},
		&&
		\mC_\upsilon=\mathfrak M\delta_q^{\frac12},
		&&
		\lambda_f=\lambda_\upsilon=\lambda_q,&&
		\mu_f=\widetilde\mu_f=\mu_\upsilon=\mathfrak M\Lambda_q.
	\end{aligned}
	$$
	Also take $K_\upsilon=1$ and $\tilde\mu_\upsilon=\mu_\upsilon$.
	Taking $M_2=0$ in \eqref{proof:DtL-unmollified-stress} and using the
	velocity-gradient estimate verifies \eqref{ass:appendix-material-mollification}
	and \eqref{ass:appendix-extended-velocity}. After increasing $a$,
	$\lambda_*=\lambda_q$, $\ell_q\lambda_q\leq1$, and
	$\mu_{\upsilon,*}=\tilde\mu_{\upsilon,*}=\mu_*=\widetilde\mu_*=\mathfrak M\Lambda_q\geq1$.
	By \eqref{est:appendix-extended-spatial-mollification}, we have
	$$ 
	\|\DTL^M R_{q,0}^{\mathrm{sp}}\|_N 
	\lesssim 
	\lambda_q^{-8\alpha}\delta_{q+1}\Lambda_q^M 
	\MM{N,N_R-M,\lambda_q,\ell_q^{-1}}, 
	\qquad 
	0\leq M\leq N_R,\quad 
	M+N\leq N_{\mathrm{max}}. 
	$$ 
	We obtain \eqref{est:uniform-Rq0}.
	
	For $M+N\leq N_R$, set 
	$L:=N_R-M-N$. Then 
	$0\leq L\leq N_R\leq N_{\mathrm{mom}}$ and 
	$M+N+L=N_R=K^*$. The third estimate in \eqref{est:appendix-material-mollification} therefore gives
	$$ 
	\|\DTL^M(R_q-R_{q,0}^{\mathrm{sp}})\|_N 
	\lesssim 
	\lambda_q^{-8\alpha}\delta_{q+1}\Lambda_q^M 
	\ell_q^{N_R-M-N}\lambda_q^{N_R-M}, 
	\qquad 
	M+N\leq N_R. 
	$$ 
	Thus \eqref{est:uniform-R-minus-Rq0} holds, which completes the proof of the
	$C^N$ estimates.
	
	For the H{\"o}lder estimates, we use
	$$
	\|f\|_{N+\alpha}
	\lesssim
	\|f\|_N^{1-\alpha}\|f\|_{N+1}^{\alpha}
	$$
	and the identity
	$$
	\MM{N,K,\lambda_q,\ell_q^{-1}}^{1-\alpha}
	\MM{N+1,K,\lambda_q,\ell_q^{-1}}^\alpha
	=
	\La{N,K,\lambda_q,\ell_q^{-1}}\MM{N,K,\lambda_q,\ell_q^{-1}}.
	$$
	Interpolating the corresponding estimates in $C^N$
	and $C^{N+1}$ proves \eqref{est:holder-orho}--\eqref{est:holder-om} when $M+N\geq1$,
	\eqref{est:holder-Rq0} throughout its stated range, and
	\eqref{est:holder-om-over-orho} when $M+N\geq1$. Applying the same argument to the
	gradient and divergence estimates proves
	\eqref{est:holder-grad-om-over-orho}--%
	\eqref{est:holder-div-om-over-orho}.
	
	At $M=N=0$, interpolation between the $C^0$ and $C^1$ estimates gives
	$$
	\begin{aligned}
		\|\orho_q-\rho^*\|_\alpha
		&\lesssim
		\|\orho_q-\rho^*\|_0^{1-\alpha}
		\|\orho_q-\rho^*\|_1^\alpha
		\lesssim
		\delta_q^\alpha\lambda_q^\alpha,\\
		\|\om_q\|_\alpha
		&\lesssim
		\|\om_q\|_0^{1-\alpha}\|\om_q\|_1^\alpha
		\lesssim
		\delta_q^{\frac\alpha2}\lambda_q^\alpha,\\
		\|\overline{u}_q\|_\alpha
		&\lesssim
		\|\overline{u}_q\|_0^{1-\alpha}
		\|\overline{u}_q\|_1^\alpha
		\lesssim
		\delta_q^{\frac\alpha2}\lambda_q^\alpha.
	\end{aligned}
	$$
	Interpolation of the estimates for the mollification differences uses
	$$
	\left(\ell_q^{N_{\mathrm{ind}}-M-N}\lambda_q^{N_{\mathrm{ind}}-M}\right)^{1-\alpha}
	\left(\ell_q^{N_{\mathrm{ind}}-M-N-1}\lambda_q^{N_{\mathrm{ind}}-M}\right)^\alpha
	=
	\ell_q^{N_{\mathrm{ind}}-M-N-\alpha}\lambda_q^{N_{\mathrm{ind}}-M}.
	$$
	For the stress, replace $N_{\mathrm{ind}}$ by $N_R$.
	Interpolating \eqref{est:uniform-rho-minus-orho}--\eqref{est:uniform-velocity-difference} at spatial orders $N$ and $N+1$
	proves \eqref{est:holder-rho-minus-orho}--\eqref{est:holder-velocity-difference}.
\end{proof}

\subsection{Cutoffs}
\label{subsec:cutoffs}

The spatial and temporal cutoffs below will localize the perturbations in the Newton step.
We define the
space-time localization scale $\tau_q>0$ by
\begin{equation}
	\tau_q^{-1}
	:=
	2\left\lceil
	\lambda_{q+1}^{4\alpha}
	\lambda_q\delta_q^{\frac12}/2
	\right\rceil
	=
	2\left\lceil
	\lambda_{q+1}^{4\alpha}\Lambda_q/2
	\right\rceil
	\in2\N.
	\label{def:tau-q}
\end{equation}
By definition,
$
\lambda_{q+1}^{4\alpha}\Lambda_q
\leq
\tau_q^{-1}
<
\lambda_{q+1}^{4\alpha}\Lambda_q+2.
$
Consequently,
$
\tau_q^{-1}
\sim
\lambda_{q+1}^{4\alpha}\lambda_q^{1-\beta}
\sim
\lambda_q^{1-\beta+4b\alpha}\leq\lambda_q^{1-\frac23\beta}.
$
For sufficiently large $a$, the following scale comparisons hold uniformly in $q$.
Since $4b\alpha<(b-1)/6<\beta$, we have
$\tau_q^{-1}\ll\lambda_q\leq\ell_q^{-1}$.
The same parameter inequalities give
\begin{equation}
	\ell_{t,q}\ll\tau_q,
	\qquad
	\ell_{t,q}^{-1}
	=\delta_q^{\frac12}\lambda_q^{\frac12}\lambda_{q+1}^{\frac12}
	\ll\Lambda_{q+1}\ll\lambda_{q+1}.
	\label{eq:material-mollification-hierarchy}
\end{equation}
Indeed,
$$
\frac{\ell_{t,q}}{\tau_q}\sim\lambda_{q+1}^{4\alpha}\left(\frac{\lambda_q}{\lambda_{q+1}}\right)^{\frac12}\ll1,
\qquad
\frac{\ell_{t,q}^{-1}}{\Lambda_{q+1}}=\left(\frac{\lambda_q}{\lambda_{q+1}}\right)^{\frac12-\beta}\ll1.
$$
We also use
\begin{equation}\label{eq:material-mollification-spatial-comparisons}
	\delta_q^{\frac12}\ell_q^{-1}\ell_{t,q}=\left(\frac{\lambda_q}{\lambda_{q+1}}\right)^{\frac14}\ll1,
	\qquad
	\frac{\ell_{t,q}\ell_q^{-1}}{\tau_q\lambda_q}\sim\lambda_{q+1}^{4\alpha}\left(\frac{\lambda_q}{\lambda_{q+1}}\right)^{\frac14}\ll1.
\end{equation}
Here both inequalities follow from $4b\alpha<(b-1)/6$.
For $a>a_{\mathrm{mol}}$, Proposition~\ref{prop:uniform-mollification-estimates} and
\eqref{def:tau-q} give
\begin{equation}
	\tau_q\|\nabla\overline{u}_q\|_N
	\leq C
	\lambda_{q+1}^{-4\alpha}\lambda_q^N
	\leq
	\lambda_q^{N-4\alpha},
	\qquad
	0\leq N\leq N_{\mathrm{ind}}-1.
	\label{eq:geometry-condition-cutoff}
\end{equation}
Here $C$ is at least as large as the implicit constant in \eqref{est:uniform-grad-om-over-orho} with $M=0$ in the stated range of $N$.
By \eqref{def:lambda-delta} and $b>1$, increase $a_{\mathrm{mol}}$ so that $C(\lambda_q/\lambda_{q+1})^{4\alpha}\leq1$ for every $q$.

The spatial partition is adapted from the partitions of unity in \cite{DK22,GK22},
while the temporal partition and enlarged temporal cutoffs are adapted from
\cite{GR24}. The enlarged temporal cutoffs are used to
glue the temporally localized perturbations in the Newton step. The Nash step uses a
separate spatial partition.

For $x\in\R^3$ and $r>0$, let
$
Q(x,r):=\{y\in\R^3:|y-x|_\infty<r\},
$
where $|\cdot|_\infty$ denotes the supremum norm on $\R^3$. Choose a
nonnegative function $\chi_0\in C_c^\infty(\R^3)$ supported in
$Q(0,9\pi/8)$ and equal to $1$ on
$\overline{Q(0,7\pi/8)}$, such that
$
\sum_{\kappa\in\Z^3}\chi_0^2(x-2\pi\kappa)=1
$
for every $x\in\R^3$. For each $\kappa\in\Z^3$, define
$
x_\kappa:=2\pi\tau_q\kappa
$
and
$
\chi_\kappa(x):=\chi_0(\tau_q^{-1}x-2\pi\kappa).
$
Then
$$
\sum_{\kappa\in\Z^3}\chi_\kappa^2(x)=1,
\qquad
x\in\R^3.
$$
Since $\tau_q^{-1}\in 2\N$, translations by elements of
$2\pi\Z^3$ permute the family $\{\chi_\kappa\}_{\kappa\in\Z^3}$, so this family
defines a partition on the periodic domain.
Moreover,
$$
\supp\chi_\kappa
\subset
Q\left(x_\kappa,\frac{9\pi}{8}\tau_q\right),
\qquad
\chi_\kappa=1
\quad\text{on}\quad
\overline{Q\left(x_\kappa,\frac{7\pi}{8}\tau_q\right)}.
$$
The family $\{\chi_\kappa\}_{\kappa\in\Z^3}$ is locally finite and has
uniformly bounded overlap. For every integer $N\geq0$,
$$
\|\chi_\kappa\|_N
\lesssim_N
\tau_q^{-N},
\qquad
\kappa\in\Z^3.
$$

For the temporal partition, choose a nonnegative function
$\theta_0\in C_c^\infty(\R)$ satisfying
$$
\theta_0=1
\quad\text{on}\quad
\left[-\frac13,\frac13\right],
\qquad
\supp\theta_0
\subset
\left(-\frac23,\frac23\right),
$$
and normalized so that
$
\sum_{u\in\Z}\theta_0^2(t-u)=1
$
for every $t\in\R$. For each $u\in\Z$, define
$
t_u:=u\tau_q
$
and
$
\theta_u(t):=\theta_0(\tau_q^{-1}t-u).
$
Then
$$
\sum_{u\in\Z}\theta_u^2(t)=1,
\qquad
t\in\R.
$$
The support properties are
$$
\supp\theta_u
\subset
\left(t_u-\frac23\tau_q,\,
t_u+\frac23\tau_q\right),
\qquad
\theta_u=1
\quad\text{on}\quad
\left[t_u-\frac13\tau_q,\,
t_u+\frac13\tau_q\right].
$$
In particular,
$
\supp\theta_{u-1}\cap\supp\theta_{u+1}=\varnothing
$
for every $u\in\Z$. Moreover, for every integer $N\geq0$,
\begin{equation}\label{eq:time-cutoff-derivative-estimate}
\|\partial_t^N\theta_u\|_0
\lesssim_N
\tau_q^{-N},
\qquad
u\in\Z.
\end{equation}

For the enlarged temporal cutoffs, let
$\tilde\theta_0\in C_c^\infty(\R)$ be nonnegative and satisfy
$$
\tilde\theta_0=1
\quad\text{on}\quad
\left[-\frac23,\frac23\right],
\qquad
\supp\tilde\theta_0
\subset
\left(-\frac{17}{24},\frac{17}{24}\right),
$$
and set
$
\tilde\theta_u(t):=\tilde\theta_0(\tau_q^{-1}t-u).
$
Since $\theta_0$ is supported in $\left(-2/3,2/3\right)$, one has
$
\theta_u\tilde\theta_u=\theta_u
$
for every $u\in\Z$. Furthermore,
$$
\supp\tilde\theta_u
\subset
\left(t_u-\frac{17}{24}\tau_q,\,
t_u+\frac{17}{24}\tau_q\right),
\qquad
\tilde\theta_u=1
\quad\text{on}\quad
\left[t_u-\frac23\tau_q,\,
t_u+\frac23\tau_q\right].
$$
In particular,
$
\supp\partial_t\tilde\theta_u\cap\supp\theta_u=\varnothing.
$
For every integer $N\geq0$,
$$
\|\partial_t^N\tilde\theta_u\|_0
\lesssim_N
\tau_q^{-N},
\qquad
u\in\Z.
$$
The cutoff estimates are uniform in $q$, $u$, and
$\kappa$.

Finally, we define the space-time index set
$
\mathscr I:=\Z\times\Z^3.
$
The spatial and temporal partitions give
$$
\sum_{(u,\kappa)\in\mathscr I}
\theta_u^2(t)\chi_\kappa^2(x)=1,
\qquad
(t,x)\in\R\times\R^3.
$$
The dependence of $\chi_\kappa$, $\theta_u$, $\tilde\theta_u$,
$x_\kappa$, and $t_u$ on the fixed stage $q$ is suppressed from the
notation.

\subsection{Forward and backward flow maps}
\label{subsec:forward-backward-flows}
Let $\upsilon$ be a smooth spatially periodic velocity field on
$\mcI\times\R^3$. Its forward flow map $X_\upsilon(s,y;t)$ is defined by
\begin{equation}\label{eq of forward flow}
	\left\lbrace
	\begin{aligned}
		&\partial_sX_\upsilon(s,y;t)
		=
		\upsilon(s,X_\upsilon(s,y;t)),\\
		&X_\upsilon(t,y;t)=y.
	\end{aligned}
	\right.
\end{equation}
The corresponding backward flow map $\Phi_\upsilon(s,x;t)$ is
characterized by
\begin{equation}\label{backward flow}
	\left\lbrace
	\begin{aligned}
		&\partial_s\Phi_\upsilon(s,x;t)
		+
		\upsilon(s,x)\cdot\nabla\Phi_\upsilon(s,x;t)
		=0,\\
		&\Phi_\upsilon(t,x;t)=x.
	\end{aligned}
	\right.
\end{equation}
Both maps are defined for $s,t\in\mcI$, and
$X_\upsilon(s,\cdot;t)$ and $\Phi_\upsilon(s,\cdot;t)$ are inverse maps. Thus,
\begin{equation}\label{for and back}
	X_\upsilon(s,\Phi_\upsilon(s,x;t);t)=x,
	\qquad
	\Phi_\upsilon(s,X_\upsilon(s,y;t);t)=y.
\end{equation}
For each time center $t_u=u\tau_q$ from
Subsection~\ref{subsec:cutoffs} such that
$[t_u-\tau_q,t_u+\tau_q]\subset\mcI$, we write
\begin{equation}
	\label{def of Xu and Phiu}
	X_u(t,x)
	:=
	X_{\overline{u}_q}(t,x;t_u),
	\qquad
	\Phi_u(t,x)
	:=
	\Phi_{\overline{u}_q}(t,x;t_u)
	=
	X_u(t,\cdot)^{-1}(x).
\end{equation}
For a velocity field $\upsilon$ and a fixed base time $t_0$, set
$$
\JX(t,\cdot)
:=
\det\nabla X_\upsilon(t,\cdot;t_0),
\qquad
\JPhi(t,\cdot)
:=
\det\nabla\Phi_\upsilon(t,\cdot;t_0).
$$
Their dependence on $\upsilon$ and $t_0$ is suppressed below. We also write
$
D_{t,\upsilon}:=\partial_t+\upsilon\cdot\nabla.
$

\begin{lem}
	\label{lem:flow-estimates-upsilon}
	Let $\mcI\subset\mathbb R$ be an interval and
	$\upsilon:\mcI\times\T^3\to\R^3$ a smooth vector field.
	Fix $K^*,K^\sharp\in\mathbb N$, $N_t\in\mathbb N_0$, and $0<\alpha<1$.
	Let $\mC_\upsilon,\mC_{\Div,\upsilon}\geq0$ and choose 
	$
	1\leq\lambda_\upsilon\leq\tilde\lambda_\upsilon,
	\
	1\leq\mu_\upsilon\leq\tilde\mu_\upsilon.
	$
	For $M,N\in\mathbb N_0$ in the indicated ranges, assume
	\begin{align}
		\nrm{D_{t,\upsilon}^{M}\nabla\upsilon}_{N,\mcI}
		&\lesssim \mC_\upsilon\MM{N+1,K^*-M,\lambda_\upsilon,\tilde\lambda_\upsilon} \MM{M,N_t,\mu_\upsilon,\tilde\mu_\upsilon},
		&& 0\leq M\leq K^*,\hspace{18pt}  M+N\leq K^\sharp-1,
		\label{eq:nabla-upsilon-material-assumption}\\
		\nrm{D_{t,\upsilon}^{M}\Div\upsilon}_{N,\mcI}
		&\lesssim \mC_{\Div,\upsilon}\MM{N+1,K^*-M,\lambda_\upsilon,\tilde\lambda_\upsilon}\MM{M,N_t,\mu_\upsilon,\tilde\mu_\upsilon},
		&& 0\leq M\leq K^*-1,\  M+N\leq K^\sharp-1.
		\label{eq:div-upsilon-material-assumption}
	\end{align}
	Assume also the corresponding H{\"o}lder estimates
	\begin{align}
		\nrm{D_{t,\upsilon}^{M}\nabla\upsilon}_{N+\alpha,\mcI}
		&\lesssim \mC_\upsilon\La{N+1,K^*-M,\lambda_\upsilon,\tilde\lambda_\upsilon}  \notag\\
		&\quad\cdot \MM{N+1,K^*-M,\lambda_\upsilon,\tilde\lambda_\upsilon}\MM{M,N_t,\mu_\upsilon,\tilde\mu_\upsilon},
		\quad 0\leq M\leq K^*,\hspace{18pt} M+N\leq K^\sharp-2,
		\label{eq:nabla-upsilon-material-assumption-Calpha}\\
		\nrm{D_{t,\upsilon}^{M}\Div\upsilon}_{N+\alpha,\mcI}
		&\lesssim \mC_{\Div,\upsilon}\La{N+1,K^*-M,\lambda_\upsilon,\tilde\lambda_\upsilon}  \notag\\
		&\quad\cdot \MM{N+1,K^*-M,\lambda_\upsilon,\tilde\lambda_\upsilon}\MM{M,N_t,\mu_\upsilon,\tilde\mu_\upsilon},
		\quad 0\leq M\leq K^*-1,\ M+N\leq K^\sharp-2.
		\label{eq:div-upsilon-material-assumption-Calpha}
	\end{align}
	Set
	$
	\mu_{\upsilon,*}
	:=
	\max\left\{
	\mu_\upsilon,
	\mC_\upsilon\lambda_\upsilon
	\right\},
	\
	\tilde\mu_{\upsilon,*}
	:=
	\max\left\{
	\tilde\mu_\upsilon,
	\mC_\upsilon\lambda_\upsilon
	\right\}.
	$
	Assume moreover that, for some $t_0\in\mcI$ and $\tau>0$,
	$$
	[t_0-\tau,t_0+\tau]\subset\mcI,
	\qquad
	\tau\mC_\upsilon\lambda_\upsilon\leq c_0,
	$$
	where $c_0>0$ is sufficiently small in terms of the implicit constant
	in \eqref{eq:nabla-upsilon-material-assumption} with $M=N=0$. Let
	$X_\upsilon(t,y;t_0)$ and $\Phi_\upsilon(t,x;t_0)$ be defined by
	\eqref{eq of forward flow} and \eqref{backward flow}. Then, for
	$|t-t_0|\leq\tau$, one has
	\begin{align}
		\nrm{\nabla\Phi_\upsilon(t,\cdot;t_0)-\Id}_{0}
		+
		\nrm{\nabla X_\upsilon(t,\cdot;t_0)-\Id}_{0}
		&\lesssim
		\tau\mC_\upsilon\lambda_\upsilon,
		\label{eq:flow-gradient-close-id}
		\\
		\nrm{\JPhi(t,\cdot)-1}_{0}
		+
		\nrm{\JX(t,\cdot)-1}_{0}
		&\lesssim
		\tau\mC_{\Div,\upsilon}\lambda_\upsilon.
		\label{eq:flow-jacobian-close-one}
	\end{align}
	The constant $c_0$ is chosen so that the left-hand side of
	\eqref{eq:flow-gradient-close-id} is at most $1/2$.
	Moreover, one has
	\begin{align}
		&\nrm{\nabla\Phi_\upsilon(t,\cdot;t_0)}_N+\nrm{\bigl(\nabla\Phi_\upsilon\bigr)^{-1}(t,\cdot;t_0)}_N+\nrm{\nabla X_\upsilon(t,\cdot;t_0)}_N \lesssim\MM{N,K^*-1,\lambda_\upsilon,\tilde\lambda_\upsilon},
		&& 0\leq N\leq K^\sharp-1,
		\label{eq:flow-gradient-CN}\\
		&\nrm{\JPhi(t,\cdot)}_N+\nrm{(\JPhi)^{-1}(t,\cdot)}_N+\nrm{\JX(t,\cdot)}_N+\nrm{(\JX)^{-1}(t,\cdot)}_N\lesssim\MM{N,K^*-1,\lambda_\upsilon,\tilde\lambda_\upsilon},
		&& 0\leq N\leq K^\sharp-1.
		\label{eq:flow-jacobian-CN}
	\end{align}
	The corresponding H{\"o}lder estimates are
	\begin{align}
		&\nrm{\nabla\Phi_\upsilon(t,\cdot;t_0)}_{N+\alpha}+\nrm{\bigl(\nabla\Phi_\upsilon\bigr)^{-1}(t,\cdot;t_0)}_{N+\alpha}+\nrm{\nabla X_\upsilon(t,\cdot;t_0)}_{N+\alpha} \notag\\
		&\quad\lesssim\La{N+1,K^*,\lambda_\upsilon,\tilde\lambda_\upsilon}\MM{N,K^*-1,\lambda_\upsilon,\tilde\lambda_\upsilon},
		&& 0\leq N\leq K^\sharp-2,
		\label{eq:flow-gradient-Calpha}\\
		&\nrm{\JPhi(t,\cdot)}_{N+\alpha}+\nrm{(\JPhi)^{-1}(t,\cdot)}_{N+\alpha}+\nrm{\JX(t,\cdot)}_{N+\alpha}+\nrm{(\JX)^{-1}(t,\cdot)}_{N+\alpha} \notag\\
		&\quad\lesssim\La{N+1,K^*,\lambda_\upsilon,\tilde\lambda_\upsilon}\MM{N,K^*-1,\lambda_\upsilon,\tilde\lambda_\upsilon},
		&& 0\leq N\leq K^\sharp-2.
		\label{eq:flow-jacobian-Calpha}
	\end{align}
	The following derivative estimates hold for $\sigma\in\mathbb R$.  
	\begin{align}
		&\nrm{D_{t,\upsilon}^M\nabla\Phi_\upsilon(t,\cdot;t_0)}_{N}+\nrm{D_{t,\upsilon}^M(\nabla\Phi_\upsilon)^{-1}(t,\cdot;t_0)}_{N}+\nrm{\partial_t^M\nabla X_\upsilon(t,\cdot;t_0)}_{N} \notag\\
		&\quad\lesssim \mC_\upsilon\MM{N+1,K^*-M+1,\lambda_\upsilon,\tilde\lambda_\upsilon} \MM{M-1,N_t,\mu_{\upsilon,*},\tilde\mu_{\upsilon,*}},
		&& 1\leq M\leq K^*+1,\  M+N\leq K^\sharp,
		\label{eq:flow-gradient-material}\\
		&\nrm{D_{t,\upsilon}^M\JPhi^\sigma(t,\cdot)}_{N}+\nrm{\partial_t^M\JX^\sigma(t,\cdot)}_{N} \notag\\
		&\quad\lesssim \mC_{\Div,\upsilon}\MM{N+1,K^*-M+1,\lambda_\upsilon,\tilde\lambda_\upsilon} \MM{M-1,N_t,\mu_{\upsilon,*},\tilde\mu_{\upsilon,*}},
		&& 1\leq M\leq K^*,\quad M+N\leq K^\sharp.
		\label{eq:flow-Jacobian-power-material}
	\end{align}
	The corresponding H{\"o}lder estimates are
	\begin{align}
		&\nrm{D_{t,\upsilon}^M\nabla\Phi_\upsilon(t,\cdot;t_0)}_{N+\alpha}+\nrm{D_{t,\upsilon}^M(\nabla\Phi_\upsilon)^{-1}(t,\cdot;t_0)}_{N+\alpha}+\nrm{\partial_t^M\nabla X_\upsilon(t,\cdot;t_0)}_{N+\alpha} \notag\\
		&\quad\lesssim \mC_\upsilon\La{N+1,K^*-M+1,\lambda_\upsilon,\tilde\lambda_\upsilon}  \notag\\
		&\hspace{30pt}\cdot \MM{N+1,K^*-M+1,\lambda_\upsilon,\tilde\lambda_\upsilon}\MM{M-1,N_t,\mu_{\upsilon,*},\tilde\mu_{\upsilon,*}},
		\quad 1\leq M\leq K^*+1,\quad M+N\leq K^\sharp-1,
		\label{eq:flow-gradient-material-Calpha}\\
		&\nrm{D_{t,\upsilon}^M\JPhi^\sigma(t,\cdot)}_{N+\alpha}+\nrm{\partial_t^M\JX^\sigma(t,\cdot)}_{N+\alpha} \notag\\
		&\quad\lesssim \mC_{\Div,\upsilon}\La{N+1,K^*-M+1,\lambda_\upsilon,\tilde\lambda_\upsilon} \notag\\
		&\hspace{30pt}\cdot \MM{N+1,K^*-M+1,\lambda_\upsilon,\tilde\lambda_\upsilon} \MM{M-1,N_t,\mu_{\upsilon,*},\tilde\mu_{\upsilon,*}},
		\quad 1\leq M\leq K^*,\hspace{24pt}M+N\leq K^\sharp-1.
		\label{eq:flow-Jacobian-power-material-Calpha}
	\end{align}
	Estimates \eqref{eq:flow-Jacobian-power-material} and
	\eqref{eq:flow-Jacobian-power-material-Calpha} also hold for
	$1\leq M\leq K^*+1$ with $\mC_{\Div,\upsilon}$ replaced by
	$\mC_\upsilon$, for $M+N\leq K^\sharp$ and
	$M+N\leq K^\sharp-1$, respectively.
	The implicit constants depend only on $K^*,K^\sharp$,
	$\sigma$, $\alpha$, and the constants in
	\eqref{eq:nabla-upsilon-material-assumption}--\eqref{eq:div-upsilon-material-assumption-Calpha}.
\end{lem}

\begin{proof}
	For compositions, use Lemma~\ref{lem:BDSV-composition} for the $C^N$ estimates when $N\geq1$ and Lemma~\ref{lem:holder-composition-general} for the H{\"o}lder estimates, including $N=0$. 
	For a vector-valued function $f$, the commutator identity is
	$
	[D_{t,\upsilon},\nabla]f=-\nabla f\nabla\upsilon
	$.
	Since $D_{t,\upsilon}\Phi_\upsilon=0$, the commutator identity implies
	$
	D_{t,\upsilon}\nabla\Phi_\upsilon
	=
	-\nabla\Phi_\upsilon\nabla\upsilon.
	$
	Applying \eqref{transport equation 0} to this equation and using
	\eqref{eq:nabla-upsilon-material-assumption} with $M=N=0$, we obtain
	$$
	\nrm{\nabla\Phi_\upsilon(t,\cdot;t_0)}_0
	\leq
	1+
	C\int_{\min\{t,t_0\}}^{\max\{t,t_0\}}
	\nrm{\nabla\Phi_\upsilon(s,\cdot;t_0)}_0
	\nrm{\nabla\upsilon(s,\cdot)}_0\,\rd s.
	$$
	Gronwall's inequality and
	$\tau\mC_\upsilon\lambda_\upsilon\leq c_0$ imply
	$$
	\nrm{\nabla\Phi_\upsilon(t,\cdot;t_0)}_0
	\lesssim
	1.
	$$
	Subtracting $\Id$ from the integral equation for
	$\nabla\Phi_\upsilon$ yields
	$$
	\nrm{\nabla\Phi_\upsilon(t,\cdot;t_0)-\Id}_0
	\lesssim
	\int_{\min\{t,t_0\}}^{\max\{t,t_0\}}
	\nrm{\nabla\Phi_\upsilon(s,\cdot;t_0)}_0
	\nrm{\nabla\upsilon(s,\cdot)}_0\,\rd s
	\lesssim
	\tau\mC_\upsilon\lambda_\upsilon.
	$$
	By choosing $c_0$ sufficiently small, we ensure that
	$\nabla\Phi_\upsilon$ is uniformly invertible and
	$$
	\nrm{(\nabla\Phi_\upsilon)^{-1}(t,\cdot;t_0)}_0
	\lesssim
	1,
	\qquad
	\nrm{(\nabla\Phi_\upsilon)^{-1}(t,\cdot;t_0)-\Id}_0
	\lesssim
	\tau\mC_\upsilon\lambda_\upsilon.
	$$
	Also,
	$
	D_{t,\upsilon}(\nabla\Phi_\upsilon)^{-1}
	=
	\nabla\upsilon(\nabla\Phi_\upsilon)^{-1}.
	$
	
	For $1\leq N\leq K^\sharp-1$, we apply the $C^N$ estimate in
	Lemma~\ref{estimate of the transport equation} to
	$
	D_{t,\upsilon}\nabla\Phi_\upsilon
	=
	-\nabla\Phi_\upsilon\nabla\upsilon
	$
	to obtain
	$$
	\begin{aligned}
		\nrm{\nabla\Phi_\upsilon(t,\cdot;t_0)}_N
		\lesssim 1+\int_{\min\{t,t_0\}}^{\max\{t,t_0\}}\nrm{\nabla\upsilon(s,\cdot)}_0\nrm{\nabla\Phi_\upsilon(s,\cdot;t_0)}_N\,\rd s+\sum_{N_1=1}^{N}\int_{\min\{t,t_0\}}^{\max\{t,t_0\}}\nrm{\nabla\upsilon(s,\cdot)}_{N_1}\nrm{\nabla\Phi_\upsilon(s,\cdot;t_0)}_{N-N_1}\,\rd s.
	\end{aligned}
	$$
	Using \eqref{eq:nabla-upsilon-material-assumption}, the estimate for $\nabla\Phi_\upsilon$ in
	\eqref{eq:flow-gradient-CN} at $N=0$, proved above, and induction on $N$, we
	obtain
	$$
	\nrm{\nabla\Phi_\upsilon(t,\cdot;t_0)}_N
	\lesssim 1+
	\int_{\min\{t,t_0\}}^{\max\{t,t_0\}}
	\mC_\upsilon\lambda_\upsilon
	\nrm{\nabla\Phi_\upsilon(s,\cdot;t_0)}_N\,\rd s
	+
	\tau\mC_\upsilon
	\MM{N+1,K^*,\lambda_\upsilon,\tilde\lambda_\upsilon}.
	$$
	Since $K^*\geq1$,
	$\MM{N+1,K^*,\lambda_\upsilon,\tilde\lambda_\upsilon}=\lambda_\upsilon\MM{N,K^*-1,\lambda_\upsilon,\tilde\lambda_\upsilon}$.
	The condition $\tau\mC_\upsilon\lambda_\upsilon\leq c_0$ and
	Gronwall's inequality therefore give
	$$
	\nrm{\nabla\Phi_\upsilon(t,\cdot;t_0)}_N
	\lesssim
	\MM{N,K^*-1,\lambda_\upsilon,\tilde\lambda_\upsilon}.
	$$
	Applying the same argument to
	$
	D_{t,\upsilon}(\nabla\Phi_\upsilon)^{-1}
	=
	\nabla\upsilon(\nabla\Phi_\upsilon)^{-1}
	$
	gives
	$$
	\nrm{(\nabla\Phi_\upsilon)^{-1}(t,\cdot;t_0)}_N
	\lesssim
	\MM{N,K^*-1,\lambda_\upsilon,\tilde\lambda_\upsilon}.
	$$
	For $0\leq N\leq K^\sharp-2$, the corresponding $C^{N+\alpha}$ estimates
	follow from the same argument, using
	\eqref{eq:nabla-upsilon-material-assumption-Calpha}:
	$$
	\begin{aligned}
		\nrm{\nabla\Phi_\upsilon(t,\cdot;t_0)}_{N+\alpha}
		+\nrm{(\nabla\Phi_\upsilon)^{-1}(t,\cdot;t_0)}_{N+\alpha}\lesssim\La{N+1,K^*,\lambda_\upsilon,\tilde\lambda_\upsilon}\MM{N,K^*-1,\lambda_\upsilon,\tilde\lambda_\upsilon}.
	\end{aligned}
	$$

	For $1\leq M\leq K^*+1$, the transport equation gives
	$$
	D_{t,\upsilon}^{M}\nabla\Phi_\upsilon
	=
	-\sum_{r=0}^{M-1}\binom{M-1}{r}
	\bigl(D_{t,\upsilon}^{r}\nabla\Phi_\upsilon\bigr)
	\bigl(D_{t,\upsilon}^{M-1-r}\nabla\upsilon\bigr).
	$$
	After taking $N$ spatial derivatives, each velocity-gradient factor has material order at most $M-1\leq K^*$. The sum of this order and the number of additional spatial derivatives is at most $M+N-1$. Thus
	\eqref{eq:nabla-upsilon-material-assumption} applies when
	$M+N\leq K^\sharp$. In the H{\"o}lder range
	$M+N\leq K^\sharp-1$, the corresponding sum for each velocity-gradient factor is at most $M+N-1\leq K^\sharp-2$, so that
	\eqref{eq:nabla-upsilon-material-assumption-Calpha} applies.
	
	We argue by induction on $M$. For $M=1$, the recurrence and the spatial
	estimates for $\nabla\Phi_\upsilon$ give
	$$
	\nrm{D_{t,\upsilon}\nabla\Phi_\upsilon}_N
	\lesssim
	\mC_\upsilon
	\MM{N+1,K^*,\lambda_\upsilon,\tilde\lambda_\upsilon}.
	$$
	Let $2\leq M\leq K^*+1$ and assume that
	$$
	\nrm{D_{t,\upsilon}^{r}\nabla\Phi_\upsilon}_N
	\lesssim
	\mC_\upsilon
	\MM{N+1,K^*-r+1,\lambda_\upsilon,\tilde\lambda_\upsilon}
	\MM{r-1,N_t,\mu_{\upsilon,*},\tilde\mu_{\upsilon,*}}
	$$
	for $1\leq r<M$ and $N\in\mathbb N_0$ with $r+N\leq K^\sharp$.
	For $r=0$, \eqref{eq:nabla-upsilon-material-assumption} and the spatial estimates give
	$$
	\nrm{\nabla\Phi_\upsilon
		D_{t,\upsilon}^{M-1}\nabla\upsilon}_N
	\lesssim
	\mC_\upsilon
	\MM{N+1,K^*-M+1,\lambda_\upsilon,\tilde\lambda_\upsilon}
	\MM{M-1,N_t,\mu_{\upsilon,*},\tilde\mu_{\upsilon,*}}.
	$$
	For $1\leq r\leq M-1$,
	$
	\min\{K^*-r+1,K^*-M+r+1\}
	\geq K^*-M+2.
	$
	Using the induction hypothesis,
	\eqref{eq:nabla-upsilon-material-assumption}, and \eqref{eq:holder-product-tame} with $\alpha=0$,
	we obtain
	$$
	\begin{aligned}
		&\nrm{
			\bigl(D_{t,\upsilon}^{r}\nabla\Phi_\upsilon\bigr)
			\bigl(D_{t,\upsilon}^{M-1-r}\nabla\upsilon\bigr)
		}_N\\
		&\quad\lesssim
		\mC_\upsilon^2\lambda_\upsilon
		\MM{N+1,K^*-M+1,\lambda_\upsilon,\tilde\lambda_\upsilon}
		\MM{r-1,N_t,\mu_{\upsilon,*},\tilde\mu_{\upsilon,*}}
		\MM{M-1-r,N_t,\mu_{\upsilon,*},\tilde\mu_{\upsilon,*}}.
	\end{aligned}
	$$
	By \eqref{eq:appendix-derivative-weight-product} and the definition of $\mu_{\upsilon,*}$,
	$$
	\begin{aligned}
		\nrm{
			\bigl(D_{t,\upsilon}^{r}\nabla\Phi_\upsilon\bigr)
			\bigl(D_{t,\upsilon}^{M-1-r}\nabla\upsilon\bigr)
		}_N\lesssim
		\mC_\upsilon
		\MM{N+1,K^*-M+1,\lambda_\upsilon,\tilde\lambda_\upsilon}
		\MM{M-1,N_t,\mu_{\upsilon,*},\tilde\mu_{\upsilon,*}}.
	\end{aligned}
	$$
	Here we used $ \MM{N+2,K^*-M+2,\lambda_\upsilon,\tilde\lambda_\upsilon} = \lambda_\upsilon \MM{N+1,K^*-M+1,\lambda_\upsilon,\tilde\lambda_\upsilon}. $
	Summing over $r$ yields
	$$
	\nrm{D_{t,\upsilon}^{M}\nabla\Phi_\upsilon(t,\cdot;t_0)}_N
	\lesssim
	\mC_\upsilon
	\MM{N+1,K^*-M+1,\lambda_\upsilon,\tilde\lambda_\upsilon}
	\MM{M-1,N_t,\mu_{\upsilon,*},\tilde\mu_{\upsilon,*}}.
	$$
	The identity
	$
	D_{t,\upsilon}(\nabla\Phi_\upsilon)^{-1}
	=
	\nabla\upsilon(\nabla\Phi_\upsilon)^{-1}
	$
	gives the estimate for $(\nabla\Phi_\upsilon)^{-1}$ by the same induction.
	With \eqref{eq:nabla-upsilon-material-assumption-Calpha} and \eqref{eq:holder-product-finite-Cnalpha}, the induction also proves
	$$
	\begin{aligned}
		&\nrm{D_{t,\upsilon}^{M}\nabla\Phi_\upsilon(t,\cdot;t_0)}_{N+\alpha}
		+
		\nrm{D_{t,\upsilon}^{M}
			(\nabla\Phi_\upsilon)^{-1}(t,\cdot;t_0)}_{N+\alpha}\\
		&\quad\lesssim
		\mC_\upsilon
		\La{N+1,K^*-M+1,\lambda_\upsilon,\tilde\lambda_\upsilon}
		\MM{N+1,K^*-M+1,\lambda_\upsilon,\tilde\lambda_\upsilon}
		\MM{M-1,N_t,\mu_{\upsilon,*},\tilde\mu_{\upsilon,*}}.
	\end{aligned}
	$$

	For the forward flow, differentiating
	$
	\partial_tX_\upsilon(t,y;t_0)
	=
	\upsilon(t,X_\upsilon(t,y;t_0))
	$
	with respect to $y$ yields
	\begin{align}
		\partial_t\nabla X_\upsilon(t,y;t_0)
		&=
		(\nabla\upsilon)(t,X_\upsilon(t,y;t_0))
		\nabla X_\upsilon(t,y;t_0),
		\qquad
		\nabla X_\upsilon(t_0,y;t_0)=\Id.
		\label{pt nabla X}
	\end{align}
	Hence
	$$
	\nrm{\nabla X_\upsilon(t,\cdot;t_0)}_0
	\leq
	1+
	C\int_{\min\{t,t_0\}}^{\max\{t,t_0\}}
	\nrm{\nabla\upsilon(s,\cdot)}_0
	\nrm{\nabla X_\upsilon(s,\cdot;t_0)}_0\,\rd s,
	$$
	and Gronwall's inequality gives
	$
	\nrm{\nabla X_\upsilon(t,\cdot;t_0)}_0
	\lesssim
	1.
	$
	Subtracting $\Id$ from the integral form of \eqref{pt nabla X},
	we obtain
	$$
	\nrm{\nabla X_\upsilon(t,\cdot;t_0)-\Id}_0
	\lesssim
	\int_{\min\{t,t_0\}}^{\max\{t,t_0\}}
	\nrm{\nabla\upsilon(s,\cdot)}_0
	\nrm{\nabla X_\upsilon(s,\cdot;t_0)}_0\,\rd s
	\lesssim
	\tau\mC_\upsilon\lambda_\upsilon.
	$$
	For $1\leq N\leq K^\sharp-1$, let $|\zeta|=N$. Applying $\partial^\zeta$
	to \eqref{pt nabla X} gives
	$$
	\begin{aligned}
		\partial_t\partial^\zeta\nabla X_\upsilon
		=
		(\nabla\upsilon)(t,X_\upsilon)
		\partial^\zeta\nabla X_\upsilon+
		\sum_{\substack{\zeta_1+\zeta_2=\zeta,\zeta_1\neq0}}
		c_{\zeta_1,\zeta_2}
		\partial^{\zeta_1}
		\bigl((\nabla\upsilon)(t,X_\upsilon)\bigr)
		\partial^{\zeta_2}\nabla X_\upsilon.
	\end{aligned}
	$$
	We argue by induction on $N$. In each term of the sum, all
	derivatives of $\nabla X_\upsilon$ have order at most $N-1$.
	By \eqref{eq:holder-product-tame} with $\alpha=0$, \eqref{est:holder-norm-composition},
	and the induction hypothesis, we have
	$$
	\begin{aligned}
		\nrm{\nabla X_\upsilon(t,\cdot;t_0)}_N
		&\lesssim 1+
		\int_{\min\{t,t_0\}}^{\max\{t,t_0\}}\nrm{\nabla\upsilon(s,\cdot)}_0\nrm{\nabla X_\upsilon(s,\cdot;t_0)}_N\,\rd s\\
		&\quad+\sum_{N_1=1}^{N}
		\int_{\min\{t,t_0\}}^{\max\{t,t_0\}}\nrm{\nabla\upsilon(s,\cdot)}_{N_1}\MM{N-N_1,K^*-1,\lambda_\upsilon,\tilde\lambda_\upsilon}\,\rd s.
	\end{aligned}
	$$
	Using \eqref{eq:nabla-upsilon-material-assumption},
	\eqref{eq:appendix-derivative-weight-product}, and
	$\tau\mC_\upsilon\lambda_\upsilon\leq c_0$, Gronwall's inequality
	gives
	\begin{equation}\label{eq:forward-flow-gradient-integer-estimate}
	\nrm{\nabla X_\upsilon(t,\cdot;t_0)}_N
	\lesssim
	\MM{N,K^*-1,\lambda_\upsilon,\tilde\lambda_\upsilon}.
	\end{equation}
	The same argument in $C^{N+\alpha}$, using
	\eqref{eq:nabla-upsilon-material-assumption-Calpha}, gives
	\begin{equation}\label{eq:forward-flow-gradient-holder-estimate}
	\nrm{\nabla X_\upsilon(t,\cdot;t_0)}_{N+\alpha}
	\lesssim
	\La{N+1,K^*,\lambda_\upsilon,\tilde\lambda_\upsilon}
	\MM{N,K^*-1,\lambda_\upsilon,\tilde\lambda_\upsilon}.
	\end{equation}
	
	By the inverse relation between $X_\upsilon$ and $\Phi_\upsilon$, we have
	$
	\nabla X_\upsilon(t,y;t_0)
	=
	(\nabla\Phi_\upsilon)^{-1}
	(t,X_\upsilon(t,y;t_0);t_0).
	$
	Hence, for every $M\geq1$,
	\begin{equation}\label{eq:pullback-flow-gradient-material}
		\partial_t^M\nabla X_\upsilon(t,y;t_0)
		=
		\bigl(D_{t,\upsilon}^{M}
		(\nabla\Phi_\upsilon)^{-1}\bigr)
		(t,X_\upsilon(t,y;t_0);t_0).
	\end{equation}
	Applying Lemmas~\ref{lem:BDSV-composition} and \ref{lem:holder-composition-general} to
	\eqref{eq:pullback-flow-gradient-material} and using \eqref{eq:forward-flow-gradient-integer-estimate}--\eqref{eq:forward-flow-gradient-holder-estimate} for
	$X_\upsilon$, we obtain
	$$
	\begin{aligned}
		\nrm{\partial_t^M\nabla X_\upsilon(t,\cdot;t_0)}_N
		&\lesssim
		\mC_\upsilon
		\MM{N+1,K^*-M+1,\lambda_\upsilon,\tilde\lambda_\upsilon}
		\MM{M-1,N_t,\mu_{\upsilon,*},\tilde\mu_{\upsilon,*}},
		\\
		\nrm{\partial_t^M\nabla X_\upsilon(t,\cdot;t_0)}_{N+\alpha}
		&\lesssim
		\mC_\upsilon
		\La{N+1,K^*-M+1,\lambda_\upsilon,\tilde\lambda_\upsilon}\\
		&\quad\cdot 
		\MM{N+1,K^*-M+1,\lambda_\upsilon,\tilde\lambda_\upsilon}
		\MM{M-1,N_t,\mu_{\upsilon,*},\tilde\mu_{\upsilon,*}}.
	\end{aligned}
	$$
	
	For the Jacobians, Jacobi's formula gives
	$
	D_{t,\upsilon}\JPhi
	=
	\JPhi\tr\bigl(
	(\nabla\Phi_\upsilon)^{-1}
	D_{t,\upsilon}\nabla\Phi_\upsilon
	\bigr)
	=
	-(\Div\upsilon)\JPhi.
	$
	Since $\JPhi(t_0,\cdot)=1$, we integrate along trajectories to obtain
	$$
	\JPhi(t,X_\upsilon(t,y;t_0))
	=\exp\left(-\int_{t_0}^t(\Div\upsilon)(s,X_\upsilon(s,y;t_0))\,\rd s\right).
	$$
	In particular,
	$
	\JPhi(t,\cdot)>0
	$
	for every $|t-t_0|\leq\tau$. Hence all real powers of $\JPhi$ are
	well-defined.
	For every fixed $\sigma\in\mathbb R$, it follows that
	\begin{align}
		D_{t,\upsilon}\JPhi^\sigma
		=
		\sigma\JPhi^{\sigma-1}D_{t,\upsilon}\JPhi
		=
		-\sigma(\Div\upsilon)\JPhi^\sigma.
		\label{Dt Jm}
	\end{align}
	Applying \eqref{transport equation 0} to \eqref{Dt Jm} with
	$\sigma=1$ gives
	$$
	\nrm{\JPhi(t,\cdot)}_0
	\leq
	1+
	C\int_{\min\{t,t_0\}}^{\max\{t,t_0\}}
	\nrm{\Div\upsilon(s,\cdot)}_0
	\nrm{\JPhi(s,\cdot)}_0\,\rd s.
	$$
	Since $\nrm{\Div\upsilon}_0\lesssim\nrm{\nabla\upsilon}_0$, Gronwall's
	inequality and the gradient smallness condition give
	$
	\nrm{\JPhi(t,\cdot)}_0\lesssim1.
	$
	Subtracting the initial value from the integral equation yields
	$$
	\nrm{\JPhi(t,\cdot)-1}_0
	\lesssim
	\int_{\min\{t,t_0\}}^{\max\{t,t_0\}}
	\nrm{\Div\upsilon(s,\cdot)}_0
	\nrm{\JPhi(s,\cdot)}_0\,\rd s
	\lesssim
	\tau\mC_{\Div,\upsilon}\lambda_\upsilon.
	$$
	The exponential formula and the gradient smallness condition give
	$1/2\leq\JPhi\leq3/2$ for sufficiently small $c_0$.
	Hence, for every fixed
	$\sigma\in\mathbb R$,
	$$
	\nrm{\JPhi^\sigma(t,\cdot)-1}_0
	\lesssim
	\nrm{\JPhi(t,\cdot)-1}_0
	\lesssim
	\tau\mC_{\Div,\upsilon}\lambda_\upsilon.
	$$
	
	For $1\leq N\leq K^\sharp-1$, we apply the $C^N$ estimate in
	Lemma~\ref{estimate of the transport equation} to
	\eqref{Dt Jm} and estimate the divergence by the gradient to obtain
	$$
	\begin{aligned}
		\nrm{\JPhi^\sigma(t,\cdot)}_N
		&\lesssim1
		+\sum_{N_1=0}^{N}
		\int_{\min\{t,t_0\}}^{\max\{t,t_0\}}
		\nrm{\nabla\upsilon(s,\cdot)}_{N_1}
		\nrm{\JPhi^\sigma(s,\cdot)}_{N-N_1}\,\rd s.
	\end{aligned}
	$$
	Using \eqref{eq:nabla-upsilon-material-assumption} and $\tau\mC_\upsilon\lambda_\upsilon\leq c_0$, induction
	on the spatial derivative order and Gronwall's inequality give
	$$
	\nrm{\JPhi^\sigma(t,\cdot)}_N
	\lesssim
	\MM{N,K^*-1,\lambda_\upsilon,\tilde\lambda_\upsilon}.
	$$
	For $0\leq N\leq K^\sharp-2$, the corresponding $C^{N+\alpha}$ estimate
	follows from the same argument, using
	\eqref{eq:nabla-upsilon-material-assumption-Calpha}:
	$$
	\nrm{\JPhi^\sigma(t,\cdot)}_{N+\alpha}
	\lesssim
	\La{N+1,K^*,\lambda_\upsilon,\tilde\lambda_\upsilon}
	\MM{N,K^*-1,\lambda_\upsilon,\tilde\lambda_\upsilon}.
	$$
	
	For the material derivatives of $\JPhi^\sigma$, applying
	$D_{t,\upsilon}^{M-1}$ to \eqref{Dt Jm} gives
	$$
	D_{t,\upsilon}^{M}\JPhi^\sigma
	=
	-\sigma
	\sum_{M_1+M_2=M-1}\binom{M-1}{M_1}
	\bigl(D_{t,\upsilon}^{M_1}\Div\upsilon\bigr)
	\bigl(D_{t,\upsilon}^{M_2}\JPhi^\sigma\bigr).
	$$
	We proceed by induction on $M$, using the spatial estimate for $\JPhi^\sigma$ when $M=1$.
	Each coefficient $D_{t,\upsilon}^m\Div\upsilon$ is used with
	$m\leq M-1\leq K^*-1$ and $n$ additional spatial derivatives, where
	$m+n\leq M+N-1$. Thus the divergence estimates apply in the stated ranges.
	When $M_2=0$, use the estimate for $D_{t,\upsilon}^{M_1}\Div\upsilon$.
	When $M_2\geq1$, the induction hypothesis supplies the factor $\mC_{\Div,\upsilon}$.
	Use the estimate for $D_{t,\upsilon}^{M_1}\nabla\upsilon$ to control $D_{t,\upsilon}^{M_1}\Div\upsilon$.
	Thus
	$$
	\nrm{D_{t,\upsilon}^{M}\JPhi^\sigma(t,\cdot)}_N
	\lesssim
	\mC_{\Div,\upsilon}
	\MM{N+1,K^*-M+1,\lambda_\upsilon,\tilde\lambda_\upsilon}\allowbreak
	\MM{M-1,N_t,\mu_{\upsilon,*},\tilde\mu_{\upsilon,*}}.
	$$
	The $C^{N+\alpha}$ estimate follows by the same argument, using
	\eqref{eq:nabla-upsilon-material-assumption-Calpha}--%
	\eqref{eq:div-upsilon-material-assumption-Calpha}:
	$$
	\begin{aligned}
		&\nrm{D_{t,\upsilon}^{M}\JPhi^\sigma(t,\cdot)}_{N+\alpha}
		\\
		&\quad\lesssim
		\mC_{\Div,\upsilon}
		\La{N+1,K^*-M+1,\lambda_\upsilon,\tilde\lambda_\upsilon}
		\MM{N+1,K^*-M+1,\lambda_\upsilon,\tilde\lambda_\upsilon}
		\MM{M-1,N_t,\mu_{\upsilon,*},\tilde\mu_{\upsilon,*}}.
	\end{aligned}
	$$
	
	For every fixed $\sigma\in\mathbb R$, the forward and backward Jacobians
	satisfy
	$
	\JX^\sigma(t,y)
	=
	\JPhi^{-\sigma}(t,X_\upsilon(t,y;t_0)).
	$
	Moreover, for every $M\geq0$,
	\begin{equation}
		\label{eq:pullback-Jacobian-material}
		\partial_t^M\JX^\sigma(t,y)
		=
		\bigl(D_{t,\upsilon}^{M}(\JPhi^{-\sigma})\bigr)
		(t,X_\upsilon(t,y;t_0)).
	\end{equation}
	For $M=0$, Lemmas~\ref{lem:BDSV-composition} and \ref{lem:holder-composition-general}, together with the elementary $C^0$ composition estimate and the estimates already
	proved for $\JPhi^{-\sigma}$ and $X_\upsilon$, give
	\begin{align*}
		\nrm{\JX^\sigma(t,\cdot)}_N
		\lesssim
		\MM{N,K^*-1,\lambda_\upsilon,\tilde\lambda_\upsilon},
		\quad
		\nrm{\JX^\sigma(t,\cdot)}_{N+\alpha}
		\lesssim
		\La{N+1,K^*,\lambda_\upsilon,\tilde\lambda_\upsilon}
		\MM{N,K^*-1,\lambda_\upsilon,\tilde\lambda_\upsilon}.
	\end{align*}
	In particular,
	$$
	\nrm{\JX(t,\cdot)-1}_0
	=
	\nrm{
		(\JPhi)^{-1}(t,X_\upsilon(t,\cdot;t_0))-1
	}_0
	\lesssim
	\tau\mC_{\Div,\upsilon}\lambda_\upsilon.
	$$
	For $1\leq M\leq K^*$, apply the same composition estimates to
	\eqref{eq:pullback-Jacobian-material}. The estimates for $D_{t,\upsilon}^M\JPhi^{-\sigma}$ already established in
	\eqref{eq:flow-Jacobian-power-material} and \eqref{eq:flow-Jacobian-power-material-Calpha} give
	$$
	\begin{aligned}
		\nrm{\partial_t^{M}\JX^\sigma(t,\cdot)}_N
		&\lesssim
		\mC_{\Div,\upsilon}
		\MM{N+1,K^*-M+1,\lambda_\upsilon,\tilde\lambda_\upsilon}
		\MM{M-1,N_t,\mu_{\upsilon,*},\tilde\mu_{\upsilon,*}},
		\\
		\nrm{\partial_t^{M}\JX^\sigma(t,\cdot)}_{N+\alpha}
		&\lesssim
		\mC_{\Div,\upsilon}
		\La{N+1,K^*-M+1,\lambda_\upsilon,\tilde\lambda_\upsilon} 
		\MM{N+1,K^*-M+1,\lambda_\upsilon,\tilde\lambda_\upsilon}
		\MM{M-1,N_t,\mu_{\upsilon,*},\tilde\mu_{\upsilon,*}}.
	\end{aligned}
	$$
	Finally, using
	$D_{t,\upsilon}^{r}\Div\upsilon
	=\operatorname{tr}(D_{t,\upsilon}^{r}\nabla\upsilon)$
	for each material order $r$ throughout the Jacobian argument replaces $\mC_{\Div,\upsilon}$ by
	$\mC_\upsilon$ and allows $M-1\leq K^*$.
	This proves the additional Jacobian estimates through $M=K^*+1$.
\end{proof}
\begin{co}
	\label{cor:mollified-velocity-flow-estimates}
	Assume the hypotheses of Proposition~\ref{prop:uniform-mollification-estimates}.
	Let $\mcI\subset\R$ be an interval, let $t_0\in\mcI$, and let $h>0$ satisfy
	$[t_0-h,t_0+h]\subset\mcI$ and
	$h\leq\tau_q$. The flows $X_{\overline u_q}$ and $\Phi_{\overline u_q}$ have
	base time $t_0$, and all norms include the supremum over $|t-t_0|\leq h$.
	There exists $a_{\mathrm{flow}}=a_{\mathrm{flow}}(\beta,\allowbreak b,\allowbreak \alpha,\allowbreak \rho^*,\allowbreak \mathfrak M,\allowbreak N_R,\allowbreak N_{\mathrm{ind}},\allowbreak N_{\mathrm{max}})\geq a_{\mathrm{mol}}$.
	For $a>a_{\mathrm{flow}}$ and every fixed $\sigma\in\R$, one has
	\begin{align}
		&\nrm{\nabla\Phi_{\overline u_q}(t,\cdot;t_0)-\Id}_{N+\alpha}+\nrm{\bigl(\nabla\Phi_{\overline u_q}(t,\cdot;t_0)\bigr)^{-1}-\Id}_{N+\alpha}+\nrm{\nabla X_{\overline u_q}(t,\cdot;t_0)-\Id}_{N+\alpha}+\nrm{\operatorname{cof}\nabla\Phi_{\overline u_q}(t,\cdot;t_0)-\Id}_{N+\alpha} \notag\\
		&\quad\lesssim h\Lambda_q\La{N,N_{\mathrm{ind}}-1,\lambda_q,\ell_q^{-1}}\MM{N,N_{\mathrm{ind}}-1,\lambda_q,\ell_q^{-1}} \notag\\
		&\quad\lesssim\lambda_{q+1}^{-4\alpha}\La{N,N_{\mathrm{ind}}-1,\lambda_q,\ell_q^{-1}} \MM{N,N_{\mathrm{ind}}-1,\lambda_q,\ell_q^{-1}},
		\hspace{90pt} 0\leq N\leq N_{\mathrm{max}}-2,
		\label{eq:mollified-flow-close-identity}\\
		&\nrm{J_{\Phi_{\overline u_q}}^\sigma(t,\cdot)-1}_{N+\alpha}+\nrm{J_{X_{\overline u_q}}^\sigma(t,\cdot)-1}_{N+\alpha} \notag\\
		&\quad\lesssim h\delta_q\Lambda_q\La{N,N_{\mathrm{ind}}-1,\lambda_q,\ell_q^{-1}}\MM{N,N_{\mathrm{ind}}-1,\lambda_q,\ell_q^{-1}} \notag\\
		&\quad\lesssim\lambda_{q+1}^{-4\alpha}\delta_q\La{N,N_{\mathrm{ind}}-1,\lambda_q,\ell_q^{-1}} \MM{N,N_{\mathrm{ind}}-1,\lambda_q,\ell_q^{-1}},
		\hspace{80pt} 0\leq N\leq N_{\mathrm{max}}-2.
		\label{eq:mollified-flow-jacobian-close-one}
	\end{align}
	For positive material orders, one has
	\begin{align}
		&\nrm{\DTL^M\nabla\Phi_{\overline u_q}(t,\cdot;t_0)}_{N+\alpha}+\nrm{\DTL^M\bigl(\nabla\Phi_{\overline u_q}(t,\cdot;t_0)\bigr)^{-1}}_{N+\alpha}+\nrm{\partial_t^M\nabla X_{\overline u_q}(t,\cdot;t_0)}_{N+\alpha}+\nrm{\DTL^M\bigl(\operatorname{cof}\nabla\Phi_{\overline u_q}(t,\cdot;t_0)\bigr)}_{N+\alpha} \notag\\
		&\quad\lesssim \Lambda_q^M\lambda_q^{-1}\La{N+1,N_{\mathrm{ind}}-M+1,\lambda_q,\ell_q^{-1}}\notag\\
		&\qquad\cdot \MM{N+1,N_{\mathrm{ind}}-M+1,\lambda_q,\ell_q^{-1}},
		\hspace{99pt}1\leq M\leq N_{\mathrm{ind}}+1,\ \ M+N\leq N_{\mathrm{max}}-1,
		\label{eq:mollified-flow-material-estimates}\\
		&\nrm{\DTL^MJ_{\Phi_{\overline u_q}}^\sigma(t,\cdot)}_{N+\alpha}+\nrm{\partial_t^MJ_{X_{\overline u_q}}^\sigma(t,\cdot)}_{N+\alpha} \notag\\
		&\quad\lesssim \delta_q\Lambda_q^M\La{N,N_{\mathrm{ind}}-M,\lambda_q,\ell_q^{-1}}  \MM{N,N_{\mathrm{ind}}-M,\lambda_q,\ell_q^{-1}},
		\ 1\leq M\leq N_{\mathrm{ind}},\ \ M+N\leq N_{\mathrm{max}}-1.
		\label{eq:mollified-flow-jacobian-material-estimates}
	\end{align}
	For $M=N_{\mathrm{ind}}+1$ and $M+N\leq N_{\mathrm{max}}-1$,
	the sum of the norms in
	\eqref{eq:mollified-flow-jacobian-material-estimates} is bounded by
	the right-hand side of \eqref{eq:mollified-flow-material-estimates}.
	All these estimates also hold in $C^N$, with the factors
	$\La{\cdot}$ omitted and the same ranges of $M$, for
	$N\leq N_{\mathrm{max}}-1$ at $M=0$ and
	$M+N\leq N_{\mathrm{max}}$ at $M\geq1$.
	The H{\"o}lder estimates also hold with $\alpha$ replaced by any fixed $\gamma\in(0,\alpha)$.
	The implicit constants may depend on $\sigma,\rho^*,\mathfrak M,N_{\mathrm{max}},\alpha$
	and, for the estimates with exponent $\gamma$, on $\gamma$.
	They are independent of $q,a,h,t_0$.
\end{co}

\begin{proof}
	Apply Lemma~\ref{lem:flow-estimates-upsilon} with $\upsilon=\overline u_q$, $\tau=h$, and
	$$
	\begin{aligned}
		K^*=N_{\mathrm{ind}},\ \ K^\sharp=N_{\mathrm{max}},\ \ N_t=N_{\mathrm{ind}},\ \
		\lambda_\upsilon=\lambda_q,\ \ \tilde\lambda_\upsilon=\ell_q^{-1},\ \  \mu_\upsilon=\tilde\mu_\upsilon=\Lambda_q,\ \
		\mC_\upsilon=C\delta_q^{\frac12},\ \ \mC_{\Div,\upsilon}=C\delta_q^{\frac32}.
	\end{aligned}
	$$
	Choose a fixed $C\geq1$ at least as large as the constants in
	\eqref{est:uniform-grad-om-over-orho}--\eqref{est:uniform-div-om-over-orho} and \eqref{est:holder-grad-om-over-orho}--\eqref{est:holder-div-om-over-orho}.
	These estimates verify the coefficient hypotheses.
	Since $\delta_q\leq1$ and $C\geq1$,
	$\mu_{\upsilon,*}=\tilde\mu_{\upsilon,*}=C\Lambda_q$.
	By \eqref{def:tau-q}, choose $a_{\mathrm{flow}}\geq a_{\mathrm{mol}}$ so that, for $a>a_{\mathrm{flow}}$,
	$$
	h\mC_\upsilon\lambda_\upsilon=hC\Lambda_q
	\leq C\lambda_{q+1}^{-4\alpha}\leq c_0,\qquad q\in\mathbb N_0.
	$$
	
	For the estimates with a factor $h$, we combine \eqref{est:uniform-grad-om-over-orho} with
	\eqref{PhiN} to obtain
	$$
	\begin{aligned}
		\nrm{\nabla\Phi_{\overline u_q}-\Id}_N
		&\lesssim h\nrm{\nabla\overline u_q}_N
		\leq Ch\Lambda_q\MM{N,N_{\mathrm{ind}}-1,\lambda_q,\ell_q^{-1}},
		\qquad 0\leq N\leq N_{\mathrm{max}}-1.
	\end{aligned}
	$$
	The identities
	$$
	\begin{aligned}
		\bigl(\nabla\Phi_{\overline u_q}\bigr)^{-1}-\Id=-\bigl(\nabla\Phi_{\overline u_q}\bigr)^{-1}\bigl(\nabla\Phi_{\overline u_q}-\Id\bigr),\quad
		\operatorname{cof}\nabla\Phi_{\overline u_q}=J_{\Phi_{\overline u_q}}\bigl(\nabla\Phi_{\overline u_q}\bigr)^{-T},
	\end{aligned}
	$$
	together with
	$X_{\overline u_q}(t,\cdot;t_0)=\Phi_{\overline u_q}(t_0,\cdot;t)$,
	give the other terms in \eqref{eq:mollified-flow-close-identity}.
	Equation~\eqref{eq:flow-jacobian-CN} and composition with
	$s\mapsto s^\sigma$ on $[1/2,3/2]$ give
	$$\nrm{J_{\Phi_{\overline u_q}}^\sigma}_N\lesssim\MM{N,N_{\mathrm{ind}}-1,\lambda_q,\ell_q^{-1}}.$$
	By \eqref{Dt Jm}, the equation for $J_{\Phi_{\overline u_q}}^\sigma-1$ has source
	$-\sigma(\Div\overline u_q)J_{\Phi_{\overline u_q}}^\sigma$ and zero initial data.
	The product $(\Div\overline u_q)J_{\Phi_{\overline u_q}}^\sigma$ satisfies
	$$
	\begin{aligned}
		\nrm{(\Div\overline u_q)J_{\Phi_{\overline u_q}}^\sigma}_N
		\lesssim\sum_{n=0}^N\nrm{\Div\overline u_q}_n\nrm{J_{\Phi_{\overline u_q}}^\sigma}_{N-n}
		\lesssim C\delta_q\Lambda_q\MM{N,N_{\mathrm{ind}}-1,\lambda_q,\ell_q^{-1}}.
	\end{aligned}
	$$
	Here we used \eqref{est:uniform-div-om-over-orho} and the displayed estimate for $J_{\Phi_{\overline u_q}}^\sigma$.
	Applying the $C^N$ estimate in Lemma~\ref{estimate of the transport equation} and using $hC\Lambda_q\leq1$, we obtain
	$$
	\begin{aligned}
		\nrm{J_{\Phi_{\overline u_q}}^\sigma-1}_N
		&\lesssim Ch\delta_q\Lambda_q\MM{N,N_{\mathrm{ind}}-1,\lambda_q,\ell_q^{-1}},
		\qquad 0\leq N\leq N_{\mathrm{max}}-1.
	\end{aligned}
	$$
	Formula~\eqref{eq:pullback-Jacobian-material} at $M=0$ and the composition
	estimates give the corresponding estimate for $J_{X_{\overline u_q}}^\sigma-1$.
	Interpolation between spatial orders $N$ and $N+1$ proves the H{\"o}lder estimates with the same factors $h$ and $\delta_q$.
	
	For positive material orders, \eqref{eq:flow-gradient-material} and
	\eqref{eq:flow-gradient-material-Calpha} give the derivative estimates for $\nabla\Phi_{\overline u_q}$, $(\nabla\Phi_{\overline u_q})^{-1}$, and $\nabla X_{\overline u_q}$ in \eqref{eq:mollified-flow-material-estimates}.
	The estimates \eqref{eq:flow-Jacobian-power-material} and
	\eqref{eq:flow-Jacobian-power-material-Calpha} yield
	\eqref{eq:mollified-flow-jacobian-material-estimates}.
	The cofactor identity and \eqref{eq:holder-product-finite-Cnalpha} prove the cofactor
	estimate in \eqref{eq:mollified-flow-material-estimates} for
	$1\leq M\leq N_{\mathrm{ind}}+1$. At $M=N_{\mathrm{ind}}+1$, the material derivatives of
	$J_{\Phi_{\overline u_q}}$ satisfy the additional estimates in
	Lemma~\ref{lem:flow-estimates-upsilon}, with $\mC_\upsilon$ in place of
	$\mC_{\Div,\upsilon}$.
\end{proof}

\subsection{Mollification along the flow for the Newton iteration} 
\label{subsec:initial-Newton-flow-mollification} 

Let $\upsilon$ be a smooth spatially periodic velocity field and $f$ a smooth field on $\mcI\times\T^3$.
As in \cite{GK22,GR24}, we mollify $f$ along the flow of $\upsilon$. With the
temporal kernel fixed above, define
$$
(\eta_{\ell_t}^{\mathrm t}\ast_{X_\upsilon}f)(t,x):=\int_{\mathbb R}f(t+s,X_\upsilon(t+s,x;t))\eta_{\ell_t}^{\mathrm t}(s)\,\rd s,
$$
whenever $[t-\ell_t,t+\ell_t]\subset\mcI$. For the initial
stress, we define
\begin{equation} 
	R_{q,0}(t,x) 
	:=\eta_{\ell_{t,q}}^{\mathrm t}\ast_{X_{\overline u_q}}R_{q,0}^{\mathrm{sp}}=\int_{\mathbb R}R_{q,0}^{\mathrm{sp}}(t+s,X_{\overline u_q}(t+s,x;t))\eta_{\ell_{t,q}}^{\mathrm t}(s)\,\rd s, 
	\qquad t\in\mathbb R. 
	\label{def:initial-Newton-flow-mollified-stress} 
\end{equation} 

\begin{lem} 
	\label{lem:initial-Newton-flow-mollification} 
	Assume the hypotheses of Proposition~\ref{prop:uniform-mollification-estimates}.
	There exists $a_{\mathrm{init}}=a_{\mathrm{init}}(\beta,\allowbreak b,\allowbreak \alpha,\allowbreak \rho^*,\allowbreak \mathfrak M,\allowbreak N_R,\allowbreak N_{\mathrm{ind}},\allowbreak N_{\mathrm{max}})\geq a_{\mathrm{flow}}$.
	For $a>a_{\mathrm{init}}$, the stress $R_{q,0}$ defined in
	\eqref{def:initial-Newton-flow-mollified-stress} satisfies 
	\begin{equation} 
		\supp_tR_{q,0} 
		\subset 
		\mN_{\ell_{t,q}}(\supp_tR_q). 
		\label{eq:initial-Newton-flow-mollified-support} 
	\end{equation} 
	Moreover, for every $M,N\in\mathbb N_0$ satisfying 
	$ 
	\min\{M,N_R\}+N\leq N_{\mathrm{max}}-1, 
	$
	one has 
	\begin{equation} 
		\begin{aligned} 
			\Nrm{\DTL^MR_{q,0}}_{N+\alpha} \lesssim
			\lambda_q^{-8\alpha}\delta_{q+1}\MM{M,N_R,\Lambda_q,\ell_{t,q}^{-1}}\La{N,(N_R-M)_+,\lambda_q,\ell_q^{-1}}
			\MM{N,(N_R-M)_+,\lambda_q,\ell_q^{-1}}. 
		\end{aligned} 
		\label{eq:initial-Newton-flow-mollified-stress-estimate} 
	\end{equation} 
	
	The mollification error satisfies the following estimate for
	$L\in\mathbb N$, $M+L\leq N_R$, and
	$M+N+L\leq N_{\mathrm{max}}-1$.
	\begin{equation} 
		\Nrm{\DTL^M(R_{q,0}^{\mathrm{sp}}-R_{q,0})}_{N+\alpha} 
		\lesssim
		\lambda_q^{-8\alpha}(\ell_{t,q}\Lambda_q)^L\delta_{q+1} 
		\Lambda_q^M
		\La{N,N_R-M-L,\lambda_q,\ell_q^{-1}} 
		\MM{N,N_R-M-L,\lambda_q,\ell_q^{-1}}. 
		\label{eq:initial-Newton-flow-mollification-remainder-estimate} 
	\end{equation} 
	At the endpoint $M+N=N_R$, one has 
	\begin{equation} 
		\Nrm{\DTL^M(R_{q,0}^{\mathrm{sp}}-R_{q,0})}_{N+\alpha} 
		\lesssim
		\lambda_q^{-8\alpha}\delta_{q+1} 
		\Lambda_q^M 
		\La{N,N_R-M,\lambda_q,\ell_q^{-1}} 
		\MM{N,N_R-M,\lambda_q,\ell_q^{-1}}. 
		\label{eq:initial-Newton-flow-mollification-endpoint-estimate} 
	\end{equation} 
	The estimates also hold with $\alpha$ replaced by any fixed $\gamma\in(0,\alpha)$,
	with the stress amplitude unchanged, and in $C^N$ with the factors
	$\La{\cdot}$ omitted.
	The implicit constants may depend on $M,\rho^*,\mathfrak M,N_{\mathrm{max}},\alpha$
	and, for the estimates with exponent $\gamma$, on $\gamma$.
	They are independent of $q$ and $a$.
\end{lem} 

\begin{proof} 
	The support inclusion \eqref{eq:initial-Newton-flow-mollified-support}
	follows from $\supp_tR_{q,0}^{\mathrm{sp}}\subset\supp_tR_q$ and
	$\supp\eta_{\ell_{t,q}}^{\mathrm t}\subset(-\ell_{t,q},\ell_{t,q})$.
	By \eqref{eq:material-mollification-hierarchy}, choose $a_{\mathrm{init}}\geq a_{\mathrm{flow}}$ so that
	$\ell_{t,q}<\tau_q$ for $a>a_{\mathrm{init}}$.
	Thus Corollary~\ref{cor:mollified-velocity-flow-estimates} applies with $h=\ell_{t,q}$, and
	\eqref{eq:mollified-flow-close-identity} holds uniformly for $|s|\leq\ell_{t,q}$.
	
	By \eqref{eq of forward flow} and \eqref{for and back}, material differentiation
	commutes with pullback along $X_{\overline u_q}$. Hence it may be passed under
	the integral in \eqref{def:initial-Newton-flow-mollified-stress}.
	Integration by parts in $s$ gives, for
	$0\leq r\leq\min\{M,N_R\}$ and $r+N\leq N_{\mathrm{max}}-1$,
	\begin{equation} 
		\begin{aligned} 
			\DTL^MR_{q,0}(t,x) 
			&=\int_{\mathbb R}(\DTL^MR_{q,0}^{\mathrm{sp}})(t+s,X_{\overline u_q}(t+s,x;t))\eta_{\ell_{t,q}}^{\mathrm t}(s)\,\rd s\\ 
			&=(-1)^{M-r}\int_{\mathbb R}(\DTL^rR_{q,0}^{\mathrm{sp}})(t+s,X_{\overline u_q}(t+s,x;t))\partial_s^{M-r}\eta_{\ell_{t,q}}^{\mathrm t}(s)\,\rd s. 
		\end{aligned} 
		\label{eq:initial-Newton-flow-material-integration-by-parts} 
	\end{equation} 
	The kernel satisfies
	$$
	\|\partial_s^{M-r}\eta_{\ell_{t,q}}^{\mathrm t}\|_{L^1(\mathbb R)} 
	=\ell_{t,q}^{-(M-r)}\|\partial_s^{M-r}\eta^{\mathrm t}\|_{L^1(\mathbb R)}
	\lesssim 
	\ell_{t,q}^{-(M-r)}. 
	$$
	For $0\leq r\leq N_R$ and $r+N\leq N_{\mathrm{max}}-1$,
	we have $N_R-r\leq N_{\mathrm{ind}}-1$. We use the estimate
	\eqref{est:holder-Rq0} and the flow estimates from
	Corollary~\ref{cor:mollified-velocity-flow-estimates}. The second inequality in
	Lemma~\ref{lem:BDSV-composition} gives the $C^N$ composition estimate for $N\geq1$,
	and Lemma~\ref{lem:holder-composition-general} proves the H{\"o}lder
	seminorm estimate, including $N=0$. For $N\geq1$, at most $N-1$ spatial derivatives of
	$\nabla X_{\overline u_q}$ are needed, so the required spatial derivative estimate remains valid
	throughout the stated range. Thus
	\begin{equation} 
		\begin{aligned} 
			\Nrm{ 
				(\DTL^rR_{q,0}^{\mathrm{sp}}) 
				(t+s,X_{\overline u_q}(t+s,\cdot;t)) 
			}_{N+\alpha}\lesssim 
			\lambda_q^{-8\alpha}\delta_{q+1} 
			\Lambda_q^r 
			\La{N,N_R-r,\lambda_q,\ell_q^{-1}} 
			\MM{N,N_R-r,\lambda_q,\ell_q^{-1}}, 
		\end{aligned} 
		\label{eq:initial-Newton-flow-pullback-estimate} 
	\end{equation} 
	uniformly for $|s|\leq\ell_{t,q}$. With $r=\min\{M,N_R\}$,
	\eqref{eq:initial-Newton-flow-material-integration-by-parts} yields
	\begin{equation} 
		\begin{aligned} 
			\Nrm{\DTL^MR_{q,0}}_{N+\alpha} 
			&\lesssim 
			\lambda_q^{-8\alpha}\delta_{q+1} 
			\ell_{t,q}^{-(M-r)}\Lambda_q^r 
			\La{N,N_R-r,\lambda_q,\ell_q^{-1}} 
			\MM{N,N_R-r,\lambda_q,\ell_q^{-1}}. 
		\end{aligned} 
		\label{eq:initial-Newton-flow-mollification-gain} 
	\end{equation} 
	Here $N_R-r=(N_R-M)_+$ and
	$ 
	\ell_{t,q}^{-(M-r)}\Lambda_q^r 
	= 
	\MM{M,N_R,\Lambda_q,\ell_{t,q}^{-1}}.
	$ 
	Hence \eqref{eq:initial-Newton-flow-mollified-stress-estimate} holds.
	
	For $L\geq1$, $M+L\leq N_R$, and $M+N+L\leq N_{\mathrm{max}}-1$,
	\eqref{eq of forward flow} gives, for $0\leq k\leq L$,
	$$
	\partial_s^k\left[(\DTL^MR_{q,0}^{\mathrm{sp}})(t+s,X_{\overline u_q}(t+s,x;t))\right]=(\DTL^{M+k}R_{q,0}^{\mathrm{sp}})(t+s,X_{\overline u_q}(t+s,x;t)).
	$$
	Taylor's formula at $s=0$ gives
	$$
	\begin{aligned}
		(\DTL^MR_{q,0}^{\mathrm{sp}})(t+s,X_{\overline u_q}(t+s,x;t))&=\sum_{k=0}^{L-1}\frac{s^k}{k!}(\DTL^{M+k}R_{q,0}^{\mathrm{sp}})(t,x)\\
		&\quad+\frac{s^L}{(L-1)!}\int_0^1(1-\theta)^{L-1}(\DTL^{M+L}R_{q,0}^{\mathrm{sp}})(t+\theta s,X_{\overline u_q}(t+\theta s,x;t))\,\rd\theta.
	\end{aligned}
	$$
	Since $L\leq N_R-M<N_{\mathrm{ind}}$, the moment conditions give
	$\int_{\mathbb R}s^k\eta_{\ell_{t,q}}^{\mathrm t}(s)\,\rd s=0$ for $1\leq k<L$.
	By the normalization of $\eta^{\mathrm t}$ and
	\eqref{def:material-mollification-scale}, we have
	$$
	(\DTL^MR_{q,0}^{\mathrm{sp}})(t,x)
	-\int_{\mathbb R}
	(\DTL^MR_{q,0}^{\mathrm{sp}})(t,x)
	\eta_{\ell_{t,q}}^{\mathrm t}(s)\,\rd s=0.
	$$
	Integrating the Taylor remainder against the kernel and using
	\eqref{eq:initial-Newton-flow-material-integration-by-parts} with $r=M$ gives
	\begin{equation}
		\DTL^M(R_{q,0}^{\mathrm{sp}}-R_{q,0})(t,x)
		=-\frac{1}{(L-1)!}\int_{\mathbb R}s^L\eta_{\ell_{t,q}}^{\mathrm t}(s)\int_0^1(1-\theta)^{L-1}(\DTL^{M+L}R_{q,0}^{\mathrm{sp}})(t+\theta s,X_{\overline u_q}(t+\theta s,x;t))\,\rd\theta\,\rd s.
		\label{eq:initial-Newton-flow-mollification-remainder-identity}
	\end{equation}
	Since
	$\int_{\mathbb R}|s|^L|\eta_{\ell_{t,q}}^{\mathrm t}(s)|\,\rd s=\ell_{t,q}^L\int_{\mathbb R}|s|^L|\eta^{\mathrm t}(s)|\,\rd s\lesssim\ell_{t,q}^L$,
	we have
	$$
	\Nrm{\DTL^M(R_{q,0}^{\mathrm{sp}}-R_{q,0})}_{N+\alpha}\lesssim\ell_{t,q}^L\sup_{|s|\leq\ell_{t,q}}\Nrm{(\DTL^{M+L}R_{q,0}^{\mathrm{sp}})(t+s,X_{\overline u_q}(t+s,\cdot;t))}_{N+\alpha}.
	$$
	Applying \eqref{eq:initial-Newton-flow-pullback-estimate} with $r=M+L$ and using
	$
	\ell_{t,q}^L\Lambda_q^{M+L}=\Lambda_q^M(\ell_{t,q}\Lambda_q)^L,
	$
	we obtain \eqref{eq:initial-Newton-flow-mollification-remainder-estimate}.
	
	For $M+N=N_R$, \eqref{est:holder-Rq0} and
	\eqref{eq:initial-Newton-flow-mollified-stress-estimate}
	give \eqref{eq:initial-Newton-flow-mollification-endpoint-estimate}.
	The $C^N$ estimates follow from \eqref{est:uniform-Rq0} by the same
	argument.
\end{proof} 

\section{A local right inverse of divergence}
\label{sec:local-inverse-divergence}

Before describing the perturbation in the Newton step, we introduce a local
right inverse of the divergence operator, with no symmetry constraint on the
resulting matrix field. Its role is to represent a compactly supported error field as the
divergence of a compactly supported matrix field. More precisely, at each fixed
time, we seek a matrix field $R[U]$ satisfying
\begin{equation}
	\Div R[U]=U,
	\qquad
	R[U]:\R^d\to\R^{d\times d}.
	\label{eq:nonsymmetric-divergence-equation}
\end{equation}
Here $U:\R^d\to\R^d$ is compactly supported and satisfies the zero-mean
condition
\begin{equation}
	\int_{\R^d} U_i(x)\,\rd x=0,
	\qquad
	i=1,\ldots,d.
	\label{eq:zero-mean-condition-for-U}
\end{equation}
This is the compatibility condition for solving
\eqref{eq:nonsymmetric-divergence-equation} with a compactly supported matrix
field.

The Newton step permits a nonsymmetric right inverse. If $A=(A_{ij})$ is
antisymmetric, then
$
\partial_i\partial_j A_{ij}=0,
$
and hence $\Div A$ is divergence-free. For a general matrix, write
$
R_{ij}
=
\frac{R_{ij}+R_{ji}}{2}
+
\frac{R_{ij}-R_{ji}}{2},
$
where the antisymmetric part produces a divergence-free source that can be treated by the
incompressible correction. Thus it suffices to solve the nonsymmetric divergence
equation $\partial_jR_{ij}[U]=U_i$ for the compressible part. For cutoff residuals, only
the symmetric part of $R[U]$ enters the gluing stress.

We use the translation-averaging argument of Isett--Oh~\cite{Isett16} to
construct a local right inverse. For the nonsymmetric divergence equation,
a first-order Taylor expansion of each component in Fourier space gives a
solution operator that preserves support in the given ball, with a possible singularity
at the center of the ball. Conjugating this operator by translations and taking
a smooth average yields a compactly supported Bogovskii-type operator.
Isett--Oh use a second-order expansion for the
symmetric divergence equation, which also requires the angular-momentum
compatibility condition.

We also need higher-order material derivative estimates. To obtain
them, we define a time-independent Bogovskii operator in Lagrangian coordinates,
which commutes with ordinary time derivatives, and then return to Eulerian coordinates
through the Piola transform.

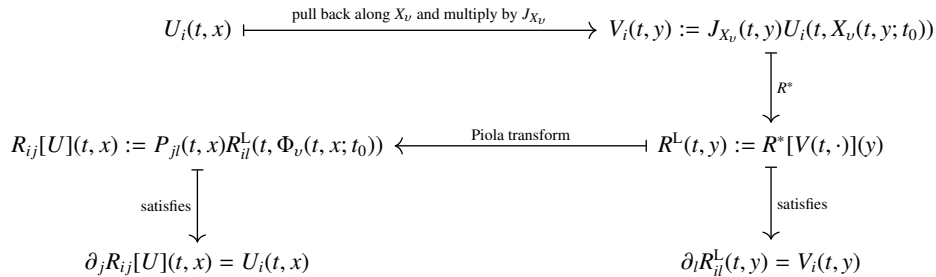
\begin{figure}[H]
	\centering
	\small
	\begin{tikzcd}[ampersand replacement=\&, row sep=3em, column sep=8.5em]
		U_i(t,x)
		\arrow[r, maps to, "{\text{pull back along }X_\upsilon\text{ and multiply by }\JX}"]
		\&
		V_i(t,y):=\JX(t,y)U_i(t,X_\upsilon(t,y;t_0))
		\arrow[d, maps to, "{R^*}"]
		\\
		R_{ij}[U](t,x):=P_{jl}(t,x)R^{\mathrm L}_{il}(t,\Phi_\upsilon(t,x;t_0))
		\arrow[d, maps to, "\text{satisfies}"']
		\&
		R^{\mathrm L}(t,y):=R^*[V(t,\cdot)](y)
		\arrow[l, maps to, "{\text{Piola transform}}"']
		\arrow[d, maps to, "\text{satisfies}"]
		\\
		\partial_jR_{ij}[U](t,x)=U_i(t,x)
		\&
		\partial_lR^{\mathrm L}_{il}(t,y)=V_i(t,y)
	\end{tikzcd}
	\caption{Construction of the local inverse-divergence operator}
	\label{fig:flow-coordinate-R}
\end{figure}

Let $y$ denote the Lagrangian label and $x$ the corresponding spatial variable,
so that $x=X_\upsilon(t,y;t_0)$ and $y=\Phi_\upsilon(t,x;t_0)$. More precisely,
we first pull the source $U$ back to the Lagrangian label $y$ and multiply it by
the Jacobian determinant of the forward flow map to obtain
$
V_i(t,y):=\JX(t,y)U_i(t,X_\upsilon(t,y;t_0)).
$
We then apply the local inverse-divergence operator $R^*$ in the $y$-variable and
set
$
R^{\mathrm L}(t,y):=R^*[V(t,\cdot)](y).
$
Thus $R^{\mathrm L}$ satisfies
$
\partial_lR^{\mathrm L}_{il}(t,y)=V_i(t,y),
$
where the derivative is taken with respect to $y$. Finally, we push this matrix
field forward to the spatial variable $x$ by the Piola transform to obtain
$
R_{ij}[U](t,x):=P_{jl}(t,x)R^{\mathrm L}_{il}(t,\Phi_\upsilon(t,x;t_0)).
$
The matrix $P$ is the transposed cofactor matrix associated with the backward flow map. By the Piola
identity, the resulting matrix field satisfies
$
\partial_jR_{ij}[U](t,x)=U_i(t,x).
$
The construction is summarized in Figure~\ref{fig:flow-coordinate-R}.

\subsection{Adapted cylinders}
\label{subsec:adapted-cylinders}
We use the adapted Eulerian and Lagrangian cylinders of
Isett--Oh~\cite{Isett16} to track supports.

\begin{df}
	Let $X_\VV(s,x;t)$ and $\Phi_\VV(s,x;t)$ be the forward and backward flow
	maps associated with the vector field $\VV$. Given $\tau,\varepsilon>0$ and
	a point $(t_0,x_0)\in\R\times\R^d$, we define the $\VV$-adapted Lagrangian
	cylinder by
	$$
	\hat\Gamma_\VV(\tau,\varepsilon;t_0,x_0)
	:=
	\left\{(s,X_\VV(s,y;t_0)): |s-t_0|\leq\tau,\ |y-x_0|<\varepsilon\right\}.
	$$
	Equivalently,
	$$
	\hat\Gamma_\VV(\tau,\varepsilon;t_0,x_0)
	=
	\left\{(s,x)\in\R\times\R^d: |s-t_0|\leq\tau,\ \Phi_\VV(s,x;t_0)\in B(x_0,\varepsilon)\right\}.
	$$
	Thus $\hat\Gamma_\VV(\tau,\varepsilon;t_0,x_0)$ is obtained by
	transporting the initial ball $B(x_0,\varepsilon)$ along the flow of $\VV$.
	We also define the $\VV$-adapted Eulerian cylinder by
	$$
	\hat C_\VV(\tau,\varepsilon;t_0,x_0)
	:=
	\left\{(s,x)\in\R\times\R^d: |s-t_0|\leq\tau,\ |x-X_\VV(s,x_0;t_0)|<\varepsilon\right\}.
	$$
	Thus $\hat C_\VV(\tau,\varepsilon;t_0,x_0)$ is the Eulerian tube of radius
	$\varepsilon$ around the center trajectory $X_\VV(s,x_0;t_0)$.
	
	If $S\subseteq\R\times\R^d$, we define
	$$
	\begin{aligned}
		\hat\Gamma_\VV(\tau,\varepsilon;S)
		:=
		\bigcup_{(t_0,x_0)\in S}\hat\Gamma_\VV(\tau,\varepsilon;t_0,x_0),\qquad
		\hat C_\VV(\tau,\varepsilon;S)
		:=
		\bigcup_{(t_0,x_0)\in S}\hat C_\VV(\tau,\varepsilon;t_0,x_0).
	\end{aligned}
	$$
\end{df}

In the cylinder inclusions and support estimates below, the $C^0$ norm of a velocity gradient is taken in the Euclidean operator norm.
The following lemma compares these cylinders.

\begin{lem}
	\label{lem:comparison-Gamma-C}
	Fix $\tau,\varepsilon>0$ and $(t_0,x_0)\in\R\times\R^d$,
	and let $\VV$ be a smooth vector field. Set
	$$
	\Omega_{\VV}(t_0,\tau)
	:=
	\int_{t_0-\tau}^{t_0+\tau}\|\nabla\VV(s,\cdot)\|_{0}\,\rd s .
	$$
	Assume that $\Omega_{\VV}(t_0,\tau)<\infty$.
	Then
	\begin{align}
		\label{eq:Gamma-contained-in-C}
		\hat\Gamma_{\VV}(\tau,\varepsilon;t_0,x_0)
		&\subseteq
		\hat C_{\VV}
		\left(
		\tau,
		e^{\Omega_{\VV}(t_0,\tau)}\varepsilon;
		t_0,x_0
		\right),\\
		\label{eq:C-contained-in-Gamma}
		\hat C_{\VV}(\tau,\varepsilon;t_0,x_0)
		&\subseteq
		\hat\Gamma_{\VV}
		\left(
		\tau,
		e^{\Omega_{\VV}(t_0,\tau)}\varepsilon;
		t_0,x_0
		\right).
	\end{align}
\end{lem}

\begin{proof}
	We prove \eqref{eq:Gamma-contained-in-C}. Let
	$(t,x)\in\hat\Gamma_{\VV}(\tau,\varepsilon;t_0,x_0)$. Then there exists
	$y\in B(x_0,\varepsilon)$ such that $x=X_{\VV}(t,y;t_0)$ and
	$|t-t_0|\leq\tau$. By the Lipschitz estimate for the flow map, we have
	$$
	\left|x-X_{\VV}(t,x_0;t_0)\right|
	=
	\left|X_{\VV}(t,y;t_0)-X_{\VV}(t,x_0;t_0)\right|
	\leq
	e^{\Omega_{\VV}(t_0,\tau)}|y-x_0|
	<
	e^{\Omega_{\VV}(t_0,\tau)}\varepsilon .
	$$
	This proves \eqref{eq:Gamma-contained-in-C}.
	
	For the inclusion \eqref{eq:C-contained-in-Gamma}, let
	$(t,x)\in\hat C_{\VV}(\tau,\varepsilon;t_0,x_0)$. Then $|t-t_0|\leq\tau$
	and
	$
	\left|x-X_{\VV}(t,x_0;t_0)\right|<\varepsilon.
	$
	Set $y:=\Phi_{\VV}(t,x;t_0)$. Then $X_{\VV}(t,y;t_0)=x$ and
	$
	\Phi_{\VV}(t,X_{\VV}(t,x_0;t_0);t_0)=x_0.
	$
	Applying the Lipschitz estimate to the backward flow gives
	$$
	|y-x_0|
	=
	\left|
	\Phi_{\VV}(t,x;t_0)
	-
	\Phi_{\VV}(t,X_{\VV}(t,x_0;t_0);t_0)
	\right|
	\leq
	e^{\Omega_{\VV}(t_0,\tau)}
	\left|
	x-X_{\VV}(t,x_0;t_0)
	\right|
	<
	e^{\Omega_{\VV}(t_0,\tau)}\varepsilon .
	$$
	Hence $(t,x)$ belongs to
	$
	\hat\Gamma_{\VV}
	\left(
	\tau,
	e^{\Omega_{\VV}(t_0,\tau)}\varepsilon;
	t_0,x_0
	\right),
	$
	which proves \eqref{eq:C-contained-in-Gamma}.
\end{proof}

Throughout the rest of this section, fix a sufficiently large universal
constant $C_{\mathrm{loc}}\geq4$. Whenever the adapted cylinder
$\hat C_{\upsilon}(\tau,C_{\mathrm{loc}}\varepsilon;t_0,x_0)$ is fixed, we use
the following spatial norms, taken uniformly in time. For $N\in\mathbb N_0$ and
$0<\alpha<1$, set
$$
\begin{aligned}
	\nrm{f}_{N}
	:=
	&
	\sum_{|\beta|\leq N}
	\sup_{(t,x)\in\hat C_{\upsilon}(\tau,C_{\mathrm{loc}}\varepsilon;t_0,x_0)}
	|\partial^\beta f(t,x)|,\\
	\nrm{f}_{N+\alpha}
	:=
	&
	\sum_{|\beta|\leq N}
	\sup_{(t,x)\in\hat C_{\upsilon}(\tau,C_{\mathrm{loc}}\varepsilon;t_0,x_0)}
	|\partial^\beta f(t,x)|
	+
	\sum_{|\beta|=N}
	\sup_{\substack{(t,x),(t,z)\in\hat C_{\upsilon}(\tau,C_{\mathrm{loc}}\varepsilon;t_0,x_0)\\ x\neq z}}
	\frac{|\partial^\beta f(t,x)-\partial^\beta f(t,z)|}{|x-z|^\alpha}.
\end{aligned}
$$
Here and below, $\partial^\beta$ denotes a spatial derivative. For compactly supported functions
of the Lagrangian variable, the same notation denotes the usual $C^N(\R^d)$ or
$C^{N+\alpha}(\R^d)$ norm, uniformly in time. For quantities associated with the flow, norms are taken on $B(x_0,\varepsilon)$ in Lagrangian coordinates
and on $X_\upsilon(t,B(x_0,\varepsilon);t_0)$ in Eulerian coordinates.
If $e^{2\tau\nrm{\nabla\upsilon}_0}<2$, Gronwall's inequality up to the first exit from
$\hat C_\upsilon(\tau,2\varepsilon;t_0,x_0)$ rules out such an exit for the forward and backward trajectories used in Lemma~\ref{lem:comparison-Gamma-C}.
Thus the same cylinder inclusions hold with the local gradient norm.

We now define the local inverse-divergence operator on the adapted cylinder.
Let $X_\upsilon(t,\cdot;t_0)$ be the forward flow map
generated by $\upsilon$ from the reference time $t_0$, and let
$\Phi_\upsilon(t,\cdot;t_0)$ be the corresponding backward flow map. Let $x=X_\upsilon(t,y;t_0)$ and
$y=\Phi_\upsilon(t,x;t_0)$. We write
$$
\JX(t,y):=\det\nabla X_\upsilon(t,y;t_0),
\qquad
\JPhi(t,x):=\det\nabla\Phi_\upsilon(t,x;t_0).
$$
Then
$$
\JPhi(t,X_\upsilon(t,y;t_0))\JX(t,y)=1.
$$
We denote the transposed cofactor matrix associated with the backward flow map by
$$
P_{ij}(t,x):=\left[(\cof\nabla\Phi_\upsilon(t,x;t_0))^T\right]_{ij}
=\JPhi(t,x)\left[\bigl(\nabla\Phi_\upsilon\bigr)^{-1}(t,x;t_0)\right]_{ij}.
$$
It satisfies the Piola identity
$
\partial_jP_{jl}(t,x)=0.
$

For any smooth matrix field $R^{\mathrm L}=(R^{\mathrm L}_{ij})$, define its Piola transform by
$$
R(t,x):=R^{\mathrm L}(t,\Phi_\upsilon(t,x;t_0))P^T(t,x).
$$
Using the Piola identity, we have
$$
\begin{aligned}
	\partial_jR_{ij}
	&=
	P_{jl}(t,x)\partial_j
	\left(R^{\mathrm L}_{il}(t,\Phi_\upsilon(t,x;t_0))\right)=
	P_{jl}(t,x)(\partial_mR^{\mathrm L}_{il})(t,\Phi_\upsilon(t,x;t_0))
	\partial_j(\Phi_\upsilon)_m(t,x;t_0).
\end{aligned}
$$
Using the definition of $P_{jl}$, we obtain
$$
\begin{aligned}
	P_{jl}(t,x)\partial_j(\Phi_\upsilon)_m(t,x;t_0)
	&=
	\JPhi(t,x)\left[\bigl(\nabla\Phi_\upsilon\bigr)^{-1}(t,x;t_0)\right]_{jl}
	\partial_j(\Phi_\upsilon)_m(t,x;t_0)=
	\JPhi(t,x)\delta_{ml}.
\end{aligned}
$$
Consequently,
$$
\Div\bigl((R^{\mathrm L}\circ\Phi_\upsilon)P^T\bigr)
=
\JPhi(\Div R^{\mathrm L})\circ\Phi_\upsilon.
$$

We denote by $R^*$ a Bogovskii-type local inverse-divergence operator. More
precisely, if $V\in C_c^\infty(B(x_0,\bar\varepsilon);\R^d)$ satisfies
$
\int_{\R^d}V(y)\,\rd y=0,
$
then $R^*[V]$ is a smooth matrix field satisfying
$
\partial_l(R^*[V])_{il}=V_i
$
and
$
\supp R^*[V]\subset B(x_0,\bar\varepsilon).
$

Suppose that
$
\supp U\subset\hat\Gamma_\upsilon(\bar\tau,\bar\varepsilon;t_0,x_0)
$
and that the zero-mean compatibility condition
$
\int_{\R^d}U(t,x)\,\rd x=0
$
holds for every $t$. Define
$$
V(t,y):=\JX(t,y)U(t,X_\upsilon(t,y;t_0)).
$$
By the support assumption, $\supp V(t,\cdot)\subset B(x_0,\bar\varepsilon)$
for $|t-t_0|\leq\bar\tau$. Moreover, the change of variables
$x=X_\upsilon(t,y;t_0)$ gives
$$
\int_{B(x_0,\bar\varepsilon)}V(t,y)\,\rd y=0,
\qquad
\forall |t-t_0|\leq\bar\tau.
$$
Therefore, setting
$
R^{\mathrm L}(t,y):=R^*[V(t,\cdot)](y),
$
we have
$$
\partial_lR^{\mathrm L}_{il}(t,y)=V_i(t,y),
\qquad
\supp R^{\mathrm L}(t,\cdot)\subset B(x_0,\bar\varepsilon).
$$

Finally, define
$$
R_{ij}[U](t,x):=P_{jl}(t,x)R^{\mathrm L}_{il}(t,\Phi_\upsilon(t,x;t_0)).
$$
Using the Piola identity above, we obtain
$$
\begin{aligned}
	\partial_jR_{ij}[U](t,x)
	=\JPhi(t,x)(\partial_lR^{\mathrm L}_{il})(t,\Phi_\upsilon(t,x;t_0))
	=\JPhi(t,x)V_i(t,\Phi_\upsilon(t,x;t_0))
	=U_i(t,x).
\end{aligned}
$$
Thus
$$
\Div R[U]=U,
\qquad
\supp R[U]\subset\hat\Gamma_\upsilon(\bar\tau,\bar\varepsilon;t_0,x_0).
$$
If $\Omega_\upsilon(t_0,\bar\tau)<\infty$, Lemma~\ref{lem:comparison-Gamma-C} gives the support inclusion
$$
\supp R[U]\subset
\hat C_\upsilon\left(\bar\tau,e^{\Omega_\upsilon(t_0,\bar\tau)}\bar\varepsilon;t_0,x_0\right).
$$

For the periodic construction in the Newton step, the local inverse is applied to each
localized source on $\R^3$. Lemma~\ref{lem support of R under cutoffs}
below shows that the sum of the resulting matrix fields is locally finite
and, under the translation covariance condition on the source fields, periodic.

\subsection{Main theorem}
\label{subsec:main-result-nonsymmetric-divergence}

For $s\in\mathbb N$, write
$\Lk{N,K,\lambda,\tilde\lambda}{s}
:=\La{N,K,\lambda,\tilde\lambda}^{s}$.
We now state the main theorem.

\begin{thm}
	\label{thm:inv4divEq}
	Fix $0<\alpha<1$, $K^*,K^\sharp\in\mathbb N$, and
	$K_U,M_*,N_t\in\mathbb N_0$ with $K_U\leq K^*$ and $K^\sharp\geq K^*+2$.
	If $M_*\geq K^*+1$, assume also $N_t\leq K^*$.
	Let $\upsilon=(\upsilon_1,\ldots,\upsilon_d)$ be a smooth vector field on
	$\R\times\R^d$. Set
	$
	D_{t,\upsilon}:=\partial_t+\upsilon\cdot\nabla .
	$
	Fix $(t_0,x_0)\in\R\times\R^d$ and $\bar\tau>0$ such that the flow
	maps used below are defined for $|t-t_0|\leq\bar\tau$. Fix a localization radius
	satisfying
	\begin{equation}
		\label{eq:epsilon-material-smallness}
		0<\bar\varepsilon\leq1.
	\end{equation}
	The norms of $\upsilon$, $U$, and $R[U]$ below are taken over
	$\hat C_{\upsilon}(\bar\tau,C_{\mathrm{loc}}\bar\varepsilon;t_0,x_0)$.
	The norms of the flow maps are taken on the domains specified in Subsection~\ref{subsec:adapted-cylinders}.
	Let
	$
	1\leq\lambda_{\upsilon}\leq\tilde\lambda_{\upsilon},
	\
	1\leq\mu_{\upsilon}\leq\tilde\mu_{\upsilon},
	\
	\mC_{\upsilon}\geq0.
	$
	Assume that, for all $M,N\in\mathbb N_0$ in the ranges displayed below,
	\begin{align}
		\nrm{D_{t,\upsilon}^{M}\nabla\upsilon}_{N}
		&\leq
		\mC_{\upsilon}
		\MM{N+1,K^*-M,\lambda_{\upsilon},\tilde\lambda_{\upsilon}}
		\MM{M,N_t,\mu_{\upsilon},\tilde\mu_{\upsilon}},\label{eq:inv4divEq:est4nbv-material-CN-MM}
		&&0\leq M\leq K^*,\qquad M+N\leq K^\sharp-1,
		\\
		\nrm{D_{t,\upsilon}^{M}\nabla\upsilon}_{N+\alpha}
		&\leq
		\mC_{\upsilon}
		\La{N+1,K^*-M,\lambda_{\upsilon},\tilde\lambda_{\upsilon}}
		\notag\\
		\label{eq:inv4divEq:est4nbv-material-Calpha-MM}
		&\quad\cdot
		\MM{N+1,K^*-M,\lambda_{\upsilon},\tilde\lambda_{\upsilon}}\MM{M,N_t,\mu_{\upsilon},\tilde\mu_{\upsilon}},
		&&0\leq M\leq K^*,\qquad M+N\leq K^\sharp-2.
	\end{align}
	Assume also the local flow smallness condition
	\begin{equation}
		\label{eq:tau-local-smallness}
		\bar\tau\mC_{\upsilon}\lambda_{\upsilon}\leq c_0.
	\end{equation}
	Here $c_0>0$ depends only on the dimension $d$ and is chosen sufficiently small so that
	$\nrm{\nabla\Phi_\upsilon(t,\cdot;t_0)-\Id}_0+\nrm{\nabla X_\upsilon(t,\cdot;t_0)-\Id}_0\leq1/2$
	for $|t-t_0|\leq\bar\tau$.
	
	Fix $s\in\mathbb N$. Let
	$
	1\leq\lambda_U\leq\tilde\lambda_U,
	\
	1\leq\mu_U\leq\tilde\mu_U,
	\
	\mC_U\geq0,
	$
	and let $U$ be a smooth vector field satisfying the zero-mean condition
	\begin{equation}
		\label{eq:inv4divEq:hyp:1}
		\int_{\R^d}U_i(t,x)\,\rd x=0,
		\qquad
		i=1,\ldots,d,
	\end{equation}
	for every $t$. Assume moreover that
	\begin{equation}
		\label{eq:inv4divEq:supp4U}
		\supp U
		\subseteq
		\hat\Gamma_{\upsilon}(\bar\tau,\bar\varepsilon;t_0,x_0),
	\end{equation}
	and that, for all $M,N\in\mathbb N_0$ satisfying
	$0\leq M\leq M_*$ and $M+N\leq K^\sharp-1$,
	\begin{equation}
		\label{eq:inv4divEq:hyp:3-material-Calpha-MM}
		\begin{aligned}
			\nrm{D_{t,\upsilon}^{M}U}_{N+\alpha}
			&\leq
			\mC_U
			\Lk{N,(K_U-M)_+,\lambda_U,\tilde\lambda_U}{s}
			\MM{N,(K_U-M)_+,\lambda_U,\tilde\lambda_U}
			\MM{M,N_t,\mu_U,\tilde\mu_U}.
		\end{aligned}
	\end{equation}
	Set $K_0:=\min\{K_U,K^*-1\}$. Define
	\begin{align}
		\label{eq:definition-lambda-RUV}
		\lambda_{R,U,\upsilon}
		&:=\max\{\bar\varepsilon^{-1},\lambda_U,\lambda_{\upsilon}\},
		&
		\tilde\lambda_{R,U,\upsilon}
		&:=\max\{\bar\varepsilon^{-1},\tilde\lambda_U,\tilde\lambda_{\upsilon}\},
		\\
		\label{eq:definition-mu-RUV}
		\mu_{R,U,\upsilon}
		&:=\max\{\mu_U,\mu_{\upsilon},\mC_{\upsilon}\lambda_{\upsilon}\},
		&
		\tilde\mu_{R,U,\upsilon}
		&:=\max\{\tilde\mu_U,\tilde\mu_{\upsilon},\mC_{\upsilon}\tilde\lambda_{\upsilon}\}.
	\end{align}
	Then there exists a smooth matrix field
	$
	R[U]=\bigl(R_{ij}[U]\bigr)_{i,j=1}^{d}
	$
	depending linearly on $U$ and solving the nonsymmetric divergence equation
	\begin{equation}
		\label{eq:inv4divEq:nonsymmetric-divergence}
		\partial_jR_{ij}[U]=U_i,
		\qquad
		i=1,\ldots,d.
	\end{equation}
	The field $R[U]$ has the following properties.
	
	\begin{enumerate}
		\item The matrix $R[U]$ satisfies the support inclusion
		\begin{equation}
			\label{eq:inv4divEq:supp4R}
			\supp R[U]
			\subseteq
			\hat\Gamma_{\upsilon}(\bar\tau,\bar\varepsilon;t_0,x_0).
		\end{equation}
		In particular, by Lemma~\ref{lem:comparison-Gamma-C}, the corresponding
		Eulerian support inclusion is
		\begin{equation}
			\label{eq:inv4divEq:supp4R-eulerian}
			\supp R[U]
			\subseteq
			\hat C_{\upsilon}
			\left(
			\bar\tau,
			e^{2\bar\tau\nrm{\nabla\upsilon}_0}\bar\varepsilon;
			t_0,x_0
			\right).
		\end{equation}
		
		\item The time support satisfies
		\begin{equation}
			\label{eq:time-supp-localR}
			\supp_t R[U]\subseteq\supp_t U.
		\end{equation}
		
		\item For $M=0$ and every $0\leq N\leq K^\sharp-2$,
		the H{\"o}lder estimate is
		\begin{equation}
			\begin{aligned}
				\nrm{R[U]}_{N+\alpha}
				\lesssim
				\mC_U\bar\varepsilon^{1+\alpha}
				\Lk{N,K_0,\lambda_{R,U,\upsilon},\tilde\lambda_{R,U,\upsilon}}{s}\La{N+1,K_0,\lambda_{R,U,\upsilon},\tilde\lambda_{R,U,\upsilon}}
				\MM{N,K_0,\lambda_{R,U,\upsilon},\tilde\lambda_{R,U,\upsilon}}.
			\end{aligned}\label{eq:inv4divEq:Calpha-material-zero-MM}
		\end{equation}
		
		\item For all $M,N\in\mathbb N_0$ satisfying
		$1\leq M\leq \min\{M_*,K^*+1\}$ and $M+N\leq K^\sharp-1$, the H{\"o}lder estimate is
		\begin{equation}
			\begin{aligned}
				\nrm{D_{t,\upsilon}^{M}R[U]}_{N+\alpha}
				&\lesssim
				\mC_U\bar\varepsilon^{1+\alpha}
				\Lk{N,(K_U-M)_+,\lambda_{R,U,\upsilon},\tilde\lambda_{R,U,\upsilon}}{s}\La{N+1,(K_U-M)_+,\lambda_{R,U,\upsilon},\tilde\lambda_{R,U,\upsilon}}\\
				&\quad\cdot
				\MM{N,(K_U-M)_+,\lambda_{R,U,\upsilon},\tilde\lambda_{R,U,\upsilon}}
				\MM{M,N_t,\mu_{R,U,\upsilon},\tilde\mu_{R,U,\upsilon}}.
			\end{aligned}\label{eq:inv4divEq:Calpha-material-MM}
		\end{equation}
		
		\item For $M=0$ and every $0\leq N\leq K^\sharp-1$,
		the corresponding $C^N$ estimate is
		\begin{equation}
			\label{eq:inv4divEq:CN-material-zero-MM}
			\begin{aligned}
				\nrm{R[U]}_{N}
				&\lesssim
				\mC_U\bar\varepsilon^{1+\alpha}
				\Lk{N,K_0,\lambda_{R,U,\upsilon},\tilde\lambda_{R,U,\upsilon}}{s}
				\MM{N,K_0,\lambda_{R,U,\upsilon},\tilde\lambda_{R,U,\upsilon}}.
			\end{aligned}
		\end{equation}
		
		\item For all $M,N\in\mathbb N_0$ satisfying
		$1\leq M\leq \min\{M_*,K^*+1\}$ and $M+N\leq K^\sharp-1$, the $C^N$ estimate is
		\begin{equation}
			\label{eq:inv4divEq:CN-material-MM-weak}
			\begin{aligned}
				\nrm{D_{t,\upsilon}^{M}R[U]}_{N}
				&\lesssim
				\mC_U\bar\varepsilon^{1+\alpha}
				\Lk{N,(K_U-M)_+,\lambda_{R,U,\upsilon},\tilde\lambda_{R,U,\upsilon}}{s}\\
				&\quad\cdot
				\MM{N,(K_U-M)_+,\lambda_{R,U,\upsilon},\tilde\lambda_{R,U,\upsilon}}
				\MM{M,N_t,\mu_{R,U,\upsilon},\tilde\mu_{R,U,\upsilon}}.
			\end{aligned}
		\end{equation}
	\end{enumerate}

	If $U$ instead satisfies the following $C^N$ estimates
	\begin{equation}\label{eq:inv4divEq:integer-source-assumption}
		\nrm{D_{t,\upsilon}^{M}U}_{N}
		\leq\mC_U\MM{N,(K_U-M)_+,\lambda_U,\tilde\lambda_U}
		\MM{M,N_t,\mu_U,\tilde\mu_U},
	\end{equation}
	for $M,N\in\mathbb N_0$ with $0\leq M\leq M_*$ and $M+N\leq K^\sharp-1$,
	then, in the respective ranges of \eqref{eq:inv4divEq:CN-material-zero-MM} and \eqref{eq:inv4divEq:CN-material-MM-weak}, we have
	\begin{align}
		\nrm{R[U]}_{N}
		&\lesssim\mC_U\bar\varepsilon\MM{N,K_0,\lambda_{R,U,\upsilon},\tilde\lambda_{R,U,\upsilon}},
		\label{eq:inv4divEq:CN-from-integer-source-zero}\\
		\nrm{D_{t,\upsilon}^{M}R[U]}_{N}
		&\lesssim\mC_U\bar\varepsilon\MM{N,(K_U-M)_+,\lambda_{R,U,\upsilon},\tilde\lambda_{R,U,\upsilon}}
		\MM{M,N_t,\mu_{R,U,\upsilon},\tilde\mu_{R,U,\upsilon}}.
		\label{eq:inv4divEq:CN-from-integer-source}
	\end{align}
	The implicit constant may depend on $d,K^*,K^\sharp,\alpha$.
\end{thm}

\begin{rem}
	In the applications below, the source $U$ in Theorem~\ref{thm:inv4divEq}
	is often obtained by multiplying a given vector field $\widetilde U$ by a
	space-time cutoff of the form
	$$
	U(t,x)
	=
	\theta(\bar\tau^{-1}(t-t_0))
	\chi(\varepsilon^{-1}(\Phi_\upsilon(t,x;t_0)-x_0))
	\widetilde U(t,x).
	$$
	Assume that $\theta\in C_c^\infty(\R)$, $\chi\in C_c^\infty(\R^d)$,
	$\supp\theta\subseteq[-c_\theta,c_\theta]$ and
	$\supp\chi\subseteq B(0,c_\chi)$. Then
	$$
	\supp U
	\subseteq
	\supp\bigl(\theta(\bar\tau^{-1}(t-t_0))
	\chi(\varepsilon^{-1}(\Phi_\upsilon(t,x;t_0)-x_0))\bigr)
	\subseteq
	\{(t,x):|t-t_0|\leq c_\theta\bar\tau,\ 
	|\Phi_\upsilon(t,x;t_0)-x_0|<c_\chi\varepsilon\}.
	$$
	In particular,
	$
	\supp U
	\subseteq
	\hat\Gamma_{\upsilon}
	(c_\theta\bar\tau,c_\chi\varepsilon;t_0,x_0).
	$
	Theorem~\ref{thm:inv4divEq} applies provided all its hypotheses hold with
	these localization radii. The zero-mean compatibility condition
	must be verified separately.
\end{rem}
\begin{rem}\label{rem:local-inverse-holder-endpoint}
	If \eqref{eq:inv4divEq:est4nbv-material-Calpha-MM} also holds for
	$M=0$ and $N=K^\sharp-1$, then
	\eqref{eq:inv4divEq:Calpha-material-zero-MM} extends to
	$0\leq N\leq K^\sharp-1$. Together with
	\eqref{eq:inv4divEq:Calpha-material-MM}, we obtain the
	H{\"o}lder estimates for
	$
	0\leq M\leq \min\{M_*,K^*+1\},
	$
	$
	M+N\leq K^\sharp-1.
	$
	The threshold for spatial derivatives at $M=0$ remains $K_0$.
\end{rem}
\begin{rem}
	The construction of $R$ uses the flow of $\upsilon$ from $t_0$
	and the averaging kernel $\zeta_{x_0,\bar\varepsilon}$.
	The time radius $\bar\tau$ specifies the interval of application.
	We suppress these parameters when they are fixed.
	The localization radii are specified in each application.
	When $(t_0,x_0)=(t_u,x_\kappa)$, we write $R^{(\upsilon)}_{u,\kappa}$
	if this dependence needs to be specified. In the Newton step,
	$\upsilon=\overline u_q$, $t_u=u\tau_q$, and
	$x_\kappa=2\pi\tau_q\kappa$, so the corresponding operator is denoted
	by $R^{(\overline u_q)}_{u,\kappa}$.
\end{rem}

\subsection{Construction and basic properties of \texorpdfstring{$R^*$}{R*}}
\label{subsec:derivation-Rstar}
\label{subsec:formula-Rstar-basic-properties}

A first-order Taylor expansion in Fourier space gives a solution operator
that preserves support in the given ball, but the resulting field may be singular at the origin.
Following Isett--Oh~\cite{Isett16}, we conjugate the operator by translations
and take a smooth average to obtain $R^*$.

For simplicity, we omit the time variable and consider a smooth vector field
$V:\R^d\to\R^d$ compactly supported in $B(0,\bar\varepsilon)$. We write $x$
for the spatial variable, $\eta$ for the averaging parameter, and $\xi$ for
the Fourier variable.
By the first-order Taylor expansion of $\widehat V_i(\xi)$ around
$\xi=0$, we have
$$
\widehat V_i(\xi)
=
\widehat V_i(0)
+
\left(
\int_0^1 \partial_l\widehat V_i(\sgm\xi)\,\rd\sgm
\right)\xi_l .
$$
Under the zero-mean condition
\begin{equation}
	\label{eq:zeroMomentum}
	\widehat V_i(0)
	=
	\int_{\R^d}V_i(x)\,\rd x
	=
	0,
	\qquad
	i=1,\ldots,d,
\end{equation}
this becomes
$$
\widehat V_i(\xi)
=
i\xi_l
\left(
\int_0^1 \frac{1}{i}\partial_l\widehat V_i(\sgm\xi)\,\rd\sgm
\right).
$$
Taking the inverse Fourier transform gives the formal expression
for $r[V]=(r_{il}[V])$:
\begin{equation}
	\label{eq:rjl0:formal}
	\begin{aligned}
		r_{il}[V](x)
		&:=
		\cF^{-1}
		\left[
		\int_0^1 \frac{1}{i}\partial_l\widehat V_i(\sgm\xi)\,\rd\sgm
		\right](x)
		=
		-\int_0^1\frac{x_l}{\sgm}V_i\left(\frac{x}{\sgm}\right)\frac{\rd\sgm}{\sgm^d}.
	\end{aligned}
\end{equation}
The integral in \eqref{eq:rjl0:formal} is interpreted as a distribution.
After the change of variables $x=\sgm y$, its pairing with a test function
is integrable in $\sgm$. Its divergence satisfies
$$
\partial_lr_{il}[V]
=
V_i
-
\left(
\int_{\R^d}V_i(x)\,\rd x
\right)\dlt_0
=
V_i,
$$
where the last equality uses \eqref{eq:zeroMomentum}. Thus $r[V]$ is a
distributional solution of $\partial_lr_{il}[V]=V_i$. 

For $x\ne0$, \eqref{eq:rjl0:formal} represents a weighted integral of
$V_i$ along the ray starting at $x$ and pointing away from the origin. Since
$\supp V\subset B(0,\bar\varepsilon)$, we obtain
\begin{equation*}
	\supp r_{il}[V]\subset B(0,\bar\varepsilon).
\end{equation*}
The distribution $r[V]$ may be singular at $x=0$. Using translation
invariance of the divergence equation, we conjugate $r$ by translations
and take a smooth average of the resulting operators.
Let $T_\eta f(x):=f(x-\eta)$ denote the
translation operator. For each $\eta\in\R^d$, define the conjugated operator by
$
{}^{(\eta)}r_{il}[V]:=T_\eta\bigl(r_{il}[T_{-\eta}V]\bigr).
$
Equivalently,
\begin{equation}
	\label{eq:translated-ray-operator}
	{}^{(\eta)}r_{il}[V](x)
	=
	-\int_0^1\frac{(x-\eta)_l}{\sgm}V_i\left(\frac{x-\eta}{\sgm}+\eta\right)\frac{\rd\sgm}{\sgm^d}.
\end{equation}
By translation invariance and the zero-mean condition
\eqref{eq:zeroMomentum}, we have
$$
\partial_l{}^{(\eta)}r_{il}[V]=V_i.
$$

For $\eta\in B(0,\bar\varepsilon)$, the distribution
${}^{(\eta)}r[V]$ satisfies the same support property. Indeed, if the integrand
in \eqref{eq:translated-ray-operator} is nonzero, then
$$
x=\sgm z+(1-\sgm)\eta,
\qquad
z\in\supp V,\quad 0\leq\sgm\leq1.
$$
By convexity of $B(0,\bar\varepsilon)$, we obtain
$
x\in B(0,\bar\varepsilon).
$
Hence
$$
\supp{}^{(\eta)}r[V]\subset B(0,\bar\varepsilon)
\qquad
\text{for every }\eta\in B(0,\bar\varepsilon).
$$

Let $\zeta\in C_c^\infty(B(0,\bar\varepsilon))$ satisfy
$
\int_{\R^d}\zeta(\eta)\,\rd\eta=1.
$
Define
\begin{equation}
	\label{eq:Rstar-formula-o}
	\begin{aligned}
		R^*_{il}[V](x)
		&:=
		\int_{\R^d}{}^{(\eta)}r_{il}[V](x)\zeta(\eta)\,\rd\eta=
		-\int_{\R^d}
		\int_0^1
		\zeta(\eta)\frac{(x-\eta)_l}{\sgm}
		V_i\left(\frac{x-\eta}{\sgm}+\eta\right)
		\frac{\rd\sgm}{\sgm^d}\,\rd\eta.
	\end{aligned}
\end{equation}
Averaging the divergence identity and using $\int_{\R^d}\zeta(\eta)\,\rd\eta=1$ yields
$$
\partial_lR^*_{il}[V]=V_i.
$$
Since $\supp\zeta$ and $\supp V$ are compact subsets of the convex ball
$B(0,\bar\varepsilon)$, the averaged operator $R^*$ satisfies the support property
\begin{equation*}
	\supp R^*[V]\subseteq B(0,\bar\varepsilon).
\end{equation*}
By translating the entire construction, the same operator can be defined on any
ball $B(x_0,\bar\varepsilon)$, with the support of $R^*[V]$ remaining inside
that ball whenever $\supp V\subset B(x_0,\bar\varepsilon)$.

Choose once and for all a nonnegative function
$\zeta_0\in C_c^\infty(B(0,1))$ such that
$\int_{\R^d}\zeta_0(x)\,\rd x=1$, and set
$$
\zeta_{x_0,\bar\varepsilon}(\eta)
:=
\bar\varepsilon^{-d}
\zeta_0\left(\frac{\eta-x_0}{\bar\varepsilon}\right).
$$
For brevity, write $\zeta:=\zeta_{x_0,\bar\varepsilon}$. Then
\begin{equation}
	\label{eq:bnd4zeta}
	\int_{\R^d}\zeta(\eta)\,\rd\eta=1,
	\qquad
	\nrm{\partial^\gamma\zeta}_0
	\lesssim_N
	\bar\varepsilon^{-d-N},
	\qquad
	|\gamma|=N,\quad N\geq0,
\end{equation}
and, for every $0<\alpha<1$,
\begin{equation*}
	\nrm{\partial^\gamma\zeta}_0
	+
	\bar\varepsilon^\alpha[\partial^\gamma\zeta]_\alpha
	\lesssim_{N,\alpha}
	\bar\varepsilon^{-d-N},
	\qquad
	|\gamma|=N,\quad N\geq0.
\end{equation*}
For a smooth vector field $V:\R^d\to\R^d$ supported in
$B(x_0,\bar\varepsilon)$, define
\begin{equation}
	\label{eq:Rstar-formula}
	R^*_{il}[V](x)
	:=
	\int_{\R^d}{}^{(\eta)}r_{il}[V](x)\zeta(\eta)\,\rd\eta
	=
	-\int_{\R^d}\int_0^1 \zeta(\eta)\frac{(x-\eta)_l}{\sgm}
	V_i\left(\frac{x-\eta}{\sgm}+\eta\right)\frac{\rd\sgm}{\sgm^d}\,\rd\eta.
\end{equation}
The next proposition records the divergence identity, support and derivative properties, and commutation of $R^*$ with time derivatives.

\begin{prop}
	\label{prop:Rstar}
	Let $V$ be a smooth vector field satisfying
	$$
	\supp V\subset B(x_0,\bar\varepsilon),
	\qquad
	\int_{\R^d}V_i(x)\,\rd x=0,
	\qquad
	i=1,\ldots,d.
	$$
	Then the operator $R^*$ defined by \eqref{eq:Rstar-formula} has the following
	properties.
	\begin{enumerate}
		\item $R^*[V]$ is linear in $V$.
		
		\item $R^*[V]$ solves the nonsymmetric divergence equation
		\begin{equation}
			\label{eq:Rstar-div}
			\partial_lR^*_{il}[V]=V_i.
		\end{equation}
		
		\item The matrix $R^*[V]$ satisfies the support inclusion
		\begin{equation}
			\label{eq:supp4Rstar}
			\supp R^*[V]\subset B(x_0,\bar\varepsilon).
		\end{equation}
		
		\item For every multi-index $\gamma$, one has
		\begin{equation}
			\label{eq:dRstar}
			\begin{aligned}
				\pa^\gamma R^*_{il}[V](x)
				=&
				-\sum_{\gamma_1+\gamma_2=\gamma}\binom{\gamma}{\gamma_1}
				\int_{\R^d}\int_0^1
				(\pa^{\gamma_1}\zeta)(\eta)\frac{(x-\eta)_l}{\sgm}
				(\pa^{\gamma_2}V_i)
				\left(\frac{x-\eta}{\sgm}+\eta\right)
				\frac{\rd\sgm}{\sgm^d}\,\rd\eta .
			\end{aligned}
		\end{equation}
		
		\item If $V=V(t,\cdot)$ depends smoothly on time, the time-independent
		operator $R^*$ commutes with ordinary time derivatives, so
		\begin{equation}
			\label{eq:Rstar-time-commute}
			\partial_t^M R^*[V(t,\cdot)]
			=
			R^*[\partial_t^M V(t,\cdot)],
			\qquad M\in\mathbb N_0.
		\end{equation}
	\end{enumerate}
\end{prop}

\begin{proof}
	Linearity follows from \eqref{eq:Rstar-formula}. We first
	prove the divergence identity \eqref{eq:Rstar-div}. For each fixed
	$\eta\in\R^d$, translation invariance and the zero-mean assumption
	on $V$ give, in the sense of distributions on $\R^d$,
	$
	\partial_l{}^{(\eta)}r_{il}[V]=V_i.
	$
	Averaging this identity against $\zeta(\eta)\,\rd\eta$, we obtain
	the distributional identity
	$$
	\partial_lR^*_{il}[V]
	=
	\int_{\R^d}\partial_l{}^{(\eta)}r_{il}[V]\zeta(\eta)\,\rd\eta
	=
	V_i\int_{\R^d}\zeta(\eta)\,\rd\eta
	=
	V_i.
	$$
	Lemma~\ref{lem:keyLp} below justifies differentiation under the integral sign
	and passage to the limit as the lower integration cutoff $\delta\rightarrow0$. In particular, $R^*[V]$
	is smooth and the distributional identity holds classically.
	The support property \eqref{eq:supp4Rstar} follows from the convexity
	argument above after translation to $B(x_0,\bar\varepsilon)$.
	
	Fix $m\in\{1,\ldots,d\}$. We write
	$$
	\partial_mR^*_{il}[V](x)
	=
	\left.
	\frac{\partial}{\partial z^m}
	R^*_{il}[V](x+z)
	\right|_{z=0}.
	$$
	Using \eqref{eq:Rstar-formula}, we have
	$$
	R^*_{il}[V](x+z)
	=
	-\int_{\R^d}
	\int_0^1
	\zeta(\eta)\frac{(x+z-\eta)_l}{\sgm}
	V_i\left(\frac{x+z-\eta}{\sgm}+\eta\right)
	\frac{\rd\sgm}{\sgm^d}\,\rd\eta .
	$$
	Make the change of variables $\widetilde\eta=\eta-z$. Then
	$$
	R^*_{il}[V](x+z)
	=
	-\int_{\R^d}
	\int_0^1
	\zeta(\widetilde\eta+z)\frac{(x-\widetilde\eta)_l}{\sgm}
	V_i\left(\frac{x-\widetilde\eta}{\sgm}+\widetilde\eta+z\right)
	\frac{\rd\sgm}{\sgm^d}\,\rd\widetilde\eta .
	$$
	Differentiating with respect to $z_m$
	at $z=0$ gives
	$$
	\begin{aligned}
		\partial_mR^*_{il}[V](x)
		&=
		-\int_{\R^d}
		\int_0^1
		(\partial_m\zeta)(\eta)\frac{(x-\eta)_l}{\sgm}
		V_i\left(\frac{x-\eta}{\sgm}+\eta\right)
		\frac{\rd\sgm}{\sgm^d}\,\rd\eta
		-\int_{\R^d}
		\int_0^1
		\zeta(\eta)\frac{(x-\eta)_l}{\sgm}
		(\partial_mV_i)\left(\frac{x-\eta}{\sgm}+\eta\right)
		\frac{\rd\sgm}{\sgm^d}\,\rd\eta .
	\end{aligned}
	$$
	After the change of variables, $(x-\eta)_l/\sgm$ is independent of
	$z$. Differentiation in $z_m$ therefore acts only on $\zeta$ and $V_i$.
	Repeating the calculation inductively yields
	$$
	\begin{aligned}
		\pa^\gamma R^*_{il}[V](x)
		&=
		-\sum_{\gamma_1+\gamma_2=\gamma}\binom{\gamma}{\gamma_1}
		\int_{\R^d}
		\int_0^1
		(\pa^{\gamma_1}\zeta)(\eta)\frac{(x-\eta)_l}{\sgm}
		(\pa^{\gamma_2}V_i)
		\left(\frac{x-\eta}{\sgm}+\eta\right)
		\frac{\rd\sgm}{\sgm^d}\,\rd\eta .
	\end{aligned}
	$$
	Thus \eqref{eq:dRstar} holds.
	
	Since $R^*$ is independent of time,
	$
	\partial_tR^*_{il}[V(t,\cdot)]
	=
	R^*_{il}[\partial_tV(t,\cdot)].
	$
	Iterating this identity, we have
	$$
	\partial_t^M R^*[V(t,\cdot)]
	=
	R^*[\partial_t^M V(t,\cdot)],
	\qquad M\in\mathbb N_0.
	$$
	This proves \eqref{eq:Rstar-time-commute}.
\end{proof}

\subsection{Estimates for the operator \texorpdfstring{$R^*$}{R*}}

We first establish the integral estimates for the spatial derivatives of $R^*$.

\begin{lem}
	\label{lem:keyLp}
	Fix $\bar\varepsilon>0$, $x_0\in\R^d$, $N_1,N_2\in\mathbb N_0$,
	and a multi-index $\kappa\in\mathbb N_0^d$ with $|\kappa|=N_2$.
	Let $\widetilde{\zeta}\in C_c^\infty(\R^d)$ satisfy
	$
	\supp \widetilde{\zeta}\subseteq B(x_0,\bar\varepsilon).
	$
	Assume that there are constants $A_{N_1},A_{N_1+1}\geq0$ such that,
	for $N=N_1,N_1+1$ and every multi-index
	$\gamma\in\mathbb N_0^d$ with $|\gamma|=N$,
	\begin{equation}
		\label{eq:keyLp:hyp-N}
		\nrm{\partial^\gamma\widetilde{\zeta}}_0
		\leq
		A_N\bar\varepsilon^{-d-N}.
	\end{equation}
	Then the following estimates hold.
	
	\begin{enumerate}
		\item For any $f\in L_x^\infty$ supported in
		$B(x_0,\bar\varepsilon)$, and for every $|\gamma|=N_1$, one has
		\begin{equation}
			\label{eq:keyLp:1}
			\operatorname*{ess\,sup}_{x\in\R^d,\ \sgm\in(0,1]}
			\left|
			\int_{\R^d}
			\partial^\gamma\widetilde{\zeta}(\eta)
			\left(\frac{x-\eta}{\sgm}\right)^\kappa
			f\left(\frac{x-\eta}{\sgm}+\eta\right)
			\frac{\rd\eta}{\sgm^d}
			\right|
			\leq
			C_{N_1,N_2}A_{N_1}
			\bar\varepsilon^{N_2-N_1}
			\|f\|_{L_x^\infty}.
		\end{equation}
		
		\item Consequently, for any $f\in L_x^\infty$ supported in
		$B(x_0,\bar\varepsilon)$, and for every $|\gamma|=N_1$, one has
		\begin{equation}
			\label{eq:keyLp:2}
			\left\|
			\int_0^1
			\int_{\R^d}
			\partial^\gamma\widetilde{\zeta}(\eta)
			\left(\frac{x-\eta}{\sgm}\right)^\kappa
			f\left(\frac{x-\eta}{\sgm}+\eta\right)
			\frac{\rd\eta}{\sgm^d}\,\rd\sgm
			\right\|_{L_x^\infty}
			\leq
			C_{N_1,N_2}A_{N_1}
			\bar\varepsilon^{N_2-N_1}
			\|f\|_{L_x^\infty}.
		\end{equation}
		
		\item Let $0<\alpha<1$. Then, for any $f\in C^\alpha(\R^d)$ supported in
		$B(x_0,\bar\varepsilon)$, and for every $|\gamma|=N_1$, one has
		\begin{equation}
			\label{eq:keyLp-holder}
			\begin{aligned}
				&
				\left[
				\int_0^1
				\int_{\R^d}
				\partial^\gamma\widetilde{\zeta}(\eta)
				\left(\frac{x-\eta}{\sgm}\right)^\kappa
				f\left(\frac{x-\eta}{\sgm}+\eta\right)
				\frac{\rd\eta}{\sgm^d}\,\rd\sgm
				\right]_{\alpha}\\
				\leq&
				C_{N_1,N_2,\alpha}
				A_{N_1}
				\bar\varepsilon^{N_2-N_1}
				\nrm{f}_{\alpha}
				+
				C_{N_1,N_2,\alpha}
				(A_{N_1}+A_{N_1+1})
				\bar\varepsilon^{N_2-N_1-\alpha}
				\nrm{f}_{0}.
			\end{aligned}
		\end{equation}
		Consequently,
		\begin{equation}
			\label{eq:keyLp-holder-norm}
			\begin{aligned}
				&
				\left\|
				\int_0^1
				\int_{\R^d}
				\partial^\gamma\widetilde{\zeta}(\eta)
				\left(\frac{x-\eta}{\sgm}\right)^\kappa
				f\left(\frac{x-\eta}{\sgm}+\eta\right)
				\frac{\rd\eta}{\sgm^d}\,\rd\sgm
				\right\|_{\alpha}\\
				\leq&
				C_{N_1,N_2,\alpha}
				A_{N_1}
				\bar\varepsilon^{N_2-N_1}
				\nrm{f}_{\alpha}
				+
				C_{N_1,N_2,\alpha}
				(A_{N_1}+A_{N_1+1})
				\bar\varepsilon^{N_2-N_1-\alpha}
				\nrm{f}_{0}.
			\end{aligned}
		\end{equation}
	\end{enumerate}

	The constants $C_{N_1,N_2}$ depend only on $d,N_1,N_2$.
	The constants $C_{N_1,N_2,\alpha}$ may also depend on $\alpha$.
	They are independent of $x_0,\bar\varepsilon,f,\widetilde\zeta,A_{N_1},A_{N_1+1}$.
\end{lem}

\begin{proof}
	We establish \eqref{eq:keyLp:1}. Fix $x\in\R^d$ and
	$\sgm\in(0,1]$, and set
	$
	z=(x-\eta)/\sgm+\eta.
	$
	On the support of the integrand, both $\eta$ and $z$ belong to
	$B(x_0,\bar\varepsilon)$. Thus
	$
	x=\sgm z+(1-\sgm)\eta\in B(x_0,\bar\varepsilon)
	$
	and
	$$
	\left|
	\frac{x-\eta}{\sgm}
	\right|
	=
	|z-\eta|
	\leq
	2\bar\varepsilon.
	$$
	In particular, since $|\kappa|=N_2$,
	$$
	\left|
	\left(\frac{x-\eta}{\sgm}\right)^\kappa
	\right|
	\leq
	C_{N_2}\bar\varepsilon^{N_2}.
	$$
	
	For $0<\sgm\leq1/2$, we estimate the integral after the change
	of variables
	$
	z=(x-\eta)/\sgm+\eta.
	$
	Since
	$
	\eta=(x-\sgm z)/(1-\sgm)
	$
	and
	$
	\rd\eta/\sgm^d=\rd z/(1-\sgm)^d,
	$
	we have $(z-x)/(1-\sgm)=z-\eta$. On the support of the transformed
	integrand, $z,\eta\in B(x_0,\bar\varepsilon)$. Thus
	$$
	\begin{aligned}
		\left|
		\int_{\R^d}
		\partial^\gamma\widetilde{\zeta}(\eta)
		\left(\frac{x-\eta}{\sgm}\right)^\kappa
		f\left(\frac{x-\eta}{\sgm}+\eta\right)
		\frac{\rd\eta}{\sgm^d}
		\right|
		&=
		\left|
		\int_{\R^d}
		\partial^\gamma\widetilde{\zeta}
		\left(\frac{x-\sgm z}{1-\sgm}\right)
		\left(\frac{z-x}{1-\sgm}\right)^\kappa
		f(z)\frac{\rd z}{(1-\sgm)^d}
		\right|\\
		&\leq
		C_{N_2}A_{N_1}\bar\varepsilon^{-d-N_1}
		(2\bar\varepsilon)^{N_2}(1-\sgm)^{-d}
		\int_{B(x_0,\bar\varepsilon)}|f(z)|\,\rd z\\
		&\leq
		C_{N_1,N_2}
		A_{N_1}\bar\varepsilon^{N_2-N_1}\|f\|_{L_x^\infty},
	\end{aligned}
	$$
	where we used $\supp f\subset B(x_0,\bar\varepsilon)$ and
	$1-\sgm\geq1/2$.
	For $1/2\leq\sgm\leq1$, we estimate directly in the $\eta$ variable.
	Since $\sgm^{-d}\leq2^d$, the same support argument gives
	$$
	\begin{aligned}
		\left|
		\int_{\R^d}
		\partial^\gamma\widetilde{\zeta}(\eta)
		\left(\frac{x-\eta}{\sgm}\right)^\kappa
		f\left(\frac{x-\eta}{\sgm}+\eta\right)
		\frac{\rd\eta}{\sgm^d}
		\right|
		\leq
		C_{N_2}A_{N_1}\bar\varepsilon^{-d-N_1}
		(2\bar\varepsilon)^{N_2}\sgm^{-d}
		\|f\|_{L_x^\infty}|B(x_0,\bar\varepsilon)|
		\leq
		C_{N_1,N_2}
		A_{N_1}\bar\varepsilon^{N_2-N_1}\|f\|_{L_x^\infty}.
	\end{aligned}
	$$
	Since $|B(x_0,\bar\varepsilon)|=|B(0,1)|\bar\varepsilon^d$, the two
	estimates prove \eqref{eq:keyLp:1}. Integration in $\sgm\in(0,1)$ gives
	\eqref{eq:keyLp:2}.
	
	We next prove \eqref{eq:keyLp-holder}. The assumptions
	\eqref{eq:keyLp:hyp-N}, together with the mean value theorem, imply
	$$
	[\partial^\gamma\widetilde{\zeta}]_\alpha
	\leq
	C_{N_1,\alpha}
	(A_{N_1}+A_{N_1+1})
	\bar\varepsilon^{-d-N_1-\alpha},
	\qquad
	|\gamma|=N_1.
	$$
	Indeed, the estimate for derivatives of order $N_1+1$ applies when
	$|h|\leq\bar\varepsilon$, while the estimate for derivatives of order $N_1$ applies when
	$|h|>\bar\varepsilon$.
	
	For a fixed multi-index $\gamma$ with $|\gamma|=N_1$, write
	$$
	Tf(x)
	:=
	\int_0^1\int_{\R^d}
	\partial^\gamma\widetilde{\zeta}(\eta)
	\left(\frac{x-\eta}{\sgm}\right)^\kappa
	f\left(\frac{x-\eta}{\sgm}+\eta\right)
	\frac{\rd\eta}{\sgm^d}\,\rd\sgm.
	$$
	If $|h|>\bar\varepsilon$, then \eqref{eq:keyLp:2} gives
	$$
	\begin{aligned}
		|Tf(x+h)-Tf(x)|
		\leq
		2\nrm{Tf}_0
		\leq
		C_{N_1,N_2}
		A_{N_1}\bar\varepsilon^{N_2-N_1}\nrm{f}_0
		\leq
		C_{N_1,N_2}
		A_{N_1}\bar\varepsilon^{N_2-N_1-\alpha}
		\nrm{f}_0|h|^\alpha.
	\end{aligned}
	$$
	
	Consider $|h|\leq\bar\varepsilon$. In the formula for
	$Tf(x+h)$, make the change of variables $\eta\mapsto\eta+h$. This
	translation leaves the factor $((x-\eta)/\sgm)^\kappa$ unchanged:
	$$
	\begin{aligned}
		Tf(x+h)
		=
		\int_0^1\int_{\R^d}
		\partial^\gamma\widetilde{\zeta}(\eta+h)
		\left(\frac{x-\eta}{\sgm}\right)^\kappa
		f\left(\frac{x-\eta}{\sgm}+\eta+h\right)
		\frac{\rd\eta}{\sgm^d}\,\rd\sgm.
	\end{aligned}
	$$
	Subtracting the formula for $Tf(x)$,
	we obtain
	$$
	\begin{aligned}
		Tf(x+h)-Tf(x)
		={}&
		\int_0^1\int_{\R^d}
		\left[
		\partial^\gamma\widetilde{\zeta}(\eta+h)
		-
		\partial^\gamma\widetilde{\zeta}(\eta)
		\right]
		\left(\frac{x-\eta}{\sgm}\right)^\kappa
		f\left(\frac{x-\eta}{\sgm}+\eta+h\right)
		\frac{\rd\eta}{\sgm^d}\,\rd\sgm\\
		&+
		\int_0^1\int_{\R^d}
		\partial^\gamma\widetilde{\zeta}(\eta)
		\left(\frac{x-\eta}{\sgm}\right)^\kappa\
		\left[
		f\left(\frac{x-\eta}{\sgm}+\eta+h\right)
		-
		f\left(\frac{x-\eta}{\sgm}+\eta\right)
		\right]
		\frac{\rd\eta}{\sgm^d}\,\rd\sgm.
	\end{aligned}
	$$
	Whenever either integrand is nonzero, the averaging variable and
	the relevant argument of $f$ lie in $B(x_0,2\bar\varepsilon)$. Moreover,
	$$
	\begin{aligned}
		\left\|
		\partial^\gamma\widetilde{\zeta}(\cdot+h)
		-
		\partial^\gamma\widetilde{\zeta}
		\right\|_0
		\leq
		C_{N_1,\alpha}
		(A_{N_1}+A_{N_1+1})
		\bar\varepsilon^{-d-N_1-\alpha}|h|^\alpha,\quad
		\nrm{f(\cdot+h)-f}_0
		\leq
		[f]_\alpha|h|^\alpha.
	\end{aligned}
	$$
	Apply the argument used in the proof of \eqref{eq:keyLp:2} separately to
	the two terms, with $\bar\varepsilon$ replaced by $2\bar\varepsilon$. We obtain
	$$
	\begin{aligned}
		|Tf(x+h)-Tf(x)|
		\leq
		C_{N_1,N_2,\alpha}
		A_{N_1}\bar\varepsilon^{N_2-N_1}
		\nrm{f}_\alpha|h|^\alpha
		+
		C_{N_1,N_2,\alpha}
		(A_{N_1}+A_{N_1+1})
		\bar\varepsilon^{N_2-N_1-\alpha}
		\nrm{f}_0|h|^\alpha.
	\end{aligned}
	$$
	The H{\"o}lder seminorm estimate
	\eqref{eq:keyLp-holder} follows. Combining it with
	\eqref{eq:keyLp:2} gives \eqref{eq:keyLp-holder-norm}.
\end{proof}

\eqref{eq:dRstar} and Lemma~\ref{lem:keyLp} give the following estimates.

\begin{lem}
	\label{lem:Rstar-spatial-holder}
	Fix $N_*\in\N$, $x_0\in\R^d$, and $0<\bar\varepsilon\leq1$.
	Let $V$ be a smooth vector field supported in
	$B(x_0,\bar\varepsilon)$ with zero mean in each component.
	Then, for every $0\leq N\leq N_*$,
	one has
	\begin{equation}
		\label{eq:Rstar-CN-basic}
		\nrm{R^*[V]}_{N}
		\leq
		C\bar\varepsilon
		\sum_{r=0}^{N}\bar\varepsilon^{-(N-r)}\nrm{V}_{r}.
	\end{equation}
	In particular,
	\begin{equation}
		\label{eq:Rstar-C0-basic}
		\nrm{R^*[V]}_{0}
		\leq
		C\bar\varepsilon\nrm{V}_{0}.
	\end{equation}
	Moreover, for every $0\leq N\leq N_*$ and $0<\alpha<1$, one has
	\begin{equation}
		\label{eq:Rstar-Calpha-basic}
		\nrm{R^*[V]}_{N+\alpha}
		\leq
		C\bar\varepsilon
		\sum_{r=0}^{N}\bar\varepsilon^{-(N-r)}\nrm{V}_{r+\alpha}.
	\end{equation}
	In particular,
	\begin{equation}
		\label{eq:Rstar-alpha-basic}
		\nrm{R^*[V]}_{\alpha}
		\leq
		C\bar\varepsilon\nrm{V}_{\alpha}.
	\end{equation}
	
	Fix $K\in\N_0$, $s\in\N$, and $0<\alpha<1$. Suppose, in addition, that for some
	$\mC_V\geq0$ and $1\leq\lambda_V\leq\tilde\lambda_V$,
	\begin{equation}
		\label{eq:Rstar-input-Calpha-MM}
		\nrm{V}_{N+\alpha}
		\leq
		\mC_V
		\Lk{N,K,\lambda_V,\tilde\lambda_V}{s}
		\MM{N,K,\lambda_V,\tilde\lambda_V}
	\end{equation}
	for every $0\leq N\leq N_*$. Define
	$
	\lambda_{R^*,V}:=\max\{\bar\varepsilon^{-1},\lambda_V\},
	\
	\tilde\lambda_{R^*,V}:=\max\{\bar\varepsilon^{-1},\tilde\lambda_V\}.
	$
	Then, for every $0\leq N\leq N_*$,
	\begin{equation}
		\label{eq:Rstar-Calpha-MM}
		\nrm{R^*[V]}_{N+\alpha}
		\lesssim
		\mC_V\bar\varepsilon
		\Lk{N,K,\lambda_{R^*,V},\tilde\lambda_{R^*,V}}{s}
		\MM{N,K,\lambda_{R^*,V},\tilde\lambda_{R^*,V}}.
	\end{equation}
	Furthermore,
	\begin{equation}
		\label{eq:Rstar-CN-MM}
		\nrm{R^*[V]}_{N}
		\lesssim
		\mC_V\bar\varepsilon^{1+\alpha}
		\Lk{N,K,\lambda_{R^*,V},\tilde\lambda_{R^*,V}}{s}
		\MM{N,K,\lambda_{R^*,V},\tilde\lambda_{R^*,V}}.
	\end{equation}

	The constants may depend on $d,N_*$ and, in the H{\"o}lder estimates,
	on $\alpha$.
	They are independent of $x_0,\bar\varepsilon,V$ and the amplitude and frequency scales.
\end{lem}

\begin{proof}
	Fix $0\leq M\leq N$ and a multi-index $\gamma$ with $|\gamma|=M$.
	By \eqref{eq:dRstar}, up to fixed signs and combinatorial factors,
	every term in $\partial^\gamma R^*[V]$ is of the form
	$$
	\int_{\R^d}\int_0^1
	\partial^{\gamma_1}\zeta(\eta)
	\left(\frac{x-\eta}{\sgm}\right)_l
	\partial^{\gamma_2}V_i\left(\frac{x-\eta}{\sgm}+\eta\right)
	\frac{\rd\sgm}{\sgm^d}\,\rd\eta,
	$$
	where $\gamma_1+\gamma_2=\gamma$.
	Using \eqref{eq:bnd4zeta} and \eqref{eq:keyLp:2}
	with $\kappa=e_l$, $f=\partial^{\gamma_2}V_i$, $N_1=|\gamma_1|$, and $N_2=1$,
	we obtain
	\begin{equation*}
		\begin{aligned}
			\left\|
			\int_{\R^d}\int_0^1
			\partial^{\gamma_1}\zeta(\eta)
			\left(\frac{x-\eta}{\sgm}\right)_l
			\partial^{\gamma_2}V_i\left(\frac{x-\eta}{\sgm}+\eta\right)
			\frac{\rd\sgm}{\sgm^d}\,\rd\eta
			\right\|_0
			\leq
			C\bar\varepsilon^{1-|\gamma_1|}
			\nrm{\partial^{\gamma_2}V_i}_0
			=
			C\bar\varepsilon^{1+|\gamma_2|-M}
			\nrm{\partial^{\gamma_2}V_i}_0.
		\end{aligned}
	\end{equation*}
	Grouping the terms by $r=|\gamma_2|$ gives
	\begin{equation*}
		\nrm{\nabla^MR^*[V]}_0
		\leq
		C\bar\varepsilon
		\sum_{r=0}^{M}\bar\varepsilon^{-(M-r)}\nrm{V}_{r}.
	\end{equation*}
	Summing over $0\leq M\leq N$ and using
	$0<\bar\varepsilon\leq1$ proves \eqref{eq:Rstar-CN-basic}. Taking
	$N=0$ gives \eqref{eq:Rstar-C0-basic}.
	
	Apply \eqref{eq:keyLp-holder-norm} to the same terms. With
	$f=\partial^{\gamma_2}V_i$, $N_1=|\gamma_1|$, and $N_2=1$, we obtain
	\begin{equation}
		\label{eq:Rstar-kernel-Calpha-bound}
		\begin{aligned}
			\left\|
			\int_{\R^d}\int_0^1
			\partial^{\gamma_1}\zeta(\eta)
			\left(\frac{x-\eta}{\sgm}\right)_l
			\partial^{\gamma_2}V_i\left(\frac{x-\eta}{\sgm}+\eta\right)
			\frac{\rd\sgm}{\sgm^d}\,\rd\eta
			\right\|_{\alpha}
			\leq
			C\bar\varepsilon^{1-|\gamma_1|}
			\nrm{\partial^{\gamma_2}V_i}_{\alpha}
			+
			C\bar\varepsilon^{1-|\gamma_1|-\alpha}
			\nrm{\partial^{\gamma_2}V_i}_{0}.
		\end{aligned}
	\end{equation}
	Here the second term is absorbed using the support of $V$. Indeed, for every
	$x\in\supp\partial^{\gamma_2}V_i$, one may choose
	$x'\notin B(x_0,\bar\varepsilon)$ such that
	$|x-x'|\lesssim\bar\varepsilon$. Since
	$\partial^{\gamma_2}V_i(x')=0$, it follows that
	$$
	|\partial^{\gamma_2}V_i(x)|
	\leq
	[\partial^{\gamma_2}V_i]_{\alpha}|x-x'|^\alpha
	\lesssim
	\bar\varepsilon^\alpha[\partial^{\gamma_2}V_i]_{\alpha}.
	$$
	Consequently,
	\begin{equation}
		\label{eq:Rstar-small-support-holder-absorb}
		\bar\varepsilon^{-\alpha}
		\nrm{\partial^{\gamma_2}V_i}_{0}
		\lesssim
		\nrm{V}_{|\gamma_2|+\alpha}.
	\end{equation}
	Combining \eqref{eq:Rstar-kernel-Calpha-bound}--\eqref{eq:Rstar-small-support-holder-absorb}, we obtain
	\begin{equation*}
		\begin{aligned}
			\left\|
			\int_{\R^d}\int_0^1
			\partial^{\gamma_1}\zeta(\eta)
			\left(\frac{x-\eta}{\sgm}\right)_l
			\partial^{\gamma_2}V_i\left(\frac{x-\eta}{\sgm}+\eta\right)
			\frac{\rd\sgm}{\sgm^d}\,\rd\eta
			\right\|_{\alpha}
			\leq
			C\bar\varepsilon^{1-|\gamma_1|}
			\nrm{V}_{|\gamma_2|+\alpha}
			=
			C\bar\varepsilon^{1+|\gamma_2|-M}
			\nrm{V}_{|\gamma_2|+\alpha}.
		\end{aligned}
	\end{equation*}
	The sum of the terms in $\partial^\gamma R^*[V]$ satisfies
	\begin{equation*}
		\nrm{\nabla^MR^*[V]}_{\alpha}
		\leq
		C\bar\varepsilon
		\sum_{r=0}^{M}\bar\varepsilon^{-(M-r)}\nrm{V}_{r+\alpha}.
	\end{equation*}
	Summing over $0\leq M\leq N$ and using
	$0<\bar\varepsilon\leq1$ proves \eqref{eq:Rstar-Calpha-basic}. Taking
	$N=0$ yields \eqref{eq:Rstar-alpha-basic}.
	
	For the two-scale estimates, use
	\eqref{eq:Rstar-Calpha-basic} with
	\eqref{eq:Rstar-input-Calpha-MM} to write
	\begin{equation}
		\label{eq:Rstar-parameter-sum}
		\begin{aligned}
			\nrm{R^*[V]}_{N+\alpha}
			\leq
			C\bar\varepsilon
			\sum_{r=0}^{N}\bar\varepsilon^{-(N-r)}
			\nrm{V}_{r+\alpha}
			\lesssim
			\mC_V\bar\varepsilon
			\sum_{r=0}^{N}\bar\varepsilon^{-(N-r)}
			\Lk{r,K,\lambda_V,\tilde\lambda_V}{s}
			\MM{r,K,\lambda_V,\tilde\lambda_V}.
		\end{aligned}
	\end{equation}
	By the definitions of $\lambda_{R^*,V}$ and
	$\tilde\lambda_{R^*,V}$, both scales are at least
	$\bar\varepsilon^{-1}$ and dominate the corresponding scales of $V$.
	Hence, for every $0\leq r\leq N$,
	$$
	\begin{aligned}
		\bar\varepsilon^{-(N-r)}
		\Lk{r,K,\lambda_V,\tilde\lambda_V}{s}
		\MM{r,K,\lambda_V,\tilde\lambda_V}
		\leq
		\Lk{N,K,\lambda_{R^*,V},\tilde\lambda_{R^*,V}}{s}
		\MM{N,K,\lambda_{R^*,V},\tilde\lambda_{R^*,V}}.
	\end{aligned}
	$$
	Applying this inequality to each term in \eqref{eq:Rstar-parameter-sum} proves
	\eqref{eq:Rstar-Calpha-MM}.
	
	Finally, applying \eqref{eq:Rstar-small-support-holder-absorb} to each derivative
	of order at most $r$ gives, for $0\leq r\leq N$,
	$$
	\nrm{V}_{r}
	\lesssim
	\bar\varepsilon^\alpha\nrm{V}_{r+\alpha}.
	$$
	Substitution into \eqref{eq:Rstar-CN-basic} proves \eqref{eq:Rstar-CN-MM}.
\end{proof}

\subsection{Construction and estimates for the matrix field \texorpdfstring{$R[U]$}{R[U]}}
\label{subsec:construction-Eulerian-solution-operator}
\label{subsec:estimates-Eulerian-R-proof-main-theorem}

Let $U$ satisfy the hypotheses of
Theorem~\ref{thm:inv4divEq}. We construct $V$, $R^{\mathrm L}$, and $R[U]$ for
$|t-t_0|\leq\bar\tau$ and extend them smoothly by zero, using
\eqref{eq:inv4divEq:supp4U}. We first pull $U$ back to the fixed Lagrangian
coordinates by setting
\begin{equation}
	\label{eq:Rjl0}
	V_i(t,y)
	:=
	\JX(t,y)U_i(t,X_\upsilon(t,y;t_0)).
\end{equation}
Under the change of variables $x=X_\upsilon(t,y;t_0)$, one has
$\rd x=\JX(t,y)\,\rd y$. Therefore,
$$
\int_{\R^d}V_i(t,y)\,\rd y
=
\int_{\R^d}\JX(t,y)U_i(t,X_\upsilon(t,y;t_0))\,\rd y
=
\int_{\R^d}U_i(t,x)\,\rd x
=
0,
$$
where the last equality follows from \eqref{eq:inv4divEq:hyp:1}.

The pullback also places the support inside a fixed spatial ball.
Indeed, since $\JX$ is strictly positive and
$X_\upsilon(t,\cdot;t_0)$ is a diffeomorphism,
$
\supp V(t,\cdot)
=
\Phi_\upsilon(t,\supp U(t,\cdot);t_0).
$
If $y\in\supp V(t,\cdot)$, then
$(t,X_\upsilon(t,y;t_0))\in\supp U$. By
\eqref{eq:inv4divEq:supp4U},
$
(t,X_\upsilon(t,y;t_0))
\in
\hat\Gamma_\upsilon(\bar\tau,\bar\varepsilon;t_0,x_0),
$
and hence, by the definition of the adapted Lagrangian cylinder,
$$
y
=
\Phi_\upsilon(t,X_\upsilon(t,y;t_0);t_0)
\in
B(x_0,\bar\varepsilon).
$$
Thus
$$
\supp V(t,\cdot)\subset B(x_0,\bar\varepsilon).
$$

The zero-mean and support properties allow us to apply $R^*$ to
$V(t,\cdot)$. Define
\begin{equation}
	\label{eq:Rjl1}
	R^{\mathrm L}_{il}(t,y):=R^*_{il}[V(t,\cdot)](y).
\end{equation}
Then
\begin{equation}
	\label{eq:Rjl1-div}
	\partial_lR^{\mathrm L}_{il}(t,y)=V_i(t,y),
	\qquad
	\supp R^{\mathrm L}(t,\cdot)\subset B(x_0,\bar\varepsilon).
\end{equation}
Define the final operator by the Piola transform
\begin{equation}
	\label{eq:Rjl2}
	R_{ij}[U](t,x)
	:=
	P_{jl}(t,x)R^{\mathrm L}_{il}(t,\Phi_\upsilon(t,x;t_0)),
\end{equation}
where
$$
P_{ij}(t,x)
:=
\left[(\cof\nabla_x\Phi_\upsilon(t,x;t_0))^T\right]_{ij}
=
\JPhi(t,x)
\left[(\nabla_x\Phi_\upsilon)^{-1}(t,x;t_0)\right]_{ij}.
$$
The construction gives the following properties of $R[U]$.

\begin{prop}
	\label{prop:Rjl}
	Let $U$ satisfy the hypotheses of Theorem~\ref{thm:inv4divEq}, and
	define $R[U]=(R_{ij}[U])$ by \eqref{eq:Rjl0}--\eqref{eq:Rjl2}. Then
	$R[U]$ has the following properties.
	\begin{enumerate}
		\item The map $U\mapsto R[U]$ is linear.
		
		\item The matrix field $R[U]$ solves the nonsymmetric divergence equation
		\begin{equation}
			\label{eq:R-div}
			\partial_jR_{ij}[U]=U_i.
		\end{equation}
		
		\item The support of $R[U]$ satisfies
		\begin{equation}
			\label{eq:supp4Rjl}
			\supp R[U]
			\subseteq
			\hat\Gamma_\upsilon(\bar\tau,\bar\varepsilon;t_0,x_0).
		\end{equation}
		In particular, Lemma~\ref{lem:comparison-Gamma-C} gives
		\begin{equation}
			\label{eq:supp4Rjl-eulerian}
			\supp R[U]
			\subseteq
			\hat C_\upsilon
			\left(
			\bar\tau,
			e^{2\bar\tau\nrm{\nabla\upsilon}_0}\bar\varepsilon;
			t_0,x_0
			\right).
		\end{equation}
		
		\item The time support satisfies
		\begin{equation}
			\label{eq:time-supp-R}
			\supp_tR[U]\subseteq\supp_tU.
		\end{equation}
		
		\item For every multi-index $\gamma$, the spatial derivatives of $R[U]$
		satisfy
		\begin{equation}
			\label{eq:dRjl0}
			\pa^\gamma R_{ij}[U]
			=
			\sum_{\gamma_1+\gamma_2=\gamma}
			\binom{\gamma}{\gamma_1}
			(\pa^{\gamma_1}P_{jl})
			\pa^{\gamma_2}
			\left(
			R^*_{il}[V(t,\cdot)]\circ\Phi_\upsilon
			\right).
		\end{equation}
		
		\item For every $M\in\mathbb N_0$, the $M$-th material derivative of
		$R[U]$ is given by
		\begin{equation}
			\label{eq:higher-material-R}
			D_{t,\upsilon}^{M}R_{ij}[U]
			=
			\sum_{k=0}^{M}
			\binom{M}{k}
			(D_{t,\upsilon}^{M-k}P_{jl})
			\left(
			R^*_{il}[\partial_t^kV(t,\cdot)]\circ\Phi_\upsilon
			\right).
		\end{equation}
		In particular,
		\begin{equation}
			\label{eq:approxMdRjl}
			\begin{aligned}
				D_{t,\upsilon}R_{ij}[U]
				={}&
				(D_{t,\upsilon}P_{jl})
				\left(R^*_{il}[V(t,\cdot)]\circ\Phi_\upsilon\right)
				+
				P_{jl}
				\left(
				R^*_{il}[\partial_tV(t,\cdot)]\circ\Phi_\upsilon
				\right).
			\end{aligned}
		\end{equation}
	\end{enumerate}
\end{prop}

\begin{proof}
	The pullback, $R^*$, and the Piola transform are linear, so $U\mapsto R[U]$ is linear.
	Using the Piola identity
	$\partial_jP_{jl}=0$, we have
	$$
	\begin{aligned}
		\partial_jR_{ij}[U]
		=
		\partial_j\left[P_{jl}(R^{\mathrm L}_{il}\circ\Phi_\upsilon)\right]
		=
		P_{jl}(\partial_mR^{\mathrm L}_{il})\circ\Phi_\upsilon\,
		\partial_j(\Phi_\upsilon)^m.
	\end{aligned}
	$$
	The definition of $P$ gives
	$
	P_{jl}\partial_j(\Phi_\upsilon)^m
	=
	\JPhi\delta_{lm}.
	$
	By \eqref{eq:Rjl1-div}, we obtain
	$$
	\begin{aligned}
		\partial_jR_{ij}[U](t,x)
		=
		\JPhi(t,x)
		(\partial_lR^{\mathrm L}_{il})(t,\Phi_\upsilon(t,x;t_0))
		=
		\JPhi(t,x)
		\JX(t,\Phi_\upsilon(t,x;t_0))
		U_i(t,x)
		=
		U_i(t,x).
	\end{aligned}
	$$
	
	By \eqref{eq:Rjl1-div},
	$\supp R^{\mathrm L}(t,\cdot)\subset B(x_0,\bar\varepsilon)$. Hence
	\eqref{eq:Rjl2} implies
	$$
	\supp R[U](t,\cdot)
	\subseteq
	X_\upsilon(t,B(x_0,\bar\varepsilon);t_0).
	$$
	If $U(t,\cdot)=0$, then $V(t,\cdot)=0$, and the linearity of $R^*$
	implies $R^{\mathrm L}(t,\cdot)=0$ and $R[U](t,\cdot)=0$. Thus
	\eqref{eq:time-supp-R} holds. Combining these observations with
	\eqref{eq:inv4divEq:supp4U} gives \eqref{eq:supp4Rjl}. The Eulerian
	support inclusion \eqref{eq:supp4Rjl-eulerian} follows from
	Lemma~\ref{lem:comparison-Gamma-C} and
	$$
	\int_{t_0-\bar\tau}^{t_0+\bar\tau}
	\nrm{\nabla\upsilon(\sigma,\cdot)}_0\,\rd \sigma
	\leq
	2\bar\tau\nrm{\nabla\upsilon}_0.
	$$
	
	Equation~\eqref{eq:Rjl2} also gives \eqref{eq:dRjl0}.
	We finish with \eqref{eq:higher-material-R}. Since
	$D_{t,\upsilon}\Phi_\upsilon=0$, one has
	$
	D_{t,\upsilon}
	\left(R^{\mathrm L}_{il}\circ\Phi_\upsilon\right)
	=
	(\partial_tR^{\mathrm L}_{il})\circ\Phi_\upsilon.
	$
	Moreover, \eqref{eq:Rstar-time-commute} gives
	$$
	\partial_t^kR^{\mathrm L}_{il}(t,\cdot)
	=
	R^*_{il}[\partial_t^kV(t,\cdot)],
	\qquad
	k\in\mathbb N_0.
	$$
	The support of $V$ remains in the fixed ball, so
	$
	\supp\partial_t^kV(t,\cdot)
	\subset B(x_0,\bar\varepsilon),
	$
	while the zero-mean identity for $V$ implies
	$
	\int_{\R^d}\partial_t^kV_i(t,y)\,\rd y=0.
	$
	Applying $D_{t,\upsilon}^{M}$ to \eqref{eq:Rjl2} and using these
	identities, we obtain
	$$
	D_{t,\upsilon}^{M}R_{ij}[U]
	=
	\sum_{k=0}^{M}
	\binom{M}{k}
	(D_{t,\upsilon}^{M-k}P_{jl})
	(\partial_t^kR^{\mathrm L}_{il})\circ\Phi_\upsilon.
	$$
	This proves \eqref{eq:higher-material-R}. Taking $M=1$ gives
	\eqref{eq:approxMdRjl}.
\end{proof}

We next combine the pullback estimates with the estimates for $R^*$.

\begin{lem}
	\label{lem:estimate-pulled-back-sources}
	Fix $s\in\mathbb N$. Let $U$ satisfy the hypotheses of
	Theorem~\ref{thm:inv4divEq} with this $s$. Define $V$ by
	\eqref{eq:Rjl0} and $R[U]$ by
	\eqref{eq:Rjl0}--\eqref{eq:Rjl2}. Then, for every
	$M\in\mathbb N_0$,
	\begin{equation}
		\label{eq:supp-time-derivatives-V-pullback}
		\supp\partial_t^MV(t,\cdot)
		\subset B(x_0,\bar\varepsilon).
	\end{equation}
	For $M=0$ and every $0\leq N\leq K^\sharp-2$,
	\begin{align}
		\nrm{V}_{N+\alpha}
		&\lesssim
		\mC_U\bar\varepsilon^\alpha
		\Lk{N,K_0,\lambda_{R,U,\upsilon},\tilde\lambda_{R,U,\upsilon}}{s}\La{N+1,K_0,\lambda_{R,U,\upsilon},\tilde\lambda_{R,U,\upsilon}}
		\MM{N,K_0,\lambda_{R,U,\upsilon},\tilde\lambda_{R,U,\upsilon}}.
		\label{eq:estimate-V-pullback-material-zero}
	\end{align}
	For every $0\leq N\leq K^\sharp-1$, the corresponding $C^N$ estimate is
	\begin{equation}
		\label{eq:estimate-V-pullback-material-zero-CN}
		\nrm{V}_{N}
		\lesssim
		\mC_U\bar\varepsilon^\alpha
		\Lk{N,K_0,\lambda_{R,U,\upsilon},\tilde\lambda_{R,U,\upsilon}}{s}
		\MM{N,K_0,\lambda_{R,U,\upsilon},\tilde\lambda_{R,U,\upsilon}}.
	\end{equation}
	For all $M,N\in\mathbb N_0$ satisfying
	$1\leq M\leq \min\{M_*,K^*+1\}$ and $M+N\leq K^\sharp-1$, we have
	\begin{align}
		\nrm{\partial_t^MV}_{N+\alpha}
		&\lesssim
		\mC_U\bar\varepsilon^\alpha
		\Lk{N,(K_U-M)_+,\lambda_{R,U,\upsilon},\tilde\lambda_{R,U,\upsilon}}{s}\La{N+1,(K_U-M)_+,\lambda_{R,U,\upsilon},\tilde\lambda_{R,U,\upsilon}}\notag\\
		\label{eq:estimate-time-derivatives-V-pullback}
		&\quad\cdot
		\MM{N,(K_U-M)_+,\lambda_{R,U,\upsilon},\tilde\lambda_{R,U,\upsilon}}
		\MM{M,N_t,\mu_{R,U,\upsilon},\tilde\mu_{R,U,\upsilon}},
		\\
		\nrm{\partial_t^MV}_{N}
		&\lesssim
		\mC_U\bar\varepsilon^\alpha
		\Lk{N,(K_U-M)_+,\lambda_{R,U,\upsilon},\tilde\lambda_{R,U,\upsilon}}{s}\notag\\
		&\quad\cdot\MM{N,(K_U-M)_+,\lambda_{R,U,\upsilon},\tilde\lambda_{R,U,\upsilon}}
		\MM{M,N_t,\mu_{R,U,\upsilon},\tilde\mu_{R,U,\upsilon}}.
		\label{eq:estimate-time-derivatives-V-pullback-CN}
	\end{align}
	In particular, for every $1\leq M\leq \min\{M_*,K^*+1\}$,
	\begin{equation}
		\label{eq:estimate-time-derivatives-V-pullback-C0-small}
		\nrm{\partial_t^MV}_{0}
		\lesssim
		\mC_U\bar\varepsilon^\alpha
		\Lk{0,(K_U-M)_+,\lambda_{R,U,\upsilon},\tilde\lambda_{R,U,\upsilon}}{s}
		\MM{M,N_t,\mu_{R,U,\upsilon},\tilde\mu_{R,U,\upsilon}}.
	\end{equation}
	
	For the solution operator, when $M=0$ and $0\leq N\leq K^\sharp-2$,
	\begin{align}
		\nrm{R[U]}_{N+\alpha}
		&\lesssim
		\mC_U\bar\varepsilon^{1+\alpha}
		\Lk{N,K_0,\lambda_{R,U,\upsilon},\tilde\lambda_{R,U,\upsilon}}{s}\La{N+1,K_0,\lambda_{R,U,\upsilon},\tilde\lambda_{R,U,\upsilon}}
		\MM{N,K_0,\lambda_{R,U,\upsilon},\tilde\lambda_{R,U,\upsilon}}.
		\label{eq:Eulerian-R-Calpha-material-zero-MM}
	\end{align}
	For every $0\leq N\leq K^\sharp-1$, the corresponding $C^N$ estimate is
	\begin{equation}
		\label{eq:Eulerian-R-CN-material-zero-MM}
		\nrm{R[U]}_{N}
		\lesssim
		\mC_U\bar\varepsilon^{1+\alpha}
		\Lk{N,K_0,\lambda_{R,U,\upsilon},\tilde\lambda_{R,U,\upsilon}}{s}
		\MM{N,K_0,\lambda_{R,U,\upsilon},\tilde\lambda_{R,U,\upsilon}}.
	\end{equation}
	For all $M,N\in\mathbb N_0$ satisfying
	$1\leq M\leq \min\{M_*,K^*+1\}$ and $M+N\leq K^\sharp-1$, we have
	\begin{align}
		\nrm{D_{t,\upsilon}^{M}R[U]}_{N+\alpha}
		&\lesssim
		\mC_U\bar\varepsilon^{1+\alpha}
		\Lk{N,(K_U-M)_+,\lambda_{R,U,\upsilon},\tilde\lambda_{R,U,\upsilon}}{s}\La{N+1,(K_U-M)_+,\lambda_{R,U,\upsilon},\tilde\lambda_{R,U,\upsilon}}\notag\\
		\label{eq:Eulerian-R-Calpha-material-MM}
		&\quad\cdot
		\MM{N,(K_U-M)_+,\lambda_{R,U,\upsilon},\tilde\lambda_{R,U,\upsilon}}
		\MM{M,N_t,\mu_{R,U,\upsilon},\tilde\mu_{R,U,\upsilon}},
		\\
		\nrm{D_{t,\upsilon}^{M}R[U]}_{N}
		&\lesssim
		\mC_U\bar\varepsilon^{1+\alpha}
		\Lk{N,(K_U-M)_+,\lambda_{R,U,\upsilon},\tilde\lambda_{R,U,\upsilon}}{s}\nonumber\\
		&\quad\cdot\MM{N,(K_U-M)_+,\lambda_{R,U,\upsilon},\tilde\lambda_{R,U,\upsilon}}
		\MM{M,N_t,\mu_{R,U,\upsilon},\tilde\mu_{R,U,\upsilon}}.
		\label{eq:Eulerian-R-CN-material-MM}
	\end{align}

	If $U$ instead satisfies \eqref{eq:inv4divEq:integer-source-assumption} for $M,N\in\mathbb N_0$ with
	$0\leq M\leq M_*$ and $M+N\leq K^\sharp-1$, the following $C^N$ estimates hold.
	For $M=0$ and $0\leq N\leq K^\sharp-1$, we have
	\begin{equation}
		\label{eq:estimate-V-pullback-integer-source-zero}
		\nrm{V}_{N}
		\lesssim\mC_U\MM{N,K_0,\lambda_{R,U,\upsilon},\tilde\lambda_{R,U,\upsilon}}.
	\end{equation}
	For $M,N\in\mathbb N_0$ with $1\leq M\leq\min\{M_*,K^*+1\}$ and $M+N\leq K^\sharp-1$,
	\begin{equation}
		\label{eq:estimate-V-pullback-integer-source}
		\nrm{\partial_t^MV}_{N}
		\lesssim\mC_U\MM{N,(K_U-M)_+,\lambda_{R,U,\upsilon},\tilde\lambda_{R,U,\upsilon}}
		\MM{M,N_t,\mu_{R,U,\upsilon},\tilde\mu_{R,U,\upsilon}}.
	\end{equation}
	For $R[U]$, when $M=0$ and $0\leq N\leq K^\sharp-1$,
	\begin{equation}
		\label{eq:Eulerian-R-CN-integer-source-zero}
		\nrm{R[U]}_{N}
		\lesssim\mC_U\bar\varepsilon\MM{N,K_0,\lambda_{R,U,\upsilon},\tilde\lambda_{R,U,\upsilon}}.
	\end{equation}
	For $M,N\in\mathbb N_0$ with $1\leq M\leq\min\{M_*,K^*+1\}$ and $M+N\leq K^\sharp-1$,
	\begin{equation}
		\label{eq:Eulerian-R-CN-integer-source}
		\nrm{D_{t,\upsilon}^{M}R[U]}_{N}
		\lesssim\mC_U\bar\varepsilon\MM{N,(K_U-M)_+,\lambda_{R,U,\upsilon},\tilde\lambda_{R,U,\upsilon}}
		\MM{M,N_t,\mu_{R,U,\upsilon},\tilde\mu_{R,U,\upsilon}}.
	\end{equation}

	The implicit constants may depend on $d,K^*,K^\sharp,\alpha$.
	They are independent of $t_0,x_0,\bar\tau,\bar\varepsilon$ and of the amplitudes and frequency scales.
\end{lem}

\begin{proof}
	By \eqref{eq:inv4divEq:supp4U} and \eqref{eq:Rjl0},
	$V(t,\cdot)$ is supported in the fixed ball $B(x_0,\bar\varepsilon)$
	for every $t$.
	The same holds for its time derivatives, proving
	\eqref{eq:supp-time-derivatives-V-pullback}.
	
	By \eqref{eq:Rjl0},
	$
	V_i(t,y)=\JX(t,y)U_i(t,X_\upsilon(t,y;t_0)).
	$
	The flow equation reads
	$
	\partial_tX_\upsilon(t,y;t_0)
	=
	\upsilon(t,X_\upsilon(t,y;t_0)),
	$
	and hence
	\begin{equation}
		\label{eq:dt-pullback-U-material}
		\partial_t^{M_2}
		\left[
		U_i(t,X_\upsilon(t,y;t_0))
		\right]
		=
		(D_{t,\upsilon}^{M_2}U_i)
		(t,X_\upsilon(t,y;t_0))
	\end{equation}
	for every $M_2\in\mathbb N_0$. Consequently,
	\begin{equation}
		\label{eq:dtM-V-expansion}
		\partial_t^MV_i(t,y)
		=
		\sum_{M_1+M_2=M}
		\binom{M}{M_1}
		(\partial_t^{M_1}\JX)(t,y)
		(D_{t,\upsilon}^{M_2}U_i)
		(t,X_\upsilon(t,y;t_0)).
	\end{equation}
	
	The assumptions \eqref{eq:inv4divEq:est4nbv-material-CN-MM}--%
	\eqref{eq:inv4divEq:est4nbv-material-Calpha-MM} also control
	$D_{t,\upsilon}^M\Div\upsilon=\operatorname{tr}(D_{t,\upsilon}^M\nabla\upsilon)$ through $M=K^*$.
	The proof of Lemma~\ref{lem:flow-estimates-upsilon} therefore applies locally in $\R^d$.
	Indeed, the trajectories starting in
	$B(x_0,\bar\varepsilon)$ remain in the enlarged Eulerian cylinder
	by the local gradient estimate and the cylinder inclusion in
	Subsection~\ref{subsec:adapted-cylinders}.
	Using \eqref{eq:flow-jacobian-CN} and
	\eqref{eq:flow-jacobian-Calpha},
	\eqref{eq:tau-local-smallness}, and the scale definitions
	\eqref{eq:definition-lambda-RUV}--\eqref{eq:definition-mu-RUV},
	we obtain, for $0\leq N\leq K^\sharp-2$,
	\begin{align*}
		\nrm{\JX}_{N+\alpha}
		&\lesssim
		\La{N,K^*-1,\lambda_{R,U,\upsilon},\tilde\lambda_{R,U,\upsilon}}
		\MM{N,K^*-1,\lambda_{R,U,\upsilon},\tilde\lambda_{R,U,\upsilon}}.
	\end{align*}
	For $0\leq N\leq K^\sharp-1$, we also have
	$$
	\nrm{\JX}_{N}
	\lesssim
	\MM{N,K^*-1,\lambda_{R,U,\upsilon},\tilde\lambda_{R,U,\upsilon}}.
	$$
	
	For $M_1,N\in\mathbb N_0$ with
	$1\leq M_1\leq\min\{M_*,K^*+1\}$ and
	$M_1+N\leq K^\sharp-1$, we absorb the extra spatial factor using
	$\mC_\upsilon\lambda_\upsilon\leq\mu_{R,U,\upsilon}$ for $M_1\leq K^*$.
	At $M_1=K^*+1$, we use $\mC_\upsilon\tilde\lambda_\upsilon\leq\tilde\mu_{R,U,\upsilon}$
	and $M_1=K^*+1>N_t$.
	The argument for \eqref{eq:flow-Jacobian-power-material} and
	\eqref{eq:flow-Jacobian-power-material-Calpha} then gives
	\begin{align*}
		\nrm{\partial_t^{M_1}\JX}_{N+\alpha}
		&\lesssim
		\La{N,(K^*-M_1)_+,\lambda_{R,U,\upsilon},\tilde\lambda_{R,U,\upsilon}}
		\MM{N,(K^*-M_1)_+,\lambda_{R,U,\upsilon},\tilde\lambda_{R,U,\upsilon}}
		\MM{M_1,N_t,\mu_{R,U,\upsilon},\tilde\mu_{R,U,\upsilon}},
		\\
		\nrm{\partial_t^{M_1}\JX}_{N}
		&\lesssim
		\MM{N,(K^*-M_1)_+,\lambda_{R,U,\upsilon},\tilde\lambda_{R,U,\upsilon}}
		\MM{M_1,N_t,\mu_{R,U,\upsilon},\tilde\mu_{R,U,\upsilon}}.
	\end{align*}
	The backward flow equations give the same estimates for $\JPhi$
	and its material derivatives.
	
	For $N,K\in\mathbb N_0$, the combined scales satisfy
	$
	\bar\varepsilon^\alpha
	\La{N,K,\lambda_{R,U,\upsilon},\tilde\lambda_{R,U,\upsilon}}
	\geq
	(\bar\varepsilon\lambda_{R,U,\upsilon})^\alpha
	\geq1.
	$
	Material and spatial differentiation preserve the support inclusion
	$
	\supp D_{t,\upsilon}^{M_2}U
	\subseteq\supp U
	\subseteq\hat\Gamma_\upsilon(\bar\tau,\bar\varepsilon;t_0,x_0).
	$
	The flow estimate gives a diameter $\lesssim\bar\varepsilon$ for
	$X_\upsilon(t,B(x_0,\bar\varepsilon);t_0)$.
	Since $C_{\mathrm{loc}}\geq4$, the small-support argument in the proof of
	Lemma~\ref{lem:Rstar-spatial-holder} applies within the enlarged Eulerian cylinder and gives
	\begin{equation}\label{eq:inv4divEq-small-support-material-estimate}
	\nrm{D_{t,\upsilon}^{M_2}U_i}_{m}
	\lesssim
	\bar\varepsilon^\alpha
	\nrm{D_{t,\upsilon}^{M_2}U_i}_{m+\alpha}.
	\end{equation}
	
	We estimate the compositions in \eqref{eq:dtM-V-expansion} using
	the Fa\`a di Bruno formula \eqref{eq:faa-di-bruno} for the $C^N$ norms with $N\geq1$ and
	the argument of Lemma~\ref{lem:holder-composition-general} for the H{\"o}lder seminorms. Here and below, the local form of \eqref{eq:flow-gradient-close-id} gives uniform Lipschitz constants for $X_\upsilon$ on the initial ball and for $\Phi_\upsilon$ on its image. The compositions vanish near the corresponding boundaries. For differences across a boundary, insert the first boundary point on the joining segment. The zero extensions satisfy the same H{\"o}lder estimates.
	For $M_2=0$ and $N\geq1$, write $N_1+\cdots+N_m=N$ with $N_\nu\geq1$.
	The spatial scales satisfy
	$
		\MM{m,K_U,\lambda_U,\tilde\lambda_U}
		\prod_{\nu=1}^{m}
		\MM{N_\nu-1,K^*-1,\lambda_\upsilon,\tilde\lambda_\upsilon}\leq
		\MM{N,K_0,\lambda_{R,U,\upsilon},\tilde\lambda_{R,U,\upsilon}}.
	$
	For $M_2\geq1$, the inequality $(K_U-M_2)_+\leq K^*-1$ allows us to
	replace the second argument $K^*-1$ in the weights for derivatives of
	$\nabla X_\upsilon$ by $(K_U-M_2)_+$, as in the source estimate.
	Assumption~\eqref{eq:inv4divEq:hyp:3-material-Calpha-MM} and \eqref{eq:inv4divEq-small-support-material-estimate} control
	the $C^{m+\alpha}$ and $C^m$ norms of $D_{t,\upsilon}^{M_2}U_i$, respectively.
	
	For $M_2=0$ and $0\leq N\leq K^\sharp-2$, we have
	\begin{align*}
		\nrm{
			U_i(t,X_\upsilon(t,\cdot;t_0))
		}_{N+\alpha}
		\lesssim
		\mC_U
		\bar\varepsilon^\alpha\Lk{N,K_0,\lambda_{R,U,\upsilon},\tilde\lambda_{R,U,\upsilon}}{s+1}
		\MM{N,K_0,\lambda_{R,U,\upsilon},\tilde\lambda_{R,U,\upsilon}}.
	\end{align*}
	For $1\leq M_2\leq \min\{M_*,K^*+1\}$ and $M_2+N\leq K^\sharp-1$, we obtain
	\begin{align*}
		\nrm{
			(D_{t,\upsilon}^{M_2}U_i)
			(t,X_\upsilon(t,\cdot;t_0))
		}_{N+\alpha}
		\lesssim&
		\mC_U
		\bar\varepsilon^\alpha\Lk{N,(K_U-M_2)_+,\lambda_{R,U,\upsilon},\tilde\lambda_{R,U,\upsilon}}{s+1}\\
		&\cdot
		\MM{N,(K_U-M_2)_+,\lambda_{R,U,\upsilon},\tilde\lambda_{R,U,\upsilon}}
		\MM{M_2,N_t,\mu_{R,U,\upsilon},\tilde\mu_{R,U,\upsilon}}.
	\end{align*}
	The same argument gives the $C^N$ estimate, for $M_2=0$ and $0\leq N\leq K^\sharp-1$,
	\begin{align}
		\nrm{
			U_i(t,X_\upsilon(t,\cdot;t_0))
		}_{N}
		&\lesssim
		\mC_U\bar\varepsilon^\alpha
		\Lk{N,K_0,\lambda_{R,U,\upsilon},\tilde\lambda_{R,U,\upsilon}}{s}
		\MM{N,K_0,\lambda_{R,U,\upsilon},\tilde\lambda_{R,U,\upsilon}},
		\label{eq:local-inverse-source-composition-CN-zero}
	\end{align}
	while, for $1\leq M_2\leq \min\{M_*,K^*+1\}$ and $M_2+N\leq K^\sharp-1$,
	\begin{align}
		\nrm{
			(D_{t,\upsilon}^{M_2}U_i)
			(t,X_\upsilon(t,\cdot;t_0))
		}_{N}
		&\lesssim
		\mC_U\bar\varepsilon^\alpha
		\Lk{N,(K_U-M_2)_+,\lambda_{R,U,\upsilon},\tilde\lambda_{R,U,\upsilon}}{s}\notag\\
		&\quad\cdot
		\MM{N,(K_U-M_2)_+,\lambda_{R,U,\upsilon},\tilde\lambda_{R,U,\upsilon}}
		\MM{M_2,N_t,\mu_{R,U,\upsilon},\tilde\mu_{R,U,\upsilon}}.
		\label{eq:local-inverse-source-composition-CN-material}
	\end{align}
	
	In \eqref{est:holder-composition-general}, the $C^{m+\alpha}$ source estimate
	gives the first term, using
	$\bar\varepsilon^\alpha\La{N,K,\lambda_{R,U,\upsilon},\tilde\lambda_{R,U,\upsilon}}\geq1$
	with $K=K_0$ when $M_2=0$ and $K=(K_U-M_2)_+$ otherwise.
	For the second term, \eqref{eq:inv4divEq-small-support-material-estimate} gives the factor $\bar\varepsilon^\alpha$,
	and the derivative of $X_\upsilon$ contributes the additional H{\"o}lder factor.
	For $N\geq1$, the $C^N$ estimates follow from \eqref{eq:inv4divEq-small-support-material-estimate} and \eqref{eq:faa-di-bruno}.
	At $N=0$, use \eqref{est:holder-composition-zero} and the $C^0$ composition estimate.
	
	We now combine the composition estimates with the estimates for $\partial_t^{M_1}\JX$
	in \eqref{eq:dtM-V-expansion}. For fixed derivative order and $1\leq\lambda\leq\tilde\lambda$, the weights
	$\mathcal M$ and $\Xi_\alpha$ are nonincreasing in their threshold parameter $K$.
	For $M_1+M_2=M\geq1$, we have
	$
	\min\{K_0,(K^*-M_1)_+,(K_U-M_2)_+\}\geq(K_U-M)_+.
	$
	The undifferentiated factors $\JX$ and $U\circ X_\upsilon$ have thresholds
	$K^*-1\geq K_0$ and $K_0$, respectively. Thus all spatial weights in
	\eqref{eq:dtM-V-expansion} can be bounded using the common threshold $(K_U-M)_+$.
	Applying the case $\alpha=0$ of \eqref{eq:holder-product-tame} and \eqref{eq:holder-product-finite-Cnalpha} to \eqref{eq:dtM-V-expansion}, we obtain
	\begin{align}
		\|\partial_t^MV\|_{N+\alpha}
		&\lesssim\sum_{M_1+M_2=M}\sum_{N_1+N_2=N}
		\Big(\|\partial_t^{M_1}\JX\|_{N_1}\|(D_{t,\upsilon}^{M_2}U)\circ X_\upsilon\|_{N_2+\alpha}+\|\partial_t^{M_1}\JX\|_{N_1+\alpha}\|(D_{t,\upsilon}^{M_2}U)\circ X_\upsilon\|_{N_2}\Big),\notag\\
		\|\partial_t^MV\|_N
		&\lesssim\sum_{M_1+M_2=M}\sum_{N_1+N_2=N}
		\|\partial_t^{M_1}\JX\|_{N_1}
		\|(D_{t,\upsilon}^{M_2}U)\circ X_\upsilon\|_{N_2}.
		\label{eq:local-inverse-pullback-CN-product}
	\end{align}
	For $N_1,N_2,K\in\mathbb N_0$ with $N_1+N_2\leq N$ and
	$1\leq\lambda\leq\tilde\lambda$, the H{\"o}lder factors satisfy
	$$
	\Lk{N_2,K,\lambda,\tilde\lambda}{s}\La{N_1+1,K,\lambda,\tilde\lambda}
	\leq\Lk{N,K,\lambda,\tilde\lambda}{s}\La{N+1,K,\lambda,\tilde\lambda}.
	$$
	Substituting the preceding $C^N$ and H{\"o}lder estimates
	gives, for $M=0$,
	\begin{align*}
		\nrm{V}_{N+\alpha}
		&\lesssim
		\mC_U\bar\varepsilon^\alpha
		\Lk{N,K_0,\lambda_{R,U,\upsilon},\tilde\lambda_{R,U,\upsilon}}{s}\La{N+1,K_0,\lambda_{R,U,\upsilon},\tilde\lambda_{R,U,\upsilon}}
		\MM{N,K_0,\lambda_{R,U,\upsilon},\tilde\lambda_{R,U,\upsilon}}.
	\end{align*}
	For $M\geq1$, the corresponding estimate is
	\begin{align*}
		\nrm{\partial_t^MV}_{N+\alpha}
		&\lesssim
		\mC_U\bar\varepsilon^\alpha
		\Lk{N,(K_U-M)_+,\lambda_{R,U,\upsilon},\tilde\lambda_{R,U,\upsilon}}{s}\La{N+1,(K_U-M)_+,\lambda_{R,U,\upsilon},\tilde\lambda_{R,U,\upsilon}}\\
		&\quad\cdot
		\MM{N,(K_U-M)_+,\lambda_{R,U,\upsilon},\tilde\lambda_{R,U,\upsilon}}
		\MM{M,N_t,\mu_{R,U,\upsilon},\tilde\mu_{R,U,\upsilon}}.
	\end{align*}
	These are \eqref{eq:estimate-V-pullback-material-zero} and
	\eqref{eq:estimate-time-derivatives-V-pullback}.
	Substituting the $C^N$ estimates into \eqref{eq:local-inverse-pullback-CN-product} gives, for $M=0$,
	$$
	\nrm{V}_{N}
	\lesssim
	\mC_U\bar\varepsilon^\alpha
	\Lk{N,K_0,\lambda_{R,U,\upsilon},\tilde\lambda_{R,U,\upsilon}}{s}
	\MM{N,K_0,\lambda_{R,U,\upsilon},\tilde\lambda_{R,U,\upsilon}},
	$$
	and, for $M\geq1$,
	$$
	\nrm{\partial_t^MV}_{N}
	\lesssim
	\mC_U\bar\varepsilon^\alpha
	\Lk{N,(K_U-M)_+,\lambda_{R,U,\upsilon},\tilde\lambda_{R,U,\upsilon}}{s}
	\MM{N,(K_U-M)_+,\lambda_{R,U,\upsilon},\tilde\lambda_{R,U,\upsilon}}
	\MM{M,N_t,\mu_{R,U,\upsilon},\tilde\mu_{R,U,\upsilon}}.
	$$
	We obtain \eqref{eq:estimate-V-pullback-material-zero-CN} and
	\eqref{eq:estimate-time-derivatives-V-pullback-CN}.
	Taking $N=0$ in \eqref{eq:estimate-time-derivatives-V-pullback-CN} gives
	\eqref{eq:estimate-time-derivatives-V-pullback-C0-small}.
	
	For $R[U]$, we first estimate $R^*[\partial_t^kV]$.
	For every $k\in\mathbb N_0$, the fixed support
	\eqref{eq:supp-time-derivatives-V-pullback} and differentiation of the
	zero-mean identity for $V$ verify the support and
	cancellation hypotheses of Lemma~\ref{lem:Rstar-spatial-holder} for
	$\partial_t^kV$. Since $0<\bar\varepsilon\leq1$ and
	$\lambda_{R,U,\upsilon}\geq\bar\varepsilon^{-1}$, the factors $\bar\varepsilon^{-(N-r)}$ in
	\eqref{eq:Rstar-CN-basic} and \eqref{eq:Rstar-Calpha-basic}
	are absorbed by the spatial weights. Applying these estimates to
	\eqref{eq:estimate-V-pullback-material-zero}--\eqref{eq:estimate-time-derivatives-V-pullback-CN}
	gives the same right-hand sides for $R^*[\partial_t^kV]$,
	with an additional factor $\bar\varepsilon$ and the same derivative ranges.
	
	For $R^*[\partial_t^kV]\circ\Phi_\upsilon$, we use \eqref{eq:faa-di-bruno} and
	the argument of Lemma~\ref{lem:holder-composition-general}.
	For $1\leq m\leq N$, the two terms in \eqref{est:holder-composition-general} use
	the $C^{m+\alpha}$ and $C^m$ estimates for $R^*[\partial_t^kV]$, respectively.
	In either term, the additional H{\"o}lder factor is bounded by the factor at order $N+1$, with threshold $K_0$ for $k=0$ and $(K_U-k)_+$ for $k\geq1$.
	Since $K_0\leq K^*-1$ and $(K_U-k)_+\leq K^*-1$ for $k\geq1$,
	the spatial scale products are controlled as for composition with $X_\upsilon$ above.
	For $N=0$, we use \eqref{est:holder-composition-zero} and the $C^0$ composition estimate.
	For $k=0$, we obtain the following H{\"o}lder and $C^N$ estimates
	for $0\leq N\leq K^\sharp-2$ and $0\leq N\leq K^\sharp-1$, respectively,
	\begin{align}
		\nrm{
			R^*[V(t,\cdot)]\circ\Phi_\upsilon
		}_{N+\alpha}
		&\lesssim
		\mC_U\bar\varepsilon^{1+\alpha}
		\Lk{N,K_0,\lambda_{R,U,\upsilon},\tilde\lambda_{R,U,\upsilon}}{s}\La{N+1,K_0,\lambda_{R,U,\upsilon},\tilde\lambda_{R,U,\upsilon}}
		\MM{N,K_0,\lambda_{R,U,\upsilon},\tilde\lambda_{R,U,\upsilon}},\notag\\
		\nrm{
			R^*[V(t,\cdot)]\circ\Phi_\upsilon
		}_{N}
		&\lesssim
		\mC_U\bar\varepsilon^{1+\alpha}
		\Lk{N,K_0,\lambda_{R,U,\upsilon},\tilde\lambda_{R,U,\upsilon}}{s}
		\MM{N,K_0,\lambda_{R,U,\upsilon},\tilde\lambda_{R,U,\upsilon}}.
		\label{eq:local-inverse-composition-CN-zero}
	\end{align}
	For $1\leq k\leq \min\{M_*,K^*+1\}$ and
	$k+N\leq K^\sharp-1$, applying the composition estimate to
	$R^*[\partial_t^kV]$ gives
	\begin{align*}
		\nrm{
			R^*[\partial_t^kV(t,\cdot)]\circ\Phi_\upsilon
		}_{N+\alpha}
		\lesssim&
		\mC_U\bar\varepsilon^{1+\alpha}
		\Lk{N,(K_U-k)_+,\lambda_{R,U,\upsilon},\tilde\lambda_{R,U,\upsilon}}{s}\La{N+1,(K_U-k)_+,\lambda_{R,U,\upsilon},\tilde\lambda_{R,U,\upsilon}}\\
		&\quad\cdot
		\MM{N,(K_U-k)_+,\lambda_{R,U,\upsilon},\tilde\lambda_{R,U,\upsilon}}
		\MM{k,N_t,\mu_{R,U,\upsilon},\tilde\mu_{R,U,\upsilon}}.
	\end{align*}
	The corresponding $C^N$ estimate is
	\begin{equation}
		\label{eq:local-inverse-composition-CN-material}
	\begin{aligned}
		&\nrm{
			R^*[\partial_t^kV(t,\cdot)]\circ\Phi_\upsilon
		}_{N}\\
		&\quad\lesssim
		\mC_U\bar\varepsilon^{1+\alpha}
		\Lk{N,(K_U-k)_+,\lambda_{R,U,\upsilon},\tilde\lambda_{R,U,\upsilon}}{s}
		\MM{N,(K_U-k)_+,\lambda_{R,U,\upsilon},\tilde\lambda_{R,U,\upsilon}}
		\MM{k,N_t,\mu_{R,U,\upsilon},\tilde\mu_{R,U,\upsilon}}.
	\end{aligned}
	\end{equation}
	
	It remains to estimate the factor $P$ and its material derivatives in
	\eqref{eq:higher-material-R}. With $P=\JPhi(\nabla\Phi_\upsilon)^{-1}$,
	\eqref{eq:flow-gradient-CN}--\eqref{eq:flow-jacobian-Calpha}
	from Lemma~\ref{lem:flow-estimates-upsilon} give, with the combined scales,
	for $0\leq N\leq K^\sharp-2$ in H{\"o}lder norms and
	$0\leq N\leq K^\sharp-1$ in $C^N$ norms,
	\begin{align*}
		\nrm{P}_{N+\alpha}
		\lesssim
		\La{N+1,K^*-1,\lambda_{R,U,\upsilon},\tilde\lambda_{R,U,\upsilon}}
		\MM{N,K^*-1,\lambda_{R,U,\upsilon},\tilde\lambda_{R,U,\upsilon}},
		\quad
		\nrm{P}_{N}
		\lesssim
		\MM{N,K^*-1,\lambda_{R,U,\upsilon},\tilde\lambda_{R,U,\upsilon}}.
	\end{align*}
	For the material derivatives, we use
	\eqref{eq:flow-gradient-material}--\eqref{eq:flow-Jacobian-power-material-Calpha},
	including the extension to order $K^*+1$ stated in
	Lemma~\ref{lem:flow-estimates-upsilon}.
	We absorb the extra spatial factor into the time scales as in the estimate
	for $\partial_t^{M_1}\JX$ above to obtain, for
	$1\leq M\leq\min\{M_*,K^*+1\}$,
	$M+N\leq K^\sharp-1$, and $0\leq k<M$,
	\begin{align*}
		\nrm{D_{t,\upsilon}^{M-k}P}_{N+\alpha}
		&\lesssim
		\La{N+1,(K^*-M+k)_+,\lambda_{R,U,\upsilon},\tilde\lambda_{R,U,\upsilon}}\MM{N,(K^*-M+k)_+,\lambda_{R,U,\upsilon},\tilde\lambda_{R,U,\upsilon}}\\
		&\quad\cdot
		\MM{M-k,N_t,\mu_{R,U,\upsilon},\tilde\mu_{R,U,\upsilon}},
		\\
		\nrm{D_{t,\upsilon}^{M-k}P}_{N}
		&\lesssim
		\MM{N,(K^*-M+k)_+,\lambda_{R,U,\upsilon},\tilde\lambda_{R,U,\upsilon}}
		\MM{M-k,N_t,\mu_{R,U,\upsilon},\tilde\mu_{R,U,\upsilon}}.
	\end{align*}
	At $M=0$, \eqref{eq:Rjl2} and the estimates for $P$ and
	$R^*[V]\circ\Phi_\upsilon$ give
	\begin{align*}
		\nrm{R[U]}_{N+\alpha}
		&\lesssim
		\mC_U\bar\varepsilon^{1+\alpha}
		\Lk{N,K_0,\lambda_{R,U,\upsilon},\tilde\lambda_{R,U,\upsilon}}{s}\La{N+1,K_0,\lambda_{R,U,\upsilon},\tilde\lambda_{R,U,\upsilon}}
		\MM{N,K_0,\lambda_{R,U,\upsilon},\tilde\lambda_{R,U,\upsilon}},\\
		\nrm{R[U]}_{N}
		&\lesssim
		\mC_U\bar\varepsilon^{1+\alpha}
		\Lk{N,K_0,\lambda_{R,U,\upsilon},\tilde\lambda_{R,U,\upsilon}}{s}
		\MM{N,K_0,\lambda_{R,U,\upsilon},\tilde\lambda_{R,U,\upsilon}}.
	\end{align*}
	Estimates~\eqref{eq:Eulerian-R-Calpha-material-zero-MM}--%
	\eqref{eq:Eulerian-R-CN-material-zero-MM} follow.
	
	Fix
	$1\leq M\leq\min\{M_*,K^*+1\}$ and $M+N\leq K^\sharp-1$.
	For $0\leq k\leq M$, we have
	$
	\min\{K_0,(K^*-M+k)_+,(K_U-k)_+\}\geq(K_U-M)_+.
	$
	The undifferentiated factors at $k=0$ and $k=M$ have thresholds
	$K_0$ and $K^*-1\geq K_0$, respectively.
	By the monotonicity of $\mathcal M$ and $\Xi_\alpha$, we may lower all spatial thresholds to $(K_U-M)_+$.
	Applying \eqref{eq:holder-product-finite-Cnalpha} to each term in
	\eqref{eq:higher-material-R} gives
	\begin{align*}
		&
		\Nrm{
			(D_{t,\upsilon}^{M-k}P_{jl})
			\left(
			R^*_{il}[\partial_t^kV(t,\cdot)]\circ\Phi_\upsilon
			\right)
		}_{N+\alpha}
		\\
		\lesssim&
		\sum_{N_1+N_2=N}\Bigl(
		\nrm{D_{t,\upsilon}^{M-k}P}_{N_1}
		\nrm{R^*[\partial_t^kV]\circ\Phi_\upsilon}_{N_2+\alpha}+
		\nrm{D_{t,\upsilon}^{M-k}P}_{N_1+\alpha}
		\nrm{R^*[\partial_t^kV]\circ\Phi_\upsilon}_{N_2}\Bigr)\\
		\lesssim&
		\mC_U\bar\varepsilon^{1+\alpha}
		\Lk{N,(K_U-M)_+,\lambda_{R,U,\upsilon},\tilde\lambda_{R,U,\upsilon}}{s}\La{N+1,(K_U-M)_+,\lambda_{R,U,\upsilon},\tilde\lambda_{R,U,\upsilon}}\\
		&\quad\cdot
		\MM{N,(K_U-M)_+,\lambda_{R,U,\upsilon},\tilde\lambda_{R,U,\upsilon}}
		\MM{M-k,N_t,\mu_{R,U,\upsilon},\tilde\mu_{R,U,\upsilon}}
		\MM{k,N_t,\mu_{R,U,\upsilon},\tilde\mu_{R,U,\upsilon}}.
	\end{align*}
	Using
	$
	\MM{M-k,N_t,\mu_{R,U,\upsilon},\tilde\mu_{R,U,\upsilon}}
	\MM{k,N_t,\mu_{R,U,\upsilon},\tilde\mu_{R,U,\upsilon}}
	\leq
	\MM{M,N_t,\mu_{R,U,\upsilon},\tilde\mu_{R,U,\upsilon}},
	$
	and summing over $0\leq k\leq M$, we obtain
	\eqref{eq:Eulerian-R-Calpha-material-MM}.
	Estimate~\eqref{eq:holder-product-tame} with $\alpha=0$ gives
	\begin{align}
		\Nrm{
			(D_{t,\upsilon}^{M-k}P_{jl})
			\left(
			R^*_{il}[\partial_t^kV(t,\cdot)]\circ\Phi_\upsilon
			\right)
		}_{N}
		\lesssim&
		\mC_U\bar\varepsilon^{1+\alpha}
		\Lk{N,(K_U-M)_+,\lambda_{R,U,\upsilon},\tilde\lambda_{R,U,\upsilon}}{s}
		\MM{N,(K_U-M)_+,\lambda_{R,U,\upsilon},\tilde\lambda_{R,U,\upsilon}}\notag\\
		&\quad\cdot
		\MM{M-k,N_t,\mu_{R,U,\upsilon},\tilde\mu_{R,U,\upsilon}}
		\MM{k,N_t,\mu_{R,U,\upsilon},\tilde\mu_{R,U,\upsilon}}.
		\label{eq:local-inverse-CN-product}
	\end{align}
	Equations~\eqref{eq:higher-material-R} and \eqref{eq:local-inverse-CN-product} prove
	\eqref{eq:Eulerian-R-CN-material-MM}.

	Finally, suppose that \eqref{eq:inv4divEq:hyp:3-material-Calpha-MM} is replaced by
	\eqref{eq:inv4divEq:integer-source-assumption}.
	For compositions with $X_\upsilon$ and $\Phi_\upsilon$, use \eqref{eq:faa-di-bruno} when $N\geq1$
	and the $C^0$ composition estimate when $N=0$.
	Apply \eqref{eq:inv4divEq:integer-source-assumption} directly in the arguments for
	\eqref{eq:local-inverse-source-composition-CN-zero} and \eqref{eq:local-inverse-source-composition-CN-material}.
	For $M_2=0$ and $0\leq N\leq K^\sharp-1$, we have
	$$
	\nrm{U_i(t,X_\upsilon(t,\cdot;t_0))}_N
	\lesssim\mC_U\MM{N,K_0,\lambda_{R,U,\upsilon},\tilde\lambda_{R,U,\upsilon}}.
	$$
	For $1\leq M_2\leq\min\{M_*,K^*+1\}$ and $M_2+N\leq K^\sharp-1$, we obtain
	$$
	\begin{aligned}
		\nrm{(D_{t,\upsilon}^{M_2}U_i)(t,X_\upsilon(t,\cdot;t_0))}_N
		\lesssim\mC_U\MM{N,(K_U-M_2)_+,\lambda_{R,U,\upsilon},\tilde\lambda_{R,U,\upsilon}}\MM{M_2,N_t,\mu_{R,U,\upsilon},\tilde\mu_{R,U,\upsilon}}.
	\end{aligned}
	$$
	Substituting these estimates into \eqref{eq:local-inverse-pullback-CN-product} proves
	\eqref{eq:estimate-V-pullback-integer-source-zero} and \eqref{eq:estimate-V-pullback-integer-source}.

	The support and zero-mean properties of $\partial_t^kV$ also remain valid.
	For $k,N\in\mathbb N_0$ with $0\leq k\leq\min\{M_*,K^*+1\}$ and $k+N\leq K^\sharp-1$,
	\eqref{eq:Rstar-CN-basic} gives
	$$
	\nrm{R^*[\partial_t^kV]}_N
	\lesssim\bar\varepsilon\sum_{r=0}^N\bar\varepsilon^{-(N-r)}\nrm{\partial_t^kV}_r.
	$$
	Using \eqref{eq:estimate-V-pullback-integer-source-zero}--\eqref{eq:estimate-V-pullback-integer-source} and
	$\lambda_{R,U,\upsilon}\geq\bar\varepsilon^{-1}$,
	repeat the argument for the $C^N$ estimates in \eqref{eq:local-inverse-composition-CN-zero} and \eqref{eq:local-inverse-composition-CN-material}.
	For $k=0$ and $0\leq N\leq K^\sharp-1$, we have
	$$
	\nrm{R^*[V]\circ\Phi_\upsilon}_N
	\lesssim\mC_U\bar\varepsilon\MM{N,K_0,\lambda_{R,U,\upsilon},\tilde\lambda_{R,U,\upsilon}}.
	$$
	For $1\leq k\leq\min\{M_*,K^*+1\}$ and $k+N\leq K^\sharp-1$,
	$$
	\nrm{R^*[\partial_t^kV]\circ\Phi_\upsilon}_N
	\lesssim\mC_U\bar\varepsilon\MM{N,(K_U-k)_+,\lambda_{R,U,\upsilon},\tilde\lambda_{R,U,\upsilon}}
	\MM{k,N_t,\mu_{R,U,\upsilon},\tilde\mu_{R,U,\upsilon}}.
	$$

	The estimates for $P$ and its material derivatives are unchanged.
	For $M=0$ and $0\leq N\leq K^\sharp-1$, \eqref{eq:Rjl2} gives
	$$
	\begin{aligned}
		\nrm{R[U]}_N
		\lesssim\sum_{N_1+N_2=N}\nrm{P}_{N_1}\nrm{R^*[V]\circ\Phi_\upsilon}_{N_2}
		\lesssim\mC_U\bar\varepsilon\MM{N,K_0,\lambda_{R,U,\upsilon},\tilde\lambda_{R,U,\upsilon}}.
	\end{aligned}
	$$
	Thus \eqref{eq:Eulerian-R-CN-integer-source-zero} holds.
	For $1\leq M\leq\min\{M_*,K^*+1\}$ and $M+N\leq K^\sharp-1$,
	use \eqref{eq:higher-material-R} and \eqref{eq:local-inverse-CN-product} with the new composition estimates.
	Lowering the spatial thresholds to $(K_U-M)_+$ gives
	$$
	\begin{aligned}
		\nrm{D_{t,\upsilon}^{M}R[U]}_N
		&\lesssim\sum_{k=0}^M\sum_{N_1+N_2=N}
		\nrm{D_{t,\upsilon}^{M-k}P}_{N_1}\nrm{R^*[\partial_t^kV]\circ\Phi_\upsilon}_{N_2}\\
		&\lesssim\mC_U\bar\varepsilon\MM{N,(K_U-M)_+,\lambda_{R,U,\upsilon},\tilde\lambda_{R,U,\upsilon}}
		\sum_{k=0}^M\MM{M-k,N_t,\mu_{R,U,\upsilon},\tilde\mu_{R,U,\upsilon}}
		\MM{k,N_t,\mu_{R,U,\upsilon},\tilde\mu_{R,U,\upsilon}}\\
		&\lesssim\mC_U\bar\varepsilon\MM{N,(K_U-M)_+,\lambda_{R,U,\upsilon},\tilde\lambda_{R,U,\upsilon}}
		\MM{M,N_t,\mu_{R,U,\upsilon},\tilde\mu_{R,U,\upsilon}}.
	\end{aligned}
	$$
	This is \eqref{eq:Eulerian-R-CN-integer-source}.
	At orders $M,N$, this argument uses \eqref{eq:inv4divEq:integer-source-assumption} only through material order $M$ and total order $M+N$.
\end{proof}

\begin{proof}[Proof of Theorem~\ref{thm:inv4divEq}]
	Define $V$ and $R[U]$ by \eqref{eq:Rjl0}--\eqref{eq:Rjl2}.
	Proposition~\ref{prop:Rjl} gives the linearity of $U\mapsto R[U]$ and the identity
	$
	\Div R[U]=U.
	$
	The support properties are
	$
	\supp R[U]
	\subseteq
	\hat\Gamma_\upsilon(\bar\tau,\bar\varepsilon;t_0,x_0),
	\
	\supp_tR[U]\subseteq\supp_tU.
	$
	Thus \eqref{eq:inv4divEq:nonsymmetric-divergence}--%
	\eqref{eq:inv4divEq:supp4R} and \eqref{eq:time-supp-localR} hold.
	The Eulerian support inclusion
	\eqref{eq:inv4divEq:supp4R-eulerian} follows from
	\eqref{eq:supp4Rjl-eulerian}.
	Lemma~\ref{lem:estimate-pulled-back-sources} proves \eqref{eq:inv4divEq:Calpha-material-zero-MM}--\eqref{eq:inv4divEq:CN-material-MM-weak} with the combined scales.
	Under \eqref{eq:inv4divEq:integer-source-assumption},
	\eqref{eq:Eulerian-R-CN-integer-source-zero} and \eqref{eq:Eulerian-R-CN-integer-source} give
	\eqref{eq:inv4divEq:CN-from-integer-source-zero} and \eqref{eq:inv4divEq:CN-from-integer-source}.
\end{proof}

\subsection{Spatially localized sources and cutoffs in the Newton step}
\label{subsec:support-properties-R-cutoffs}

In the Newton step, finite propagation enlarges the spatial support
of each localized source. We choose the radius in
Theorem~\ref{thm:inv4divEq} so that the corresponding adapted Lagrangian
cylinder contains this support.
Assume the hypotheses of Proposition~\ref{prop:uniform-mollification-estimates}.
We use the cutoff notation from Subsection~\ref{subsec:cutoffs} and
$X_u,\Phi_u$ from \eqref{def of Xu and Phiu}.
Let $U_{u,\kappa}$ be a smooth vector field. Assume that there is a constant
$C>0$, independent of $q,u$, and $\kappa$, such that
$$
\supp U_{u,\kappa}
\subset
\left\{
(t,x)\in\R\times\R^3:
|t-t_u|<\frac23\tau_q,\ 
|\Phi_u(t,x)-x_\kappa|
\leq\frac{C}{2}\tau_q
\right\}.
$$

For every $(u,\kappa)\in\mathscr I$, set
\begin{equation}
	\label{def of R support parameters under cutoffs}
	\bar\tau:=\frac23\tau_q,
	\qquad
	\bar\varepsilon:=C\tau_q.
\end{equation}
For each $(u,\kappa)$, we use $R^{(\overline{u}_q)}_{u,\kappa}$ from
Theorem~\ref{thm:inv4divEq} on $|t-t_u|\leq\bar\tau$ with these
localization radii and the translated averaging kernels
$\zeta_{x_\kappa,\bar\varepsilon}$. The flow smallness condition is verified in the proof below.
The remaining hypotheses of Theorem~\ref{thm:inv4divEq} are verified when
the operator is applied in the Newton step.
The localized matrix fields satisfy the following support and overlap properties.

\begin{lem}
	\label{lem support of R under cutoffs}
	Under the support and localization assumptions above, assume in addition
	that, for every $(u,\kappa)\in\mathscr I$ and every fixed time $t$,
	the localized source $U_{u,\kappa}$ satisfies
	$$
	\int_{\R^3}(U_{u,\kappa})_i(t,x)\,\rd x=0,
	\qquad
	i=1,2,3.
	$$
	Set $\bar\varepsilon_E:=13C\tau_q/12$ and
	$S_{u,\kappa}:=\hat C_{\overline u_q}(\bar\tau,\bar\varepsilon_E;t_u,x_\kappa)$.
	There exists $a_{\mathrm{loc}}=a_{\mathrm{loc}}(\beta,\allowbreak b,\allowbreak \alpha,\allowbreak \rho^*,\allowbreak \mathfrak M,\allowbreak C,\allowbreak N_R,\allowbreak N_{\mathrm{ind}},\allowbreak N_{\mathrm{max}})\geq a_{\mathrm{init}}$.
	For $a>a_{\mathrm{loc}}$, the following statements hold.
	\begin{enumerate}
		\item For every $(u,\kappa)\in\mathscr I$,
		\begin{equation}
			\label{supp of localized U uk}
			\supp U_{u,\kappa}
			\subseteq
			\hat\Gamma_{\overline u_q}
			(\bar\tau,\bar\varepsilon;t_u,x_\kappa).
		\end{equation}
		
		\item For every $(u,\kappa)\in\mathscr I$,
		\begin{equation}
			\label{supp of R localized U uk Gamma}
			\supp R^{(\overline{u}_q)}_{u,\kappa}[U_{u,\kappa}]
			\subseteq
			\hat\Gamma_{\overline u_q}
			(\bar\tau,\bar\varepsilon;t_u,x_\kappa).
		\end{equation}
		Moreover,
		$$
		\supp_tR^{(\overline{u}_q)}_{u,\kappa}[U_{u,\kappa}]
		\subseteq
		\supp_tU_{u,\kappa}.
		$$
		In particular,
		\begin{equation}
			\label{supp of R localized U uk C}
			\supp R^{(\overline{u}_q)}_{u,\kappa}[U_{u,\kappa}]
			\subseteq
			S_{u,\kappa}.
		\end{equation}
		
		\item The family
		$\{S_{u,\kappa}\}_{(u,\kappa)\in\mathscr I}$ is locally finite and
		has uniformly bounded overlap. More precisely,
		\begin{equation}
			\label{finite overlap S uk R}
			\#\left\{
			(u',\kappa')\in\mathscr I:
			S_{u,\kappa}\cap S_{u',\kappa'}\neq\varnothing
			\right\}
			\leq
			3\left(
			2\left\lceil\frac{13C}{11\pi}\right\rceil+1
			\right)^3.
		\end{equation}
		Consequently, the family
		$
		\left\{
		\supp R^{(\overline{u}_q)}_{u,\kappa}[U_{u,\kappa}]
		\right\}_{(u,\kappa)\in\mathscr I}
		$
		is also locally finite and has uniformly bounded overlap.
		
		\item The locally finite sum
		$
		\sum_{(u,\kappa)\in\mathscr I}
		R^{(\overline{u}_q)}_{u,\kappa}[U_{u,\kappa}]
		$
		defines a smooth matrix field on
		$\R\times\R^3$. If
		\begin{equation}
			\label{eq:covariance of Newton cutoff inputs}
			U_{u,\kappa+\tau_q^{-1}n}(t,x+2\pi n)
			=
			U_{u,\kappa}(t,x)
		\end{equation}
		for every $n\in\Z^3$, then the sum is $2\pi$-periodic in each spatial variable
		and defines a smooth matrix field on $\R\times\T^3$.
	\end{enumerate}
\end{lem}

\begin{proof}
	In this proof, spatial norms are taken over $\R^3$.
	Use $\mC_\upsilon=C_0\delta_q^{\frac12}$ and $\lambda_\upsilon=\lambda_q$ in \eqref{eq:tau-local-smallness}.
	Choose $C_0\geq1$ at least as large as the implicit constant in \eqref{est:uniform-grad-om-over-orho} at $M=N=0$.
	Then $\|\nabla\overline u_q\|_0\leq C_0\Lambda_q$.
	By \eqref{def:tau-q}, choose $a_{\mathrm{loc}}\geq a_{\mathrm{init}}$ so that, for $a>a_{\mathrm{loc}}$,
	$$
	\begin{aligned}
		\bar\varepsilon=C\tau_q\leq1,\quad
		\bar\tau\mC_\upsilon\lambda_\upsilon
		=\frac23 C_0\tau_q\Lambda_q
		\leq\frac23 C_0\lambda_{q+1}^{-4\alpha}\leq c_0.
	\end{aligned}
	$$
	By \eqref{eq:geometry-condition-cutoff} and \eqref{eq:flow-gradient-close-id}, this threshold can also ensure, for $|t-t_u|\leq\tau_q$,
	$
	e^{2\bar\tau\|\nabla\overline u_q\|_0}\leq13/12,
	\
	\|\nabla X_u-\Id\|_0\leq1/12.
	$
	Consequently,
	$
	|X_u(t,y)-X_u(t,z)|\geq(11/12)|y-z|.
	$
	
	For $(t,x)\in\supp U_{u,\kappa}$, the support assumption gives
	$|t-t_u|<\bar\tau$ and
	$|\Phi_u(t,x)-x_\kappa|\leq C\tau_q/2<\bar\varepsilon$.
	We obtain \eqref{supp of localized U uk} by the definition of
	$\hat\Gamma_{\overline u_q}$.
	The zero-mean condition and \eqref{supp of localized U uk} allow the
	construction \eqref{eq:Rjl0}--\eqref{eq:Rjl2}. Its support proof in
	Proposition~\ref{prop:Rjl} gives \eqref{supp of R localized U uk Gamma}
	and the temporal support inclusion.
	By Lemma~\ref{lem:comparison-Gamma-C}, we have
	$$
	\hat\Gamma_{\overline u_q}(\bar\tau,\bar\varepsilon;t_u,x_\kappa)
	\subseteq
	\hat C_{\overline u_q}\left(\bar\tau,
	e^{2\bar\tau\|\nabla\overline u_q\|_0}\bar\varepsilon;t_u,x_\kappa\right)
	\subseteq S_{u,\kappa},
	$$
	where the last inclusion follows from
	$
	e^{2\bar\tau\|\nabla\overline u_q\|_0}\bar\varepsilon
	\leq(13/12)C\tau_q=\bar\varepsilon_E.
	$
	The support inclusion \eqref{supp of R localized U uk C} follows.
	
	To prove \eqref{finite overlap S uk R}, fix $(u,\kappa)\in\mathscr I$
	and suppose $S_{u,\kappa}\cap S_{u',\kappa'}\neq\varnothing$.
	Choose $(t,x)$ in this intersection.
	The time intervals imply
	$$
	|u-u'|\tau_q\leq2\bar\tau=4\tau_q/3,
	$$
	so $u'\in\{u-1,u,u+1\}$. By the definition of these cylinders,
	$|x-X_u(t,x_\kappa)|\leq\bar\varepsilon_E$ and
	$|x-X_{u'}(t,x_{\kappa'})|\leq\bar\varepsilon_E$.
	Hence
	$$
	|X_u(t,x_\kappa)-X_{u'}(t,x_{\kappa'})|
	\leq2\bar\varepsilon_E.
	$$
	Since $|t_{u'}-t_u|\leq\tau_q$ and $|t-t_{u'}|\leq\bar\tau<\tau_q$,
	the flow identity
	$
	X_u(t,x_\kappa)=X_{u'}(t,X_u(t_{u'},x_\kappa))
	$
	and the lower Lipschitz bound for $X_{u'}(t,\cdot)$ give
	$$
	|X_u(t_{u'},x_\kappa)-x_{\kappa'}|
	\leq\frac{12}{11}\,2\bar\varepsilon_E
	=\frac{26C}{11}\tau_q.
	$$
	Since $x_{\kappa'}=2\pi\tau_q\kappa'$, we have
	$$
	\left|\kappa'-\frac{X_u(t_{u'},x_\kappa)}{2\pi\tau_q}\right|
	\leq\frac{13C}{11\pi}.
	$$
	For each fixed $u'$, the center depends only on $(u,\kappa)$.
	Each coordinate of $\kappa'$ therefore has at most
	$2\lceil13C/(11\pi)\rceil+1$ possible integer values.
	Multiplying the resulting spatial count by the three choices of $u'$
	proves \eqref{finite overlap S uk R}.
	
	We now prove local finiteness. Let
	$K\subset\R\times\R^3$ be compact. If
	$K\cap S_{u,\kappa}\neq\varnothing$, then the time projection of $K$
	intersects $[t_u-\bar\tau,t_u+\bar\tau]$. Since $t_u=u\tau_q$, only
	finitely many values of $u$ are possible. For each such $u$, choose
	$(t,x)\in K\cap S_{u,\kappa}$. Then
	$
	|X_u(t,x_\kappa)-x|
	\leq
	\bar\varepsilon_E.
	$
	The bi-Lipschitz estimate for the flow gives
	$$
	\frac{11}{12}|x_\kappa|
	\leq
	|X_u(t,x_\kappa)-X_u(t,0)|
	\leq
	|x|+\bar\varepsilon_E+|X_u(t,0)|.
	$$
	The right-hand side is uniformly bounded for $(t,x)\in K$ and for the
	finitely many admissible values of $u$. Hence $|x_\kappa|$ is bounded.
	Since $x_\kappa=2\pi\tau_q\kappa$, only finitely many values of $\kappa$
	are possible. Thus $\{S_{u,\kappa}\}_{(u,\kappa)\in\mathscr I}$ is locally
	finite. By \eqref{supp of R localized U uk C}, the family of supports of the local matrix fields is also locally finite and has uniformly bounded
	overlap. Consequently, their sum is locally finite and smooth.
	
	To prove periodicity, let $n\in\Z^3$. Periodicity of
	$\overline u_q$ gives
	$
	\Phi_u(t,x+2\pi n)
	=
	\Phi_u(t,x)+2\pi n.
	$
	Since $\tau_q^{-1}\in2\N$, we also have
	$
	x_{\kappa+\tau_q^{-1}n}
	=
	x_\kappa+2\pi n.
	$
	The translated averaging kernels, \eqref{eq:Rstar-formula}, and
	\eqref{eq:covariance of Newton cutoff inputs} then give
	$$
	R^{(\overline{u}_q)}_{u,\kappa+\tau_q^{-1}n}
	[U_{u,\kappa+\tau_q^{-1}n}](t,x+2\pi n)
	=
	R^{(\overline{u}_q)}_{u,\kappa}[U_{u,\kappa}](t,x).
	$$
	Reindexing the locally finite sum in $\kappa$ proves periodicity.
\end{proof}

\section{Nonlocal operators and material derivative estimates}
\label{sec:periodic-Hodge-and-Newton-stresses}

We establish spatial and material derivative estimates for the nonlocal operators
and the linearized equation used in the Newton step.
We begin with the periodic Hodge potentials in Subsection~\ref{subsec:periodic-Hodge-decomposition}.
In Subsection~\ref{subsec:symmetric-inverse-divergence-curl-sources}, we estimate
the operators $\mathcal R\curl$ and $\mathcal B\curl$.
Finally, Subsection~\ref{subsec:linearized-equation-material-estimates} applies these operator estimates
to the linearized equation with a curl source and estimates the solution, its vector potential, and the pressure.

To prove the higher material derivative estimates, we combine transport estimates
with periodic elliptic estimates and proceed by induction on the material order.
At each order, we differentiate the divergence and curl identities and estimate
the resulting commutator terms using \eqref{eq:first-order-material-spatial-commutator}
and H{\"o}lder product estimates at a finite sequence of exponents.
Since the background field need not be divergence-free, \eqref{eq:first-order-material-mean-identity}
shows that material derivatives need not preserve zero mean. We therefore compute the spatial means separately.
This div--curl argument parallels the proof of \cite[Theorem~3]{GST12}.
The identities for iterated material derivatives in that setting are given in
\cite[Proposition~6]{GST12}.

Throughout this section, $J\subset\mathbb R$ is a compact interval and
$0<\alpha<1$. For $r\in\mathbb N_0$, set $\alpha_r:=2^{-r}\alpha$, so that
$\alpha_0=\alpha$ and $2\alpha_{r+1}=\alpha_r$. For $s\in\mathbb N$ and
$0<\beta<1$, write
$\Lkb{N,K,\lambda,\tilde\lambda}{s}{\beta}:=\left(\Xi_\beta\left(N,K,\lambda,\tilde\lambda\right)\right)^s$.
Only finitely many of the exponents $\alpha_r$ are used in each argument. Let $\upsilon$ be a
smooth vector field on $J\times\T^3$, and set
$D_{t,\upsilon}:=\partial_t+\upsilon\cdot\nabla$. 

\subsection{The periodic Hodge decomposition}
\label{subsec:periodic-Hodge-decomposition}

We use the periodic Hodge decomposition to write the divergence of a
matrix field as the sum of a gradient and a curl. We first record the
commutator and mean identities used in the material derivative estimates.
For every smooth scalar, vector, or tensor field $h$, we have
\begin{equation}
	\label{eq:first-order-material-spatial-commutator}
	[D_{t,\upsilon},\partial_i]h=-(\partial_i\upsilon_j)\partial_jh.
\end{equation}
Consequently, for every smooth scalar field $\varphi$, vector field $Y$, and
matrix field $A$,
\begin{align*}
	[D_{t,\upsilon},\nabla]\varphi
	&=-(\nabla\upsilon)^T\nabla\varphi,
	\\
	[D_{t,\upsilon},\Div]Y
	&=-(\partial_i\upsilon_j)\partial_jY_i,
	\\
	\bigl([D_{t,\upsilon},\curl]Y\bigr)_i
	&=-\epsilon_{ijk}(\partial_j\upsilon_{l})\partial_{l}Y_k,
	\\
	\bigl([D_{t,\upsilon},\Div]A\bigr)_i
	&=-(\partial_j\upsilon_{k})\partial_{k}A_{ij}.
\end{align*}
Periodic integration by parts also gives
\begin{equation}
	\label{eq:first-order-material-mean-identity}
	\fint_{\T^3}D_{t,\upsilon}h\,\rd x
	=\frac{\rd}{\rd t}\fint_{\T^3}h\,\rd x-\fint_{\T^3}h\Div\upsilon\,\rd x.
\end{equation}

For a smooth periodic matrix field $T$ on $J\times\T^3$, define
\begin{equation}
	\label{def:periodic-Hodge-potentials}
	Q[T]:=\Delta^{-1}\Div\Div T,\qquad W[T]:=-\Delta^{-1}\curl\Div T.
\end{equation}
Both potentials have zero spatial mean. Since $\Div T$ also has zero mean,
the periodic Hodge decomposition gives
\begin{equation}
	\label{eq:periodic-Hodge-decomposition}
	\Div T=\nabla Q[T]+\curl W[T],\qquad \Div W[T]=0.
\end{equation}
For $N\in\mathbb N_0$ and $0<\alpha<1$, apply
\eqref{eq:app-multiplier-besov-estimate} to the zero-mean part of $T$, with
$p=q=\infty$, $s=N+\alpha$, and multiplier order $m=0$.
Using Remark~\ref{rem:app-holder-besov-equivalence}, we obtain
\begin{equation}
	\label{eq:periodic-Hodge-material-zero-estimate}
	\nrm{Q[T]}_{N+\alpha,J}+\nrm{W[T]}_{N+\alpha,J}
	\lesssim_{N,\alpha}\nrm{T}_{N+\alpha,J}.
\end{equation}
The next proposition estimates higher material derivatives of these
potentials, with separate spatial and material scales for $T$ and $\upsilon$.

\begin{prop}
	\label{prop:periodic-Hodge-material-estimates}
	Fix $K^*,K^\sharp\in\mathbb N$ and $K_T,M_*,N_t\in\mathbb N_0$
	with $K_T\leq K^*$ and $K^\sharp\geq K^*+1$.
	If $M_*\geq K^*+1$, assume also $N_t\leq K^*$.
	Let $A_\upsilon>0$ and choose scales for $\upsilon$ satisfying
	$
	1\leq\lambda_\upsilon\leq\tilde\lambda_\upsilon,\
	1\leq\mu_\upsilon\leq\tilde\mu_\upsilon.
	$
	Assume the following $C^N$ estimate for the material derivatives of $\nabla\upsilon$
	\begin{equation}
		\label{eq:periodic-Hodge-upsilon-integer-assumption}
		\begin{aligned}
			\nrm{D_{t,\upsilon}^{M}\nabla\upsilon}_{N,J}
			&\leq A_\upsilon
			\MM{N+1,K^*-M,\lambda_\upsilon,\tilde\lambda_\upsilon}
			\MM{M,N_t,\mu_\upsilon,\tilde\mu_\upsilon},
			&& 0\leq M\leq K^*,\qquad M+N\leq K^\sharp-1.
		\end{aligned}
	\end{equation}
	At each exponent
	$\alpha_0,\ldots,\alpha_{K^*+1}$, assume also the H{\"o}lder estimates
	\begin{align}
		\nrm{D_{t,\upsilon}^{M}\nabla\upsilon}_{N+\alpha_r,J}
		&\leq A_\upsilon
		\Lab{N+1,K^*-M,\lambda_\upsilon,\tilde\lambda_\upsilon}{\alpha_r}\notag\\
		&\quad\cdot
		\MM{N+1,K^*-M,\lambda_\upsilon,\tilde\lambda_\upsilon}
		\MM{M,N_t,\mu_\upsilon,\tilde\mu_\upsilon},
		&& 0\leq M\leq K^*,\qquad M+N\leq K^\sharp-2.
		\label{eq:periodic-Hodge-upsilon-Holder-assumption}
	\end{align}

	Fix $s\in\mathbb N$ and $A_T>0$. Choose scales for $T$ satisfying
	$
	1\leq\lambda_T\leq\tilde\lambda_T,\
	1\leq\mu_T\leq\tilde\mu_T.
	$
	Fix $r^*\in\mathbb N_0$ with $r^*\leq K^*+1$.
	For $M,N,r\in\mathbb N_0$ with $M+r\leq r^*$, assume
	\begin{align}
		\nrm{D_{t,\upsilon}^{M}T}_{N+\alpha_r,J}
		&\leq A_T\Lkb{N,(K_T-M)_+,\lambda_T,\tilde\lambda_T}{s}{\alpha_r}\notag\\
		&\quad\cdot
		\MM{N,(K_T-M)_+,\lambda_T,\tilde\lambda_T}
		\MM{M,N_t,\mu_T,\tilde\mu_T},
		&& 0\leq M\leq M_*,\qquad M+N\leq K^\sharp-1.
		\label{eq:periodic-Hodge-T-Holder-assumption}
	\end{align}
	
	Set
	$$
	\begin{aligned}
		\lambda:=\max\{\lambda_\upsilon,\lambda_T\},\qquad
		\tilde\lambda:=\max\{\tilde\lambda_\upsilon,\tilde\lambda_T\},\qquad
		\mu:=\max\{\mu_\upsilon,\mu_T,A_\upsilon\lambda_\upsilon\},\qquad
		\tilde\mu:=\max\{\tilde\mu_\upsilon,\tilde\mu_T,A_\upsilon\tilde\lambda_\upsilon\}.
	\end{aligned}
	$$
	Then, for every $M,N\in\mathbb N_0$ satisfying
	$1\leq M\leq \min\{M_*,r^*\}$ and $M+N\leq K^\sharp-1$,
	\begin{equation}
		\label{eq:Hodge-high-material-estimate}
		\begin{aligned}
		&\nrm{D_{t,\upsilon}^{M}Q[T]}_{N+\alpha,J}
		+\nrm{D_{t,\upsilon}^{M}W[T]}_{N+\alpha,J}\\
		&\qquad\lesssim A_T\Lk{N,(K_T-M)_+,\lambda,\tilde\lambda}{\max\{2,s\}}
		\MM{N,(K_T-M)_+,\lambda,\tilde\lambda}
		\MM{M,N_t,\mu,\tilde\mu}.
		\end{aligned}
	\end{equation}
	At material order zero, for every $0\leq N\leq K^\sharp-1$,
	\begin{equation}
		\label{eq:Hodge-material-zero-estimate}
		\nrm{Q[T]}_{N+\alpha,J}+\nrm{W[T]}_{N+\alpha,J}
		\lesssim A_T\Lk{N,K_T,\lambda_T,\tilde\lambda_T}{s}
		\MM{N,K_T,\lambda_T,\tilde\lambda_T}.
	\end{equation}

	For $r\in\mathbb N_0$ with $M+r\leq r^*$, the same estimates hold with
	$\alpha$ replaced by $\alpha_r$, in the stated derivative ranges. It suffices to assume \eqref{eq:periodic-Hodge-T-Holder-assumption}
	at material orders $k\leq M$ and exponents $\alpha_{r'}$ with $k+r'\leq M+r$.
	The implicit constants depend only on $K^*,K^\sharp,\alpha$.
	They are independent of $J,\upsilon,T$, the amplitudes, and the spatial and material scales.
\end{prop}

\begin{proof}
	The case $M=0$ follows from
	\eqref{eq:periodic-Hodge-material-zero-estimate} and
	\eqref{eq:periodic-Hodge-T-Holder-assumption}. More precisely, for
	$0\leq r\leq r^*$ and $0\leq N\leq K^\sharp-1$, we have
	$$
	\nrm{Q[T]}_{N+\alpha_r,J}+\nrm{W[T]}_{N+\alpha_r,J}
	\lesssim A_T\Lkb{N,K_T,\lambda_T,\tilde\lambda_T}{s}{\alpha_r}
	\MM{N,K_T,\lambda_T,\tilde\lambda_T}.
	$$
	For $m\geq0$, write
	$$
	T^{[m]}:=D_{t,\upsilon}^{m}T,\quad
	Q^{[m]}:=D_{t,\upsilon}^{m}Q[T],\quad
	W^{[m]}:=D_{t,\upsilon}^{m}W[T].
	$$
	Bracketed superscripts denote material derivatives, while subscripts on auxiliary
	fields denote their order in the recursive construction.
	
	We argue by induction on $M$, simultaneously at the exponents $\alpha_r$.
	For $1\leq M\leq M_*$ and $r,N\in\mathbb N_0$ with
	$M+r\leq r^*$ and $M+N\leq K^\sharp-1$, we prove
	$$
	\begin{aligned}
		\nrm{D_{t,\upsilon}^{M}Q[T]}_{N+\alpha_r,J}
		+\nrm{D_{t,\upsilon}^{M}W[T]}_{N+\alpha_r,J}\lesssim
		A_T\Lkb{N,(K_T-M)_+,\lambda,\tilde\lambda}{\max\{2,s\}}{\alpha_r}
		\MM{N,(K_T-M)_+,\lambda,\tilde\lambda}
		\MM{M,N_t,\mu,\tilde\mu}.
	\end{aligned}
	$$
	Fix $M$ and assume the estimates at all lower material orders in the
	stated ranges. We derive the divergence, curl, and mean identities at order
	$M$, then use periodic elliptic estimates to reduce the induction step to
	products involving lower material derivatives of the potentials.
	We construct a matrix field $H_m$ and a vector field
	$K_m$, $0\leq m\leq M$, such that
	\begin{equation}
		\label{eq:material-Hodge-local-system}
		\nabla Q^{[m]}+\curl W^{[m]}=\Div(T^{[m]}+H_m),\qquad \Div W^{[m]}=\Div K_m.
	\end{equation}
	For $m=0$, \eqref{eq:periodic-Hodge-decomposition} gives
	\eqref{eq:material-Hodge-local-system} with $H_0=0$ and $K_0=0$.
	For a matrix field $A$, a scalar field $\varphi$, and a vector field $Y$,
	we rewrite \eqref{eq:first-order-material-spatial-commutator} as
	\begin{equation}\label{eq:material-Hodge-divergence-form-commutators}
	\begin{aligned}
		\bigl([D_{t,\upsilon},\Div]A\bigr)_i
		&=\partial_{k}\left(-(\partial_j\upsilon_{k})A_{ij}
		+(\Div\upsilon)A_{ik}\right)
		-(\Div\upsilon)\partial_jA_{ij},\\
		\bigl([D_{t,\upsilon},\nabla]\varphi\bigr)_i
		&=\partial_{k}\left(-(\partial_i\upsilon_{k})\varphi
		+(\Div\upsilon)\varphi\delta_{ik}\right)
		-(\Div\upsilon)\partial_i\varphi,\\
		\bigl([D_{t,\upsilon},\curl]Y\bigr)_i
		&=\partial_{l}\left(-\epsilon_{ijk}(\partial_j\upsilon_{l})Y_k
		+(\Div\upsilon)\epsilon_{ilk}Y_k\right)
		-(\Div\upsilon)(\curl Y)_i,\\
		[D_{t,\upsilon},\Div]Y
		&=\partial_{k}\left(-(\partial_i\upsilon_{k})Y_i
		+(\Div\upsilon)Y_{k}\right)
		-(\Div\upsilon)\Div Y.
	\end{aligned}
	\end{equation}
	
	Suppose that \eqref{eq:material-Hodge-local-system} holds at order $m$.
	Applying $D_{t,\upsilon}$ to its first identity, we obtain
	$$
	\begin{aligned}
		\nabla Q^{[m+1]}+\curl W^{[m+1]}
		&=\Div\bigl(T^{[m+1]}+D_{t,\upsilon}H_m\bigr)
		+[D_{t,\upsilon},\Div](T^{[m]}+H_m)-[D_{t,\upsilon},\nabla]Q^{[m]}
		-[D_{t,\upsilon},\curl]W^{[m]}.
	\end{aligned}
	$$
	Substituting \eqref{eq:material-Hodge-divergence-form-commutators}, we have, componentwise,
	$$
	\begin{aligned}
		\bigl(\nabla Q^{[m+1]}+\curl W^{[m+1]}\bigr)_i
		=&\partial_{l}\Bigl((T^{[m+1]})_{il}
		+D_{t,\upsilon}(H_m)_{il}
		-(\partial_j\upsilon_{l})(T^{[m]}+H_m)_{ij}
		+(\Div\upsilon)(T^{[m]}+H_m)_{il}\\
		&\hspace{15pt}+(\partial_i\upsilon_{l})Q^{[m]}
		-(\Div\upsilon)Q^{[m]}\delta_{il}
		+\epsilon_{ijk}(\partial_j\upsilon_{l})(W^{[m]})_k
		-(\Div\upsilon)\epsilon_{ilk}(W^{[m]})_k\Bigr)\\
		&-(\Div\upsilon)\bigl(\Div(T^{[m]}+H_m)\bigr)_i
		+(\Div\upsilon)(\nabla Q^{[m]})_i
		+(\Div\upsilon)(\curl W^{[m]})_i.
	\end{aligned}
	$$
	The last three terms cancel by the first identity in
	\eqref{eq:material-Hodge-local-system}, giving
	$$
	-(\Div\upsilon)\Div(T^{[m]}+H_m)
	+(\Div\upsilon)\nabla Q^{[m]}
	+(\Div\upsilon)\curl W^{[m]}=0.
	$$
	For the second identity in \eqref{eq:material-Hodge-local-system}, we
	similarly obtain
	$$
	\begin{aligned}
		\Div W^{[m+1]}
		=D_{t,\upsilon}\Div W^{[m]}-[D_{t,\upsilon},\Div]W^{[m]}
		=D_{t,\upsilon}\Div K_m-[D_{t,\upsilon},\Div]W^{[m]}
		=\Div(D_{t,\upsilon}K_m)
		+[D_{t,\upsilon},\Div](K_m-W^{[m]}).
	\end{aligned}
	$$
	Since $\Div(K_m-W^{[m]})=0$, the last commutator identity reads
	$$
	\begin{aligned}
		\relax[D_{t,\upsilon},\Div](K_m-W^{[m]})
		=\partial_{k}\left(
		-(\partial_i\upsilon_{k})(K_m-W^{[m]})_i
		+(\Div\upsilon)(K_m-W^{[m]})_{k}
		\right).
	\end{aligned}
	$$
	Consequently,
	$$
	\begin{aligned}
		\Div W^{[m+1]}
		=\partial_{k}\Bigl(
		D_{t,\upsilon}(K_m)_{k}
		-(\partial_i\upsilon_{k})(K_m-W^{[m]})_i
		+(\Div\upsilon)(K_m-W^{[m]})_{k}
		\Bigr).
	\end{aligned}
	$$
	We therefore define
	\begin{align*}
		(H_{m+1})_{il}
		&:=D_{t,\upsilon}(H_m)_{il}
		-(\partial_j\upsilon_{l})(T^{[m]}+H_m)_{ij}
		+(\Div\upsilon)(T^{[m]}+H_m)_{il}\\
		&\quad+(\partial_i\upsilon_{l})Q^{[m]}
		-(\Div\upsilon)Q^{[m]}\delta_{il}
		+\epsilon_{ijk}(\partial_j\upsilon_{l})(W^{[m]})_k
		-(\Div\upsilon)\epsilon_{ilk}(W^{[m]})_k,
		\\
		(K_{m+1})_{k}
		&:=D_{t,\upsilon}(K_m)_{k}
		-(\partial_i\upsilon_{k})(K_m-W^{[m]})_i
		+(\Div\upsilon)(K_m-W^{[m]})_{k}.
	\end{align*}
	With these definitions,
	$$
	\nabla Q^{[m+1]}+\curl W^{[m+1]}
	=\Div(T^{[m+1]}+H_{m+1}),
	\qquad
	\Div W^{[m+1]}=\Div K_{m+1}.
	$$
	This proves \eqref{eq:material-Hodge-local-system} by induction on the
	material order.
	
	Each component of $H_M$ or $K_M$ is a finite sum of tensor contractions
	of products of the form
	\begin{equation}
		\label{eq:material-Hodge-correction-product}
		\left(
		\prod_{\nu=1}^{\kappa}
		D_{t,\upsilon}^{M_\nu}\nabla\upsilon
		\right)
		D_{t,\upsilon}^{k}Z,
		\qquad
		Z\in\{T,Q[T],W[T]\},
	\end{equation}
	where $\kappa\geq1$, $0\leq k\leq M-1$, and
	\begin{equation}
		\label{eq:material-Hodge-correction-order}
		k+\sum_{\nu=1}^{\kappa}(M_\nu+1)=M.
	\end{equation}
	Since $H_0=K_0=0$, the first corrections are
	$$
	\begin{aligned}
		(H_1)_{il}
		&=-(\partial_j\upsilon_{l})T_{ij}
		+(\Div\upsilon)T_{il}
		+(\partial_i\upsilon_{l})Q[T]
		-(\Div\upsilon)Q[T]\delta_{il}+\epsilon_{ijk}(\partial_j\upsilon_{l})W[T]_k
		-(\Div\upsilon)\epsilon_{ilk}W[T]_k,\\
		(K_1)_{k}
		&=(\partial_i\upsilon_{k})W[T]_i
		-(\Div\upsilon)W[T]_{k}.
	\end{aligned}
	$$
	These terms satisfy \eqref{eq:material-Hodge-correction-product}--\eqref{eq:material-Hodge-correction-order} with $\kappa=1$ and $M_1=k=0$.
	When $D_{t,\upsilon}$ differentiates a product, it raises one of the
	material indices by one, while multiplication by $\nabla\upsilon$ adds a
	factor of material order zero. The recursions therefore preserve the product
	form and order identity.
	
	Taking the divergence of the first identity in
	\eqref{eq:material-Hodge-local-system} gives
	$$
	\Delta Q^{[M]}=\Div\Div(T^{[M]}+H_M).
	$$
	Taking the curl of the same identity, we have
	$$
	\curl\curl W^{[M]}=\curl\Div(T^{[M]}+H_M).
	$$
	Using
	$$
	\curl\curl W^{[M]}=\nabla\Div W^{[M]}-\Delta W^{[M]},\quad\Div W^{[M]}=\Div K_M,$$ we obtain
	$$
	\Delta W^{[M]}
	=\nabla\Div K_M-\curl\Div(T^{[M]}+H_M).
	$$
	Inverting the Laplacian and restoring the spatial means, we find
	\begin{align*}
		Q^{[M]}
		&=\Delta^{-1}\Div\Div(T^{[M]}+H_M)
		+\fint_{\T^3}Q^{[M]}\,\rd x,
		\\
		W^{[M]}
		&=\nabla\Delta^{-1}\Div K_M
		-\Delta^{-1}\curl\Div(T^{[M]}+H_M)
		+\fint_{\T^3}W^{[M]}\,\rd x.
	\end{align*}
	The operators $\Delta^{-1}\Div\Div$, $\nabla\Delta^{-1}\Div$, and $\Delta^{-1}\curl\Div$ are periodic Fourier multipliers of order zero.
	Apply \eqref{eq:app-multiplier-besov-estimate} to the zero-mean parts of the fields and use Remark~\ref{rem:app-holder-besov-equivalence}. For $0\leq r\leq K^*+1$ and $N\in\mathbb N_0$, the resulting estimate is
	\begin{equation}
		\label{eq:material-Hodge-reconstruction-estimate}
		\begin{aligned}
			&\nrm{Q^{[M]}}_{N+\alpha_r,J}+\nrm{W^{[M]}}_{N+\alpha_r,J}\\
			&\quad\lesssim\nrm{T^{[M]}}_{N+\alpha_r,J}
			+\nrm{H_M}_{N+\alpha_r,J}+\nrm{K_M}_{N+\alpha_r,J}
			+\sup_{t\in J}\left|\fint_{\T^3}Q^{[M]}\,\rd x\right|
			+\sup_{t\in J}\left|\fint_{\T^3}W^{[M]}\,\rd x\right|.
		\end{aligned}
	\end{equation}
	
	In \eqref{eq:material-Hodge-reconstruction-estimate}, the $C^{N+\alpha_r}$ norm of $T^{[M]}$ is estimated by
	\eqref{eq:periodic-Hodge-T-Holder-assumption}. It remains to estimate $H_M$, $K_M$, and the two spatial means.
	For the spatial means, use
	\eqref{eq:first-order-material-mean-identity}.
	For $Z\in\{Q[T],W[T]\}$, its zero spatial mean gives
	$$
	\fint_{\T^3}D_{t,\upsilon}Z\,\rd x
	=-\fint_{\T^3}Z\Div\upsilon\,\rd x.
	$$
	For every $M\geq0$, applying
	\eqref{eq:first-order-material-mean-identity} to
	$D_{t,\upsilon}^{M}Z$ gives the recursion
	$$
	\fint_{\T^3}D_{t,\upsilon}^{M+1}Z\,\rd x
	=\frac{\rd}{\rd t}\fint_{\T^3}D_{t,\upsilon}^{M}Z\,\rd x
	-\fint_{\T^3}D_{t,\upsilon}^{M}Z\Div\upsilon\,\rd x.
	$$
	For $M=2$,
	we first compute
	$$
	\begin{aligned}
		\frac{\rd}{\rd t}\fint_{\T^3}Z\Div\upsilon\,\rd x
		&=\fint_{\T^3}D_{t,\upsilon}(Z\Div\upsilon)\,\rd x+\fint_{\T^3}Z(\Div\upsilon)^2\,\rd x\\
		&=\fint_{\T^3}(D_{t,\upsilon}\Div\upsilon)Z\,\rd x+\fint_{\T^3}(\Div\upsilon)D_{t,\upsilon}Z\,\rd x+\fint_{\T^3}Z(\Div\upsilon)^2\,\rd x,
	\end{aligned}
	$$
	so, applying \eqref{eq:first-order-material-mean-identity} to $D_{t,\upsilon}Z$, we obtain
	$$
	\begin{aligned}
		\fint_{\T^3}D_{t,\upsilon}^{2}Z\,\rd x
		=-\fint_{\T^3}(D_{t,\upsilon}\Div\upsilon)Z\,\rd x
		-2\fint_{\T^3}(\Div\upsilon)D_{t,\upsilon}Z\,\rd x
		-\fint_{\T^3}Z(\Div\upsilon)^2\,\rd x.
	\end{aligned}
	$$
	
	For every $M\geq1$, this recursion expresses the mean of $D_{t,\upsilon}^{M}Z$ as a finite linear
	combination of expressions of the form
	\begin{equation}
		\label{eq:material-Hodge-mean-product}
		\fint_{\T^3}\left(\prod_{\nu=1}^{\kappa}
		D_{t,\upsilon}^{M_\nu}\Div\upsilon\right)
		D_{t,\upsilon}^{k}Z\,\rd x,
	\end{equation}
	where $\kappa\geq1$, $0\leq k\leq M-1$, and the indices satisfy
	\eqref{eq:material-Hodge-correction-order}.
	To pass from order $M$ to $M+1$, differentiate each average using
	$$
	\frac{\rd}{\rd t}\fint_{\T^3}h\,\rd x
	=\fint_{\T^3}D_{t,\upsilon}h\,\rd x
	+\fint_{\T^3}h\Div\upsilon\,\rd x.
	$$
	The resulting terms are
	$$
	\begin{aligned}
		\sum_{i=1}^{\kappa}\fint_{\T^3}
		\left(D_{t,\upsilon}^{M_i+1}\Div\upsilon
		\prod_{\substack{\nu=1\\ \nu\neq i}}^{\kappa}
		D_{t,\upsilon}^{M_\nu}\Div\upsilon\right)
		D_{t,\upsilon}^{k}Z\,\rd x,\quad
		\fint_{\T^3}\left(\prod_{\nu=1}^{\kappa}
		D_{t,\upsilon}^{M_\nu}\Div\upsilon\right)
		D_{t,\upsilon}^{k+1}Z\,\rd x,\quad
		\fint_{\T^3}\left(\prod_{\nu=1}^{\kappa}
		D_{t,\upsilon}^{M_\nu}\Div\upsilon\right)
		D_{t,\upsilon}^{k}Z\Div\upsilon\,\rd x.
	\end{aligned}
	$$
	Each term increases the total order in \eqref{eq:material-Hodge-correction-order} by one.
	The additional term in the mean recursion
	$
	-\fint_{\T^3}D_{t,\upsilon}^{M}Z\Div\upsilon\,\rd x
	$
	has $\kappa=1$, $M_1=0$, and $k=M$, which proves the mean expansion.
	
	Fix $r,N\in\mathbb N_0$ with $M+r\leq r^*$ and $M+N\leq K^\sharp-1$.
	We first estimate the products in \eqref{eq:material-Hodge-correction-product}, then the spatial means in \eqref{eq:material-Hodge-mean-product}.
	For $0\leq k\leq M-1$,
	\eqref{eq:holder-product-finite-Cnalpha} gives
	$$
	\begin{aligned}
		&\Nrm{
			\left(\prod_{\nu=1}^{\kappa}
			D_{t,\upsilon}^{M_\nu}\nabla\upsilon\right)
			D_{t,\upsilon}^{k}Z
		}_{N+\alpha_r,J}\\
		&\quad\lesssim
		\sum_{\substack{N_0+\sum_{\nu=1}^{\kappa}N_\nu=N\\
				N_0,N_1,\ldots,N_\kappa\geq0}}
		\Biggl[
		\nrm{D_{t,\upsilon}^{k}Z}_{N_0+\alpha_r,J}
		\prod_{\nu=1}^{\kappa}
		\nrm{D_{t,\upsilon}^{M_\nu}\nabla\upsilon}_{N_\nu,J}+
		\sum_{i=1}^{\kappa}
		\nrm{D_{t,\upsilon}^{M_i}\nabla\upsilon}_{N_i+\alpha_r,J}
		\nrm{D_{t,\upsilon}^{k}Z}_{N_0,J}
		\prod_{\substack{\nu=1\\ \nu\neq i}}^{\kappa}
		\nrm{D_{t,\upsilon}^{M_\nu}\nabla\upsilon}_{N_\nu,J}
		\Biggr].
	\end{aligned}
	$$
	By \eqref{eq:material-Hodge-correction-order}, the indices in this sum satisfy
	$k+N_0\leq M+N-1$ and $M_\nu+N_\nu\leq M+N-1\leq K^\sharp-2$.
	Thus \eqref{eq:periodic-Hodge-upsilon-integer-assumption} and
	\eqref{eq:periodic-Hodge-upsilon-Holder-assumption} control the $C^{N_\nu}$ and
	$C^{N_\nu+\alpha_r}$ norms of every coefficient $D_{t,\upsilon}^{M_\nu}\nabla\upsilon$, respectively.
	When $K^\sharp=K^*+1$, the restriction $M+N\leq K^\sharp-1$ gives $M+N\leq K^*$,
	so every coefficient satisfies $M_\nu+N_\nu\leq K^*-1$ and $M_\nu<K^*$.
	It remains to estimate the $C^{N_0+\alpha_r}$ and $C^{N_0}$ norms of $D_{t,\upsilon}^kZ$.
	Assumption~\eqref{eq:periodic-Hodge-T-Holder-assumption}, the induction hypothesis for $1\leq k\leq M-1$, and
	\eqref{eq:Hodge-material-zero-estimate} for $k=0$ give, respectively,
	$$
	\begin{aligned}
		\nrm{D_{t,\upsilon}^{k}T}_{N_0+\alpha_r,J}
		&\leq A_T\Lkb{N_0,(K_T-k)_+,\lambda_T,\tilde\lambda_T}{s}{\alpha_r}
		\MM{N_0,(K_T-k)_+,\lambda_T,\tilde\lambda_T}
		\MM{k,N_t,\mu_T,\tilde\mu_T},\\
		\nrm{D_{t,\upsilon}^{k}Q[T]}_{N_0+\alpha_r,J}
		+\nrm{D_{t,\upsilon}^{k}W[T]}_{N_0+\alpha_r,J}
		&\lesssim
		A_T\Lkb{N_0,(K_T-k)_+,\lambda,\tilde\lambda}{\max\{2,s\}}{\alpha_r}
		\MM{N_0,(K_T-k)_+,\lambda,\tilde\lambda}
		\MM{k,N_t,\mu,\tilde\mu},\\
		\nrm{Q[T]}_{N_0+\alpha_r,J}
		+\nrm{W[T]}_{N_0+\alpha_r,J}
		&\lesssim
		A_T\Lkb{N_0,K_T,\lambda_T,\tilde\lambda_T}{s}{\alpha_r}
		\MM{N_0,K_T,\lambda_T,\tilde\lambda_T}.
	\end{aligned}
	$$
	For each distribution of the $N$ spatial derivatives and $Z\in\{T,Q[T],W[T]\}$, the two terms in brackets satisfy
	\begin{align*}
		&\nrm{D_{t,\upsilon}^{k}Z}_{N_0+\alpha_r,J}
		\prod_{\nu=1}^{\kappa}\nrm{D_{t,\upsilon}^{M_\nu}\nabla\upsilon}_{N_\nu,J}\\
		&\quad\lesssim A_T
	\Lkb{N_0,(K_T-k)_+,\lambda,\tilde\lambda}{\max\{2,s\}}{\alpha_r}
	\MM{N_0,(K_T-k)_+,\lambda,\tilde\lambda}\MM{k,N_t,\mu,\tilde\mu}\\
	&\qquad\cdot\prod_{\nu=1}^{\kappa}\Bigl(
	A_\upsilon\MM{N_\nu+1,K^*-M_\nu,\lambda_\upsilon,\tilde\lambda_\upsilon}
	\MM{M_\nu,N_t,\mu_\upsilon,\tilde\mu_\upsilon}\Bigr)\\
	&\quad\lesssim A_T
	\Lkb{N,(K_T-M)_+,\lambda,\tilde\lambda}{\max\{2,s\}}{\alpha_r}
	\MM{N_0,(K_T-k)_+,\lambda,\tilde\lambda}\MM{k,N_t,\mu,\tilde\mu}\\
	&\qquad\cdot\prod_{\nu=1}^{\kappa}\Bigl(
	\MM{N_\nu,(K^*-M_\nu-1)_+,\lambda,\tilde\lambda}
	\MM{M_\nu+1,N_t,\mu,\tilde\mu}\Bigr)\\
		&\quad\lesssim
		A_T\Lkb{N,(K_T-M)_+,\lambda,\tilde\lambda}{\max\{2,s\}}{\alpha_r}
		\MM{N,(K_T-M)_+,\lambda,\tilde\lambda}\MM{M,N_t,\mu,\tilde\mu},\\
		&\sum_{i=1}^{\kappa}\nrm{D_{t,\upsilon}^{M_i}\nabla\upsilon}_{N_i+\alpha_r,J}
		\nrm{D_{t,\upsilon}^{k}Z}_{N_0,J}
		\prod_{\substack{\nu=1\\\nu\neq i}}^{\kappa}\nrm{D_{t,\upsilon}^{M_\nu}\nabla\upsilon}_{N_\nu,J}\\
		&\quad\leq\sum_{i=1}^{\kappa}\nrm{D_{t,\upsilon}^{M_i}\nabla\upsilon}_{N_i+\alpha_r,J}
		\nrm{D_{t,\upsilon}^{k}Z}_{N_0+\alpha_{r+1},J}
		\prod_{\substack{\nu=1\\\nu\neq i}}^{\kappa}\nrm{D_{t,\upsilon}^{M_\nu}\nabla\upsilon}_{N_\nu,J}\\
		&\quad\lesssim A_T
	\sum_{i=1}^{\kappa}\Bigl(
	\Lab{N_i+1,K^*-M_i,\lambda_\upsilon,\tilde\lambda_\upsilon}{\alpha_r}
	\Lkb{N_0,(K_T-k)_+,\lambda,\tilde\lambda}{\max\{2,s\}}{\alpha_{r+1}}\Bigr)\\
	&\qquad\cdot\MM{N_0,(K_T-k)_+,\lambda,\tilde\lambda}\MM{k,N_t,\mu,\tilde\mu}
	\prod_{\nu=1}^{\kappa}\Bigl(
	A_\upsilon\MM{N_\nu+1,K^*-M_\nu,\lambda_\upsilon,\tilde\lambda_\upsilon}
	\MM{M_\nu,N_t,\mu_\upsilon,\tilde\mu_\upsilon}\Bigr)\\
	&\quad\lesssim A_T
	\left(\Lab{N,(K_T-M)_+,\lambda,\tilde\lambda}{\alpha_r}\right)^{1+\frac{\max\{2,s\}}{2}}
	\MM{N_0,(K_T-k)_+,\lambda,\tilde\lambda}\MM{k,N_t,\mu,\tilde\mu}\\
	&\qquad\cdot\prod_{\nu=1}^{\kappa}\Bigl(
	\MM{N_\nu,(K^*-M_\nu-1)_+,\lambda,\tilde\lambda}
	\MM{M_\nu+1,N_t,\mu,\tilde\mu}\Bigr)\\
		&\quad\lesssim
		A_T\Lkb{N,(K_T-M)_+,\lambda,\tilde\lambda}{\max\{2,s\}}{\alpha_r}
		\MM{N,(K_T-M)_+,\lambda,\tilde\lambda}\MM{M,N_t,\mu,\tilde\mu},
	\end{align*}
	where the first $\lesssim$ inequality in each estimate follows from
	\eqref{eq:periodic-Hodge-upsilon-integer-assumption}--\eqref{eq:periodic-Hodge-upsilon-Holder-assumption}
	and the three estimates above.
	The estimates with H{\"o}lder exponents $\alpha_r$ and $\alpha_{r+1}$ apply since $k+N_0\leq M+N-1$ and $k+r+1\leq M+r\leq r^*$.
	We replace the individual scales by the combined scales and use $\Xi_{\alpha_r},\Xi_{\alpha_{r+1}}\geq1$ to replace $s$ by $\max\{2,s\}$.
	Since $K_T\leq K^*$ and $k,M_i\leq M-1$, we have
	$
	(K_T-k)_+\geq(K_T-M)_+,\
	(K^*-M_i-1)_+\geq(K_T-M)_+.
	$
	These threshold inequalities control the H{\"o}lder factors.
	For the sum over $i$, we estimate the product of the two H{\"o}lder factors using $2\alpha_{r+1}=\alpha_r$ and
	$1+\max\{2,s\}/2\leq\max\{2,s\}$.
	
	For the second $\lesssim$ inequality in each estimate, we use
	$$
	\MM{N_\nu+1,K^*-M_\nu,\lambda_\upsilon,\tilde\lambda_\upsilon}
	=
	\begin{cases}
		\lambda_\upsilon\MM{N_\nu,K^*-M_\nu-1,\lambda_\upsilon,\tilde\lambda_\upsilon}, & M_\nu\leq K^*-1,\\
		\tilde\lambda_\upsilon\MM{N_\nu,0,\lambda_\upsilon,\tilde\lambda_\upsilon}, & M_\nu=K^*.
	\end{cases}
	$$
	The second case can occur only when $M=K^*+1$.
	Estimate~\eqref{eq:appendix-derivative-weight-product} and the threshold
	comparisons above give
	\begin{equation}
		\label{eq:material-Hodge-spatial-scale-product}
		\MM{N_0,(K_T-k)_+,\lambda,\tilde\lambda}
		\prod_{\nu=1}^{\kappa}
		\MM{N_\nu,(K^*-M_\nu-1)_+,\lambda,\tilde\lambda}
		\leq
		\MM{N,(K_T-M)_+,\lambda,\tilde\lambda}.
	\end{equation}
	For the material scales, we use \eqref{eq:material-Hodge-correction-order} to obtain
	\begin{equation}
		\label{eq:material-Hodge-material-scale-product}
		\MM{k,N_t,\mu,\tilde\mu}
		\prod_{\nu=1}^{\kappa}\MM{M_\nu+1,N_t,\mu,\tilde\mu}
		\leq\MM{M,N_t,\mu,\tilde\mu}.
	\end{equation}
	For $M_\nu<K^*$, this spatial factor is absorbed using
	$A_\upsilon\lambda_\upsilon\MM{M_\nu,N_t,\mu_\upsilon,\tilde\mu_\upsilon}\leq\MM{M_\nu+1,N_t,\mu,\tilde\mu}$.
	For $M_\nu=K^*$, the assumption $N_t\leq K^*$ implies
	$A_\upsilon\tilde\lambda_\upsilon\MM{K^*,N_t,\mu_\upsilon,\tilde\mu_\upsilon}\leq\MM{K^*+1,N_t,\mu,\tilde\mu}$.
	The last inequality in each estimate follows from
	\eqref{eq:material-Hodge-spatial-scale-product}--\eqref{eq:material-Hodge-material-scale-product}.
	Summing over the spatial indices gives, for every product in
	\eqref{eq:material-Hodge-correction-product},
	$$
	\begin{aligned}
		&\Nrm{
			\left(\prod_{\nu=1}^{\kappa}
			D_{t,\upsilon}^{M_\nu}\nabla\upsilon\right)
			D_{t,\upsilon}^{k}Z
		}_{N+\alpha_r,J}\lesssim
		A_T\Lkb{N,(K_T-M)_+,\lambda,\tilde\lambda}{\max\{2,s\}}{\alpha_r}
		\MM{N,(K_T-M)_+,\lambda,\tilde\lambda}
		\MM{M,N_t,\mu,\tilde\mu}.
	\end{aligned}
	$$
	
	Since $\Div\upsilon=\operatorname{tr}(\nabla\upsilon)$, each integrand in
	\eqref{eq:material-Hodge-mean-product} is a contraction of the products just estimated.
	The preceding estimate with $N=0$ gives, for every $N$ in the stated range,
	$$
	\begin{aligned}
		&\sup_{t\in J}\left|
		\fint_{\T^3}
		\left(\prod_{\nu=1}^{\kappa}
		D_{t,\upsilon}^{M_\nu}\Div\upsilon\right)
		D_{t,\upsilon}^{k}Z\,\rd x
		\right|\lesssim
		A_T\Lkb{N,(K_T-M)_+,\lambda,\tilde\lambda}{\max\{2,s\}}{\alpha_r}
		\MM{N,(K_T-M)_+,\lambda,\tilde\lambda}
		\MM{M,N_t,\mu,\tilde\mu}.
	\end{aligned}
	$$
	
	Together, these estimates give
	$$
	\begin{aligned}
		&\nrm{H_M}_{N+\alpha_r,J}+\nrm{K_M}_{N+\alpha_r,J}
		+\sup_{t\in J}\left|\fint_{\T^3}Q^{[M]}\,\rd x\right|
		+\sup_{t\in J}\left|\fint_{\T^3}W^{[M]}\,\rd x\right|\\
		&\qquad\lesssim
		A_T\Lkb{N,(K_T-M)_+,\lambda,\tilde\lambda}{\max\{2,s\}}{\alpha_r}
		\MM{N,(K_T-M)_+,\lambda,\tilde\lambda}
		\MM{M,N_t,\mu,\tilde\mu}.
	\end{aligned}
	$$
	The remaining term in \eqref{eq:material-Hodge-reconstruction-estimate} is $T^{[M]}$, for which \eqref{eq:periodic-Hodge-T-Holder-assumption} gives
	$$
	\begin{aligned}
		\nrm{T^{[M]}}_{N+\alpha_r,J}
		&\leq A_T
		\Lkb{N,(K_T-M)_+,\lambda_T,\tilde\lambda_T}{s}{\alpha_r}
		\MM{N,(K_T-M)_+,\lambda_T,\tilde\lambda_T}
		\MM{M,N_t,\mu_T,\tilde\mu_T}.
	\end{aligned}
	$$
	Substitution into \eqref{eq:material-Hodge-reconstruction-estimate}
	closes the induction simultaneously for all $r$ with $M+r\leq r^*$.
\end{proof}

\subsection{The symmetric inverse-divergence operator and curl sources}
\label{subsec:symmetric-inverse-divergence-curl-sources}

We now apply the Hodge estimates to $\mathcal R\curl$ and
$\mathcal B\curl$. In the Newton step, differentiating a time cutoff $\theta(t)$
multiplying $\omega=\curl\Psi$ produces the term
$\theta'\omega=\curl(\theta'\Psi)$. Applying the periodic inverse-divergence
operator $\mathcal R$ from \cite{BDSV19,DS13} represents this term as the
divergence of a symmetric trace-free matrix. For the two-dimensional Newton
step, estimates for the stream function are proved in \cite[Lemma~3.7]{GR24}.
The operator $\mathcal B\curl$ maps $G$ to a zero-mean
divergence-free field with the same curl, so the linearized equation below
does not require $G$ to be divergence-free.

For every smooth periodic vector field $F$, define
\begin{equation}
	\label{def:periodic-symmetric-inverse-divergence}
	\begin{aligned}
		(\mathcal RF)_{ij}
		:=-\frac12\Delta^{-2}\partial_i\partial_j\partial_kF_k
		-\frac12\delta_{ij}\Delta^{-1}\partial_kF_k+\Delta^{-1}\partial_iF_j+\Delta^{-1}\partial_jF_i.
	\end{aligned}
\end{equation}
For every smooth zero-mean vector field
$h$ on $\T^3$, define the periodic vector-potential operator by
\begin{equation}
	\label{def:periodic-vector-potential-operator}
	\mathcal B h:=-\Delta^{-1}\curl h.
\end{equation}
Then $\Div\mathcal B h=0$ and $\fint_{\T^3}\mathcal B h\,\rd x=0$. If also
$\Div h=0$, then $\curl\mathcal B h=h$ and $\mathcal B h$ is the unique
zero-mean divergence-free vector potential of $h$.

\begin{lem}
	\label{lem:periodic-symmetric-inverse-divergence}
	For every smooth periodic vector field $F$, the matrix field
	$\mathcal RF$ has zero spatial mean and satisfies
	\begin{equation}
		\label{eq:periodic-symmetric-inverse-divergence-properties}
		(\mathcal RF)^T=\mathcal RF,\qquad
		\operatorname{tr}(\mathcal RF)=0,\qquad
		\Div\mathcal RF=F-\fint_{\T^3}F\,\rd x.
	\end{equation}
	Moreover, for every $N\in\mathbb N_0$ and $0<\alpha<1$,
	\begin{equation}
		\label{eq:periodic-symmetric-inverse-divergence-Schauder}
		\nrm{\mathcal RF}_{N+1+\alpha}
		\lesssim\nrm{F}_{N+\alpha}.
	\end{equation}
	The implicit constant depends only on $N,\alpha$ and is independent of $F$.
\end{lem}

\begin{proof}
	Symmetry and zero mean follow from \eqref{def:periodic-symmetric-inverse-divergence}.
	Taking the trace and divergence gives \eqref{eq:periodic-symmetric-inverse-divergence-properties},
	using $\Delta\Delta^{-1}F=F-\fint_{\T^3}F\,\rd x$.
	Each component of $\mathcal R$ is a periodic Fourier multiplier of order minus one.
	Applying \eqref{eq:app-multiplier-besov-estimate} to the zero-mean part of $F$
	and using Remark~\ref{rem:app-holder-besov-equivalence}, we obtain
	\eqref{eq:periodic-symmetric-inverse-divergence-Schauder}.
\end{proof}

For a smooth periodic vector field $G$, $\Div\curl G=0$, so
\eqref{def:periodic-symmetric-inverse-divergence} reduces to
\begin{equation}
	\label{def:direct-symmetric-curl-inverse}
	\bigl(\mathcal R\curl G\bigr)_{ij}
	=\epsilon_{ikl}\partial_k\partial_j\Delta^{-1}G_{l}
	+\epsilon_{jkl}\partial_k\partial_i\Delta^{-1}G_{l}.
\end{equation}
Equivalently,
\begin{equation}
	\label{eq:direct-symmetric-curl-standard-form}
	\mathcal R\curl G
	=\nabla\Delta^{-1}\curl G
	+\bigl(\nabla\Delta^{-1}\curl G\bigr)^T.
\end{equation}
Similarly, the vector identity
$\curl\curl G=\nabla\Div G-\Delta G$ gives
\begin{equation}
	\label{eq:direct-vector-curl-standard-form}
	\mathcal B\curl G=G-\nabla\Delta^{-1}\Div G-\fint_{\T^3}G\,\rd x.
\end{equation}
\begin{lem}
	\label{lem:direct-symmetric-curl-inverse}
	\label{lem:direct-vector-curl-inverse}
	For every smooth periodic vector field $G$, the matrix field
	$\mathcal R\curl G$ has zero spatial mean and satisfies
	\begin{equation}
		\label{eq:direct-symmetric-curl-properties}
		(\mathcal R\curl G)^T=\mathcal R\curl G,\qquad
		\operatorname{tr}(\mathcal R\curl G)=0,\qquad
		\Div(\mathcal R\curl G)=\curl G.
	\end{equation}
	The vector field
	$\mathcal B\curl G$ has zero spatial mean and satisfies
	\begin{equation}
		\label{eq:direct-vector-curl-properties}
		\Div(\mathcal B\curl G)=0,\qquad \curl(\mathcal B\curl G)=\curl G.
	\end{equation}
	For every $N\in\mathbb N_0$ and $0<\alpha<1$, we have
	\begin{equation}
		\label{eq:direct-symmetric-curl-Schauder}
		\nrm{\mathcal R\curl G}_{N+\alpha}
		\lesssim\nrm{G}_{N+\alpha}.
	\end{equation}
	\begin{equation}
		\label{eq:direct-vector-curl-Schauder}
		\nrm{\mathcal B\curl G}_{N+\alpha}\lesssim\nrm{G}_{N+\alpha}.
	\end{equation}
	One also has
	$\mathcal R\curl(\mathcal B\curl G)=\mathcal R\curl G$.
	If $G$ is divergence-free and has zero spatial mean, then
	$\mathcal B\curl G=G$.
	The implicit constants depend only on $N,\alpha$ and are independent of $G$.
\end{lem}

\begin{proof}
	Since $\curl G$ has zero spatial mean,
	Lemma~\ref{lem:periodic-symmetric-inverse-divergence} gives
	\eqref{eq:direct-symmetric-curl-properties} and the zero mean of $\mathcal R\curl G$.
	Identity~\eqref{eq:direct-vector-curl-standard-form} proves
	\eqref{eq:direct-vector-curl-properties}, the zero mean of $\mathcal B\curl G$,
	and $\mathcal B\curl G=G$ for divergence-free, zero-mean $G$. Applying $\mathcal R$ to the curl identity gives
	$\mathcal R\curl(\mathcal B\curl G)=\mathcal R\curl G$.
	By \eqref{def:direct-symmetric-curl-inverse} and
	\eqref{eq:direct-vector-curl-standard-form}, both operators are Fourier multipliers of order zero.
	Remark~\ref{rem:app-holder-besov-equivalence}, together with
	the multiplier estimate \eqref{eq:app-multiplier-besov-estimate}, proves
	\eqref{eq:direct-symmetric-curl-Schauder}--\eqref{eq:direct-vector-curl-Schauder}.
\end{proof}

The material derivative estimates follow from
Proposition~\ref{prop:periodic-Hodge-material-estimates}. We use its
notation and assumptions for $K^*,K^\sharp,N_t,\ \alpha$ and for the
transport field $\upsilon$, its amplitude, and its scales.

\begin{lem}
	\label{lem:direct-symmetric-curl-material-estimate}
	Fix $K_G,M_*\in\mathbb N_0$ with $K_G\leq K^*$.
	If $M_*\geq K^*+1$, assume also $N_t\leq K^*$.
	Assume \eqref{eq:periodic-Hodge-upsilon-integer-assumption}--%
	\eqref{eq:periodic-Hodge-upsilon-Holder-assumption}.
	Let $G$ be a smooth
	periodic vector field on $J\times\T^3$. Fix $s\in\mathbb N$, $A_G>0$, and
	scales
	$
	1\leq\lambda_G\leq\tilde\lambda_G,\
	1\leq\mu_G\leq\tilde\mu_G.
	$
	At each exponent
	$\alpha_0,\ldots,\alpha_{K^*+1}$, assume the H{\"o}lder estimate
	\begin{equation}
		\label{eq:symmetric-curl-G-Holder-assumption}
		\begin{aligned}
			\nrm{D_{t,\upsilon}^{M}G}_{N+\alpha_r,J}
			&\leq A_G\Lkb{N,(K_G-M)_+,\lambda_G,\tilde\lambda_G}{s}{\alpha_r}\\
			&\quad\cdot\MM{N,(K_G-M)_+,\lambda_G,\tilde\lambda_G}
			\MM{M,N_t,\mu_G,\tilde\mu_G},
			\qquad 0\leq M\leq M_*,\qquad M+N\leq K^\sharp-1.
		\end{aligned}
	\end{equation}
	
	Set
	$$
	\begin{aligned}
		\lambda:=\max\{\lambda_\upsilon,\lambda_G\},\qquad
		\tilde\lambda:=\max\{\tilde\lambda_\upsilon,\tilde\lambda_G\},\qquad
		\mu:=\max\{\mu_\upsilon,\mu_G,A_\upsilon\lambda_\upsilon\},\qquad
		\tilde\mu:=\max\{\tilde\mu_\upsilon,\tilde\mu_G,A_\upsilon\tilde\lambda_\upsilon\}.
	\end{aligned}
	$$
	Then, for every $M,N\in\mathbb N_0$ satisfying
	$1\leq M\leq \min\{M_*,K^*+1\}$ and $M+N\leq K^\sharp-1$,
	\begin{equation}
		\label{eq:direct-symmetric-curl-high-material}
		\begin{aligned}
			\nrm{D_{t,\upsilon}^{M}\mathcal R\curl G}_{N+\alpha,J}
			\lesssim A_G\Lk{N,(K_G-M)_+,\lambda,\tilde\lambda}{\max\{2,s\}}
			\MM{N,(K_G-M)_+,\lambda,\tilde\lambda}
			\MM{M,N_t,\mu,\tilde\mu}.
		\end{aligned}
	\end{equation}
	At material order zero, for every $0\leq N\leq K^\sharp-1$,
	\begin{equation}
		\label{eq:direct-symmetric-curl-material-zero}
		\nrm{\mathcal R\curl G}_{N+\alpha,J}
		\lesssim A_G\Lk{N,K_G,\lambda_G,\tilde\lambda_G}{s}
		\MM{N,K_G,\lambda_G,\tilde\lambda_G}.
	\end{equation}
	For $\mathcal B\curl G$, the same estimates hold in the ranges of \eqref{eq:direct-symmetric-curl-high-material} and \eqref{eq:direct-symmetric-curl-material-zero}, respectively:
	\begin{align}
		\nrm{D_{t,\upsilon}^{M}\mathcal B\curl G}_{N+\alpha,J}
		&\lesssim A_G\Lk{N,(K_G-M)_+,\lambda,\tilde\lambda}{\max\{2,s\}}
		\MM{N,(K_G-M)_+,\lambda,\tilde\lambda}
		\MM{M,N_t,\mu,\tilde\mu},	\label{eq:direct-vector-curl-high-material}\\
		\nrm{\mathcal B\curl G}_{N+\alpha,J}
		&\lesssim A_G\Lk{N,K_G,\lambda_G,\tilde\lambda_G}{s}
		\MM{N,K_G,\lambda_G,\tilde\lambda_G}.\label{eq:direct-vector-curl-material-zero}
	\end{align}

	For $r\in\mathbb N_0$ with $M+r\leq K^*+1$, the same estimates hold with $\alpha$ replaced
	by $\alpha_r$, in the same derivative ranges. At material order $M$, the proof uses the H{\"o}lder assumptions
	for $\upsilon$ and $G$ only with H{\"o}lder exponents $\alpha_r,\ldots,\alpha_{r+M}$.
	The implicit constants depend only on $K^*,K^\sharp,\alpha$.
	They are independent of $J,\upsilon,G$, the amplitudes, and the spatial and material scales.
\end{lem}

\begin{proof}
	Let $e_1,e_2,e_3$ denote the standard basis of $\mathbb R^3$. For every
	$i,j\in\{1,2,3\}$,
	$$
	\Div(G\otimes e_j)=\partial_jG,
	\qquad
	W[G\otimes e_j]_i
	=-\epsilon_{ikl}\partial_k\partial_j\Delta^{-1}G_{l}.
	$$
	The definition \eqref{def:direct-symmetric-curl-inverse} therefore gives
	$$
	(\mathcal R\curl G)_{ij}
	=-W[G\otimes e_j]_i-W[G\otimes e_i]_j.
	$$
	Since $e_j$ is constant,
	$
	D_{t,\upsilon}^{M}(G\otimes e_j)
	=(D_{t,\upsilon}^{M}G)\otimes e_j.
	$
	For $\mathcal B\curl$, define the skew-symmetric matrix field $T$ by
	$T_{ij}:=\epsilon_{ijk}G_k$. Then
	$$
	\Div T=\curl G,\qquad W[T]=\mathcal B\curl G,\qquad
	(D_{t,\upsilon}^{M}T)_{ij}=\epsilon_{ijk}D_{t,\upsilon}^{M}G_k.
	$$
	These matrix fields satisfy \eqref{eq:periodic-Hodge-T-Holder-assumption},
	with $K_T=K_G$, the same $M_*$, $A_T$ a fixed dimensional multiple
	of $A_G$, and
	$
	\lambda_T=\lambda_G,\
	\tilde\lambda_T=\tilde\lambda_G,\
	\mu_T=\mu_G,\
	\tilde\mu_T=\tilde\mu_G.
	$
	Proposition~\ref{prop:periodic-Hodge-material-estimates}, applied with $r^*=K^*+1$, gives
	all four estimates in the stated ranges.
\end{proof}

\subsection{Material derivative estimates for the linearized Euler equations}
\label{subsec:linearized-equation-material-estimates}

We consider the linearized equation for the incompressible
correction in the Newton step, with source $\curl G$ and zero initial data. Neither $U$
nor $G$ is required to be divergence-free. By
Lemma~\ref{lem:direct-vector-curl-inverse}, replacing $G$ by
$\mathcal B\curl G$ preserves the source and gives a zero-mean
divergence-free field.
\begin{equation}\label{eq:linearized-equation 1}
	\left\{
	\begin{aligned}
		&\partial_t\omega+\Div(U\otimes\omega+\omega\otimes U)+\nabla\pi=\curl G,\\
		&\Div\omega=0,\\
		&\omega|_{t=t_0}=0.
	\end{aligned}
	\right.
\end{equation}
Here $D_{t,U}:=\partial_t+U\cdot\nabla$. Integrating the equation in space
and using the zero initial data gives $\fint_{\T^3}\omega\,\rd x=0$.
We normalize $\pi$ to have zero spatial mean and introduce the vector
potential $\Psi=\mathcal B\omega$.

Proposition~\ref{prop:linearized-equation} gives existence, uniqueness,
and the H{\"o}lder estimates for the solution and vector potential at material order zero. For smooth data,
Remark~\ref{rem:smooth-data-smooth-solution}
implies that the solution and pressure are smooth. For higher orders,
we use Lemma~\ref{lem:direct-symmetric-curl-material-estimate} to estimate
$\mathcal B\curl G$. The next proposition estimates $\omega$, $\Psi$,
and $\pi$ in terms of the derivatives of $G$.

\begin{prop} 
	\label{prop:linearized-high-material-derivatives} 
	Let $I\subset\mathbb R$ be an interval. Fix $t_0\in I$ and
	$\tau>0$ such that 
	$ 
	J:=J_\tau=[t_0-\tau,t_0+\tau]\subset I. 
	$ 
	Let $U$ and $G$ be smooth vector fields on $J\times\T^3$.
	Fix $K^*,K^\sharp\in\mathbb N$ and $K_G,M_*,N_t\in\mathbb N_0$
	with $K^\sharp\geq K^*+2$ and $K_G+1\leq K^*$.
	If $M_*\geq K^*$, assume also $N_t\leq K^*$.
	Let $A_U>0$ and choose scales
	$ 
	1\leq\lambda_U\leq\tilde\lambda_U 
	$ 
	and 
	$ 
	1\leq\mu_U\leq\tilde\mu_U 
	$ 
	for $U$.
	For $0<c_0\leq1$, assume the short-time condition 
	\begin{equation}\label{eq:linearized-short-time-assumption} 
	\tau A_U\La{1,K^*,\lambda_U,\tilde\lambda_U} 
	\MM{1,K^*,\lambda_U,\tilde\lambda_U}\leq c_0. 
	\end{equation} 
	Assume the following $C^N$ estimates for the material derivatives of $\nabla U$ 
	\begin{equation}\label{eq:linearized-U-integer-assumption} 
	\begin{aligned} 
		\sup_{t\in J}\|D_{t,U}^{M}\nabla U(t,\cdot)\|_{N} 
		\leq A_U 
		\MM{N+1,K^*-M,\lambda_U,\tilde\lambda_U} 
		\MM{M,N_t,\mu_U,\tilde\mu_U}, \
		0\leq M\leq K^*,\ M+N\leq K^\sharp-1. 
	\end{aligned} 
	\end{equation} 
	For every $0\leq r\leq K^*+1$, assume also the H{\"o}lder estimates 
	\begin{equation}\label{eq:linearized-U-Holder-assumption} 
	\begin{aligned} 
		\sup_{t\in J}\|D_{t,U}^{M}\nabla U(t,\cdot)\|_{N+\alpha_r} 
		&\leq A_U 
		\Lab{N+1,K^*-M,\lambda_U,\tilde\lambda_U}{\alpha_r}\MM{N+1,K^*-M,\lambda_U,\tilde\lambda_U} \\ 
		&\quad\cdot 
		\MM{M,N_t,\mu_U,\tilde\mu_U}, 
		\hspace{110pt}0\leq M\leq K^*,\ \ M+N\leq K^\sharp-2. 
	\end{aligned} 
	\end{equation} 
	Fix $s\in\mathbb N$ and $A_G>0$, and choose scales
	$ 
	1\leq\lambda_G\leq\tilde\lambda_G 
	$ 
	and 
	$ 
	1\leq\mu_G\leq\tilde\mu_G 
	$ 
	for $G$. Assume at the same exponents 
	\begin{equation}\label{eq:linearized-G-Holder-assumption} 
	\begin{aligned} 
		\sup_{t\in J}\|D_{t,U}^{M}G(t,\cdot)\|_{N+\alpha_r} 
		&\leq A_G 
		\Lkb{N,(K_G-M)_+,\lambda_G,\tilde\lambda_G}{s}{\alpha_r}\\ 
		&\quad\cdot 
		\MM{N,(K_G-M)_+,\lambda_G,\tilde\lambda_G} 
		\MM{M,N_t,\mu_G,\tilde\mu_G}, 
		&0\leq M\leq M_*,\ \ M+N\leq K^\sharp-1. 
	\end{aligned} 
	\end{equation} 
	
	Define the combined scales by
	$$ 
	\lambda:=\max\{\lambda_U,\lambda_G\},\qquad 
	\tilde\lambda:=\max\{\tilde\lambda_U,\tilde\lambda_G\},\qquad 
	\mu:=\max\{\mu_U,\mu_G,A_U\lambda_U\},\qquad 
	\tilde\mu:=\max\{\tilde\mu_U,\tilde\mu_G,A_U\tilde\lambda_U\}. 
	$$ 
	Let $(\omega,\pi)$ solve \eqref{eq:linearized-equation 1}, with $\pi$
	normalized to have zero spatial mean, and set $\Psi:=\mathcal B\omega$.
	Then, for every $M,N\in\mathbb N_0$ with 
	$ 
	1\leq M\leq\min\{K^*,M_*\}+1,\ M+N\leq K^\sharp-1, 
	$ 
	one has 
	\begin{equation} 
		\begin{aligned} 
			&\sup_{t\in J}\left( 
			\|D_{t,U}^{M}\omega(t,\cdot)\|_{N+\alpha} 
			+\|D_{t,U}^{M}\Psi(t,\cdot)\|_{N+1+\alpha} 
			\right)\\ 
			&\quad\lesssim 
			\Lk{N+1,(K_G-M+1)_+,\lambda,\tilde\lambda}{\max\{2,s\}} 
			\MM{N+1,(K_G-M+1)_+,\lambda,\tilde\lambda}\\ 
			&\qquad\cdot A_G\Bigl( 
			\tau\MM{M,N_t,\mu,\tilde\mu} 
			+\MM{M-1,N_t,\mu,\tilde\mu} 
			\Bigr). 
		\end{aligned}\label{eq:linearized-spatial-material-scale-conclusion-Calpha} 
	\end{equation} 
	In the same range, the potential satisfies 
	\begin{equation} 
		\begin{aligned} 
			\sup_{t\in J}\|D_{t,U}^{M}\Psi(t,\cdot)\|_{N+\alpha}
			&\lesssim 
			\Lk{N,(K_G-M+1)_+,\lambda,\tilde\lambda}{\max\{2,s\}} 
			\MM{N,(K_G-M+1)_+,\lambda,\tilde\lambda}\\
			&\quad\cdot A_G\Bigl( 
			\tau\MM{M,N_t,\mu,\tilde\mu} 
			+\MM{M-1,N_t,\mu,\tilde\mu} 
			\Bigr). 
		\end{aligned}\label{eq:potential-spatial-material-scale-conclusion} 
	\end{equation} 
	For $1\leq M\leq\min\{K^*,M_*+1\}$ and $M+N\leq K^\sharp-1$,
	the pressure satisfies 
	\begin{equation} 
		\begin{aligned} 
			\sup_{t\in J}\|D_{t,U}^{M}\pi(t,\cdot)\|_{N+\alpha}
			&\lesssim 
			A_U\Lk{N+1,(K_G-M)_+,\lambda,\tilde\lambda}{\max\{2,s\}} 
			\MM{N+1,(K_G-M)_+,\lambda,\tilde\lambda}\\
			&\quad\cdot 
			A_G\Bigl( 
			\tau\MM{M,N_t,\mu,\tilde\mu} 
			+\MM{M-1,N_t,\mu,\tilde\mu} 
			\Bigr). 
		\end{aligned}\label{eq:pressure-spatial-material-scale-conclusion} 
	\end{equation} 
	
	At material order $M=0$, for every $0\leq N\leq K^\sharp-2$, one has 
	$$ 
	\begin{aligned} 
		\sup_{t\in J}\left( 
		\|\omega(t,\cdot)\|_{N+\alpha} 
		+\|\Psi(t,\cdot)\|_{N+1+\alpha} 
		\right)\lesssim\tau A_G 
		\Lk{N+1,K_G,\lambda,\tilde\lambda}{s} 
		\MM{N+1,K_G,\lambda,\tilde\lambda}. 
	\end{aligned} 
	$$ 
	For every
	$0\leq N\leq K^\sharp-1$, the potential also satisfies 
	$$ 
	\sup_{t\in J}\|\Psi(t,\cdot)\|_{N+\alpha} 
	\lesssim 
	\tau A_G 
	\Lk{N,K_G,\lambda,\tilde\lambda}{s} 
	\MM{N,K_G,\lambda,\tilde\lambda}. 
	$$ 
	For every $0\leq N\leq K^\sharp-1$, the pressure satisfies
	$$ 
	\sup_{t\in J}\|\pi(t,\cdot)\|_{N+\alpha} 
	\lesssim 
	\tau A_UA_G 
	\Lk{N+1,K_G,\lambda,\tilde\lambda}{\max\{2,s\}} 
	\MM{N+1,K_G,\lambda,\tilde\lambda}. 
	$$ 
	The implicit constants depend only on
	$K^*,K^\sharp,\alpha$. 
\end{prop}
\begin{proof}
	Throughout this proof, write $K_{G,M}:=(K_G-M)_+$ for $M\geq0$.
	We first estimate $\mathcal B\curl G$, which appears in the equation for the vector potential.
	For $0\leq M\leq\min\{K^*,M_*\}$, $M+r\leq K^*$, and
	$M+N\leq K^\sharp-1$, we have
	\begin{equation}
		\label{eq:linearized-projected-source-material}
		\begin{aligned}
			\sup_{t\in J}\|D_{t,U}^{M}\mathcal B\curl G(t,\cdot)\|_{N+\alpha_r}
			\lesssim A_G\Lkb{N,K_{G,M},\lambda,\tilde\lambda}{\max\{2,s\}}{\alpha_r}
			\MM{N,K_{G,M},\lambda,\tilde\lambda}\MM{M,N_t,\mu,\tilde\mu}.
		\end{aligned}
	\end{equation}
	For $M\geq1$, apply \eqref{eq:direct-vector-curl-high-material} at exponent $\alpha_r$,
	as permitted by Lemma~\ref{lem:direct-symmetric-curl-material-estimate}.
	At $M=0$, \eqref{eq:direct-vector-curl-Schauder}
	gives the sharper H{\"o}lder factor
	with power $s$ from \eqref{eq:linearized-G-Holder-assumption} at all
	exponents $\alpha_0,\ldots,\alpha_{K^*+1}$.
	For the rest of the proof, we replace $G$ by $\mathcal B\curl G$ and use the same notation
	$G$. By \eqref{eq:direct-vector-curl-properties}, this replacement leaves
	the source $\curl G$ unchanged and makes $G$ divergence-free with zero spatial mean.
	
	We begin with the estimates at material order zero, which provide the base cases for the inductions below.
	Fix $0\leq r\leq K^*+1$ and $0\leq N\leq K^\sharp-2$. For $t\in J$, let
	$I(t,t_0)$ denote the interval with endpoints $t$ and $t_0$.
	By \eqref{eq:direct-vector-curl-material-zero}, \eqref{eq:linearized-short-time-assumption}, and
	\eqref{eq:linearized-U-Holder-assumption}--%
	\eqref{eq:linearized-G-Holder-assumption}, we have
	\begin{equation}\label{eq:linearized-zero-order-integral-estimates}
	\begin{aligned}
		\int_{I(t,t_0)}\|\nabla U(\sigma,\cdot)\|_{\alpha_r}\,\rd \sigma
		&\lesssim
		\tau A_U\Lab{1,K^*,\lambda_U,\tilde\lambda_U}{\alpha_r}
		\MM{1,K^*,\lambda_U,\tilde\lambda_U}
		\lesssim c_0,\\
		\int_{I(t,t_0)}\|\curl G(\sigma,\cdot)\|_{N+\alpha_r}\,\rd \sigma
		&\lesssim
		\tau A_G\Lkb{N+1,K_G,\lambda,\tilde\lambda}{s}{\alpha_r}
		\MM{N+1,K_G,\lambda,\tilde\lambda},\\
		\int_{I(t,t_0)}\|G(\sigma,\cdot)\|_{N+\alpha_r}\,\rd \sigma
		&\lesssim
		\tau A_G\Lkb{N,K_G,\lambda,\tilde\lambda}{s}{\alpha_r}
		\MM{N,K_G,\lambda,\tilde\lambda}.
	\end{aligned}
	\end{equation}
	Apply \eqref{eq:holder-low-estimate} for $N=0$ and
	\eqref{eq:holder-estimate} for $N\geq1$, at exponent $\alpha_r$ and in
	both time directions from $t_0$, with zero initial data.
	For $N\geq1$, the last integral term in \eqref{eq:holder-estimate} is controlled using
	\eqref{eq:linearized-short-time-assumption},
	\eqref{eq:linearized-U-Holder-assumption} at material order zero, and
	\eqref{eq:linearized-zero-order-integral-estimates}. Together with $K_G+1\leq K^*$
	and $\alpha_r<1$, these imply
	\begin{equation}\label{eq:linearized-zero-order-integral-product-estimate}
		\begin{aligned}
			\left(\int_{I(t,t_0)}\|\nabla U(\sigma,\cdot)\|_{N+\alpha_r}\,\rd\sigma\right)
			\left(\int_{I(t,t_0)}\|\curl G(\sigma,\cdot)\|_{\alpha_r}\,\rd\sigma\right)\lesssim c_0\tau A_G
			\Lkb{N+1,K_G,\lambda,\tilde\lambda}{s}{\alpha_r}
			\MM{N+1,K_G,\lambda,\tilde\lambda}.
		\end{aligned}
	\end{equation}
	Substituting \eqref{eq:linearized-zero-order-integral-estimates}--\eqref{eq:linearized-zero-order-integral-product-estimate} into
	\eqref{eq:holder-low-estimate}--\eqref{eq:holder-estimate}
	and using $c_0\leq1$, we obtain
	\begin{equation}
	\label{eq:linearized-solution-material-zero}
	\begin{aligned}
		\sup_{t\in J}\|\omega(t,\cdot)\|_{N+\alpha_r}
		\lesssim\tau A_G
		\Lkb{N+1,K_G,\lambda,\tilde\lambda}{s}{\alpha_r}
		\MM{N+1,K_G,\lambda,\tilde\lambda}.
	\end{aligned}
	\end{equation}
	Since $\Div\Psi=0$, $\curl\Psi=\omega$, and $\Psi$ has zero spatial mean,
	estimate~\eqref{eq:linearized-solution-material-zero}, together with
	Lemma~\ref{lem:periodic-Hodge-estimates} with $L=N+1$ and exponent $\alpha_r$, gives
	$$
	\begin{aligned}
		\sup_{t\in J}\|\Psi(t,\cdot)\|_{N+1+\alpha_r}
		&\lesssim\tau A_G
		\Lkb{N+1,K_G,\lambda,\tilde\lambda}{s}{\alpha_r}
		\MM{N+1,K_G,\lambda,\tilde\lambda}.
	\end{aligned}
	$$
	
	We next estimate $\Psi$ with one fewer spatial derivative. Let
	$1\leq N\leq K^\sharp-1$. Since
	$\Psi=\mathcal B\omega$,
	\begin{align*}
	\|\Psi(t,\cdot)\|_{N+\alpha_r}
	&\lesssim
	\|\omega(t,\cdot)\|_{N-1+\alpha_r},
	\\
	\|\curl G(t,\cdot)\|_{N-1+\alpha_r}
	&\lesssim
	\|G(t,\cdot)\|_{N+\alpha_r}.
	\end{align*}
	Using \eqref{eq:linearized-solution-material-zero} at spatial order $N-1$, we obtain
	$$
	\sup_{t\in J}\|\Psi(t,\cdot)\|_{N+\alpha_r}
	\lesssim\sup_{t\in J}\|\omega(t,\cdot)\|_{N-1+\alpha_r}
	\lesssim
	\tau A_G
	\Lkb{N,K_G,\lambda,\tilde\lambda}{s}{\alpha_r}
	\MM{N,K_G,\lambda,\tilde\lambda}.
	$$
	When $N=0$, apply \eqref{eq:besov-estimate} to $\Psi$
	with zero initial data. Since $\mathcal B\curl G=G$,
	the estimates in \eqref{eq:linearized-zero-order-integral-estimates}
	give
	$$
	\sup_{t\in J}\|\Psi(t,\cdot)\|_{\alpha_r}
	\lesssim
	\tau A_G\Lkb{0,K_G,\lambda,\tilde\lambda}{s}{\alpha_r}.
	$$
	Combining the cases $N\geq1$ and $N=0$, we have for every $0\leq r\leq K^*+1$ and
	$0\leq N\leq K^\sharp-1$,
	\begin{equation}
	\label{eq:linearized-potential-material-zero}
	\sup_{t\in J}\|\Psi(t,\cdot)\|_{N+\alpha_r}
	\lesssim
	\tau A_G
	\Lkb{N,K_G,\lambda,\tilde\lambda}{s}{\alpha_r}
	\MM{N,K_G,\lambda,\tilde\lambda}.
	\end{equation}
	Taking $r=0$ proves the estimates with no material derivatives for $\omega$ and
	$\Psi$.
	
	We now prove the estimate for $\omega$ in $C^{N+\alpha_r}$ and $\Psi$ in $C^{N+1+\alpha_r}$ at positive material orders. Set
	$$
	\omega^{[M]}:=D_{t,U}^{M}\omega,
	\qquad
	\Psi^{[M]}:=D_{t,U}^{M}\Psi
	$$
	for $M\geq0$. We argue by induction on $M$, simultaneously at the
	exponents $\alpha_r$. For $1\leq M\leq M_*+1$ and $r,N\in\mathbb N_0$
	with $M+r\leq K^*+1$ and $M+N\leq K^\sharp-1$, we prove
	\begin{equation}\label{eq:linearized-omega-potential-material-induction}
	\begin{aligned}
		&\sup_{t\in J}\left(
		\|\omega^{[M]}(t,\cdot)\|_{N+\alpha_r}
		+\|\Psi^{[M]}(t,\cdot)\|_{N+1+\alpha_r}
		\right)\\
		&\quad\lesssim
		\Lkb{N+1,K_{G,M-1},\lambda,\tilde\lambda}{\max\{2,s\}}{\alpha_r}
		\MM{N+1,K_{G,M-1},\lambda,\tilde\lambda}
		A_G\Bigl(
		\tau\MM{M,N_t,\mu,\tilde\mu}
		+\MM{M-1,N_t,\mu,\tilde\mu}
		\Bigr).
	\end{aligned}
	\end{equation}
	Fix $M$ and assume the estimates at all lower material orders in the
	stated ranges.
	
	For every smooth vector field $h$, \eqref{eq:first-order-material-spatial-commutator} gives
	$$
	[D_{t,U},\Div]h=-\partial_iU_j\,\partial_jh_i, \qquad \bigl([D_{t,U},\curl]h\bigr)_k
	=-\epsilon_{kij}\partial_iU_{l}\,\partial_{l}h_j.
	$$
	If $\Div h=0$, then
	\begin{equation}
		\label{eq:linearized-divergence-free-commutator}
	\begin{aligned}
		-[D_{t,U},\Div]h
		=\partial_iU_j\,\partial_jh_i
		=\partial_i\left(
		(\partial_jU_i)h_j-(\Div U)h_i
		\right).
	\end{aligned}
	\end{equation}
	If $\curl h=0$, then $\partial_i h_{l}=\partial_{l} h_i$, and hence
	\begin{equation}
		\label{eq:linearized-curl-free-commutator}
	\begin{aligned}
		\bigl([D_{t,U},\curl]h\bigr)_k
		=-\epsilon_{kij}\partial_iU_{l}\,\partial_{l}h_j
		=\epsilon_{kij}\partial_i\left(
		(\partial_jU_{l})h_{l}
		\right).
	\end{aligned}
	\end{equation}
	
	At material order one, expanding
	\eqref{eq:linearized-equation 1}, we have
	$$
	\omega^{[1]}
	=
	\curl G-\omega\cdot\nabla U-\omega\,\Div U-\nabla\pi.
	$$
	Taking the curl eliminates the pressure and leaves
	$$
	\curl \omega^{[1]}
	=
	\curl\left(
	\curl G-\omega\cdot\nabla U-\omega\,\Div U
	\right).
	$$
	The divergence-free constraint gives
	$$
	\Div \omega^{[1]}
	=
	\partial_iU_j\,\partial_j\omega_i
	=
	\partial_i\left(
	(\partial_jU_i)\omega_j-(\Div U)\omega_i
	\right).
	$$
	Moreover, since $\fint_{\mathbb T^3}\omega\,\rd x=0$,
	\eqref{eq:first-order-material-mean-identity} gives
	$$
	\fint_{\mathbb T^3}\omega^{[1]}\,\rd x
	=
	-\fint_{\mathbb T^3}\omega\,\Div U\,\rd x.
	$$
	
	For every $M\geq1$, we construct vector fields
	$H_M^{\Div}$ and $H_M^{\curl}$ such that
	$$
	\Div \omega^{[M]}=\Div H_M^{\Div},
	\qquad
	\curl \omega^{[M]}=\curl H_M^{\curl}.
	$$
	Each component of $H_M^{\Div}$ is a finite sum of fixed tensor
	contractions of products of the form
	$$
	\left(
	\prod_{\nu=1}^{\kappa}D_{t,U}^{M_\nu}\nabla U
	\right)\omega^{[j]},
	\qquad
	0\leq j\leq M-1,
	\qquad
	j+\sum_{\nu=1}^{\kappa}(M_\nu+1)=M.
	$$
	Each component of $H_M^{\curl}$ is a finite sum of fixed tensor contractions
	of two types of products. The first contains $\omega^{[j]}$
	as described for $H_M^{\Div}$. The second has the form
	$$
	\left(
	\prod_{\nu=1}^{\kappa}D_{t,U}^{M_\nu}\nabla U
	\right)\nabla D_{t,U}^{k}G,
	\qquad
	0\leq k\leq M-1,
	\qquad
	k+\sum_{\nu=1}^{\kappa}(M_\nu+1)=M-1.
	$$
	For products containing $\omega^{[j]}$, one has $\kappa\geq1$. For products containing $\nabla D_{t,U}^kG$,
	the empty coefficient product is one, and $\kappa=0$ gives the direct term
	with $k=M-1$.
	
	At $M=1$, we take
	\begin{align*}
	\left(H_1^{\Div}\right)_i
	&:=
	(\partial_jU_i)\omega_j-(\Div U)\omega_i,
	\\
	H_1^{\curl}
	&:=
	\curl G-\omega\cdot\nabla U-\omega\,\Div U.
	\end{align*}
	Suppose that $H_M^{\Div}$ and $H_M^{\curl}$ have been constructed. Since
	$
	\Div(H_M^{\Div}-\omega^{[M]})=0,
	$
	we have
	$$
	\begin{aligned}
		\Div \omega^{[M+1]}
		=
		D_{t,U}\Div \omega^{[M]}-[D_{t,U},\Div]\omega^{[M]}
		=
		\Div(D_{t,U}H_M^{\Div})
		+[D_{t,U},\Div](H_M^{\Div}-\omega^{[M]}).
	\end{aligned}
	$$
	Applying \eqref{eq:linearized-divergence-free-commutator} to the divergence-free
	field $H_M^{\Div}-\omega^{[M]}$, we obtain
	$$
	\begin{aligned}
		\relax[D_{t,U},\Div](H_M^{\Div}-\omega^{[M]})
		=
		-\partial_i\left(
		(\partial_jU_i)(H_M^{\Div}-\omega^{[M]})_j
		-(\Div U)(H_M^{\Div}-\omega^{[M]})_i
		\right).
	\end{aligned}
	$$
	Thus $\Div \omega^{[M+1]}=\Div H_{M+1}^{\Div}$, where
	$$
		\left(H_{M+1}^{\Div}\right)_i
		:=
		D_{t,U}\left(H_M^{\Div}\right)_i
		-(\partial_jU_i)(H_M^{\Div}-\omega^{[M]})_j
		+(\Div U)(H_M^{\Div}-\omega^{[M]})_i.
	$$
	
	Similarly,
	$
	\curl(H_M^{\curl}-\omega^{[M]})=0,
	$
	and therefore
	$$
	\begin{aligned}
		\curl \omega^{[M+1]}
		=
		D_{t,U}\curl \omega^{[M]}-[D_{t,U},\curl]\omega^{[M]}=
		\curl(D_{t,U}H_M^{\curl})
		+[D_{t,U},\curl](H_M^{\curl}-\omega^{[M]}).
	\end{aligned}
	$$
	Applying \eqref{eq:linearized-curl-free-commutator} to the curl-free field
	$H_M^{\curl}-\omega^{[M]}$, we have
	$$
	\begin{aligned}
		\bigl([D_{t,U},\curl](H_M^{\curl}-\omega^{[M]})\bigr)_k
		=
		\epsilon_{kij}\partial_i\left(
		(\partial_jU_{l})(H_M^{\curl}-\omega^{[M]})_{l}
		\right).
	\end{aligned}
	$$
	Hence $\curl \omega^{[M+1]}=\curl H_{M+1}^{\curl}$, where
	$$
	\left(H_{M+1}^{\curl}\right)_j
	:=
	D_{t,U}\left(H_M^{\curl}\right)_j
	+(\partial_jU_{k})(H_M^{\curl}-\omega^{[M]})_{k}.
	$$
	
	For products containing $\omega^{[j]}$, applying $D_{t,U}$ or multiplying by $\nabla U$
	raises $j+\sum_{\nu=1}^{\kappa}(M_\nu+1)$ by one.
	The new terms containing $\omega^{[M]}$ have $\kappa=1$, $M_1=0$, and $j=M$,
	so the same sum equals $M+1$.
	Thus the recursions preserve the stated product forms and order identities.
	For products containing $\nabla D_{t,U}^kG$,
	\eqref{eq:first-order-material-spatial-commutator} gives, componentwise,
	$$
	D_{t,U}\bigl(\partial_iD_{t,U}^kG\bigr)
	=\partial_iD_{t,U}^{k+1}G
	-(\partial_iU_l)\partial_lD_{t,U}^kG.
	$$
	Both terms increase $k+\sum_{\nu=1}^{\kappa}(M_\nu+1)$ by one,
	preserving the stated order identity at the next material order.
	
	The mean identity \eqref{eq:first-order-material-mean-identity} gives
	\begin{equation}\label{eq:linearized-omega-mean-recursion}
	\fint_{\mathbb T^3}\omega^{[M+1]}\,\rd x
	=
	\frac{\rd}{\rd t}\fint_{\mathbb T^3}\omega^{[M]}\,\rd x
	-
	\fint_{\mathbb T^3}\omega^{[M]}\,\Div U\,\rd x.
	\end{equation}
	Since $\omega^{[0]}=\omega$ has zero spatial mean, \eqref{eq:linearized-omega-mean-recursion} with $M=0$ reads
	$$
	\fint_{\mathbb T^3}\omega^{[1]}\,\rd x
	=-\fint_{\mathbb T^3}\omega\,\Div U\,\rd x.
	$$
	This has $\kappa=1$, $M_1=0$, and $j=0$, since $\Div U=\tr(\nabla U)$.
	We claim that, for every $M\geq1$, $\fint_{\mathbb T^3}\omega^{[M]}\,\rd x$
	is a finite linear combination of fixed tensor contractions of
	$$
	\fint_{\mathbb T^3}
	\left(
	\prod_{\nu=1}^{\kappa}D_{t,U}^{M_\nu}\nabla U
	\right)\omega^{[j]}\,\rd x,
	\qquad
	\kappa\geq1,
	\qquad
	0\leq j\leq M-1,
	\qquad
	j+\sum_{\nu=1}^{\kappa}(M_\nu+1)=M.
	$$
	Here $M_\nu\in\mathbb N_0$. The computation above proves the claim for $M=1$.
	Assume it holds at order $M\geq1$. For each contraction in that expansion,
	\eqref{eq:first-order-material-mean-identity} gives
	\begin{equation}\label{eq:linearized-omega-mean-product-derivative}
	\begin{aligned}
		\frac{\rd}{\rd t}\fint_{\mathbb T^3}
		\left(\prod_{\nu=1}^{\kappa}D_{t,U}^{M_\nu}\nabla U\right)\omega^{[j]}\,\rd x&=\sum_{i=1}^{\kappa}\fint_{\mathbb T^3}
		\left(D_{t,U}^{M_i+1}\nabla U
		\prod_{\substack{\nu=1\\\nu\neq i}}^{\kappa}D_{t,U}^{M_\nu}\nabla U\right)
		\omega^{[j]}\,\rd x\\
		&\quad+\fint_{\mathbb T^3}
		\left(\prod_{\nu=1}^{\kappa}D_{t,U}^{M_\nu}\nabla U\right)\omega^{[j+1]}\,\rd x+\fint_{\mathbb T^3}
		\left(\prod_{\nu=1}^{\kappa}D_{t,U}^{M_\nu}\nabla U\right)
		\omega^{[j]}\Div U\,\rd x.
	\end{aligned}
	\end{equation}
	In \eqref{eq:linearized-omega-mean-product-derivative}, the first sum increases one $M_i$ by one, and the second term increases $j$ by one.
	The last term adds a factor $\nabla U$ of material order zero through its trace $\Div U$.
	Each term therefore has total order $M+1$ and $j\leq M$.
	The last term in \eqref{eq:linearized-omega-mean-recursion},
	$
	-\fint_{\mathbb T^3}\omega^{[M]}\Div U\,\rd x,
	$
	has $\kappa=1$, $M_1=0$, and $j=M$.
	Thus the expansion and its order identity hold at order $M+1$. This argument applies to any smooth periodic field with zero spatial mean at each time.
	
	For $r,N$ in the stated ranges, we estimate separately the products containing $\omega^{[j]}$ and those containing $\nabla D_{t,U}^kG$.
	Consider first a product containing $\omega^{[j]}$ in
	$H_M^{\Div}$, $H_M^{\curl}$, or the mean expansion.
	For $j\geq1$, the induction hypothesis \eqref{eq:linearized-omega-potential-material-induction} applies with H{\"o}lder exponents $\alpha_r$ and $\alpha_{r+1}$
	since $j+r+1\leq M+r\leq K^*+1$. For $j=0$, use
	\eqref{eq:linearized-solution-material-zero}.
	Together with \eqref{eq:linearized-U-integer-assumption}--\eqref{eq:linearized-U-Holder-assumption}
	and \eqref{eq:holder-product-lower-exponent} with $\eta=\alpha_{r+1}$,
	these estimates give
	\begin{equation}\label{eq:linearized-omega-correction-product-estimate}
	\begin{aligned}
		&\sup_{t\in J}
		\left\|
		\left(
		\prod_{\nu=1}^{\kappa}D_{t,U}^{M_\nu}\nabla U
		\right)\omega^{[j]}
		\right\|_{N+\alpha_r}\\
		&\qquad\lesssim
		\sum_{N_0+\sum_{\nu=1}^{\kappa}N_\nu=N}\Biggl[
		\nrm{\omega^{[j]}}_{N_0+\alpha_r,J}
		\prod_{\nu=1}^{\kappa}\nrm{D_{t,U}^{M_\nu}\nabla U}_{N_\nu,J}\\
		&\hspace{80pt}+
		\nrm{\omega^{[j]}}_{N_0+\alpha_{r+1},J}
		\sum_{i=1}^{\kappa}\nrm{D_{t,U}^{M_i}\nabla U}_{N_i+\alpha_r,J}
		\prod_{\substack{\nu=1\\\nu\neq i}}^{\kappa}
		\nrm{D_{t,U}^{M_\nu}\nabla U}_{N_\nu,J}\Biggr]\\
		&\qquad\lesssim
		A_G\Lkb{N+1,K_{G,M-1},\lambda,\tilde\lambda}{\max\{2,s\}}{\alpha_r}
		\sum_{N_0+\sum_{\nu=1}^{\kappa}N_\nu=N}
		\MM{N_0+1,K_{G,M-1},\lambda,\tilde\lambda}
		\prod_{\nu=1}^{\kappa}\MM{N_\nu,K_{G,M-1},\lambda,\tilde\lambda}\\
		&\qquad\quad\cdot
		\sum_{h=0}^{\min\{1,j\}}\tau^{1-h}\MM{j-h,N_t,\mu,\tilde\mu}
		\prod_{\nu=1}^{\kappa}\MM{M_\nu+1,N_t,\mu,\tilde\mu}\\
		&\qquad\lesssim
		\Lkb{N+1,K_{G,M-1},\lambda,\tilde\lambda}{\max\{2,s\}}{\alpha_r}
		\MM{N+1,K_{G,M-1},\lambda,\tilde\lambda}
		A_G\Bigl(
		\tau\MM{M,N_t,\mu,\tilde\mu}
		+\MM{M-1,N_t,\mu,\tilde\mu}
		\Bigr).
	\end{aligned}
	\end{equation}
	
	For a product containing a derivative of $G$,
	$$
	\left(
	\prod_{\nu=1}^{\kappa}D_{t,U}^{M_\nu}\nabla U
	\right)\nabla D_{t,U}^{k}G,
	\qquad
	0\leq k\leq M-1,
	\qquad
	k+\sum_{\nu=1}^{\kappa}(M_\nu+1)=M-1,
	$$
	use \eqref{eq:direct-vector-curl-material-zero} when $k=0$, and
	\eqref{eq:linearized-projected-source-material} for $k\geq1$, with one
	additional spatial derivative. For $\kappa=0$, these estimates apply directly with $k=M-1$.
	For $\kappa\geq1$, use
	\eqref{eq:holder-product-lower-exponent} with $\eta=\alpha_{r+1}$.
	The estimates at exponent $\alpha_{r+1}$ apply since
	$
	k+r+1\leq M+r-1\leq K^*.
	$
	We obtain
	\begin{equation}\label{eq:linearized-gradient-source-product-estimate}
	\begin{aligned}
		&\sup_{t\in J}
		\left\|
		\left(
		\prod_{\nu=1}^{\kappa}D_{t,U}^{M_\nu}\nabla U
		\right)\nabla D_{t,U}^{k}G
		\right\|_{N+\alpha_r}\\
		&\lesssim
		\sum_{N_0+\sum_{\nu=1}^{\kappa}N_\nu=N}\Biggl[
		\nrm{\nabla D_{t,U}^{k}G}_{N_0+\alpha_r,J}
		\prod_{\nu=1}^{\kappa}\nrm{D_{t,U}^{M_\nu}\nabla U}_{N_\nu,J}\\
		&\hspace{60pt}+
		\nrm{\nabla D_{t,U}^{k}G}_{N_0+\alpha_{r+1},J}
		\sum_{i=1}^{\kappa}\nrm{D_{t,U}^{M_i}\nabla U}_{N_i+\alpha_r,J}
		\prod_{\substack{\nu=1\\\nu\neq i}}^{\kappa}
		\nrm{D_{t,U}^{M_\nu}\nabla U}_{N_\nu,J}\Biggr]\\
		&\lesssim
		A_G\Lkb{N+1,K_{G,M-1},\lambda,\tilde\lambda}{\max\{2,s\}}{\alpha_r}
		\sum_{N_0+\sum_{\nu=1}^{\kappa}N_\nu=N}
		\MM{N_0+1,K_{G,M-1},\lambda,\tilde\lambda}
		\prod_{\nu=1}^{\kappa}\MM{N_\nu,K_{G,M-1},\lambda,\tilde\lambda}\\
		&\quad\cdot\MM{k,N_t,\mu,\tilde\mu}
		\prod_{\nu=1}^{\kappa}\MM{M_\nu+1,N_t,\mu,\tilde\mu}\\
		&\lesssim
		A_G\Lkb{N+1,K_{G,M-1},\lambda,\tilde\lambda}{\max\{2,s\}}{\alpha_r}
		\MM{N+1,K_{G,M-1},\lambda,\tilde\lambda}
		\MM{M-1,N_t,\mu,\tilde\mu}.
	\end{aligned}
	\end{equation}
	For the second inequality in \eqref{eq:linearized-omega-correction-product-estimate}--\eqref{eq:linearized-gradient-source-product-estimate}, substitute
	\eqref{eq:linearized-U-integer-assumption}--\eqref{eq:linearized-U-Holder-assumption}
	and the estimates for $\omega^{[j]}$ and $\nabla D_{t,U}^kG$ cited above.
	For $M_\nu<K^*$, use
	$$
	\begin{aligned}
	\MM{N_\nu+1,K^*-M_\nu,\lambda_U,\tilde\lambda_U}
	&=\lambda_U\MM{N_\nu,K^*-M_\nu-1,\lambda_U,\tilde\lambda_U},\\
	\Lab{N_\nu+1,K^*-M_\nu,\lambda_U,\tilde\lambda_U}{\alpha_r}
	&=\Lab{N_\nu,K^*-M_\nu-1,\lambda_U,\tilde\lambda_U}{\alpha_r}.
	\end{aligned}
	$$
	Using $A_U\lambda_U\leq\mu$ in the material weight gives
	the factor $\MM{M_\nu+1,N_t,\mu,\tilde\mu}$ appearing in both estimates.
	Since $M_\nu\leq M-1$ and $K_G+1\leq K^*$,
	the remaining spatial thresholds satisfy
	\begin{equation}
		\label{eq:linearized-coefficient-threshold-comparison}
		(K^*-M_\nu-1)_+\geq(K^*-M)_+\geq(K_G-M+1)_+.
	\end{equation}
	If $M_\nu=K^*$, the order identity forces
	$M=K^*+1$, $\kappa=1$, and $j=0$. Hence $M_*\geq K^*$ and $N_t\leq K^*$.
	The spatial threshold is zero. Using $A_U\tilde\lambda_U\leq\tilde\mu$,
	we absorb the coefficient factor into the material weight at order $K^*+1$.

	The estimates for $\omega^{[j]}$ have spatial threshold $K_G$
	when $j=0$ and $(K_G-j+1)_+$ when $j\geq1$,
	while the estimate for $\nabla D_{t,U}^kG$ has threshold $(K_G-k)_+$.
	These thresholds and \eqref{eq:linearized-coefficient-threshold-comparison}
	give the common threshold $K_{G,M-1}$ used in both estimates.
	As in the proof of Proposition~\ref{prop:periodic-Hodge-material-estimates},
	$2\alpha_{r+1}=\alpha_r$ and
	$1+\max\{2,s\}/2\leq\max\{2,s\}$ bound the H{\"o}lder factors by
	$\Lkb{N+1,K_{G,M-1},\lambda,\tilde\lambda}{\max\{2,s\}}{\alpha_r}$.
	Estimate~\eqref{eq:appendix-derivative-weight-product}
	gives the final inequalities in
	\eqref{eq:linearized-omega-correction-product-estimate}--%
	\eqref{eq:linearized-gradient-source-product-estimate}.

	For each integrand $H$ in the expansion of the spatial mean of $\omega^{[M]}$, the estimate
	\eqref{eq:linearized-omega-correction-product-estimate} gives
	\begin{equation}\label{eq:linearized-omega-mean-product-estimate}
	\begin{aligned}
		&\sup_{t\in J}
		\left|
		\fint_{\mathbb T^3}H(t,x)\,\rd x
		\right|\\
		&\quad\lesssim
		\Lkb{1,K_{G,M-1},\lambda,\tilde\lambda}{\max\{2,s\}}{\alpha_r}
		\MM{1,K_{G,M-1},\lambda,\tilde\lambda}
		A_G\Bigl(
		\tau\MM{M,N_t,\mu,\tilde\mu}
		+\MM{M-1,N_t,\mu,\tilde\mu}
		\Bigr).
	\end{aligned}
	\end{equation}
	
	We use $\Div\omega^{[M]}=\Div H_M^{\Div}$ and $\curl\omega^{[M]}=\curl H_M^{\curl}$ to apply Lemma~\ref{lem:periodic-Hodge-estimates} to $\omega^{[M]}$.
	For $N\geq1$, \eqref{eq:linearized-omega-correction-product-estimate}--\eqref{eq:linearized-gradient-source-product-estimate} control $H_M^{\Div}$ and $H_M^{\curl}$ in $C^{N+\alpha_r}$. Together with \eqref{eq:linearized-omega-mean-product-estimate}, they give
	$$
	\begin{aligned}
		&\sup_{t\in J}\left(
		\|\Div \omega^{[M]}(t,\cdot)\|_{N-1+\alpha_r}
		+\|\curl \omega^{[M]}(t,\cdot)\|_{N-1+\alpha_r}
		+\left|\fint_{\mathbb T^3}\omega^{[M]}(t,x)\,\rd x\right|
		\right)\\
		&\qquad\lesssim
		\nrm{H_M^{\Div}}_{N+\alpha_r,J}
		+\nrm{H_M^{\curl}}_{N+\alpha_r,J}
		+\sup_{t\in J}\left|\fint_{\mathbb T^3}\omega^{[M]}(t,x)\,\rd x\right|\\
		&\qquad\lesssim
		\Lkb{N+1,K_{G,M-1},\lambda,\tilde\lambda}{\max\{2,s\}}{\alpha_r}
		\MM{N+1,K_{G,M-1},\lambda,\tilde\lambda}
		A_G\Bigl(
		\tau\MM{M,N_t,\mu,\tilde\mu}
		+\MM{M-1,N_t,\mu,\tilde\mu}
		\Bigr).
	\end{aligned}
	$$
	For $N=0$, the same identities, Remark~\ref{rem:app-holder-besov-equivalence},
	and \eqref{eq:app-multiplier-besov-estimate} give
	$$
	\begin{aligned}
		&\sup_{t\in J}\Bigl(
		\|\Div\omega^{[M]}(t,\cdot)\|_{B^{-1+\alpha_r}_{\infty,\infty}}
		+\|\curl\omega^{[M]}(t,\cdot)\|_{B^{-1+\alpha_r}_{\infty,\infty}}
		+\Bigl|\fint_{\mathbb T^3}\omega^{[M]}(t,x)\,\rd x\Bigr|\Bigr)\\
		&\quad=\sup_{t\in J}\Bigl(
		\|\Div H_M^{\Div}(t,\cdot)\|_{B^{-1+\alpha_r}_{\infty,\infty}}
		+\|\curl H_M^{\curl}(t,\cdot)\|_{B^{-1+\alpha_r}_{\infty,\infty}}
		+\Bigl|\fint_{\mathbb T^3}\omega^{[M]}(t,x)\,\rd x\Bigr|\Bigr)\\
		&\quad\lesssim
		\nrm{H_M^{\Div}}_{\alpha_r,J}
		+\nrm{H_M^{\curl}}_{\alpha_r,J}
		+\sup_{t\in J}\left|\fint_{\mathbb T^3}\omega^{[M]}(t,x)\,\rd x\right|\\
		&\quad\lesssim
		\Lkb{1,K_{G,M-1},\lambda,\tilde\lambda}{\max\{2,s\}}{\alpha_r}
		\MM{1,K_{G,M-1},\lambda,\tilde\lambda}
		A_G\Bigl(
		\tau\MM{M,N_t,\mu,\tilde\mu}
		+\MM{M-1,N_t,\mu,\tilde\mu}
		\Bigr),
	\end{aligned}
	$$
	where the last inequality uses \eqref{eq:linearized-omega-correction-product-estimate}--\eqref{eq:linearized-gradient-source-product-estimate} at $N=0$ and \eqref{eq:linearized-omega-mean-product-estimate}.
	Applying Lemma~\ref{lem:periodic-Hodge-estimates} to $\omega^{[M]}$,
	with $L=N$ for $N\geq1$ and its lowest-order estimate for $N=0$, gives
	$$
	\begin{aligned}
		&\sup_{t\in J}\|\omega^{[M]}(t,\cdot)\|_{N+\alpha_r}\\
		&\qquad\lesssim
		\nrm{H_M^{\Div}}_{N+\alpha_r,J}
		+\nrm{H_M^{\curl}}_{N+\alpha_r,J}
		+\sup_{t\in J}\left|\fint_{\mathbb T^3}\omega^{[M]}(t,x)\,\rd x\right|\\
		&\qquad\lesssim
		\Lkb{N+1,K_{G,M-1},\lambda,\tilde\lambda}{\max\{2,s\}}{\alpha_r}
		\MM{N+1,K_{G,M-1},\lambda,\tilde\lambda}
		A_G\Bigl(
		\tau\MM{M,N_t,\mu,\tilde\mu}
		+\MM{M-1,N_t,\mu,\tilde\mu}
		\Bigr).
	\end{aligned}
	$$
	
	To complete the induction step, we estimate the divergence, curl, and
	spatial mean of $\Psi^{[M]}$ using the estimate for $\omega^{[M]}$
	and the induction hypothesis, then apply
	Lemma~\ref{lem:periodic-Hodge-estimates}.
	Since $\Div\Psi=0$ and $\curl\Psi=\omega$,
	$
	\Div \Psi^{[0]}=0,
	\;
	\curl \Psi^{[0]}=\omega^{[0]}.
	$
	For every $M\geq0$,
	\begin{align*}
	\Div \Psi^{[M+1]}
	&=
	D_{t,U}\Div \Psi^{[M]}-[D_{t,U},\Div]\Psi^{[M]},
	\\
	\curl \Psi^{[M+1]}
	&=
	D_{t,U}\curl \Psi^{[M]}-[D_{t,U},\curl]\Psi^{[M]}.
	\end{align*}

	Since $\omega^{[M+1]}=D_{t,U}\omega^{[M]}$,
	$
	\curl \Psi^{[M+1]}-\omega^{[M+1]}
	=
	D_{t,U}(\curl \Psi^{[M]}-\omega^{[M]})-[D_{t,U},\curl]\Psi^{[M]}.
	$
	At the first material order, the commutator identities give
	$$
	\Div \Psi^{[1]}=\partial_iU_j\,\partial_j(\Psi^{[0]})_i,
	\qquad
	(\curl \Psi^{[1]}-\omega^{[1]})_k
	=\epsilon_{kij}\partial_iU_{l}\,\partial_{l}(\Psi^{[0]})_j.
	$$
	We claim that, for every $M\geq1$, each term in $\Div \Psi^{[M]}$ and
	$\curl \Psi^{[M]}-\omega^{[M]}$ is a fixed tensor contraction of a product of the form
	\begin{equation}
		\label{eq:linearized-potential-div-curl-products}
	\left(
	\prod_{\nu=1}^{\kappa}D_{t,U}^{M_\nu}\nabla U
	\right)\nabla \Psi^{[j]},
	\qquad
	0\leq j\leq M-1,
	\qquad
	j+\sum_{\nu=1}^{\kappa}(M_\nu+1)=M.
	\end{equation}
	Here $\kappa\geq1$.
	The first-order computation proves the claim for $M=1$. Suppose that it
	holds at order $M$. For a product occurring at that order,
	$$
	\begin{aligned}
		D_{t,U}\left[
		\left(\prod_{\nu=1}^{\kappa}D_{t,U}^{M_\nu}\nabla U\right)
		\nabla \Psi^{[j]}\right]&=
		\sum_{i=1}^{\kappa}
		\left(D_{t,U}^{M_i+1}\nabla U
		\prod_{\substack{\nu=1\\ \nu\neq i}}^{\kappa}
		D_{t,U}^{M_\nu}\nabla U\right)\nabla \Psi^{[j]}+\left(\prod_{\nu=1}^{\kappa}D_{t,U}^{M_\nu}\nabla U\right)
		\left(\nabla \Psi^{[j+1]}+[D_{t,U},\nabla]\Psi^{[j]}\right).
	\end{aligned}
	$$
	The gradient commutator adds one factor of $\nabla U$. Thus each term
	raises $j+\sum_{\nu=1}^{\kappa}(M_\nu+1)$ by one. The new divergence
	and curl commutators have the same product structure, which proves the claim.
	
	For the spatial mean, \eqref{eq:first-order-material-mean-identity} gives
	$$
	\fint_{\mathbb T^3}\Psi^{[M+1]}\,\rd x
	=
	\frac{\rd}{\rd t}\fint_{\mathbb T^3}\Psi^{[M]}\,\rd x
	-
	\fint_{\mathbb T^3}\Psi^{[M]}\,\Div U\,\rd x.
	$$
	Since $\Psi^{[0]}=\Psi$ has zero spatial mean,
	$
	\fint_{\mathbb T^3}\Psi^{[1]}\,\rd x
	=
	-\fint_{\mathbb T^3}\Psi\,\Div U\,\rd x.
	$
	By \eqref{eq:linearized-omega-mean-recursion}--\eqref{eq:linearized-omega-mean-product-derivative}, with $\omega$ replaced by $\Psi$,
	the mean of $\Psi^{[M]}$ is a finite linear combination
	of fixed tensor contractions of
	$$
	\fint_{\mathbb T^3}
	\left(
	\prod_{\nu=1}^{\kappa}D_{t,U}^{M_\nu}\nabla U
	\right)\Psi^{[j]}\,\rd x,
	\qquad
	\kappa\geq1,
	\qquad
	0\leq j\leq M-1,
	\qquad
	j+\sum_{\nu=1}^{\kappa}(M_\nu+1)=M.
	$$
	
	For $1\leq j\leq M-1$, the induction hypothesis \eqref{eq:linearized-omega-potential-material-induction} applies with H{\"o}lder exponents $\alpha_r$ and $\alpha_{r+1}$ since $j+r+1\leq M+r\leq K^*+1$.
	For $j=0$, use \eqref{eq:linearized-potential-material-zero} with one additional spatial derivative at both exponents.
	For $N_\nu\leq N$, the coefficient derivatives satisfy
	$M_\nu\leq M-1\leq K^*$ and $M_\nu+N_\nu\leq M+N-1\leq K^\sharp-2$.
	Applying \eqref{eq:holder-product-lower-exponent} with $\eta=\alpha_{r+1}$
	and combining these estimates with
	\eqref{eq:linearized-U-integer-assumption}--\eqref{eq:linearized-U-Holder-assumption}, we obtain
	\begin{equation}\label{eq:linearized-potential-correction-product-estimate}
	\begin{aligned}
		&\sup_{t\in J}
		\left\|
		\left(
		\prod_{\nu=1}^{\kappa}D_{t,U}^{M_\nu}\nabla U
		\right)\nabla \Psi^{[j]}
		\right\|_{N+\alpha_r}\\
		&\quad\lesssim
		\sum_{N_0+\sum_{\nu=1}^{\kappa}N_\nu=N}\Biggl[
		\nrm{\nabla \Psi^{[j]}}_{N_0+\alpha_r,J}
		\prod_{\nu=1}^{\kappa}\nrm{D_{t,U}^{M_\nu}\nabla U}_{N_\nu,J}\\
		&\hspace{70pt}+
		\nrm{\nabla \Psi^{[j]}}_{N_0+\alpha_{r+1},J}
		\sum_{i=1}^{\kappa}\nrm{D_{t,U}^{M_i}\nabla U}_{N_i+\alpha_r,J}
		\prod_{\substack{\nu=1\\\nu\neq i}}^{\kappa}
		\nrm{D_{t,U}^{M_\nu}\nabla U}_{N_\nu,J}\Biggr]\\
		&\quad\lesssim
		A_G\Lkb{N+1,K_{G,M-1},\lambda,\tilde\lambda}{\max\{2,s\}}{\alpha_r}
		\sum_{N_0+\sum_{\nu=1}^{\kappa}N_\nu=N}
		\MM{N_0+1,K_{G,M-1},\lambda,\tilde\lambda}
		\prod_{\nu=1}^{\kappa}\MM{N_\nu,K_{G,M-1},\lambda,\tilde\lambda}\\
		&\qquad\cdot\sum_{h=0}^{\min\{1,j\}}\tau^{1-h}\MM{j-h,N_t,\mu,\tilde\mu}
		\prod_{\nu=1}^{\kappa}\MM{M_\nu+1,N_t,\mu,\tilde\mu}\\
		&\quad\lesssim
		\Lkb{N+1,K_{G,M-1},\lambda,\tilde\lambda}{\max\{2,s\}}{\alpha_r}
		\MM{N+1,K_{G,M-1},\lambda,\tilde\lambda}
		A_G\Bigl(
		\tau\MM{M,N_t,\mu,\tilde\mu}
		+\MM{M-1,N_t,\mu,\tilde\mu}
		\Bigr),
	\end{aligned}
	\end{equation}
	where the estimates for $\Psi^{[j]}$ at spatial order $N_0+1$
	give the same H{\"o}lder and $\mathcal M$ factors as those for $\omega^{[j]}$ at order $N_0$.
	Thus the inequalities for the scale factors in the proof of \eqref{eq:linearized-omega-correction-product-estimate} also apply here.
	For the mean, use
	$\nrm{\Psi^{[j]}}_{0,J}\leq\nrm{\Psi^{[j]}}_{1+\alpha_r,J}$,
	together with the $C^0$ bound in \eqref{eq:linearized-U-integer-assumption}.
	
	The expansion \eqref{eq:linearized-potential-div-curl-products} and \eqref{eq:linearized-potential-correction-product-estimate} give,
	for $1\leq M\leq M_*+1$, $M+r\leq K^*+1$, and $M+N\leq K^\sharp-1$,
	\begin{equation}\label{eq:linearized-potential-div-curl-defect-estimate}
	\begin{aligned}
		&\nrm{\Div\Psi^{[M]}}_{N+\alpha_r,J}
		+\nrm{\curl\Psi^{[M]}-\omega^{[M]}}_{N+\alpha_r,J}\\
		&\ \lesssim
		\Lkb{N+1,K_{G,M-1},\lambda,\tilde\lambda}{\max\{2,s\}}{\alpha_r}
		\MM{N+1,K_{G,M-1},\lambda,\tilde\lambda}
		A_G\Bigl(
		\tau\MM{M,N_t,\mu,\tilde\mu}
		+\MM{M-1,N_t,\mu,\tilde\mu}
		\Bigr).
	\end{aligned}
	\end{equation}
	By \eqref{eq:linearized-potential-div-curl-defect-estimate} and the estimates for $\omega^{[M]}$
	and the spatial mean of $\Psi^{[M]}$, Lemma~\ref{lem:periodic-Hodge-estimates}
	at spatial order $N+1$ gives
	$$
	\begin{aligned}
		&\sup_{t\in J}\|\Psi^{[M]}(t,\cdot)\|_{N+1+\alpha_r}\\
		&\quad\lesssim
		\nrm{\Div\Psi^{[M]}}_{N+\alpha_r,J}
		+\nrm{\curl\Psi^{[M]}}_{N+\alpha_r,J}
		+\sup_{t\in J}\left|\fint_{\mathbb T^3}\Psi^{[M]}(t,x)\,\rd x\right|\\
		&\quad\lesssim
		\Lkb{N+1,K_{G,M-1},\lambda,\tilde\lambda}{\max\{2,s\}}{\alpha_r}
		\MM{N+1,K_{G,M-1},\lambda,\tilde\lambda}
		A_G\Bigl(
		\tau\MM{M,N_t,\mu,\tilde\mu}
		+\MM{M-1,N_t,\mu,\tilde\mu}
		\Bigr).
	\end{aligned}
	$$
	This closes the joint induction for $\omega^{[M]}$ and $\Psi^{[M]}$.
	Taking $r=0$ proves
	\eqref{eq:linearized-spatial-material-scale-conclusion-Calpha}.
	
	We next estimate $\Psi^{[M]}$ in $C^{N+\alpha_r}$
	using the equation for $\Psi$.
	For every smooth, divergence-free, zero-mean vector
	field $g$, define
	$$
	\begin{aligned}
		\mathcal B_Ug
		&:=\mathcal B\Div\left(
		U\otimes\curl g+\curl g\otimes U
		\right)-U\cdot\nabla g,\quad
		(\cC_Ug)_{ik}
		:=(g\times\nabla U_i)_k.
	\end{aligned}
	$$
	After the replacement $G\mapsto\mathcal B\curl G$, the field denoted by $G$
	is divergence-free and has zero spatial mean. Thus
	\eqref{eq:direct-vector-curl-standard-form} gives $\mathcal B\curl G=G$,
	and \eqref{eq:transport-vector-potential-equation} yields
	$$
	D_{t,U}\Psi+\mathcal B_U\Psi=G.
	$$
	For $g=\Psi$, the definition of $\cC_Ug$ and $\curl\Psi=\omega$ give
	$
	\Div(\cC_U\Psi)=(\omega\cdot\nabla)U.
	$
	The definition of $Q$ in \eqref{def:periodic-Hodge-potentials}
	and the pressure formula \eqref{eq:pressure-vector-potential-identity} give
	$
	\pi=-2Q[\cC_U\Psi].
	$
	
	For every zero-mean vector field $h$, the definition of $\mathcal B$
	implies
	$
	\curl\mathcal B h=h-\nabla\Delta^{-1}\Div h.
	$
	Applying this identity to $h=\Div(U\otimes\omega+\omega\otimes U)$,
	which has zero mean by periodicity, and using \eqref{eq:pressure-linearized}, we obtain
	$$
	\curl(\mathcal B_U\Psi)
	=\Div(U\otimes\omega+\omega\otimes U)
	+\nabla\pi-\curl(U\cdot\nabla\Psi).
	$$
	Since $\omega=\curl\Psi$ and $\Div\omega=0$, one has
	$$
	\begin{aligned}
		\Div(U\otimes\omega+\omega\otimes U)
		-\curl(U\cdot\nabla\Psi)=\curl((\nabla U)^T\Psi)+2(\omega\cdot\nabla)U.
	\end{aligned}
	$$
	Therefore,
	\begin{equation}\label{eq:linearized-BU-curl-identity}
	\begin{aligned}
		\curl(\mathcal B_U\Psi)
		&=\curl((\nabla U)^T\Psi)
		+2\Div(\cC_U\Psi)+\nabla\pi
		=\curl\left((\nabla U)^T\Psi+2W[\cC_U\Psi]\right).
	\end{aligned}
	\end{equation}
	Here the second equality follows from \eqref{eq:periodic-Hodge-decomposition}
	and $\pi=-2Q[\cC_U\Psi]$.
	The definition of $\mathcal B_U$, the identity
	$\Div\mathcal B h=0$, and periodic integration by parts also give
	\begin{equation}\label{eq:linearized-BU-divergence-mean-identities}
	\begin{aligned}
		\Div(\mathcal B_U\Psi)
		&=-\partial_i\left((\partial_jU_i)\Psi_j-(\Div U)\Psi_i\right),\quad
		\fint_{\mathbb T^3}\mathcal B_U\Psi\,\rd x
		=\fint_{\mathbb T^3}\Psi\Div U\,\rd x.
	\end{aligned}
	\end{equation}
	
	We again argue by induction on $M$. For $1\leq M\leq M_*+1$ and
	$r,N\in\mathbb N_0$ with $M+r\leq K^*+1$ and $M+N\leq K^\sharp-1$,
	we prove
	\begin{equation}
		\label{eq:linearized-potential-material-induction}
	\begin{aligned}
		&\sup_{t\in J}\|\Psi^{[M]}(t,\cdot)\|_{N+\alpha_r}\\
		&\qquad\lesssim
		\Lkb{N,K_{G,M-1},\lambda,\tilde\lambda}{\max\{2,s\}}{\alpha_r}
		\MM{N,K_{G,M-1},\lambda,\tilde\lambda}
		A_G\Bigl(
		\tau\MM{M,N_t,\mu,\tilde\mu}
		+\MM{M-1,N_t,\mu,\tilde\mu}
		\Bigr).
	\end{aligned}
	\end{equation}
	The induction starts from \eqref{eq:linearized-potential-material-zero}.
	Assume the estimate for $\Psi^{[j]}$ at all
	lower material orders in the stated ranges. In the identity
	$
	\Psi^{[M]}=D_{t,U}^{M-1}G-D_{t,U}^{M-1}(\mathcal B_U\Psi),
	$
	the $C^{N+\alpha_r}$ norm of $D_{t,U}^{M-1}G$ is estimated by \eqref{eq:linearized-projected-source-material} at material order $M-1$.
	It remains to estimate $D_{t,U}^{M-1}(\mathcal B_U\Psi)$ using the lower material derivatives of $\Psi$.
	
	Set
	$
	\Theta^{[m]}:=D_{t,U}^{m}(\mathcal B_U\Psi).
	$
	At material order zero, define
	$$
	\begin{aligned}
		(L_0^{\Div})_i
		:=-(\partial_jU_i)\Psi_j+(\Div U)\Psi_i,\quad
		L_0^{\curl}
		:=(\nabla U)^T\Psi+2W[\cC_U\Psi].
	\end{aligned}
	$$
	Identities~\eqref{eq:linearized-BU-curl-identity}--\eqref{eq:linearized-BU-divergence-mean-identities} give
	$$
	\Div\Theta^{[0]}=\Div L_0^{\Div},
	\qquad
	\curl\Theta^{[0]}=\curl L_0^{\curl},
	\qquad
	\fint_{\mathbb T^3}\Theta^{[0]}\,\rd x
	=\fint_{\mathbb T^3}\Psi\Div U\,\rd x.
	$$
	Suppose that $L_m^{\Div}$ and $L_m^{\curl}$ have been defined so that
	$$
	\Div\Theta^{[m]}=\Div L_m^{\Div},
	\qquad
	\curl\Theta^{[m]}=\curl L_m^{\curl}.
	$$
	Since $\Div(L_m^{\Div}-\Theta^{[m]})=0$ and
	$\curl(L_m^{\curl}-\Theta^{[m]})=0$, \eqref{eq:first-order-material-spatial-commutator} gives
	\begin{equation}\label{eq:linearized-potential-div-curl-recursion}
	\begin{aligned}
		\Div\Theta^{[m+1]}
		&=\Div(D_{t,U}L_m^{\Div})
		+[D_{t,U},\Div](L_m^{\Div}-\Theta^{[m]}),\\
		\curl\Theta^{[m+1]}
		&=\curl(D_{t,U}L_m^{\curl})
		+[D_{t,U},\curl](L_m^{\curl}-\Theta^{[m]}).
	\end{aligned}
	\end{equation}
	Apply \eqref{eq:linearized-divergence-free-commutator} and
	\eqref{eq:linearized-curl-free-commutator} to
	$L_m^{\Div}-\Theta^{[m]}$ and $L_m^{\curl}-\Theta^{[m]}$, respectively. We obtain
	$$
	\begin{aligned}
		\relax[D_{t,U},\Div](L_m^{\Div}-\Theta^{[m]})
		&=-\partial_i\Bigl(
		(\partial_jU_i)(L_m^{\Div}-\Theta^{[m]})_j
		-(\Div U)(L_m^{\Div}-\Theta^{[m]})_i\Bigr),\\
		\bigl([D_{t,U},\curl](L_m^{\curl}-\Theta^{[m]})\bigr)_k
		&=\epsilon_{kij}\partial_i\Bigl(
		(\partial_jU_l)(L_m^{\curl}-\Theta^{[m]})_l\Bigr).
	\end{aligned}
	$$
	Substituting these identities into \eqref{eq:linearized-potential-div-curl-recursion}, we define
	$$
	\begin{aligned}
		(L_{m+1}^{\Div})_i
		&:=D_{t,U}(L_m^{\Div})_i
		-(\partial_jU_i)(L_m^{\Div}-\Theta^{[m]})_j
		+(\Div U)(L_m^{\Div}-\Theta^{[m]})_i,\\
		(L_{m+1}^{\curl})_j
		&:=D_{t,U}(L_m^{\curl})_j
		+(\partial_jU_{k})(L_m^{\curl}-\Theta^{[m]})_{k}.
	\end{aligned}
	$$
	With these definitions,
	$
	\Div\Theta^{[m+1]}=\Div L_{m+1}^{\Div}
	$
	and
	$
	\curl\Theta^{[m+1]}=\curl L_{m+1}^{\curl}.
	$
	The mean satisfies
	$$
	\fint_{\mathbb T^3}\Theta^{[m+1]}\,\rd x
	=
	\frac{\rd}{\rd t}\fint_{\mathbb T^3}\Theta^{[m]}\,\rd x
	-
	\fint_{\mathbb T^3}\Theta^{[m]}\Div U\,\rd x.
	$$
	
	Set $m=M-1$ and substitute
	$\Theta^{[j]}=D_{t,U}^{j}G-\Psi^{[j+1]}$ in the recursions.
	Then $L_{M-1}^{\Div}$ and $L_{M-1}^{\curl}$ are finite linear combinations
	of contractions of the following products:
	\begin{equation}\label{eq:linearized-L-product-forms}\begin{alignedat}{4}
		&\left(\prod_{\nu=1}^{\kappa}D_{t,U}^{M_\nu}\nabla U\right)\Psi^{[j]},
		&\quad&\kappa\geq1,
		&\quad&0\leq j\leq M-1,
		&\quad&j+\sum_{\nu=1}^{\kappa}(M_\nu+1)=M,\\[3pt]
		&\left(\prod_{\nu=1}^{\kappa}D_{t,U}^{M_\nu}\nabla U\right)D_{t,U}^{k}G,
		&\quad&\kappa\geq1,
		&\quad&0\leq k\leq M-2,
		&\quad&k+\sum_{\nu=1}^{\kappa}(M_\nu+1)=M-1,\\[3pt]
		&\left(\prod_{\nu=1}^{\kappa}D_{t,U}^{M_\nu}\nabla U\right)D_{t,U}^{k}W[\cC_U\Psi],
		&\quad&\kappa\geq0,
		&\quad&0\leq k\leq M-1,
		&\quad&k+\sum_{\nu=1}^{\kappa}(M_\nu+1)=M-1.
	\end{alignedat}\end{equation}
	The second class is absent when $M=1$.
	The third occurs only in $L_{M-1}^{\curl}$, with the empty product equal to $1$.
	These product forms and order identities follow from the recursions and the substitution above,
	as in the proof of
	\eqref{eq:material-Hodge-correction-product}--\eqref{eq:material-Hodge-correction-order}.
	
	We first estimate the products in \eqref{eq:linearized-L-product-forms} containing $\Psi^{[j]}$ and $D_{t,U}^{k}G$.
	Write $N_0+\sum_{\nu=1}^{\kappa}N_\nu=N$, with $N_0$ the number of spatial derivatives assigned to the last factor.
	At positive material orders, use \eqref{eq:linearized-projected-source-material} for $D_{t,U}^{k}G$
	and the induction hypothesis \eqref{eq:linearized-potential-material-induction} for $\Psi^{[j]}$.
	The estimates with H{\"o}lder exponent $\alpha_{r+1}$ apply because
	$
	j+r+1\leq M+r\leq K^*+1,\ k+r+1\leq M+r-1\leq K^*.
	$
	At material order zero, use \eqref{eq:direct-vector-curl-material-zero} for the projected $G$ and
	\eqref{eq:linearized-potential-material-zero} for $\Psi^{[0]}$, at both exponents.
	Apply \eqref{eq:holder-product-lower-exponent} with
	$g_\nu=D_{t,U}^{M_\nu}\nabla U$, $\alpha=\alpha_r$,
	$\eta=\alpha_{r+1}$, and $F=\Psi^{[j]}$ or $D_{t,U}^{k}G$.
	We obtain the following estimates for the first and second classes.
	\begin{align}
		&\Nrm{\left(\prod_{\nu=1}^{\kappa}D_{t,U}^{M_\nu}\nabla U\right)\Psi^{[j]}}_{N+\alpha_r,J}\notag\\
		&\quad\lesssim
		\sum_{N_0+\sum_{\nu=1}^{\kappa}N_\nu=N}\Biggl[
		\nrm{\Psi^{[j]}}_{N_0+\alpha_r,J}
		\prod_{\nu=1}^{\kappa}\nrm{D_{t,U}^{M_\nu}\nabla U}_{N_\nu,J}+
		\nrm{\Psi^{[j]}}_{N_0+\alpha_{r+1},J}
		\sum_{i=1}^{\kappa}\nrm{D_{t,U}^{M_i}\nabla U}_{N_i+\alpha_r,J}
		\prod_{\substack{\nu=1\\\nu\neq i}}^{\kappa}\nrm{D_{t,U}^{M_\nu}\nabla U}_{N_\nu,J}\Biggr]\notag\\
		&\quad\lesssim
		A_G\Lkb{N,K_{G,M-1},\lambda,\tilde\lambda}{\max\{2,s\}}{\alpha_r}
		\sum_{N_0+\sum_{\nu=1}^{\kappa}N_\nu=N}
		\MM{N_0,K_{G,M-1},\lambda,\tilde\lambda}
		\prod_{\nu=1}^{\kappa}\MM{N_\nu,K_{G,M-1},\lambda,\tilde\lambda}\notag\\
		&\qquad\quad\cdot \sum_{h=0}^{\min\{1,j\}}\tau^{1-h}\MM{j-h,N_t,\mu,\tilde\mu}
		\prod_{\nu=1}^{\kappa}\MM{M_\nu+1,N_t,\mu,\tilde\mu}\notag\\
		&\quad\lesssim
		A_G\Lkb{N,K_{G,M-1},\lambda,\tilde\lambda}{\max\{2,s\}}{\alpha_r}
		\MM{N,K_{G,M-1},\lambda,\tilde\lambda}\sum_{h=0}^{\min\{1,j\}}\tau^{1-h}\MM{M-h,N_t,\mu,\tilde\mu}.
		\label{eq:linearized-potential-local-product-estimate}\\[6pt]
		&\Nrm{\left(\prod_{\nu=1}^{\kappa}D_{t,U}^{M_\nu}\nabla U\right)D_{t,U}^{k}G}_{N+\alpha_r,J}\notag\\
		&\quad\lesssim
		\sum_{N_0+\sum_{\nu=1}^{\kappa}N_\nu=N}\Biggl[
		\nrm{D_{t,U}^{k}G}_{N_0+\alpha_r,J}
		\prod_{\nu=1}^{\kappa}\nrm{D_{t,U}^{M_\nu}\nabla U}_{N_\nu,J}+
		\nrm{D_{t,U}^{k}G}_{N_0+\alpha_{r+1},J}
		\sum_{i=1}^{\kappa}\nrm{D_{t,U}^{M_i}\nabla U}_{N_i+\alpha_r,J}
		\prod_{\substack{\nu=1\\\nu\neq i}}^{\kappa}\nrm{D_{t,U}^{M_\nu}\nabla U}_{N_\nu,J}\Biggr]\notag\\
		&\quad\lesssim
		A_G\Lkb{N,K_{G,M-1},\lambda,\tilde\lambda}{\max\{2,s\}}{\alpha_r}
		\sum_{N_0+\sum_{\nu=1}^{\kappa}N_\nu=N}
		\MM{N_0,K_{G,M-1},\lambda,\tilde\lambda}
		\prod_{\nu=1}^{\kappa}\MM{N_\nu,K_{G,M-1},\lambda,\tilde\lambda}\notag\\
		&\qquad\quad\cdot \MM{k,N_t,\mu,\tilde\mu}
		\prod_{\nu=1}^{\kappa}\MM{M_\nu+1,N_t,\mu,\tilde\mu}\notag\\
		&\quad\lesssim
		A_G\Lkb{N,K_{G,M-1},\lambda,\tilde\lambda}{\max\{2,s\}}{\alpha_r}
		\MM{N,K_{G,M-1},\lambda,\tilde\lambda}\MM{M-1,N_t,\mu,\tilde\mu}.
		\label{eq:linearized-potential-source-product-estimate}
	\end{align}

	In the second inequality in each of
	\eqref{eq:linearized-potential-local-product-estimate}--\eqref{eq:linearized-potential-source-product-estimate},
	use \eqref{eq:linearized-U-integer-assumption}--\eqref{eq:linearized-U-Holder-assumption}
	for the coefficients.
	For the coefficient factors in these products, we have
	$$
	\begin{aligned}
		A_U\MM{N_\nu+1,K^*-M_\nu,\lambda_U,\tilde\lambda_U}
		\MM{M_\nu,N_t,\mu_U,\tilde\mu_U}\leq
		\MM{N_\nu,(K^*-M_\nu-1)_+,\lambda,\tilde\lambda}
		\MM{M_\nu+1,N_t,\mu,\tilde\mu}.
	\end{aligned}
	$$
	For $M_\nu<K^*$, use $A_U\lambda_U\leq\mu$.
	If $M_\nu=K^*$, the order identity forces $M=K^*+1$,
	so $M_*\geq K^*$ and $N_t\leq K^*$; use $A_U\tilde\lambda_U\leq\tilde\mu$.
	To estimate the H{\"o}lder factors in this step, write
	$$
	L_\nu:=(K^*-M_\nu-1)_+,\qquad
	L_0:=\begin{cases}
	K_G, & F=\Psi^{[0]},\\
	(K_G-j+1)_+, & F=\Psi^{[j]},\ j\geq1,\\
	(K_G-k)_+, & F=D_{t,U}^{k}G,
	\end{cases}
	$$
	for $1\leq\nu\leq\kappa$.
	By \eqref{eq:linearized-coefficient-threshold-comparison} and the stated ranges of $j$ and $k$,
	all these thresholds are at least $K_{G,M-1}$.
	Thus, for $1\leq i\leq\kappa$,
	\begin{equation*}
		\Lab{N_i,L_i,\lambda,\tilde\lambda}{\alpha_r}
		\Lkb{N_0,L_0,\lambda,\tilde\lambda}{\max\{2,s\}}{\alpha_{r+1}}
		\leq\Lkb{N,K_{G,M-1},\lambda,\tilde\lambda}{\max\{2,s\}}{\alpha_r},
	\end{equation*}
	using $2\alpha_{r+1}=\alpha_r$ and
	$1+\max\{2,s\}/2\leq\max\{2,s\}$, as in the proof of
	Proposition~\ref{prop:periodic-Hodge-material-estimates}.
	The terms with only one H{\"o}lder factor satisfy the same estimate by the definition of $\Xi_{\alpha_r}$.
	Estimate~\eqref{eq:appendix-derivative-weight-product}
	gives the final inequalities in
	\eqref{eq:linearized-potential-local-product-estimate}--%
	\eqref{eq:linearized-potential-source-product-estimate}.
	
	For the products containing $D_{t,U}^kW[\cC_U\Psi]$ in \eqref{eq:linearized-L-product-forms}, we first estimate the material derivatives of $W[\cC_U\Psi]$.
	For $0\leq m\leq M-1$, $m+r\leq K^*$, and $m+N\leq K^\sharp-2$, we claim
	\begin{equation}
		\label{eq:linearized-auxiliary-Hodge-bound}
		\begin{aligned}
			&\nrm{D_{t,U}^{m} Q[\cC_U\Psi]}_{N+\alpha_r,J}+\nrm{D_{t,U}^{m} W[\cC_U\Psi]}_{N+\alpha_r,J}\\
			&\qquad\lesssim A_G\Lkb{N,K_{G,m},\lambda,\tilde\lambda}{\max\{2,s\}}{\alpha_r}
			\cdot\MM{N,K_{G,m},\lambda,\tilde\lambda}
			\Bigl(\tau\MM{m+1,N_t,\mu,\tilde\mu}+\MM{m,N_t,\mu,\tilde\mu}\Bigr).
		\end{aligned}
	\end{equation}
	Fix $m$ in this range. Since $(\cC_U\Psi)_{ik}=(\Psi\times\nabla U_i)_k$,
	each term in $D_{t,U}^{k}(\cC_U\Psi)$ is a contraction of
	$(D_{t,U}^{k-j}\nabla U)\Psi^{[j]}$, with $0\leq j\leq k$.
	Estimates \eqref{eq:linearized-U-integer-assumption}--\eqref{eq:linearized-U-Holder-assumption},
	the induction hypothesis \eqref{eq:linearized-potential-material-induction},
	and \eqref{eq:holder-product-lower-exponent}
	give, for $k,n,r'\in\mathbb N_0$ with $k\leq m$, $k+r'\leq K^*$, and $k+n\leq K^\sharp-2$,
	\begin{equation}\label{eq:linearized-auxiliary-tensor-material-estimate}
	\begin{aligned}
	\nrm{D_{t,U}^{k}(\cC_U\Psi)}_{n+\alpha_{r'},J}
	\lesssim A_G\Lkb{n,K_{G,k},\lambda,\tilde\lambda}{\max\{2,s\}}{\alpha_{r'}}
	\MM{n,K_{G,k},\lambda,\tilde\lambda}\Bigl(\tau\MM{k+1,N_t,\mu,\tilde\mu}+\MM{k,N_t,\mu,\tilde\mu}\Bigr).
	\end{aligned}
	\end{equation}
	Here $j=0$ is covered by \eqref{eq:linearized-potential-material-zero}.
	The scale inequalities are those in \eqref{eq:linearized-potential-local-product-estimate}
	with $M=k+1$ and $\kappa=1$. Only lower material orders $j\leq k\leq m<M$
	of $\Psi$ occur, and $j+r'+1\leq k+r'+1\leq K^*+1$ ensures that
	the estimates with H{\"o}lder exponent $\alpha_{r'+1}$ apply.
	Choose
	$$
	A_T:=C A_G\left(1+\tau
	\frac{\MM{m+1,N_t,\mu,\tilde\mu}}{\MM{m,N_t,\mu,\tilde\mu}}\right),
	$$
	where $C$ is the implicit constant in \eqref{eq:linearized-auxiliary-tensor-material-estimate},
	taken uniformly over the stated ranges of $k,n,r'$.
	For $0\leq k\leq m$, the definition of
	$\mathcal M$ and $\mu\leq\tilde\mu$ give
	$$
	\frac{\MM{k+1,N_t,\mu,\tilde\mu}}{\MM{k,N_t,\mu,\tilde\mu}}
	=
	\begin{cases}
		\mu, & k<N_t,\\
		\tilde\mu, & k\geq N_t,
	\end{cases}
	\leq
	\frac{\MM{m+1,N_t,\mu,\tilde\mu}}{\MM{m,N_t,\mu,\tilde\mu}}.
	$$
	Hence
	$$
	C A_G\Bigl(\tau\MM{k+1,N_t,\mu,\tilde\mu}+\MM{k,N_t,\mu,\tilde\mu}\Bigr)
	\leq A_T\MM{k,N_t,\mu,\tilde\mu}.
	$$
	This verifies \eqref{eq:periodic-Hodge-T-Holder-assumption} for $T=\cC_U\Psi$,
	with $K_T=K_G$, $s$ replaced by $\max\{2,s\}$, and scales $(\lambda,\tilde\lambda,\mu,\tilde\mu)$.
	We apply Proposition~\ref{prop:periodic-Hodge-material-estimates} with
	$\upsilon=U$, $A_\upsilon=A_U$, $M_*=m$, $r^*=K^*$,
	and $K^\sharp$ replaced by $K^\sharp-1$, with $K^*$ and $N_t$ unchanged.
	Its coefficient assumptions follow from
	\eqref{eq:linearized-U-integer-assumption}--\eqref{eq:linearized-U-Holder-assumption},
	and the parameter condition holds since $K^\sharp-1\geq K^*+1$.
	The combined scales remain $(\lambda,\tilde\lambda,\mu,\tilde\mu)$.
	Substituting $A_T$ into the conclusion at material order $m$ gives
	\eqref{eq:linearized-auxiliary-Hodge-bound}, using \eqref{eq:Hodge-material-zero-estimate} when $m=0$.

	For the products in \eqref{eq:linearized-L-product-forms} containing $D_{t,U}^{k}W[\cC_U\Psi]$ with $\kappa\geq1$, apply
	\eqref{eq:linearized-auxiliary-Hodge-bound} and \eqref{eq:holder-product-lower-exponent} to obtain
	\begin{equation}
		\label{eq:linearized-potential-Hodge-product-estimate}
		\begin{aligned}
		&\Nrm{\left(\prod_{\nu=1}^{\kappa}D_{t,U}^{M_\nu}\nabla U\right)D_{t,U}^{k}W[\cC_U\Psi]}_{N+\alpha_r,J}\\
		&\quad\lesssim
		\sum_{N_0+\sum_{\nu=1}^{\kappa}N_\nu=N}\Biggl[
		\nrm{D_{t,U}^{k}W[\cC_U\Psi]}_{N_0+\alpha_r,J}
		\prod_{\nu=1}^{\kappa}\nrm{D_{t,U}^{M_\nu}\nabla U}_{N_\nu,J}\\
		&\hspace{70pt}+
		\nrm{D_{t,U}^{k}W[\cC_U\Psi]}_{N_0+\alpha_{r+1},J}
		\sum_{i=1}^{\kappa}\nrm{D_{t,U}^{M_i}\nabla U}_{N_i+\alpha_r,J}
		\prod_{\substack{\nu=1\\\nu\neq i}}^{\kappa}\nrm{D_{t,U}^{M_\nu}\nabla U}_{N_\nu,J}\Biggr]\\
		&\quad\lesssim
		A_G\Lkb{N,K_{G,M-1},\lambda,\tilde\lambda}{\max\{2,s\}}{\alpha_r}
		\sum_{N_0+\sum_{\nu=1}^{\kappa}N_\nu=N}
		\MM{N_0,K_{G,M-1},\lambda,\tilde\lambda}
		\prod_{\nu=1}^{\kappa}\MM{N_\nu,K_{G,M-1},\lambda,\tilde\lambda}\\
		&\qquad\quad\cdot \Bigl(\tau\MM{k+1,N_t,\mu,\tilde\mu}+\MM{k,N_t,\mu,\tilde\mu}\Bigr)
		\prod_{\nu=1}^{\kappa}\MM{M_\nu+1,N_t,\mu,\tilde\mu}\\
		&\quad\lesssim
		A_G\Lkb{N,K_{G,M-1},\lambda,\tilde\lambda}{\max\{2,s\}}{\alpha_r}
		\MM{N,K_{G,M-1},\lambda,\tilde\lambda} \Bigl(\tau\MM{M,N_t,\mu,\tilde\mu}+\MM{M-1,N_t,\mu,\tilde\mu}\Bigr).
		\end{aligned}
	\end{equation}
	For the second inequality in \eqref{eq:linearized-potential-Hodge-product-estimate}, we used
	\eqref{eq:linearized-U-integer-assumption}--\eqref{eq:linearized-U-Holder-assumption}
	and \eqref{eq:linearized-auxiliary-Hodge-bound} with H{\"o}lder exponents $\alpha_r$ and $\alpha_{r+1}$.
	The auxiliary Hodge estimate applies since $k+r+1\leq M+r-1\leq K^*$ and
	$k+N_0\leq M+N-2\leq K^\sharp-3$.
	The coefficient scales and H{\"o}lder factors are treated as in the proof of
	\eqref{eq:linearized-potential-local-product-estimate}, using $K_{G,k}\geq K_{G,M-1}$.
	The last inequality follows from
	the order identities in \eqref{eq:linearized-L-product-forms}
	and \eqref{eq:appendix-derivative-weight-product}.
	When $\kappa=0$, $k=M-1$, and the same final estimate follows directly from
	\eqref{eq:linearized-auxiliary-Hodge-bound} at exponent $\alpha_r$.

	It remains to estimate the spatial mean of $\Theta^{[M-1]}$.
	Starting from $\fint_{\mathbb T^3}\Theta^{[0]}\,\rd x=\fint_{\mathbb T^3}\Psi\Div U\,\rd x,$
	iterate \eqref{eq:linearized-omega-mean-recursion} with $\omega$ replaced by $\Theta^{[0]}$ and substitute
	$\Theta^{[j]}=D_{t,U}^{j}G-\Psi^{[j+1]}$.
	By \eqref{eq:linearized-omega-mean-product-derivative}, the mean
	$\fint_{\mathbb T^3}\Theta^{[M-1]}\,\rd x$ is a finite linear combination of terms $\fint_{\mathbb T^3}H\,\rd x$,
	where $H$ is a fixed tensor contraction of a product from the first two classes in \eqref{eq:linearized-L-product-forms}.
	Applying \eqref{eq:linearized-potential-local-product-estimate}--%
	\eqref{eq:linearized-potential-source-product-estimate} to each integrand at spatial order zero gives
	\begin{equation}
		\label{eq:linearized-Theta-mean-estimate}
		\begin{aligned}
		&\sup_{t\in J}\left|\fint_{\mathbb T^3}\Theta^{[M-1]}(t,x)\,\rd x\right|\lesssim
		A_G\Lkb{0,K_{G,M-1},\lambda,\tilde\lambda}{\max\{2,s\}}{\alpha_r}
		\Bigl(\tau\MM{M,N_t,\mu,\tilde\mu}+\MM{M-1,N_t,\mu,\tilde\mu}\Bigr).
		\end{aligned}
	\end{equation}
	Applying \eqref{eq:linearized-potential-local-product-estimate}--\eqref{eq:linearized-potential-Hodge-product-estimate} to the products in \eqref{eq:linearized-L-product-forms}
	controls $L_{M-1}^{\Div}$ and $L_{M-1}^{\curl}$ in $C^{N+\alpha_r}$.
	Using $\Div\Theta^{[M-1]}=\Div L_{M-1}^{\Div}$ and
	$\curl\Theta^{[M-1]}=\curl L_{M-1}^{\curl}$, apply
	Lemma~\ref{lem:periodic-Hodge-estimates} as for $\omega^{[M]}$,
	with $L=N$ for $N\geq1$ and the divergence and curl estimated in
	$B^{-1+\alpha_r}_{\infty,\infty}$ for $N=0$. Together with
	\eqref{eq:linearized-Theta-mean-estimate}, we obtain
	$$
	\begin{aligned}
		&\sup_{t\in J}\|D_{t,U}^{M-1}(\mathcal B_U\Psi)(t,\cdot)\|_{N+\alpha_r}\\
		&\qquad\lesssim
		\nrm{L_{M-1}^{\Div}}_{N+\alpha_r,J}
		+\nrm{L_{M-1}^{\curl}}_{N+\alpha_r,J}
		+\sup_{t\in J}\left|\fint_{\mathbb T^3}\Theta^{[M-1]}(t,x)\,\rd x\right|\\
		&\qquad\lesssim
		\Lkb{N,K_{G,M-1},\lambda,\tilde\lambda}{\max\{2,s\}}{\alpha_r}
		\MM{N,K_{G,M-1},\lambda,\tilde\lambda}
		A_G\Bigl(
		\tau\MM{M,N_t,\mu,\tilde\mu}
		+\MM{M-1,N_t,\mu,\tilde\mu}
		\Bigr).
	\end{aligned}
	$$
	Combining this estimate for $\Theta^{[M-1]}$ with \eqref{eq:linearized-projected-source-material} in
	$\Psi^{[M]}=D_{t,U}^{M-1}G-\Theta^{[M-1]}$ closes the induction.
	Taking $r=0$ proves \eqref{eq:potential-spatial-material-scale-conclusion}.
	
	For the pressure estimate, we first establish a joint estimate for $\omega^{[M]}$ and $\Psi^{[M]}$.
	For
	$1\leq M\leq M_*+1$,
	$M+r\leq K^*+1$, $N\geq1$, and $M+N\leq K^\sharp-1$, we have
	\begin{equation}\label{eq:linearized-joint-lower-spatial-estimate}
	\begin{aligned}
		&\sup_{t\in J}\left(
		\|\omega^{[M]}(t,\cdot)\|_{N-1+\alpha_r}
		+\|\Psi^{[M]}(t,\cdot)\|_{N+\alpha_r}
		\right)\\
		&\qquad\leq\nrm{\Psi^{[M]}}_{N+\alpha_r,J}
		+\nrm{\curl\Psi^{[M]}}_{N-1+\alpha_r,J}+\nrm{\curl\Psi^{[M]}-\omega^{[M]}}_{N-1+\alpha_r,J}\\
		&\qquad\lesssim\nrm{\Psi^{[M]}}_{N+\alpha_r,J}
		+\nrm{\curl\Psi^{[M]}-\omega^{[M]}}_{N-1+\alpha_r,J}\\
		&\qquad\lesssim
		\Lkb{N,K_{G,M-1},\lambda,\tilde\lambda}{\max\{2,s\}}{\alpha_r}
		\MM{N,K_{G,M-1},\lambda,\tilde\lambda}
		A_G\Bigl(
		\tau\MM{M,N_t,\mu,\tilde\mu}
		+\MM{M-1,N_t,\mu,\tilde\mu}
		\Bigr).
	\end{aligned}
	\end{equation}
	For the last inequality in \eqref{eq:linearized-joint-lower-spatial-estimate},
	use \eqref{eq:linearized-potential-div-curl-defect-estimate} at spatial order $N-1$
	and \eqref{eq:linearized-potential-material-induction} at spatial order $N$.
	At material order zero, the required estimate follows from \eqref{eq:linearized-potential-material-zero} and $\omega^{[0]}=\curl \Psi^{[0]}$.
	
	Set
	$
	\pi^{[M]}:=D_{t,U}^{M}\pi,
	\;
	Y_M:=\nabla \pi^{[M]}.
	$
	Equation~\eqref{eq:pressure-linearized}, $\Div\omega^{[0]}=0$, and periodicity give
	$$
	\Div Y_0=-2\Div((\omega^{[0]}\cdot\nabla)U),
	\qquad
	\curl Y_0=0,
	\qquad
	\fint_{\mathbb T^3}Y_0\,\rd x=0.
	$$
	Thus, with
	$
	F_{\pi,0}:=-2(\omega^{[0]}\cdot\nabla)U,
	$
	one has $\Div Y_0=\Div F_{\pi,0}$. Since $\omega^{[0]}=\curl \Psi^{[0]}$, antisymmetry of the
	Levi--Civita symbol also gives
	\begin{equation}
		\pi^{[0]}
		=-2\partial_i\partial_k\Delta^{-1}
		\left(\epsilon_{jkl}(\Psi^{[0]})_l\partial_jU_i\right).
		\label{eq:pressure-vector-potential-material-identity}
	\end{equation}
	This is the component form of $\pi=-2Q[\cC_U\Psi]$.
	
	For $N\geq1$, apply \eqref{eq:holder-product-lower-exponent} to obtain
	\begin{equation}\label{eq:linearized-pressure-source-zero-product-estimate}
	\begin{aligned}
	\|F_{\pi,0}\|_{N-1+\alpha_r,J}
	\lesssim\sum_{N_0+N_1=N-1}\Big(&\|\omega^{[0]}\|_{N_0+\alpha_r,J}\|\nabla U\|_{N_1,J}
	+\|\omega^{[0]}\|_{N_0+\alpha_{r+1},J}\|\nabla U\|_{N_1+\alpha_r,J}\Big).
	\end{aligned}
	\end{equation}
	At $N=0$, the Hodge potential estimate \eqref{eq:periodic-Hodge-material-zero-estimate} applied to $T=\cC_U\Psi$ and
	the pressure identity \eqref{eq:pressure-vector-potential-material-identity} give
	\begin{equation}\label{eq:linearized-pressure-low-zero-product-estimate}
	\|\pi^{[0]}\|_{\alpha_r,J}
	\lesssim\|\Psi^{[0]}\|_{\alpha_r,J}\|\nabla U\|_{0,J}
	+\|\Psi^{[0]}\|_{\alpha_{r+1},J}\|\nabla U\|_{\alpha_r,J}.
	\end{equation}
	In these product estimates, use the assumptions on $\nabla U$ in
	\eqref{eq:linearized-U-integer-assumption}--\eqref{eq:linearized-U-Holder-assumption}
	and the estimates for $\omega^{[0]}$ and $\Psi^{[0]}$ in
	\eqref{eq:linearized-solution-material-zero}--\eqref{eq:linearized-potential-material-zero},
	combining the H{\"o}lder factors with $2\alpha_{r+1}=\alpha_r$ and $1+s/2\leq\max\{2,s\}$.
	For $N\geq1$, the zero mean of $\pi^{[0]}$ and Lemma~\ref{lem:periodic-Hodge-estimates} give
	$$
	\nrm{\pi^{[0]}}_{N+\alpha_r,J}
	\lesssim\nrm{Y_0}_{N-1+\alpha_r,J}
	\lesssim\nrm{F_{\pi,0}}_{N-1+\alpha_r,J}.
	$$
	For $N\geq2$, the second inequality uses $\Div Y_0=\Div F_{\pi,0}$ and the zero curl and mean of $Y_0$.
	For $N=1$, use the lowest-order Hodge estimate, Remark~\ref{rem:app-holder-besov-equivalence},
	and \eqref{eq:app-multiplier-besov-estimate} to obtain
	$$
	\begin{aligned}
	\nrm{Y_0}_{\alpha_r,J}
	&\lesssim\sup_{t\in J}\|\Div F_{\pi,0}(t,\cdot)\|_{B^{-1+\alpha_r}_{\infty,\infty}}
	\lesssim\nrm{F_{\pi,0}}_{\alpha_r,J}.
	\end{aligned}
	$$
	Using \eqref{eq:linearized-pressure-source-zero-product-estimate}--\eqref{eq:linearized-pressure-low-zero-product-estimate}, we have, for
	$0\leq r\leq K^*$ and $0\leq N\leq K^\sharp-1$,
	\begin{equation}
	\label{eq:linearized-pressure-material-zero}
	\begin{aligned}
		\sup_{t\in J}\|\pi^{[0]}(t,\cdot)\|_{N+\alpha_r}
		&\lesssim
		\tau A_UA_G
		\Lkb{N+1,K_G,\lambda,\tilde\lambda}{\max\{2,s\}}{\alpha_r}
		\MM{N+1,K_G,\lambda,\tilde\lambda}.
	\end{aligned}
	\end{equation}
	Applying \eqref{eq:periodic-Hodge-material-zero-estimate} to $T=\cC_U\Psi$
	and using \eqref{eq:holder-product-lower-exponent} gives, for $0\leq r\leq K^*$,
	\begin{equation}\label{eq:linearized-pressure-auxiliary-material-zero}
	\begin{aligned}
		\sup_{t\in J}\|W[\cC_U\Psi](t,\cdot)\|_{\alpha_r}
		&\lesssim\nrm{\Psi^{[0]}}_{\alpha_r,J}\nrm{\nabla U}_{0,J}
		+\nrm{\Psi^{[0]}}_{\alpha_{r+1},J}\nrm{\nabla U}_{\alpha_r,J}\\
		&\lesssim\tau A_UA_G
		\Lkb{1,K_G,\lambda,\tilde\lambda}{\max\{2,s\}}{\alpha_r}
		\MM{1,K_G,\lambda,\tilde\lambda}.
	\end{aligned}
	\end{equation}

	Since $Y_M=\nabla \pi^{[M]}$, periodicity gives, for every $M\geq0$,
	$
	\curl Y_M=0,
	\fint_{\mathbb T^3}Y_M\,\rd x=0.
	$
	Hence only the divergence requires a recursive representation.
	By \eqref{eq:first-order-material-spatial-commutator}, we have
	$
	Y_{M+1}=D_{t,U}Y_M+(\nabla U)^TY_M.
	$
	Suppose that $\Div Y_M=\Div F_{\pi,M}$. Then
	$$
	\begin{aligned}
		\Div Y_{M+1}
		&=\Div\left(D_{t,U}F_{\pi,M}+(\nabla U)^TY_M\right)
		+[D_{t,U},\Div](F_{\pi,M}-Y_M).
	\end{aligned}
	$$
	Since $\Div(F_{\pi,M}-Y_M)=0$, \eqref{eq:linearized-divergence-free-commutator} gives
	$\Div Y_{M+1}=\Div F_{\pi,M+1}$, where
	$$
	\begin{aligned}
		(F_{\pi,M+1})_i
		&:=D_{t,U}(F_{\pi,M})_i
		+(\partial_iU_j)(Y_M)_j
		-(\partial_jU_i)(F_{\pi,M}-Y_M)_j
		+(\Div U)(F_{\pi,M}-Y_M)_i.
	\end{aligned}
	$$
	As in the proof of \eqref{eq:material-Hodge-correction-product}--\eqref{eq:material-Hodge-correction-order},
	the recursion expresses $F_{\pi,M}$ as a finite sum of contractions
	of products of the forms
	\begin{equation}\label{eq:linearized-pressure-divergence-products}
	\begin{alignedat}{3}
		&\left(\prod_{\nu=1}^{\kappa}D_{t,U}^{M_\nu}\nabla U\right)\omega^{[j]},
		&\quad&\kappa\geq1,\ 0\leq j\leq M,
		&\quad&j+\sum_{\nu=1}^{\kappa}(M_\nu+1)=M+1,\\[3pt]
		&\left(\prod_{\nu=1}^{\kappa}D_{t,U}^{M_\nu}\nabla U\right)Y_j,
		&\quad&\kappa\geq1,\ 0\leq j\leq M-1,
		&\quad&j+\sum_{\nu=1}^{\kappa}(M_\nu+1)=M.
	\end{alignedat}
	\end{equation}
	
	For the spatial mean, \eqref{eq:first-order-material-mean-identity} gives
	$$
	\fint_{\mathbb T^3}\pi^{[M+1]}\,\rd x
	=
	\frac{\rd}{\rd t}\fint_{\mathbb T^3}\pi^{[M]}\,\rd x
	-
	\fint_{\mathbb T^3}\pi^{[M]}\Div U\,\rd x.
	$$
	Since $\pi^{[0]}$ has zero spatial mean, apply \eqref{eq:linearized-omega-mean-recursion}--\eqref{eq:linearized-omega-mean-product-derivative} with $\omega$ replaced by $\pi$.
	Thus $\fint_{\mathbb T^3}\pi^{[M]}\,\rd x$ is a finite linear combination of fixed tensor contractions of
	$$
	\fint_{\mathbb T^3}\left(
	\prod_{\nu=1}^{\kappa}D_{t,U}^{M_\nu}\nabla U
	\right)\pi^{[j]}\,\rd x,
	\qquad\kappa\geq1,
	0\leq j\leq M-1,
	\qquad
	j+\sum_{\nu=1}^{\kappa}(M_\nu+1)=M.
	$$
	
	We prove the pressure estimate by induction on $M$, simultaneously with
	the $C^{\alpha_r}$ estimate for
	$D_{t,U}^{M}W[\cC_U\Psi]$, with the same right-hand side evaluated at $N=0$.
	Assume both estimates at lower material orders in their stated ranges.
	For $1\leq M\leq M_*+1$, $M+r+1\leq K^*+1$, and
	$M+N\leq K^\sharp-1$, we prove
	\begin{equation}\label{eq:linearized-pressure-material-induction}
	\begin{aligned}
		&\sup_{t\in J}\|\pi^{[M]}(t,\cdot)\|_{N+\alpha_r}\\
		&\lesssim
		A_U\Lkb{N+1,K_{G,M},\lambda,\tilde\lambda}{\max\{2,s\}}{\alpha_r}\cdot
		\MM{N+1,K_{G,M},\lambda,\tilde\lambda}\cdot
		A_G\Bigl(
		\tau\MM{M,N_t,\mu,\tilde\mu}
		+\MM{M-1,N_t,\mu,\tilde\mu}
		\Bigr).
	\end{aligned}
	\end{equation}
	
	For $N\geq1$, the product forms in \eqref{eq:linearized-pressure-divergence-products} and
	\eqref{eq:holder-product-lower-exponent}, with
	$N_0+\sum_{\nu=1}^{\kappa}N_\nu=N-1$, give
	\begin{equation}\label{eq:linearized-pressure-divergence-estimate}
	\begin{aligned}
		&\sup_{t\in J}\|F_{\pi,M}(t,\cdot)\|_{N-1+\alpha_r}\\
		&\qquad\lesssim
		\sum_{\substack{\kappa\geq1,\ 0\leq j\leq M\\
		j+\sum_{\nu=1}^{\kappa}(M_\nu+1)=M+1\\
		N_0+\sum_{\nu=1}^{\kappa}N_\nu=N-1}}\Biggl[
		\nrm{\omega^{[j]}}_{N_0+\alpha_r,J}
		\prod_{\nu=1}^{\kappa}\nrm{D_{t,U}^{M_\nu}\nabla U}_{N_\nu,J}\\
		&\quad\hspace{100pt}+
		\nrm{\omega^{[j]}}_{N_0+\alpha_{r+1},J}
		\sum_{i=1}^{\kappa}\nrm{D_{t,U}^{M_i}\nabla U}_{N_i+\alpha_r,J}
		\prod_{\substack{\nu=1\\\nu\neq i}}^{\kappa}\nrm{D_{t,U}^{M_\nu}\nabla U}_{N_\nu,J}\Biggr]\\
		&\qquad\quad+
		\sum_{\substack{\kappa\geq1,\ 0\leq j\leq M-1\\
		j+\sum_{\nu=1}^{\kappa}(M_\nu+1)=M\\
		N_0+\sum_{\nu=1}^{\kappa}N_\nu=N-1}}\Biggl[
		\nrm{Y_j}_{N_0+\alpha_r,J}
		\prod_{\nu=1}^{\kappa}\nrm{D_{t,U}^{M_\nu}\nabla U}_{N_\nu,J}\\
		&\quad\hspace{100pt}+
		\nrm{Y_j}_{N_0+\alpha_{r+1},J}
		\sum_{i=1}^{\kappa}\nrm{D_{t,U}^{M_i}\nabla U}_{N_i+\alpha_r,J}
		\prod_{\substack{\nu=1\\\nu\neq i}}^{\kappa}\nrm{D_{t,U}^{M_\nu}\nabla U}_{N_\nu,J}\Biggr]\\
		&\qquad\lesssim
		A_U\Lkb{N+1,K_{G,M},\lambda,\tilde\lambda}{\max\{2,s\}}{\alpha_r}
		\MM{N+1,K_{G,M},\lambda,\tilde\lambda}
		A_G\Bigl(
		\tau\MM{M,N_t,\mu,\tilde\mu}
		+\MM{M-1,N_t,\mu,\tilde\mu}
		\Bigr).
	\end{aligned}
	\end{equation}
	For the last inequality in \eqref{eq:linearized-pressure-divergence-estimate}, estimate the coefficients by
	\eqref{eq:linearized-U-integer-assumption}--\eqref{eq:linearized-U-Holder-assumption}.
	For $j\geq1$, estimate $\omega^{[j]}$ in the first sum by \eqref{eq:linearized-joint-lower-spatial-estimate}
	and $Y_j=\nabla\pi^{[j]}$ in the second sum by the induction hypothesis \eqref{eq:linearized-pressure-material-induction},
	both used at spatial order $N_0+1$.
	These estimates apply with H{\"o}lder exponent $\alpha_{r+1}$ because
	$$
	\begin{aligned}
		j+r+1&\leq M+r+1\leq K^*+1, &&1\leq j\leq M \text{ in } \omega^{[j]},\\
		j+r+2&\leq M+r+1\leq K^*+1, &&1\leq j<M \text{ in } Y_j.
	\end{aligned}
	$$
	At material order zero, use \eqref{eq:linearized-potential-material-zero} with
	$\omega^{[0]}=\curl\Psi^{[0]}$ and \eqref{eq:linearized-pressure-material-zero} with
	$Y_0=\nabla\pi^{[0]}$, at both exponents.
	The spatial orders are covered since $j+N_0+1\leq M+N\leq K^\sharp-1$
	and $M_\nu+N_\nu\leq M+N-1\leq K^\sharp-2$.
	
	In the first sum, keep the amplitude $A_U$ and the spatial weight of order $N_1+1$ from the estimate for $D_{t,U}^{M_1}\nabla U$.
	Absorb $A_U\lambda_U$ from each other coefficient into its material scale. The order identities give
	\begin{equation*}
		\begin{aligned}
			j+M_1+\sum_{\nu=2}^{\kappa}(M_\nu+1)=M,\quad
			(N_0+1)+(N_1+1)+\sum_{\nu=2}^{\kappa}N_\nu=N+1.
		\end{aligned}
	\end{equation*}
	In the second sum, the estimate for $Y_j$ has the factor $A_U$ and contributes spatial scale factors of order $N_0+2$.
	Absorbing $A_U\lambda_U$ from every coefficient into its material scale gives
	$(N_0+2)+\sum_\nu N_\nu=N+1$ and $j+\sum_\nu(M_\nu+1)=M$.
	In both sums, each coefficient contributing the factor $\lambda_U$ satisfies
	$M_\nu\leq M-1\leq K^*-1$ and $K^*-M_\nu-1\geq K_{G,M}$.
	The coefficient $D_{t,U}^{M_1}\nabla U$ in the first sum has $K^*-M_1\geq K_{G,M}$,
	and the thresholds for $\omega^{[j]}$ and $Y_j$ are also at least $K_{G,M}$.

	At this common threshold, the H{\"o}lder factors combine as in the proof of
	Proposition~\ref{prop:periodic-Hodge-material-estimates}, giving the displayed power $\max\{2,s\}$.
	Estimate~\eqref{eq:appendix-derivative-weight-product} gives the spatial and material weights
	in \eqref{eq:linearized-pressure-divergence-estimate}.
	
	For the spatial mean of $\pi^{[M]}$, use
	\eqref{eq:linearized-U-integer-assumption}--\eqref{eq:linearized-U-Holder-assumption},
	\eqref{eq:linearized-pressure-material-zero} for $j=0$, the induction hypothesis \eqref{eq:linearized-pressure-material-induction} for $j\geq1$,
	and \eqref{eq:holder-product-lower-exponent} at spatial order zero
	to obtain
	\begin{equation}\label{eq:linearized-pressure-material-mean-estimate}
	\begin{aligned}
		&\sup_{t\in J}\left|\fint_{\mathbb T^3}\pi^{[M]}(t,x)\,\rd x\right|\\
		&\quad\lesssim A_U
		\Lkb{1,K_{G,M},\lambda,\tilde\lambda}{\max\{2,s\}}{\alpha_r}
		\MM{1,K_{G,M},\lambda,\tilde\lambda}A_G\Bigl(
		\tau\MM{M,N_t,\mu,\tilde\mu}
		+\MM{M-1,N_t,\mu,\tilde\mu}
		\Bigr).
	\end{aligned}
	\end{equation}
	Apply Lemma~\ref{lem:periodic-Hodge-estimates} to $Y_M$ at spatial
	order $N-1$, using its lowest-order estimate and
	\eqref{eq:app-multiplier-besov-estimate} for the divergence when $N=1$.
	We have
	$$
	\begin{aligned}
	\|\pi^{[M]}\|_{N+\alpha_r}
	\lesssim
	\|Y_M\|_{N-1+\alpha_r}
	+
	\left|\fint_{\mathbb T^3}\pi^{[M]}\,\rd x\right|
	\lesssim
	\|F_{\pi,M}\|_{N-1+\alpha_r}
	+\left|\fint_{\mathbb T^3}\pi^{[M]}\,\rd x\right|,
	\end{aligned}
	$$
	which, together with \eqref{eq:linearized-pressure-divergence-estimate}--\eqref{eq:linearized-pressure-material-mean-estimate}, proves \eqref{eq:linearized-pressure-material-induction} for $N\geq1$.
	
	For $N=0$, set $T=\cC_U\Psi$ and use
	$
	\pi^{[M]}=-2D_{t,U}^{M}Q[\cC_U\Psi].
	$
	The terms in $D_{t,U}^{M}T$ and the recursive products with $Z=T$
	in \eqref{eq:material-Hodge-correction-product} are contractions of
	$$
	\left(
	\prod_{\nu=1}^{\kappa}D_{t,U}^{M_\nu}\nabla U
	\right)\Psi^{[j]},
	\qquad
	\kappa\geq1,
	\qquad
	0\leq j\leq M,
	\qquad
	j+\sum_{\nu=1}^{\kappa}(M_\nu+1)=M+1.
	$$
	Estimate $\Psi^{[j]}$ by \eqref{eq:linearized-potential-material-zero} for $j=0$
	and \eqref{eq:linearized-potential-material-induction} for $j\geq1$.
	For $Z=Q[T]$ or $W[T]$, \eqref{eq:material-Hodge-correction-order} gives $k<M$.
	Use the simultaneous induction hypothesis for $k\geq1$ and
	\eqref{eq:linearized-pressure-material-zero}--\eqref{eq:linearized-pressure-auxiliary-material-zero} for $k=0$,
	with $D_{t,U}^{k}Q[T]=-\pi^{[k]}/2$.
	These estimates apply with H{\"o}lder exponent $\alpha_{r+1}$ because
	$
	k+r+2\leq M+r+1\leq K^*+1,\ j+r+1\leq M+r+1\leq K^*+1,\ 2\alpha_{r+1}=\alpha_r.
	$
	Applying \eqref{eq:holder-product-lower-exponent} and the estimates for $\Psi^{[j]}$ and $D_{t,U}^{M_\nu}\nabla U$, we obtain
	\begin{align*}
	&\Nrm{\left(\prod_{\nu=1}^{\kappa}D_{t,U}^{M_\nu}\nabla U\right)\Psi^{[j]}}_{\alpha_r,J}\\
	&\quad\lesssim
	\nrm{\Psi^{[j]}}_{\alpha_r,J}
	\prod_{\nu=1}^{\kappa}\nrm{D_{t,U}^{M_\nu}\nabla U}_{0,J}+
	\nrm{\Psi^{[j]}}_{\alpha_{r+1},J}
	\sum_{i=1}^{\kappa}\nrm{D_{t,U}^{M_i}\nabla U}_{\alpha_r,J}
	\prod_{\substack{\nu=1\\\nu\neq i}}^{\kappa}\nrm{D_{t,U}^{M_\nu}\nabla U}_{0,J}\\
	&\quad\lesssim
	A_UA_G\Lkb{1,K_{G,M},\lambda,\tilde\lambda}{\max\{2,s\}}{\alpha_r}
	\MM{1,K_{G,M},\lambda,\tilde\lambda}\\
	&\qquad\cdot\sum_{h=0}^{\min\{1,j\}}\tau^{1-h}
	\MM{j-h,N_t,\mu,\tilde\mu}\MM{M_1,N_t,\mu,\tilde\mu}
	\prod_{\nu=2}^{\kappa}\MM{M_\nu+1,N_t,\mu,\tilde\mu}\\
	&\quad\lesssim
	A_UA_G\Lkb{1,K_{G,M},\lambda,\tilde\lambda}{\max\{2,s\}}{\alpha_r}
	\MM{1,K_{G,M},\lambda,\tilde\lambda}
	\sum_{h=0}^{\min\{1,j\}}\tau^{1-h}\MM{M-h,N_t,\mu,\tilde\mu}.
	\end{align*}
	For $Z=Q[T]$ or $W[T]$, the recursive products satisfy
	$k+\sum_{\nu=1}^{\kappa}(M_\nu+1)=M$. Similarly,
	\begin{align*}
	&\Nrm{\left(\prod_{\nu=1}^{\kappa}D_{t,U}^{M_\nu}\nabla U\right)D_{t,U}^kZ}_{\alpha_r,J}\\
	&\quad\lesssim
	\nrm{D_{t,U}^kZ}_{\alpha_r,J}
	\prod_{\nu=1}^{\kappa}\nrm{D_{t,U}^{M_\nu}\nabla U}_{0,J}+
	\nrm{D_{t,U}^kZ}_{\alpha_{r+1},J}
	\sum_{i=1}^{\kappa}\nrm{D_{t,U}^{M_i}\nabla U}_{\alpha_r,J}
	\prod_{\substack{\nu=1\\\nu\neq i}}^{\kappa}\nrm{D_{t,U}^{M_\nu}\nabla U}_{0,J}\\
	&\quad\lesssim
	A_UA_G\Lkb{1,K_{G,M},\lambda,\tilde\lambda}{\max\{2,s\}}{\alpha_r}
	\MM{1,K_{G,M},\lambda,\tilde\lambda}\sum_{h=0}^{\min\{1,k\}}\tau^{1-h}
	\MM{k-h,N_t,\mu,\tilde\mu}
	\prod_{\nu=1}^{\kappa}\MM{M_\nu+1,N_t,\mu,\tilde\mu}\\
	&\quad\lesssim
	A_UA_G\Lkb{1,K_{G,M},\lambda,\tilde\lambda}{\max\{2,s\}}{\alpha_r}
	\MM{1,K_{G,M},\lambda,\tilde\lambda}
	\sum_{h=0}^{\min\{1,k\}}\tau^{1-h}\MM{M-h,N_t,\mu,\tilde\mu}.
	\end{align*}
	For the intermediate inequalities, estimate the H{\"o}lder and spatial factors as in
	\eqref{eq:linearized-pressure-divergence-estimate}.
	In the first estimate, use the factor $A_U\MM{1,K^*-M_1,\lambda_U,\tilde\lambda_U}$
	from the first coefficient. In the second, induction gives
	$A_U\MM{1,K_{G,k},\lambda,\tilde\lambda}$, and $K_{G,k}\geq K_{G,M}$.
	Absorb the remaining coefficient factors into the material scales.
	The final inequalities
	follow from
	\eqref{eq:appendix-derivative-weight-product}.
	The mean terms in \eqref{eq:material-Hodge-mean-product} are estimated in the same way.
	Applying \eqref{eq:material-Hodge-reconstruction-estimate} at $N=0$, we obtain
	$$
	\begin{aligned}
		&\sup_{t\in J}\Bigl(
		\|\pi^{[M]}(t,\cdot)\|_{\alpha_r}
		+\|D_{t,U}^{M}W[\cC_U\Psi](t,\cdot)\|_{\alpha_r}
		\Bigr)\\
		&\quad\lesssim{}A_U
		\Lkb{1,K_{G,M},\lambda,\tilde\lambda}{\max\{2,s\}}{\alpha_r}
		\MM{1,K_{G,M},\lambda,\tilde\lambda}A_G\Bigl(
		\tau\MM{M,N_t,\mu,\tilde\mu}
		+\MM{M-1,N_t,\mu,\tilde\mu}
		\Bigr).
	\end{aligned}
	$$

	This closes the simultaneous induction at $N=0$.
	Taking $r=0$ proves
	\eqref{eq:pressure-spatial-material-scale-conclusion} in the stated range.
	For $M\geq1$ at exponent $\alpha_0$, the estimates for $\omega$ and $\Psi$
	use exponents through $\alpha_M$, and the pressure estimate through $\alpha_{M+1}$.
	All these exponents are included in the hypotheses for the stated range of material orders.
\end{proof}

\section{The Newton steps}\label{sec:Newton steps}

We construct and estimate the compressible and incompressible corrections
in a finite Newton iteration for the compressible Euler system,
adapting the scheme of Giri--Radu~\cite{GR24}.
Section~\ref{sec:local-inverse-divergence} provides the local right inverse of divergence
for the compressible residuals. Section~\ref{sec:periodic-Hodge-and-Newton-stresses}
supplies the Hodge decomposition and the estimates for the nonlocal operators
and the linearized equation.

We begin with the scheme and parameters in Subsection~\ref{subsec:Newton-scheme-parameters}.
In Subsection~\ref{subsec:Newton-compressible-perturbations}, we construct and
estimate the compressible corrections using Kirchhoff's formula in
Lagrangian coordinates and the Piola transform.
Subsection~\ref{subsec:Newton-incompressible-perturbation} constructs, glues,
and estimates the local incompressible corrections.
Finally, Subsection~\ref{subsec:Newton-stress-state-estimates} estimates the
updated state, stress tensors, and terminal velocity in the Newton step, and proves the
temporal support inclusions.

The Newton steps use the time scale $\tau_q$ and the flow of $\overline u_q$.
The Nash step corrects the gluing stress $S_{q,J_*}$. The flux remainders
and the two source mollification stresses enter $R_{q+1}$ through $P_{q+1,J_*}$.

Throughout this section, fix $q\in\mathbb N_0$ and a smooth triple $(\rho_q,m_q,R_q)$ satisfying the hypotheses of Proposition~\ref{prop:uniform-mollification-estimates}.
We write
$\overline{u}_q=\om_q/\orho_q$ and
$\DTL=\partial_t+\overline{u}_q\cdot\nabla$.
The cutoffs $\theta_u,\tilde\theta_u,\chi_\kappa$, the index set
$\mathscr I$, and the flow maps $X_u,\Phi_u$ are those introduced in
Section~\ref{Mollification}. The mollified fields $\orho_q$ and $\om_q$
remain fixed throughout the Newton iteration.

\subsection{The Newton scheme and parameters}
\label{subsec:Newton-scheme-parameters}

Recall from \eqref{def:Newton-stopping-depth} that
\begin{equation*}
	J_*=\left\lceil\frac{b-1}{4\alpha}\right\rceil+1.
\end{equation*}

Let $R_{q,0}$ be the stress defined in
\eqref{def:initial-Newton-flow-mollified-stress}. We initialize the Newton
iteration by
$$
\rho_{q,0}:=\rho_q,\qquad m_{q,0}:=m_q,\qquad S_{q,0}:=0,
\qquad P_{q+1,0}:=R_q-R_{q,0}.
$$
For $0\leq j\leq J_*$, let $\rho_{q,j}$ and $m_{q,j}$ denote the
density and momentum after $j$ substeps. We construct these fields and
the stresses $R_{q,j}$, $S_{q,j}$, and $P_{q+1,j}$ to satisfy
\begin{equation}\label{eq:Newton-substep-Euler-Reynolds-system}
	\left\{
	\begin{aligned}
		&\partial_t\rho_{q,j}+\Div m_{q,j}=0,\\
		&\partial_t m_{q,j}
		+\Div\left(\frac{m_{q,j}\otimes m_{q,j}}{\rho_{q,j}}\right)
		+\nabla p(\rho_{q,j})=\Div R_{q,j}+\Div S_{q,j}
		+\Div P_{q+1,j}.
	\end{aligned}
	\right.
\end{equation}
Here $S_{q,j}$ is the sum of the symmetric gluing stresses from the
first $j$ Newton substeps. The terminal sum $S_{q,J_*}$ is corrected in
the Nash step. The stress $P_{q+1,j}$ contains $R_q-R_{q,0}$, the source
mollification stresses, and the flux remainders from these substeps.
Its terminal value $P_{q+1,J_*}$ is included in
the next-stage Reynolds stress $R_{q+1}$.
The interior stress $R_{q,j}$ need not be symmetric. For
$0\leq j<J_*$, \eqref{def:periodic-Hodge-potentials} gives
$$
\Div R_{q,j}=\nabla Q[R_{q,j}]+\curl W[R_{q,j}],
\qquad \Div W[R_{q,j}]=0.
$$
We apply the Hodge decomposition before localization and place
the cutoffs outside the along-flow convolutions.

For a smooth vector potential $V$, define
\begin{equation}\label{def:Newton-perturbations-from-potential}
	\varrho[V]:=\Div V,\qquad
	w_i[V]:=-\partial_tV_i+\partial_k\left(
	(\overline{u}_{q})_iV_k-(\overline{u}_{q})_kV_i\right).
\end{equation}
Write
$
J_u:=\det\nabla\Phi_u,\ \mathbf C_u:=\cof\nabla\Phi_u=J_u(\nabla\Phi_u)^{-T}.
$
Using the kernel fixed in Section~\ref{Mollification}, define for
$(u,\kappa)\in\mathscr I$
\begin{equation}\label{def:Newton-sharp-acoustic-input}
	G_{q,j,u,\kappa}^{\mathrm{mol}}:=\theta_u\chi_\kappa^2(\Phi_u)\bigl(\eta_{\ell_{t,q}}^{\mathrm t}\ast_{X_{\overline u_q}}(J_u^{-1}Q[R_{q,j}])\bigr).
\end{equation}
For the vector potential, define
\begin{equation}\label{def:Newton-flow-mollified-vector-source}
	\begin{aligned}
		W_{q,j}^{\mathrm{mol}}(t,x):=&\bigl(\eta_{\ell_{t,q}}^{\mathrm t}\ast_{X_{\overline u_q}}W[R_{q,j}]\bigr)(t,x).
	\end{aligned}
\end{equation}
Define the two source mollification stresses included in $P_{q+1,j+1}$ by
\begin{equation}\label{def:Newton-substep-flow-mollification}
	\begin{aligned}
		P_{q+1,j+1}^{\mathrm{mol,comp}}:=\Bigl(Q[R_{q,j}]-\sum_{(u,\kappa)\in\mathscr I}\theta_uJ_uG_{q,j,u,\kappa}^{\mathrm{mol}}\Bigr)\Id,\quad
		P_{q+1,j+1}^{\mathrm{mol,inc}}:=\mathcal R\curl\bigl(W[R_{q,j}]-W_{q,j}^{\mathrm{mol}}\bigr).
	\end{aligned}
\end{equation}
Both stresses are symmetric.

Let $V_{q,j+1,u,\kappa}^{\circ}$ be the potential with zero spatial integral defined in \eqref{def:Newton-mean-corrected-Piola-potential}.
We localize it in time and set
$$
V_{q,j+1,u,\kappa}\comp
:=\theta_uV_{q,j+1,u,\kappa}^{\circ},
\qquad
(\trho_{q,j+1,u,\kappa}\comp,\tm_{q,j+1,u,\kappa}\comp)
:=(\varrho[V_{q,j+1,u,\kappa}\comp],w[V_{q,j+1,u,\kappa}\comp]).
$$
The definitions give $\partial_t\varrho[V]+\Div w[V]=0$ for every $V$.
Since $\theta_u$ depends only on time,
$$
\begin{aligned}
	\varrho[\theta_uV_{q,j+1,u,\kappa}^{\circ}]
	=\theta_u\varrho[V_{q,j+1,u,\kappa}^{\circ}],\quad
	w[\theta_uV_{q,j+1,u,\kappa}^{\circ}]
	=\theta_uw[V_{q,j+1,u,\kappa}^{\circ}]
	-(\partial_t\theta_u)V_{q,j+1,u,\kappa}^{\circ}.
\end{aligned}
$$

Differentiating the time cutoff $\theta_u$ in
$V_{q,j+1,u,\kappa}\comp=\theta_uV_{q,j+1,u,\kappa}^{\circ}$ produces the cutoff residual
$$
F_{q,j+1,u,\kappa}^{\mathrm{cut}}
:=-2(\partial_t\theta_u)
\left(\DTL V_{q,j+1,u,\kappa}^{\circ}
+(\Div\overline{u}_q)V_{q,j+1,u,\kappa}^{\circ}\right)
-(\partial_t^2\theta_u)V_{q,j+1,u,\kappa}^{\circ}.
$$
The terms in which no derivative falls on $\theta_u$ consist of
$-\theta_u\nabla(J_uG_{q,j,u,\kappa}^{\mathrm{mol}})$ and the interior
residual $F_{q,j+1,u,\kappa}\comp$. The local fields therefore satisfy
$$
\left\{
\begin{aligned}
	&\partial_t\trho_{q,j+1,u,\kappa}\comp
	+\Div\tm_{q,j+1,u,\kappa}\comp=0,\\
	&\partial_t\tm_{q,j+1,u,\kappa}\comp
	+\Div\left(\overline{u}_q\otimes\tm_{q,j+1,u,\kappa}\comp
	+\tm_{q,j+1,u,\kappa}\comp\otimes\overline{u}_q
	-\trho_{q,j+1,u,\kappa}\comp\overline{u}_q\otimes\overline{u}_q\right)+\nabla\left(p'(\orho_q)\trho_{q,j+1,u,\kappa}\comp\right)\\
	&\hspace{98pt}
	=-\theta_u\nabla(J_uG_{q,j,u,\kappa}^{\mathrm{mol}})
	+F_{q,j+1,u,\kappa}\comp
	+F_{q,j+1,u,\kappa}^{\mathrm{cut}}.
\end{aligned}
\right.
$$
The temporal supports of the perturbations and residuals satisfy
$$
\supp_t\trho_{q,j+1,u,\kappa}\comp
\cup\supp_t\tm_{q,j+1,u,\kappa}\comp
\cup\supp_tF_{q,j+1,u,\kappa}\comp
\cup\supp_tF_{q,j+1,u,\kappa}^{\mathrm{cut}}
\subset\supp\theta_u.
$$

We define the total potential by
$$
V_{q,j+1}\comp
:=\sum_{(u,\kappa)\in\mathscr I}
	\theta_uV_{q,j+1,u,\kappa}^{\circ}
=\sum_{(u,\kappa)\in\mathscr I}
V_{q,j+1,u,\kappa}\comp.
$$
The corresponding density and momentum perturbations are
$$
\begin{aligned}
	\trho_{q,j+1}\comp
	:=\varrho[V_{q,j+1}\comp]
	=\sum_{(u,\kappa)\in\mathscr I}
	\theta_u\varrho[V_{q,j+1,u,\kappa}^{\circ}],\quad
	\tm_{q,j+1}\comp
	:=w[V_{q,j+1}\comp]
	=\sum_{(u,\kappa)\in\mathscr I}
	\left(
	\theta_uw[V_{q,j+1,u,\kappa}^{\circ}]
	-(\partial_t\theta_u)V_{q,j+1,u,\kappa}^{\circ}
	\right).
\end{aligned}
$$
Define the corresponding residual sums by
$$
F_{q,j+1}\comp
:=\sum_{(u,\kappa)\in\mathscr I}F_{q,j+1,u,\kappa}\comp,
\qquad
F_{q,j+1}^{\mathrm{cut}}
:=\sum_{(u,\kappa)\in\mathscr I}F_{q,j+1,u,\kappa}^{\mathrm{cut}}.
$$
The blockwise equations sum to
\begin{equation}\label{eq:Newton-approximate-compressible-linearized-system}
	\left\{
	\begin{aligned}
		&\partial_t\trho_{q,j+1}\comp
		+\Div\tm_{q,j+1}\comp=0,\\
		&\partial_t\tm_{q,j+1}\comp
		+\Div\left(\overline{u}_q\otimes\tm_{q,j+1}\comp
		+\tm_{q,j+1}\comp\otimes\overline{u}_q
		-\trho_{q,j+1}\comp\overline{u}_q\otimes\overline{u}_q\right)+\nabla\left(p'(\orho_q)\trho_{q,j+1}\comp\right)\\
		&\hspace{88pt}=-\sum_{(u,\kappa)\in\mathscr I}\theta_u\nabla(J_uG_{q,j,u,\kappa}^{\mathrm{mol}})+F_{q,j+1}\comp+F_{q,j+1}^{\mathrm{cut}}.
	\end{aligned}
	\right.
\end{equation}

On each block, write the explicit gradient and divergence terms in $F_{q,j+1,u,\kappa}\comp$ as matrix divergences. Apply the local right inverse from Theorem~\ref{thm:inv4divEq} to the remaining zero-mean source $H_{q,j+1,u,\kappa}^{\mathrm{int}}$ defined in \eqref{def:Newton-interior-vector-remainder}. Summing over the blocks, we obtain an interior matrix $R_{q,j+1}\comp$ satisfying
$$
F_{q,j+1}\comp=\Div R_{q,j+1}\comp.
$$

Each term in the cutoff residual contains $\partial_t\theta_u$ or
$\partial_t^2\theta_u$. Its temporal support is therefore contained in
$\bigcup_u\supp\partial_t\theta_u$.
We group it by time block as
$$
F_{q,j+1}^{\mathrm{cut}}
=\sum_uF_{q,j+1,u}^{\mathrm{cut}},
\qquad
F_{q,j+1,u}^{\mathrm{cut}}
:=\sum_{\kappa\in\mathbb Z^3}F_{q,j+1,u,\kappa}^{\mathrm{cut}}.
$$
Each localized cutoff source also has zero spatial integral by the
mean correction. Apply Theorem~\ref{thm:inv4divEq} to each localized source and sum
in $\kappa$ to obtain matrices $S_{q,j+1,u}^{\mathrm{cut}}$.
The local right inverse does not enlarge temporal support, so
$\supp_tS_{q,j+1,u}^{\mathrm{cut}}\subset\supp\partial_t\theta_u$. They satisfy
$$
F_{q,j+1,u}^{\mathrm{cut}}
=\Div S_{q,j+1,u}^{\mathrm{cut}}.
$$
These matrices need not be symmetric. For each $u$, we therefore write
$$
S_{q,j+1,u}^{\mathrm{cut}}
=S_{q,j+1,u}^{\mathrm{cut,sym}}
+S_{q,j+1,u}^{\mathrm{cut,skew}},
\qquad
(S_{q,j+1,u}^{\mathrm{cut,sym}})^T
=S_{q,j+1,u}^{\mathrm{cut,sym}},
\qquad
(S_{q,j+1,u}^{\mathrm{cut,skew}})^T
=-S_{q,j+1,u}^{\mathrm{cut,skew}}.
$$
The sum over the time blocks satisfies
$$
F_{q,j+1}^{\mathrm{cut}}
=\sum_u\Div S_{q,j+1,u}^{\mathrm{cut,sym}}
+\sum_u\Div S_{q,j+1,u}^{\mathrm{cut,skew}}.
$$
The symmetric matrices are included in the gluing stress
$R_{q,j+1}^{\mathrm{gluing}}$. The incompressible correction cancels the divergences of
the antisymmetric matrices.

For any periodic antisymmetric matrix $A$, antisymmetry and
\eqref{def:periodic-Hodge-potentials} give
$$
\Div\Div A=0,
\qquad
Q[A]=0,
\qquad
\Div A=\curl W[A].
$$
The incompressible perturbation changes only the momentum.
For each time block $u$, let $(\omega_{q,j,u},\pi_{q,j,u})$ solve
the linearized system with zero initial data
$$
\left\{
\begin{aligned}
	&\partial_t\omega_{q,j,u}
	+\Div\left(\overline{u}_q\otimes\omega_{q,j,u}
	+\omega_{q,j,u}\otimes\overline{u}_q\right)
	+\nabla\pi_{q,j,u}\hspace{-2pt}=\hspace{-2pt}-\theta_u^2\curl W_{q,j}^{\mathrm{mol}}
	-\curl W\left[S_{q,j+1,u}^{\mathrm{cut,skew}}\right],\\
	&\Div\omega_{q,j,u}=0,\\
	&\omega_{q,j,u}|_{t=t_u}=0,
	\qquad
	\fint_{\T^3}\pi_{q,j,u}\,\rd x=0.
\end{aligned}
\right.
$$
The background and forcing are smooth. Since $\theta_u$ depends only on time,
the forcing is a periodic curl and is therefore divergence-free with zero spatial mean.
Apply Proposition~\ref{prop:linearized-equation} with
$I=[t_u-\tau_q,t_u+\tau_q]$, $U=\overline u_q$,
$F=-\theta_u^2\curl W_{q,j}^{\mathrm{mol}}-\curl W[S_{q,j+1,u}^{\mathrm{cut,skew}}]$,
$t_0=t_u$, and $\omega_0=0$. We obtain existence and uniqueness on $I$,
which contains $\supp\tilde\theta_u$. Smoothness follows from
Remark~\ref{rem:smooth-data-smooth-solution}.
Integrating the equation in space and using the zero initial data shows
that $\omega_{q,j,u}$ has zero spatial mean.

With $\mathcal B$ defined in \eqref{eq:vector-potential-operator}, let
$
\Psi_{q,j,u}:=\mathcal B\omega_{q,j,u},
\curl\Psi_{q,j,u}=\omega_{q,j,u},
$
and define the glued momentum perturbation and the corresponding pressure by
$$
	\tm_{q,j+1}\inc:=\sum_u\tilde\theta_u\omega_{q,j,u},
\qquad
\pi_{q,j+1}\inc:=\sum_u\tilde\theta_u\pi_{q,j,u}.
$$
Multiplying the local equations by $\tilde\theta_u$ and summing in $u$, we obtain
$$
\left\{
\begin{aligned}
	&\partial_t\tm_{q,j+1}\inc
	+\Div\left(\overline{u}_q\otimes\tm_{q,j+1}\inc
	+\tm_{q,j+1}\inc\otimes\overline{u}_q\right)
	+\!\nabla\pi_{q,j+1}\inc=\curl\left(
	\sum_u(\partial_t\tilde\theta_u)\Psi_{q,j,u}
	\right)-\!\curl W_{q,j}^{\mathrm{mol}}-\!\sum_u\curl W\left[
	S_{q,j+1,u}^{\mathrm{cut,skew}}
	\right],\\
	&\Div\tm_{q,j+1}\inc=0,\\
	&\fint_{\T^3}\tm_{q,j+1}\inc\,\rd x=0,\qquad
	\fint_{\T^3}\pi_{q,j+1}\inc\,\rd x=0.
\end{aligned}
\right.
$$
Here we used $\sum_u\theta_u^2=1$ and
$\tilde\theta_u\theta_u^2=\theta_u^2$.
Since $W$ acts only in space,
$W[S_{q,j+1,u}^{\mathrm{cut,skew}}]$ has temporal support in
$\supp\partial_t\theta_u$, where $\tilde\theta_u=1$.
The terms $-\curl W[S_{q,j+1,u}^{\mathrm{cut,skew}}]$ cancel the divergences of the antisymmetric cutoff
matrices in the compressible equation.
The term containing $\partial_t\tilde\theta_u$ is the incompressible
gluing residual.
Since $\Div\tm_{q,j+1}\inc=0$, adding this momentum correction
preserves the continuity equation. The gluing residual is a curl
and therefore has zero spatial mean. Applying the periodic symmetric
right inverse $\mathcal R$ of $\Div$ to it and adding the symmetric
cutoff matrices gives the gluing stress $R_{q,j+1}^{\mathrm{gluing}}$ defined by
$$
\sum_uS_{q,j+1,u}^{\mathrm{cut,sym}}
+\mathcal R\curl\left(
\sum_u(\partial_t\tilde\theta_u)\Psi_{q,j,u}
\right).
$$
We add this stress to $S_{q,j}$ to obtain $S_{q,j+1}$.

The pressure contribution on the right-hand side satisfies
$
-\nabla\pi_{q,j+1}\inc
=\Div\left(-\pi_{q,j+1}\inc\Id\right),
$
so we define the next interior stress by
\begin{equation}\label{def:Newton-unsmoothed-interior-output}
	\begin{aligned}
		R_{q,j+1}:=R_{q,j+1}\comp-\pi_{q,j+1}\inc\Id.
	\end{aligned}
\end{equation}
The interior stress $R_{q,j+1}$ is used in the next Newton substep when
$j+1<J_*$. At $j+1=J_*$, the terminal interior stress $R_{q,J_*}$ is included
in $R_{q+1}$.

We update the density and momentum by
$$
\rho_{q,j+1}:=\rho_{q,j}+\trho_{q,j+1}\comp,\qquad m_{q,j+1}:=m_{q,j}+\tm_{q,j+1}\comp+\tm_{q,j+1}\inc.
$$
From the exact changes in the convective and pressure fluxes, subtract the terms
linear in the new perturbations, with coefficients evaluated at $(\orho_q,\om_q)$.
Together with the two source mollification stresses, this remainder defines
$P_{q+1,j+1}-P_{q+1,j}$ in \eqref{def:Newton-perturbative-tensor}.
Adding the correction equations to \eqref{eq:Newton-substep-Euler-Reynolds-system}
gives the system at substep $j+1$.

Along-flow mollification enlarges temporal support by at most
$\ell_{t,q}$. The local right inverse of $\Div$ and the periodic
operators $Q$, $W$, $\mathcal B$, and $\mathcal R$ act at each fixed time, so they do not enlarge temporal support. For $0\leq j<J_*$, a correction on block $u$ can be nonzero only if
$\supp\theta_u\cap\mN_{\ell_{t,q}}(\supp_tR_{q,j})\neq\varnothing$. The cutoff supports and the definitions of
the source mollification stresses then give
\begin{equation}\label{eq:Newton-substep-temporal-support}
	\begin{aligned}
		&\supp_t\trho_{q,j+1}\comp\cup\supp_t\tm_{q,j+1}\comp\cup\supp_t\tm_{q,j+1}\inc\cup\supp_tR_{q,j+1}\cup\supp_tR_{q,j+1}^{\mathrm{gluing}}\subset\mN_{2\tau_q}(\supp_tR_{q,j}),\\
		&\supp_tP_{q+1,j+1}^{\mathrm{mol,comp}}\cup\supp_tP_{q+1,j+1}^{\mathrm{mol,inc}}\subset\mN_{\ell_{t,q}}(\supp_tR_{q,j}).
	\end{aligned}
\end{equation}
Indeed, on a block satisfying the intersection condition, the distance
from any point of $\supp\tilde\theta_u$ to $\supp_tR_{q,j}$ is at most
$(17/24+2/3)\tau_q+\ell_{t,q}<2\tau_q$, since
$\ell_{t,q}<\tau_q/6$. The source mollification stresses vanish
outside $\mN_{\ell_{t,q}}(\supp_tR_{q,j})$ by
\eqref{def:Newton-substep-flow-mollification}.

Since $P_{q+1,0}=R_q-R_{q,0}$, \eqref{eq:initial-Newton-flow-mollified-support} gives
$\supp_tR_{q,0}\cup\supp_tP_{q+1,0}\subset\mN_{\ell_{t,q}}(\supp_tR_q)$.
Using $\ell_{t,q}<\tau_q$ and iterating
\eqref{eq:Newton-substep-temporal-support}, we obtain
\begin{equation}\label{eq:Newton-iterated-temporal-support}
	\begin{aligned}
		\supp_tR_{q,j}&\subset\mN_{(2j+1)\tau_q}(\supp_tR_q),\qquad 0\leq j\leq J_*,\\
		\supp_tP_{q+1,j}^{\mathrm{mol,comp}}\cup\supp_tP_{q+1,j}^{\mathrm{mol,inc}}&\subset\mN_{(2j+1)\tau_q}(\supp_tR_q),\qquad 1\leq j\leq J_*.
	\end{aligned}
\end{equation}

The flux remainders in \eqref{def:Newton-perturbative-tensor} vanish outside the temporal supports of
$\trho_{q,j+1}\comp$, $\tm_{q,j+1}\comp$, and $\tm_{q,j+1}\inc$. Combining this with
\eqref{eq:Newton-substep-temporal-support}--%
\eqref{eq:Newton-iterated-temporal-support}, we have, for $0\leq j<J_*$,
\begin{equation}\label{eq:Newton-perturbative-increment-temporal-support}
	\supp_t(P_{q+1,j+1}-P_{q+1,j})\subset\mN_{(2j+3)\tau_q}(\supp_tR_q).
\end{equation}
Using \eqref{eq:Newton-substep-temporal-support}--\eqref{eq:Newton-perturbative-increment-temporal-support}, sum the perturbations, gluing stresses, and differences $P_{q+1,\nu+1}-P_{q+1,\nu}$ over $0\leq\nu<j$.
The initial values then give
\begin{equation}\label{eq:Newton-partial-state-stress-temporal-support}
	\begin{aligned}
		\supp_t(\rho_{q,j}-\rho_q)\cup\supp_t(m_{q,j}-m_q)\cup\supp_tS_{q,j}\cup\supp_tP_{q+1,j}\subset\mN_{(2j+1)\tau_q}(\supp_tR_q),\qquad 0\leq j\leq J_*.
	\end{aligned}
\end{equation}
We now state the substep estimates and fix the remaining parameters. Take
$$
N_R\geq4J_*+12,\qquad
N_{\mathrm{ind}}\geq N_R+4J_*+12,\qquad
N_{\mathrm{max}}\geq N_{\mathrm{ind}}+4J_*+12.
$$
We also require $N_{\mathrm{ind}}$ to satisfy
\begin{equation}\label{eq:Newton-top-material-reserve}
	(b-1)\left(\frac12-\beta\right)(N_{\mathrm{ind}}-1)>(2b-1)\beta.
\end{equation}
Set
\begin{equation}\label{def:Newton-derivative-bookkeeping}
	\begin{aligned}
		K_j^{(s)}:=N_{\mathrm{max}}-3j-2-s,
		&& 
		N_{R,j}^{(s)}(M):=\max\{N_R-3j-M-s,0\},
		&& 0\leq s\leq3.
	\end{aligned}
\end{equation}
Here $M+N\leq K_j^{(s)}$ is the total derivative-order range. The index $N_{R,j}^{(s)}(M)$ is the threshold
in the spatial factor $\MM{N,N_{R,j}^{(s)}(M),\lambda_q,\ell_q^{-1}}$.
For $r\in\mathbb N_0$, set $\alpha_r:=2^{-r}\alpha$, so that $2\alpha_{r+1}=\alpha_r$. The number of substeps and the derivative ranges are fixed, so only finitely many exponents are needed. The background estimates in Proposition~\ref{prop:uniform-mollification-estimates} also hold at these fixed exponents.

The implicit constants in this section are independent of $q$ and $a$.
Unless stated otherwise, they may depend on
$\rho^*,\mathfrak M,N_{\mathrm{max}},\alpha,\alpha_r$.
Fix these parameters and then choose $C_*$.
Choose $a_{\mathrm{Newton}}=a_{\mathrm{Newton}}(\beta,\allowbreak b,\allowbreak \alpha,\allowbreak \rho^*,\allowbreak \mathfrak M,\allowbreak C_*,\allowbreak N_R,\allowbreak N_{\mathrm{ind}},\allowbreak N_{\mathrm{max}})\geq a_{\mathrm{loc}}$
large enough for all estimates and support conditions in this section.
Throughout this section, assume $a>a_{\mathrm{Newton}}$.

By \eqref{eq:material-mollification-spatial-comparisons}, for $K\in\mathbb N$ and $M,N\in\mathbb N_0$,
\begin{equation}\label{eq:Newton-one-threshold-weight-shift}
	\begin{aligned}
		\frac{\MM{M,K,\tau_q^{-1},\ell_{t,q}^{-1}}
			\MM{N+1,(K-M)_+,\lambda_q,\ell_q^{-1}}}
		{\MM{M,K-1,\tau_q^{-1},\ell_{t,q}^{-1}}
			\MM{N,(K-1-M)_+,\lambda_q,\ell_q^{-1}}}=
		\begin{cases}
			\lambda_q,&M<K,\\
			\tau_q^{-1}\ell_{t,q}\ell_q^{-1},&M\geq K,
		\end{cases}
		\leq\lambda_q.
	\end{aligned}
\end{equation}
For $a>a_{\mathrm{Newton}}$, $\ell_{t,q}<\tau_q/6$ and
\begin{equation}\label{eq:Newton-scale-separation}
	\tau_q^{-1}\leq\lambda_q\leq\ell_q^{-1},\qquad
	\tau_q\delta_q^{\frac12}\lambda_q^{1+\alpha}\leq\lambda_{q+1}^{-3\alpha},\qquad
	\tau_q\Lambda_q\leq\lambda_{q+1}^{-4\alpha}.
\end{equation}
The first relation implies $\tau_q\ell_q^{-1}\geq1$. The last two estimates follow from \eqref{def:tau-q}.

At substep $0\leq j<J_*$, assume that
\begin{equation}\label{eq:Newton-interior-input-estimate}
	\begin{aligned}
		\nrm{\DTL^MR_{q,j}}_{N+\alpha}\leq A_{q,j}
		\MM{M,N_R-3j,\tau_q^{-1},\ell_{t,q}^{-1}}
		\La{N,N_{R,j}^{(0)}(M),\lambda_q,\ell_q^{-1}}
		\MM{N,N_{R,j}^{(0)}(M),\lambda_q,\ell_q^{-1}},
	\end{aligned}
\end{equation}
for $
0\leq M\leq N_{\mathrm{ind}}-1
$ and $
M+N\leq K_j^{(0)}=N_{\mathrm{max}}-3j-2,
$ where the stress amplitudes $A_{q,j}>0$ will be fixed below. We will prove that the new interior stress satisfies
\begin{equation}\label{eq:Newton-final-interior-stress-estimate}
	\begin{aligned}
		&\nrm{\DTL^MR_{q,j+1}}_{N+\alpha}\\
		&\ \ \lesssim\tau_q\Lambda_q
		(\tau_q\ell_q^{-1})^\alpha A_{q,j}\MM{M,N_R-3(j+1),\tau_q^{-1},\ell_{t,q}^{-1}}
		\Lat{N,N_{R,j+1}^{(0)}(M),\lambda_q,\ell_q^{-1}}
		\MM{N,N_{R,j+1}^{(0)}(M),\lambda_q,\ell_q^{-1}},
	\end{aligned}
\end{equation}
for $
0\leq M\leq N_{\mathrm{ind}}-1
$ and $
M+N\leq K_j^{(3)}=K_{j+1}^{(0)}.
$
Choose $C_*>1$ to dominate the constants in \eqref{eq:initial-Newton-flow-mollified-stress-estimate} and
\eqref{eq:Newton-final-interior-stress-estimate}, for all substeps, stated derivative orders, and fixed H{\"o}lder exponents. The uniformity of these estimates in $q$ and $a$ allows $C_*$ to be chosen independently of both. Define the amplitudes recursively by
$$
A_{q,0}:=C_*\lambda_q^{-8\alpha}\delta_{q+1},
\qquad
A_{q,j+1}:=C_*\tau_q\Lambda_q
(\tau_q\ell_q^{-1})^\alpha\ell_q^{-\alpha}A_{q,j}.
$$
Use $\La{N,N_{R,j+1}^{(0)}(M),\lambda_q,\ell_q^{-1}}\leq\ell_q^{-\alpha}$ in \eqref{eq:Newton-final-interior-stress-estimate}. The amplitude recursion then gives
\eqref{eq:Newton-interior-input-estimate} at substep $j+1$. At $j=0$, this estimate follows from
\eqref{eq:initial-Newton-flow-mollified-stress-estimate}. At each smaller fixed H{\"o}lder exponent, we repeat the source and perturbation estimates. The same amplitude recursion applies, with the parameters unchanged, because the additional H{\"o}lder and localization factors are no larger than those for H{\"o}lder exponent $\alpha$. Each substep preserves the range of $M$ and lowers the maximum total derivative order by three.

Iterating the recursion, we have
$$
A_{q,j}=A_{q,0}\left[
C_*\tau_q\Lambda_q
(\tau_q\ell_q^{-1})^\alpha\ell_q^{-\alpha}
\right]^j.
$$
By \eqref{def:parameter-mollification} and \eqref{def:tau-q}, $
	\tau_q\ell_q^{-1}
\sim
\lambda_q^{\frac{b-1}{4}+\beta-4b\alpha}\leq\lambda_q^{\frac32\beta}
$ and
$$
	C_*\tau_q\Lambda_q
(\tau_q\ell_q^{-1})^\alpha\ell_q^{-\alpha}
\lesssim_{C_*}
\lambda_q^{-\alpha\left(\frac{7b-1}{2}-\beta+4b\alpha\right)}\leq\lambda_q^{-\frac83\alpha}\leq1.
$$
Since the exponent is negative, take $a_{\mathrm{Newton}}$ sufficiently large to absorb the implicit constant.

At the end of the Newton step, the amplitude formula and the estimate for the recursive factor, together with
$\lambda_{q+1}\sim\lambda_q^b$, $\delta_q=\lambda_q^{-2\beta}$, and
$\ell_q^{-1}\sim\lambda_q^{\frac{3+b}{4}}$, give
\begin{equation}\label{eq:Newton-terminal-amplitude-ratio}
\frac{A_{q,J_*}\ell_q^{-\alpha}}
{\lambda_{q+1}^{-8\alpha}\delta_{q+2}}
\lesssim
\lambda_q^{-J_*\alpha\left(\frac{7b-1}{2}-\beta+4b\alpha\right)
	+(b-1)(2b\beta+8\alpha)+\frac{3+b}{4}\alpha}\leq\lambda_q^{-\alpha}.
\end{equation}
By \eqref{def:admissible-alpha-losses},
$
2b\beta+8\alpha<\beta+\frac13<\frac23,
\
\frac{7b-1}{2}-\beta+4b\alpha>\frac83>\frac{7+b}{4}.
$
Since $J_*\alpha\geq(b-1)/4+\alpha$, we obtain
$$
J_*\alpha\left(\frac{7b-1}{2}-\beta+4b\alpha\right)
>(b-1)(2b\beta+8\alpha)+\frac{7+b}{4}\alpha.
$$
The exponent in \eqref{eq:Newton-terminal-amplitude-ratio} is less than $-\alpha$. For sufficiently large $a$, we have
\begin{equation}\label{eq:Newton-terminal-amplitude-bound}
	A_{q,J_*}\ell_q^{-\alpha}
	\leq\lambda_q^{-\alpha}\lambda_{q+1}^{-8\alpha}\delta_{q+2}.
\end{equation}

Under \eqref{eq:Newton-interior-input-estimate} at the fixed H{\"o}lder exponents, the scalar and vector potentials satisfy, for $0\leq j<J_*$, $0\leq M\leq N_{\mathrm{ind}}-1$, and $M+N\leq K_j^{(0)}$,
$$
\begin{aligned}
	&\nrm{\DTL^MQ[R_{q,j}]}_{N+\alpha}
	+\nrm{\DTL^MW[R_{q,j}]}_{N+\alpha}\\&\quad\lesssim A_{q,j}
	\MM{M,N_R-3j,\tau_q^{-1},\ell_{t,q}^{-1}}
	\Lat{N,N_{R,j}^{(0)}(M),\lambda_q,\ell_q^{-1}}
	\MM{N,N_{R,j}^{(0)}(M),\lambda_q,\ell_q^{-1}}.
\end{aligned}
$$
For $0\leq j<J_*$, the mollified sources defined in \eqref{def:Newton-sharp-acoustic-input}--\eqref{def:Newton-flow-mollified-vector-source} satisfy
$$
\begin{aligned}
	&\nrm{\DTL^MG_{q,j,u,\kappa}^{\mathrm{mol}}}_{N+\alpha}+\nrm{\DTL^MW_{q,j}^{\mathrm{mol}}}_{N+\alpha}\\
	&\quad\lesssim A_{q,j}\MM{M,N_R-3j,\tau_q^{-1},\ell_{t,q}^{-1}}\Lat{N,N_{R,j}^{(0)}(M),\lambda_q,\ell_q^{-1}}\MM{N,N_{R,j}^{(0)}(M),\lambda_q,\ell_q^{-1}},
\end{aligned}
$$
for $0\leq M\leq N_{\mathrm{max}}$ and $\min\{M,N_R-3j\}+N\leq K_j^{(0)}$. Along-flow mollification provides estimates for the higher material derivatives needed in the compressible construction.
The estimates for these potentials and the mollified sources are proved in Lemma~\ref{lem:Newton-source-flow-mollification}.

The compressible perturbation and stress estimates, with their derivative ranges, are stated in
Proposition~\ref{prop:Newton-three-dimensional-acoustic-estimates}.
The incompressible correction and gluing-stress estimates are stated in
Proposition~\ref{prop:Newton-incompressible-flow-gluing-estimates}.

For $0\leq j<J_*$, $M\leq N_{\mathrm{ind}}-1$, and $M+N\leq K_j^{(0)}$, set $L:=N_{R,j}^{(0)}(M+N)$. The two source mollification stresses satisfy
$$
\begin{aligned}
	\nrm{\DTL^MP_{q+1,j+1}^{\mathrm{mol,comp}}}_N+\nrm{\DTL^MP_{q+1,j+1}^{\mathrm{mol,inc}}}_N
	&\lesssim A_{q,j}\ell_q^{-\alpha}(\ell_{t,q}\tau_q^{-1})^L
	\MM{M,N_R-3j,\tau_q^{-1},\ell_{t,q}^{-1}}\ell_q^{-N}.
\end{aligned}
$$
Estimates~\eqref{eq:Newton-substep-mollification-remainder-estimate} and \eqref{eq:Newton-source-mollification-two-scale-estimate} imply, in the same range,
$$
\begin{aligned}
	\nrm{\DTL^MP_{q+1,j+1}^{\mathrm{mol,comp}}}_N+\nrm{\DTL^MP_{q+1,j+1}^{\mathrm{mol,inc}}}_N
	&\lesssim A_{q,j}\ell_q^{-\alpha}(\ell_q\lambda_q)^{\frac{N_R-3j}{2}}\ell_{t,q}^{-M}\ell_q^{-N}.
\end{aligned}
$$
To obtain a factor $\lambda_q^{-\alpha}$ in the estimates for $P_{q+1,j}$, we use
the last inequality in \eqref{def:admissible-alpha-losses}. The estimates for
$R_{q,j+1}$ and the source mollification stresses are proved in Subsection~\ref{subsec:Newton-stress-state-estimates}. The estimates for
$P_{q+1,j}$ along the background flow are proved in Proposition~\ref{prop:Newton-background-flow-stresses}.

\subsection{The compressible perturbations}
\label{subsec:Newton-compressible-perturbations}

\subsubsection{Kirchhoff's formula in Lagrangian coordinates}

Fix $0\leq j<J_*$. We construct $\widetilde V_{q,j+1,u,\kappa}$ in flow coordinates, then correct its spatial mean and localize it in time. Let $X_u(s,y)$ be the flow of $\overline u_q$ based at $t_u$, and let $\Phi_u(t,\cdot)$ be
its spatial inverse. On the block $(u,\kappa)$, we freeze the sound speed at $(t_u,x_\kappa)$ by setting
$$
c_{q,u,\kappa}^2:=p'(\orho_q)(t_u,x_\kappa),
\qquad
c_{q,u,\kappa}:=\sqrt{p'(\orho_q)(t_u,x_\kappa)}.
$$
Since $\orho_q\in I_{\rho^*}$ and $p'>0$ on $I_{\rho^*}$, $c_{q,u,\kappa}$ is uniformly bounded above and away from zero. For $\xi\in\mathbb S^2$, define
$$
Y_{q,u,\kappa}(t,s,x,\xi):=X_u\left(s,\Phi_u(t,x)+c_{q,u,\kappa}(t-s)\xi\right).
$$
We construct the vector potential from a constant-coefficient wave equation in Lagrangian coordinates, with sound speed $c_{q,u,\kappa}$ and zero initial data at
$t_{u-1}$. In Eulerian coordinates, its solution operator is the three-dimensional Kirchhoff operator
\begin{equation}\label{def:Newton-three-dimensional-cone-operator}
	\begin{aligned}
		I_{q,u,\kappa}[H](t,x)
		&:=\int_{t_{u-1}}^t\frac{t-s}{4\pi}
		\int_{\mathbb S^2}
		H(s,Y_{q,u,\kappa}(t,s,x,\xi))\,\rd S_\xi\,\rd s.
	\end{aligned}
\end{equation}
In Lagrangian coordinates, we write
\begin{equation}\label{def:Newton-Lagrangian-cone-operator}
	\begin{aligned}
		\widehat I_{q,u,\kappa}[H](t,y)
		&:=\int_{t_{u-1}}^t\frac{t-s}{4\pi}
		\int_{\mathbb S^2}
		H\left(s,X_u\left(s,y+c_{q,u,\kappa}(t-s)\xi\right)\right)
		\,\rd S_\xi\,\rd s.
	\end{aligned}
\end{equation}
Thus
$
I_{q,u,\kappa}[H](t,x)=\widehat I_{q,u,\kappa}[H](t,\Phi_u(t,x)).
$
The factor $\theta_u$ makes $G_{q,j,u,\kappa}^{\mathrm{mol}}$ and all its material derivatives vanish near the lower endpoint $t_{u-1}$. Thus,
integration by parts in time produces no boundary terms at $t_{u-1}$.

\begin{lem}\label{lem:Newton-three-dimensional-cone-identities}
	Fix $M\in\mathbb N$ and let $H$ be smooth. All statements below concern $t\geq t_{u-1}$. Assume that $
	(\DTL^mH)(t_{u-1},\cdot)=0
	$
	for every $0\leq m\leq M$. Then, for $1\leq r\leq M+1$,
	\begin{equation}\label{eq:Newton-cone-material-identities}
		\DTL^rI_{q,u,\kappa}[H]=I_{q,u,\kappa}[\DTL^rH],
	\end{equation}
	and
	\begin{equation}\label{eq:Newton-Kirchhoff-identity}
		c_{q,u,\kappa}^2\bigl(\Delta_y\widehat I_{q,u,\kappa}[H]\bigr)(t,\Phi_u(t,x))
		=I_{q,u,\kappa}[\DTL^2H](t,x)-H(t,x).
	\end{equation}
	If $H$ is scalar, then
	\begin{equation}\label{eq:Newton-cone-gradient-no-loss}
		\begin{aligned}
			\nabla_y\widehat I_{q,u,\kappa}[H](t,y)
			=\frac{1}{4\pi c_{q,u,\kappa}}
			\int_{t_{u-1}}^t\int_{\mathbb S^2}\xi
			\left[H(s,Y_y)+(t-s)(\DTL H)(s,Y_y)\right]\,\rd S_\xi\,\rd s,
		\end{aligned}
	\end{equation}
	where $
	Y_y=X_u(s,y+c_{q,u,\kappa}(t-s)\xi).
	$ For every $0\leq m\leq M$,
	\begin{equation}\label{eq:Newton-commuted-cone-gradient}
		\begin{aligned}
			\partial_t^m\nabla_y\widehat I_{q,u,\kappa}[H](t,y)
			=\frac{1}{4\pi c_{q,u,\kappa}}
			\int_{t_{u-1}}^t\int_{\mathbb S^2}\xi
			\left[
			(\DTL^mH)(s,Y_y)
			+(t-s)(\DTL^{m+1}H)(s,Y_y)
			\right]\,\rd S_\xi\,\rd s.
		\end{aligned}
	\end{equation}
	For every $0\leq m\leq M$, $N\in\mathbb N_0$, and $0<\alpha<1$,
	\begin{equation}\label{eq:Newton-Lagrangian-cone-holder-estimates}
		\begin{aligned}
			\nrm{\partial_t^m\widehat I_{q,u,\kappa}[H](t,\cdot)}_{N+\alpha}
			&\lesssim
			(t-t_{u-1})^2
			\sup_{t_{u-1}\leq s\leq t}
			\nrm{(\DTL^mH)(s,X_u(s,\cdot))}_{N+\alpha},\\
			\nrm{\partial_t^m\nabla_y\widehat I_{q,u,\kappa}[H](t,\cdot)}_{N+\alpha}
			&\lesssim
			c_{q,u,\kappa}^{-1}(t-t_{u-1})^2
			\sup_{t_{u-1}\leq s\leq t}
			\nrm{(\DTL^{m+1}H)(s,X_u(s,\cdot))}_{N+\alpha}.
		\end{aligned}
	\end{equation}
	
	Suppose, in addition, that for some $\bar\tau,\bar\varepsilon>0$ and $x_0\in\mathbb R^3$,
	\begin{equation}\label{eq:Newton-cone-input-support}
		\supp H
		\subseteq
		\hat\Gamma_{\overline u_q}
		(\bar\tau,\bar\varepsilon;t_u,x_0).
	\end{equation}
	Then, for every $t\geq t_{u-1}$,
	\begin{equation}\label{eq:Newton-cone-slice-support}
		\begin{aligned}
			\supp_y\widehat I_{q,u,\kappa}[H](t,\cdot)
			&\subseteq
			\overline B\left(
			x_0,\bar\varepsilon
			+c_{q,u,\kappa}(t-t_{u-1})
			\right),\\
			\supp_xI_{q,u,\kappa}[H](t,\cdot)
			&\subseteq
			X_u\left(
			t,\overline B\left(
			x_0,\bar\varepsilon
			+c_{q,u,\kappa}(t-t_{u-1})
			\right)
			\right).
		\end{aligned}
	\end{equation}
	For $\bar\tau\geq2\tau_q/3$, the time-localized potential satisfies
	\begin{equation}\label{eq:Newton-cone-cylinder-support}
		\supp\left(\theta_u I_{q,u,\kappa}[H]\right)
		\subseteq
		\hat\Gamma_{\overline u_q}\left(
		\bar\tau,
		\bar\varepsilon+2c_{q,u,\kappa}\bar\tau;
		t_u,x_0
		\right).
	\end{equation}
	Here the product is extended by zero for $t<t_{u-1}$.
	For every fixed $K\in\mathbb N_0$, the constants in \eqref{eq:Newton-Lagrangian-cone-holder-estimates} for $m+N\leq K$ depend only on $K,\alpha$.

\end{lem}

\begin{proof}
	We first compute $\DTL Y_{q,u,\kappa}$ along the integration trajectories. Since $\DTL\Phi_u=0$,
	$$
	\DTL Y_{q,u,\kappa}=c_{q,u,\kappa}(\nabla_yX_u)\left(s,\Phi_u(t,x)+c_{q,u,\kappa}(t-s)\xi\right)\xi,
	$$
	whereas $
	\partial_sY_{q,u,\kappa}=\overline{u}_q(s,Y_{q,u,\kappa})-c_{q,u,\kappa}(\nabla_yX_u)\left(s,\Phi_u(t,x)+c_{q,u,\kappa}(t-s)\xi\right)\xi.
	$ Hence
	$
	\DTL Y_{q,u,\kappa}=\overline{u}_q(s,Y_{q,u,\kappa})-\partial_sY_{q,u,\kappa}.
	$
	For any scalar or vector field $H$,
	\begin{equation}\label{eq:Newton-material-chain-rule-on-cone}
		\DTL\bigl(H(s,Y_{q,u,\kappa})\bigr)=(\DTL H)(s,Y_{q,u,\kappa})-\partial_s\bigl(H(s,Y_{q,u,\kappa})\bigr),
	\end{equation}
	where $\DTL H$ is taken before composition with $Y_{q,u,\kappa}$.
	From \eqref{def:Newton-three-dimensional-cone-operator} and \eqref{eq:Newton-material-chain-rule-on-cone},
	$$
	\begin{aligned}
		\DTL I_{q,u,\kappa}[H]
		=&\int_{t_{u-1}}^t\frac{1}{4\pi}\int_{\mathbb S^2}H(s,Y_{q,u,\kappa}(t,s,x,\xi))\,\rd S_\xi\,\rd s\\
		&+\int_{t_{u-1}}^t\frac{t-s}{4\pi}\int_{\mathbb S^2}\left[(\DTL H)(s,Y_{q,u,\kappa}(t,s,x,\xi))-\partial_s(H(s,Y_{q,u,\kappa}(t,s,x,\xi)))\right]\,\rd S_\xi\,\rd s.
	\end{aligned}
	$$
	Integration by parts in $s$ gives
	$$
	\begin{aligned}
		&-\int_{t_{u-1}}^t\frac{t-s}{4\pi}\int_{\mathbb S^2}\partial_s(H(s,Y_{q,u,\kappa}(t,s,x,\xi)))\,\rd S_\xi\,\rd s\\
		=&-\left[\frac{t-s}{4\pi}\int_{\mathbb S^2}H(s,Y_{q,u,\kappa}(t,s,x,\xi))\,\rd S_\xi\right]_{s=t_{u-1}}^{s=t}-\int_{t_{u-1}}^t\frac{1}{4\pi}\int_{\mathbb S^2}H(s,Y_{q,u,\kappa}(t,s,x,\xi))\,\rd S_\xi\,\rd s.
	\end{aligned}
	$$
	The boundary terms vanish since $t-s=0$ at $s=t$ and $H(t_{u-1},\cdot)=0$. The remaining integral cancels the first term in
	the formula for $\DTL I_{q,u,\kappa}[H]$. Hence $
	\DTL I_{q,u,\kappa}[H]=I_{q,u,\kappa}[\DTL H].
	$ Applying this identity successively to $H,\DTL H,\ldots,\DTL^MH$ proves
	\eqref{eq:Newton-cone-material-identities}.
	
	Since $
	I_{q,u,\kappa}[H](t,x)=\widehat I_{q,u,\kappa}[H](t,\Phi_u(t,x))
	$ and $\DTL\Phi_u=0$,
	$$
	\DTL^rI_{q,u,\kappa}[H](t,x)
	=
	\bigl(\partial_t^r\widehat I_{q,u,\kappa}[H]\bigr)(t,\Phi_u(t,x)).
	$$
	Equation \eqref{eq:Newton-cone-material-identities} and the invertibility of $\Phi_u(t,\cdot)$ imply
	\begin{equation}\label{eq:Newton-Lagrangian-material-commutation}
		\partial_t^r\widehat I_{q,u,\kappa}[H](t,y)
		=\widehat I_{q,u,\kappa}[\DTL^rH](t,y),
		\qquad
		0\leq r\leq M+1.
	\end{equation}
	In particular,
	$$
	\DTL^m\left(
	\nabla_y\widehat I_{q,u,\kappa}[H](t,\Phi_u(t,x))
	\right)
	=
	\left(
	\partial_t^m\nabla_y\widehat I_{q,u,\kappa}[H]
	\right)(t,\Phi_u(t,x)).
	$$
	
	To prove \eqref{eq:Newton-Kirchhoff-identity}, set
	$$
	h_s(y):=H(s,X_u(s,y)),
	\qquad
	M_\sigma h_s(y):=\frac{1}{4\pi}
	\int_{\mathbb S^2}h_s(y+\sigma\xi)\,\rd S_\xi.
	$$
	The three-dimensional spherical-mean identity is
	$$
	\partial_\sigma^2\bigl(\sigma M_\sigma h_s(y)\bigr)
	=\sigma M_\sigma(\Delta_yh_s)(y).
	$$
	Since $c_{q,u,\kappa}$ is constant on the block,
	$$
	\left(\partial_\sigma^2-c_{q,u,\kappa}^2\Delta_y\right)
	\bigl(\sigma M_{c_{q,u,\kappa}\sigma}h_s(y)\bigr)=0.
	$$
	At $\sigma=0$, the function $\sigma M_{c_{q,u,\kappa}\sigma}h_s(y)$ vanishes and its $\sigma$-derivative equals $h_s(y)$. Moreover,
	$$
	\widehat I_{q,u,\kappa}[H](t,y)
	=\int_{t_{u-1}}^t(t-s)M_{c_{q,u,\kappa}(t-s)}h_s(y)\,\rd s.
	$$
	Therefore,
	\begin{align*}
	\partial_t\widehat I_{q,u,\kappa}[H](t,y)
	&=\int_{t_{u-1}}^t
	\left.\partial_\sigma\bigl(\sigma M_{c_{q,u,\kappa}\sigma}h_s(y)\bigr)
	\right|_{\sigma=t-s}\,\rd s,\\
	\partial_t^2\widehat I_{q,u,\kappa}[H](t,y)
	&=h_t(y)+\int_{t_{u-1}}^t
	\left.\partial_\sigma^2\bigl(\sigma M_{c_{q,u,\kappa}\sigma}h_s(y)\bigr)
	\right|_{\sigma=t-s}\,\rd s.
	\end{align*}
	It follows that
	$$
	\left(\partial_t^2-c_{q,u,\kappa}^2\Delta_y\right)\widehat I_{q,u,\kappa}[H](t,y)=h_t(y).
	$$
	The commutation identity \eqref{eq:Newton-Lagrangian-material-commutation} then gives
	$$
	c_{q,u,\kappa}^2\bigl(\Delta_y\widehat I_{q,u,\kappa}[H]\bigr)(t,\Phi_u(t,x))
	=I_{q,u,\kappa}[\DTL^2H](t,x)-H(t,x).
	$$
	
	To prove \eqref{eq:Newton-cone-gradient-no-loss}, set $
	K(s,y,R,\xi):=H(s,X_u(s,y+R\xi)).
	$ Here $\nabla_{\mathbb S^2}$ denotes the tangential gradient with respect to $\xi$. For $R>0$,
	$$
	\partial_RK=\xi\cdot\nabla_yK,
	\qquad
	\nabla_{\mathbb S^2}K
	=R(\Id-\xi\otimes\xi)\nabla_yK,
	$$
	and consequently
	$
	\nabla_yK
	=\xi\partial_RK+R^{-1}\nabla_{\mathbb S^2}K.
	$
	For every constant vector $e\in\mathbb R^3$, $
	\Div_{\mathbb S^2}\bigl((\Id-\xi\otimes\xi)e\bigr)
	=-2e\cdot\xi.
	$ Since $\nabla_{\mathbb S^2}K$ is tangent to $\mathbb S^2$, we integrate by parts on $\mathbb S^2$ to obtain
	$$
	\begin{aligned}
		e\cdot\int_{\mathbb S^2}\nabla_{\mathbb S^2}K\,\rd S_\xi
		=\int_{\mathbb S^2}
		(\Id-\xi\otimes\xi)e\cdot\nabla_{\mathbb S^2}K\,\rd S_\xi
		=2e\cdot\int_{\mathbb S^2}\xi K\,\rd S_\xi.
	\end{aligned}
	$$
	Since $e$ is arbitrary,
	$$
	\int_{\mathbb S^2}\nabla_{\mathbb S^2}K\,\rd S_\xi
	=2\int_{\mathbb S^2}\xi K\,\rd S_\xi.
	$$
	Differentiate \eqref{def:Newton-Lagrangian-cone-operator} in space. Using the identities for $\nabla_yK$ and $\int_{\mathbb S^2}\nabla_{\mathbb S^2}K\,\rd S_\xi$, we obtain
	$$
	\begin{aligned}
		\nabla_y\widehat I_{q,u,\kappa}[H](t,y)
		&=\frac{1}{4\pi}
		\int_{t_{u-1}}^t\int_{\mathbb S^2}
		(t-s)\left(
		\xi\partial_RK+\frac1R\nabla_{\mathbb S^2}K
		\right)\Big|_{R=c_{q,u,\kappa}(t-s)}\,\rd S_\xi\,\rd s\\
		&=\frac{1}{4\pi c_{q,u,\kappa}}
		\int_{t_{u-1}}^t\int_{\mathbb S^2}
		\xi\left[c_{q,u,\kappa}(t-s)\partial_RK+2K\right]\Big|_{R=c_{q,u,\kappa}(t-s)}
		\,\rd S_\xi\,\rd s.
	\end{aligned}
	$$
	For $
	Y_y=X_u(s,y+c_{q,u,\kappa}(t-s)\xi),
	$ we have
	$$
	c_{q,u,\kappa}\partial_RK(s,y,c_{q,u,\kappa}(t-s),\xi)
	=(\DTL H)(s,Y_y)-\partial_s(H(s,Y_y)).
	$$
	Therefore,
	$$
	\begin{aligned}
		\nabla_y\widehat I_{q,u,\kappa}[H](t,y)
		=\frac{1}{4\pi c_{q,u,\kappa}}\int_{t_{u-1}}^t
		\int_{\mathbb S^2}\xi\Big[(t-s)(\DTL H)(s,Y_y)
		-(t-s)\partial_s(H(s,Y_y))+2H(s,Y_y)\Big]
		\,\rd S_\xi\,\rd s.
	\end{aligned}
	$$
	Moreover,
	$$
	\begin{aligned}
		-\int_{t_{u-1}}^t
		(t-s)\partial_s(H(s,Y_y))\,\rd s
		=-\left[(t-s)H(s,Y_y)\right]_{s=t_{u-1}}^{s=t}
		-\int_{t_{u-1}}^tH(s,Y_y)\,\rd s
		=-\int_{t_{u-1}}^tH(s,Y_y)\,\rd s.
	\end{aligned}
	$$
	The boundary terms vanish because $t-s=0$ at $s=t$ and $H(t_{u-1},\cdot)=0$. Combining these two identities proves
	\eqref{eq:Newton-cone-gradient-no-loss}.
	
	By \eqref{eq:Newton-Lagrangian-material-commutation}, we have
	$
	\partial_t^m\nabla_y\widehat I_{q,u,\kappa}[H]
	=
	\nabla_y\widehat I_{q,u,\kappa}[\DTL^mH].
	$
	Applying \eqref{eq:Newton-cone-gradient-no-loss} to $\DTL^mH$ proves \eqref{eq:Newton-commuted-cone-gradient}. Since translations and normalized spherical
	averaging do not increase the $C^{N+\alpha}$ norm,
	$$
	\nrm{\partial_t^m\nabla_y\widehat I_{q,u,\kappa}[H](t,\cdot)}_{N+\alpha}
	\lesssim c_{q,u,\kappa}^{-1}(t-t_{u-1})
	\left[
	\sup_{t_{u-1}\leq s\leq t}\nrm{(\DTL^mH)(s,X_u(s,\cdot))}_{N+\alpha}
	+(t-t_{u-1})\sup_{t_{u-1}\leq s\leq t}\nrm{(\DTL^{m+1}H)(s,X_u(s,\cdot))}_{N+\alpha}
	\right].
	$$
	For $t_{u-1}\leq s\leq t$,
	$$
	\begin{aligned}
		(\DTL^mH)(s,X_u(s,\cdot))
		=\int_{t_{u-1}}^s
		(\DTL^{m+1}H)(\sigma,X_u(\sigma,\cdot))
		\,\rd\sigma,
	\end{aligned}
	$$
	because $(\DTL^mH)(t_{u-1},\cdot)=0$. Hence
	$$
	\begin{aligned}
		\sup_{t_{u-1}\leq s\leq t}
		\nrm{(\DTL^mH)(s,X_u(s,\cdot))}_{N+\alpha}
		\leq
		(t-t_{u-1})
		\sup_{t_{u-1}\leq s\leq t}
		\nrm{(\DTL^{m+1}H)(s,X_u(s,\cdot))}_{N+\alpha}.
	\end{aligned}
	$$
	This proves the second estimate in \eqref{eq:Newton-Lagrangian-cone-holder-estimates}. For the first estimate, \eqref{def:Newton-Lagrangian-cone-operator} and
	\eqref{eq:Newton-Lagrangian-material-commutation} give
	$$
	\begin{aligned}
		\partial_t^m\widehat I_{q,u,\kappa}[H](t,y)
		=\int_{t_{u-1}}^t\frac{t-s}{4\pi}
		\int_{\mathbb S^2}
		(\DTL^mH)\left(
		s,X_u\left(s,y+c_{q,u,\kappa}(t-s)\xi\right)
		\right)\,\rd S_\xi\,\rd s.
	\end{aligned}
	$$
	Therefore,
	$$
	\begin{aligned}
		\nrm{\partial_t^m\widehat I_{q,u,\kappa}[H](t,\cdot)}_{N+\alpha}
		\lesssim
		\int_{t_{u-1}}^t(t-s)
		\nrm{(\DTL^mH)(s,X_u(s,\cdot))}_{N+\alpha}\,\rd s
		\lesssim
		(t-t_{u-1})^2
		\sup_{t_{u-1}\leq s\leq t}
		\nrm{(\DTL^mH)(s,X_u(s,\cdot))}_{N+\alpha}.
	\end{aligned}
	$$
	
	To prove \eqref{eq:Newton-cone-slice-support}--\eqref{eq:Newton-cone-cylinder-support}, suppose that $\widehat I_{q,u,\kappa}[H](t,y)\neq0$. By \eqref{def:Newton-Lagrangian-cone-operator}, there exist $t_{u-1}\leq s\leq t$ and $\xi\in\mathbb S^2$ such that
	$H\left(s,X_u\left(s,y+c_{q,u,\kappa}(t-s)\xi\right)\right)\neq0$. The support assumption \eqref{eq:Newton-cone-input-support} gives
	$$
	|y-x_0|
	<
	\bar\varepsilon+c_{q,u,\kappa}(t-s)
	\leq
	\bar\varepsilon+c_{q,u,\kappa}(t-t_{u-1}).
	$$
	Taking the closure gives the inclusion for $\widehat I_{q,u,\kappa}[H]$ in \eqref{eq:Newton-cone-slice-support}. The inclusion for $I_{q,u,\kappa}[H]$ follows from
	$I_{q,u,\kappa}[H](t,x)=\widehat I_{q,u,\kappa}[H](t,\Phi_u(t,x))$.
	If $\bar\tau\geq2\tau_q/3$ and $t\in\supp\theta_u$, then $t$ and every contributing source time $s$ lie in $[t_u-\bar\tau,t_u+\bar\tau]$, so $0\leq t-s\leq2\bar\tau$.
	The spatial projection of the pulled-back support of $H$ is a compact subset of $B(x_0,\bar\varepsilon)$, so the radius in \eqref{eq:Newton-cone-input-support} may be chosen strictly smaller than $\bar\varepsilon$.
	The time support of $\theta_u$ then gives \eqref{eq:Newton-cone-cylinder-support}.
\end{proof}

\subsubsection{Construction and localization of the potential}

We now take the gradient of the Kirchhoff solution in flow coordinates and apply the Piola transform. For the mollified source in
\eqref{def:Newton-sharp-acoustic-input}, define the local potential by
\begin{equation}\label{def:Newton-raw-Piola-potential}
	(\widetilde V_{q,j+1,u,\kappa})_i
	:=(\mathbf C_u)_{n i}\partial_{y_n}
	\widehat I_{q,u,\kappa}[G_{q,j,u,\kappa}^{\mathrm{mol}}](t,\Phi_u).
\end{equation}
Since the source vanishes near $t_{u-1}$, the Kirchhoff solution extends smoothly by zero for $t<t_{u-1}$.

For every smooth vector field $Z=Z(t,y)$, the Piola identity gives
\begin{equation}\label{eq:Newton-Piola-transform}
	\Div_x\bigl(\mathbf C_u^TZ(t,\Phi_u)\bigr)=J_u(\Div_yZ)(t,\Phi_u).
\end{equation}
This follows from $\partial_i(\mathbf C_u)_{n i}=0$ and $(\mathbf C_u)_{n i}\partial_i(\Phi_{u})_r=J_u\delta_{nr}$. Applying Lemma~\ref{lem:Newton-three-dimensional-cone-identities} and
\eqref{eq:Newton-Piola-transform}, we obtain
\begin{equation}\label{eq:Newton-sharp-density-identity}
	\varrho[\widetilde V_{q,j+1,u,\kappa}]
	=\frac{J_u}{c_{q,u,\kappa}^2}\left(
	I_{q,u,\kappa}[\DTL^2G_{q,j,u,\kappa}^{\mathrm{mol}}]
	-G_{q,j,u,\kappa}^{\mathrm{mol}}
	\right).
\end{equation}
It follows that
\begin{equation}\label{eq:Newton-sharp-pressure-cancellation}
	c_{q,u,\kappa}^2\varrho[\widetilde V_{q,j+1,u,\kappa}]+J_uG_{q,j,u,\kappa}^{\mathrm{mol}}=J_uI_{q,u,\kappa}[\DTL^2G_{q,j,u,\kappa}^{\mathrm{mol}}].
\end{equation}
The factor $J_u^{-1}$ in \eqref{def:Newton-sharp-acoustic-input} cancels the Jacobian factor in the unmollified source.
The along-flow mollification error is included in $P_{q+1,j+1}^{\mathrm{mol,comp}}$ in \eqref{def:Newton-substep-flow-mollification}.

For a smooth vector field $V$, \eqref{def:Newton-perturbations-from-potential} can be written in terms of $\DTL$ as
\begin{equation}\label{eq:Newton-momentum-in-material-form}
	w_i[V]
	=-\partial_tV_i+\partial_k\left((\overline{u}_{q})_iV_k-(\overline{u}_{q})_kV_i\right)
	=-\DTL V_i+(\overline{u}_{q})_i\partial_kV_k+\partial_k(\overline{u}_{q})_iV_k-(\Div\overline{u}_q)V_i.
\end{equation}
For $\widetilde V_{q,j+1,u,\kappa}$, the identity 
$$
\DTL\mathbf C_u=\mathbf C_u(\nabla\overline u_q)^T
-(\Div\overline u_q)\mathbf C_u,
$$
together with \eqref{eq:Newton-Lagrangian-material-commutation} and \eqref{eq:Newton-momentum-in-material-form}, gives
\begin{equation}\label{eq:Newton-raw-momentum-Piola-cancellation}
	w_i[\widetilde V_{q,j+1,u,\kappa}]=(\overline u_{q})_i\varrho[\widetilde V_{q,j+1,u,\kappa}]-(\mathbf C_u)_{ni}\partial_{y_n}\widehat I_{q,u,\kappa}[\DTL G_{q,j,u,\kappa}^{\mathrm{mol}}](t,\Phi_u).
\end{equation}

For the local right inverse used below, we also need the support of $\theta_u\widetilde V_{q,j+1,u,\kappa}$. By \eqref{def:Newton-sharp-acoustic-input} and the cutoff support inclusions,
$$
\supp G_{q,j,u,\kappa}^{\mathrm{mol}}
\subseteq
\hat\Gamma_{\overline u_q}\left(
\frac23\tau_q,
\frac{9\sqrt3\pi}{8}\tau_q;
t_u,x_\kappa
\right).
$$
Choose a fixed constant $C_{\mathrm{supp}}$ such that $C_{\mathrm{supp}}/2>9\sqrt3\pi/8+
(5/3)\max_{\rho\in I_{\rho^*}}\sqrt{p'(\rho)}$. Applying \eqref{eq:Newton-cone-cylinder-support} with $\bar\tau=2\tau_q/3$,
$\bar\varepsilon=9\sqrt3\pi\tau_q/8$, and $x_0=x_\kappa$ shows that the spatial radius is at most $C_{\mathrm{supp}}\tau_q/2$. In particular,
\begin{equation}\label{eq:Newton-localized-Piola-potential-support}
	\supp\left(
	\theta_u\widetilde V_{q,j+1,u,\kappa}
	\right)
	\subseteq
	\hat\Gamma_{\overline u_q}
	\left(
	\frac23\tau_q,C_{\mathrm{supp}}\tau_q;
	t_u,x_\kappa
	\right).
\end{equation}
Multiplication by $\mathbf C_u^T$ and differentiation in $y$ do not enlarge the support. The same inclusion holds with
$\theta_u$ replaced by any of its time derivatives.

For the local right inverse of $\Div$, we correct the spatial mean of the potential so that the interior and cutoff residuals are separately mean-free. Choose $\phi\in C_c^\infty\left(B\left(0,\frac13\right)\right)$ with $\int_{\mathbb R^3}\phi\,\rd x=1$ and define
\begin{equation}\label{def:Newton-mean-profile}
\phi_{u,\kappa}(t,x):=\tau_q^{-3}J_u(t,x)
\phi\left(\frac{\Phi_u(t,x)-x_\kappa}{\tau_q}\right).
\end{equation}
Changing variables $y=\Phi_u(t,x)$ gives $\int_{\mathbb R^3}\phi_{u,\kappa}\,\rd x=1$. The identities $\DTL\Phi_u=0$ and $\DTL J_u=-(\Div\overline u_q)J_u$ also imply
\begin{equation}\label{eq:Newton-transported-mean-profile}
	\partial_t\phi_{u,\kappa}+\Div(\overline{u}_q\phi_{u,\kappa})=0.
\end{equation}
Moreover,
$$
\supp\left(\theta_u\phi_{u,\kappa}\right)
\subseteq
\hat\Gamma_{\overline u_q}\left(
\frac23\tau_q,\frac13\tau_q;
t_u,x_\kappa
\right).
$$
Set
\begin{equation}\label{def:Newton-mean-corrected-Piola-potential}
	\begin{aligned}
		\overline V_{q,j+1,u,\kappa}(t):=\int_{\mathbb R^3}\widetilde V_{q,j+1,u,\kappa}(t,x)\,\rd x,\quad
		V_{q,j+1,u,\kappa}^{\circ}
		:=\widetilde V_{q,j+1,u,\kappa}
		-\overline V_{q,j+1,u,\kappa}(t)\phi_{u,\kappa}(t,x).
	\end{aligned}
\end{equation}
Then $
\int_{\mathbb R^3}V_{q,j+1,u,\kappa}^{\circ}\,\rd x=0.
$ Since $\overline V_{q,j+1,u,\kappa}$ depends only on time, \eqref{def:Newton-perturbations-from-potential},
\eqref{eq:Newton-momentum-in-material-form}, and \eqref{eq:Newton-transported-mean-profile} give
\begin{equation}\label{eq:Newton-mean-perturbation-identities}
	\begin{aligned}
		\varrho[\overline V_{q,j+1,u,\kappa}\phi_{u,\kappa}]
		&=\overline V_{q,j+1,u,\kappa}\cdot\nabla\phi_{u,\kappa},\\
		w[\overline V_{q,j+1,u,\kappa}\phi_{u,\kappa}]
		&=-\left(\partial_t\overline V_{q,j+1,u,\kappa}\right)\phi_{u,\kappa}
		+\overline u_q\varrho[\overline V_{q,j+1,u,\kappa}\phi_{u,\kappa}]+(\nabla\overline u_q)(\overline V_{q,j+1,u,\kappa}\phi_{u,\kappa}).
	\end{aligned}
\end{equation}
Subtracting these expressions from \eqref{eq:Newton-sharp-density-identity} and \eqref{eq:Newton-raw-momentum-Piola-cancellation} gives the density and momentum associated with $V_{q,j+1,u,\kappa}^{\circ}$.

The Kirchhoff solution need not vanish after its source vanishes. We therefore localize the corrected potential in time by setting
$
V_{q,j+1,u,\kappa}\comp:=\theta_uV_{q,j+1,u,\kappa}^{\circ}.
$
The inclusion \eqref{eq:Newton-localized-Piola-potential-support}, together with the support of $\theta_u\phi_{u,\kappa}$, implies
$$
\supp V_{q,j+1,u,\kappa}\comp
\subseteq
\hat\Gamma_{\overline u_q}
\left(
\frac23\tau_q,C_{\mathrm{supp}}\tau_q;
t_u,x_\kappa
\right).
$$
The local density and momentum perturbations are
$$
	\trho_{q,j+1,u,\kappa}\comp:=\varrho[V_{q,j+1,u,\kappa}\comp],
\qquad
	\tm_{q,j+1,u,\kappa}\comp:=w[V_{q,j+1,u,\kappa}\comp].
$$
Since $\theta_u$ depends only on time,
\begin{equation}\label{eq:Newton-local-compressible-cutoff-formulas}
	\trho_{q,j+1,u,\kappa}\comp
	=\theta_u\varrho[V_{q,j+1,u,\kappa}^{\circ}],
	\qquad
	\tm_{q,j+1,u,\kappa}\comp
	=\theta_uw[V_{q,j+1,u,\kappa}^{\circ}]
	-(\partial_t\theta_u)V_{q,j+1,u,\kappa}^{\circ}.
\end{equation}
For any smooth vector field $V$,
$
\partial_i\partial_k\left((\overline{u}_{q})_iV_k-(\overline{u}_{q})_kV_i\right)=0,
$
so \eqref{def:Newton-perturbations-from-potential} gives
\begin{equation}\label{eq:Newton-exact-continuity-from-potential}
	\partial_t\varrho[V]+\Div w[V]=0.
\end{equation}
Applying \eqref{eq:Newton-exact-continuity-from-potential} to $V_{q,j+1,u,\kappa}\comp$, we obtain
$$
\partial_t\trho_{q,j+1,u,\kappa}\comp
+\Div\tm_{q,j+1,u,\kappa}\comp=0.
$$

We first compute the momentum residual for $\widetilde V_{q,j+1,u,\kappa}$. We then include the mean correction $\overline V_{q,j+1,u,\kappa}\phi_{u,\kappa}$ and the time cutoff $\theta_u$. The first two material derivatives of $\widetilde V_{q,j+1,u,\kappa}$ are
\begin{align}
	\DTL(\widetilde V_{q,j+1,u,\kappa})_i
	&=(\DTL\mathbf C_u)_{n i}\partial_{y_n}\widehat I_{q,u,\kappa}[G_{q,j,u,\kappa}^{\mathrm{mol}}](t,\Phi_u)+(\mathbf C_u)_{n i}\partial_{y_n}\widehat I_{q,u,\kappa}[\DTL G_{q,j,u,\kappa}^{\mathrm{mol}}](t,\Phi_u),
	\notag\\
	\DTL^2(\widetilde V_{q,j+1,u,\kappa})_i
	&=(\DTL^2\mathbf C_u)_{n i}\partial_{y_n}\widehat I_{q,u,\kappa}[G_{q,j,u,\kappa}^{\mathrm{mol}}](t,\Phi_u)+2(\DTL\mathbf C_u)_{n i}\partial_{y_n}\widehat I_{q,u,\kappa}[\DTL G_{q,j,u,\kappa}^{\mathrm{mol}}](t,\Phi_u)\label{eq:Newton-second-material-Piola}\\
	&\quad+(\mathbf C_u)_{n i}\partial_{y_n}\widehat I_{q,u,\kappa}[\DTL^2G_{q,j,u,\kappa}^{\mathrm{mol}}](t,\Phi_u).\notag	
\end{align}
\begin{samepage}
	Combining \eqref{eq:Newton-sharp-pressure-cancellation} with \eqref{eq:Newton-second-material-Piola}, we have
	\begin{equation}\label{eq:Newton-sharp-acoustic-local-identity}
		\begin{aligned}
			-\DTL^2(\widetilde V_{q,j+1,u,\kappa})_i
			+\partial_i\left(p'(\orho_q)\partial_k(\widetilde V_{q,j+1,u,\kappa})_k\right)
			+\partial_i(J_uG_{q,j,u,\kappa}^{\mathrm{mol}})=\partial_i\left[
			\left(p'(\orho_q)-c_{q,u,\kappa}^2\right)\varrho[\widetilde V_{q,j+1,u,\kappa}]\right]
			+(\mathcal E_{q,j,u,\kappa}^{\sharp})_i,
		\end{aligned}
	\end{equation}
	where
	\begin{equation}\label{def:Newton-sharp-vector-remainder}
		\begin{aligned}
			(\mathcal E_{q,j,u,\kappa}^{\sharp})_i
			:=&\partial_iJ_u\,I_{q,u,\kappa}[\DTL^2G_{q,j,u,\kappa}^{\mathrm{mol}}]
			+\left[J_u\partial_i(\Phi_{u})_n-(\mathbf C_u)_{n i}\right]
			\partial_{y_n}\widehat I_{q,u,\kappa}[\DTL^2G_{q,j,u,\kappa}^{\mathrm{mol}}](t,\Phi_u)\\
			&-2(\DTL\mathbf C_u)_{n i}
			\partial_{y_n}\widehat I_{q,u,\kappa}[\DTL G_{q,j,u,\kappa}^{\mathrm{mol}}](t,\Phi_u)
			-(\DTL^2\mathbf C_u)_{n i}
			\partial_{y_n}\widehat I_{q,u,\kappa}[G_{q,j,u,\kappa}^{\mathrm{mol}}](t,\Phi_u).
		\end{aligned}
	\end{equation}
\end{samepage}
Indeed,
$$
\partial_i\left(J_uI_{q,u,\kappa}[\DTL^2G_{q,j,u,\kappa}^{\mathrm{mol}}]\right)
=\partial_iJ_uI_{q,u,\kappa}[\DTL^2G_{q,j,u,\kappa}^{\mathrm{mol}}]
+J_u\partial_i(\Phi_{u})_n
\partial_{y_n}\widehat I_{q,u,\kappa}[\DTL^2G_{q,j,u,\kappa}^{\mathrm{mol}}](t,\Phi_u).
$$
Subtracting the formula for the second material derivative in \eqref{eq:Newton-second-material-Piola} gives \eqref{eq:Newton-sharp-acoustic-local-identity}.

For a general potential $V$, substituting \eqref{eq:Newton-momentum-in-material-form}, we obtain
\begin{equation}\label{eq:Newton-linearized-momentum-from-potential}
	\begin{aligned}
		&\partial_tw_i[V]
		+\partial_k\left((\overline{u}_{q})_kw_i[V]+w_k[V](\overline{u}_{q})_i
		-\varrho[V](\overline{u}_{q})_i(\overline{u}_{q})_k\right)
		+\partial_i\left(p'(\orho_q)\varrho[V]\right)\\
		=&-\DTL^2V_i
		+\partial_i\left(p'(\orho_q)\partial_kV_k\right)
		+\partial_k(\DTL(\overline{u}_{q})_iV_k)
		-\DTL\left((\Div\overline{u}_q)V_i\right)-(\Div\overline{u}_q)\DTL V_i-(\Div\overline{u}_q)^2V_i.
	\end{aligned}
\end{equation}
Indeed, by \eqref{eq:Newton-exact-continuity-from-potential}, we have
$$
\begin{aligned}
	&\partial_tw_i[V]+\partial_k\left((\overline{u}_{q})_kw_i[V]+w_k[V](\overline{u}_{q})_i-\varrho[V](\overline{u}_{q})_i(\overline{u}_{q})_k\right)\\
	=&\DTL\left(w_i[V]-(\overline{u}_{q})_i\varrho[V]\right)+(\Div\overline u_q)\left(w_i[V]-(\overline{u}_{q})_i\varrho[V]\right)+\partial_k(\overline{u}_{q})_i\left(w_k[V]-(\overline{u}_{q})_k\varrho[V]\right)+(\DTL\overline{u}_{q})_i\varrho[V].
\end{aligned}
$$
By \eqref{eq:Newton-momentum-in-material-form}, the terms $\partial_k(\overline{u}_{q})_i\DTL V_k$ cancel in pairs, as do the terms $(\Div\overline u_q)\partial_k(\overline{u}_{q})_iV_k$. The commutator $[\DTL,\partial_k]=-(\partial_k\overline{u}_{q})_r\partial_r$ gives
$$
(\DTL\partial_k\overline{u}_{q})_iV_k+\partial_k(\overline{u}_{q})_i\partial_r(\overline{u}_{q})_kV_r
=\partial_k\DTL(\overline{u}_{q})_iV_k.
$$
Since $\varrho[V]=\partial_kV_k$, this term combines with $(\DTL\overline{u}_{q})_i\varrho[V]$ to give $\partial_k((\DTL\overline{u}_{q})_iV_k)$. Adding the pressure gradient yields
\eqref{eq:Newton-linearized-momentum-from-potential}.
Applying \eqref{eq:Newton-linearized-momentum-from-potential} to $\widetilde V_{q,j+1,u,\kappa}$ and using \eqref{eq:Newton-sharp-acoustic-local-identity}, we obtain
\begin{equation}\label{eq:Newton-raw-Piola-momentum-equation}
\begin{aligned}
	&\partial_tw_i[\widetilde V_{q,j+1,u,\kappa}]
	+\partial_i\left(p'(\orho_q)\varrho[\widetilde V_{q,j+1,u,\kappa}]\right)
	+\partial_i(J_uG_{q,j,u,\kappa}^{\mathrm{mol}})\\
	&\quad+\partial_k\left((\overline{u}_{q})_kw_i[\widetilde V_{q,j+1,u,\kappa}]
	+w_k[\widetilde V_{q,j+1,u,\kappa}](\overline{u}_{q})_i
	-\varrho[\widetilde V_{q,j+1,u,\kappa}](\overline{u}_{q})_i(\overline{u}_{q})_k\right)\\
	=&\partial_i\left[
	\left(p'(\orho_q)-c_{q,u,\kappa}^2\right)\varrho[\widetilde V_{q,j+1,u,\kappa}]\right]
	+(\mathcal E_{q,j,u,\kappa}^{\sharp})_i
	+\partial_k(\DTL(\overline{u}_{q})_i(\widetilde V_{q,j+1,u,\kappa})_k)\\
	&-\DTL\left((\Div\overline{u}_q)(\widetilde V_{q,j+1,u,\kappa})_i\right)-(\Div\overline{u}_q)\DTL(\widetilde V_{q,j+1,u,\kappa})_i
	-(\Div\overline{u}_q)^2(\widetilde V_{q,j+1,u,\kappa})_i.
\end{aligned}
\end{equation}
We next compute the contribution of $\overline V_{q,j+1,u,\kappa}\phi_{u,\kappa}$. Since $\overline V_{q,j+1,u,\kappa}$ depends only on time, \eqref{eq:Newton-transported-mean-profile} gives
$$
\DTL(\overline V_{q,j+1,u,\kappa}\phi_{u,\kappa})+(\Div\overline u_q)\overline V_{q,j+1,u,\kappa}\phi_{u,\kappa}=\left(\partial_t\overline V_{q,j+1,u,\kappa}\right)\phi_{u,\kappa}.
$$
Applying $\DTL+\Div\overline u_q$ once more, we have
$$
(\DTL+\Div\overline u_q)^2(\overline V_{q,j+1,u,\kappa}\phi_{u,\kappa})=\left(\partial_t^{2}\overline V_{q,j+1,u,\kappa}\right)\phi_{u,\kappa}.
$$
In \eqref{eq:Newton-linearized-momentum-from-potential}, the first term and the last three terms on the right combine to give $-(\DTL+\Div\overline u_q)^2V$. Set $V=\overline V_{q,j+1,u,\kappa}\phi_{u,\kappa}$. Since $\varrho[V]=\Div V$, the pressure gradients on both sides cancel, giving
$$
\begin{aligned}
	&\partial_tw[\overline V_{q,j+1,u,\kappa}\phi_{u,\kappa}]+\Div\left(\overline u_q\otimes w[\overline V_{q,j+1,u,\kappa}\phi_{u,\kappa}]+w[\overline V_{q,j+1,u,\kappa}\phi_{u,\kappa}]\otimes\overline u_q\right)-\Div\left(\varrho[\overline V_{q,j+1,u,\kappa}\phi_{u,\kappa}]\overline u_q\otimes\overline u_q\right)\\=&-\left(\partial_t^{2}\overline V_{q,j+1,u,\kappa}\right)\phi_{u,\kappa}
	+\Div\left((\DTL\overline u_q)\otimes(\overline V_{q,j+1,u,\kappa}\phi_{u,\kappa})\right).
\end{aligned}
$$
We subtract this identity and the pressure contribution of the mean correction from \eqref{eq:Newton-raw-Piola-momentum-equation}. Multiplying the resulting
residual by $\theta_u$, we define
\begin{equation}\label{eq:Newton-explicit-local-compressible-residual}
	\begin{aligned}
		(F_{q,j+1,u,\kappa}\comp)_i
		:=&\theta_u\partial_i\Bigl[\left(p'(\orho_q)-c_{q,u,\kappa}^2\right)\varrho[\widetilde V_{q,j+1,u,\kappa}]-p'(\orho_q)\varrho[\overline V_{q,j+1,u,\kappa}\phi_{u,\kappa}]\Bigr]+\theta_u\partial_k\left((\DTL\overline u_{q})_i(V_{q,j+1,u,\kappa}^{\circ})_k\right)\\
		&+\theta_u\Bigl[(\mathcal E_{q,j,u,\kappa}^{\sharp})_i-\DTL\left((\Div\overline u_q)(\widetilde V_{q,j+1,u,\kappa})_i\right)-(\Div\overline u_q)\DTL(\widetilde V_{q,j+1,u,\kappa})_i-(\Div\overline u_q)^2(\widetilde V_{q,j+1,u,\kappa})_i\Bigr]\\
		&+\theta_u\left(\partial_t^{2}\overline V_{q,j+1,u,\kappa}\right)_i\phi_{u,\kappa}.
	\end{aligned}
\end{equation}

The remaining residual comes from differentiating the time cutoff $\theta_u$. For a smooth scalar function $f=f(t)$, the definitions give $\varrho[fV]=f\varrho[V]$ and $w[fV]=fw[V]-(\partial_tf)V$. Using
\eqref{eq:Newton-linearized-momentum-from-potential}, we find 
\begin{equation}\label{eq:Newton-temporal-cutoff-momentum-identity}
\begin{aligned}
	&\partial_tw_i[fV]
	+\partial_k\left((\overline{u}_{q})_kw_i[fV]+w_k[fV](\overline{u}_{q})_i
	-\varrho[fV](\overline{u}_{q})_i(\overline{u}_{q})_k\right)+\partial_i\left(p'(\orho_q)\varrho[fV]\right)\\
	=&f\Bigl[\partial_tw_i[V]
	+\partial_k\left((\overline{u}_{q})_kw_i[V]+w_k[V](\overline{u}_{q})_i
	-\varrho[V](\overline{u}_{q})_i(\overline{u}_{q})_k\right)
	+\partial_i\left(p'(\orho_q)\varrho[V]\right)\Bigr]\\
	&-2(\partial_tf)\left(\DTL V_i+(\Div\overline{u}_q)V_i\right)
	-(\partial_t^2f)V_i.
\end{aligned}
\end{equation}
With $f=\theta_u$ and $V=V_{q,j+1,u,\kappa}^{\circ}$, the terms containing $\partial_t\theta_u$ or $\partial_t^2\theta_u$ form the residual
\begin{equation}\label{def:Newton-compressible-cutoff-residual}
	\begin{aligned}
		F_{q,j+1,u,\kappa}^{\mathrm{cut}}
		&:=-2(\partial_t\theta_u)
		\left(\DTL V_{q,j+1,u,\kappa}^{\circ}
		+(\Div\overline{u}_q)V_{q,j+1,u,\kappa}^{\circ}\right)-(\partial_t^2\theta_u)
		V_{q,j+1,u,\kappa}^{\circ},\\
		F_{q,j+1,u}^{\mathrm{cut}}
		&:=\sum_{\kappa\in\mathbb Z^3}F_{q,j+1,u,\kappa}^{\mathrm{cut}},\\
		F_{q,j+1}^{\mathrm{cut}}
		&:=\sum_uF_{q,j+1,u}^{\mathrm{cut}}=\sum_u\sum_{\kappa\in\mathbb Z^3}F_{q,j+1,u,\kappa}^{\mathrm{cut}}.
	\end{aligned}
\end{equation}
Combining \eqref{eq:Newton-explicit-local-compressible-residual}--\eqref{def:Newton-compressible-cutoff-residual}, we have
\begin{equation}\label{eq:Newton-local-compressible-block-equation}
	\left\{
	\begin{aligned}
		&\partial_t\trho_{q,j+1,u,\kappa}\comp
		+\Div\tm_{q,j+1,u,\kappa}\comp=0,\\
		&\partial_t\tm_{q,j+1,u,\kappa}\comp
		+\Div\left(\overline{u}_q\otimes\tm_{q,j+1,u,\kappa}\comp
		+\tm_{q,j+1,u,\kappa}\comp\otimes\overline{u}_q
		-\trho_{q,j+1,u,\kappa}\comp\overline{u}_q\otimes\overline{u}_q\right)+\nabla\left(p'(\orho_q)\trho_{q,j+1,u,\kappa}\comp\right)\\
		&\hspace{98pt}
		=-\theta_u\nabla(J_uG_{q,j,u,\kappa}^{\mathrm{mol}})
		+F_{q,j+1,u,\kappa}\comp
		+F_{q,j+1,u,\kappa}^{\mathrm{cut}}.
	\end{aligned}
	\right.
\end{equation}
The temporal supports satisfy
$$
\supp_t\trho_{q,j+1,u,\kappa}\comp
\cup\supp_t\tm_{q,j+1,u,\kappa}\comp
\cup\supp_tF_{q,j+1,u,\kappa}\comp
\cup\supp_tF_{q,j+1,u,\kappa}^{\mathrm{cut}}
\subset\supp\theta_u.
$$
Moreover, all four fields are supported in the adapted cylinder
$
\hat\Gamma_{\overline u_q}
\left(
\frac23\tau_q,C_{\mathrm{supp}}\tau_q;
t_u,x_\kappa
\right).
$

Define the total compressible perturbations by
\begin{equation}\label{def:Newton-total-compressible-perturbation}
	\begin{aligned}
		V_{q,j+1}\comp
		&:=\sum_{(u,\kappa)\in\mathscr I}V_{q,j+1,u,\kappa}\comp
		=\sum_{(u,\kappa)\in\mathscr I}
		\theta_uV_{q,j+1,u,\kappa}^{\circ},\\
		\trho_{q,j+1}\comp
		&:=\varrho[V_{q,j+1}\comp]
		=\sum_{(u,\kappa)\in\mathscr I}
		\theta_u\varrho[V_{q,j+1,u,\kappa}^{\circ}],\\
		\tm_{q,j+1}\comp
		&:=w[V_{q,j+1}\comp]
		=\sum_{(u,\kappa)\in\mathscr I}
		\left(
		\theta_uw[V_{q,j+1,u,\kappa}^{\circ}]
		-(\partial_t\theta_u)V_{q,j+1,u,\kappa}^{\circ}
		\right).
	\end{aligned}
\end{equation}
Define the total interior residual by
$$
F_{q,j+1}\comp
:=\sum_{(u,\kappa)\in\mathscr I}F_{q,j+1,u,\kappa}\comp.
$$
Summing \eqref{eq:Newton-local-compressible-block-equation} proves \eqref{eq:Newton-approximate-compressible-linearized-system}.

We next verify that $F_{q,j+1,u,\kappa}\comp$ and $F_{q,j+1,u,\kappa}^{\mathrm{cut}}$ have zero spatial integral, as required for the local right inverse. The mean correction gives $\int_{\mathbb R^3}V_{q,j+1,u,\kappa}\comp\,\rd x=0$ at every time. Integrating \eqref{def:Newton-perturbations-from-potential} then yields $\int_{\mathbb R^3}w[V_{q,j+1,u,\kappa}\comp]\,\rd x=0$.
Integrating \eqref{eq:Newton-local-compressible-block-equation} in space and using compact support gives
$$
\int_{\mathbb R^3}\left(
F_{q,j+1,u,\kappa}\comp
+F_{q,j+1,u,\kappa}^{\mathrm{cut}}
\right)\,\rd x=0.
$$
Moreover,
$$
\int_{\mathbb R^3}\left(
\DTL V_{q,j+1,u,\kappa}^{\circ}
+(\Div\overline{u}_q)V_{q,j+1,u,\kappa}^{\circ}
\right)\,\rd x
=\partial_t\int_{\mathbb R^3}V_{q,j+1,u,\kappa}^{\circ}\,\rd x=0.
$$
Together with $\int_{\mathbb R^3}V_{q,j+1,u,\kappa}^{\circ}\,\rd x=0$, this shows that every component of $F_{q,j+1,u,\kappa}^{\mathrm{cut}}$ has zero integral. It follows that
$$
\int_{\mathbb R^3}F_{q,j+1,u,\kappa}\comp\,\rd x=0.
$$

Since $\theta_u$ depends only on time, the explicit gradient and divergence terms in \eqref{eq:Newton-explicit-local-compressible-residual} equal the divergence of the matrix
\begin{equation}\label{def:Newton-explicit-local-acoustic-stress}
	\begin{aligned}
		\mathcal S_{q,j+1,u,\kappa}
		:=&\theta_u\Bigl[
		\left(p'(\orho_q)-c_{q,u,\kappa}^2\right)\varrho[\widetilde V_{q,j+1,u,\kappa}]
		-p'(\orho_q)\varrho[\overline V_{q,j+1,u,\kappa}\phi_{u,\kappa}]\Bigr]\Id+\theta_u(\DTL\overline{u}_q)\otimes V_{q,j+1,u,\kappa}^{\circ}.
	\end{aligned}
\end{equation}
We include this matrix in $R_{q,j+1}\comp$ and apply the local right inverse only to the remaining vector field, defined by
\begin{equation}\label{def:Newton-interior-vector-remainder}
	\begin{aligned}
		H_{q,j+1,u,\kappa}^{\mathrm{int}}\hspace{-2pt}
		:=\theta_u\Biggl[\mathcal E_{q,j,u,\kappa}^{\sharp}
		-\DTL\left((\Div\overline u_q)\widetilde V_{q,j+1,u,\kappa}\right)-(\Div\overline u_q)\DTL\widetilde V_{q,j+1,u,\kappa}
		-(\Div\overline u_q)^2\widetilde V_{q,j+1,u,\kappa}
		+\left(\partial_t^{2}\overline V_{q,j+1,u,\kappa}\right)\phi_{u,\kappa}\Biggr].
	\end{aligned}
\end{equation}
By construction, $F_{q,j+1,u,\kappa}\comp
=\Div\mathcal S_{q,j+1,u,\kappa}+H_{q,j+1,u,\kappa}^{\mathrm{int}}$. Since $F_{q,j+1,u,\kappa}\comp$ has zero spatial integral and $\mathcal S_{q,j+1,u,\kappa}$ is compactly supported,
$$
\int_{\mathbb R^3}H_{q,j+1,u,\kappa}^{\mathrm{int}}\,\rd x=0.
$$
Both $H_{q,j+1,u,\kappa}^{\mathrm{int}}$ and $F_{q,j+1,u,\kappa}^{\mathrm{cut}}$ therefore satisfy the support and zero-mean assumptions of Theorem~\ref{thm:inv4divEq}. For sufficiently large
$a$, $C_{\mathrm{supp}}\tau_q\leq1$. Estimate~\eqref{est:uniform-grad-om-over-orho} and \eqref{eq:Newton-scale-separation} imply the local short-time condition \eqref{eq:tau-local-smallness}. For each
$(u,\kappa)$, let $R^{(\overline{u}_q)}_{u,\kappa}$ be the local right inverse from Theorem~\ref{thm:inv4divEq}, with transport field $\overline u_q$,
point $(t_u,x_\kappa)$, time radius $2\tau_q/3$, and spatial radius $C_{\mathrm{supp}}\tau_q$. Define
\begin{equation}\label{def:Newton-local-compressible-residual-matrix}
	R_{q,j+1}\comp
	:=\sum_{(u,\kappa)\in\mathscr I}\left\{
	\mathcal S_{q,j+1,u,\kappa}
	+R^{(\overline{u}_q)}_{u,\kappa}\left[H_{q,j+1,u,\kappa}^{\mathrm{int}}\right]\right\}.
\end{equation}
The translation properties of the cutoffs and the periodicity of the background fields imply that $\mathcal S_{q,j+1,u,\kappa}$, $H_{q,j+1,u,\kappa}^{\mathrm{int}}$, and $F_{q,j+1,u,\kappa}^{\mathrm{cut}}$ satisfy
\eqref{eq:covariance of Newton cutoff inputs}. The sum is therefore smooth and periodic by Lemma~\ref{lem support of R under cutoffs}. The right-inverse identity gives
\begin{equation}\label{eq:Newton-compressible-residual-divergence-representation}
	F_{q,j+1}\comp=\Div R_{q,j+1}\comp.
\end{equation}

For the cutoff residual, use these local right inverses and set
\begin{equation}\label{def:Newton-cutoff-Piola-stress}
	S_{q,j+1,u}^{\mathrm{cut}}
	:=\sum_{\kappa\in\mathbb Z^3}
	R^{(\overline{u}_q)}_{u,\kappa}\left[F_{q,j+1,u,\kappa}^{\mathrm{cut}}\right].
\end{equation}
Lemma~\ref{lem support of R under cutoffs} and the covariance relation \eqref{eq:covariance of Newton cutoff inputs} imply that $S_{q,j+1,u}^{\mathrm{cut}}$ is smooth and periodic. It satisfies
\begin{equation}\label{eq:Newton-cutoff-Piola-divergence}
	\Div S_{q,j+1,u}^{\mathrm{cut}}
	=F_{q,j+1,u}^{\mathrm{cut}}.
\end{equation}
The local right inverse does not enlarge temporal support, so $S_{q,j+1,u}^{\mathrm{cut}}$ has temporal support in $\supp\partial_t\theta_u$. For each time block, write
\begin{align*}
S_{q,j+1,u}^{\mathrm{cut,sym}}
&:=\frac12\left(
S_{q,j+1,u}^{\mathrm{cut}}
+(S_{q,j+1,u}^{\mathrm{cut}})^T
\right),
\\
S_{q,j+1,u}^{\mathrm{cut,skew}}
&:=\frac12\left(
S_{q,j+1,u}^{\mathrm{cut}}
-(S_{q,j+1,u}^{\mathrm{cut}})^T
\right).
\end{align*}
The symmetric part enters $R_{q,j+1}^{\mathrm{gluing}}$. The divergence of the antisymmetric part is a curl, so it can be canceled by the incompressible correction in Subsection~\ref{subsec:Newton-incompressible-perturbation}. Indeed,
$$
\curl W\left[S_{q,j+1,u}^{\mathrm{cut,skew}}\right]
=\Div S_{q,j+1,u}^{\mathrm{cut,skew}}.
$$
Since $W$ acts only in space, both sides have temporal support in $\supp\partial_t\theta_u$.

\subsubsection{Estimates for the compressible perturbation}
\label{subsec:Estimates for the compressible perturbation}
\begin{lem}
	\label{lem:Newton-source-flow-mollification}
	Fix an integer $0\leq j<J_*$. Assume the background estimates in
	Proposition~\ref{prop:uniform-mollification-estimates} and \eqref{eq:Newton-interior-input-estimate} at the fixed H{\"o}lder exponents. For $a>a_{\mathrm{Newton}}$, the following estimates hold.
	For $M,N\in\mathbb N_0$ with $M\leq N_{\mathrm{ind}}-1$ and $M+N\leq K_j^{(0)}$, the scalar and vector potentials satisfy
	\begin{equation}\label{eq:Newton-compressible-input-estimate}
		\begin{aligned}
			&\nrm{\DTL^MQ[R_{q,j}]}_{N+\alpha}
			+\nrm{\DTL^MW[R_{q,j}]}_{N+\alpha}\\
			&\qquad\lesssim A_{q,j}
			\MM{M,N_R-3j,\tau_q^{-1},\ell_{t,q}^{-1}}
			\Lat{N,N_{R,j}^{(0)}(M),\lambda_q,\ell_q^{-1}}
			\MM{N,N_{R,j}^{(0)}(M),\lambda_q,\ell_q^{-1}}.
		\end{aligned}
	\end{equation}
	The mollified sources satisfy
	\begin{equation}\label{eq:Newton-mollified-source-estimate}
		\begin{aligned}
			&\nrm{\DTL^MG_{q,j,u,\kappa}^{\mathrm{mol}}}_{N+\alpha}+\nrm{\DTL^MW_{q,j}^{\mathrm{mol}}}_{N+\alpha}\\
			&\quad\lesssim A_{q,j}\MM{M,N_R-3j,\tau_q^{-1},\ell_{t,q}^{-1}}\Lat{N,N_{R,j}^{(0)}(M),\lambda_q,\ell_q^{-1}}\MM{N,N_{R,j}^{(0)}(M),\lambda_q,\ell_q^{-1}},
		\end{aligned}
	\end{equation}
	for $0\leq M\leq N_{\mathrm{max}}$ and $\min\{M,N_R-3j\}+N\leq K_j^{(0)}$.
	Moreover,
	\begin{equation}\label{eq:Newton-mollified-source-support}
		\supp G_{q,j,u,\kappa}^{\mathrm{mol}}
		\subseteq\hat\Gamma_{\overline u_q}\left(
		\frac23\tau_q,\frac{9\sqrt3\pi}{8}\tau_q;t_u,x_\kappa
		\right).
	\end{equation}
	The implicit constants may depend on $\rho^*,\mathfrak M,N_{\mathrm{max}},\alpha,\alpha_r$.
	They are independent of $q,j,u,\kappa,a,C_*$, and $A_{q,j}$.
\end{lem}

\begin{proof}
	To prove \eqref{eq:Newton-compressible-input-estimate}, apply Proposition~\ref{prop:periodic-Hodge-material-estimates} with
	$$
	\begin{aligned}
		&K^*=N_{\mathrm{ind}}-2,\quad M_*=N_{\mathrm{ind}}-1,&& K_T=N_t=N_R-3j,&& K^\sharp=K_j^{(0)}+1,&& A_\upsilon=C\delta_q^{\frac12},&& A_T=A_{q,j},\quad s=1,\\
		&\lambda_\upsilon=\lambda_T=\lambda_q,\quad r^*=K^*+1,&& \tilde\lambda_\upsilon=\tilde\lambda_T=\ell_q^{-1},&& \mu_\upsilon=\tilde\mu_\upsilon=\Lambda_q,&& \mu_T=\tau_q^{-1},&& \tilde\mu_T=\ell_{t,q}^{-1}.
	\end{aligned}
	$$
	Choose $C\geq1$ to cover the constants in \eqref{est:uniform-grad-om-over-orho}
	and \eqref{est:holder-grad-om-over-orho}. These estimates verify the coefficient assumptions, while the choice of derivative parameters gives
	$K_T\leq K^*$ and $K^\sharp\geq K^*+2$. For sufficiently large $a$,
	$A_\upsilon\lambda_\upsilon=C\Lambda_q\leq\tau_q^{-1}$,
	and
	$A_\upsilon\tilde\lambda_\upsilon=C\delta_q^{\frac12}\ell_q^{-1}\leq\ell_{t,q}^{-1}$.
	Estimates~\eqref{eq:Hodge-high-material-estimate}--\eqref{eq:Hodge-material-zero-estimate}, with
	$s=1$ and $\Xi_\alpha\geq1$, give
	\eqref{eq:Newton-compressible-input-estimate}.
	
	The time supremum in \eqref{eq:Newton-Jacobian-material-estimates} is over $|t-t_u|\leq\tau_q$. Apply Corollary~\ref{cor:mollified-velocity-flow-estimates} with $h=\tau_q$ and $\sigma=\pm1$. The second estimate holds for $1\leq r\leq N_{\mathrm{ind}}$.
	\begin{equation}\label{eq:Newton-Jacobian-material-estimates}
	\begin{alignedat}{2}
		&\nrm{J_u^\sigma-1}_{N+\alpha}\lesssim\tau_q\delta_q\Lambda_q\La{N,N_{\mathrm{ind}}-1,\lambda_q,\ell_q^{-1}}\MM{N,N_{\mathrm{ind}}-1,\lambda_q,\ell_q^{-1}},
		&&\qquad N\leq N_{\mathrm{max}}-2,\\[2pt]
		&\nrm{\DTL^rJ_u^\sigma}_{N+\alpha}\lesssim\delta_q\Lambda_q^r\La{N,N_{\mathrm{ind}}-r,\lambda_q,\ell_q^{-1}}\MM{N,N_{\mathrm{ind}}-r,\lambda_q,\ell_q^{-1}},
		&&\qquad r+N\leq N_{\mathrm{max}}-1,\\[2pt]
		&\nrm{\DTL^{N_{\mathrm{ind}}+1}J_u^\sigma}_{N+\alpha}\lesssim\delta_q^{\frac12}\Lambda_q^{N_{\mathrm{ind}}}\ell_q^{-N-1-\alpha},
		&&\qquad N_{\mathrm{ind}}+1+N\leq N_{\mathrm{max}}-1.
	\end{alignedat}
	\end{equation}
	The $C^N$ versions replace $N+\alpha$ by $N$, omit the $\La{\cdot}$ factors and the last factor $\ell_q^{-\alpha}$, and use the same ranges with $N_{\mathrm{max}}+1$ in place of $N_{\mathrm{max}}$.

	For $r,N\in\mathbb N_0$ with $r\leq N_R-3j$ and $r+N\leq K_j^{(0)}$, apply \eqref{eq:holder-product-lower-exponent} with $\eta=\alpha_1$ on $|t-t_u|\leq\tau_q$. We obtain
	\begin{equation}\label{eq:Newton-Jacobian-weighted-source-estimate}
	\begin{aligned}
		&\Nrm{\DTL^r\bigl(J_u^{-1}Q[R_{q,j}]\bigr)}_{N+\alpha}\\
		&\quad\lesssim\sum_{M=0}^r\binom{r}{M}\sum_{N_1+N_2=N}\Bigl(
		\nrm{\DTL^MJ_u^{-1}}_{N_1+\alpha}
		\nrm{\DTL^{r-M}Q[R_{q,j}]}_{N_2+\alpha_1}+\nrm{\DTL^MJ_u^{-1}}_{N_1}
		\nrm{\DTL^{r-M}Q[R_{q,j}]}_{N_2+\alpha}\Bigr)\\
		&\quad\lesssim A_{q,j}\tau_q^{-r}
		\Lat{N,N_{R,j}^{(0)}(r),\lambda_q,\ell_q^{-1}}
		\MM{N,N_{R,j}^{(0)}(r),\lambda_q,\ell_q^{-1}}.
	\end{aligned}
	\end{equation}
	Here the last inequality uses \eqref{eq:Newton-scale-separation}, \eqref{eq:Newton-compressible-input-estimate} with H{\"o}lder exponents $\alpha_1$ and $\alpha$, and \eqref{eq:Newton-Jacobian-material-estimates} with $\sigma=-1$, together with $2\alpha_1=\alpha$.

	For the mollified sources, use \eqref{eq:Newton-compressible-input-estimate} for $W[R_{q,j}]$ and \eqref{eq:Newton-Jacobian-weighted-source-estimate} for $J_u^{-1}Q[R_{q,j}]$.
	Set $r:=\min\{M,N_R-3j\}$. Repeat the argument for
	\eqref{eq:initial-Newton-flow-pullback-estimate}--\eqref{eq:initial-Newton-flow-mollification-gain}
	with $F=J_u^{-1}Q[R_{q,j}]$ or $F=W[R_{q,j}]$.
	For the scalar source, take $|t-t_u|\leq2\tau_q/3$. For the vector source, take $t\in\mathbb R$.
	In the scalar case, $\ell_{t,q}<\tau_q/6$ ensures $|t+s-t_u|<\tau_q$ for $|s|\leq\ell_{t,q}$.
	Integration by parts as in \eqref{eq:initial-Newton-flow-material-integration-by-parts} gives
	\begin{equation}\label{eq:Newton-source-flow-mollification-gain}
	\begin{aligned}
		&\nrm{\DTL^M\bigl(\eta_{\ell_{t,q}}^{\mathrm t}\ast_{X_{\overline u_q}}F\bigr)(t,\cdot)}_{N+\alpha}\\
		&\quad\leq\|\partial_s^{M-r}\eta_{\ell_{t,q}}^{\mathrm t}\|_{L^1(\mathbb R)}
		\sup_{|s|\leq\ell_{t,q}}
		\Nrm{(\DTL^rF)(t+s,X_{\overline u_q}(t+s,\cdot;t))}_{N+\alpha}\\
		&\quad\lesssim A_{q,j}\ell_{t,q}^{-(M-r)}\tau_q^{-r}
		\Lat{N,N_R-3j-r,\lambda_q,\ell_q^{-1}}
		\MM{N,N_R-3j-r,\lambda_q,\ell_q^{-1}}\\
		&\quad=A_{q,j}\MM{M,N_R-3j,\tau_q^{-1},\ell_{t,q}^{-1}}
		\Lat{N,N_{R,j}^{(0)}(M),\lambda_q,\ell_q^{-1}}
		\MM{N,N_{R,j}^{(0)}(M),\lambda_q,\ell_q^{-1}}.
	\end{aligned}
	\end{equation}
	Here $N_R-3j-r=N_{R,j}^{(0)}(M)$.
	The source and flow estimates apply since $r\leq N_R-3j\leq N_{\mathrm{ind}}-1$, $r+N\leq K_j^{(0)}$, and $K_j^{(0)}+2\leq N_{\mathrm{max}}$. Combining \eqref{eq:time-cutoff-derivative-estimate} with \eqref{eq:Newton-source-flow-mollification-gain}, $\|\chi_\kappa\|_N\lesssim\tau_q^{-N}$, and
	$\DTL(\chi_\kappa(\Phi_u))=0$, we obtain \eqref{eq:Newton-mollified-source-estimate}.
	The cutoff factors $\theta_u$ and $\chi_\kappa^2(\Phi_u)$ give the support inclusion \eqref{eq:Newton-mollified-source-support}.
\end{proof}
\begin{prop}
	\label{prop:Newton-three-dimensional-acoustic-estimates}
	Fix $0\leq j<J_*$. Assume the background estimates in Proposition~\ref{prop:uniform-mollification-estimates}, \eqref{eq:Newton-scale-separation},
	and the source estimates in \eqref{eq:Newton-mollified-source-estimate} at the fixed H{\"o}lder exponents. Set $\mcI_u:=[t_{u-1},t_u+2\tau_q/3]$.
	The time supremum is taken over the indicated interval, or over $\mathbb R$ if no interval is specified. For $a>a_{\mathrm{Newton}}$ and
	$M,N\in\mathbb N_0$ in the indicated ranges, we have
	\begin{alignat}{2}
		\nrm{\DTL^M\widetilde V_{q,j+1,u,\kappa}}_{N+\alpha,\mcI_u}
		&\lesssim\mathrlap{\tau_qA_{q,j}\MM{M,N_R-3j-1,\tau_q^{-1},\ell_{t,q}^{-1}}\Lat{N,N_{R,j}^{(1)}(M),\lambda_q,\ell_q^{-1}}}\notag\\*
		&\quad\cdot\MM{N,N_{R,j}^{(1)}(M),\lambda_q,\ell_q^{-1}},
		&&\qquad 0\leq M\leq N_{\mathrm{ind}}+1,\quad M+N\leq K_j^{(1)}.
		\label{eq:Newton-raw-Piola-estimate}\\*[2pt]
		\nrm{\DTL^M\varrho[\widetilde V_{q,j+1,u,\kappa}]}_{N+\alpha,\mcI_u}
		&\lesssim\mathrlap{A_{q,j}\MM{M,N_R-3j-2,\tau_q^{-1},\ell_{t,q}^{-1}}\Lat{N,N_{R,j}^{(2)}(M),\lambda_q,\ell_q^{-1}}}\notag\\*
		&\quad\cdot\MM{N,N_{R,j}^{(2)}(M),\lambda_q,\ell_q^{-1}},
		&&\qquad 0\leq M\leq N_{\mathrm{ind}},\quad M+N\leq K_j^{(2)}.
		\label{eq:Newton-raw-density-estimate}\\*[2pt]
		\nrm{\DTL^Mw[\widetilde V_{q,j+1,u,\kappa}]}_{N+\alpha,\mcI_u}
		&\lesssim\mathrlap{A_{q,j}\MM{M,N_R-3j-2,\tau_q^{-1},\ell_{t,q}^{-1}}\Lat{N,N_{R,j}^{(2)}(M),\lambda_q,\ell_q^{-1}}}\notag\\*
		&\quad\cdot\MM{N,N_{R,j}^{(2)}(M),\lambda_q,\ell_q^{-1}},
		&&\qquad 0\leq M\leq N_{\mathrm{ind}},\quad M+N\leq K_j^{(2)}.
		\label{eq:Newton-raw-momentum-estimate}
	\end{alignat}
	The mean correction satisfies
	\begin{alignat}{2}
		\nrm{\DTL^M(\overline V_{q,j+1,u,\kappa}\phi_{u,\kappa})}_{N+\alpha,\mcI_u}
		&\lesssim\mathrlap{\tau_q^2\Lambda_q A_{q,j}\MM{M,N_R-3j-1,\tau_q^{-1},\ell_{t,q}^{-1}}\Lat{N,N_{R,j}^{(1)}(M),\lambda_q,\ell_q^{-1}}}\notag\\*
		&\quad\cdot\MM{N,N_{R,j}^{(1)}(M),\lambda_q,\ell_q^{-1}},
		&&\qquad 0\leq M\leq N_{\mathrm{ind}}+1,\quad M+N\leq K_j^{(1)}.
		\label{eq:Newton-mean-potential-estimate}\\*[2pt]
		\nrm{\DTL^M\varrho[\overline V_{q,j+1,u,\kappa}\phi_{u,\kappa}]}_{N+\alpha,\mcI_u}
		&\lesssim\mathrlap{\tau_q\Lambda_q A_{q,j}\MM{M,N_R-3j-1,\tau_q^{-1},\ell_{t,q}^{-1}}\Lat{N,N_{R,j}^{(1)}(M),\lambda_q,\ell_q^{-1}}}\notag\\*
		&\quad\cdot\MM{N,N_{R,j}^{(1)}(M),\lambda_q,\ell_q^{-1}},
		&&\qquad 0\leq M\leq N_{\mathrm{ind}},\quad M+N\leq K_j^{(1)}.
		\label{eq:Newton-mean-density-estimate}\\*[2pt]
		\nrm{\DTL^Mw[\overline V_{q,j+1,u,\kappa}\phi_{u,\kappa}]}_{N+\alpha,\mcI_u}
		&\lesssim\mathrlap{\tau_q\Lambda_q A_{q,j}\MM{M,N_R-3j-2,\tau_q^{-1},\ell_{t,q}^{-1}}\Lat{N,N_{R,j}^{(2)}(M),\lambda_q,\ell_q^{-1}}}\notag\\*
		&\quad\cdot\MM{N,N_{R,j}^{(2)}(M),\lambda_q,\ell_q^{-1}},
		&&\qquad 0\leq M\leq N_{\mathrm{ind}},\quad M+N\leq K_j^{(2)}.
		\label{eq:Newton-mean-momentum-estimate}
	\end{alignat}
	At $M=0$, \eqref{eq:Newton-mean-potential-estimate} also holds for $N\leq N_{\mathrm{max}}-2$, and \eqref{eq:Newton-mean-density-estimate}--\eqref{eq:Newton-mean-momentum-estimate} for $N\leq N_{\mathrm{max}}-3$.
	For $M\geq1$, \eqref{eq:Newton-mean-potential-estimate} also holds for $M+N\leq N_{\mathrm{max}}-1$, and \eqref{eq:Newton-mean-density-estimate}--\eqref{eq:Newton-mean-momentum-estimate} for $M+N\leq N_{\mathrm{max}}-2$, with the same material-order ranges.
	The corrected and total potentials satisfy
	\begin{equation}\label{eq:Newton-total-compressible-potential-estimate}
		\begin{aligned}
			\nrm{\DTL^MV_{q,j+1}\comp}_{N+\alpha}
			&\lesssim A_{q,j}\tau_q\MM{M,N_R-3j-1,\tau_q^{-1},\ell_{t,q}^{-1}}
			\Lat{N,N_{R,j}^{(1)}(M),\lambda_q,\ell_q^{-1}}
			\MM{N,N_{R,j}^{(1)}(M),\lambda_q,\ell_q^{-1}},
		\end{aligned}
	\end{equation}
	for $0\leq M\leq N_{\mathrm{ind}}+1$ and $M+N\leq K_j^{(1)}$. Estimate~\eqref{eq:Newton-total-compressible-potential-estimate} also holds for $V_{q,j+1,u,\kappa}^{\circ}$ on $\mcI_u$ and for $V_{q,j+1,u,\kappa}\comp$
	over all times. 
	
	For $0\leq M\leq N_{\mathrm{ind}}$ and $M+N\leq K_j^{(2)}$, the total density and momentum perturbations satisfy
	\begin{align}
		\nrm{\DTL^M\trho_{q,j+1}\comp}_{N+\alpha}
		&\lesssim A_{q,j}\MM{M,N_R-3j-2,\tau_q^{-1},\ell_{t,q}^{-1}}
		\Lat{N,N_{R,j}^{(2)}(M),\lambda_q,\ell_q^{-1}}\MM{N,N_{R,j}^{(2)}(M),\lambda_q,\ell_q^{-1}},
		\label{eq:Newton-total-compressible-density-estimate}\\
		\nrm{\DTL^M\tm_{q,j+1}\comp}_{N+\alpha}
		&\lesssim A_{q,j}\MM{M,N_R-3j-2,\tau_q^{-1},\ell_{t,q}^{-1}}
		\Lat{N,N_{R,j}^{(2)}(M),\lambda_q,\ell_q^{-1}}\MM{N,N_{R,j}^{(2)}(M),\lambda_q,\ell_q^{-1}}.
		\label{eq:Newton-total-compressible-momentum-estimate}
	\end{align}
	The same estimates hold for the corresponding local perturbations $\trho_{q,j+1,u,\kappa}\comp$ and $\tm_{q,j+1,u,\kappa}\comp$.	
	
	The  matrix $\mathcal S_{q,j+1,u,\kappa}$ and the remaining vector field $H_{q,j+1,u,\kappa}^{\mathrm{int}}$ satisfy
	\begin{alignat}{2}
		\nrm{\DTL^M\mathcal S_{q,j+1,u,\kappa}}_{N+\alpha}
		&\lesssim\mathrlap{\tau_q\Lambda_q A_{q,j}\MM{M,N_R-3j-2,\tau_q^{-1},\ell_{t,q}^{-1}}\Lat{N,N_{R,j}^{(2)}(M),\lambda_q,\ell_q^{-1}}}\notag\\
		&\quad\cdot\MM{N,N_{R,j}^{(2)}(M),\lambda_q,\ell_q^{-1}},
		&&\qquad 0\leq M\leq N_{\mathrm{ind}}-1,\quad M+N\leq K_j^{(2)},\label{eq:Newton-interior-acoustic-residual-estimate}\\[2pt]
		\nrm{\DTL^MH_{q,j+1,u,\kappa}^{\mathrm{int}}}_{N+\alpha}
		&\lesssim\mathrlap{\Lambda_q A_{q,j}\MM{M,N_R-3j-3,\tau_q^{-1},\ell_{t,q}^{-1}}\Lat{N,N_{R,j}^{(3)}(M),\lambda_q,\ell_q^{-1}}}\notag\\
		&\quad\cdot\MM{N,N_{R,j}^{(3)}(M),\lambda_q,\ell_q^{-1}},
		&&\qquad 0\leq M\leq N_{\mathrm{ind}}-1,\quad M+N\leq K_j^{(3)}.\label{eq:Newton-interior-vector-source-estimate}
	\end{alignat}
	Each cutoff residual satisfies
	\begin{equation}\label{eq:Newton-compressible-cutoff-residual-estimate}
		\begin{alignedat}{2}
			\nrm{\DTL^MF_{q,j+1,u,\kappa}^{\mathrm{cut}}}_{N+\alpha}
			&\lesssim\mathrlap{A_{q,j}\tau_q^{-1}\MM{M,N_R-3j-2,\tau_q^{-1},\ell_{t,q}^{-1}}\Lat{N,N_{R,j}^{(2)}(M),\lambda_q,\ell_q^{-1}}}\\
			&\quad\cdot\MM{N,N_{R,j}^{(2)}(M),\lambda_q,\ell_q^{-1}},
			&&\qquad 0\leq M\leq N_{\mathrm{ind}},\quad M+N\leq K_j^{(2)}.
		\end{alignedat}
	\end{equation}
	Moreover, $F_{q,j+1,u,\kappa}^{\mathrm{cut}}$ satisfies $\supp_tF_{q,j+1,u,\kappa}^{\mathrm{cut}}\subset\supp\partial_t\theta_u$.
	
	The interior matrix satisfies
	\begin{equation}\label{eq:Newton-interior-stress-potential-estimate}
		\begin{alignedat}{2}
			\nrm{\DTL^MR_{q,j+1}\comp}_{N+\alpha}
			&\lesssim\mathrlap{\tau_q\Lambda_q(\tau_q\ell_q^{-1})^\alpha A_{q,j}\MM{M,N_R-3j-3,\tau_q^{-1},\ell_{t,q}^{-1}}\Lat{N,N_{R,j+1}^{(0)}(M),\lambda_q,\ell_q^{-1}}}\\
			&\quad\cdot\MM{N,N_{R,j+1}^{(0)}(M),\lambda_q,\ell_q^{-1}},
			&&\qquad 0\leq M\leq N_{\mathrm{ind}}-1,\quad M+N\leq K_j^{(3)}.
		\end{alignedat}
	\end{equation}
	The cutoff stresses and the vector potential $W[S_{q,j+1,u}^{\mathrm{cut,skew}}]$ satisfy
	\begin{equation}\label{eq:Newton-cutoff-stress-potential-estimate}
		\begin{alignedat}{2}
			&\sup_u\Bigl(\nrm{\DTL^MS_{q,j+1,u}^{\mathrm{cut}}}_{N+\alpha}
			+\nrm{\DTL^MS_{q,j+1,u}^{\mathrm{cut,sym}}}_{N+\alpha}
			+\nrm{\DTL^MS_{q,j+1,u}^{\mathrm{cut,skew}}}_{N+\alpha}
			+\Nrm{\DTL^MW\left[S_{q,j+1,u}^{\mathrm{cut,skew}}\right]}_{N+\alpha}\Bigr)\\
			&\qquad\lesssim A_{q,j}(\tau_q\ell_q^{-1})^\alpha\MM{M,N_R-3j-2,\tau_q^{-1},\ell_{t,q}^{-1}}\Lat{N,N_{R,j}^{(2)}(M),\lambda_q,\ell_q^{-1}}\\
			&\qquad\quad\cdot\MM{N,N_{R,j}^{(2)}(M),\lambda_q,\ell_q^{-1}},
			\hspace{80pt} 0\leq M\leq N_{\mathrm{ind}},\quad M+N\leq K_j^{(2)}.
		\end{alignedat}
	\end{equation}
	For each fixed exponent $\alpha_r$, the same estimates hold with $\alpha$ replaced by $\alpha_r$.
	The implicit constants may depend on $\rho^*,\mathfrak M,N_{\mathrm{max}},\alpha,\alpha_r$.
	They are independent of $q,j,u,\kappa,a,C_*$ and $A_{q,j}$.
\end{prop}

\begin{proof}
	We first estimate $\widetilde V_{q,j+1,u,\kappa}$ and the perturbations $\varrho[\widetilde V_{q,j+1,u,\kappa}]$ and $w[\widetilde V_{q,j+1,u,\kappa}]$.
	We then estimate $\overline V_{q,j+1,u,\kappa}\phi_{u,\kappa}$, followed by the stresses $R_{q,j+1}\comp$ and $S_{q,j+1,u}^{\mathrm{cut}}$.
	We work on $\mcI_u$ and omit the interval subscript. For fields supported in $\supp\theta_u$, the norms agree with those over $\mathbb R$.
	Since $a>a_{\mathrm{Newton}}\geq a_{\mathrm{flow}}$ and $|t-t_u|\leq\tau_q$ on $\mcI_u$,
	apply Corollary~\ref{cor:mollified-velocity-flow-estimates} with $\mcI=\mathbb R$, $t_0=t_u$, and $h=\tau_q$.
	The flow estimates cover all required derivatives since $K_j^{(0)}+2\leq N_{\mathrm{max}}$. For $J_u$, \eqref{eq:Newton-Jacobian-material-estimates} with $\sigma=1$ and $N=0$ gives $1/2\leq J_u\leq2$ for sufficiently large $a$.
	The frozen sound speed is constant in $(t,x)$ and satisfies
	$$
	0<\min_{\rho\in I_{\rho^*}}p'(\rho)
	\leq c_{q,u,\kappa}^2
	\leq\max_{\rho\in I_{\rho^*}}p'(\rho).
	$$

	To prove \eqref{eq:Newton-cone-derivative-estimate}--\eqref{eq:Newton-cone-gradient-derivative-estimate}, estimate the Kirchhoff solution $I_{q,u,\kappa}[G_{q,j,u,\kappa}^{\mathrm{mol}}]$ and its gradient in flow coordinates using \eqref{eq:Newton-mollified-source-estimate}.
	All material derivatives of $G_{q,j,u,\kappa}^{\mathrm{mol}}$ vanish near $t_{u-1}$. For $t\in\mcI_u$, apply \eqref{eq:Newton-Lagrangian-cone-holder-estimates} with $0\leq t-t_{u-1}\leq5\tau_q/3$.
	For composition with $X_u$ and $\Phi_u$, use Lemmas~\ref{lem:BDSV-composition} and \ref{lem:holder-composition-general}. Together with Corollary~\ref{cor:mollified-velocity-flow-estimates}, these give the first inequalities in \eqref{eq:Newton-cone-derivative-estimate}--\eqref{eq:Newton-cone-gradient-derivative-estimate}.
	\begin{align}
		&\nrm{\DTL^MI_{q,u,\kappa}[G_{q,j,u,\kappa}^{\mathrm{mol}}]}_{N+\alpha}\notag\\
		&\quad\lesssim\tau_q^2\sum_{n=0}^N
		\MM{N-n,N_{R,j}^{(0)}(M),\lambda_q,\ell_q^{-1}}
		\Bigl[\nrm{\DTL^MG_{q,j,u,\kappa}^{\mathrm{mol}}}_{n+\alpha}+
		\La{N,N_{R,j}^{(0)}(M),\lambda_q,\ell_q^{-1}}
		\nrm{\DTL^MG_{q,j,u,\kappa}^{\mathrm{mol}}}_{n+\alpha_1}\Bigr]\notag\\
		&\quad\lesssim A_{q,j}\tau_q^{2}\MM{M,N_R-3j,\tau_q^{-1},\ell_{t,q}^{-1}}
		\Lat{N,N_{R,j}^{(0)}(M),\lambda_q,\ell_q^{-1}}\MM{N,N_{R,j}^{(0)}(M),\lambda_q,\ell_q^{-1}},\label{eq:Newton-cone-derivative-estimate}\\
		&\nrm{\DTL^M\nabla_y
		\widehat I_{q,u,\kappa}[G_{q,j,u,\kappa}^{\mathrm{mol}}](t,\Phi_u)}_{N+\alpha}
			\notag\\
		&\quad\lesssim\tau_q^2\sum_{n=0}^N
		\MM{N-n,N_{R,j}^{(1)}(M),\lambda_q,\ell_q^{-1}}
		\Bigl[\nrm{\DTL^{M+1}G_{q,j,u,\kappa}^{\mathrm{mol}}}_{n+\alpha}+
		\La{N,N_{R,j}^{(1)}(M),\lambda_q,\ell_q^{-1}}
		\nrm{\DTL^{M+1}G_{q,j,u,\kappa}^{\mathrm{mol}}}_{n+\alpha_1}\Bigr]\notag\\
		&\quad\lesssim A_{q,j}\tau_q^2\MM{M+1,N_R-3j,\tau_q^{-1},\ell_{t,q}^{-1}}
		\Lat{N,N_{R,j}^{(0)}(M+1),\lambda_q,\ell_q^{-1}}\MM{N,N_{R,j}^{(0)}(M+1),\lambda_q,\ell_q^{-1}}\notag\\
		&\quad= A_{q,j}\tau_q\MM{M,N_R-3j-1,\tau_q^{-1},\ell_{t,q}^{-1}}
		\Lat{N,N_{R,j}^{(1)}(M),\lambda_q,\ell_q^{-1}}\MM{N,N_{R,j}^{(1)}(M),\lambda_q,\ell_q^{-1}}.\label{eq:Newton-cone-gradient-derivative-estimate}
	\end{align}
	For both sums in \eqref{eq:Newton-cone-derivative-estimate}--\eqref{eq:Newton-cone-gradient-derivative-estimate}, we used \eqref{eq:Newton-mollified-source-estimate} at exponents $\alpha$ and $\alpha_1$, with $2\alpha_1=\alpha$.
	Estimate \eqref{eq:Newton-cone-derivative-estimate} holds for $0\leq M\leq N_{\mathrm{max}}$ and $\min\{M,N_R-3j\}+N\leq K_j^{(0)}$. Estimate \eqref{eq:Newton-cone-gradient-derivative-estimate} holds for $0\leq M\leq N_{\mathrm{max}}-1$ and $\min\{M+1,N_R-3j\}+N\leq K_j^{(0)}$.
	
	For $0\leq M\leq N_{\mathrm{ind}}+1$ and $M+N\leq K_j^{(1)}$, \eqref{def:Newton-raw-Piola-potential} gives
	$$
	\begin{aligned}
		&\nrm{\DTL^M\widetilde V_{q,j+1,u,\kappa}}_{N+\alpha}\\
		&\quad\lesssim
		\sum_{r=0}^M\sum_{N_1+N_2=N}\Bigl[
		\nrm{\DTL^r\mathbf C_u^T}_{N_1}
		\nrm{\DTL^{M-r}\bigl(\nabla_y\widehat I_{q,u,\kappa}[G_{q,j,u,\kappa}^{\mathrm{mol}}](t,\Phi_u)\bigr)}_{N_2+\alpha}+
		\nrm{\DTL^r\mathbf C_u^T}_{N_1+\alpha}
		\nrm{\DTL^{M-r}\bigl(\nabla_y\widehat I_{q,u,\kappa}[G_{q,j,u,\kappa}^{\mathrm{mol}}](t,\Phi_u)\bigr)}_{N_2+\alpha_1}\Bigr]\\
		&\quad\lesssim\tau_q A_{q,j}
		\MM{M,N_R-3j-1,\tau_q^{-1},\ell_{t,q}^{-1}}
		\Lat{N,N_{R,j}^{(1)}(M),\lambda_q,\ell_q^{-1}}
		\MM{N,N_{R,j}^{(1)}(M),\lambda_q,\ell_q^{-1}}.
	\end{aligned}
	$$
	For the second inequality, we used the gradient estimate in \eqref{eq:Newton-cone-gradient-derivative-estimate}
	and the following cofactor estimates. For $1\leq M\leq N_{\mathrm{ind}}+1$, \eqref{eq:mollified-flow-material-estimates} gives
	\begin{equation}\label{eq:Newton-raw-cofactor-material-estimates}
	\begin{aligned}
		\nrm{\DTL^M\mathbf C_u}_N
		&\lesssim\lambda_q^{-1}\Lambda_q^M
		\MM{N+1,N_{\mathrm{ind}}-M+1,\lambda_q,\ell_q^{-1}}\\
		&\leq\MM{M,N_R-3j-1,\tau_q^{-1},\ell_{t,q}^{-1}}
		\MM{N,N_{R,j}^{(1)}(M),\lambda_q,\ell_q^{-1}},\\
		\nrm{\DTL^M\mathbf C_u}_{N+\alpha}
		&\lesssim\lambda_q^{-1}\Lambda_q^M
		\La{N+1,N_{\mathrm{ind}}-M+1,\lambda_q,\ell_q^{-1}}
		\MM{N+1,N_{\mathrm{ind}}-M+1,\lambda_q,\ell_q^{-1}}\\
		&\leq\MM{M,N_R-3j-1,\tau_q^{-1},\ell_{t,q}^{-1}}
		\La{N,N_{R,j}^{(1)}(M),\lambda_q,\ell_q^{-1}}\MM{N,N_{R,j}^{(1)}(M),\lambda_q,\ell_q^{-1}}.
	\end{aligned}
	\end{equation}
	The $C^N$ estimate holds for $M+N\leq N_{\mathrm{max}}$, and the H{\"o}lder estimate for $M+N\leq N_{\mathrm{max}}-1$.
	At $M=N_{\mathrm{ind}}+1$, we use $\lambda_q^{-1}\Lambda_q\ell_q^{-1}=\delta_q^{\frac12}\ell_q^{-1}\leq\ell_{t,q}^{-1}$ from \eqref{eq:material-mollification-spatial-comparisons} and $M>N_R-3j-1$.
	The H{\"o}lder estimate also holds at the fixed smaller exponents. 	For the cofactor terms with $r=0$, use \eqref{eq:mollified-flow-close-identity}.
	Applying the gradient estimate in \eqref{eq:Newton-cone-gradient-derivative-estimate} and
	\eqref{eq:Newton-raw-cofactor-material-estimates} to the preceding sum proves
	\eqref{eq:Newton-raw-Piola-estimate} through $M=N_{\mathrm{ind}}+1$.
	
	For $\varrho[\widetilde V_{q,j+1,u,\kappa}]$, \eqref{eq:Newton-cone-material-identities} and \eqref{eq:Newton-sharp-density-identity} give
	$$
	\DTL^M\varrho[\widetilde V_{q,j+1,u,\kappa}]
	=c_{q,u,\kappa}^{-2}\sum_{r=0}^M\binom Mr(\DTL^rJ_u)\left(\DTL^{M-r+2}I_{q,u,\kappa}[G_{q,j,u,\kappa}^{\mathrm{mol}}]-\DTL^{M-r}G_{q,j,u,\kappa}^{\mathrm{mol}}\right).
	$$
	We use $
	\tau_q^2\MM{M+2,N_R-3j,\tau_q^{-1},\ell_{t,q}^{-1}}
	=\MM{M,N_R-3j-2,\tau_q^{-1},\ell_{t,q}^{-1}}$, and 
	$N_{R,j}^{(0)}(M+2)=N_{R,j}^{(2)}(M).
	$
	For $0\leq M\leq N_{\mathrm{ind}}$ and $M+N\leq K_j^{(2)}$,
	\eqref{eq:Newton-mollified-source-estimate} applies through material order $M+2$ and spatial order $N$, since $M+2\leq N_{\mathrm{max}}$
	and $M+N+2\leq K_j^{(0)}$.
	Applying \eqref{eq:holder-product-lower-exponent} with $\eta=\alpha_1$, we have
	$$
	\begin{aligned}
		\nrm{\DTL^M\varrho[\widetilde V_{q,j+1,u,\kappa}]}_{N+\alpha}&\lesssim\sum_{r=0}^M\binom Mr\Bigl(
		\nrm{(\DTL^rJ_u)\DTL^{M-r+2}I_{q,u,\kappa}[G_{q,j,u,\kappa}^{\mathrm{mol}}]}_{N+\alpha}+\nrm{(\DTL^rJ_u)\DTL^{M-r}G_{q,j,u,\kappa}^{\mathrm{mol}}}_{N+\alpha}\Bigr)\\
		&\lesssim\sum_{r=0}^M\sum_{N_1+N_2=N}\Bigl[
		\nrm{\DTL^rJ_u}_{N_1}
		\Bigl(\nrm{\DTL^{M-r+2}I_{q,u,\kappa}[G_{q,j,u,\kappa}^{\mathrm{mol}}]}_{N_2+\alpha}
		+\nrm{\DTL^{M-r}G_{q,j,u,\kappa}^{\mathrm{mol}}}_{N_2+\alpha}\Bigr)\\
		&\hspace{50pt}\quad+
		\nrm{\DTL^rJ_u}_{N_1+\alpha}
		\Bigl(\nrm{\DTL^{M-r+2}I_{q,u,\kappa}[G_{q,j,u,\kappa}^{\mathrm{mol}}]}_{N_2+\alpha_1}
		+\nrm{\DTL^{M-r}G_{q,j,u,\kappa}^{\mathrm{mol}}}_{N_2+\alpha_1}\Bigr)\Bigr]\\
		&\lesssim A_{q,j}\sum_{r=0}^M
		\MM{r,N_R-3j-2,\tau_q^{-1},\ell_{t,q}^{-1}}
		\MM{M-r,N_R-3j-2,\tau_q^{-1},\ell_{t,q}^{-1}}\\
		&\qquad\quad\cdot\Lat{N,N_{R,j}^{(2)}(M),\lambda_q,\ell_q^{-1}}
		\MM{N,N_{R,j}^{(2)}(M),\lambda_q,\ell_q^{-1}}\\
		&\lesssim A_{q,j}\MM{M,N_R-3j-2,\tau_q^{-1},\ell_{t,q}^{-1}}\Lat{N,N_{R,j}^{(2)}(M),\lambda_q,\ell_q^{-1}}
		\MM{N,N_{R,j}^{(2)}(M),\lambda_q,\ell_q^{-1}}.
	\end{aligned}
	$$
	In the third inequality, we used \eqref{eq:Newton-mollified-source-estimate} for $G_{q,j,u,\kappa}^{\mathrm{mol}}$ and \eqref{eq:Newton-cone-derivative-estimate} for $I_{q,u,\kappa}[G_{q,j,u,\kappa}^{\mathrm{mol}}]$.
	We estimated $J_u$ and its material derivatives by \eqref{eq:Newton-Jacobian-material-estimates} with $\sigma=1$.
	Since $2\alpha_1=\alpha$, the product of the H{\"o}lder factors in each mixed term is at most $\Lat{N,N_{R,j}^{(2)}(M),\lambda_q,\ell_q^{-1}}$.
	Thus \eqref{eq:Newton-raw-density-estimate} holds.
	
	For $w[\widetilde V_{q,j+1,u,\kappa}]$, \eqref{eq:Newton-raw-momentum-Piola-cancellation} gives
	$$
	\begin{aligned}
		&\nrm{\DTL^Mw[\widetilde V_{q,j+1,u,\kappa}]}_{N+\alpha}\\
		&\quad\leq\nrm{\DTL^M\bigl(\overline u_q\varrho[\widetilde V_{q,j+1,u,\kappa}]\bigr)}_{N+\alpha}+\nrm{\DTL^M\bigl(\mathbf C_u^T\nabla_y\widehat I_{q,u,\kappa}[\DTL G_{q,j,u,\kappa}^{\mathrm{mol}}](t,\Phi_u)\bigr)}_{N+\alpha}\\
		&\quad\lesssim\sum_{r=0}^M\sum_{N_1+N_2=N}\Bigl[
		\Nrm{\DTL^r\overline u_q}_{N_1}\Nrm{\DTL^{M-r}\varrho[\widetilde V_{q,j+1,u,\kappa}]}_{N_2+\alpha}+\Nrm{\DTL^r\overline u_q}_{N_1+\alpha}\Nrm{\DTL^{M-r}\varrho[\widetilde V_{q,j+1,u,\kappa}]}_{N_2+\alpha_1}\\
		&\qquad+\Nrm{\DTL^r\mathbf C_u^T}_{N_1}\Nrm{\DTL^{M-r}\bigl(\nabla_y\widehat I_{q,u,\kappa}[\DTL G_{q,j,u,\kappa}^{\mathrm{mol}}](t,\Phi_u)\bigr)}_{N_2+\alpha}+\Nrm{\DTL^r\mathbf C_u^T}_{N_1+\alpha}\Nrm{\DTL^{M-r}\bigl(\nabla_y\widehat I_{q,u,\kappa}[\DTL G_{q,j,u,\kappa}^{\mathrm{mol}}](t,\Phi_u)\bigr)}_{N_2+\alpha_1}\Bigr]\\
		&\lesssim A_{q,j}\MM{M,N_R-3j-2,\tau_q^{-1},\ell_{t,q}^{-1}}\Lat{N,N_{R,j}^{(2)}(M),\lambda_q,\ell_q^{-1}}\MM{N,N_{R,j}^{(2)}(M),\lambda_q,\ell_q^{-1}}.
	\end{aligned}
	$$
	For the last inequality, we used \eqref{eq:Newton-raw-density-estimate} for the density factors.
	For the Kirchhoff gradient, we used \eqref{eq:Newton-commuted-cone-gradient} and \eqref{eq:Newton-cone-gradient-derivative-estimate} at material order $M-r+1$.
	For $\overline u_q$, we used \eqref{est:uniform-zero-velocity}--\eqref{est:uniform-om-over-orho} and \eqref{est:holder-om-over-orho}.
	For $\mathbf C_u^T$, we applied \eqref{eq:mollified-flow-close-identity} at material order zero and \eqref{eq:Newton-raw-cofactor-material-estimates} at positive material orders.
	The gradient term requires source estimates through material order $M+2$ and flow coefficient estimates through order $M\leq N_{\mathrm{ind}}$.
	These estimates give \eqref{eq:Newton-raw-momentum-estimate} for $0\leq M\leq N_{\mathrm{ind}}$ and $M+N\leq K_j^{(2)}$.
	
	We next estimate the coefficients in \eqref{eq:Newton-frozen-flow-coefficient-estimates}, used in the mean correction and the interior and cutoff stresses. Except for $p'(\orho_q)-c_{q,u,\kappa}^2$ at material order zero, spatial norms are taken over $\mathbb R^3$.
	For that coefficient, take the spatial norms over $X_u(t,B(x_\kappa,C_{\mathrm{supp}}\tau_q))$, which contains the supports of $\widetilde V_{q,j+1,u,\kappa}$
	and $\phi_{u,\kappa}$ on $\mcI_u$.
	
	The coefficient $p'(\orho_q)-c_{q,u,\kappa}^2$ satisfies
	$$
	\begin{aligned}
		p'(\orho_q(t,x))-c_{q,u,\kappa}^2
		&=p'(\orho_q(t,x))-p'(\orho_q(t_u,\Phi_u(t,x)))+p'(\orho_q(t_u,\Phi_u(t,x)))-p'(\orho_q(t_u,x_\kappa)),\\
		\nrm{p'(\orho_q)-c_{q,u,\kappa}^2}_0
		&\lesssim\tau_q\left(\nrm{\DTL\orho_q}_0+\nrm{\nabla\orho_q}_0\right)
		\lesssim\tau_q\delta_q(\Lambda_q+\lambda_q)
		\lesssim\tau_q\Lambda_q.
	\end{aligned}
	$$
	In the second inequality, we used \eqref{est:uniform-orho} for $\DTL\orho_q$ and $\nabla\orho_q$.
	For the positive-order spatial derivatives of $p'(\orho_q)$, apply
	\eqref{eq:faa-di-bruno} and use \eqref{est:uniform-orho}.
	Since $4b\alpha<\beta$, \eqref{def:lambda-delta} and \eqref{def:tau-q} give $\delta_q\leq\tau_q\Lambda_q$ after increasing $a$.
	Interpolate the resulting $C^N$ and $C^{N+1}$ estimates. For $J_u-1$, use \eqref{eq:Newton-Jacobian-material-estimates} with $\sigma=1$.
	For $\nabla\Phi_u-\Id$ and $\mathbf C_u-\Id$, use \eqref{eq:mollified-flow-close-identity}. For $0\leq N\leq N_{\mathrm{max}}-2$, we obtain
	\begin{equation}\label{eq:Newton-frozen-flow-coefficient-estimates}
		\begin{aligned}
			&\nrm{p'(\orho_q)-c_{q,u,\kappa}^2}_{N+\alpha}
			+\nrm{J_u-1}_{N+\alpha}
			+\nrm{\nabla\Phi_u-\Id}_{N+\alpha}
			+\nrm{\mathbf C_u-\Id}_{N+\alpha}\\
			&\quad\lesssim\tau_q\Lambda_q
			\La{N,N_{R,j}^{(0)}(0),\lambda_q,\ell_q^{-1}}\MM{N,N_{R,j}^{(0)}(0),\lambda_q,\ell_q^{-1}}.
		\end{aligned}
	\end{equation}
	Using
	$
	J_u\nabla\Phi_u-\mathbf C_u
	=J_u(\nabla\Phi_u-\Id)+(J_u-1)\Id-(\mathbf C_u-\Id)
	$
	and \eqref{eq:Newton-frozen-flow-coefficient-estimates}, we obtain
	$$
	\nrm{J_u\nabla\Phi_u-\mathbf C_u}_{N+\alpha}
	\lesssim\tau_q\Lambda_q\La{N,N_{R,j}^{(0)}(0),\lambda_q,\ell_q^{-1}}\MM{N,N_{R,j}^{(0)}(0),\lambda_q,\ell_q^{-1}}.
	$$
	
	For $M\geq1$, material differentiation gives
	$$
	\DTL^M(J_u\nabla\Phi_u-\mathbf C_u)
	=\sum_{r=0}^M\binom Mr(\DTL^rJ_u)(\DTL^{M-r}\nabla\Phi_u)-\DTL^M\mathbf C_u.
	$$
	For $\nabla\Phi_u$, use \eqref{eq:mollified-flow-close-identity} and \eqref{eq:mollified-flow-material-estimates}.
	Together with \eqref{eq:Newton-Jacobian-material-estimates} with $\sigma=1$, \eqref{eq:Newton-raw-cofactor-material-estimates}, and \eqref{eq:holder-product-tame}, these give
	\begin{equation}\label{eq:Newton-cofactor-material-estimates}
		\begin{aligned}
			\nrm{\DTL^M\mathbf C_u}_N+\nrm{\DTL^M(J_u\nabla\Phi_u-\mathbf C_u)}_N
			&\lesssim\Lambda_q^M\MM{N,N_{\mathrm{ind}}-M,\lambda_q,\ell_q^{-1}},\\
			\nrm{\DTL^M\mathbf C_u}_{N+\alpha}+\nrm{\DTL^M(J_u\nabla\Phi_u-\mathbf C_u)}_{N+\alpha}
			&\lesssim\Lambda_q^M\La{N,N_{\mathrm{ind}}-M,\lambda_q,\ell_q^{-1}}\MM{N,N_{\mathrm{ind}}-M,\lambda_q,\ell_q^{-1}}.
		\end{aligned}
	\end{equation}
	Both estimates hold for $1\leq M\leq N_{\mathrm{ind}}$. The $C^N$ estimate is valid for $M+N\leq N_{\mathrm{max}}$, while the H{\"o}lder estimate is valid for $M+N\leq N_{\mathrm{max}}-1$.
	
	For $\nabla J_u$ and $\nabla\DTL^rJ_u$, use the $C^{N+1}$ versions of \eqref{eq:Newton-Jacobian-material-estimates}.
	For $1\leq r\leq N_{\mathrm{ind}}$, apply Lemma~\ref{lem:appendix-scale-dependent-commuted-gradient} at order $r-1$ with
	$$
	\begin{aligned}
		G&=\DTL J_u,& \upsilon&=\overline u_q,& K_G=N_t&=N_{\mathrm{ind}}-1,& K_\upsilon&=N_{\mathrm{ind}},& K^\sharp&=N_{\mathrm{max}}-1,\\
		\lambda&=\lambda_q,& \widetilde\lambda&=\ell_q^{-1},& \mu=\tilde\mu&=C\Lambda_q,& \mC_G&=\delta_q\Lambda_q,& \mC_\upsilon&=C\delta_q^{\frac12}.
	\end{aligned}
	$$
	Choose $C\geq1$ to cover \eqref{est:uniform-grad-om-over-orho}. The velocity estimate and the $C^N$ estimates in \eqref{eq:Newton-Jacobian-material-estimates} verify the hypotheses, since $\mC_\upsilon\lambda=\mu$.
	For $r+N\leq N_{\mathrm{max}}-1$, the lemma gives
	$$
	\nrm{\DTL^{r-1}\nabla(\DTL J_u)}_N
	\lesssim\delta_q\Lambda_q^r
	\MM{N+1,N_{\mathrm{ind}}-r,\lambda_q,\ell_q^{-1}}.
	$$
	We then use
	$$
	\DTL^r\nabla J_u
	=\DTL^{r-1}\nabla(\DTL J_u)
	-\DTL^{r-1}\bigl((\nabla\overline u_q)^T\nabla J_u\bigr).
	$$
	By induction on $r$, we have
	$$
	\begin{aligned}
		\nrm{\DTL^r\nabla J_u}_N
		\leq\nrm{\DTL^{r-1}\nabla(\DTL J_u)}_N
		+\nrm{\DTL^{r-1}\bigl((\nabla\overline u_q)^T\nabla J_u\bigr)}_N\
		\lesssim\delta_q\Lambda_q^r(1+\tau_q\Lambda_q)
		\MM{N+1,N_{\mathrm{ind}}-r,\lambda_q,\ell_q^{-1}}.
	\end{aligned}
	$$
	The product term uses \eqref{est:uniform-grad-om-over-orho} and the induction hypothesis.
	Using $\delta_q^{\frac12}\lambda_q=\Lambda_q$ and $\tau_q\Lambda_q\leq1$, we obtain
	\begin{equation}\label{eq:Newton-Jacobian-gradient-CN}
	\begin{aligned}
		\nrm{\nabla J_u}_{N}&\lesssim\tau_q\delta_q\Lambda_q\lambda_q
		\MM{N,N_{\mathrm{ind}}-2,\lambda_q,\ell_q^{-1}},
		&&0\leq N\leq N_{\mathrm{max}}-2,\\
		\nrm{\nabla\DTL^rJ_u}_{N}+\nrm{\DTL^r\nabla J_u}_{N}&\lesssim\delta_q\Lambda_q^r
		\MM{N+1,N_{\mathrm{ind}}-r,\lambda_q,\ell_q^{-1}},
		&&1\leq r\leq N_{\mathrm{ind}},\quad r+N\leq N_{\mathrm{max}}-1.
	\end{aligned}
	\end{equation}
	Spatial interpolation at orders $N$ and $N+1$ gives
	\begin{equation}\label{eq:Newton-Jacobian-gradient-Holder}
	\begin{aligned}
		\nrm{\nabla J_u}_{N+\alpha}&\lesssim\tau_q\delta_q\Lambda_q\lambda_q\La{N,N_{\mathrm{ind}}-2,\lambda_q,\ell_q^{-1}}\MM{N,N_{\mathrm{ind}}-2,\lambda_q,\ell_q^{-1}},\\
		\nrm{\nabla\DTL^rJ_u}_{N+\alpha}+\nrm{\DTL^r\nabla J_u}_{N+\alpha}&\lesssim\delta_q\Lambda_q^r\La{N+1,N_{\mathrm{ind}}-r,\lambda_q,\ell_q^{-1}}\MM{N+1,N_{\mathrm{ind}}-r,\lambda_q,\ell_q^{-1}}.
	\end{aligned}
	\end{equation}
	The first estimate holds for $N\leq N_{\mathrm{max}}-3$. The second holds for $1\leq r\leq N_{\mathrm{ind}}$ and $r+N\leq N_{\mathrm{max}}-2$.
	
	We next estimate the time derivatives of $\overline V_{q,j+1,u,\kappa}$. By \eqref{def:Newton-perturbations-from-potential}, \eqref{def:Newton-mean-corrected-Piola-potential}, and compact support,
	$$
	\int_{\mathbb R^3}w[\widetilde V_{q,j+1,u,\kappa}]\,\rd x=-\partial_t\overline V_{q,j+1,u,\kappa}.
	$$ Integrating \eqref{eq:Newton-raw-Piola-momentum-equation} over $\mathbb R^3$ gives
	\begin{equation}\label{eq:Newton-mean-Piola-identity}
		\begin{aligned}
			-\partial_t^{2}(\overline V_{q,j+1,u,\kappa})_i
			=\int_{\mathbb R^3}\Bigl[
			(\mathcal E_{q,j,u,\kappa}^{\sharp})_i
			-\DTL\left((\Div\overline{u}_q)(\widetilde V_{q,j+1,u,\kappa})_i\right)-(\Div\overline{u}_q)\DTL(\widetilde V_{q,j+1,u,\kappa})_i
			-(\Div\overline{u}_q)^2(\widetilde V_{q,j+1,u,\kappa})_i\Bigr]\,\rd x.
		\end{aligned}
	\end{equation}
	The potential $\widetilde V_{q,j+1,u,\kappa}$ vanishes near $t_{u-1}$, so $
	\overline V_{q,j+1,u,\kappa}(t_{u-1})
	=\partial_t\overline V_{q,j+1,u,\kappa}(t_{u-1})=0.
	$
	
	We first estimate $\mathcal E_{q,j,u,\kappa}^{\sharp}$ using \eqref{def:Newton-sharp-vector-remainder}. For $0\leq M\leq N_{\mathrm{ind}}-2$ and $M+N\leq K_j^{(3)}$, the term with coefficient
	$J_u\partial_i(\Phi_u)_n-(\mathbf C_u)_{ni}$ satisfies
	$$
	\begin{aligned}
		&\Nrm{\DTL^M\Bigl[\bigl(J_u\partial_i(\Phi_u)_n-(\mathbf C_u)_{ni}\bigr)\partial_{y_n}\widehat I_{q,u,\kappa}[\DTL^2G_{q,j,u,\kappa}^{\mathrm{mol}}](t,\Phi_u)\Bigr]}_{N+\alpha}\\
		&\quad\lesssim\sum_{r=0}^M\sum_{N_1+N_2=N}\Bigl[
		\Nrm{\DTL^r\bigl(J_u\partial_i(\Phi_u)_n-(\mathbf C_u)_{ni}\bigr)}_{N_1}\Nrm{\DTL^{M-r}\Bigl(\partial_{y_n}\widehat I_{q,u,\kappa}[\DTL^2G_{q,j,u,\kappa}^{\mathrm{mol}}](t,\Phi_u)\Bigr)}_{N_2+\alpha}\\
		&\hspace{70pt}+\Nrm{\DTL^r\bigl(J_u\partial_i(\Phi_u)_n-(\mathbf C_u)_{ni}\bigr)}_{N_1+\alpha}\Nrm{\DTL^{M-r}\Bigl(\partial_{y_n}\widehat I_{q,u,\kappa}[\DTL^2G_{q,j,u,\kappa}^{\mathrm{mol}}](t,\Phi_u)\Bigr)}_{N_2+\alpha_1}\Bigr]\\
		&\quad\lesssim A_{q,j}\left[\tau_q^3\Lambda_q\MM{M+3,N_R-3j,\tau_q^{-1},\ell_{t,q}^{-1}}+\tau_q^2\sum_{r=1}^M\Lambda_q^r\MM{M-r+3,N_R-3j,\tau_q^{-1},\ell_{t,q}^{-1}}\right]\\
		&\qquad\cdot\Lat{N,N_{R,j}^{(3)}(M),\lambda_q,\ell_q^{-1}}\MM{N,N_{R,j}^{(3)}(M),\lambda_q,\ell_q^{-1}}\\
		&\quad\lesssim\Lambda_q A_{q,j}\MM{M,N_R-3j-3,\tau_q^{-1},\ell_{t,q}^{-1}}\Lat{N,N_{R,j}^{(3)}(M),\lambda_q,\ell_q^{-1}}\MM{N,N_{R,j}^{(3)}(M),\lambda_q,\ell_q^{-1}}.
	\end{aligned}
	$$
	For the second inequality, we used \eqref{eq:Newton-commuted-cone-gradient} and the gradient estimate in \eqref{eq:Newton-cone-gradient-derivative-estimate} at material order $M-r+2$.
	We also used \eqref{eq:Newton-frozen-flow-coefficient-estimates} for $r=0$ and \eqref{eq:Newton-cofactor-material-estimates} for $r\geq1$. The source estimates apply since
	$
	M+3\leq N_{\mathrm{ind}}+1\leq N_{\mathrm{max}},\ M+N+3\leq K_j^{(0)}.
	$
	
	For the terms containing $\DTL\mathbf C_u$ and $\DTL^2\mathbf C_u$ in \eqref{def:Newton-sharp-vector-remainder}, let $h=1,2$. In the same range of $M,N$, \eqref{eq:holder-product-lower-exponent} gives
	$$
	\begin{aligned}
		&\Nrm{\DTL^M\Bigl[(\DTL^h\mathbf C_u)^T
		\nabla_y\widehat I_{q,u,\kappa}[\DTL^{2-h}G_{q,j,u,\kappa}^{\mathrm{mol}}](t,\Phi_u)\Bigr]}_{N+\alpha}\\
		&\quad\lesssim\sum_{r=0}^M\sum_{N_1+N_2=N}\Bigl[
		\nrm{\DTL^{r+h}\mathbf C_u}_{N_1}
		\Nrm{\DTL^{M-r}\bigl(\nabla_y\widehat I_{q,u,\kappa}[\DTL^{2-h}G_{q,j,u,\kappa}^{\mathrm{mol}}](t,\Phi_u)\bigr)}_{N_2+\alpha}\\
		&\hspace{70pt}+
		\nrm{\DTL^{r+h}\mathbf C_u}_{N_1+\alpha}
		\Nrm{\DTL^{M-r}\bigl(\nabla_y\widehat I_{q,u,\kappa}[\DTL^{2-h}G_{q,j,u,\kappa}^{\mathrm{mol}}](t,\Phi_u)\bigr)}_{N_2+\alpha_1}\Bigr]\\
		&\quad\lesssim\tau_q^2\Lambda_q^h A_{q,j}
		\sum_{r=0}^M\Lambda_q^r\MM{M-r+3-h,N_R-3j,\tau_q^{-1},\ell_{t,q}^{-1}}\Lat{N,N_{R,j}^{(3)}(M),\lambda_q,\ell_q^{-1}}
		\MM{N,N_{R,j}^{(3)}(M),\lambda_q,\ell_q^{-1}}\\
		&\quad\lesssim\Lambda_q A_{q,j}
		\MM{M,N_R-3j-3,\tau_q^{-1},\ell_{t,q}^{-1}}
		\Lat{N,N_{R,j}^{(3)}(M),\lambda_q,\ell_q^{-1}}
		\MM{N,N_{R,j}^{(3)}(M),\lambda_q,\ell_q^{-1}}.
	\end{aligned}
	$$
	Here the second inequality uses \eqref{eq:Newton-cofactor-material-estimates} at material order $r+h$.
	We also used \eqref{eq:Newton-commuted-cone-gradient} and the gradient estimate in \eqref{eq:Newton-cone-gradient-derivative-estimate} at material order $M-r+2-h$.
	The required material orders satisfy $r+h\leq M+2\leq N_{\mathrm{ind}}$ and $M-r+3-h\leq M+2$.
	The last inequality also uses
	$$
	\Lambda_q^r\leq\tau_q^{-r},\qquad
	\tau_q^{h-1}\Lambda_q^h=\Lambda_q(\tau_q\Lambda_q)^{h-1}\leq\Lambda_q,\qquad h=1,2.
	$$
	For the term containing $\nabla J_u$, combine \eqref{eq:Newton-cone-material-identities} and \eqref{eq:Newton-cone-derivative-estimate} at material order $M-r+2$
	with \eqref{eq:Newton-Jacobian-gradient-CN}--\eqref{eq:Newton-Jacobian-gradient-Holder} to obtain, in the same range,
	$$
	\begin{aligned}
		&\Nrm{\DTL^M\Bigl(\nabla J_u\,I_{q,u,\kappa}[\DTL^2G_{q,j,u,\kappa}^{\mathrm{mol}}]\Bigr)}_{N+\alpha}\\
		&\quad\lesssim\sum_{r=0}^M\sum_{N_1+N_2=N}\Bigl[
		\nrm{\DTL^r\nabla J_u}_{N_1}
		\nrm{\DTL^{M-r}I_{q,u,\kappa}[\DTL^2G_{q,j,u,\kappa}^{\mathrm{mol}}]}_{N_2+\alpha}
		+\nrm{\DTL^r\nabla J_u}_{N_1+\alpha}
		\nrm{\DTL^{M-r}I_{q,u,\kappa}[\DTL^2G_{q,j,u,\kappa}^{\mathrm{mol}}]}_{N_2}\Bigr]\\
		&\quad\lesssim\tau_q\delta_q\lambda_q\Lambda_q A_{q,j}\MM{M,N_R-3j-3,\tau_q^{-1},\ell_{t,q}^{-1}}\Lat{N,N_{R,j}^{(3)}(M),\lambda_q,\ell_q^{-1}}\MM{N,N_{R,j}^{(3)}(M),\lambda_q,\ell_q^{-1}}.
	\end{aligned}
	$$
	Since $\tau_q\delta_q\lambda_q\leq1$ by \eqref{eq:Newton-scale-separation}, this completes the estimate for $\mathcal E_{q,j,u,\kappa}^{\sharp}$ in the stated range.

	For the remaining products in \eqref{eq:Newton-mean-Piola-identity}, use \eqref{est:uniform-div-om-over-orho} and \eqref{est:holder-div-om-over-orho} for the $\Div\overline u_q$ factors, and \eqref{eq:Newton-raw-Piola-estimate} for the $\widetilde V_{q,j+1,u,\kappa}$ factors. For $M\leq N_{\mathrm{ind}}-2$, the required divergence estimates hold through material order $M+1\leq N_{\mathrm{ind}}-1$.
	At $M=N_{\mathrm{ind}}-1$, \eqref{eq:Newton-mollified-source-estimate} applies through order $M+3=N_{\mathrm{ind}}+2\leq N_{\mathrm{max}}$, since $M+N+3\leq K_j^{(0)}$.
	The only coefficient factors requiring a separate estimate are
	$\DTL^{N_{\mathrm{ind}}+1}\mathbf C_u$ and $\DTL^{N_{\mathrm{ind}}}\Div\overline u_q$. They arise from the $\DTL^2\mathbf C_u$ term in $\mathcal E_{q,j,u,\kappa}^{\sharp}$ and the differentiated divergence factor, respectively.
	For these coefficients, use \eqref{eq:Newton-raw-cofactor-material-estimates} at $M=N_{\mathrm{ind}}+1$ and the estimates \eqref{est:uniform-grad-om-over-orho} and \eqref{est:holder-grad-om-over-orho}, respectively, with $\Div\overline u_q=\operatorname{tr}(\nabla\overline u_q)$. For $N_{\mathrm{ind}}+N+2\leq N_{\mathrm{max}}$, we have
	$$
	\nrm{\DTL^{N_{\mathrm{ind}}+1}\mathbf C_u}_{N+\alpha}+\nrm{\DTL^{N_{\mathrm{ind}}}\Div\overline u_q}_{N+\alpha}\lesssim\delta_q^{\frac12}\Lambda_q^{N_{\mathrm{ind}}}\ell_q^{-N-1-\alpha}.
	$$
	For the potential factors with no material derivatives, use \eqref{eq:Newton-raw-Piola-estimate} and the gradient estimate in \eqref{eq:Newton-cone-gradient-derivative-estimate}. Both have amplitude $\tau_q A_{q,j}$.
	By \eqref{eq:material-mollification-hierarchy}--\eqref{eq:material-mollification-spatial-comparisons} and \eqref{eq:Newton-scale-separation},
	$$
	\begin{aligned}
		\tau_q\delta_q^{\frac12}\Lambda_q^{N_{\mathrm{ind}}}\ell_q^{-1}
		&\leq\Lambda_q
		\MM{N_{\mathrm{ind}}-1,N_R-3j-3,\tau_q^{-1},\ell_{t,q}^{-1}}.
	\end{aligned}
	$$
	At $M=N_{\mathrm{ind}}-1$, $N_{R,j}^{(3)}(M)=0$. Apply
	\eqref{eq:holder-product-lower-exponent} using the $C^N$ and $C^{N+\alpha}$ estimates for these coefficients
	and the $C^{N+\alpha}$ and $C^{N+\alpha_1}$ estimates for the potential
	factors with no material derivatives.
	The estimates for $\mathcal E_{q,j,u,\kappa}^{\sharp}$ and the other terms in \eqref{eq:Newton-mean-Piola-identity} give, for $0\leq M\leq N_{\mathrm{ind}}-1$ and $M+N\leq K_j^{(3)}$,
	\begin{align}
		&\nrm{\DTL^M\mathcal E_{q,j,u,\kappa}^{\sharp}}_{N+\alpha}
		+\nrm{\DTL^{M+1}\bigl((\Div\overline u_q)\widetilde V_{q,j+1,u,\kappa}\bigr)}_{N+\alpha}+\nrm{\DTL^M\bigl((\Div\overline u_q)\DTL\widetilde V_{q,j+1,u,\kappa}\bigr)}_{N+\alpha}
		+\nrm{\DTL^M\bigl((\Div\overline u_q)^2\widetilde V_{q,j+1,u,\kappa}\bigr)}_{N+\alpha}\nonumber\\
		&\quad\lesssim\Lambda_q A_{q,j}\MM{M,N_R-3j-3,\tau_q^{-1},\ell_{t,q}^{-1}}\Lat{N,N_{R,j}^{(3)}(M),\lambda_q,\ell_q^{-1}}\MM{N,N_{R,j}^{(3)}(M),\lambda_q,\ell_q^{-1}}.\label{eq:Newton-mean-integrand-estimate}
	\end{align}
	
	For $M=0,1$, use \eqref{eq:Newton-mean-integrand-estimate} at material and spatial order zero with exponent $\alpha_1$. Integrate \eqref{eq:Newton-mean-Piola-identity} twice and once, respectively. The zero initial data and the interval length
	$O(\tau_q)$ give factors $\tau_q^2$ and $\tau_q$. For $M\geq2$, differentiate $M-2$ times in $t$ and repeatedly
	use the following identity for spatially compactly supported $h$.
	$$
	\partial_t\int_{\mathbb R^3}h\,\rd x=\int_{\mathbb R^3}(\DTL h+(\Div\overline u_q)h)\,\rd x.
	$$ 
	For the additional factors $\DTL^r\Div\overline u_q$, \eqref{est:uniform-div-om-over-orho} and \eqref{def:tau-q} give
	$$
	\nrm{\DTL^r\Div\overline u_q}_0
	\lesssim\delta_q\Lambda_q^{r+1}
	=\delta_q(\tau_q\Lambda_q)^{r+1}\tau_q^{-r-1}
	\leq\tau_q^{-r-1},
	\qquad 0\leq r\leq N_{\mathrm{ind}}-2.
	$$
	For $2\leq M\leq N_{\mathrm{ind}}+1$, apply \eqref{eq:Newton-mean-integrand-estimate} at $N=0$ and exponent $\alpha_1$, with material indices at most $M-2\leq N_{\mathrm{ind}}-1$.

	By \eqref{eq:Newton-cone-slice-support} and
	$J_u^{-1}\lesssim1$ from \eqref{eq:Newton-Jacobian-material-estimates} with $\sigma=-1$, the differentiated integrands in \eqref{eq:Newton-mean-Piola-identity} have spatial support of volume $O(\tau_q^3)$
	on $\mcI_u$. For
	$M\geq2$, $N_{R,j}^{(3)}(M-2)=N_{R,j}^{(1)}(M)$ and $
	\tau_q^3\MM{M-2,N_R-3j-3,\tau_q^{-1},\ell_{t,q}^{-1}}=\tau_q^5\MM{M,N_R-3j-1,\tau_q^{-1},\ell_{t,q}^{-1}}.
	$
	Thus, for $0\leq M\leq N_{\mathrm{ind}}+1$,
	\begin{equation}\label{eq:Newton-mean-time-derivative-estimate}
	\Nrm{\partial_t^{M}\overline V_{q,j+1,u,\kappa}}_0
	\lesssim\Lambda_qA_{q,j}\tau_q^{5}\MM{M,N_R-3j-1,\tau_q^{-1},\ell_{t,q}^{-1}}\La{0,N_{R,j}^{(1)}(M),\lambda_q,\ell_q^{-1}}.
	\end{equation}
	
	By \eqref{def:Newton-mean-profile}, we have
\begin{equation}\label{eq:Newton-mean-profile-gradient}
	\begin{aligned}
		\nabla\phi_{u,\kappa}(t,x)
		&=\tau_q^{-3}\nabla J_u(t,x)
		\phi\left(\tau_q^{-1}(\Phi_u-x_\kappa)\right)+\tau_q^{-4}J_u(t,x)(\nabla\Phi_u(t,x))^T
		(\nabla\phi)\left(\tau_q^{-1}(\Phi_u-x_\kappa)\right).
	\end{aligned}
\end{equation}
	Since $\DTL\Phi_u=0$, \eqref{def:Newton-mean-profile}, \eqref{eq:Newton-mean-profile-gradient}, and \eqref{eq:holder-product-finite-Cnalpha} give
	\begin{equation}\label{eq:Newton-mean-profile-Leibniz-estimates}
	\begin{aligned}
		\nrm{\DTL^r\phi_{u,\kappa}}_{N+\alpha}
		&\lesssim\tau_q^{-3}\sum_{N_1+N_2=N}\Bigl[
		\nrm{\DTL^rJ_u}_{N_1}\nrm{\phi\left(\tau_q^{-1}(\Phi_u-x_\kappa)\right)}_{N_2+\alpha}
		+\nrm{\DTL^rJ_u}_{N_1+\alpha}\nrm{\phi\left(\tau_q^{-1}(\Phi_u-x_\kappa)\right)}_{N_2}\Bigr],\\
		\nrm{\DTL^r\nabla\phi_{u,\kappa}}_{N+\alpha}
		&\lesssim\tau_q^{-3}\sum_{N_1+N_2=N}\Bigl[
		\nrm{\DTL^r\nabla J_u}_{N_1}\nrm{\phi\left(\tau_q^{-1}(\Phi_u-x_\kappa)\right)}_{N_2+\alpha}
		+\nrm{\DTL^r\nabla J_u}_{N_1+\alpha}\nrm{\phi\left(\tau_q^{-1}(\Phi_u-x_\kappa)\right)}_{N_2}\Bigr]\\
		&\qquad+\tau_q^{-4}\sum_{N_1+N_2=N}\Bigl[
		\nrm{\DTL^r(J_u\nabla\Phi_u)}_{N_1}\nrm{(\nabla\phi)\left(\tau_q^{-1}(\Phi_u-x_\kappa)\right)}_{N_2+\alpha}\\
		&\hspace{80pt}+\nrm{\DTL^r(J_u\nabla\Phi_u)}_{N_1+\alpha}\nrm{(\nabla\phi)\left(\tau_q^{-1}(\Phi_u-x_\kappa)\right)}_{N_2}\Bigr].
	\end{aligned}
	\end{equation}
	For the composed profiles, apply Lemmas~\ref{lem:BDSV-composition} and \ref{lem:holder-composition-general} with \eqref{eq:mollified-flow-close-identity} and $\tau_q^{-1}\leq\lambda_q$.
	Use \eqref{eq:Newton-Jacobian-material-estimates} with $\sigma=1$ for $J_u$ and \eqref{eq:Newton-Jacobian-gradient-CN}--\eqref{eq:Newton-Jacobian-gradient-Holder} for $\nabla J_u$.
	For $J_u\nabla\Phi_u$, combine the estimate for $J_u$ with \eqref{eq:mollified-flow-close-identity} at material order zero. At positive material orders, use \eqref{eq:Newton-cofactor-material-estimates}.
	The spatial estimates for $\phi_{u,\kappa}$ and $\nabla\phi_{u,\kappa}$ are
	\begin{equation}\label{eq:Newton-mean-profile-spatial-estimates}
	\begin{aligned}
		\nrm{\phi_{u,\kappa}}_{N+\alpha}
		&\lesssim\tau_q^{-3}
		\La{N,N_{R,j}^{(0)}(0),\lambda_q,\ell_q^{-1}}\MM{N,N_{R,j}^{(0)}(0),\lambda_q,\ell_q^{-1}},
		&&0\leq N\leq N_{\mathrm{max}}-2,\\
		\nrm{\nabla\phi_{u,\kappa}}_{N+\alpha}
		&\lesssim\tau_q^{-4}
		\La{N,N_{R,j}^{(0)}(0),\lambda_q,\ell_q^{-1}}\MM{N,N_{R,j}^{(0)}(0),\lambda_q,\ell_q^{-1}},
		&&0\leq N\leq N_{\mathrm{max}}-3.
	\end{aligned}
	\end{equation}
	At positive material orders, apply \eqref{eq:Newton-Jacobian-material-estimates} and \eqref{eq:Newton-cofactor-material-estimates}--\eqref{eq:Newton-Jacobian-gradient-Holder} in \eqref{eq:Newton-mean-profile-Leibniz-estimates}.
	These estimates prove the profile estimate in \eqref{eq:Newton-mean-profile-material-estimates} and the gradient estimate for $r<N_{\mathrm{ind}}$.
	The gradient estimate at $r=N_{\mathrm{ind}}$ is proved next.
	\begin{equation}\label{eq:Newton-mean-profile-material-estimates}
	\begin{aligned}
		&\nrm{\DTL^r\phi_{u,\kappa}}_{N+\alpha}\\
		&\quad\lesssim\tau_q^{-3-r}
		\La{N,N_{R,j}^{(0)}(r),\lambda_q,\ell_q^{-1}}\MM{N,N_{R,j}^{(0)}(r),\lambda_q,\ell_q^{-1}},
		&&1\leq r\leq N_{\mathrm{ind}},\quad r+N\leq N_{\mathrm{max}}-1,\\
		&\nrm{\DTL^r\nabla\phi_{u,\kappa}}_{N+\alpha}\\
		&\quad\lesssim\tau_q^{-4-r}
		\La{N,N_{R,j}^{(0)}(r),\lambda_q,\ell_q^{-1}}\MM{N,N_{R,j}^{(0)}(r),\lambda_q,\ell_q^{-1}},
		&&1\leq r\leq N_{\mathrm{ind}},\quad r+N\leq N_{\mathrm{max}}-2.
	\end{aligned}
	\end{equation}
	At $r=N_{\mathrm{ind}}$, use the second estimate in \eqref{eq:Newton-mean-profile-Leibniz-estimates}.
	Estimate $\DTL^r\nabla J_u$ by \eqref{eq:Newton-Jacobian-gradient-CN}--\eqref{eq:Newton-Jacobian-gradient-Holder}, and $\DTL^r(J_u\nabla\Phi_u)$ by \eqref{eq:Newton-cofactor-material-estimates}.
	Using $\tau_q\Lambda_q\leq1$, we have
	\begin{equation}\label{eq:Newton-mean-profile-gradient-endpoint}
	\begin{aligned}
		\nrm{\DTL^{N_{\mathrm{ind}}}\nabla\phi_{u,\kappa}}_{N+\alpha}
		&\lesssim\left(\tau_q^{-4-N_{\mathrm{ind}}}+\tau_q^{-3}\delta_q\Lambda_q^{N_{\mathrm{ind}}}\ell_q^{-1}\right)\ell_q^{-N-\alpha},
		&&N_{\mathrm{ind}}+N\leq N_{\mathrm{max}}-2.
	\end{aligned}
	\end{equation}
	Since $b<1+2\beta$, definitions~\eqref{def:lambda-delta}, \eqref{def:parameter-mollification}, and \eqref{def:tau-q} give, for sufficiently large $a$,
	$
	\delta_q\tau_q\ell_q^{-1}
	\lesssim\lambda_q^{-\beta+\frac{b-1}{4}}\leq\lambda_q^{-\frac\beta2}\leq1.
	$
	Together with $\tau_q\Lambda_q\leq1$, this proves the second estimate in \eqref{eq:Newton-mean-profile-material-estimates} at $r=N_{\mathrm{ind}}$.
	At $r=N_{\mathrm{ind}}+1$, apply \eqref{eq:Newton-Jacobian-material-estimates} with $\sigma=1$ to \eqref{def:Newton-mean-profile}, using $\DTL\Phi_u=0$ and $\delta_q^{\frac12}\ell_q^{-1}\leq\ell_{t,q}^{-1}$.
	\begin{equation}\label{eq:Newton-mean-profile-material-endpoint}
	\begin{aligned}
		\nrm{\DTL^{N_{\mathrm{ind}}+1}\phi_{u,\kappa}}_{N+\alpha}
		&\lesssim\tau_q^{-3}\MM{N_{\mathrm{ind}}+1,N_R-3j-1,\tau_q^{-1},\ell_{t,q}^{-1}}\ell_q^{-N-\alpha},
		&&N_{\mathrm{ind}}+1+N\leq N_{\mathrm{max}}-1.
	\end{aligned}
	\end{equation}
	
	For $0\leq M\leq N_{\mathrm{ind}}+1$, with $N\leq N_{\mathrm{max}}-2$ when $M=0$ and $M+N\leq N_{\mathrm{max}}-1$ when $M\geq1$, we obtain
	$$
	\begin{aligned}
		&\nrm{\DTL^M(\overline V_{q,j+1,u,\kappa}\phi_{u,\kappa})}_{N+\alpha}\\
		&\quad\lesssim\sum_{r=0}^M\binom Mr
		\Nrm{\partial_t^{r}\overline V_{q,j+1,u,\kappa}}_0
		\nrm{\DTL^{M-r}\phi_{u,\kappa}}_{N+\alpha}\\
		&\quad\lesssim\tau_q^2\Lambda_q A_{q,j}
		\Lat{N,N_{R,j}^{(1)}(M),\lambda_q,\ell_q^{-1}}
		\MM{N,N_{R,j}^{(1)}(M),\lambda_q,\ell_q^{-1}}\\
		&\qquad\cdot\sum_{r=0}^M\binom Mr
		\MM{r,N_R-3j-1,\tau_q^{-1},\ell_{t,q}^{-1}}
		\MM{M-r,N_R-3j-1,\tau_q^{-1},\ell_{t,q}^{-1}}\\
		&\quad\lesssim\tau_q^2\Lambda_q A_{q,j}\MM{M,N_R-3j-1,\tau_q^{-1},\ell_{t,q}^{-1}}\Lat{N,N_{R,j}^{(1)}(M),\lambda_q,\ell_q^{-1}}\MM{N,N_{R,j}^{(1)}(M),\lambda_q,\ell_q^{-1}}.
	\end{aligned}
	$$
	Here the second inequality uses \eqref{eq:Newton-mean-time-derivative-estimate} and the first estimates in \eqref{eq:Newton-mean-profile-spatial-estimates}--\eqref{eq:Newton-mean-profile-material-estimates}. When $M-r=N_{\mathrm{ind}}+1$, we also used \eqref{eq:Newton-mean-profile-material-endpoint}. At the common spatial threshold $N_{R,j}^{(1)}(M)$, each H{\"o}lder factor is at most $\La{N,N_{R,j}^{(1)}(M),\lambda_q,\ell_q^{-1}}$. Estimate~\eqref{eq:Newton-mean-potential-estimate} follows.
	
	For the density and momentum perturbations from the mean correction, let $0\leq M\leq N_{\mathrm{ind}}$, with $N\leq N_{\mathrm{max}}-3$ at $M=0$ and $M+N\leq N_{\mathrm{max}}-2$ at $M\geq1$.
	The first identity in \eqref{eq:Newton-mean-perturbation-identities} gives
	$$
	\begin{aligned}
		&\nrm{\DTL^M\varrho[\overline V_{q,j+1,u,\kappa}\phi_{u,\kappa}]}_{N+\alpha}\\
		&\quad\lesssim\sum_{r=0}^M\binom Mr
		\Nrm{\partial_t^{r}\overline V_{q,j+1,u,\kappa}}_0\nrm{\DTL^{M-r}\nabla\phi_{u,\kappa}}_{N+\alpha}\\
		&\quad\lesssim\tau_q\Lambda_q A_{q,j}
		\Lat{N,N_{R,j}^{(1)}(M),\lambda_q,\ell_q^{-1}}
		\MM{N,N_{R,j}^{(1)}(M),\lambda_q,\ell_q^{-1}}\\
		&\qquad\cdot\sum_{r=0}^M\binom Mr
		\MM{r,N_R-3j-1,\tau_q^{-1},\ell_{t,q}^{-1}}
		\MM{M-r,N_R-3j-1,\tau_q^{-1},\ell_{t,q}^{-1}}\\
		&\quad\lesssim\tau_q\Lambda_q A_{q,j}\MM{M,N_R-3j-1,\tau_q^{-1},\ell_{t,q}^{-1}}\Lat{N,N_{R,j}^{(1)}(M),\lambda_q,\ell_q^{-1}}\MM{N,N_{R,j}^{(1)}(M),\lambda_q,\ell_q^{-1}}.
	\end{aligned}
	$$
	For the second inequality, we used \eqref{eq:Newton-mean-time-derivative-estimate} and the second estimates in \eqref{eq:Newton-mean-profile-spatial-estimates}--\eqref{eq:Newton-mean-profile-material-estimates}. Thus \eqref{eq:Newton-mean-density-estimate} holds.
	
	For the first momentum term in \eqref{eq:Newton-mean-perturbation-identities}, use \eqref{eq:Newton-mean-time-derivative-estimate} and the first estimates in \eqref{eq:Newton-mean-profile-spatial-estimates}--\eqref{eq:Newton-mean-profile-material-estimates}. Then
	$$
	\begin{aligned}
		&\Nrm{\DTL^M\left(\left(\partial_t\overline V_{q,j+1,u,\kappa}\right)\phi_{u,\kappa}\right)}_{N+\alpha}\\
		&\quad\lesssim\sum_{r=0}^M\binom Mr
		\Nrm{\partial_t^{r+1}\overline V_{q,j+1,u,\kappa}}_0
		\nrm{\DTL^{M-r}\phi_{u,\kappa}}_{N+\alpha}\\
		&\quad\lesssim\tau_q^2\Lambda_q A_{q,j}
		\Lat{N,N_{R,j}^{(2)}(M),\lambda_q,\ell_q^{-1}}
		\MM{N,N_{R,j}^{(2)}(M),\lambda_q,\ell_q^{-1}}\\
		&\qquad\cdot\sum_{r=0}^M\binom Mr
		\MM{r+1,N_R-3j-1,\tau_q^{-1},\ell_{t,q}^{-1}}
		\MM{M-r,N_R-3j-1,\tau_q^{-1},\ell_{t,q}^{-1}}\\
		&\quad\lesssim\tau_q\Lambda_q A_{q,j}
		\MM{M,N_R-3j-2,\tau_q^{-1},\ell_{t,q}^{-1}}
		\Lat{N,N_{R,j}^{(2)}(M),\lambda_q,\ell_q^{-1}}
		\MM{N,N_{R,j}^{(2)}(M),\lambda_q,\ell_q^{-1}}.
	\end{aligned}
	$$
	For the velocity and gradient factors in \eqref{eq:Newton-mean-perturbation-identities}, use \eqref{est:uniform-zero-velocity}--\eqref{est:uniform-grad-om-over-orho} and \eqref{est:holder-om-over-orho}--\eqref{est:holder-grad-om-over-orho}.
	Combine these with \eqref{eq:Newton-mean-potential-estimate} for the product containing $\nabla\overline u_q$
	and \eqref{eq:Newton-mean-density-estimate} for the product containing $\overline u_q$. At $M=N_{\mathrm{ind}}$, the term with all material derivatives on $\nabla\overline u_q$ uses
	\eqref{est:holder-grad-om-over-orho}. Since $N_{\mathrm{ind}}>N_R-3j-2$, the ratio satisfies
	\begin{equation}\label{eq:Newton-mean-momentum-endpoint-coefficient}
	\begin{aligned}
		\frac{\delta_q^{\frac12}\Lambda_q^{N_{\mathrm{ind}}}\ell_q^{-1}}
		{\tau_q^{-1}\MM{N_{\mathrm{ind}},N_R-3j-2,\tau_q^{-1},\ell_{t,q}^{-1}}}=(\delta_q^{\frac12}\ell_{t,q}\ell_q^{-1})(\tau_q\Lambda_q)^{N_R-3j-1}
		(\ell_{t,q}\Lambda_q)^{N_{\mathrm{ind}}-N_R+3j+1}\leq1.
	\end{aligned}
	\end{equation}
	In the last inequality, we used \eqref{eq:material-mollification-spatial-comparisons}, $\tau_q\Lambda_q\leq1$, and $\ell_{t,q}\leq\tau_q$.
	We obtain \eqref{eq:Newton-mean-momentum-estimate}.
	
	By \eqref{def:Newton-mean-corrected-Piola-potential} and \eqref{eq:Newton-mean-potential-estimate}, \eqref{eq:Newton-raw-Piola-estimate} also holds with
	$V_{q,j+1,u,\kappa}^{\circ}$ in place of $\widetilde V_{q,j+1,u,\kappa}$, for $0\leq M\leq N_{\mathrm{ind}}+1$ and $M+N\leq K_j^{(1)}$. Use \eqref{eq:time-cutoff-derivative-estimate}. At each point, only a fixed number of local terms contribute to
	\eqref{def:Newton-total-compressible-perturbation}. Hence \eqref{eq:Newton-total-compressible-potential-estimate} holds.
	
	For $\trho_{q,j+1}\comp$, use the first identity in \eqref{eq:Newton-local-compressible-cutoff-formulas}, \eqref{eq:Newton-raw-density-estimate}, and \eqref{eq:Newton-mean-density-estimate}. The cutoff estimate \eqref{eq:time-cutoff-derivative-estimate} and finite overlap give \eqref{eq:Newton-total-compressible-density-estimate} for $0\leq M\leq N_{\mathrm{ind}}$ and $M+N\leq K_j^{(2)}$.
	For $\tm_{q,j+1}\comp$, estimate $w[V_{q,j+1,u,\kappa}^{\circ}]$ in \eqref{eq:Newton-local-compressible-cutoff-formulas} by \eqref{eq:Newton-raw-momentum-estimate} and \eqref{eq:Newton-mean-momentum-estimate}.
	For the term containing $\partial_t\theta_u$, use \eqref{eq:Newton-total-compressible-potential-estimate}. For $0\leq M\leq N_{\mathrm{ind}}$ and $M+N\leq K_j^{(2)}$, we obtain
	$$
	\begin{aligned}
		&\Nrm{\DTL^M\left((\partial_t\theta_u)V_{q,j+1,u,\kappa}^{\circ}\right)}_{N+\alpha}\\
		&\quad\lesssim\sum_{r=0}^M\tau_q^{-r-1}
		\nrm{\DTL^{M-r}V_{q,j+1,u,\kappa}^{\circ}}_{N+\alpha}\\
		&\quad\lesssim A_{q,j}\Lat{N,N_{R,j}^{(2)}(M),\lambda_q,\ell_q^{-1}}
		\MM{N,N_{R,j}^{(2)}(M),\lambda_q,\ell_q^{-1}}
		\sum_{r=0}^M\tau_q^{-r}
		\MM{M-r,N_R-3j-1,\tau_q^{-1},\ell_{t,q}^{-1}}\\
		&\quad\lesssim A_{q,j}\MM{M,N_R-3j-2,\tau_q^{-1},\ell_{t,q}^{-1}}
		\Lat{N,N_{R,j}^{(2)}(M),\lambda_q,\ell_q^{-1}}
		\MM{N,N_{R,j}^{(2)}(M),\lambda_q,\ell_q^{-1}}.
	\end{aligned}
	$$
	For the last inequality, we lowered the material threshold from $N_R-3j-1$ to $N_R-3j-2$. Summing with finite overlap proves \eqref{eq:Newton-total-compressible-momentum-estimate}.
	
	To estimate $H_{q,j+1,u,\kappa}^{\mathrm{int}}$, use \eqref{def:Newton-interior-vector-remainder}.
	The first four terms inside the brackets satisfy \eqref{eq:Newton-mean-integrand-estimate}.
	For $(\partial_t^2\overline V_{q,j+1,u,\kappa})\phi_{u,\kappa}$, let $0\leq M\leq N_{\mathrm{ind}}-1$ and $M+N\leq K_j^{(3)}$. Then
	$$
	\begin{aligned}
		&\Nrm{\DTL^M\left(\left(\partial_t^{2}\overline V_{q,j+1,u,\kappa}\right)\phi_{u,\kappa}\right)}_{N+\alpha}\\
		&\quad\lesssim\sum_{r=0}^M\binom Mr
		\Nrm{\partial_t^{r+2}\overline V_{q,j+1,u,\kappa}}_0
		\nrm{\DTL^{M-r}\phi_{u,\kappa}}_{N+\alpha}\\
		&\quad\lesssim\tau_q^2\Lambda_qA_{q,j}
		\Lat{N,N_{R,j}^{(3)}(M),\lambda_q,\ell_q^{-1}}
		\MM{N,N_{R,j}^{(3)}(M),\lambda_q,\ell_q^{-1}}\sum_{r=0}^M\binom Mr\tau_q^{-(M-r)}
		\MM{r+2,N_R-3j-1,\tau_q^{-1},\ell_{t,q}^{-1}}\\
		&\quad\lesssim\Lambda_qA_{q,j}\MM{M,N_R-3j-3,\tau_q^{-1},\ell_{t,q}^{-1}}\Lat{N,N_{R,j}^{(3)}(M),\lambda_q,\ell_q^{-1}}\MM{N,N_{R,j}^{(3)}(M),\lambda_q,\ell_q^{-1}}.
	\end{aligned}
	$$
	For the second inequality, we used \eqref{eq:Newton-mean-time-derivative-estimate} and the first estimates in \eqref{eq:Newton-mean-profile-spatial-estimates}--\eqref{eq:Newton-mean-profile-material-estimates}. The time derivative orders satisfy $r+2\leq M+2\leq N_{\mathrm{ind}}+1$.
	The spatial factors use $N_{R,j}^{(1)}(r+2)=N_{R,j}^{(3)}(r)\geq N_{R,j}^{(3)}(M)$.
	Combining these estimates with \eqref{eq:time-cutoff-derivative-estimate} proves
	\eqref{eq:Newton-interior-vector-source-estimate}.
	
	For $\mathcal S_{q,j+1,u,\kappa}$ in \eqref{def:Newton-explicit-local-acoustic-stress}, we first estimate the material derivatives of the pressure coefficient.
	Since $c_{q,u,\kappa}$ is constant, applying \eqref{eq:faa-di-bruno} along the flow of $\overline u_q$ gives, for $r\geq1$,
	$$
	\DTL^r\bigl(p'(\orho_q)-c_{q,u,\kappa}^2\bigr)
	=\sum_{k=1}^r\frac{p^{(k+1)}(\orho_q)}{k!}\sum_{\substack{r_1+\cdots+r_k=r\\r_1,\ldots,r_k\geq1}}\binom{r}{r_1,\ldots,r_k}\prod_{\nu=1}^k\DTL^{r_\nu}\orho_q.
	$$
	Using \eqref{est:uniform-orho} in this expansion, we obtain the amplitude
	$\delta_q\Lambda_q^r$. The differentiated density factors have spatial thresholds at least
	$N_{\mathrm{ind}}-r$.
	Interpolation between the $C^N$ and $C^{N+1}$ estimates gives
	\begin{equation}\label{eq:Newton-frozen-pressure-material-estimate}
		\begin{aligned}
			\nrm{\DTL^r\bigl(p'(\orho_q)-c_{q,u,\kappa}^2\bigr)}_{N+\alpha}
			&\lesssim\delta_q\Lambda_q^r\La{N,N_{\mathrm{ind}}-r,\lambda_q,\ell_q^{-1}}\MM{N,N_{\mathrm{ind}}-r,\lambda_q,\ell_q^{-1}},
		\end{aligned}
	\end{equation}
	for $1\leq r\leq N_{\mathrm{ind}}$ and $r+N\leq N_{\mathrm{max}}-1$. The corresponding $C^N$ estimate is
	$$
	\nrm{\DTL^r\bigl(p'(\orho_q)-c_{q,u,\kappa}^2\bigr)}_{N}
	\lesssim\delta_q\Lambda_q^r\MM{N,N_{\mathrm{ind}}-r,\lambda_q,\ell_q^{-1}},
	$$
	for $1\leq r\leq N_{\mathrm{ind}}$ and $r+N\leq N_{\mathrm{max}}$.
	
	Fix $0\leq M\leq N_{\mathrm{ind}}-1$ and $M+N\leq K_j^{(2)}$. Applying \eqref{eq:holder-product-lower-exponent} with $\eta=\alpha_1$ to the three terms in \eqref{def:Newton-explicit-local-acoustic-stress}, we have
	$$
	\begin{aligned}
		&\Nrm{\DTL^M\Bigl[\bigl(p'(\orho_q)-c_{q,u,\kappa}^2\bigr)\varrho[\widetilde V_{q,j+1,u,\kappa}]\Bigr]}_{N+\alpha}
		+\Nrm{\DTL^M\Bigl[p'(\orho_q)\varrho[\overline V_{q,j+1,u,\kappa}\phi_{u,\kappa}]\Bigr]}_{N+\alpha}+\Nrm{\DTL^M\Bigl[(\DTL\overline u_q)\otimes V_{q,j+1,u,\kappa}^{\circ}\Bigr]}_{N+\alpha}\\
		&\quad\lesssim\sum_{r=0}^M\sum_{N_1+N_2=N}\Bigl[
		\Nrm{\DTL^r\bigl(p'(\orho_q)-c_{q,u,\kappa}^2\bigr)}_{N_1}\Nrm{\DTL^{M-r}\varrho[\widetilde V_{q,j+1,u,\kappa}]}_{N_2+\alpha}+\Nrm{\DTL^r\bigl(p'(\orho_q)-c_{q,u,\kappa}^2\bigr)}_{N_1+\alpha}\Nrm{\DTL^{M-r}\varrho[\widetilde V_{q,j+1,u,\kappa}]}_{N_2+\alpha_1}\\
		&\hspace{70pt}+\Nrm{\DTL^rp'(\orho_q)}_{N_1}\Nrm{\DTL^{M-r}\varrho[\overline V_{q,j+1,u,\kappa}\phi_{u,\kappa}]}_{N_2+\alpha}+\Nrm{\DTL^rp'(\orho_q)}_{N_1+\alpha}\Nrm{\DTL^{M-r}\varrho[\overline V_{q,j+1,u,\kappa}\phi_{u,\kappa}]}_{N_2+\alpha_1}\\
		&\hspace{70pt}+\Nrm{\DTL^{r+1}\overline u_q}_{N_1}\Nrm{\DTL^{M-r}V_{q,j+1,u,\kappa}^{\circ}}_{N_2+\alpha}+\Nrm{\DTL^{r+1}\overline u_q}_{N_1+\alpha}\Nrm{\DTL^{M-r}V_{q,j+1,u,\kappa}^{\circ}}_{N_2+\alpha_1}\Bigr]\\
		&\quad\lesssim\tau_q\Lambda_q A_{q,j}\MM{M,N_R-3j-2,\tau_q^{-1},\ell_{t,q}^{-1}}\Lat{N,N_{R,j}^{(2)}(M),\lambda_q,\ell_q^{-1}}\MM{N,N_{R,j}^{(2)}(M),\lambda_q,\ell_q^{-1}}.
	\end{aligned}
	$$
	For the first product, we combined \eqref{eq:Newton-raw-density-estimate}, \eqref{eq:Newton-frozen-flow-coefficient-estimates}, and \eqref{eq:Newton-frozen-pressure-material-estimate}.
	Here we used $\delta_q\Lambda_q^r\leq\tau_q\Lambda_q\tau_q^{-r}$ for $r\geq1$.
	For the second product, we used \eqref{eq:Newton-frozen-pressure-material-estimate} at positive material orders and \eqref{eq:faa-di-bruno} in space with \eqref{est:uniform-orho}.
	We combined them with \eqref{eq:Newton-mean-density-estimate}, which gives the factor $\tau_q\Lambda_q$.
	For the third product, we combined \eqref{est:uniform-om-over-orho} and \eqref{est:holder-om-over-orho} at material order $r+1$ with \eqref{eq:Newton-total-compressible-potential-estimate}. The potential contributes $\tau_q$, giving the factor $\tau_q\Lambda_q$.
	The coefficient estimates apply since $M+1\leq N_{\mathrm{ind}}$ and $M+N+1\leq K_j^{(2)}+1\leq N_{\mathrm{max}}-1$.
	The same argument applies in $C^{N+\alpha_r}$ using \eqref{eq:holder-product-lower-exponent}.
	Lower the spatial thresholds to $N_{R,j}^{(2)}(M)$ and use spatial order $N$ in the H{\"o}lder factors.
	These factors then combine by the following identity for $N,K\in\mathbb N_0$.
	$$
	\Lab{N,K,\lambda_q,\ell_q^{-1}}{\alpha_r}
	\Lkb{N,K,\lambda_q,\ell_q^{-1}}{2}{\alpha_{r+1}}
	=\Lkb{N,K,\lambda_q,\ell_q^{-1}}{2}{\alpha_r},
	\ 2\alpha_{r+1}=\alpha_r.
	$$
	Using this identity at exponent $\alpha$ and \eqref{eq:time-cutoff-derivative-estimate} proves \eqref{eq:Newton-interior-acoustic-residual-estimate}.
	
	We next estimate the cutoff residual in \eqref{def:Newton-compressible-cutoff-residual}. Fix $0\leq M\leq N_{\mathrm{ind}}$ and $M+N\leq K_j^{(2)}$.
	For $\DTL V_{q,j+1,u,\kappa}^{\circ}$, use \eqref{eq:Newton-total-compressible-potential-estimate} with $M$ replaced by $M+1$.
	The identities $\tau_q\MM{M+1,N_R-3j-1,\tau_q^{-1},\ell_{t,q}^{-1}}=\MM{M,N_R-3j-2,\tau_q^{-1},\ell_{t,q}^{-1}}$ and $N_{R,j}^{(1)}(M+1)=N_{R,j}^{(2)}(M)$ give the estimate for $\DTL V_{q,j+1,u,\kappa}^{\circ}$.
	Estimate $\nabla\overline u_q$ and $\Div\overline u_q$ using \eqref{est:uniform-grad-om-over-orho} and \eqref{est:holder-grad-om-over-orho}.
	At $M=N_{\mathrm{ind}}$, also use \eqref{eq:Newton-mean-momentum-endpoint-coefficient}. Estimate~\eqref{eq:holder-product-lower-exponent} gives
	$$
	\begin{aligned}
		&\nrm{\DTL^M\bigl(\DTL V_{q,j+1,u,\kappa}^{\circ}+(\Div\overline u_q)V_{q,j+1,u,\kappa}^{\circ}\bigr)}_{N+\alpha}
		+\nrm{\DTL^M\bigl((\nabla\overline u_q)V_{q,j+1,u,\kappa}^{\circ}\bigr)}_{N+\alpha}
		+\tau_q^{-1}\nrm{\DTL^MV_{q,j+1,u,\kappa}^{\circ}}_{N+\alpha}\\
		&\qquad\lesssim\nrm{\DTL^{M+1}V_{q,j+1,u,\kappa}^{\circ}}_{N+\alpha}
		+\tau_q^{-1}\nrm{\DTL^MV_{q,j+1,u,\kappa}^{\circ}}_{N+\alpha}+\sum_{r=0}^M\sum_{N_1+N_2=N}\Bigl[
		\bigl(\nrm{\DTL^r\Div\overline u_q}_{N_1}
		+\nrm{\DTL^r\nabla\overline u_q}_{N_1}\bigr)
		\nrm{\DTL^{M-r}V_{q,j+1,u,\kappa}^{\circ}}_{N_2+\alpha}\\
		&\hspace{230pt}+
		\bigl(\nrm{\DTL^r\Div\overline u_q}_{N_1+\alpha}
		+\nrm{\DTL^r\nabla\overline u_q}_{N_1+\alpha}\bigr)
		\nrm{\DTL^{M-r}V_{q,j+1,u,\kappa}^{\circ}}_{N_2+\alpha_1}\Bigr]\\
		&\qquad\lesssim A_{q,j}\MM{M,N_R-3j-2,\tau_q^{-1},\ell_{t,q}^{-1}}\Lat{N,N_{R,j}^{(2)}(M),\lambda_q,\ell_q^{-1}}\MM{N,N_{R,j}^{(2)}(M),\lambda_q,\ell_q^{-1}}.
	\end{aligned}
	$$
	By \eqref{def:Newton-compressible-cutoff-residual},
	$$
	\begin{aligned}
		\nrm{\DTL^MF_{q,j+1,u,\kappa}^{\mathrm{cut}}}_{N+\alpha}
		&\lesssim\sum_{r=0}^M\tau_q^{-r-1}\Bigl(\nrm{\DTL^{M-r}\bigl(\DTL V_{q,j+1,u,\kappa}^{\circ}+(\Div\overline u_q)V_{q,j+1,u,\kappa}^{\circ}\bigr)}_{N+\alpha}
		+\tau_q^{-1}\nrm{\DTL^{M-r}V_{q,j+1,u,\kappa}^{\circ}}_{N+\alpha}\Bigr)\\
		&\lesssim A_{q,j}\tau_q^{-1}\MM{M,N_R-3j-2,\tau_q^{-1},\ell_{t,q}^{-1}}\Lat{N,N_{R,j}^{(2)}(M),\lambda_q,\ell_q^{-1}}\MM{N,N_{R,j}^{(2)}(M),\lambda_q,\ell_q^{-1}}.
	\end{aligned}
	$$
	Thus \eqref{eq:Newton-compressible-cutoff-residual-estimate} holds. The temporal support of $F_{q,j+1,u,\kappa}^{\mathrm{cut}}$ follows from \eqref{def:Newton-compressible-cutoff-residual}.
	
	Apply Theorem~\ref{thm:inv4divEq} with $\upsilon=\overline u_q$ and $s=2$ to $H_{q,j+1,u,\kappa}^{\mathrm{int}}$ and $F_{q,j+1,u,\kappa}^{\mathrm{cut}}$, using \eqref{eq:Newton-interior-vector-source-estimate} and \eqref{eq:Newton-compressible-cutoff-residual-estimate}, respectively.
	The required support and zero-mean hypotheses were verified before \eqref{def:Newton-local-compressible-residual-matrix}.
	Take
	$$
	K^*=N_{\mathrm{ind}}-1,\quad \mC_\upsilon=C\delta_q^{\frac12},\quad \lambda_\upsilon=\lambda_U=\lambda_q,\quad \tilde\lambda_\upsilon=\tilde\lambda_U=\ell_q^{-1},\quad
		\mu_\upsilon=\tilde\mu_\upsilon=\Lambda_q,\quad \mu_U=\tau_q^{-1},\quad \tilde\mu_U=\ell_{t,q}^{-1}.
	$$
	For the interior and cutoff sources, respectively, take
	$$
	\begin{aligned}
		K_U=N_t&=N_R-3j-3,
		&M_*&=N_{\mathrm{ind}}-1,
		&K^\sharp&=K_j^{(3)}+1,
		&\mC_U&=C\Lambda_qA_{q,j},\\
		K_U=N_t&=N_R-3j-2,
		&M_*&=N_{\mathrm{ind}},
		&K^\sharp&=K_j^{(2)}+1,
		&\mC_U&=CA_{q,j}\tau_q^{-1}.
	\end{aligned}
	$$
	Here $C$ covers the implicit constants in \eqref{est:uniform-grad-om-over-orho}, \eqref{est:holder-grad-om-over-orho}, and \eqref{eq:Newton-interior-vector-source-estimate}--\eqref{eq:Newton-compressible-cutoff-residual-estimate}, uniformly over the required derivative orders and fixed H{\"o}lder exponents.
	These choices satisfy $K_U=N_t\leq K^*$, $M_*\leq K^*+1$, and $K^\sharp\geq K^*+2$.
	At each fixed exponent $\alpha_r$, the source estimates verify the theorem's H{\"o}lder hypothesis with $s=2$.
	The velocity hypotheses follow from \eqref{est:uniform-grad-om-over-orho} and \eqref{est:holder-grad-om-over-orho}. The short-time condition \eqref{eq:tau-local-smallness} follows from \eqref{eq:Newton-scale-separation}.
	The theorem gives the $C^{N+\alpha_r}$ estimates for $1\leq M\leq M_*$ and $M+N\leq K^\sharp-1$.
	At $M=0$, the range is $0\leq N\leq K^\sharp-2$.
	
	At $M=0$ and $N=K^\sharp-1$, use Remark~\ref{rem:local-inverse-holder-endpoint}. Its additional velocity estimate follows from \eqref{est:holder-grad-om-over-orho}, since $K^\sharp-1\leq K_j^{(2)}\leq N_{\mathrm{max}}-4$.
	For sufficiently large $a$, \eqref{def:tau-q} and \eqref{eq:material-mollification-spatial-comparisons} give
	\begin{equation}\label{eq:Newton-acoustic-combined-material-scales}
	\begin{aligned}
		\max\{\mu_U,\mu_\upsilon,\mC_\upsilon\lambda_\upsilon\}=\tau_q^{-1},\quad
		\max\{\tilde\mu_U,\tilde\mu_\upsilon,\mC_\upsilon\tilde\lambda_\upsilon\}=\ell_{t,q}^{-1}.
	\end{aligned}
	\end{equation}
	The spatial scales are $\lambda_q$ and $\ell_q^{-1}$, since $(C_{\mathrm{supp}}\tau_q)^{-1}\leq\tau_q^{-1}\leq\lambda_q$.
	For $N,k\in\mathbb N_0$, the localization factor satisfies
	$$
	(C_{\mathrm{supp}}\tau_q)^{1+\alpha_r}\Lab{N+1,k,\lambda_q,\ell_q^{-1}}{\alpha_r}
	\lesssim\tau_q(\tau_q\ell_q^{-1})^{\alpha_r}
	\leq\tau_q(\tau_q\ell_q^{-1})^\alpha.
	$$
	For the $C^{N+\alpha_r}$ estimates, use the factor $\tau_q(\tau_q\ell_q^{-1})^{\alpha_r}$.
	The material weights satisfy, for $0\leq m\leq M$ and $k\in\mathbb N_0$,
	$$
	\begin{aligned}
	\Lambda_q^{M-m}\MM{m,k,\tau_q^{-1},\ell_{t,q}^{-1}}
	=(\tau_q\Lambda_q)^{M-m}\tau_q^{-(M-m)}\MM{m,k,\tau_q^{-1},\ell_{t,q}^{-1}}
	\leq(\tau_q\Lambda_q)^{M-m}\MM{M,k,\tau_q^{-1},\ell_{t,q}^{-1}}.
	\end{aligned}
	$$
	
	By \eqref{eq:Newton-interior-vector-source-estimate}, the material threshold is $N_R-3j-3$, and the source and flow terms have spatial thresholds at least $N_{R,j}^{(3)}(M)$. At the common spatial threshold $N_{R,j+1}^{(0)}(M)=N_{R,j}^{(3)}(M)$, the local inverse satisfies
	\begin{equation}\label{eq:Newton-interior-local-inverse-estimate}
	\begin{aligned}
		&\nrm{\DTL^MR^{(\overline u_q)}_{u,\kappa}[H_{q,j+1,u,\kappa}^{\mathrm{int}}]}_{N+\alpha}\\
		&\quad\lesssim\tau_q\Lambda_q(\tau_q\ell_q^{-1})^\alpha A_{q,j}\MM{M,N_R-3j-3,\tau_q^{-1},\ell_{t,q}^{-1}}\Lat{N,N_{R,j+1}^{(0)}(M),\lambda_q,\ell_q^{-1}}\MM{N,N_{R,j+1}^{(0)}(M),\lambda_q,\ell_q^{-1}}.
	\end{aligned}
	\end{equation}
	For $0\leq M\leq N_{\mathrm{ind}}-1$ and $M+N\leq K_j^{(3)}$, lower the material and spatial thresholds in \eqref{eq:Newton-interior-acoustic-residual-estimate} to $N_R-3j-3$ and $N_{R,j+1}^{(0)}(M)=N_{R,j}^{(3)}(M)$, respectively. We obtain
	$$
	\begin{aligned}
		&\nrm{\DTL^M\mathcal S_{q,j+1,u,\kappa}}_{N+\alpha}\\
		&\quad\lesssim\tau_q\Lambda_q A_{q,j}\MM{M,N_R-3j-3,\tau_q^{-1},\ell_{t,q}^{-1}}\Lat{N,N_{R,j+1}^{(0)}(M),\lambda_q,\ell_q^{-1}}
		\MM{N,N_{R,j+1}^{(0)}(M),\lambda_q,\ell_q^{-1}}.
	\end{aligned}
	$$
	Combining this estimate with \eqref{eq:Newton-interior-local-inverse-estimate} and using \eqref{def:Newton-local-compressible-residual-matrix} and finite overlap proves \eqref{eq:Newton-interior-stress-potential-estimate}.
	
	For $F_{q,j+1,u,\kappa}^{\mathrm{cut}}$, fix $0\leq M\leq N_{\mathrm{ind}}$ and $M+N\leq K_j^{(2)}$. The local inverse contributes $\tau_q$, which cancels the $\tau_q^{-1}$ from the cutoff source. At spatial threshold $N_{R,j}^{(2)}(M)$, we have
	\begin{equation}\label{eq:Newton-cutoff-local-inverse-estimate}
	\begin{aligned}
		&\nrm{\DTL^MR^{(\overline u_q)}_{u,\kappa}[F_{q,j+1,u,\kappa}^{\mathrm{cut}}]}_{N+\alpha}\\
		&\quad\lesssim A_{q,j}(\tau_q\ell_q^{-1})^\alpha\MM{M,N_R-3j-2,\tau_q^{-1},\ell_{t,q}^{-1}}\Lat{N,N_{R,j}^{(2)}(M),\lambda_q,\ell_q^{-1}}\MM{N,N_{R,j}^{(2)}(M),\lambda_q,\ell_q^{-1}}.
	\end{aligned}
	\end{equation}
	Summing \eqref{eq:Newton-cutoff-local-inverse-estimate} over $\kappa$ in \eqref{def:Newton-cutoff-Piola-stress} and using finite spatial overlap, we obtain
	\begin{equation}\label{eq:Newton-cutoff-stress-estimate}
	\begin{aligned}
		&\sup_u\Bigl(\nrm{\DTL^MS_{q,j+1,u}^{\mathrm{cut}}}_{N+\alpha}
		+\nrm{\DTL^MS_{q,j+1,u}^{\mathrm{cut,sym}}}_{N+\alpha}
		+\nrm{\DTL^MS_{q,j+1,u}^{\mathrm{cut,skew}}}_{N+\alpha}\Bigr)\\
		&\quad\lesssim A_{q,j}(\tau_q\ell_q^{-1})^\alpha\MM{M,N_R-3j-2,\tau_q^{-1},\ell_{t,q}^{-1}}\Lat{N,N_{R,j}^{(2)}(M),\lambda_q,\ell_q^{-1}}\MM{N,N_{R,j}^{(2)}(M),\lambda_q,\ell_q^{-1}}.
	\end{aligned}
	\end{equation}
	
	To estimate $W[S_{q,j+1,u}^{\mathrm{cut,skew}}]$, apply Proposition~\ref{prop:periodic-Hodge-material-estimates} with $\upsilon=\overline u_q$, $T=S_{q,j+1,u}^{\mathrm{cut,skew}}$, and the parameters
	$$
	\begin{aligned}
		&K^*=N_{\mathrm{ind}}-1,&& K_T=N_t=N_R-3j-2,&& K^\sharp=K_j^{(2)}+1,&& M_*=r^*=N_{\mathrm{ind}},&& \mu_T=\tau_q^{-1},&& \tilde\mu_T=\ell_{t,q}^{-1},\\
		&A_\upsilon=C\delta_q^{\frac12},&& A_T=C A_{q,j}(\tau_q\ell_q^{-1})^\alpha,&& \lambda_\upsilon=\lambda_T=\lambda_q,&& \tilde\lambda_\upsilon=\tilde\lambda_T=\ell_q^{-1},&& \mu_\upsilon=\tilde\mu_\upsilon=\Lambda_q,&& s=2.
	\end{aligned}
	$$
	Choose $C$ to cover the implicit constants in \eqref{est:uniform-grad-om-over-orho}, \eqref{est:holder-grad-om-over-orho}, and \eqref{eq:Newton-cutoff-stress-estimate} at the required derivative orders and fixed H{\"o}lder exponents.
	The source estimates at the smaller fixed H{\"o}lder exponents also satisfy
	\eqref{eq:periodic-Hodge-T-Holder-assumption} with the same $A_T$.
	The coefficient hypotheses follow from \eqref{est:uniform-grad-om-over-orho} and \eqref{est:holder-grad-om-over-orho}.
	Here $K_T=N_t\leq K^*$, $M_*=r^*=K^*+1$, and $K^\sharp\geq K^*+1$.
	By \eqref{eq:Newton-acoustic-combined-material-scales}, $\mu=\tau_q^{-1}$ and $\tilde\mu=\ell_{t,q}^{-1}$.
	For $W[S_{q,j+1,u}^{\mathrm{cut,skew}}]$, use \eqref{eq:Hodge-high-material-estimate} when $M\geq1$ and \eqref{eq:Hodge-material-zero-estimate} when $M=0$.
	Together with \eqref{eq:Newton-cutoff-stress-estimate}, these estimates give \eqref{eq:Newton-cutoff-stress-potential-estimate}.
	For the $C^{N+\alpha_r}$ estimate of $\DTL^MW[S_{q,j+1,u}^{\mathrm{cut,skew}}]$, apply Proposition~\ref{prop:periodic-Hodge-material-estimates} with $\alpha$ replaced by $\alpha_r$ and $A_T=C A_{q,j}(\tau_q\ell_q^{-1})^{\alpha_r}$.
\end{proof}

\subsection{The incompressible perturbation}
\label{subsec:Newton-incompressible-perturbation}

The incompressible correction $\tm_{q,j+1}\inc$ cancels $\curl W_{q,j}^{\mathrm{mol}}$ and the sum of $\Div S_{q,j+1,u}^{\mathrm{cut,skew}}$ over $u$. It changes only the momentum, with $\Div\tm_{q,j+1}\inc=0$.

For every $u\in\mathbb Z$, let $(\omega_{q,j,u},\pi_{q,j,u})$ solve the system below on $[t_u-\tau_q,t_u+\tau_q]$. This interval contains $\supp\tilde\theta_u$.
\begin{equation}\label{eq:Newton-local-incompressible-linearized-system}
	\left\{
	\begin{aligned}
		&\partial_t\omega_{q,j,u}
		+\Div\left(\overline{u}_q\otimes\omega_{q,j,u}
		+\omega_{q,j,u}\otimes\overline{u}_q\right)
		+\nabla\pi_{q,j,u}=-\theta_u^2\curl W_{q,j}^{\mathrm{mol}}-\curl W\left[S_{q,j+1,u}^{\mathrm{cut,skew}}\right],\\
		&\Div\omega_{q,j,u}=0,\\
		&\omega_{q,j,u}|_{t=t_u}=0,\qquad
		\fint_{\T^3}\pi_{q,j,u}\,\rd x=0.
	\end{aligned}
	\right.
\end{equation}
These are the unique local solutions given by Proposition~\ref{prop:linearized-equation}, with zero-mean pressure.
Their zero spatial mean and the identity $\Div\omega_{q,j,u}=0$ give
$$
\Psi_{q,j,u}:=\mathcal B\omega_{q,j,u},\qquad
\curl\Psi_{q,j,u}=\omega_{q,j,u}.
$$

Define the glued momentum perturbation and the corresponding pressure by
\begin{equation}\label{def:Newton-glued-incompressible-perturbation}
	\tm_{q,j+1}\inc:=\sum_u\tilde\theta_u\omega_{q,j,u},\qquad
	\pi_{q,j+1}\inc:=\sum_u\tilde\theta_u\pi_{q,j,u}.
\end{equation}
The summands extend smoothly by zero outside the local time intervals because $\tilde\theta_u$ vanishes near their endpoints.
Multiplying \eqref{eq:Newton-local-incompressible-linearized-system} by $\tilde\theta_u$, summing in $u$, and using
$\tilde\theta_u\theta_u^2=\theta_u^2$ and $\sum_u\theta_u^2=1$, we obtain
\begin{equation}\label{eq:Newton-glued-incompressible-equation}
	\begin{aligned}
		&\partial_t\tm_{q,j+1}\inc
		+\Div\left(\overline{u}_q\otimes\tm_{q,j+1}\inc
		+\tm_{q,j+1}\inc\otimes\overline{u}_q\right)
		+\nabla\pi_{q,j+1}\inc\\
		&\quad=-\curl W_{q,j}^{\mathrm{mol}}-\sum_u\curl W\left[
		S_{q,j+1,u}^{\mathrm{cut,skew}}
		\right]+\curl\left(\sum_u(\partial_t\tilde\theta_u)
		\Psi_{q,j,u}\right).
	\end{aligned}
\end{equation}
Here $(\partial_t\tilde\theta_u)\omega_{q,j,u}
=\curl((\partial_t\tilde\theta_u)\Psi_{q,j,u})$. Since $W$ acts only in space, $W[S_{q,j+1,u}^{\mathrm{cut,skew}}]$ has temporal support contained in that of
$S_{q,j+1,u}^{\mathrm{cut,skew}}$, and $\tilde\theta_u=1$ on this support.

Writing the pressure term from \eqref{eq:Newton-glued-incompressible-equation} as $-\nabla\pi_{q,j+1}\inc=\Div(-\pi_{q,j+1}\inc\Id)$, we define the interior stress by
\begin{equation}\label{def:Newton-interior-stress}
	R_{q,j+1}
	=R_{q,j+1}\comp-\pi_{q,j+1}\inc\Id.
\end{equation}
Define the gluing stress by
\begin{equation}\label{def:Newton-gluing-stress}
	\begin{aligned}
		R_{q,j+1}^{\mathrm{gluing}}:=&
		\sum_uS_{q,j+1,u}^{\mathrm{cut,sym}}+\mathcal R\curl\left(
		\sum_u(\partial_t\tilde\theta_u)\Psi_{q,j,u}\right).
	\end{aligned}
\end{equation}
Both terms are symmetric, and \eqref{eq:direct-symmetric-curl-properties} gives
$$
\Div R_{q,j+1}^{\mathrm{gluing}}
=\sum_u\Div S_{q,j+1,u}^{\mathrm{cut,sym}}
+\curl\left(\sum_u(\partial_t\tilde\theta_u)\Psi_{q,j,u}\right).
$$
We set
$$
S_{q,j+1}:=S_{q,j}+R_{q,j+1}^{\mathrm{gluing}}.
$$

Update the density and momentum by
\begin{equation}\label{def:Newton-final-state}
	\rho_{q,j+1}:=\rho_{q,j}+\trho_{q,j+1}\comp,\qquad
	m_{q,j+1}:=m_{q,j}+\tm_{q,j+1}\comp+\tm_{q,j+1}\inc.
\end{equation}
Since $\theta_u$ depends only on time, the Hodge decomposition and
\eqref{def:Newton-substep-flow-mollification} give
$$
\begin{aligned}
	-\sum_{(u,\kappa)\in\mathscr I}\theta_u\nabla(J_uG_{q,j,u,\kappa}^{\mathrm{mol}})-\curl W_{q,j}^{\mathrm{mol}}
	&=-\Div R_{q,j}+\Div P_{q+1,j+1}^{\mathrm{mol,comp}}+\Div P_{q+1,j+1}^{\mathrm{mol,inc}}.
\end{aligned}
$$
The correction equations \eqref{eq:Newton-approximate-compressible-linearized-system} and \eqref{eq:Newton-glued-incompressible-equation} are linearized at
$(\orho_q,\om_q)$. The flux remainder is the full flux change from $(\rho_{q,j},m_{q,j})$ to $(\rho_{q,j+1},m_{q,j+1})$ minus the linear flux terms in these equations.
Proposition~\ref{prop:Newton-updated-state-estimates} gives $\rho_{q,j},\rho_{q,j+1}\in[\rho^*/2,3\rho^*/2]$. For $0\leq j<J_*$, define
\begin{equation}\label{def:Newton-perturbative-tensor}
	\begin{aligned}
		P_{q+1,j+1}:=&P_{q+1,j}+P_{q+1,j+1}^{\mathrm{mol,comp}}+P_{q+1,j+1}^{\mathrm{mol,inc}}+\left[p(\rho_{q,j+1})-p(\rho_{q,j})
		-p'(\orho_q)\trho_{q,j+1}\comp\right]\Id
		\\
		&+\frac{m_{q,j+1}\otimes m_{q,j+1}}{\rho_{q,j+1}}
		-\frac{m_{q,j}\otimes m_{q,j}}{\rho_{q,j}}\\
		&-\Bigl[\overline{u}_q\otimes
		\left(\tm_{q,j+1}\comp+\tm_{q,j+1}\inc\right)
		+\left(\tm_{q,j+1}\comp+\tm_{q,j+1}\inc\right)
		\otimes\overline{u}_q
		-\trho_{q,j+1}\comp\overline{u}_q\otimes\overline{u}_q\Bigr].
	\end{aligned}
\end{equation}
The flux remainder and the two source mollification stresses in \eqref{def:Newton-perturbative-tensor} are symmetric. Adding \eqref{eq:Newton-approximate-compressible-linearized-system} and \eqref{eq:Newton-glued-incompressible-equation} cancels the divergences of the antisymmetric cutoff matrices. The remaining linear residuals are the divergences of the stresses defined in
\eqref{def:Newton-interior-stress}--\eqref{def:Newton-gluing-stress}. Together with \eqref{eq:Newton-substep-Euler-Reynolds-system} and
\eqref{def:Newton-perturbative-tensor}, we find
\begin{equation}\label{eq:Newton-final-Euler-Reynolds-system}
	\left\{
	\begin{aligned}
		&\partial_t\rho_{q,j+1}+\Div m_{q,j+1}=0,\\
		&\partial_t m_{q,j+1}
		+\Div\left(\frac{m_{q,j+1}\otimes m_{q,j+1}}{\rho_{q,j+1}}\right)
		+\nabla p(\rho_{q,j+1})=\Div R_{q,j+1}
		+\Div S_{q,j+1}
		+\Div P_{q+1,j+1}.
	\end{aligned}
	\right.
\end{equation}

\begin{prop}
	\label{prop:Newton-incompressible-flow-gluing-estimates}
	Fix $0\leq j<J_*$ and assume the source estimates \eqref{eq:Newton-mollified-source-estimate} at the fixed H{\"o}lder exponents. For $a>a_{\mathrm{Newton}}$, the following estimates hold.
	The time supremum in local norms is taken over $[t_u-\tau_q,t_u+\tau_q]$, while for the glued fields it is taken over all times. For $M,N\in\mathbb N_0$ with $M\leq N_{\mathrm{ind}}+1$ and $M+N\leq K_j^{(2)}$, we have
	\begin{equation}\label{eq:Newton-local-vector-potential-estimate}
		\begin{aligned}
			&\nrm{\DTL^M\Psi_{q,j,u}}_{N+\alpha}
			\\
			&\quad\lesssim A_{q,j}
			(\tau_q\ell_q^{-1})^\alpha\tau_q\MM{M,N_R-3j-2,\tau_q^{-1},\ell_{t,q}^{-1}}\Lat{N,N_{R,j}^{(2)}(M),\lambda_q,\ell_q^{-1}}
			\MM{N,N_{R,j}^{(2)}(M),\lambda_q,\ell_q^{-1}}.
		\end{aligned}
	\end{equation}
	Take $M,N\in\mathbb N_0$ with $M\leq N_{\mathrm{ind}}+1$. For $M=0$, assume $N\leq K_j^{(3)}$. For $M\geq1$, assume $M+N\leq K_j^{(2)}$. Then
	\begin{align}
		\nrm{\DTL^M\Psi_{q,j,u}}_{N+1+\alpha}
		&\lesssim A_{q,j}
		(\tau_q\ell_q^{-1})^\alpha\tau_q\MM{M,N_R-3j-2,\tau_q^{-1},\ell_{t,q}^{-1}}\notag\\
		&\quad\cdot\Lat{N+1,N_{R,j}^{(2)}(M),\lambda_q,\ell_q^{-1}}
		\MM{N+1,N_{R,j}^{(2)}(M),\lambda_q,\ell_q^{-1}},
		\label{eq:Newton-local-vector-potential-auxiliary-estimate}\\
		\nrm{\DTL^M\omega_{q,j,u}}_{N+\alpha}
		&\lesssim A_{q,j}
		(\tau_q\ell_q^{-1})^\alpha\tau_q\MM{M,N_R-3j-2,\tau_q^{-1},\ell_{t,q}^{-1}}\notag\\
		&\quad\cdot\Lat{N+1,N_{R,j}^{(2)}(M),\lambda_q,\ell_q^{-1}}
		\MM{N+1,N_{R,j}^{(2)}(M),\lambda_q,\ell_q^{-1}},
		\label{eq:Newton-local-incompressible-solution-estimate}\\
		\nrm{\DTL^M\tm_{q,j+1}\inc}_{N+\alpha}
		&\lesssim A_{q,j}
		(\tau_q\ell_q^{-1})^\alpha\tau_q\MM{M,N_R-3j-2,\tau_q^{-1},\ell_{t,q}^{-1}}\notag\\
		&\quad\cdot\Lat{N+1,N_{R,j}^{(2)}(M),\lambda_q,\ell_q^{-1}}
		\MM{N+1,N_{R,j}^{(2)}(M),\lambda_q,\ell_q^{-1}}.
		\label{eq:Newton-incompressible-momentum-estimate}
	\end{align}
	For $M\leq N_{\mathrm{ind}}$ and $M+N\leq K_j^{(2)}$, we have
	\begin{equation}\label{eq:Newton-pressure-potential-estimate}
		\begin{aligned}
			&\nrm{\DTL^M\pi_{q,j,u}}_{N+\alpha}
			+\nrm{\DTL^M\pi_{q,j+1}\inc}_{N+\alpha}\\
			&\qquad\lesssim\Lambda_q A_{q,j}
			(\tau_q\ell_q^{-1})^\alpha\tau_q
			\MM{M,N_R-3j-3,\tau_q^{-1},\ell_{t,q}^{-1}}
			\Lat{N,N_{R,j}^{(3)}(M),\lambda_q,\ell_q^{-1}}
			\MM{N,N_{R,j}^{(3)}(M),\lambda_q,\ell_q^{-1}}.
		\end{aligned}
	\end{equation}
	In the same range, the gluing stress satisfies
	\begin{equation}
		\begin{aligned}
			&\nrm{\DTL^MR_{q,j+1}^{\mathrm{gluing}}}_{N+\alpha}\\
			&\qquad\lesssim A_{q,j}
			(\tau_q\ell_q^{-1})^\alpha\MM{M,N_R-3j-2,\tau_q^{-1},\ell_{t,q}^{-1}}
			\Lat{N,N_{R,j}^{(2)}(M),\lambda_q,\ell_q^{-1}}\MM{N,N_{R,j}^{(2)}(M),\lambda_q,\ell_q^{-1}}.
		\end{aligned}
		\label{eq:Newton-gluing-stress-estimate}
	\end{equation}
	Moreover,
	$$
	\supp_t R_{q,j+1}^{\mathrm{gluing}}
	\subset\bigcup_u\left(
	\supp\partial_t\theta_u
	\cup\supp\partial_t\tilde\theta_u
	\right).
	$$
	For each fixed exponent $\alpha_r$, the same estimates hold with $\alpha$ replaced by $\alpha_r$.
	The implicit constants may depend on $\rho^*,\mathfrak M,N_{\mathrm{max}},\alpha,\alpha_r$.
	They are independent of $q,j,u,a,C_*$ and $A_{q,j}$.
\end{prop}

\begin{proof}
	We first estimate $\Psi_{q,j,u}$, $\omega_{q,j,u}$, and $\pi_{q,j,u}$. We then use these estimates in
	\eqref{def:Newton-glued-incompressible-perturbation} and \eqref{def:Newton-gluing-stress}.
	Apply Proposition~\ref{prop:linearized-high-material-derivatives} with 
	$$
	t_0=t_u,\qquad\tau=\tau_q,\qquad s=2,\qquad U=\overline u_q, \qquad G=-\theta_u^2W_{q,j}^{\mathrm{mol}}-W\left[S_{q,j+1,u}^{\mathrm{cut,skew}}\right].
	$$
	The initial data in
	\eqref{eq:Newton-local-incompressible-linearized-system} are zero. The source $G$ satisfies
	\begin{equation}\label{eq:Newton-local-incompressible-source-estimate}
	\begin{aligned}
		\nrm{\DTL^MG}_{N+\alpha}
		&\lesssim\sum_{r=0}^M\|\partial_t^r(\theta_u^2)\|_0
		\nrm{\DTL^{M-r}W_{q,j}^{\mathrm{mol}}}_{N+\alpha}
		+\nrm{\DTL^MW[S_{q,j+1,u}^{\mathrm{cut,skew}}]}_{N+\alpha}\\
		&\lesssim A_{q,j}(\tau_q\ell_q^{-1})^\alpha
		\MM{M,N_R-3j-2,\tau_q^{-1},\ell_{t,q}^{-1}}\Lat{N,N_{R,j}^{(2)}(M),\lambda_q,\ell_q^{-1}}\MM{N,N_{R,j}^{(2)}(M),\lambda_q,\ell_q^{-1}},
	\end{aligned}
	\end{equation}
	for $M\leq N_{\mathrm{ind}}$ and $M+N\leq K_j^{(2)}$. The last inequality uses \eqref{eq:Newton-mollified-source-estimate}, \eqref{eq:Newton-cutoff-stress-potential-estimate}, and $\|\partial_t^r(\theta_u^2)\|_0\lesssim\tau_q^{-r}$.
	Take
	$$
	\begin{aligned}
		&K^*=M_*=N_{\mathrm{ind}},&& K_G=N_t=N_R-3j-2,&& K^\sharp=K_j^{(2)}+1,&& A_U=C\delta_q^{\frac12},&& A_G=C A_{q,j}(\tau_q\ell_q^{-1})^\alpha,\\
		&\lambda_U=\lambda_G=\lambda_q,&& \tilde\lambda_U=\tilde\lambda_G=\ell_q^{-1},&& \mu_U=\tilde\mu_U=\Lambda_q,&& \mu_G=\tau_q^{-1},&& \tilde\mu_G=\ell_{t,q}^{-1}.
	\end{aligned}
	$$
	Here $C$ covers the implicit constants in \eqref{est:uniform-grad-om-over-orho}, \eqref{est:holder-grad-om-over-orho}, and \eqref{eq:Newton-local-incompressible-source-estimate} at the required derivative orders and fixed H{\"o}lder exponents.
	The derivative parameters satisfy $K^\sharp\geq K^*+2$ and $K_G+1=N_R-3j-1\leq K^*$.
	The coefficient estimates and short-time condition follow from \eqref{est:uniform-grad-om-over-orho},
	\eqref{est:holder-grad-om-over-orho}, and \eqref{eq:Newton-scale-separation}.
	By \eqref{eq:Newton-acoustic-combined-material-scales},
	the combined material scales are $\tau_q^{-1}$ and $\ell_{t,q}^{-1}$.
	The same $A_G$ works at the smaller fixed H{\"o}lder exponents.
	For $M\geq1$, use
	\begin{equation}\label{eq:Newton-linearized-material-shift}\MM{M-1,N_R-3j-2,\tau_q^{-1},\ell_{t,q}^{-1}}\leq\tau_q\MM{M,N_R-3j-2,\tau_q^{-1},\ell_{t,q}^{-1}}.\end{equation}
	
	At $M=0$, use \eqref{eq:linearized-solution-material-zero}--\eqref{eq:linearized-potential-material-zero}.
	We obtain the estimate for $\Psi_{q,j,u}$ when $N\leq K^\sharp-1=K_j^{(2)}$ and the joint estimate for $\omega_{q,j,u}$ and $\Psi_{q,j,u}$ when $N\leq K^\sharp-2=K_j^{(3)}$.
	For $1\leq M\leq N_{\mathrm{ind}}+1$ and $M+N\leq K^\sharp-1=K_j^{(2)}$,
	\eqref{eq:potential-spatial-material-scale-conclusion} and \eqref{eq:Newton-linearized-material-shift} give
	$$
	\begin{aligned}
	\nrm{\DTL^M\Psi_{q,j,u}}_{N+\alpha}
	&\lesssim\tau_q A_G\MM{M,K_G,\tau_q^{-1},\ell_{t,q}^{-1}}\Lat{N,(K_G-M+1)_+,\lambda_q,\ell_q^{-1}}
	\MM{N,(K_G-M+1)_+,\lambda_q,\ell_q^{-1}}\\
	&\lesssim\tau_q A_{q,j}(\tau_q\ell_q^{-1})^\alpha
	\MM{M,K_G,\tau_q^{-1},\ell_{t,q}^{-1}}\Lat{N,N_{R,j}^{(2)}(M),\lambda_q,\ell_q^{-1}}
	\MM{N,N_{R,j}^{(2)}(M),\lambda_q,\ell_q^{-1}}.
	\end{aligned}
	$$
	Here the last inequality uses $(K_G-M+1)_+\geq N_{R,j}^{(2)}(M)$.
	We obtain \eqref{eq:Newton-local-vector-potential-estimate}.
	For $1\leq M\leq N_{\mathrm{ind}}+1$ and $M+N\leq K_j^{(2)}$, \eqref{eq:linearized-spatial-material-scale-conclusion-Calpha} gives
	$$
	\begin{aligned}
		&\nrm{\DTL^M\omega_{q,j,u}}_{N+\alpha}
		+\nrm{\DTL^M\Psi_{q,j,u}}_{N+1+\alpha}\\
		&\quad\lesssim A_G\Bigl[
		\tau_q\MM{M,K_G,\tau_q^{-1},\ell_{t,q}^{-1}}
		+\MM{M-1,K_G,\tau_q^{-1},\ell_{t,q}^{-1}}\Bigr]\\
		&\qquad\cdot\Lat{N+1,(K_G-M+1)_+,\lambda_q,\ell_q^{-1}}
		\MM{N+1,(K_G-M+1)_+,\lambda_q,\ell_q^{-1}}\\
		&\quad\lesssim\tau_q A_G\MM{M,K_G,\tau_q^{-1},\ell_{t,q}^{-1}}\Lat{N+1,N_{R,j}^{(2)}(M),\lambda_q,\ell_q^{-1}}
		\MM{N+1,N_{R,j}^{(2)}(M),\lambda_q,\ell_q^{-1}}.
	\end{aligned}
	$$
	For the last inequality, we used \eqref{eq:Newton-linearized-material-shift}. The estimates \eqref{eq:Newton-local-vector-potential-auxiliary-estimate}--\eqref{eq:Newton-local-incompressible-solution-estimate} follow.
	
	For $\pi_{q,j,u}$, apply \eqref{eq:linearized-pressure-material-zero} when $M=0$.
	For $M\geq1$, use \eqref{eq:pressure-spatial-material-scale-conclusion} and \eqref{eq:Newton-linearized-material-shift}.
	We have, for $0\leq M\leq N_{\mathrm{ind}}$ and $M+N\leq K_j^{(2)}$,
	$$
	\begin{aligned}
		\nrm{\DTL^M\pi_{q,j,u}}_{N+\alpha}
		&\lesssim\tau_q A_UA_G\MM{M,K_G,\tau_q^{-1},\ell_{t,q}^{-1}}\Lat{N+1,(K_G-M)_+,\lambda_q,\ell_q^{-1}}
		\MM{N+1,(K_G-M)_+,\lambda_q,\ell_q^{-1}}\\
		&\lesssim\tau_q\Lambda_q A_{q,j}(\tau_q\ell_q^{-1})^\alpha
		\MM{M,K_G-1,\tau_q^{-1},\ell_{t,q}^{-1}}\Lat{N,N_{R,j}^{(3)}(M),\lambda_q,\ell_q^{-1}}
		\MM{N,N_{R,j}^{(3)}(M),\lambda_q,\ell_q^{-1}}.
	\end{aligned}
	$$
	In the last inequality, we used \eqref{eq:Newton-one-threshold-weight-shift}.
	By the local estimate for $\pi_{q,j,u}$, $\|\partial_t^r\tilde\theta_u\|_0\lesssim\tau_q^{-r}$,
	and the fact that at most two indices contribute to the sums in \eqref{def:Newton-glued-incompressible-perturbation} at each time,
	we obtain \eqref{eq:Newton-pressure-potential-estimate}.
	For $\tm_{q,j+1}\inc$, use \eqref{eq:Newton-local-incompressible-solution-estimate} in \eqref{def:Newton-glued-incompressible-perturbation}.
	When all material derivatives fall on $\tilde\theta_u$, the factor $\omega_{q,j,u}$ has material order zero.
	For $M\geq1$, the required spatial range follows from $N\leq K_j^{(2)}-M\leq K_j^{(3)}$.
	Thus \eqref{eq:Newton-incompressible-momentum-estimate} holds in its stated range.
	
	For $R_{q,j+1}^{\mathrm{gluing}}$, apply Lemma~\ref{lem:direct-symmetric-curl-material-estimate} with
	$\upsilon=\overline u_q$ and $G=\Psi_{q,j,u}$. Take
	$$
	\begin{aligned}
		&K^*=N_{\mathrm{ind}}-1,&& M_*=N_{\mathrm{ind}},&& K_G=N_t=N_R-3j-2,&& K^\sharp=K_j^{(2)}+1,&& A_\upsilon=C\delta_q^{\frac12},&& s=2,\\
		&\lambda_\upsilon=\lambda_G=\lambda_q,&& \tilde\lambda_\upsilon=\tilde\lambda_G=\ell_q^{-1},&& \mu_\upsilon=\tilde\mu_\upsilon=\Lambda_q,&& \mu_G=\tau_q^{-1}, \tilde\mu_G=\ell_{t,q}^{-1},
		&&A_G=C\tau_q A_{q,j}(\tau_q\ell_q^{-1})^\alpha.
	\end{aligned}
	$$
	Estimate~\eqref{eq:Newton-local-vector-potential-estimate} verifies \eqref{eq:symmetric-curl-G-Holder-assumption} at the fixed exponents,
	with $C$ covering its implicit constants and those in \eqref{est:uniform-grad-om-over-orho} and \eqref{est:holder-grad-om-over-orho} over the required derivative orders and fixed H{\"o}lder exponents.
	The coefficient hypotheses follow from \eqref{est:uniform-grad-om-over-orho} and \eqref{est:holder-grad-om-over-orho}.
	The combined material scales are again given by \eqref{eq:Newton-acoustic-combined-material-scales}.
	The estimate $\|\partial_t^{r+1}\tilde\theta_u\|_0\lesssim\tau_q^{-r-1}$ and
	\eqref{eq:direct-symmetric-curl-high-material}--\eqref{eq:direct-symmetric-curl-material-zero} give
	$$
	\begin{aligned}
		&\Nrm{\DTL^M\left((\partial_t\tilde\theta_u)\mathcal R\curl\Psi_{q,j,u}\right)}_{N+\alpha}\\
		&\quad\lesssim A_{q,j}(\tau_q\ell_q^{-1})^\alpha
		\Lat{N,N_{R,j}^{(2)}(M),\lambda_q,\ell_q^{-1}}\MM{N,N_{R,j}^{(2)}(M),\lambda_q,\ell_q^{-1}}
		\sum_{r=0}^M\tau_q^{-r}\MM{M-r,N_R-3j-2,\tau_q^{-1},\ell_{t,q}^{-1}}\\
		&\quad\lesssim A_{q,j}(\tau_q\ell_q^{-1})^\alpha
		\MM{M,N_R-3j-2,\tau_q^{-1},\ell_{t,q}^{-1}}\Lat{N,N_{R,j}^{(2)}(M),\lambda_q,\ell_q^{-1}}\MM{N,N_{R,j}^{(2)}(M),\lambda_q,\ell_q^{-1}}.
	\end{aligned}
	$$
	For $0\leq M\leq N_{\mathrm{ind}}$ and $M+N\leq K_j^{(2)}$, sum the preceding estimate over $u$. At each time, at most two indices contribute.
	Combining this sum with \eqref{eq:Newton-cutoff-stress-estimate} for $S_{q,j+1,u}^{\mathrm{cut,sym}}$ proves \eqref{eq:Newton-gluing-stress-estimate}.
	The temporal support inclusion for $R_{q,j+1}^{\mathrm{gluing}}$ follows from \eqref{def:Newton-gluing-stress}.
	$S_{q,j+1,u}^{\mathrm{cut,sym}}$ is supported in time in $\supp\partial_t\theta_u$,
	and $\mathcal R\curl$ does not enlarge temporal support.
	Repeating the proof with $\alpha_r$ in place of $\alpha$ gives the estimates at each fixed exponent.
\end{proof}

\subsection{Estimates for the updated state and the stress tensors}
\label{subsec:Newton-stress-state-estimates}

We estimate the interior stress $R_{q,j}$, the updated state $(\rho_{q,j},m_{q,j})$, and the stresses $S_{q,j}$ and $P_{q+1,j}$ along the flow of $\overline u_q$. We also estimate $u_{q,J_*}-\overline u_q$ and prove the terminal temporal support inclusion.

\begin{prop}
	\label{prop:Newton-interior-stress-estimates}
	Assume the pressure-law conditions in Section~\ref{sec:main-iteration}
	and the hypotheses of Proposition~\ref{prop:uniform-mollification-estimates}.
	For $a>a_{\mathrm{Newton}}$, the following estimates hold at each of the fixed H{\"o}lder exponents.
	For $0\leq j\leq J_*$ and $M,N\in\mathbb N_0$ with $M\leq N_{\mathrm{ind}}-1$ and $M+N\leq K_j^{(0)}$, one has
	\begin{equation}\label{eq:Newton-interior-stress-estimate}
		\nrm{\DTL^MR_{q,j}}_{N+\alpha}
		\leq A_{q,j}\MM{M,N_R-3j,\tau_q^{-1},\ell_{t,q}^{-1}}
		\La{N,N_{R,j}^{(0)}(M),\lambda_q,\ell_q^{-1}}\MM{N,N_{R,j}^{(0)}(M),\lambda_q,\ell_q^{-1}}.
	\end{equation}
	At $j=0$, \eqref{eq:Newton-interior-stress-estimate} also holds for $0\leq M\leq N_{\mathrm{max}}$ and $\min\{M,N_R\}+N\leq N_{\mathrm{max}}-1$.
	At $j=J_*$, we have
	\begin{equation}\label{eq:Newton-terminal-interior-stress-estimate}
		\nrm{\DTL^MR_{q,J_*}}_{N+\alpha}
		\leq A_{q,J_*}\MM{M,N_R-3J_*,\tau_q^{-1},\ell_{t,q}^{-1}}
		\La{N,N_{R,J_*}^{(0)}(M),\lambda_q,\ell_q^{-1}}\MM{N,N_{R,J_*}^{(0)}(M),\lambda_q,\ell_q^{-1}}.
	\end{equation}
	For $M\leq N_{\mathrm{ind}}-1$ and $M+N\leq K_{J_*}^{(0)}$, one also has
	\begin{equation}\label{eq:Newton-terminal-interior-stress-CN}
		\begin{aligned}
			\nrm{\DTL^MR_{q,J_*}}_{N}
			&\leq A_{q,J_*}\ell_q^{-\alpha}
			\MM{M,N_R-3J_*,\tau_q^{-1},\ell_{t,q}^{-1}}
			\MM{N,N_{R,J_*}^{(0)}(M),\lambda_q,\ell_q^{-1}}.
		\end{aligned}
	\end{equation}
\end{prop}

\begin{proof}
	The case $j=0$ follows from \eqref{eq:initial-Newton-flow-mollified-stress-estimate}, the definition of $A_{q,0}$, and $\Lambda_q\leq\tau_q^{-1}$. We prove \eqref{eq:Newton-interior-stress-estimate} for $j\geq1$ by induction. For $0\leq j<J_*$, assume it holds at substep $j$. This verifies \eqref{eq:Newton-interior-input-estimate}, so Lemma~\ref{lem:Newton-source-flow-mollification} and the perturbation estimates apply.
	Since $N_{R,j}^{(3)}(M)=N_{R,j+1}^{(0)}(M)$,
	\eqref{eq:Newton-interior-stress-potential-estimate}, \eqref{def:Newton-interior-stress}, and \eqref{eq:Newton-pressure-potential-estimate} give \eqref{eq:Newton-final-interior-stress-estimate}.
	For $0\leq M\leq N_{\mathrm{ind}}-1$ and $M+N\leq K_{j+1}^{(0)}=K_j^{(3)}$, the choice of $C_*$ gives
		\begin{align}
			&\nrm{\DTL^MR_{q,j+1}}_{N+\alpha}\nonumber\\
			&\quad\leq C_*\tau_q\Lambda_q(\tau_q\ell_q^{-1})^\alpha A_{q,j}
			\MM{M,N_R-3(j+1),\tau_q^{-1},\ell_{t,q}^{-1}}\Lat{N,N_{R,j+1}^{(0)}(M),\lambda_q,\ell_q^{-1}}\MM{N,N_{R,j+1}^{(0)}(M),\lambda_q,\ell_q^{-1}}\nonumber\\
			&\quad\leq A_{q,j+1}
			\MM{M,N_R-3(j+1),\tau_q^{-1},\ell_{t,q}^{-1}}\La{N,N_{R,j+1}^{(0)}(M),\lambda_q,\ell_q^{-1}}\MM{N,N_{R,j+1}^{(0)}(M),\lambda_q,\ell_q^{-1}}.\label{eq:Newton-substep-restored-material-range}
		\end{align}
	This completes the induction for \eqref{eq:Newton-interior-stress-estimate}.
	The same induction applies at the other fixed H{\"o}lder exponents.
	Taking $j=J_*$ in \eqref{eq:Newton-interior-stress-estimate} proves \eqref{eq:Newton-terminal-interior-stress-estimate}.
	Using $\|\cdot\|_N\leq\|\cdot\|_{N+\alpha}$ and $\Xi_\alpha\leq\ell_q^{-\alpha}$ in \eqref{eq:Newton-terminal-interior-stress-estimate} gives \eqref{eq:Newton-terminal-interior-stress-CN}.
\end{proof}

\begin{prop}
	\label{prop:Newton-updated-state-estimates}
	Assume the pressure-law conditions in Section~\ref{sec:main-iteration}
	and the hypotheses of Proposition~\ref{prop:uniform-mollification-estimates}.
	For $a>a_{\mathrm{Newton}}$ and $0\leq j\leq J_*$, set
	$u_{q,j}:=m_{q,j}/\rho_{q,j}$.
	The density, momentum, and velocity satisfy
	\begin{equation}\label{eq:Newton-state-zero-order-bounds}
	\begin{aligned}
		\nrm{\rho_{q,j}-\rho^*}_0
		\leq\mathfrak M\sum_{i=1}^{q+1}\delta_i
		\leq\frac{\rho^*}{2},\quad
		\nrm{m_{q,j}}_0\leq\mathfrak M\left(1-\frac12\delta_q^{\frac12}\right),\qquad
		\nrm{u_{q,j}}_0\leq\frac{2\mathfrak M}{\rho^*}\left(1-\frac12\delta_q^{\frac12}\right).
	\end{aligned}
	\end{equation}
	In particular, $\rho_{q,j}\in[\rho^*/2,3\rho^*/2]$ at each substep.
	For $0\leq j\leq J_*$ and $M,N\in\mathbb N_0$ with $1\leq M+N\leq N_{\mathrm{ind}}$, one has
	\begin{equation}\label{eq:Newton-sharp-integer-state-coefficients}
		\begin{aligned}
			\|\DTL^M(\rho_{q,j}-\rho^*)\|_N+\|\DTL^Mm_{q,j}\|_N+\|\DTL^Mu_{q,j}\|_N\lesssim\delta_q^{\frac12}\MM{M,N_R-3j,\tau_q^{-1},\ell_{t,q}^{-1}}\ell_q^{-N}.
		\end{aligned}
	\end{equation}
	Spatial interpolation then gives, for $1\leq M+N\leq N_{\mathrm{ind}}-1$,
	\begin{equation}\label{eq:Newton-updated-state-material-estimate}
		\begin{aligned}
			&\nrm{\DTL^M(\rho_{q,j}-\rho^*)}_{N+\alpha}
			+\nrm{\DTL^Mm_{q,j}}_{N+\alpha}
			+\nrm{\DTL^Mu_{q,j}}_{N+\alpha}\lesssim
			\delta_q^{\frac12}
			\MM{M,N_R-3j,\tau_q^{-1},\ell_{t,q}^{-1}}
			\ell_q^{-N-\alpha}.
		\end{aligned}
	\end{equation}

	For $0\leq\nu\leq j<J_*$ and $M,N\in\mathbb N_0$ with $M+N\leq N_{\mathrm{ind}}$, the individual perturbations satisfy
	\begin{equation}\label{eq:Newton-quadratic-input-material-estimates}
	\begin{aligned}
		&\nrm{\DTL^M\trho_{q,\nu+1}\comp}_N+\nrm{\DTL^M\tm_{q,\nu+1}\comp}_N
		\lesssim A_{q,\nu}\ell_q^{-2\alpha_2}
		\MM{M,N_R-3j-2,\tau_q^{-1},\ell_{t,q}^{-1}}
		\MM{N,N_{R,j}^{(2)}(M),\lambda_q,\ell_q^{-1}},\\
		&\nrm{\DTL^M\tm_{q,\nu+1}\inc}_N\lesssim A_{q,\nu}\tau_q(\tau_q\ell_q^{-1})^{\alpha_2}
		\ell_q^{-2\alpha_2}
		\MM{M,N_R-3j-2,\tau_q^{-1},\ell_{t,q}^{-1}}
		\MM{N+1,N_{R,j}^{(2)}(M),\lambda_q,\ell_q^{-1}}.
	\end{aligned}
	\end{equation}
	For $0\leq\nu<J_*$ and $M,N\in\mathbb N_0$ with $M+N\leq N_{\mathrm{ind}}$, we also have
	\begin{equation}\label{eq:Newton-uniform-perturbation-material-estimate}
	\begin{aligned}
		\nrm{\DTL^M\trho_{q,\nu+1}\comp}_N+\nrm{\DTL^M\tm_{q,\nu+1}\comp}_N+\nrm{\DTL^M\tm_{q,\nu+1}\inc}_N
		&\lesssim A_{q,\nu}\delta_q^{-\frac12}\lambda_{q+1}^{-3\alpha}\ell_{t,q}^{-M}\ell_q^{-N}.
	\end{aligned}
	\end{equation}
	For $0\leq j\leq J_*$ and $M,N\in\mathbb N_0$ with $M\leq N_{\mathrm{ind}}$, assume $N\leq K_j^{(0)}$ when $M=0$, and $M+N\leq K_j^{(0)}+1$ when $M\geq1$. The following estimate also holds with $m_{q,j}-m_q$ in place of $\rho_{q,j}-\rho_q$. 
	\begin{equation}\label{eq:Newton-accumulated-state-increments-integer}
	\begin{aligned}
		\nrm{\DTL^M(\rho_{q,j}-\rho_q)}_N
		&\leq\sum_{\nu=0}^{j-1}\left(
		\nrm{\DTL^M\trho_{q,\nu+1}\comp}_N
		+\nrm{\DTL^M\tm_{q,\nu+1}\comp}_N
		+\nrm{\DTL^M\tm_{q,\nu+1}\inc}_N\right)\\
		&\lesssim A_{q,0}\delta_q^{-\frac12}\lambda_{q+1}^{-3\alpha}
		\MM{M,N_R-3j,\tau_q^{-1},\ell_{t,q}^{-1}}
		\MM{N,N_{R,j}^{(0)}(M),\lambda_q,\ell_q^{-1}}.
	\end{aligned}
	\end{equation}
	The implicit constants may depend on $\rho^*,\mathfrak M,C_*,N_{\mathrm{max}},\alpha,\alpha_r$.
	They are independent of $q$, $j$, and $a$.
\end{prop}

\begin{proof}
	We first prove the density estimate in \eqref{eq:Newton-state-zero-order-bounds} and deduce $\rho_{q,j}\in[\rho^*/2,3\rho^*/2]$. At $j=0$, use \eqref{ind:rho-zero} and \eqref{ind:rho-gap}.
	We then estimate the perturbation sums in \eqref{eq:Newton-accumulated-state-increments-integer} and the state $(\rho_{q,j},m_{q,j},u_{q,j})$.
	Finally, we prove \eqref{eq:Newton-quadratic-input-material-estimates}--\eqref{eq:Newton-uniform-perturbation-material-estimate} for the stress estimates below.
	Proposition~\ref{prop:Newton-interior-stress-estimates} gives \eqref{eq:Newton-interior-stress-estimate} at every substep. For $0\leq j<J_*$, apply \eqref{eq:Newton-total-compressible-density-estimate} with H{\"o}lder exponent $\alpha_2$. The amplitude recursion and the definition of $A_{q,0}$ give
	$$
	\begin{aligned}
		\nrm{\rho_{q,j+1}-\rho_q}_0
		\leq\sum_{\nu=0}^{j}\nrm{\trho_{q,\nu+1}\comp}_0
		\lesssim\lambda_q^{\frac{\alpha}{2}}\sum_{\nu=0}^{j}A_{q,\nu}\lesssim\lambda_q^{-\frac{15\alpha}{2}}\delta_{q+1}\leq\mathfrak M\delta_{q+1}.
	\end{aligned}
	$$
	For sufficiently large $a_{\mathrm{Newton}}$, the last inequality holds.
	Then \eqref{ind:rho-zero} and \eqref{ind:rho-gap} give the density estimate in \eqref{eq:Newton-state-zero-order-bounds} at substep $j+1$.

	For $0\leq j<J_*$, apply \eqref{eq:Newton-total-compressible-density-estimate}--\eqref{eq:Newton-total-compressible-momentum-estimate} and \eqref{eq:Newton-incompressible-momentum-estimate} at $M=N=0$ with exponent $\alpha_2$ to obtain
	\begin{equation}\label{eq:Newton-zero-order-perturbation-amplitudes}
	\begin{aligned}
	&\nrm{\trho_{q,j+1}\comp}_{\alpha_2}+\nrm{\tm_{q,j+1}\comp}_{\alpha_2}+\nrm{\tm_{q,j+1}\inc}_{\alpha_2}\lesssim_{\mathfrak M}A_{q,j}\tau_q\lambda_q(\tau_q\ell_q^{-1})^{\alpha_2}\lambda_q^{2\alpha_2}\lesssim A_{q,j}\delta_q^{-\frac12}\lambda_{q+1}^{-3\alpha}.
	\end{aligned}
	\end{equation}
	For the second inequality, use the scale definitions to write
	$$
	\frac{\tau_q\lambda_q(\tau_q\ell_q^{-1})^{\alpha_2}\lambda_q^{2\alpha_2}}
	{\delta_q^{-\frac12}\lambda_{q+1}^{-3\alpha}}
	\lesssim\lambda_q^{\alpha_2(\frac{b+7}{4}+\beta-4b\alpha)-b\alpha}\leq\lambda_q^{-\frac{5\alpha}{12}}.
	$$
	Since $\alpha_2=\alpha/4$, $b>1$, and $\beta<1/3$, the exponent is less than $-5\alpha/12$. As $\lambda_q\geq a$, one choice of $a_{\mathrm{Newton}}$ makes the second inequality in \eqref{eq:Newton-zero-order-perturbation-amplitudes} valid for every $q$.
	Summing \eqref{eq:Newton-zero-order-perturbation-amplitudes} and using $j\leq J_*\leq N_{\mathrm{max}}$, we obtain, for $0\leq j\leq J_*$,
	\begin{equation}\label{eq:Newton-zero-order-accumulated-perturbations}
	\begin{aligned}
		&\sum_{\nu=0}^{j-1}\left(
		\nrm{\trho_{q,\nu+1}\comp}_{\alpha_2}
		+\nrm{\tm_{q,\nu+1}\comp}_{\alpha_2}
		+\nrm{\tm_{q,\nu+1}\inc}_{\alpha_2}\right)\lesssim{}A_{q,0}\delta_q^{-\frac12}
		\lambda_{q+1}^{-3\alpha}.
	\end{aligned}
	\end{equation}

	By \eqref{def:lambda-delta}, \eqref{eq:Newton-zero-order-accumulated-perturbations}, and the definition of $A_{q,0}$, we have
	$$
	\begin{aligned}
		\nrm{m_{q,j}-m_q}_0
		&\lesssim\lambda_q^{-8\alpha}
		\left(\frac{\lambda_q}{\lambda_{q+1}}\right)^{2\beta}
		\lambda_{q+1}^{-3\alpha}\delta_q^{\frac12}
		\leq\frac12\delta_q^{\frac12}.
	\end{aligned}
	$$
	Here the last inequality holds for sufficiently large $a_{\mathrm{Newton}}$.
	The momentum increment estimate and \eqref{ind:m-zero} give
	$$
	\begin{aligned}
		\nrm{m_{q,j}}_0
		\leq\nrm{m_q}_0+\nrm{m_{q,j}-m_q}_0
		\leq\mathfrak M(1-\delta_q^{\frac12})
		+\nrm{m_{q,j}-m_q}_0
		\leq\mathfrak M\left(1-\frac12\delta_q^{\frac12}\right).
	\end{aligned}
	$$
	The lower bound $\rho_{q,j}\geq\rho^*/2$ gives the velocity estimate for $u_{q,j}$ in \eqref{eq:Newton-state-zero-order-bounds}.

	The following inequalities for the derivative weights will be used repeatedly below. Let $K,M,N\in\mathbb N_0$ and $k\in\mathbb N$. For $M_i,N_i\in\mathbb N_0$ with $M_1+\cdots+M_k=M$ and $N_1+\cdots+N_k=N$, one has
	\begin{equation}\label{eq:Newton-quadratic-weight-product}
	\begin{aligned}
		\prod_{i=1}^k\Bigl(
		\MM{M_i,K,\tau_q^{-1},\ell_{t,q}^{-1}}
		\MM{N_i,(K-M_i)_+,\lambda_q,\ell_q^{-1}}\Bigr)
		\leq\MM{M,K,\tau_q^{-1},\ell_{t,q}^{-1}}
		\MM{N,(K-M)_+,\lambda_q,\ell_q^{-1}}.
	\end{aligned}
	\end{equation}
	This is \eqref{eq:appendix-material-mollification-weight-product} with $K_f=K$, $\ell=\ell_q$, $\lambda=\lambda_q$, $\mu=\tau_q^{-1}$, and $\tilde\mu=\ell_{t,q}^{-1}$.
	For the same $K,M,N$, \eqref{eq:material-mollification-spatial-comparisons} gives
	\begin{equation}\label{eq:Newton-two-scale-weight-comparison}
	\begin{aligned}
		\MM{M,K,\tau_q^{-1},\ell_{t,q}^{-1}}
		\MM{N,(K-M)_+,\lambda_q,\ell_q^{-1}}
		\leq(\ell_q\lambda_q)^{\min\{M+N,K\}}
		\ell_{t,q}^{-M}\ell_q^{-N}.
	\end{aligned}
	\end{equation}

	Fix $1\leq j\leq J_*$. Take $M,N\in\mathbb N_0$ in the derivative ranges stated for \eqref{eq:Newton-accumulated-state-increments-integer}.
	Apply \eqref{eq:Newton-total-compressible-density-estimate}--\eqref{eq:Newton-total-compressible-momentum-estimate} and \eqref{eq:Newton-incompressible-momentum-estimate} with exponent $\alpha_2$ for $0\leq\nu<j$.
	At $M=0$, the range follows from $K_\nu^{(3)}=K_{\nu+1}^{(0)}\geq K_j^{(0)}$. For $M\geq1$, use $K_\nu^{(2)}=K_{\nu+1}^{(0)}+1\geq K_j^{(0)}+1$.
	For $\tm_{q,\nu+1}\inc$, \eqref{eq:Newton-one-threshold-weight-shift} with $K=N_R-3\nu-2$ gives
	$$
	\frac{\MM{M,N_R-3\nu-2,\tau_q^{-1},\ell_{t,q}^{-1}}
	\MM{N+1,N_{R,\nu}^{(2)}(M),\lambda_q,\ell_q^{-1}}}
	{\MM{M,N_R-3\nu-3,\tau_q^{-1},\ell_{t,q}^{-1}}
	\MM{N,N_{R,\nu}^{(3)}(M),\lambda_q,\ell_q^{-1}}}
	\leq\lambda_q.
	$$
	Since $\nu<j$, lower the resulting thresholds to $N_R-3j$ and $N_{R,j}^{(0)}(M)$.
	The scale definitions give
	\begin{equation}\label{eq:Newton-quadratic-Holder-amplitude-factors}
	\begin{aligned}
		&\ell_q^{-2\alpha_2}
		\lesssim\delta_q^{-\frac12}\lambda_{q+1}^{-3\alpha},\quad
		\tau_q\lambda_q(\tau_q\ell_q^{-1})^{\alpha_2}\ell_q^{-2\alpha_2}
		\lesssim\delta_q^{-\frac12}\lambda_{q+1}^{-4\alpha}
		\lambda_q^{\alpha_2\left(\frac{3b+5}{4}+\beta-4b\alpha\right)}\lesssim\delta_q^{-\frac12}\lambda_{q+1}^{-3\alpha}\lambda_q^{-\frac{5\alpha}{12}}
		\lesssim\delta_q^{-\frac12}\lambda_{q+1}^{-3\alpha}.
	\end{aligned}
	\end{equation}
	Here the first inequality uses $(25b+3)\alpha/8<4b\alpha<\beta$.
	The factor $\lambda_q^{-5\alpha/12}$ follows from $\alpha_2=\alpha/4$, $b>1$, $\beta<1/3$, and $\lambda_{q+1}\sim\lambda_q^b$.
	Using $\Xi_{\alpha_2}^2\leq\ell_q^{-2\alpha_2}$ and summing over $0\leq\nu<j$, we obtain
	$$
	\begin{aligned}
		\nrm{\DTL^M(\rho_{q,j}-\rho_q)}_N&\leqslant\sum_{\nu=0}^{j-1}\left(
		\nrm{\DTL^M\trho_{q,\nu+1}\comp}_N
		+\nrm{\DTL^M\tm_{q,\nu+1}\comp}_N
		+\nrm{\DTL^M\tm_{q,\nu+1}\inc}_N\right)\\
		&\lesssim\left(\sum_{\nu=0}^{j-1}A_{q,\nu}\right)
		\left[1+\tau_q\lambda_q(\tau_q\ell_q^{-1})^{\alpha_2}\right]\ell_q^{-2\alpha_2}
		\MM{M,N_R-3j,\tau_q^{-1},\ell_{t,q}^{-1}}
		\MM{N,N_{R,j}^{(0)}(M),\lambda_q,\ell_q^{-1}}\\
		&\lesssim A_{q,0}\delta_q^{-\frac12}\lambda_{q+1}^{-3\alpha}
		\MM{M,N_R-3j,\tau_q^{-1},\ell_{t,q}^{-1}}
		\MM{N,N_{R,j}^{(0)}(M),\lambda_q,\ell_q^{-1}}.
	\end{aligned}
	$$
	For the last inequality, we used \eqref{eq:Newton-quadratic-Holder-amplitude-factors} and $\sum_{\nu=0}^{j-1}A_{q,\nu}\leq J_*A_{q,0}$. The first inequality in \eqref{eq:Newton-accumulated-state-increments-integer} follows from $\rho_{q,j}-\rho_q=\sum_{\nu=0}^{j-1}\trho_{q,\nu+1}\comp$. It also holds for the momentum increment, since $m_{q,j}-m_q=\sum_{\nu=0}^{j-1}(\tm_{q,\nu+1}\comp+\tm_{q,\nu+1}\inc)$. At $j=0$, both increments and the sum vanish.
	For $K,M,N\in\mathbb N_0$ with $M+N\leq N_{\mathrm{ind}}$, the initial spatial mollification differences satisfy
	\begin{equation}\label{eq:Newton-initial-state-difference-material-estimate}
	\begin{aligned}
	&\|\DTL^M(\rho_q-\orho_q)\|_N+\|\DTL^M(m_q-\om_q)\|_N\\
		&\quad\lesssim(\ell_q\lambda_q)^{N_{\mathrm{ind}}}\delta_q^{\frac12}(\delta_q^{\frac12}\ell_q^{-1})^M\ell_q^{-N}
		\leq(\ell_q\lambda_q)^{N_{\mathrm{ind}}-\min\{M+N,K\}}\delta_q^{\frac12}
		\MM{M,K,\tau_q^{-1},\ell_{t,q}^{-1}}\MM{N,(K-M)_+,\lambda_q,\ell_q^{-1}}.
	\end{aligned}
	\end{equation}
	In the first inequality, we used \eqref{est:uniform-rho-minus-orho}--\eqref{est:uniform-m-minus-om}. For the second, \eqref{eq:material-mollification-spatial-comparisons} and \eqref{eq:Newton-scale-separation} give
	$\delta_q^{\frac12}\ell_q^{-1}\ell_{t,q}=\ell_q\lambda_q$ and $(\ell_q\lambda_q)^2\leq\ell_{t,q}\tau_q^{-1}$. Thus
	$$
		(\delta_q^{\frac12}\ell_q^{-1})^M\ell_q^{-N}\leq(\ell_q\lambda_q)^{-\min\{M+N,K\}}\MM{M,K,\tau_q^{-1},\ell_{t,q}^{-1}}\MM{N,(K-M)_+,\lambda_q,\ell_q^{-1}}.
	$$
	
	To prove \eqref{eq:Newton-sharp-integer-state-coefficients}, apply \eqref{eq:Newton-initial-state-difference-material-estimate} with $K=N_R-3j$. At $j=0$, the perturbation sums are empty. The choice of $a_{\mathrm{Newton}}$ above ensures $A_{q,0}\delta_q^{-\frac12}\lambda_{q+1}^{-3\alpha}\leq\delta_q^{\frac12}$.
	Combining these estimates with \eqref{est:uniform-orho}--\eqref{est:uniform-om} and \eqref{eq:Newton-accumulated-state-increments-integer}, we have, for $1\leq M+N\leq N_{\mathrm{ind}}$,
	\begin{equation}\label{eq:Newton-updated-density-momentum-integer-estimate}
	\begin{aligned}
		\|\DTL^M(\rho_{q,j}-\rho^*)\|_N+\|\DTL^Mm_{q,j}\|_N
		&\lesssim\delta_q^{\frac12}\MM{M,N_R-3j,\tau_q^{-1},\ell_{t,q}^{-1}}\ell_q^{-N}.
	\end{aligned}
	\end{equation}
	Using $u_{q,j}=m_{q,j}\rho_{q,j}^{-1}$, the lower bound $\rho_{q,j}\geq\rho^*/2$ from \eqref{eq:Newton-state-zero-order-bounds}, and \eqref{eq:faa-di-bruno}, we obtain
	\begin{equation}\label{eq:Newton-updated-velocity-product-estimate}
	\begin{aligned}
		\|\DTL^Mu_{q,j}\|_N
		\lesssim\|\DTL^Mm_{q,j}\|_N+\sum_{k=1}^{M+N}
		\sum_{M_0+\cdots+M_k=M}
		\sum_{\substack{N_0+\cdots+N_k\leq N\\M_i+N_i\geq1,\ 1\leq i\leq k}}
		\|\nabla^{N_0}\DTL^{M_0}m_{q,j}\|_0
		\prod_{i=1}^k\|\nabla^{N_i}\DTL^{M_i}\rho_{q,j}\|_0.
	\end{aligned}
	\end{equation}
	Use \eqref{eq:Newton-updated-density-momentum-integer-estimate} for the differentiated factors. For $M_0+N_0=0$, use \eqref{eq:Newton-state-zero-order-bounds}.
	Since $k\geq1$, each product in \eqref{eq:Newton-updated-velocity-product-estimate} satisfies
	$$
	\begin{aligned}
		\|\nabla^{N_0}\DTL^{M_0}m_{q,j}\|_0
		\prod_{i=1}^k\|\nabla^{N_i}\DTL^{M_i}\rho_{q,j}\|_0
		&\lesssim\delta_q^{\frac k2}
		\prod_{i=0}^k\MM{M_i,N_R-3j,\tau_q^{-1},\ell_{t,q}^{-1}}
		\ell_q^{-\sum_{i=0}^k N_i}\\
		&\leq\delta_q^{\frac12}
		\MM{M,N_R-3j,\tau_q^{-1},\ell_{t,q}^{-1}}\ell_q^{-N}.
	\end{aligned}
	$$
	Here the last inequality uses \eqref{eq:Newton-quadratic-weight-product}. We obtain \eqref{eq:Newton-sharp-integer-state-coefficients}.
	For $1\leq M+N\leq N_{\mathrm{ind}}-1$, spatial interpolation of \eqref{eq:Newton-sharp-integer-state-coefficients} at orders $N$ and $N+1$ gives \eqref{eq:Newton-updated-state-material-estimate}.

	It remains to prove \eqref{eq:Newton-quadratic-input-material-estimates}--\eqref{eq:Newton-uniform-perturbation-material-estimate}.
	For $0\leq\nu\leq j<J_*$ and $M+N\leq N_{\mathrm{ind}}$, apply \eqref{eq:Newton-total-compressible-density-estimate}--\eqref{eq:Newton-total-compressible-momentum-estimate} and \eqref{eq:Newton-incompressible-momentum-estimate} with exponent $\alpha_2$.
	Their derivative ranges apply since $M+N\leq N_{\mathrm{ind}}\leq K_j^{(3)}\leq K_\nu^{(3)}$. Use $\Xi_{\alpha_2}^2\leq\ell_q^{-2\alpha_2}$ and lower the material and spatial thresholds to $N_R-3j-2$ and $N_{R,j}^{(2)}(M)$, respectively, to obtain \eqref{eq:Newton-quadratic-input-material-estimates}.

	For $0\leq\nu<J_*$ and $M+N\leq N_{\mathrm{ind}}$, estimate the three perturbations in \eqref{eq:Newton-uniform-perturbation-material-estimate} at the mollification scales.
	Use \eqref{eq:Newton-quadratic-input-material-estimates} with $j=\nu$.
	Apply \eqref{eq:Newton-two-scale-weight-comparison} with $K=N_R-3\nu-2$ at spatial orders $N$ and $N+1$.
	Since $K\geq2$, we have
	$$
	\begin{aligned}
		\MM{M,K,\tau_q^{-1},\ell_{t,q}^{-1}}
		\MM{N+1,(K-M)_+,\lambda_q,\ell_q^{-1}}\leq
		(\ell_q\lambda_q)^{\min\{M+N+1,K\}}
		\ell_{t,q}^{-M}\ell_q^{-N-1}
		\leq\lambda_q\ell_{t,q}^{-M}\ell_q^{-N}.
	\end{aligned}
	$$
	Substituting this comparison into \eqref{eq:Newton-quadratic-input-material-estimates} with $j=\nu$, we obtain
	$$
	\begin{aligned}
		\nrm{\DTL^M\trho_{q,\nu+1}\comp}_N+\nrm{\DTL^M\tm_{q,\nu+1}\comp}_N+\nrm{\DTL^M\tm_{q,\nu+1}\inc}_N
		&\lesssim_{\mathfrak M}A_{q,\nu}\left[1+\tau_q\lambda_q(\tau_q\ell_q^{-1})^{\alpha_2}\right]\ell_q^{-2\alpha_2}\ell_{t,q}^{-M}\ell_q^{-N}.
	\end{aligned}
	$$
	Use \eqref{eq:Newton-quadratic-Holder-amplitude-factors} and choose $a_{\mathrm{Newton}}$ large enough to absorb the fixed implicit constant.
	The uniform estimate \eqref{eq:Newton-uniform-perturbation-material-estimate} follows.
\end{proof}

\begin{prop}
	\label{prop:Newton-background-flow-stresses}
	Assume the pressure-law conditions in Section~\ref{sec:main-iteration}
	and the hypotheses of Proposition~\ref{prop:uniform-mollification-estimates}.
	For $a>a_{\mathrm{Newton}}$, the following estimates hold.
	For $M,N\in\mathbb N_0$ with $M\leq N_{\mathrm{ind}}$ and $M+N\leq K_{J_*}^{(0)}$, the density perturbation at the end of the Newton step satisfies
	\begin{equation}\label{eq:Newton-terminal-density-increment-material}
	\begin{aligned}
    \nrm{\DTL^M(\rho_{q,J_*}-\rho_q)}_N
    \lesssim
    A_{q,0}\ell_q^{-\alpha/2}
    \MM{M,N_R-3J_*,\tau_q^{-1},\ell_{t,q}^{-1}}
    \MM{N,N_{R,J_*}^{(0)}(M),\lambda_q,\ell_q^{-1}}.
	\end{aligned}
	\end{equation}
	For $M,N\in\mathbb N_0$ with $M+N\leq N_{\mathrm{ind}}$, the density difference satisfies
	\begin{equation}\label{eq:Newton-terminal-density-difference-two-scale}
	\begin{aligned}
		&\|\DTL^M(\rho_{q,J_*}-\orho_q)\|_N\\
		&\quad\lesssim\left[A_{q,0}\ell_q^{-2\alpha_2}
		+(\ell_q\lambda_q)^{N_{\mathrm{ind}}-\min\{M+N,N_R-3J_*\}}\delta_q\right]\MM{M,N_R-3J_*,\tau_q^{-1},\ell_{t,q}^{-1}}
		\MM{N,N_{R,J_*}^{(0)}(M),\lambda_q,\ell_q^{-1}}.
	\end{aligned}
	\end{equation}
	In the same range, the momentum and velocity differences satisfy
	\begin{equation}\label{eq:Newton-terminal-velocity-difference-two-scale}
		\begin{aligned}
			&\|\DTL^M(m_{q,J_*}-\om_q)\|_N+\|\DTL^M(u_{q,J_*}-\overline u_q)\|_N\\
			&\quad\lesssim\left[A_{q,0}\delta_q^{-\frac12}\lambda_{q+1}^{-3\alpha}+(\ell_q\lambda_q)^{N_{\mathrm{ind}}-N_R+3J_*}\delta_q^{\frac12}\right]\MM{M,N_R-3J_*,\tau_q^{-1},\ell_{t,q}^{-1}}\MM{N,N_{R,J_*}^{(0)}(M),\lambda_q,\ell_q^{-1}}.
		\end{aligned}
	\end{equation}

	For $M,N\in\mathbb N_0$ with $1\leq M+N\leq N_{\mathrm{ind}}$, the terminal density in the Newton step satisfies
	\begin{equation}\label{eq:Newton-terminal-density-material-estimate}
	\|\DTL^M\rho_{q,J_*}\|_N\lesssim\delta_q\MM{M,N_R-3J_*,\tau_q^{-1},\ell_{t,q}^{-1}}\MM{N,N_{R,J_*}^{(0)}(M),\lambda_q,\ell_q^{-1}}.
	\end{equation}
	In the same range, the terminal momentum in the Newton step satisfies
	\begin{equation}\label{eq:Newton-terminal-momentum-material-estimate}
	\begin{aligned}
		\nrm{\DTL^Mm_{q,J_*}}_N
		&\lesssim\delta_q^{\frac12}
		\MM{M,N_R-3J_*,\tau_q^{-1},\ell_{t,q}^{-1}}
		\MM{N,N_{R,J_*}^{(0)}(M),\lambda_q,\ell_q^{-1}}.
	\end{aligned}
	\end{equation}

	For $1\leq j\leq J_*$, $0\leq M\leq N_{\mathrm{ind}}$, and $M+N\leq K_j^{(0)}+1$, the gluing stress $S_{q,j}$ satisfies
	\begin{equation}\label{eq:Newton-accumulated-gluing-stress-background}
		\begin{aligned}
			&\nrm{\DTL^MS_{q,j}}_{N+\alpha}\\
			&\quad\lesssim A_{q,0}
			(\tau_q\ell_q^{-1})^\alpha\MM{M,N_R-3j+1,\tau_q^{-1},\ell_{t,q}^{-1}}
			\Lat{N,N_{R,j-1}^{(2)}(M),\lambda_q,\ell_q^{-1}}
			\MM{N,N_{R,j-1}^{(2)}(M),\lambda_q,\ell_q^{-1}}.
		\end{aligned}
	\end{equation}
	The H{\"o}lder estimate \eqref{eq:Newton-accumulated-gluing-stress-background} also holds at each fixed exponent used in the construction.
	
	For $0\leq j<J_*$, $M\leq N_{\mathrm{ind}}-1$, and $M+N\leq K_j^{(0)}$, set $L=N_{R,j}^{(0)}(M+N)$. The two source mollification stresses satisfy
	\begin{equation}\label{eq:Newton-substep-mollification-remainder-estimate}
		\begin{aligned}
			\nrm{\DTL^MP_{q+1,j+1}^{\mathrm{mol,comp}}}_N+\nrm{\DTL^MP_{q+1,j+1}^{\mathrm{mol,inc}}}_N
			&\lesssim A_{q,j}\ell_q^{-\alpha}(\ell_{t,q}\tau_q^{-1})^L
			\MM{M,N_R-3j,\tau_q^{-1},\ell_{t,q}^{-1}}\ell_q^{-N}.
		\end{aligned}
	\end{equation}
	In the same range, one also has
	\begin{equation}\label{eq:Newton-substep-mollification-Newton-scales}
		\begin{aligned}
			\nrm{\DTL^MP_{q+1,j+1}^{\mathrm{mol,comp}}}_N+\nrm{\DTL^MP_{q+1,j+1}^{\mathrm{mol,inc}}}_N
			&\lesssim A_{q,j}\ell_q^{-\alpha}(\ell_q\lambda_q)^{\frac{N_R-3j}{2}}\ell_{t,q}^{-M}\ell_q^{-N}.
		\end{aligned}
	\end{equation}
	For $M+N\leq N_R$, the stress $P_{q+1,j}$ satisfies
	\begin{equation}\label{eq:Newton-substep-perturbative-tensor-detailed-estimate}
		\begin{aligned}
			\sup_{0\leq j\leq J_*}\nrm{\DTL^MP_{q+1,j}}_N
			&\lesssim{}A_{q,0}^2\delta_q^{-1}\lambda_{q+1}^{-6\alpha}\ell_{t,q}^{-M}\ell_q^{-N}.
		\end{aligned}
	\end{equation}
	The implicit constants may depend on $\beta,b,\alpha,N_{\mathrm{max}},\rho^*,\mathfrak M,C_*,\alpha_r$.
	They are independent of $q$, $a$, and $j$.

	Finally, the temporal support inclusion below leaves room for the support enlargement in the Nash step. If \eqref{ind:R-support} holds, the perturbations and stresses at the end of the Newton step satisfy
\begin{equation}\label{eq:Newton-terminal-temporal-support-next-scale}
	\begin{aligned}
		&\supp_t(\rho_{q,J_*}-\rho_q)\cup\supp_t(m_{q,J_*}-m_q)\cup\supp_tR_{q,J_*}\cup\supp_tS_{q,J_*}\cup\supp_tP_{q+1,J_*}\\
		&\qquad\subset[-2+\tfrac12\Lambda_q^{-1},-1-\tfrac12\Lambda_q^{-1}]
		\cup[1+\tfrac12\Lambda_q^{-1},2-\tfrac12\Lambda_q^{-1}].
	\end{aligned}
\end{equation}
\end{prop}

\begin{proof}
	We first prove \eqref{eq:Newton-terminal-density-increment-material}--\eqref{eq:Newton-terminal-momentum-material-estimate}.
	We then estimate the gluing stresses $S_{q,j}$ and the source mollification stresses $P_{q+1,j+1}^{\mathrm{mol,comp}}$ and $P_{q+1,j+1}^{\mathrm{mol,inc}}$.
	We use the remainders in \eqref{eq:Newton-convective-flux-remainder-identity}--\eqref{eq:Newton-pressure-flux-remainder-identity} to estimate $P_{q+1,j}$.
	Finally, we prove \eqref{eq:Newton-terminal-temporal-support-next-scale}.
	The source and perturbation estimates apply by \eqref{eq:Newton-interior-stress-estimate}.
	For $0\leq j\leq J_*$, the state differences satisfy
	\begin{equation}\label{eq:Newton-state-difference-decomposition}
	\begin{aligned}
		\rho_{q,j}-\orho_q
		=(\rho_q-\orho_q)
		+\sum_{\nu=0}^{j-1}\trho_{q,\nu+1}\comp,\qquad
		m_{q,j}-\om_q
		=(m_q-\om_q)
		+\sum_{\nu=0}^{j-1}
		\left(\tm_{q,\nu+1}\comp+\tm_{q,\nu+1}\inc\right).
	\end{aligned}
	\end{equation}
	For \eqref{eq:Newton-terminal-density-increment-material}, use the density terms in the proof of \eqref{eq:Newton-accumulated-state-increments-integer} with $j=J_*$ and $\alpha_2=\alpha/4$.
	For $M\leq N_{\mathrm{ind}}$ and $M+N\leq K_{J_*}^{(0)}$, \eqref{eq:Newton-total-compressible-density-estimate} applies at every $\nu<J_*$, since $K_{J_*}^{(0)}\leq K_\nu^{(2)}$.
	Use $\Xi_{\alpha_2}^2\leq\ell_q^{-2\alpha_2}=\ell_q^{-\alpha/2}$ and lower the material and spatial thresholds to $N_R-3J_*$ and $N_{R,J_*}^{(0)}(M)$, respectively. We obtain
	$$
	\begin{aligned}
	    \nrm{\DTL^M(\rho_{q,J_*}-\rho_q)}_N
	    &\lesssim\left(\sum_{\nu=0}^{J_*-1}A_{q,\nu}\right)\ell_q^{-\alpha/2}
	    \MM{M,N_R-3J_*,\tau_q^{-1},\ell_{t,q}^{-1}}
	    \MM{N,N_{R,J_*}^{(0)}(M),\lambda_q,\ell_q^{-1}}\\
	    &\lesssim A_{q,0}\ell_q^{-\alpha/2}
	    \MM{M,N_R-3J_*,\tau_q^{-1},\ell_{t,q}^{-1}}
	    \MM{N,N_{R,J_*}^{(0)}(M),\lambda_q,\ell_q^{-1}}.
	\end{aligned}
	$$
	The last inequality uses $\sum_{\nu=0}^{J_*-1}A_{q,\nu}\leq J_*A_{q,0}$ and $J_*\leq N_{\mathrm{max}}$.

	To prove \eqref{eq:Newton-terminal-density-difference-two-scale}, take $j=J_*$ in \eqref{eq:Newton-state-difference-decomposition}. For $M+N\leq N_{\mathrm{ind}}$,
	$$
	\begin{aligned}
	\|\DTL^M(\rho_{q,J_*}-\orho_q)\|_N
	&\leq\|\DTL^M(\rho_q-\orho_q)\|_N+\sum_{\nu=0}^{J_*-1}\|\DTL^M\trho_{q,\nu+1}\comp\|_N\\
	&\lesssim\Bigl[(\ell_q\lambda_q)^{N_{\mathrm{ind}}-\min\{M+N,N_R-3J_*\}}\delta_q
	+A_{q,0}\ell_q^{-2\alpha_2}\Bigr]\MM{M,N_R-3J_*,\tau_q^{-1},\ell_{t,q}^{-1}}
	\MM{N,N_{R,J_*}^{(0)}(M),\lambda_q,\ell_q^{-1}}.
	\end{aligned}
	$$
	For the perturbation sum, we used \eqref{eq:Newton-terminal-density-increment-material}, with $\ell_q^{-\alpha/2}=\ell_q^{-2\alpha_2}$.
	For $\rho_q-\orho_q$, we used \eqref{est:uniform-rho-minus-orho} as in the proof of \eqref{eq:Newton-initial-state-difference-material-estimate} with $K=N_R-3J_*$.

	The definition of $A_{q,0}$, \eqref{def:parameter-mollification}, and $\lambda_{q+1}\sim\lambda_q^b$ from \eqref{def:lambda-delta} give the following amplitude comparisons.
	We use $b<5/3$ and $N_{\mathrm{ind}}-N_R+3J_*\geq7J_*+12>16+8\beta$.
	After increasing $a_{\mathrm{Newton}}$, we have
	\begin{equation}\label{eq:Newton-terminal-velocity-difference-amplitudes}
	\begin{aligned}
		\frac{A_{q,0}\delta_q^{-\frac12}\lambda_{q+1}^{-3\alpha}}{\delta_{q+1}^{\frac12}}
		\lesssim\lambda_q^{-(8+3b)\alpha-(b-1)\beta}\leq\lambda_q^{-11\alpha},\quad
		\frac{(\ell_q\lambda_q)^{N_{\mathrm{ind}}-N_R+3J_*}\delta_q^{\frac12}}{\delta_{q+1}^{\frac12}}
		=(\lambda_q/\lambda_{q+1})^{\frac{N_{\mathrm{ind}}-N_R+3J_*}{4}-\beta}\leq(\lambda_q/\lambda_{q+1})^4.
	\end{aligned}
	\end{equation}
	For the density difference, we also have
	$
	A_{q,0}\ell_q^{-2\alpha_2}
	\lesssim_{C_*}\lambda_q^{-\frac{61-b}{8}\alpha}\delta_{q+1}\leq\lambda_q^{-7\alpha}\delta_{q+1}
	\leq\delta_{q+1}.
	$
	For $M+N\leq N_{\mathrm{ind}}$, the same relations give
	$$
	\frac{(\ell_q\lambda_q)^{N_{\mathrm{ind}}-\min\{M+N,N_R-3J_*\}}\delta_q}{\delta_{q+1}}
	\leq(\lambda_q/\lambda_{q+1})^{\frac{N_{\mathrm{ind}}-N_R+3J_*}{4}-2\beta}\leq(\lambda_q/\lambda_{q+1})^4.
	$$
	Substituting these estimates into \eqref{eq:Newton-terminal-density-difference-two-scale} gives, for $M+N\leq N_{\mathrm{ind}}$,
	$$
	\|\DTL^M(\rho_{q,J_*}-\orho_q)\|_N
	\lesssim\delta_{q+1}\MM{M,N_R-3J_*,\tau_q^{-1},\ell_{t,q}^{-1}}
	\MM{N,N_{R,J_*}^{(0)}(M),\lambda_q,\ell_q^{-1}}.
	$$
	This estimate and \eqref{est:uniform-orho} prove \eqref{eq:Newton-terminal-density-material-estimate} for $1\leq M+N\leq N_{\mathrm{ind}}$.
	For $M=0<N$, the $C^0$ term from \eqref{eq:Newton-state-zero-order-bounds} is covered by $\delta_q\lambda_q\geq1$, which follows from \eqref{def:lambda-delta}.
	
	For $M+N\leq N_{\mathrm{ind}}$, \eqref{eq:Newton-state-difference-decomposition} with $j=J_*$ gives
	$$
	\begin{aligned}
	\|\DTL^M(m_{q,J_*}-\om_q)\|_N
	&\leq\|\DTL^M(m_q-\om_q)\|_N+\sum_{\nu=0}^{J_*-1}\left(
	\|\DTL^M\tm_{q,\nu+1}\comp\|_N
	+\|\DTL^M\tm_{q,\nu+1}\inc\|_N\right)\\
	&\lesssim\Bigl[(\ell_q\lambda_q)^{N_{\mathrm{ind}}-\min\{M+N,N_R-3J_*\}}\delta_q^{\frac12}
	\hspace{-2pt}+\hspace{-2pt}A_{q,0}\delta_q^{-\frac12}\lambda_{q+1}^{-3\alpha}\Bigr]\MM{M,N_R-3J_*,\tau_q^{-1},\ell_{t,q}^{-1}}
	\MM{N,N_{R,J_*}^{(0)}(M),\lambda_q,\ell_q^{-1}}.
	\end{aligned}
	$$
	For the perturbation sum, we used \eqref{eq:Newton-accumulated-state-increments-integer} with $j=J_*$.
	For $m_q-\om_q$, we used \eqref{eq:Newton-initial-state-difference-material-estimate} with $K=N_R-3J_*$. By \eqref{est:uniform-om} and \eqref{eq:Newton-terminal-velocity-difference-amplitudes}, the preceding estimate for $m_{q,J_*}-\om_q$ proves \eqref{eq:Newton-terminal-momentum-material-estimate}.
	
	For $u_{q,J_*}-\overline u_q$, write
	$u_{q,J_*}-\overline u_q=\rho_{q,J_*}^{-1}\bigl(m_{q,J_*}-\om_q-\overline u_q(\rho_{q,J_*}-\orho_q)\bigr)$.
	At $M=N=0$, use \eqref{est:uniform-zero-velocity} and \eqref{eq:Newton-state-zero-order-bounds}.
	At positive total orders, use \eqref{eq:Newton-terminal-density-material-estimate} for $\rho_{q,J_*}$.
	Apply \eqref{eq:Newton-quadratic-weight-product} and \eqref{eq:faa-di-bruno} to $\rho_{q,J_*}^{-1}$, and use \eqref{est:uniform-om-over-orho} for $\overline u_q$.
	For $M+N\leq N_{\mathrm{ind}}$, we have
	\begin{align*}
	\|\DTL^M(\rho_{q,J_*}^{-1})\|_N+\|\DTL^M\overline u_q\|_N\lesssim\MM{M,N_R-3J_*,\tau_q^{-1},\ell_{t,q}^{-1}}\MM{N,N_{R,J_*}^{(0)}(M),\lambda_q,\ell_q^{-1}}.
	\end{align*}
	Using these coefficient estimates, we obtain, in the same range,
	\begin{align*}
	&\|\DTL^M(u_{q,J_*}-\overline u_q)\|_N\\
	&\quad\lesssim\sum_{\substack{M_1+M_2=M\\N_1+N_2=N}}
	\|\DTL^{M_1}(\rho_{q,J_*}^{-1})\|_{N_1}
	\|\DTL^{M_2}(m_{q,J_*}-\om_q)\|_{N_2}+\sum_{\substack{M_1+M_2+M_3=M\\N_1+N_2+N_3=N}}
	\|\DTL^{M_1}(\rho_{q,J_*}^{-1})\|_{N_1}
	\|\DTL^{M_2}\overline u_q\|_{N_2}
	\|\DTL^{M_3}(\rho_{q,J_*}-\orho_q)\|_{N_3}\\
	&\quad\lesssim\Bigl[A_{q,0}\delta_q^{-\frac12}\lambda_{q+1}^{-3\alpha}
	+(\ell_q\lambda_q)^{N_{\mathrm{ind}}-N_R+3J_*}\delta_q^{\frac12}\Bigr]\MM{M,N_R-3J_*,\tau_q^{-1},\ell_{t,q}^{-1}}
	\MM{N,N_{R,J_*}^{(0)}(M),\lambda_q,\ell_q^{-1}}.
	\end{align*}
	For the last inequality, we used \eqref{eq:Newton-quadratic-weight-product}, \eqref{eq:Newton-quadratic-Holder-amplitude-factors}, \eqref{eq:Newton-terminal-density-difference-two-scale},
	and the estimate for $m_{q,J_*}-\om_q$ above.
	Estimate~\eqref{eq:Newton-terminal-velocity-difference-two-scale} follows.

	We next prove \eqref{eq:Newton-accumulated-gluing-stress-background}.
	Fix $1\leq j\leq J_*$, $0\leq M\leq N_{\mathrm{ind}}$, and $M+N\leq K_j^{(0)}+1$.
	Apply \eqref{eq:Newton-gluing-stress-estimate} with $j$ replaced by $\nu-1$ for each $1\leq\nu\leq j$. The derivative ranges and thresholds satisfy
	$
		K_{\nu-1}^{(2)}=K_\nu^{(0)}+1\geq K_j^{(0)}+1,\
		N_R-3\nu+1\geq N_R-3j+1,\
		N_{R,\nu-1}^{(2)}(M)\geq N_{R,j-1}^{(2)}(M)=\max\{N_R-3j+1-M,0\}.
	$
	With these common thresholds, we have
	\begin{align*}
		\nrm{\DTL^MS_{q,j}}_{N+\alpha}
		\lesssim\sum_{\nu=1}^{j}A_{q,\nu-1}
		(\tau_q\ell_q^{-1})^\alpha\MM{M,N_R-3j+1,\tau_q^{-1},\ell_{t,q}^{-1}}\Lat{N,N_{R,j-1}^{(2)}(M),\lambda_q,\ell_q^{-1}}
		\MM{N,N_{R,j-1}^{(2)}(M),\lambda_q,\ell_q^{-1}}.
	\end{align*}
	Since $\sum_{\nu=1}^{j}A_{q,\nu-1}\leq jA_{q,0}\leq J_*A_{q,0}$,
	we obtain \eqref{eq:Newton-accumulated-gluing-stress-background}.

	We next estimate $P_{q+1,j+1}^{\mathrm{mol,comp}}$ and $P_{q+1,j+1}^{\mathrm{mol,inc}}$, defined in \eqref{def:Newton-substep-flow-mollification}. For $0\leq j<J_*$, $M\leq N_{\mathrm{ind}}-1$, and $M+N\leq K_j^{(0)}$, set $L:=N_{R,j}^{(0)}(M+N)$. Since
	$\sum_u\theta_u^2=1$ and $\sum_\kappa\chi_\kappa^2=1$, substituting \eqref{def:Newton-sharp-acoustic-input} gives
	$$
		P_{q+1,j+1}^{\mathrm{mol,comp}}
		=\sum_u\theta_u^2J_u\Bigl[J_u^{-1}Q[R_{q,j}]-\eta_{\ell_{t,q}}^{\mathrm t}\ast_{X_{\overline u_q}}(J_u^{-1}Q[R_{q,j}])\Bigr]\Id.
	$$
	For each fixed H{\"o}lder exponent $\alpha_r\leq\alpha_2$ and integers $0\leq m\leq M$ and $n\geq0$ with $m+n+L\leq K_j^{(0)}$, we have
	\begin{equation}\label{eq:Newton-source-mollification-differences}
		\begin{aligned}
			\nrm{\DTL^mP_{q+1,j+1}^{\mathrm{mol,comp}}}_{n+\alpha_r}
			+\nrm{\DTL^m(W[R_{q,j}]-W_{q,j}^{\mathrm{mol}})}_{n+\alpha_r}\lesssim A_{q,j}(\ell_{t,q}\tau_q^{-1})^L
			\MM{m,N_R-3j,\tau_q^{-1},\ell_{t,q}^{-1}}\ell_q^{-n-2\alpha_r}.
		\end{aligned}
	\end{equation}
	For $L\geq1$, we applied \eqref{eq:initial-Newton-flow-mollification-remainder-identity} to $J_u^{-1}Q[R_{q,j}]$ and $W[R_{q,j}]$, using \eqref{eq:Newton-compressible-input-estimate} and \eqref{eq:Newton-Jacobian-weighted-source-estimate}. Their derivative ranges cover $m+L\leq M+L\leq N_R-3j$ and $m+n+L\leq K_j^{(0)}$. We estimated $\theta_u^2J_u$ by \eqref{eq:time-cutoff-derivative-estimate} and \eqref{eq:Newton-Jacobian-material-estimates}. Since at most two indices contribute at each time, \eqref{eq:Newton-source-mollification-differences} follows. For $L=0$, we used \eqref{eq:Newton-compressible-input-estimate}--\eqref{eq:Newton-mollified-source-estimate}.
	
	For $L\geq1$, estimate $P_{q+1,j+1}^{\mathrm{mol,inc}}$ by applying Lemma~\ref{lem:direct-symmetric-curl-material-estimate} to $W[R_{q,j}]-W_{q,j}^{\mathrm{mol}}$ with H{\"o}lder exponent $\alpha_2$, $s=2$, and parameters 
	$$
	\begin{aligned}
		K^\sharp&=K_j^{(0)}-L+1,&& K^*=M+1,&& K_G=0,&& M_*=M,&& N_t=N_R-3j,&&\mu_G=\tau_q^{-1},\\
		\lambda_\upsilon&=\tilde\lambda_\upsilon=\ell_q^{-1},&& \lambda_G=\tilde\lambda_G=\ell_q^{-1},&& \mu_\upsilon=\tilde\mu_\upsilon=\Lambda_q,&&
		A_\upsilon=C\Lambda_q\ell_q,&& A_G=CA_{q,j}(\ell_{t,q}\tau_q^{-1})^L,&&\tilde\mu_G=\ell_{t,q}^{-1}.
	\end{aligned}
	$$
	The parameters satisfy $M+N\leq K^\sharp-1$, $K^\sharp\geq K^*+2$, $K^*<N_{\mathrm{ind}}$, and $M+2=K^*+1$, as required for the estimate with H{\"o}lder exponent $\alpha_2$. The coefficient hypotheses follow from \eqref{est:uniform-grad-om-over-orho} and \eqref{est:holder-grad-om-over-orho}. The source hypotheses follow from \eqref{eq:Newton-source-mollification-differences}, with $C$ chosen to cover the implicit constants in these three estimates over the required derivative ranges. For sufficiently large $a$, $A_\upsilon\ell_q^{-1}=C\Lambda_q\leq\tau_q^{-1}\leq\ell_{t,q}^{-1}$. The lemma gives
	$$
	\nrm{\DTL^MP_{q+1,j+1}^{\mathrm{mol,inc}}}_{N+\alpha_2}
	\lesssim A_{q,j}(\ell_{t,q}\tau_q^{-1})^L
	\MM{M,N_R-3j,\tau_q^{-1},\ell_{t,q}^{-1}}\ell_q^{-N-2\alpha_2}.
	$$
	Combine this with \eqref{eq:Newton-source-mollification-differences} with H{\"o}lder exponent $\alpha_2$. Passing to $C^N$ proves \eqref{eq:Newton-substep-mollification-remainder-estimate} for $L\geq1$.
	
	For $M\leq N_{\mathrm{ind}}-1$ and $M+N\leq K_j^{(0)}$, we have
	\begin{equation}\label{eq:Newton-source-mollification-two-scale-estimate}
	\begin{aligned}
		\nrm{\DTL^MP_{q+1,j+1}^{\mathrm{mol,comp}}}_N+\nrm{\DTL^MP_{q+1,j+1}^{\mathrm{mol,inc}}}_N\lesssim A_{q,j}\ell_q^{-\alpha}\MM{M,N_R-3j,\tau_q^{-1},\ell_{t,q}^{-1}}\MM{N,N_{R,j}^{(0)}(M),\lambda_q,\ell_q^{-1}}.
	\end{aligned}
	\end{equation}
	For $P_{q+1,j+1}^{\mathrm{mol,comp}}$, we combined the cutoff estimate \eqref{eq:time-cutoff-derivative-estimate}, \eqref{eq:Newton-compressible-input-estimate}--\eqref{eq:Newton-mollified-source-estimate}, and the Jacobian estimates \eqref{eq:Newton-Jacobian-material-estimates}. For $P_{q+1,j+1}^{\mathrm{mol,inc}}$, we applied Lemma~\ref{lem:direct-symmetric-curl-material-estimate} to $W[R_{q,j}]-W_{q,j}^{\mathrm{mol}}$ with H{\"o}lder exponent $\alpha_2$, $s=2$, and parameters 
	$$
	\begin{aligned}
		K^\sharp&=K_j^{(0)}+1,&& K^*=N_{\mathrm{ind}},&& K_G=N_t=N_R-3j,&& M_*=N_{\mathrm{ind}}-1,\\
		\lambda_\upsilon&=\lambda_G=\lambda_q,&& \tilde\lambda_\upsilon=\tilde\lambda_G=\ell_q^{-1},&& \mu_\upsilon=\tilde\mu_\upsilon=\Lambda_q,  &&
		A_\upsilon=C\delta_q^{\frac12},&& A_G=CA_{q,j},&& \mu_G=\tau_q^{-1},\quad \tilde\mu_G=\ell_{t,q}^{-1}.
	\end{aligned}
	$$
	We used \eqref{est:uniform-grad-om-over-orho} and \eqref{est:holder-grad-om-over-orho} for the coefficients, and \eqref{eq:Newton-compressible-input-estimate}--\eqref{eq:Newton-mollified-source-estimate} for the source. Here $C$ covers the coefficient and source constants. The derivative ranges satisfy $M+N\leq K_j^{(0)}=K^\sharp-1$, $K^\sharp\geq K^*+2$, and $M+2\leq N_{\mathrm{ind}}+1=K^*+1$. With the material scales from \eqref{eq:Newton-acoustic-combined-material-scales}, the lemma proves \eqref{eq:Newton-source-mollification-two-scale-estimate} for $P_{q+1,j+1}^{\mathrm{mol,inc}}$. At $L=0$, its spatial weight is at most $\ell_q^{-N}$, so \eqref{eq:Newton-substep-mollification-remainder-estimate} follows as well.
	Estimates~\eqref{eq:Newton-substep-mollification-remainder-estimate} and \eqref{eq:Newton-source-mollification-two-scale-estimate} imply
	$$
	\begin{aligned}
		\nrm{\DTL^MP_{q+1,j+1}^{\mathrm{mol,comp}}}_N+\nrm{\DTL^MP_{q+1,j+1}^{\mathrm{mol,inc}}}_N
		&\lesssim A_{q,j}\ell_q^{-\alpha}\ell_{t,q}^{-M}\ell_q^{-N}
		\begin{cases}
			(\ell_{t,q}\tau_q^{-1})^{N_R-3j-M-N},&2(M+N)\leq N_R-3j,\\
			(\ell_q\lambda_q)^{\min\{M+N,N_R-3j\}},&2(M+N)>N_R-3j,
		\end{cases}\\
	&\leq A_{q,j}\ell_q^{-\alpha}(\ell_q\lambda_q)^{\frac{N_R-3j}{2}}\ell_{t,q}^{-M}\ell_q^{-N}.
	\end{aligned}
	$$
	The first case uses \eqref{eq:Newton-substep-mollification-remainder-estimate} with $L=N_R-3j-M-N$.
	For the second case, we applied \eqref{eq:Newton-two-scale-weight-comparison} with $K=N_R-3j$ to \eqref{eq:Newton-source-mollification-two-scale-estimate}.
	In each case, the exponent is at least $(N_R-3j)/2$. Using $\ell_{t,q}\tau_q^{-1}\leq\ell_q\lambda_q\leq1$ from \eqref{eq:material-mollification-spatial-comparisons} and \eqref{eq:Newton-scale-separation} gives \eqref{eq:Newton-substep-mollification-Newton-scales}.
	
	To prove \eqref{eq:Newton-substep-perturbative-tensor-detailed-estimate}, we estimate the initial value $P_{q+1,0}$ and the differences $P_{q+1,j+1}-P_{q+1,j}$, then sum over the substeps.
	Begin with $P_{q+1,0}=R_q-R_{q,0}$.
	For $M+N\leq N_R$, apply \eqref{est:uniform-R-minus-Rq0} to $R_q-R_{q,0}^{\mathrm{sp}}$. For $R_{q,0}^{\mathrm{sp}}-R_{q,0}$, use the $C^N$ version of \eqref{eq:initial-Newton-flow-mollification-remainder-estimate} with $L=N_R-M-N$ when $M+N<N_R$. At $M+N=N_R$, use the $C^N$ version of \eqref{eq:initial-Newton-flow-mollification-endpoint-estimate}. We have
	\begin{equation}\label{eq:Newton-initial-perturbative-tensor-material-estimate}
	\begin{aligned}
		\nrm{\DTL^MP_{q+1,0}}_N
		&\leq\nrm{\DTL^M(R_q-R_{q,0}^{\mathrm{sp}})}_N
		+\nrm{\DTL^M(R_{q,0}^{\mathrm{sp}}-R_{q,0})}_N\\
		&\lesssim A_{q,0}\Bigl[(\ell_q\lambda_q)^{N_R}
		\left(\delta_q^{\frac12}\ell_{t,q}\ell_q^{-1}\right)^M
		+(\ell_{t,q}\Lambda_q)^{N_R-N}(\ell_q\lambda_q)^N\Bigr]\ell_{t,q}^{-M}\ell_q^{-N}\\
		&\lesssim A_{q,0}(\ell_q\lambda_q)^{N_R}\ell_{t,q}^{-M}\ell_q^{-N}.
	\end{aligned}
	\end{equation}
	Here the last inequality uses $\delta_q^{\frac12}\ell_{t,q}\ell_q^{-1}\leq1$ and
	$\ell_{t,q}\Lambda_q=(\ell_q\lambda_q)^2\leq\ell_q\lambda_q$.
	
	For $0\leq j<J_*$, the convective flux remainder satisfies
	\begin{equation}\label{eq:Newton-convective-flux-remainder-identity}
	\begin{aligned}
		&\frac{m_{q,j+1}\otimes m_{q,j+1}}{\rho_{q,j+1}}
		-\frac{m_{q,j}\otimes m_{q,j}}{\rho_{q,j}}-\overline u_q\otimes(\tm_{q,j+1}\comp+\tm_{q,j+1}\inc)
		-(\tm_{q,j+1}\comp+\tm_{q,j+1}\inc)\otimes\overline u_q
		+\trho_{q,j+1}\comp\overline u_q\otimes\overline u_q\\
		=&(u_{q,j}-\overline u_q)\otimes(\tm_{q,j+1}\comp+\tm_{q,j+1}\inc)
		+(\tm_{q,j+1}\comp+\tm_{q,j+1}\inc)\otimes(u_{q,j}-\overline u_q)\\
		&-\trho_{q,j+1}\comp\left(u_{q,j}\otimes u_{q,j}-\overline u_q\otimes\overline u_q\right)\\
		&+\frac{(\tm_{q,j+1}\comp+\tm_{q,j+1}\inc-\trho_{q,j+1}\comp u_{q,j})\otimes(\tm_{q,j+1}\comp+\tm_{q,j+1}\inc-\trho_{q,j+1}\comp u_{q,j})}{\rho_{q,j}+\trho_{q,j+1}\comp},
	\end{aligned}
	\end{equation}
	and Taylor's formula gives the pressure flux remainder
	\begin{equation}\label{eq:Newton-pressure-flux-remainder-identity}
	\begin{aligned}
		&p(\rho_{q,j}+\trho_{q,j+1}\comp)-p(\rho_{q,j})-p'(\orho_q)\trho_{q,j+1}\comp\\
		&\quad=\left(p'(\rho_{q,j})-p'(\orho_q)\right)\trho_{q,j+1}\comp
		+(\trho_{q,j+1}\comp)^2\int_0^1(1-\sigma)p''(\rho_{q,j}+\sigma\trho_{q,j+1}\comp)\,\rd\sigma.
	\end{aligned}
	\end{equation}
	Use the identity
	$u_{q,j}-\overline u_q=\rho_{q,j}^{-1}\bigl(m_{q,j}-\om_q-\overline u_q(\rho_{q,j}-\orho_q)\bigr)$.
	Each flux remainder is a sum of products with smooth coefficients. Each product contains a new perturbation and either another new perturbation or one of the differences $\rho_{q,j}-\orho_q$ and $m_{q,j}-\om_q$.

	Fix $0\leq j<J_*$ and $M,N\in\mathbb N_0$ with $M+N\leq N_{\mathrm{ind}}$. We use \eqref{eq:Newton-quadratic-input-material-estimates} for each $0\leq\nu\leq j$.
	For $j\geq1$, substitute \eqref{eq:Newton-accumulated-state-increments-integer} and \eqref{eq:Newton-initial-state-difference-material-estimate} with $K=N_R-3j-2$ into \eqref{eq:Newton-state-difference-decomposition}, lowering the material threshold for the perturbations from $N_R-3j$ to $N_R-3j-2$.
	For $j=0$, the perturbation sum is empty. We therefore obtain, for $M+N\leq N_{\mathrm{ind}}$,
	\begin{equation}\label{eq:Newton-accumulated-state-difference-material-estimate}
	\begin{aligned}
		&\nrm{\DTL^M(\rho_{q,j}-\orho_q)}_N+\nrm{\DTL^M(m_{q,j}-\om_q)}_N\\
		&\quad\lesssim\left[(\ell_q\lambda_q)^{N_{\mathrm{ind}}-\min\{M+N,N_R-3j-2\}}\delta_q^{\frac12}+A_{q,0}\delta_q^{-\frac12}\lambda_{q+1}^{-3\alpha}\right]
		\MM{M,N_R-3j-2,\tau_q^{-1},\ell_{t,q}^{-1}}\MM{N,N_{R,j}^{(2)}(M),\lambda_q,\ell_q^{-1}}\\
		&\quad\lesssim\delta_q^{\frac12}\MM{M,N_R-3j-2,\tau_q^{-1},\ell_{t,q}^{-1}}\MM{N,N_{R,j}^{(2)}(M),\lambda_q,\ell_q^{-1}}.
	\end{aligned}
	\end{equation}
	Since $N_R-3j-2\leq N_{\mathrm{ind}}$, the first term in brackets is at most $\delta_q^{\frac12}$.
	Using $A_{q,0}=C_*\lambda_q^{-8\alpha}\delta_{q+1}$, we also have $A_{q,0}\delta_q^{-\frac12}\lambda_{q+1}^{-3\alpha}\leq\delta_q^{\frac12}$ for sufficiently large $a$.

	To estimate the coefficients in the two flux remainders, define the following interpolation segments for $0\leq\sigma\leq1$
	\begin{equation}\label{eq:Newton-flux-interpolation-segments}
	\begin{aligned}
		(\rho_\sigma^{(1)},m_\sigma^{(1)})&=\bigl(\orho_q+\sigma(\rho_{q,j}-\orho_q),\,
		\om_q+\sigma(m_{q,j}-\om_q)\bigr),\\
		(\rho_\sigma^{(2)},m_\sigma^{(2)})&=\bigl(\rho_{q,j}+\sigma\trho_{q,j+1}\comp,\,
		m_{q,j}+\sigma(\tm_{q,j+1}\comp+\tm_{q,j+1}\inc)\bigr).
	\end{aligned}
	\end{equation}
	Proposition~\ref{prop:uniform-mollification-estimates} and \eqref{eq:Newton-state-zero-order-bounds} give $\rho_\sigma^{(s)}\in[\rho^*/2,3\rho^*/2]$.
	Both momentum segments are convex combinations of their endpoints. By \eqref{est:uniform-zero-velocity} and \eqref{eq:Newton-state-zero-order-bounds}, we have $\nrm{m_\sigma^{(s)}}_0\leq\mathfrak M\left(1-\frac12\delta_q^{\frac12}\right).$
	The definition of $A_{q,0}$, $A_{q,j}\leq A_{q,0}$, and $b<1+2\beta$ give
	$$
	\begin{aligned}
		A_{q,j}\ell_q^{-2\alpha_2}\left[1+\tau_q(\tau_q\ell_q^{-1})^{\alpha_2}
		\frac{\MM{N+1,N_{R,j}^{(2)}(M),\lambda_q,\ell_q^{-1}}}{\MM{N,N_{R,j}^{(2)}(M),\lambda_q,\ell_q^{-1}}}\right]&\leq2A_{q,j}(\tau_q\ell_q^{-1})^{1+\alpha_2}\ell_q^{-2\alpha_2}\\
		&\lesssim\lambda_q^{-(2b-1)\beta+\frac{b-1}{4}}\leq\lambda_q^{-\frac\beta2}\leq1.
	\end{aligned}
	$$
	Since $b<1+2\beta$, take $a_{\mathrm{Newton}}$ sufficiently large to absorb the implicit constant.
	Combining this with \eqref{est:uniform-orho}--\eqref{est:uniform-om}, \eqref{eq:Newton-quadratic-input-material-estimates}, and \eqref{eq:Newton-accumulated-state-difference-material-estimate}, we obtain, uniformly for $s\in\{1,2\}$ and $\sigma\in[0,1]$,
	\begin{equation}\label{eq:Newton-interpolated-state-material-estimate}
		\nrm{\DTL^M\rho_\sigma^{(s)}}_N+\nrm{\DTL^M m_\sigma^{(s)}}_N
		\lesssim\MM{M,N_R-3j-2,\tau_q^{-1},\ell_{t,q}^{-1}}
		\MM{N,N_{R,j}^{(2)}(M),\lambda_q,\ell_q^{-1}},
		\qquad 1\leq M+N\leq N_{\mathrm{ind}}.
	\end{equation}

	For the second derivatives of $m\otimes m/\rho$ and $p(\rho)$ with respect to $(\rho,m)$ and $1\leq M+N\leq N_{\mathrm{ind}}$, we have 
	\begin{equation}\label{eq:Newton-flux-coefficient-material-estimate}
	\begin{aligned}
		&\nrm{\DTL^M\bigl((\rho_\sigma^{(s)})^{-1}\bigr)}_N
		+\Nrm{\DTL^M\left(\frac{m_\sigma^{(s)}}{(\rho_\sigma^{(s)})^2}\right)}_N
		+\Nrm{\DTL^M\left(\frac{m_\sigma^{(s)}\otimes m_\sigma^{(s)}}{(\rho_\sigma^{(s)})^3}\right)}_N
		+\nrm{\DTL^M p''(\rho_\sigma^{(s)})}_N\\
		&\quad\lesssim\sum_{k=1}^{M+N}
		\sum_{\substack{M_1+\cdots+M_k=M\\N_1+\cdots+N_k=N\\M_i+N_i\geq1}}
		\prod_{i=1}^{k}\Bigl(\nrm{\DTL^{M_i}\rho_\sigma^{(s)}}_{N_i}+\nrm{\DTL^{M_i}m_\sigma^{(s)}}_{N_i}\Bigr)\\
		&\quad\lesssim\sum_{k=1}^{M+N}
		\sum_{\substack{M_1+\cdots+M_k=M\\N_1+\cdots+N_k=N\\M_i+N_i\geq1}}
		\prod_{i=1}^{k}\Bigl(\MM{M_i,N_R-3j-2,\tau_q^{-1},\ell_{t,q}^{-1}}
		\MM{N_i,N_{R,j}^{(2)}(M_i),\lambda_q,\ell_q^{-1}}\Bigr)\\
		&\quad\lesssim
		\MM{M,N_R-3j-2,\tau_q^{-1},\ell_{t,q}^{-1}}
		\MM{N,N_{R,j}^{(2)}(M),\lambda_q,\ell_q^{-1}}.
	\end{aligned}
	\end{equation}
	Here the three inequalities use \eqref{eq:faa-di-bruno}, \eqref{eq:Newton-interpolated-state-material-estimate}, and \eqref{eq:Newton-quadratic-weight-product}, respectively. For $M=N=0$, the final estimate follows from $\rho_\sigma^{(s)}\geq\rho^*/2$ and the $C^0$ estimates on both segments.
	
	For the term containing $p'(\rho_{q,j})-p'(\orho_q)$ in \eqref{eq:Newton-pressure-flux-remainder-identity}, the mean value formula along the first segment in \eqref{eq:Newton-flux-interpolation-segments} gives
	\begin{equation}\label{eq:Newton-pressure-mixed-Leibniz-estimate}
	\begin{aligned}
		&\Nrm{\DTL^M\Bigl((p'(\rho_{q,j})-p'(\orho_q))\trho_{q,j+1}\comp\Bigr)}_N\\
		&\quad\lesssim\sum_{M_1+M_2+M_3=M}\int_0^1
		\Nrm{\DTL^{M_1}\Bigl[p''(\rho_\sigma^{(1)})\Bigr]
		\DTL^{M_2}(\rho_{q,j}-\orho_q)\DTL^{M_3}\trho_{q,j+1}\comp}_N\,\rd\sigma\\
		&\quad\lesssim\sum_{\substack{M_1+M_2+M_3=M\\N_1+N_2+N_3=N}}
		\sup_{0\leq\sigma\leq1}\Nrm{\DTL^{M_1}\Bigl[p''(\rho_\sigma^{(1)})\Bigr]}_{N_1}\nrm{\DTL^{M_2}(\rho_{q,j}-\orho_q)}_{N_2}\nrm{\DTL^{M_3}\trho_{q,j+1}\comp}_{N_3}.
	\end{aligned}
	\end{equation}
	For the quadratic term in \eqref{eq:Newton-pressure-flux-remainder-identity}, use the second segment in \eqref{eq:Newton-flux-interpolation-segments}.
	Estimate its derivatives as in \eqref{eq:Newton-pressure-mixed-Leibniz-estimate}, with $p''(\rho_\sigma^{(2)})$ and two factors of $\trho_{q,j+1}\comp$.

	For the convective flux $m\otimes m/\rho$, compare its first derivatives at $(\orho_q,\om_q)$ and $(\rho_{q,j},m_{q,j})$ using the mean value formula on the first segment in \eqref{eq:Newton-flux-interpolation-segments}. Taylor's formula on the second segment gives the quadratic remainder in \eqref{eq:Newton-convective-flux-remainder-identity}. Estimate the resulting products as in \eqref{eq:Newton-pressure-mixed-Leibniz-estimate}.
	The definition \eqref{def:Newton-perturbative-tensor} and estimate \eqref{eq:Newton-flux-coefficient-material-estimate} give, for $M+N\leq N_{\mathrm{ind}}$,
	\begin{equation}\label{eq:Newton-perturbative-tensor-increment-estimate}
	\begin{aligned}
		&\nrm{\DTL^M(P_{q+1,j+1}-P_{q+1,j}-P_{q+1,j+1}^{\mathrm{mol,comp}}-P_{q+1,j+1}^{\mathrm{mol,inc}})}_N\\
		&\quad\lesssim
		\sum_{\substack{M_1+M_2+M_3=M\\N_1+N_2+N_3=N}}
		\sup_{\substack{s\in\{1,2\}\\0\leq\sigma\leq1}}\left(
		\nrm{\DTL^{M_1}\bigl((\rho_\sigma^{(s)})^{-1}\bigr)}_{N_1}
		+\Nrm{\DTL^{M_1}\left(\frac{m_\sigma^{(s)}}{(\rho_\sigma^{(s)})^2}\right)}_{N_1}
		+\Nrm{\DTL^{M_1}\left(\frac{m_\sigma^{(s)}\otimes m_\sigma^{(s)}}{(\rho_\sigma^{(s)})^3}\right)}_{N_1}
		+\nrm{\DTL^{M_1}p''(\rho_\sigma^{(s)})}_{N_1}\right)\\
		&\qquad\cdot\Bigl(
		\nrm{\DTL^{M_2}(\rho_{q,j}-\orho_q)}_{N_2}
		+\nrm{\DTL^{M_2}(m_{q,j}-\om_q)}_{N_2}
		+\nrm{\DTL^{M_2}\trho_{q,j+1}\comp}_{N_2}
		+\nrm{\DTL^{M_2}\tm_{q,j+1}\comp}_{N_2}
		+\nrm{\DTL^{M_2}\tm_{q,j+1}\inc}_{N_2}\Bigr)\\
		&\qquad\cdot\Bigl(
		\nrm{\DTL^{M_3}\trho_{q,j+1}\comp}_{N_3}
		+\nrm{\DTL^{M_3}\tm_{q,j+1}\comp}_{N_3}
		+\nrm{\DTL^{M_3}\tm_{q,j+1}\inc}_{N_3}\Bigr)\\
		&\quad\lesssim
		\sum_{\substack{M_1+M_2+M_3=M\\N_1+N_2+N_3=N}}
		\MM{M_1,N_R-3j-2,\tau_q^{-1},\ell_{t,q}^{-1}}
		\MM{N_1,N_{R,j}^{(2)}(M_1),\lambda_q,\ell_q^{-1}}\\
		&\qquad\cdot\Bigl(
		\nrm{\DTL^{M_2}(\rho_{q,j}-\orho_q)}_{N_2}
		+\nrm{\DTL^{M_2}(m_{q,j}-\om_q)}_{N_2}
		+\nrm{\DTL^{M_2}\trho_{q,j+1}\comp}_{N_2}
		+\nrm{\DTL^{M_2}\tm_{q,j+1}\comp}_{N_2}
		+\nrm{\DTL^{M_2}\tm_{q,j+1}\inc}_{N_2}\Bigr)\\
		&\qquad\cdot\Bigl(
		\nrm{\DTL^{M_3}\trho_{q,j+1}\comp}_{N_3}
		+\nrm{\DTL^{M_3}\tm_{q,j+1}\comp}_{N_3}
		+\nrm{\DTL^{M_3}\tm_{q,j+1}\inc}_{N_3}\Bigr).
	\end{aligned}
	\end{equation}
	
	Use \eqref{eq:Newton-state-difference-decomposition} in \eqref{eq:Newton-perturbative-tensor-increment-estimate}.
	For $M,N\in\mathbb N_0$ with $M+N\leq N_{\mathrm{ind}}$, apply \eqref{eq:Newton-uniform-perturbation-material-estimate} to the two perturbation factors to obtain
	\begin{align}
		&\sum_{\nu=0}^{j}
		\sum_{\substack{M_1+M_2+M_3=M\\N_1+N_2+N_3=N}}
		\MM{M_1,N_R-3j-2,\tau_q^{-1},\ell_{t,q}^{-1}}
		\MM{N_1,N_{R,j}^{(2)}(M_1),\lambda_q,\ell_q^{-1}}
		\nrm{\DTL^{M_2}\tm_{q,j+1}\inc}_{N_2}
		\nrm{\DTL^{M_3}\tm_{q,\nu+1}\inc}_{N_3}\nonumber\\
		&\ \lesssim{}A_{q,j}\left(\sum_{\nu=0}^{j}A_{q,\nu}\right)
		\delta_q^{-1}\lambda_{q+1}^{-6\alpha}
		\sum_{\substack{M_1+M_2+M_3=M\\N_1+N_2+N_3=N}}\MM{M_1,N_R-3j-2,\tau_q^{-1},\ell_{t,q}^{-1}}
		\MM{N_1,N_{R,j}^{(2)}(M_1),\lambda_q,\ell_q^{-1}}
		\ell_{t,q}^{-(M_2+M_3)}\ell_q^{-(N_2+N_3)}\nonumber\\
		&\ \lesssim{}A_{q,j}\left(\sum_{\nu=0}^{j}A_{q,\nu}\right)
		\delta_q^{-1}\lambda_{q+1}^{-6\alpha}
		\sum_{\substack{M_1+M_2+M_3=M\\N_1+N_2+N_3=N}}
		\prod_{i=1}^{3}\ell_{t,q}^{-M_i}\ell_q^{-N_i}\lesssim A_{q,j}A_{q,0}\delta_q^{-1}\lambda_{q+1}^{-6\alpha}
		\ell_{t,q}^{-M}\ell_q^{-N}.	\label{eq:Newton-quadratic-product-material-estimate}
	\end{align}
	The term $\nu=j$ accounts for products of two new perturbations.  By \eqref{eq:Newton-uniform-perturbation-material-estimate}, the same estimate holds when either or both perturbations are replaced by density or compressible momentum perturbations.

	For products containing $m_q-\om_q$, use \eqref{eq:Newton-uniform-perturbation-material-estimate} and \eqref{eq:Newton-initial-state-difference-material-estimate} with $K=0$. Apply \eqref{eq:Newton-two-scale-weight-comparison} to the coefficient weight. We have
	\begin{align}
		&\sum_{\substack{M_1+M_2+M_3=M\\N_1+N_2+N_3=N}}
		\MM{M_1,N_R-3j-2,\tau_q^{-1},\ell_{t,q}^{-1}}\MM{N_1,N_{R,j}^{(2)}(M_1),\lambda_q,\ell_q^{-1}}
		\nrm{\DTL^{M_2}(m_q-\om_q)}_{N_2}\nrm{\DTL^{M_3}\tm_{q,j+1}\inc}_{N_3}\nonumber\\
		&\ \lesssim{}A_{q,j}\lambda_{q+1}^{-3\alpha}(\ell_q\lambda_q)^{N_{\mathrm{ind}}}
		\sum_{\substack{M_1+M_2+M_3=M\\N_1+N_2+N_3=N}}
		\MM{M_1,N_R-3j-2,\tau_q^{-1},\ell_{t,q}^{-1}}\MM{N_1,N_{R,j}^{(2)}(M_1),\lambda_q,\ell_q^{-1}}\ell_{t,q}^{-M_2}\ell_q^{-N_2}
		\ell_{t,q}^{-M_3}\ell_q^{-N_3}\nonumber\\
		&\ \lesssim{}A_{q,j}\lambda_{q+1}^{-3\alpha}(\ell_q\lambda_q)^{N_{\mathrm{ind}}}
		\sum_{\substack{M_1+M_2+M_3=M\\N_1+N_2+N_3=N}}
		\prod_{i=1}^{3}\ell_{t,q}^{-M_i}\ell_q^{-N_i}\lesssim A_{q,j}\lambda_{q+1}^{-3\alpha}(\ell_q\lambda_q)^{N_{\mathrm{ind}}}
		\ell_{t,q}^{-M}\ell_q^{-N}.\label{eq:Newton-initial-mixed-coefficient-material-estimate}
	\end{align}
	Estimate~\eqref{eq:Newton-initial-mixed-coefficient-material-estimate} also holds with $\rho_q-\orho_q$ in place of $m_q-\om_q$. One may also replace $\tm_{q,j+1}\inc$ by $\trho_{q,j+1}\comp$ or $\tm_{q,j+1}\comp$.

	Apply \eqref{eq:Newton-quadratic-product-material-estimate}--\eqref{eq:Newton-initial-mixed-coefficient-material-estimate} to the two types of products in \eqref{eq:Newton-perturbative-tensor-increment-estimate}. For $0\leq j<J_*$ and $M+N\leq N_{\mathrm{ind}}$, we have
	\begin{equation}\label{eq:Newton-perturbative-tensor-material-increment-estimate}
		\begin{aligned}
			\nrm{\DTL^M(P_{q+1,j+1}-P_{q+1,j}-P_{q+1,j+1}^{\mathrm{mol,comp}}-P_{q+1,j+1}^{\mathrm{mol,inc}})}_N\lesssim\left[A_{q,j}A_{q,0}\delta_q^{-1}\lambda_{q+1}^{-6\alpha}+A_{q,j}\lambda_{q+1}^{-3\alpha}(\ell_q\lambda_q)^{N_{\mathrm{ind}}}\right]\ell_{t,q}^{-M}\ell_q^{-N}.
		\end{aligned}
	\end{equation}	
	For $0\leq j\leq J_*$ and $M+N\leq N_R$, summing \eqref{eq:Newton-perturbative-tensor-material-increment-estimate} gives
	\begin{equation}\label{eq:Newton-accumulated-perturbative-tensor-material-estimate}
		\begin{aligned}
			\nrm{\DTL^MP_{q+1,j}}_N
			&\leq\nrm{\DTL^MP_{q+1,0}}_N+\sum_{\nu=0}^{j-1}\nrm{\DTL^M(P_{q+1,\nu+1}-P_{q+1,\nu})}_N\\
			&\lesssim\nrm{\DTL^MP_{q+1,0}}_N+\sum_{\nu=1}^{j}\left(\nrm{\DTL^MP_{q+1,\nu}^{\mathrm{mol,comp}}}_N+\nrm{\DTL^MP_{q+1,\nu}^{\mathrm{mol,inc}}}_N\right)\\
			&\quad+\left[A_{q,0}\delta_q^{-1}\lambda_{q+1}^{-6\alpha}+\lambda_{q+1}^{-3\alpha}(\ell_q\lambda_q)^{N_{\mathrm{ind}}}\right]\left(\sum_{\nu=0}^{j-1}A_{q,\nu}\right)\ell_{t,q}^{-M}\ell_q^{-N}.
		\end{aligned}
	\end{equation}
	By \eqref{def:admissible-alpha-losses}, $b>1$, and $\beta<1/3$,
	$(8+6b)\alpha+2\beta(b-1)+(3+b)\alpha/4<2(b-1)/3$.
	 Using $N_R\geq4J_*+12$, we obtain
	\begin{equation}\label{eq:Newton-small-error-allocation}
	\begin{aligned}
		&\frac{(\ell_q\lambda_q)^{N_R}+\ell_q^{-\alpha}(\ell_q\lambda_q)^{\frac{N_R-3J_*}{2}}+\lambda_{q+1}^{-3\alpha}(\ell_q\lambda_q)^{N_{\mathrm{ind}}}}
		{A_{q,0}\delta_q^{-1}\lambda_{q+1}^{-6\alpha}}
		\lesssim\lambda_q^{(8+6b)\alpha+2\beta(b-1)+\frac{(3+b)\alpha}{4}-\frac{(b-1)(N_R-3J_*)}{8}}\leq\lambda_q^{-\frac56(b-1)},
	\end{aligned}
	\end{equation}
	for sufficiently large $a$. 
	Substituting \eqref{eq:Newton-substep-mollification-Newton-scales} and \eqref{eq:Newton-initial-perturbative-tensor-material-estimate} into \eqref{eq:Newton-accumulated-perturbative-tensor-material-estimate}, we find
	$$
	\begin{aligned}
		\sup_{0\leq j\leq J_*}\nrm{\DTL^MP_{q+1,j}}_N
		&\lesssim A_{q,0}\Bigl[
		(\ell_q\lambda_q)^{N_R}+\ell_q^{-\alpha}(\ell_q\lambda_q)^{\frac{N_R-3J_*}{2}}
		+\lambda_{q+1}^{-3\alpha}(\ell_q\lambda_q)^{N_{\mathrm{ind}}}
		+A_{q,0}\delta_q^{-1}\lambda_{q+1}^{-6\alpha}\Bigr]
		\ell_{t,q}^{-M}\ell_q^{-N}
		\lesssim A_{q,0}^2\delta_q^{-1}\lambda_{q+1}^{-6\alpha}\ell_{t,q}^{-M}\ell_q^{-N}.
	\end{aligned}
	$$
	For the last inequality, we used \eqref{eq:Newton-small-error-allocation}.
	Thus \eqref{eq:Newton-substep-perturbative-tensor-detailed-estimate} holds.
	It remains to prove \eqref{eq:Newton-terminal-temporal-support-next-scale}. Taking $j=J_*$ in \eqref{eq:Newton-iterated-temporal-support} and \eqref{eq:Newton-partial-state-stress-temporal-support}, we obtain
\begin{equation}\label{eq:Newton-terminal-temporal-support}
		\supp_t(\rho_{q,J_*}-\rho_q)\cup\supp_t(m_{q,J_*}-m_q)
		\cup\supp_tR_{q,J_*}\cup\supp_tS_{q,J_*}\cup\supp_tP_{q+1,J_*}\subset\mN_{(2J_*+1)\tau_q}(\supp_tR_q).
\end{equation}
Since $J_*$ is fixed, \eqref{def:tau-q} gives, for sufficiently large $a$,
$
(2J_*+1)\tau_q\leq\frac12\Lambda_q^{-1}.
$
	Combining \eqref{ind:R-support} with \eqref{eq:Newton-terminal-temporal-support} proves \eqref{eq:Newton-terminal-temporal-support-next-scale}.
\end{proof}
\section{The Nash step}
\label{sec:Nash-step}

We start from the terminal state of the Newton step $(\rho_{q,J_*},m_{q,J_*})$ and construct a divergence-free momentum perturbation $w_{q+1}$.
We set $\rho_{q+1}=\rho_{q,J_*}$ and $m_{q+1}=m_{q,J_*}+w_{q+1}$.
Since $\Div w_{q+1}=0$, the continuity equation is preserved. The perturbation is built from spatially separated Mikado flows transported by the spatial mollification $\bar u_{q,J_*}$ of the terminal velocity in the Newton step.
We choose the amplitudes so that the mollified gluing stress can be canceled, see \eqref{def:Nash-pulled-back-matrix} and \eqref{eq:Nash-principal-quadratic-average}.

In Subsection~\ref{subsec:Nash-mollification-estimates}, we estimate the density, momentum, velocity, and stresses at the end of the Newton step with the material derivative
generated by $\bar u_{q,J_*}$ and define the flow maps and space--time cutoffs.
We define the Mikado profiles and $w_{q+1}$ in Subsection~\ref{subsec:Nash-separated-Mikado}.
In Subsection~\ref{subsec:Nash-new-Reynolds-stress}, we solve the divergence equations for localized oscillatory sources to construct the transport and oscillation errors and define the new Reynolds stress $R_{q+1}$.
Subsection~\ref{subsec:Nash-estimates} estimates the perturbation and the transport, oscillation, mollification, and corrector errors in $C^N$.
In Subsection~\ref{subsec:closure-density-momentum}, we estimate $\rho_{q+1}-\rho_q$, $m_{q+1}-m_q$, $\rho_{q+1}$, $m_{q+1}$, and $R_{q+1}$ with $D_{t,q+1}$.

We use the parameters from Section~\ref{sec:Newton steps}.
The kernels and periodic profiles used below are fixed independently of $q$ and $a$.
Choose a common threshold
$
    a_{\mathrm{Nash}}=a_{\mathrm{Nash}}(\beta,b,\alpha,\rho^*,\mathfrak M,C_*,N_R,N_{\mathrm{ind}},N_{\mathrm{max}})
    \geq a_{\mathrm{Newton}}
$
for the construction, estimates, and support conditions below, and assume
$a>a_{\mathrm{Nash}}$ throughout this section.

\subsection{Mollification, flow maps, and cutoffs in the Nash step}
\label{subsec:Nash-mollification-estimates}

\subsubsection{Mollification and stress estimates}

Using the spatial mollification from Section~\ref{Mollification}, set
\begin{equation}\label{def:Nash-mollified-terminal-fields}
	\bar u_{q,J_*}:=P_{\leq\ell_q^{-1}}u_{q,J_*},
	\qquad
	\bar\rho_{q,J_*}:=P_{\leq\ell_q^{-1}}\rho_{q,J_*}.
\end{equation}
The material derivative generated by $\bar u_{q,J_*}$ is
$\widetilde D_{t,q}:=\partial_t+\bar u_{q,J_*}\cdot\nabla$.
The following lemma converts estimates with $\DTL$ into estimates with $D_{t,\upsilon}$.
We will use it with $\upsilon=\bar u_{q,J_*}$.

\begin{lem}
	\label{lem:Newton-two-scale-material-change}
	Assume the background estimate \eqref{est:uniform-grad-om-over-orho} and the scale relations \eqref{eq:material-mollification-spatial-comparisons} and \eqref{eq:Newton-scale-separation}. Fix $J,K\in\mathbb N$ with $1\leq J\leq\min\{K,N_{\mathrm{ind}}+1\}$ and $K\leq N_{\mathrm{max}}+1$.
	Let $C_0\geq1$ and let $\upsilon$ be a smooth velocity. Set
	$D_{t,\upsilon}:=\partial_t+\upsilon\cdot\nabla$. Suppose that, for $M,N\in\mathbb N_0$ with $0\leq M<J$ and $M+N\leq K-1$,
	$$
	\|\DTL^M(\upsilon-\overline u_q)\|_N
	\leq C_0\delta_q^{\frac12}\MM{M,N_R-3J_*,\tau_q^{-1},\ell_{t,q}^{-1}}
	\MM{N,N_{R,J_*}^{(0)}(M),\lambda_q,\ell_q^{-1}}.
	$$
	Let $A>0$. If a smooth scalar, vector, or tensor field $F$ satisfies, for $M,N\in\mathbb N_0$,
	$$
	\|\DTL^MF\|_N
	\leq A\MM{M,N_R-3J_*,\tau_q^{-1},\ell_{t,q}^{-1}}
	\MM{N,N_{R,J_*}^{(0)}(M),\lambda_q,\ell_q^{-1}},
	\qquad 0\leq M\leq J,\quad M+N\leq K,
	$$
	then, in the same range,
	\begin{equation}\label{eq:Newton-two-scale-material-change}
		\|D_{t,\upsilon}^MF\|_N
		\lesssim A\MM{M,N_R-3J_*,\tau_q^{-1},\ell_{t,q}^{-1}}
		\MM{N,N_{R,J_*}^{(0)}(M),\lambda_q,\ell_q^{-1}}.
	\end{equation}
	The implicit constant may depend on $J,K,C_0,\rho^*,\mathfrak M,N_{\mathrm{max}}$.
	It is independent of $A$, $\upsilon$, $F$, $q$, and $a$.
\end{lem}

\begin{proof}
	The case $M=0$ follows from the hypothesis on $F$.
	For $M\geq1$, apply Lemma~\ref{lem:source-material-derivative-change} with background velocity $\overline u_q$, velocity difference $\upsilon-\overline u_q$, and
	$$
	\begin{aligned}
		\lambda=\lambda_q,\quad \tilde\lambda=\ell_q^{-1},\quad \mu=\tau_q^{-1},\quad \tilde\mu=\ell_{t,q}^{-1},\quad
		N_t=N_R-3J_*,\quad \mC_\omega=\delta_q^{\frac12},\quad \mC_F=A.
	\end{aligned}
	$$
	Choose $C\geq C_0$ larger than the implicit constant in
	\eqref{est:uniform-grad-om-over-orho} over its stated range.
	In Lemma~\ref{lem:source-material-derivative-change}, take $\widetilde{\mC}_0=C$.
	For fixed $C_0$, this choice is independent of $A$, $\upsilon$, $F$, $q$, and $a$.
	Definition~\eqref{def:Newton-derivative-bookkeeping} gives $N_{R,J_*}^{(0)}(M)=(N_t-M)_+$.
	The derivative parameters in Subsection~\ref{subsec:Newton-scheme-parameters} give $0\leq N_t\leq N_{\mathrm{ind}}-1$.
	The scale conditions follow from \eqref{eq:material-mollification-spatial-comparisons} and
	\eqref{eq:Newton-scale-separation}. In particular,
	$
		\mC_\omega\lambda
		=\delta_q^{\frac12}\lambda_q=\Lambda_q\leq\mu,\
		\mC_\omega\tilde\lambda
		=\delta_q^{\frac12}\ell_q^{-1}\leq\ell_{t,q}^{-1}=\tilde\mu.
	$
	The hypotheses on $\upsilon-\overline u_q$ and $F$ verify the first and third
	estimates in \eqref{ass:source-two-scale-material-change}.

	It remains to verify the estimate for $\DTL^{M_0}\nabla\overline u_q$ in
	\eqref{ass:source-two-scale-material-change}, required only for $K\geq2$.
	For $0\leq M_0<J$ and $M_0+N_0\leq K-2$, we have
	$M_0\leq N_{\mathrm{ind}}$ and $M_0+N_0\leq N_{\mathrm{max}}-1$, so \eqref{est:uniform-grad-om-over-orho} gives
	$$
	\begin{aligned}
		\|\DTL^{M_0}\nabla\overline u_q\|_{N_0}
		&\lesssim\delta_q^{\frac12}\Lambda_q^{M_0}
		\MM{N_0+1,N_{\mathrm{ind}}-M_0,\lambda_q,\ell_q^{-1}}\\
		&\leq\Lambda_q
		\MM{M_0,N_{\mathrm{ind}}-1,\mu,\tilde\mu}
		\MM{N_0,(N_{\mathrm{ind}}-1-M_0)_+,\lambda_q,\ell_q^{-1}}\\
		&\leq\mu
		\MM{M_0,N_t,\mu,\tilde\mu}
		\MM{N_0,(N_t-M_0)_+,\lambda_q,\ell_q^{-1}}.
	\end{aligned}
	$$
	Here the second inequality uses \eqref{eq:Newton-one-threshold-weight-shift} with threshold $N_{\mathrm{ind}}$ and $\delta_q^{\frac12}\lambda_q=\Lambda_q\leq\mu$.
	At $M_0=N_{\mathrm{ind}}$, this shift uses $\mu\ell_{t,q}\ell_q^{-1}\leq\lambda_q$ from \eqref{eq:material-mollification-spatial-comparisons}.
	With $\widetilde{\mC}_0=C$, this verifies the remaining hypothesis.
	Lemma~\ref{lem:source-material-derivative-change} now gives
	\eqref{eq:Newton-two-scale-material-change}.
\end{proof}

We now estimate $\bar u_{q,J_*}$ and $\bar\rho_{q,J_*}$ with $\widetilde D_{t,q}$, together with $u_{q,J_*}-\bar u_{q,J_*}$ and $\rho_{q,J_*}-\bar\rho_{q,J_*}$.
We also estimate $m_{q,J_*}$, the density and momentum perturbations in the Newton step, and the stress tensors with the same material derivative.

\begin{prop}
	\label{prop:Nash-terminal-Newton-velocity-estimates}
	Fix $q\in\mathbb N_0$. Assume the pressure-law conditions in Section~\ref{sec:main-iteration}
	and the hypotheses of Proposition~\ref{prop:uniform-mollification-estimates}.
	Let $(\rho_{q,J_*},m_{q,J_*})$ be the terminal state of the Newton scheme
	with the parameters in Subsection~\ref{subsec:Newton-scheme-parameters}.
	Let $R_{q,J_*}$, $S_{q,J_*}$, and $P_{q+1,j}$, $0\leq j\leq J_*$, be the stress tensors from the Newton step.
	For $a>a_{\mathrm{Nash}}$, we have
	\begin{equation}\label{eq:Nash-mollified-terminal-zero-order-bounds}
	\nrm{\bar u_{q,J_*}}_0\leq\frac{2\mathfrak M}{\rho^*}\left(1-\frac14\delta_q^{\frac12}\right),\qquad\frac{\rho^*}{2}\leq\bar\rho_{q,J_*}\leq\frac{3\rho^*}{2}.
	\end{equation}
	For $M,N\in\mathbb N_0$ with $0\leq M\leq N_{\mathrm{ind}}$ and
	$1\leq M+N\leq N_{\mathrm{max}}$, the mollified velocity $\bar u_{q,J_*}$ satisfies
	\begin{equation}\label{eq:Nash-mollified-terminal-velocity-integer-estimate}
	\begin{aligned}
		\nrm{\widetilde D_{t,q}^M\bar u_{q,J_*}}_N
		&\lesssim\delta_q^{\frac12}
		\MM{M,N_R-3J_*,\tau_q^{-1},\ell_{t,q}^{-1}}
		\MM{N,N_{R,J_*}^{(0)}(M),\lambda_q,\ell_q^{-1}}.
	\end{aligned}
	\end{equation}
	For $M,N\in\mathbb N_0$ with $0\leq M\leq N_{\mathrm{ind}}$
	and $1\leq M+N\leq N_{\mathrm{max}}-1$, spatial interpolation gives
	\begin{equation}\label{eq:Nash-mollified-terminal-velocity-estimate}
	\begin{aligned}
		\nrm{\widetilde D_{t,q}^M\bar u_{q,J_*}}_{N+\alpha}
		\lesssim\delta_q^{\frac12}
		\MM{M,N_R-3J_*,\tau_q^{-1},\ell_{t,q}^{-1}}
		\La{N,N_{R,J_*}^{(0)}(M),\lambda_q,\ell_q^{-1}}
		\MM{N,N_{R,J_*}^{(0)}(M),\lambda_q,\ell_q^{-1}}.
	\end{aligned}
	\end{equation}

	For each fixed $s\in\mathbb R$ and $M,N\in\mathbb N_0$ with $1\leq M+N\leq N_{\mathrm{ind}}$, we have
	\begin{equation}\label{eq:Nash-density-coefficient-material-estimate}
		\Nrm{\widetilde D_{t,q}^M[(\rho_{q,J_*})^s]}_N\lesssim_s\delta_q\MM{M,N_R-3J_*,\tau_q^{-1},\ell_{t,q}^{-1}}\MM{N,N_{R,J_*}^{(0)}(M),\lambda_q,\ell_q^{-1}}.
	\end{equation}
	For $M,N\in\mathbb N_0$ with $1\leq M+N\leq N_{\mathrm{ind}}$, the terminal momentum in the Newton step satisfies
	\begin{equation}\label{eq:Nash-terminal-momentum-material-estimate}
	\begin{aligned}
		\nrm{\widetilde D_{t,q}^Mm_{q,J_*}}_N
		&\lesssim\delta_q^{\frac12}
		\MM{M,N_R-3J_*,\tau_q^{-1},\ell_{t,q}^{-1}}
		\MM{N,N_{R,J_*}^{(0)}(M),\lambda_q,\ell_q^{-1}}.
	\end{aligned}
	\end{equation}
	For each fixed $s\in\mathbb R$ and $M,N\in\mathbb N_0$ with $M\leq N_{\mathrm{ind}}$ and $1\leq M+N\leq N_{\mathrm{max}}$, we have
	\begin{equation}\label{eq:Nash-mollified-density-material-estimate}
	\begin{aligned}
		\Nrm{\widetilde D_{t,q}^M[(\bar\rho_{q,J_*})^s]}_N
		&\lesssim_s\delta_q\MM{M,N_R-3J_*,\tau_q^{-1},\ell_{t,q}^{-1}}
		\MM{N,N_{R,J_*}^{(0)}(M),\lambda_q,\ell_q^{-1}}.
	\end{aligned}
	\end{equation}

	For the gradient estimates, set
\begin{equation}
	\label{def:Nash-terminal-material-threshold}
	N_t^{\mathrm N}:=N_R-3J_*-4.
\end{equation}
	For $M,N\in\mathbb N_0$ with $0\leq M\leq N_{\mathrm{ind}}$ and
	$M+N\leq N_{\mathrm{max}}-1$, the velocity gradient $\nabla\bar u_{q,J_*}$ satisfies
	\begin{equation}\label{eq:Nash-terminal-Newton-velocity-integer-estimate}
	\begin{aligned}
		\nrm{\widetilde D_{t,q}^M\nabla\bar u_{q,J_*}}_N
		&\lesssim
		\tau_q^{-1}\MM{M,N_t^{\mathrm N}-1,\tau_q^{-1},\ell_{t,q}^{-1}}
		\MM{N,N_{R,J_*}^{(1)}(M),\lambda_q,\ell_q^{-1}}.
	\end{aligned}
	\end{equation}
	Moreover, for $M,N\in\mathbb N_0$ with
	$0\leq M\leq N_{\mathrm{ind}}$ and
	$M+N\leq N_{\mathrm{max}}-2$,
	\begin{equation}
		\begin{aligned}
			\nrm{\widetilde D_{t,q}^M\nabla \bar u_{q,J_*}}_{N+\alpha}
			&\lesssim
			\tau_q^{-1}\MM{M,N_t^{\mathrm N}-1,\tau_q^{-1},\ell_{t,q}^{-1}}
			\La{N,N_{R,J_*}^{(1)}(M),\lambda_q,\ell_q^{-1}}
			\MM{N,N_{R,J_*}^{(1)}(M),\lambda_q,\ell_q^{-1}}.
		\end{aligned}
		\label{eq:Nash-terminal-Newton-velocity-estimate}
	\end{equation}

	The spatial mollification differences satisfy, for
	$M,N,L\in\mathbb N_0$ with $M+N+L\leq N_{\mathrm{ind}}$,
\begin{equation}\label{eq:Nash-drift-mollification-difference}
	\begin{aligned}
		\|\widetilde D_{t,q}^M(u_{q,J_*}-\bar u_{q,J_*})\|_N
		&\lesssim\ell_q^L\delta_q^{\frac12}\MM{M,N_t^{\mathrm N}+4,\tau_q^{-1},\ell_{t,q}^{-1}}\MM{N+L,N_{R,J_*}^{(1)}(M),\lambda_q,\ell_q^{-1}}.
	\end{aligned}
\end{equation}
	\begin{equation}\label{eq:Nash-density-mollification-difference}
	\|\widetilde D_{t,q}^M(\rho_{q,J_*}-\bar\rho_{q,J_*})\|_N
	\lesssim\ell_q^L\delta_q\MM{M,N_R-3J_*,\tau_q^{-1},\ell_{t,q}^{-1}}
	\MM{N+L,N_{R,J_*}^{(1)}(M),\lambda_q,\ell_q^{-1}}.
	\end{equation}
	For $M,N\in\mathbb N_0$ with $M\leq N_{\mathrm{ind}}$ and $M+N\leq K_{J_*}^{(0)}$, the momentum perturbation in the Newton step satisfies
	\begin{equation}\label{eq:Nash-terminal-Newton-momentum-increment-material}
	\begin{aligned}
		\nrm{\widetilde D_{t,q}^M(m_{q,J_*}-m_q)}_N
		&\lesssim A_{q,0}\delta_q^{-\frac12}\lambda_{q+1}^{-3\alpha}
		\MM{M,N_R-3J_*,\tau_q^{-1},\ell_{t,q}^{-1}}
		\MM{N,N_{R,J_*}^{(0)}(M),\lambda_q,\ell_q^{-1}}.
	\end{aligned}
	\end{equation}
	In the same derivative range, the density perturbation in the Newton step satisfies
\begin{equation}\label{eq:Nash-terminal-Newton-density-increment-material}
\begin{aligned}
    \nrm{\widetilde D_{t,q}^M(\rho_{q,J_*}-\rho_q)}_N
    &\lesssim A_{q,0}\ell_q^{-\alpha/2}
    \MM{M,N_R-3J_*,\tau_q^{-1},\ell_{t,q}^{-1}}
    \MM{N,N_{R,J_*}^{(0)}(M),\lambda_q,\ell_q^{-1}}.
\end{aligned}
\end{equation}

	The stress tensors satisfy the following estimates for $M,N\in\mathbb N_0$.
	For $M\leq N_{\mathrm{ind}}$ and $M+N\leq K_{J_*}^{(0)}+1$,
	\begin{equation}\label{eq:Nash-terminal-S-material-estimate}
		\|\widetilde D_{t,q}^MS_{q,J_*}\|_N\lesssim A_{q,0}(\tau_q\ell_q^{-1})^\alpha\MM{M,N_R-3J_*,\tau_q^{-1},\ell_{t,q}^{-1}}\MM{N,N_{R,J_*}^{(0)}(M),\lambda_q,\ell_q^{-1}}.
	\end{equation}
	For $M+N\leq N_R$, we have
	\begin{equation}\label{eq:Nash-retained-Newton-stresses-material}
		\sup_{0\leq j\leq J_*}\|\widetilde D_{t,q}^MP_{q+1,j}\|_N+\|\widetilde D_{t,q}^MR_{q,J_*}\|_N\leq\lambda_q^{-\alpha}\lambda_{q+1}^{-8\alpha}\delta_{q+2}\ell_{t,q}^{-M}\ell_q^{-N}.
	\end{equation}
	The implicit constants may depend on $\beta,b,\alpha,N_{\mathrm{max}},\rho^*,\mathfrak M,C_*$, unless stated otherwise.
	They are independent of $q$ and $a$.
	The implicit constants in \eqref{eq:Nash-terminal-S-material-estimate} depend only on $N_{\mathrm{max}}$ and $\alpha$.
\end{prop}

\begin{proof}
	To pass from $\DTL$ to $\widetilde D_{t,q}$, we need an estimate for $\bar u_{q,J_*}-\overline u_q$.
	We use the decomposition
	$$
	\bar u_{q,J_*}-\overline u_q=P_{\leq\ell_q^{-1}}(u_{q,J_*}-\overline u_q)-P_{>\ell_q^{-1}}\overline u_q.
	$$
	We first estimate $u_{q,J_*}-\overline u_q$ and $u_{q,J_*}$ with $\DTL$, then the spatial mollifications.
	Substituting \eqref{eq:Newton-terminal-velocity-difference-amplitudes} into \eqref{eq:Newton-terminal-velocity-difference-two-scale}, we obtain, for
	$M,N\in\mathbb N_0$ with $M+N\leq N_{\mathrm{ind}}$,
	$$
	\|\DTL^M(u_{q,J_*}-\overline u_q)\|_N
	\lesssim\delta_{q+1}^{\frac12}
	\MM{M,N_R-3J_*,\tau_q^{-1},\ell_{t,q}^{-1}}
	\MM{N,N_{R,J_*}^{(0)}(M),\lambda_q,\ell_q^{-1}}.
	$$
	This estimate and \eqref{est:uniform-om-over-orho} imply, for
	$1\leq M+N\leq N_{\mathrm{ind}}$,
	\begin{equation}\label{eq:Nash-terminal-velocity-input-estimate}
	\begin{aligned}
		\nrm{\DTL^Mu_{q,J_*}}_N
		&\lesssim\delta_q^{\frac12}
		\MM{M,N_R-3J_*,\tau_q^{-1},\ell_{t,q}^{-1}}
		\MM{N,N_{R,J_*}^{(0)}(M),\lambda_q,\ell_q^{-1}}.
	\end{aligned}
	\end{equation}
	For $\rho_{q,J_*}$, choose $C\geq1$ to cover the $C^0$ bound $3\rho^*/2$ from \eqref{eq:Newton-state-zero-order-bounds} and the implicit constant in \eqref{eq:Newton-terminal-density-material-estimate}.

	For $\bar\rho_{q,J_*}$, \eqref{eq:Newton-state-zero-order-bounds}
	and \eqref{est:appendix-spatial-mollification} with $N=0$ imply
	$$
	\|\bar\rho_{q,J_*}-\rho^*\|_0
	\leq\|\eta\|_{L^1(\mathbb R^3)}\|\rho_{q,J_*}-\rho^*\|_0
	\leq\|\eta\|_{L^1(\mathbb R^3)}\mathfrak M\sum_{i=1}^{q+1}\delta_i\leq\frac{\rho^*}{2}.
	$$
	Here the last inequality follows from the choice of $a_{\mathrm{mol}}$ used in \eqref{eq:mollification-density-smallness},
	since $a_{\mathrm{Nash}}\geq a_{\mathrm{Newton}}\geq a_{\mathrm{mol}}$.
	This proves the second estimate in \eqref{eq:Nash-mollified-terminal-zero-order-bounds}.
	For the mollified velocity, \eqref{eq:Nash-terminal-velocity-input-estimate} with $M=0$ and $N=1$
	and \eqref{est:appendix-spatial-mollification} with $N=0$ and $L=1$ imply
	$$
	\nrm{u_{q,J_*}-\bar u_{q,J_*}}_0
	\lesssim\ell_q\nrm{u_{q,J_*}}_1
	\lesssim\delta_q^{\frac12}\ell_q\lambda_q.
	$$
	Enlarge $C$ to cover the estimate for $\nrm{u_{q,J_*}-\bar u_{q,J_*}}_0$.
	By \eqref{eq:Newton-state-zero-order-bounds} with $j=J_*$, we have
	$$
	\begin{aligned}
		\nrm{\bar u_{q,J_*}}_0
		\leq\nrm{u_{q,J_*}}_0+\nrm{u_{q,J_*}-\bar u_{q,J_*}}_0
		\leq\frac{2\mathfrak M}{\rho^*}\left(1-\frac12\delta_q^{\frac12}\right)
		+C\delta_q^{\frac12}\ell_q\lambda_q
		\leq\frac{2\mathfrak M}{\rho^*}\left(1-\frac14\delta_q^{\frac12}\right).
	\end{aligned}
	$$
	For the last inequality, choose $a_{\mathrm{Nash}}$ so that
	$
		C\ell_q\lambda_q\leq(2\rho^*)^{-1}.
	$
	By \eqref{def:parameter-mollification}, we can make this choice and obtain the first estimate in \eqref{eq:Nash-mollified-terminal-zero-order-bounds}.

	We now estimate the positive-order derivatives of $\bar u_{q,J_*}$ and $\bar\rho_{q,J_*}$ with $\DTL$.
	Apply \eqref{est:appendix-two-scale-extended-spatial-mollification}
	in Proposition~\ref{prop:appendix-material-mollification} with $\upsilon=\overline u_q$.
	Take $f=u_{q,J_*}$ with $\mC_f=\delta_q^{\frac12}$ or $f=\rho_{q,J_*}$ with $\mC_f=\delta_q$, and set
	$$
	\begin{aligned}
		\ell&=\ell_q,\quad K_f=N_R-3J_*,\quad K^*=K_\upsilon=N_{\mathrm{ind}},\quad K^\sharp=N_{\mathrm{max}},\quad \lambda=\lambda_q,\quad \mu=\tau_q^{-1},\quad \tilde\mu=\ell_{t,q}^{-1},\quad \mC_\upsilon=C\delta_q^{\frac12}.
	\end{aligned}
	$$
	The source assumptions follow from \eqref{eq:Nash-terminal-velocity-input-estimate} and
	\eqref{eq:Newton-terminal-density-material-estimate}.
	Enlarge $C$ so that \eqref{est:uniform-om-over-orho}--%
	\eqref{est:uniform-grad-om-over-orho} verify the velocity assumptions.
	The scale conditions follow from \eqref{eq:material-mollification-spatial-comparisons}
	and \eqref{eq:Newton-scale-separation}, after increasing $a_{\mathrm{Nash}}$.
	For $0\leq M\leq N_{\mathrm{ind}}$ and $1\leq M+N\leq N_{\mathrm{max}}$, we find
	\begin{equation}\label{eq:Nash-mollified-velocity-background-estimate}
	\|\DTL^M\bar u_{q,J_*}\|_N
	\lesssim\delta_q^{\frac12}
	\MM{M,N_R-3J_*,\tau_q^{-1},\ell_{t,q}^{-1}}
	\MM{N,N_{R,J_*}^{(0)}(M),\lambda_q,\ell_q^{-1}}.
	\end{equation}
	In the same range, the mollified density $\bar\rho_{q,J_*}$ satisfies
	\begin{equation}\label{eq:Nash-mollified-density-background-estimate}
	\|\DTL^M\bar\rho_{q,J_*}\|_N
	\lesssim\delta_q\MM{M,N_R-3J_*,\tau_q^{-1},\ell_{t,q}^{-1}}
	\MM{N,N_{R,J_*}^{(0)}(M),\lambda_q,\ell_q^{-1}}.
	\end{equation}
	At $M=N=0$, use \eqref{eq:Nash-mollified-terminal-zero-order-bounds}.

	Next we estimate $\bar u_{q,J_*}-\overline u_q$.
	First take $0\leq M\leq N_{\mathrm{ind}}$ and $1\leq M+N\leq N_{\mathrm{max}}$.
	In both applications of Proposition~\ref{prop:appendix-material-mollification}, take
	$$
		\upsilon=\overline u_q,\quad \ell=\ell_q,\quad K_f=N_R-3J_*,\quad K^*=K_\upsilon=N_{\mathrm{ind}},\quad K^\sharp=N_{\mathrm{max}},\quad
		\lambda=\lambda_q,\quad \mu=\tau_q^{-1},\quad \tilde\mu=\ell_{t,q}^{-1},\quad \mC_\upsilon=C\delta_q^{\frac12}.
	$$
	The choice of $C$, the velocity assumptions, and the scale conditions are as in the previous application.
	For $P_{\leq\ell_q^{-1}}(u_{q,J_*}-\overline u_q)$, set $f=u_{q,J_*}-\overline u_q$ and $\mC_f=\delta_q^{\frac12}$
	in \eqref{est:appendix-two-scale-extended-spatial-mollification}.
	The source assumptions follow from \eqref{eq:Newton-terminal-velocity-difference-two-scale}
	and \eqref{eq:Newton-terminal-velocity-difference-amplitudes}.
	For $0\leq M\leq N_{\mathrm{ind}}$ and $1\leq M+N\leq N_{\mathrm{max}}$, we obtain
	$$
	\begin{aligned}
		\|\DTL^MP_{\leq\ell_q^{-1}}(u_{q,J_*}-\overline u_q)\|_N
		&\lesssim\delta_q^{\frac12}\MM{M,N_R-3J_*,\tau_q^{-1},\ell_{t,q}^{-1}}
		\MM{N,N_{R,J_*}^{(0)}(M),\lambda_q,\ell_q^{-1}}.
	\end{aligned}
	$$
	For $P_{>\ell_q^{-1}}\overline u_q$, apply Proposition~\ref{prop:appendix-material-mollification}
	with $f=\overline u_q$ and $\mC_f=\delta_q^{\frac12}$, leaving the displayed parameters unchanged.
	The source assumptions follow from \eqref{est:uniform-om-over-orho}, since $\Lambda_q\leq\tau_q^{-1}$ and $K_f\leq N_{\mathrm{ind}}$.
	For $0\leq M\leq N_{\mathrm{ind}}$ and $1\leq M+N\leq N_{\mathrm{max}}$,
	\eqref{est:appendix-two-scale-extended-spatial-mollification} gives
	$$
	\begin{aligned}
		\|\DTL^MP_{\leq\ell_q^{-1}}\overline u_q\|_N
		&\lesssim\delta_q^{\frac12}\MM{M,N_R-3J_*,\tau_q^{-1},\ell_{t,q}^{-1}}
		\MM{N,N_{R,J_*}^{(0)}(M),\lambda_q,\ell_q^{-1}}.
	\end{aligned}
	$$
	This estimate and \eqref{est:uniform-om-over-orho} imply,
	for $0\leq M\leq N_{\mathrm{ind}}$ and $1\leq M+N\leq N_{\mathrm{max}}$,
	$$
	\begin{aligned}
		\|\DTL^MP_{>\ell_q^{-1}}\overline u_q\|_N
		\leq\|\DTL^M\overline u_q\|_N+\|\DTL^MP_{\leq\ell_q^{-1}}\overline u_q\|_N
		\lesssim\delta_q^{\frac12}\MM{M,N_R-3J_*,\tau_q^{-1},\ell_{t,q}^{-1}}
		\MM{N,N_{R,J_*}^{(0)}(M),\lambda_q,\ell_q^{-1}}.
	\end{aligned}
	$$

	At $M=N=0$, Proposition~\ref{prop:appendix-spatial-mollification}
	with $\ell=\ell_q$ and $N=0$, using $L=1$ for the mollification difference, gives
	$$
	\begin{aligned}
		\|P_{\leq\ell_q^{-1}}(u_{q,J_*}-\overline u_q)\|_0
		+\|P_{>\ell_q^{-1}}\overline u_q\|_0\lesssim\|u_{q,J_*}-\overline u_q\|_0
		+\ell_q\|\overline u_q\|_1
		\lesssim\delta_q^{\frac12}.
	\end{aligned}
	$$
	For the last inequality, we used \eqref{est:uniform-om-over-orho} with $M=0$ and $N=1$,
	\eqref{eq:Newton-terminal-velocity-difference-two-scale} with $M=N=0$, and
	\eqref{eq:Newton-terminal-velocity-difference-amplitudes}.
	Combining the two terms of the decomposition, we have, for $0\leq M\leq N_{\mathrm{ind}}$ and $M+N\leq N_{\mathrm{max}}$,
	\begin{equation}\label{eq:Nash-drift-background-difference-two-scale-estimate}
	\begin{aligned}
		\|\DTL^M(\bar u_{q,J_*}-\overline u_q)\|_N\
		&\leq\|\DTL^MP_{\leq\ell_q^{-1}}(u_{q,J_*}-\overline u_q)\|_N
		+\|\DTL^MP_{>\ell_q^{-1}}\overline u_q\|_N\\
		&\lesssim\delta_q^{\frac12}\MM{M,N_R-3J_*,\tau_q^{-1},\ell_{t,q}^{-1}}
		\MM{N,N_{R,J_*}^{(0)}(M),\lambda_q,\ell_q^{-1}}.
	\end{aligned}
	\end{equation}
	Choose $\widetilde{\mC}_0\geq1$ larger than the constant in \eqref{eq:Nash-drift-background-difference-two-scale-estimate}.

	Before changing the material derivative, we estimate $\DTL^M\nabla\bar u_{q,J_*}$.
	For $0\leq M\leq N_{\mathrm{ind}}$ and $M+N\leq N_{\mathrm{max}}-1$,
	apply Lemma~\ref{lem:appendix-scale-dependent-commuted-gradient} with $\upsilon=\overline u_q$, $G=\bar u_{q,J_*}$, and
	$$
	\begin{aligned}
		K_\upsilon&=N_{\mathrm{ind}},\quad K_G=N_t=N_R-3J_*,\quad K^\sharp=N_{\mathrm{max}},\quad \mC_\upsilon=\mC_G=\delta_q^{\frac12},\quad \lambda=\lambda_q,\quad \widetilde\lambda=\ell_q^{-1},\quad \mu=\tau_q^{-1},\quad \tilde\mu=\ell_{t,q}^{-1}.
	\end{aligned}
	$$
	The hypotheses follow from \eqref{est:uniform-grad-om-over-orho} and \eqref{eq:Nash-mollified-velocity-background-estimate}. The lemma gives
	\begin{equation}\label{eq:Nash-gradient-threshold-shift}
	\begin{aligned}
		\|\DTL^M\nabla\bar u_{q,J_*}\|_N+\|\nabla\DTL^M\bar u_{q,J_*}\|_N
		&\lesssim\delta_q^{\frac12}\MM{M,N_R-3J_*,\tau_q^{-1},\ell_{t,q}^{-1}}
		\MM{N+1,N_{R,J_*}^{(0)}(M),\lambda_q,\ell_q^{-1}}\\
		&\leq\Lambda_q\MM{M,N_R-3J_*-1,\tau_q^{-1},\ell_{t,q}^{-1}}
		\MM{N,N_{R,J_*}^{(1)}(M),\lambda_q,\ell_q^{-1}}\\
		&\leq\tau_q^{-1}\MM{M,N_t^{\mathrm N}-1,\tau_q^{-1},\ell_{t,q}^{-1}}
		\MM{N,N_{R,J_*}^{(1)}(M),\lambda_q,\ell_q^{-1}}.
	\end{aligned}
	\end{equation}
	Here the second inequality uses \eqref{eq:Newton-one-threshold-weight-shift} with $K=N_R-3J_*$.
	The last inequality follows from \eqref{eq:Newton-scale-separation}
	and \eqref{def:Nash-terminal-material-threshold}.

	We now pass from $\DTL$ to $\widetilde D_{t,q}$, first for $\nabla\bar u_{q,J_*}$, then for $\bar u_{q,J_*},u_{q,J_*}$, $\rho_{q,J_*},\bar\rho_{q,J_*},\rho_{q,J_*}-\rho_q$, and $m_{q,J_*}$.
	Apply Lemma~\ref{lem:source-material-derivative-change} with $\upsilon=\overline u_q$, $\omega=\bar u_{q,J_*}-\overline u_q$, and
	$$
	\begin{aligned}
		\mcI&=\R,\quad J=N_{\mathrm{ind}},\quad \mC_\omega=\delta_q^{\frac12},\quad \lambda=\lambda_q,\quad \tilde\lambda=\ell_q^{-1},\quad \mu=\tau_q^{-1},\quad \tilde\mu=\ell_{t,q}^{-1}.
	\end{aligned}
	$$
	Enlarge $C$ to cover the implicit constants in
	\eqref{eq:Newton-terminal-density-material-estimate}--\eqref{eq:Newton-terminal-momentum-material-estimate},
	\eqref{eq:Nash-terminal-velocity-input-estimate},
	\eqref{eq:Nash-mollified-velocity-background-estimate}--\eqref{eq:Nash-mollified-density-background-estimate},
	and \eqref{eq:Nash-gradient-threshold-shift}, over the stated derivative ranges.
	Enlarge $\widetilde{\mC}_0$ to cover the estimate for $\DTL^{M_0}\nabla\overline u_q$ in \eqref{est:uniform-grad-om-over-orho}.
	The coefficient assumptions require only $M_0<J=N_{\mathrm{ind}}$.
	They follow from \eqref{est:uniform-grad-om-over-orho} and \eqref{eq:Nash-drift-background-difference-two-scale-estimate},
	using $N_t\leq N_R-3J_*\leq N_{\mathrm{ind}}-1$.
	The scale conditions follow from \eqref{eq:material-mollification-spatial-comparisons} and \eqref{eq:Newton-scale-separation}.
	
	First take $F=\nabla\bar u_{q,J_*}$, with
	$\mC_F=C\Lambda_q$, $N_t=N_R-3J_*-1$, and $K=N_{\mathrm{max}}-1$.
	Estimate~\eqref{eq:Nash-gradient-threshold-shift} verifies the source assumption, including $M=N=0$.
	Lemma~\ref{lem:source-material-derivative-change} gives, for $0\leq M\leq N_{\mathrm{ind}}$ and $M+N\leq N_{\mathrm{max}}-1$,
	\begin{equation}\label{eq:Nash-terminal-Newton-velocity-sharp-gradient}
	\begin{aligned}
		\|\widetilde D_{t,q}^M\nabla\bar u_{q,J_*}\|_N
		\lesssim\Lambda_q\MM{M,N_R-3J_*-1,\tau_q^{-1},\ell_{t,q}^{-1}}
		\MM{N,N_{R,J_*}^{(1)}(M),\lambda_q,\ell_q^{-1}}.
	\end{aligned}
	\end{equation}
	The last inequality in \eqref{eq:Nash-gradient-threshold-shift} then gives \eqref{eq:Nash-terminal-Newton-velocity-integer-estimate}.
	Next take $F=\bar u_{q,J_*}$, with
	$\mC_F=C\delta_q^{\frac12}$, $N_t=N_R-3J_*$, and $K=N_{\mathrm{max}}$.
	Estimate~\eqref{eq:Nash-mollified-velocity-background-estimate} verifies the source assumption at positive total orders.
	Estimate~\eqref{eq:source-material-derivative-change} proves \eqref{eq:Nash-mollified-terminal-velocity-integer-estimate}.
	Spatial interpolation at orders $N$ and $N+1$ in
	\eqref{eq:Nash-mollified-terminal-velocity-integer-estimate} and
	\eqref{eq:Nash-terminal-Newton-velocity-integer-estimate} gives
	\eqref{eq:Nash-mollified-terminal-velocity-estimate} and
	\eqref{eq:Nash-terminal-Newton-velocity-estimate}, respectively.
	With $F=u_{q,J_*}$ and $K=N_{\mathrm{ind}}$, the same lemma and
	\eqref{eq:Nash-terminal-velocity-input-estimate} give, for $1\leq M+N\leq N_{\mathrm{ind}}$,
	\begin{equation}\label{eq:Nash-terminal-velocity-material-estimate}
	\begin{aligned}
		\|\widetilde D_{t,q}^Mu_{q,J_*}\|_N
		&\lesssim\delta_q^{\frac12}
		\MM{M,N_R-3J_*,\tau_q^{-1},\ell_{t,q}^{-1}}
		\MM{N,N_{R,J_*}^{(0)}(M),\lambda_q,\ell_q^{-1}}.
	\end{aligned}
	\end{equation}

	For $F=\rho_{q,J_*}$, apply Lemma~\ref{lem:source-material-derivative-change} with $K=N_{\mathrm{ind}}$ and $\mC_F=C\delta_q$.
	Estimate~\eqref{eq:Newton-terminal-density-material-estimate} verifies the source assumption at positive total orders.
	We obtain \eqref{eq:Nash-density-coefficient-material-estimate} for $s=1$.
	For fixed $s\in\mathbb R$ and $1\leq M+N\leq N_{\mathrm{ind}}$, apply \eqref{eq:faa-di-bruno}, as in the proof of \eqref{eq:Newton-flux-coefficient-material-estimate}, now with $\widetilde D_{t,q}$.
	For positive spatial orders, use the second inequality in \eqref{est:holder-norm-composition} for the resulting functions of the density.
	For $M_\nu,N_\nu\in\mathbb N_0$, we have
	$$
	\begin{aligned}
	    \Nrm{\nabla^N\widetilde D_{t,q}^M[(\rho_{q,J_*})^s]}_0
	    &\lesssim_s\sum_{k=1}^{M+N}
	    \sum_{\substack{M_1+\cdots+M_k=M\\N_1+\cdots+N_k=N\\M_\nu+N_\nu\geq1}}
	    \prod_{\nu=1}^k\Nrm{\widetilde D_{t,q}^{M_\nu}\rho_{q,J_*}}_{N_\nu}
	    \lesssim_s\delta_q\MM{M,N_R-3J_*,\tau_q^{-1},\ell_{t,q}^{-1}}
	    \MM{N,N_{R,J_*}^{(0)}(M),\lambda_q,\ell_q^{-1}}.
	\end{aligned}
	$$
	Here we used the density bounds in \eqref{eq:Newton-state-zero-order-bounds}, \eqref{eq:Newton-quadratic-weight-product}, and \eqref{eq:Nash-density-coefficient-material-estimate} with $s=1$.
	For $M=0<N$, the $C^0$ term is covered by \eqref{eq:Newton-state-zero-order-bounds} and $\delta_q\lambda_q\geq1$ from \eqref{def:lambda-delta}.
	The general estimate \eqref{eq:Nash-density-coefficient-material-estimate} follows.
	With $F=\bar\rho_{q,J_*}$, $K=N_{\mathrm{max}}$, and $\mC_F=C\delta_q$, Lemma~\ref{lem:source-material-derivative-change} gives \eqref{eq:Nash-mollified-density-material-estimate} for $s=1$.
	Repeating the calculation above with \eqref{eq:Nash-mollified-terminal-zero-order-bounds} proves \eqref{eq:Nash-mollified-density-material-estimate} for each fixed $s\in\mathbb R$.
	For $M\leq N_{\mathrm{ind}}$ and $M+N\leq N_{\mathrm{max}}$, combining this estimate with the zero-order bounds and lowering the thresholds to
	$N_t^{\mathrm N}$ and $N_{R,J_*}^{(1)}(M)$ gives
	\begin{equation}\label{eq:Nash-mollified-density-coefficient-estimate}
	\begin{aligned}
	\Nrm{\widetilde D_{t,q}^M[(\bar\rho_{q,J_*})^s]}_N
	&\lesssim_s\MM{M,N_t^{\mathrm N},\tau_q^{-1},\ell_{t,q}^{-1}}
	\MM{N,N_{R,J_*}^{(1)}(M),\lambda_q,\ell_q^{-1}}.
	\end{aligned}
	\end{equation}
	    Apply Lemma~\ref{lem:Newton-two-scale-material-change} with
    $\upsilon=\bar u_{q,J_*}$, $F=\rho_{q,J_*}-\rho_q$, $J=N_{\mathrm{ind}}$,
    $K=K_{J_*}^{(0)}$, $C_0=\widetilde{\mC}_0$, and $A=CA_{q,0}\ell_q^{-\alpha/2}$,
    where $C$ covers the implicit constant in \eqref{eq:Newton-terminal-density-increment-material}.
    The velocity-difference and source assumptions follow from \eqref{eq:Nash-drift-background-difference-two-scale-estimate}
    and \eqref{eq:Newton-terminal-density-increment-material}, respectively.
    We obtain \eqref{eq:Nash-terminal-Newton-density-increment-material}.

	Finally take $F=m_{q,J_*}$, $N_t=N_R-3J_*$, $K=N_{\mathrm{ind}}$, and $\mC_F=C\delta_q^{\frac12}$ in Lemma~\ref{lem:source-material-derivative-change}.
	The source assumption follows from \eqref{eq:Newton-terminal-momentum-material-estimate}.
	Thus \eqref{eq:Nash-terminal-momentum-material-estimate} holds.

	We next estimate $u_{q,J_*}-\bar u_{q,J_*}$ and $\rho_{q,J_*}-\bar\rho_{q,J_*}$.
	We then estimate the momentum perturbation $m_{q,J_*}-m_q$ in the Newton step.
	For $1\leq L\leq N_{\mathrm{ind}}$, apply
	\eqref{est:appendix-two-scale-material-mollification} in Proposition~\ref{prop:appendix-material-mollification} with $\upsilon=\bar u_{q,J_*}$.
	Take $f=u_{q,J_*}$ with $\mC_f=\delta_q^{\frac12}$ or $f=\rho_{q,J_*}$ with $\mC_f=\delta_q$, and set
	$$
	\begin{aligned}
		\ell&=\ell_q,\quad K^*=N_{\mathrm{ind}},\quad K_\upsilon=K_f=N_R-3J_*,\quad \mC_\upsilon=C\delta_q^{\frac12},\quad \lambda=\lambda_q,\quad \mu=\tau_q^{-1},\quad \tilde\mu=\ell_{t,q}^{-1}.
	\end{aligned}
	$$
	The scale conditions follow from \eqref{eq:material-mollification-spatial-comparisons} and \eqref{eq:Newton-scale-separation} for sufficiently large $a$.
	Use \eqref{eq:Nash-terminal-Newton-velocity-sharp-gradient} for $\widetilde D_{t,q}^M\nabla\bar u_{q,J_*}$.
	For $0\leq M\leq N_{\mathrm{ind}}$ and $M+N\leq N_{\mathrm{max}}-1$, we also have
	$$
	\begin{aligned}
		\|\nabla\widetilde D_{t,q}^M\bar u_{q,J_*}\|_N
		\leq\|\widetilde D_{t,q}^M\bar u_{q,J_*}\|_{N+1}
		&\lesssim\delta_q^{\frac12}\MM{M,N_R-3J_*,\tau_q^{-1},\ell_{t,q}^{-1}}
		\MM{N+1,N_{R,J_*}^{(0)}(M),\lambda_q,\ell_q^{-1}}\\
		&\leq\Lambda_q\MM{M,N_R-3J_*-1,\tau_q^{-1},\ell_{t,q}^{-1}}
		\MM{N,N_{R,J_*}^{(1)}(M),\lambda_q,\ell_q^{-1}}.
	\end{aligned}
	$$
	Here we used \eqref{eq:Newton-one-threshold-weight-shift} with $K=N_R-3J_*$ and \eqref{eq:Nash-mollified-terminal-velocity-integer-estimate}.
	These estimates verify the velocity assumption in \eqref{ass:appendix-two-scale-material-mollification}.
	For $f=u_{q,J_*}$, use \eqref{eq:Nash-terminal-velocity-material-estimate}.
	For $f=\rho_{q,J_*}$, use \eqref{eq:Nash-density-coefficient-material-estimate} with $s=1$.
	Since $L\leq N_{\mathrm{mom}}=N_{\mathrm{ind}}$, \eqref{est:appendix-two-scale-material-mollification} applies for $M+N+L\leq N_{\mathrm{ind}}$.
	Lowering the spatial threshold to $N_{R,J_*}^{(1)}(M)$ gives
	\eqref{eq:Nash-drift-mollification-difference}--%
	\eqref{eq:Nash-density-mollification-difference}, since $N_t^{\mathrm N}+4=N_R-3J_*$.
	For $L=0$ and $1\leq M+N\leq N_{\mathrm{ind}}$, use the velocity estimates
	\eqref{eq:Nash-mollified-terminal-velocity-integer-estimate} and \eqref{eq:Nash-terminal-velocity-material-estimate},
	and the density estimates   \eqref{eq:Nash-density-coefficient-material-estimate} and \eqref{eq:Nash-mollified-density-material-estimate} with $s=1$.
	At $M=N=0$, use the estimates with $L=1$.

	For the momentum perturbation in the Newton step, \eqref{eq:Newton-accumulated-state-increments-integer} with $j=J_*$ gives
	\begin{equation}\label{eq:Nash-Newton-momentum-increment-input-estimate}
	\begin{aligned}
		\nrm{\DTL^M(m_{q,J_*}-m_q)}_N
		&\lesssim A_{q,0}\delta_q^{-\frac12}\lambda_{q+1}^{-3\alpha}
		\MM{M,N_R-3J_*,\tau_q^{-1},\ell_{t,q}^{-1}}
		\MM{N,N_{R,J_*}^{(0)}(M),\lambda_q,\ell_q^{-1}}.
	\end{aligned}
	\end{equation}
	This estimate holds for $M,N\in\mathbb N_0$ with $M\leq N_{\mathrm{ind}}$, $N\leq K_{J_*}^{(0)}$ if $M=0$, and $M+N\leq K_{J_*}^{(0)}+1$ if $M\geq1$.
	Apply Lemma~\ref{lem:Newton-two-scale-material-change} with
	$\upsilon=\bar u_{q,J_*}$, $F=m_{q,J_*}-m_q$, $J=N_{\mathrm{ind}}$, $K=K_{J_*}^{(0)}$, and $C_0=\widetilde{\mC}_0$.
	Take $A=CA_{q,0}\delta_q^{-\frac12}\lambda_{q+1}^{-3\alpha}$, with $C$ enlarged to cover the implicit constant in \eqref{eq:Nash-Newton-momentum-increment-input-estimate}.
	The velocity-difference and source assumptions follow from \eqref{eq:Nash-drift-background-difference-two-scale-estimate}
	and \eqref{eq:Nash-Newton-momentum-increment-input-estimate}, respectively.
	We obtain \eqref{eq:Nash-terminal-Newton-momentum-increment-material}.

	We now estimate the stress tensors, starting with $S_{q,J_*}$.
	We then estimate $P_{q+1,j}$ and $R_{q,J_*}$ at the derivative scales in the Newton step.
	Take $r=\lceil\log_2(3/\beta)\rceil$.
	Estimate~\eqref{eq:Newton-accumulated-gluing-stress-background} with $j=J_*$ and exponent $\alpha_r$ holds for $M\leq N_{\mathrm{ind}}$ and $M+N\leq K_{J_*}^{(0)}+1$.
	Use $\Xi_{\alpha_r}^2\leq\ell_q^{-2\alpha_r}$ and lower the thresholds to $N_R-3J_*$ and $N_{R,J_*}^{(0)}(M)$.
	Apply Lemma~\ref{lem:Newton-two-scale-material-change} with $\upsilon=\bar u_{q,J_*}$, $F=S_{q,J_*}$, and
	$$
	\begin{aligned}
		J&=N_{\mathrm{ind}},\quad K=K_{J_*}^{(0)}+1,\quad C_0=\widetilde{\mC}_0,\quad A=CA_{q,0}(\tau_q\ell_q^{-1})^{\alpha_r}\ell_q^{-2\alpha_r}.
	\end{aligned}
	$$
	The assumption on $\bar u_{q,J_*}-\overline u_q$ follows from \eqref{eq:Nash-drift-background-difference-two-scale-estimate}, which covers total orders through $K_{J_*}^{(0)}$.
	The choice of $r$, $b<5/3$, and \eqref{def:admissible-alpha-losses} give
	$
	\alpha_r\left(\frac{3b+5}{4}+\beta-4b\alpha\right)
	<\alpha\left(\frac{b-1}{4}+\beta-4b\alpha\right).
	$
	By \eqref{def:parameter-mollification} and \eqref{def:tau-q}, we have $\tau_q\ell_q^{-1}\sim\lambda_q^{\frac{b-1}{4}+\beta-4b\alpha}$.
	Hence $(\tau_q\ell_q^{-1})^{\alpha_r-\alpha}\ell_q^{-2\alpha_r}$ is a negative power of $\lambda_q$, up to a fixed factor.
	More precisely, $\alpha_r\leq\alpha\beta/3$ and $4b\alpha<(b-1)/6$ give
	$$
	\begin{aligned}
	(\tau_q\ell_q^{-1})^{\alpha_r-\alpha}\ell_q^{-2\alpha_r}
	\lesssim\lambda_q^{\alpha_r\left(\frac{3b+5}{4}+\beta-4b\alpha\right)-\alpha\left(\frac{b-1}{4}+\beta-4b\alpha\right)}
	\leq\lambda_q^{-\frac29\alpha\beta}.
	\end{aligned}
	$$
	After increasing $a_{\mathrm{Nash}}$ to absorb the fixed implicit constant, we obtain, throughout the stated derivative range,
	\begin{equation}\label{eq:Nash-terminal-S-material-reserve}
	\begin{aligned}
		\|\widetilde D_{t,q}^MS_{q,J_*}\|_N
		&\lesssim_{\beta,\mathfrak M,C_*}A_{q,0}(\tau_q\ell_q^{-1})^{\alpha_r}\ell_q^{-2\alpha_r}
		\MM{M,N_R-3J_*,\tau_q^{-1},\ell_{t,q}^{-1}}\MM{N,N_{R,J_*}^{(0)}(M),\lambda_q,\ell_q^{-1}}\\
		&\lesssim A_{q,0}(\tau_q\ell_q^{-1})^\alpha
		\MM{M,N_R-3J_*,\tau_q^{-1},\ell_{t,q}^{-1}}\MM{N,N_{R,J_*}^{(0)}(M),\lambda_q,\ell_q^{-1}}.
	\end{aligned}
	\end{equation}
	In the first inequality, the implicit constant may depend on $\beta$ through $\alpha_r$.
	This dependence is absorbed into $a_{\mathrm{Nash}}$, so \eqref{eq:Nash-terminal-S-material-estimate} holds with an implicit constant depending only on $N_{\mathrm{max}}$ and $\alpha$.
	
	For $F=P_{q+1,j}$, $0\leq j\leq J_*$, or $F=R_{q,J_*}$, apply Lemma~\ref{lem:source-material-derivative-change}
	with $\upsilon=\overline u_q$, $\omega=\bar u_{q,J_*}-\overline u_q$, and
	$$
	\begin{aligned}
		J&=K=N_R,\quad \lambda=\tilde\lambda=\ell_q^{-1},\quad \mu=\tilde\mu=\ell_{t,q}^{-1},\quad N_t=0,\quad \mC_\omega=\mu/\lambda.
	\end{aligned}
	$$
	For $\nabla\upsilon$ and $\omega$, use \eqref{est:uniform-grad-om-over-orho} and \eqref{eq:Nash-drift-background-difference-two-scale-estimate}, respectively. We compare the scales using \eqref{eq:material-mollification-spatial-comparisons}.
	These verify the first two estimates in \eqref{ass:source-two-scale-material-change} through total orders $N_R-1$ and $N_R-2$.
	Enlarge $\widetilde{\mC}_0$ to cover both implicit constants.
	Choose $\mC_F$ as a sufficiently large fixed multiple of
	$
	A_{q,0}^2\delta_q^{-1}\lambda_{q+1}^{-6\alpha}+A_{q,J_*}\ell_q^{-\alpha}.
	$
	The source assumption holds for $M+N\leq N_R$ by \eqref{eq:Newton-terminal-interior-stress-CN} and \eqref{eq:Newton-substep-perturbative-tensor-detailed-estimate}, with \eqref{eq:Newton-two-scale-weight-comparison} for $R_{q,J_*}$.
	By \eqref{def:lambda-delta}, the last inequality in \eqref{def:admissible-alpha-losses}, and the definition of $A_{q,0}$, we have, for sufficiently large $a$,
	\begin{equation}\label{eq:Newton-quadratic-amplitude-ratio}
	\begin{aligned}
		\frac{A_{q,0}^2\delta_q^{-1}\lambda_{q+1}^{-6\alpha}}
		{\lambda_{q+1}^{-8\alpha}\delta_{q+2}}
		=C_*^2\lambda_q^{-16\alpha+2\beta}
		\lambda_{q+1}^{2\alpha-4\beta}\lambda_{q+2}^{2\beta}
		\lesssim\lambda_q^{-2\left((8-b)\alpha-\beta(b-1)^2\right)}\leq\lambda_q^{-2\alpha}.
	\end{aligned}
	\end{equation}
	Taking the supremum over $j$ and adding the estimate for $R_{q,J_*}$, we obtain, for $M+N\leq N_R$,
	$$
	\begin{aligned}
		\sup_{0\leq j\leq J_*}\|\widetilde D_{t,q}^MP_{q+1,j}\|_N+\|\widetilde D_{t,q}^MR_{q,J_*}\|_N
		\leq C\left[A_{q,0}^2\delta_q^{-1}\lambda_{q+1}^{-6\alpha}+A_{q,J_*}\ell_q^{-\alpha}\right]\ell_{t,q}^{-M}\ell_q^{-N}
		\leq\lambda_q^{-\alpha}\lambda_{q+1}^{-8\alpha}\delta_{q+2}\ell_{t,q}^{-M}\ell_q^{-N}.
	\end{aligned}
	$$
	In the last inequality, we used \eqref{def:admissible-alpha-losses}, \eqref{eq:Newton-terminal-amplitude-ratio}, and \eqref{eq:Newton-quadratic-amplitude-ratio} after increasing $a_{\mathrm{Nash}}$.
	This proves \eqref{eq:Nash-retained-Newton-stresses-material}.
\end{proof}

The following lemma estimates ordered mixed spatial and material derivatives.
\begin{lem}\label{lem:Nash-coefficient-mixed-derivatives}
	Assume the hypotheses of Proposition~\ref{prop:Nash-terminal-Newton-velocity-estimates}.
	Let $\mcI$ be a time interval.
	Let $f$ be a smooth scalar, vector, or matrix field on $\mcI\times\mathbb T^3$ or $\mcI\times\mathbb R^3$, with finite norms.
	Fix $s\in\{0,1\}$ and integers $0\leq M\leq N_{\mathrm{ind}}+1$ and $M\leq N_*\leq N_{\mathrm{max}}$, with $N_*\geq1$.
	Fix $\mC_f\geq0$. Suppose that
	$$
		\|\widetilde D_{t,q}^{M_0}f\|_{N,\mcI}
		\lesssim\mC_f\MM{M_0,N_t^{\mathrm N},\tau_q^{-1},\ell_{t,q}^{-1}}
		\MM{N,N_{R,J_*}^{(s)}(M_0),\lambda_q,\ell_q^{-1}}
	$$
	for integers $0\leq M_0\leq M$ and $N\geq0$ with $1\leq M_0+N\leq N_*$.
	For $k\geq1$, write $\vartheta=(\vartheta_1,\ldots,\vartheta_k)$ and $\sigma=(\sigma_1,\ldots,\sigma_k)$, with $\vartheta,\sigma\in\mathbb N_0^k$.
	Set $|\vartheta|:=\sum_{i=1}^k\vartheta_i$ and $|\sigma|:=\sum_{i=1}^k\sigma_i$.
	For $N\geq0$ with $|\sigma|\leq M$ and $1\leq N+|\vartheta|+|\sigma|\leq N_*$, we have
	\begin{equation}\label{eq:Nash-coefficient-mixed-integer-estimate}
	\begin{aligned}
		\left\|\left(\prod_{i=1}^k\nabla^{\vartheta_i}\widetilde D_{t,q}^{\sigma_i}\right)f\right\|_{N,\mcI}
		&\lesssim\mC_f\MM{|\sigma|,N_t^{\mathrm N},\tau_q^{-1},\ell_{t,q}^{-1}}
		\MM{N+|\vartheta|,N_{R,J_*}^{(s)}(|\sigma|),\lambda_q,\ell_q^{-1}}.
	\end{aligned}
	\end{equation}
	The implicit constants depend only on $N_{\mathrm{max}}$, $\alpha$, and the implicit constants in the hypothesis and \eqref{eq:Nash-terminal-Newton-velocity-integer-estimate}.
\end{lem}

\begin{proof}
	For $|\vartheta|+|\sigma|=0$, the integer estimate is the hypothesis with $M_0=0$.
	Otherwise, apply \eqref{eq:C0-iterate} in Lemma~\ref{lem:holder-cooper-C0} componentwise. Use the same $N_*$ and $\mcI$, and take $\upsilon=\bar u_{q,J_*}$, $M_*=|\sigma|$, and
	$$
	\begin{alignedat}{4}
		N_\upsilon&=N_{R,J_*}^{(1)}(|\sigma|)+1,&\quad N_f&=\widetilde N=N_{R,J_*}^{(s)}(|\sigma|),&\quad M_\upsilon&=N_t^{\mathrm N}-1,\quad M_f=\widetilde M=N_t^{\mathrm N},&\quad\mC_\upsilon&=\tau_q^{-1}\lambda_q^{-1},\\[2pt]
		\lambda_\upsilon&=\lambda_f=\lambda_q,&\quad\tilde\lambda_\upsilon&=\tilde\lambda_f=\ell_q^{-1},&\quad\mu_\upsilon&=\mu_f=\mu_*=\tau_q^{-1},&\quad\tilde\mu_\upsilon&=\tilde\mu_f=\tilde\mu_*=\ell_{t,q}^{-1}.
	\end{alignedat}
	$$
	The velocity estimate \eqref{eq:Nash-terminal-Newton-velocity-integer-estimate} and the hypothesis on $f$ verify the derivative assumptions after lowering the spatial thresholds.
	The required velocity indices satisfy $M_0<|\sigma|\leq N_{\mathrm{ind}}+1$ and $M_0+N_0\leq N_*-1\leq N_{\mathrm{max}}-1$.
	Applying \eqref{eq:C0-iterate} with up to $N$ additional spatial derivatives on the left proves \eqref{eq:Nash-coefficient-mixed-integer-estimate}.
\end{proof}

Proposition~\ref{prop:Nash-terminal-Newton-velocity-estimates} and Lemma~\ref{lem:Nash-coefficient-mixed-derivatives} give the following estimates for the mollified density and the gluing stress.
\begin{co}\label{cor:Nash-density-stress-mixed-derivatives}
	Assume the hypotheses of Proposition~\ref{prop:Nash-terminal-Newton-velocity-estimates}.
	Let $\mcI$ be a time interval. Fix $k\geq1$ and $\vartheta,\sigma\in\mathbb N_0^k$.
	Take $N\in\mathbb N_0$ with $|\sigma|\leq N_{\mathrm{ind}}$ and $N+|\vartheta|+|\sigma|\geq1$.
	For $N+|\vartheta|+|\sigma|\leq N_{\mathrm{max}}$, the mollified density satisfies
	\begin{equation}\label{eq:Nash-mollified-density-mixed-integer-estimate}
	\begin{aligned}
		\Nrm{\left(\prod_{i=1}^k\nabla^{\vartheta_i}\widetilde D_{t,q}^{\sigma_i}\right)\bar\rho_{q,J_*}}_{N,\mcI}
		&\lesssim\delta_q\MM{|\sigma|,N_t^{\mathrm N},\tau_q^{-1},\ell_{t,q}^{-1}}\MM{N+|\vartheta|,N_{R,J_*}^{(0)}(|\sigma|),\lambda_q,\ell_q^{-1}}.
	\end{aligned}
	\end{equation}
	For $N+|\vartheta|+|\sigma|\leq K_{J_*}^{(0)}+1$, the gluing stress satisfies
	\begin{equation}\label{eq:Nash-terminal-S-mixed-integer-estimate}
	\begin{aligned}
		\Nrm{\left(\prod_{i=1}^k\nabla^{\vartheta_i}\widetilde D_{t,q}^{\sigma_i}\right)S_{q,J_*}}_{N,\mcI}
		&\lesssim A_{q,0}(\tau_q\ell_q^{-1})^\alpha
		\MM{|\sigma|,N_t^{\mathrm N},\tau_q^{-1},\ell_{t,q}^{-1}}\MM{N+|\vartheta|,N_{R,J_*}^{(0)}(|\sigma|),\lambda_q,\ell_q^{-1}}.
	\end{aligned}
	\end{equation}
	The implicit constants have the dependencies stated in Lemma~\ref{lem:Nash-coefficient-mixed-derivatives}.
\end{co}

\begin{proof}
	Apply Lemma~\ref{lem:Nash-coefficient-mixed-derivatives} with $s=0$ and $M=N_{\mathrm{ind}}$.
	The choices of $N_*$ below satisfy $N_{\mathrm{ind}}\leq N_*\leq N_{\mathrm{max}}$, since $K_{J_*}^{(0)}+1=N_{\mathrm{max}}-3J_*-1\geq N_{\mathrm{ind}}$.
	For $f=\bar\rho_{q,J_*}$, take $\mC_f=\delta_q$ and $N_*=N_{\mathrm{max}}$, using \eqref{eq:Nash-mollified-density-material-estimate} with power $1$.
	For $f=S_{q,J_*}$, take $\mC_f=A_{q,0}(\tau_q\ell_q^{-1})^\alpha$ and $N_*=K_{J_*}^{(0)}+1$, using \eqref{eq:Nash-terminal-S-material-estimate}.
	Since $\ell_{t,q}\leq\tau_q$ and $N_t^{\mathrm N}=N_R-3J_*-4$, we have, for $M\in\mathbb N_0$,
	$$
	\MM{M,N_R-3J_*,\tau_q^{-1},\ell_{t,q}^{-1}}
	\leq\MM{M,N_t^{\mathrm N},\tau_q^{-1},\ell_{t,q}^{-1}}.
	$$
	Thus each cited estimate verifies the source hypothesis.
	The conclusions follow from \eqref{eq:Nash-coefficient-mixed-integer-estimate}.
\end{proof}

\subsubsection{Flow maps and space--time cutoffs}

The time cutoffs localize $S_{q,J_*}$ to disjoint intervals.
Within each interval, the spatial cutoffs are transported by the flow of $\bar u_{q,J_*}$.
For the material derivative estimates, we use $\tau_q$ and $\ell_{t,q}$ with the scale comparisons in
\eqref{eq:material-mollification-hierarchy}--\eqref{eq:material-mollification-spatial-comparisons}.

Set $\mathcal G_q:=\supp_tS_{q,J_*}$. Under \eqref{ind:R-support},
this set is compact by \eqref{eq:Newton-terminal-temporal-support-next-scale}.
By Proposition~\ref{prop:Newton-incompressible-flow-gluing-estimates}, $\mathcal G_q$ is contained in
$\bigcup_{u\in\mathbb Z}(\supp\partial_t\theta_u\cup\supp\partial_t\tilde\theta_u)$.
The cutoff definitions in Subsection~\ref{subsec:cutoffs} place these derivative supports in the union of shifted copies of a fixed compact subinterval of $(7\tau_q/24,17\tau_q/24)$.
The shifts are integer multiples of $\tau_q$.
For $u\in\mathbb Z$, set
$I_u:=\left(\left(u+\frac14\right)\tau_q,\left(u+\frac34\right)\tau_q\right)$.
Select the intervals $I_u$ that meet $\mathcal G_q$. There are finitely many since $\mathcal G_q$ is compact.
All sums over $u$ below are taken over these selected intervals. They satisfy
$$
|I_u|=\frac12\tau_q,
\qquad
\operatorname{dist}(I_u,I_{u'})\geq\frac12\tau_q
\quad\text{if }u\neq u',
$$
and $\mathcal G_q\subset\bigcup_u I_u$. They also satisfy
$\operatorname{dist}(\mathcal G_q\cap I_u,\partial I_u)\geq\tau_q/24$.
The separation from $\partial I_u$ allows us to choose
$\theta_u^{\mathrm N}\in C_c^\infty(I_u)$ with $0\leq\theta_u^{\mathrm N}\leq1$ so that
\begin{equation}
	\label{def:Nash-transition-cutoff}
	\begin{aligned}
		&\theta_u^{\mathrm N}=1\quad\text{when }\operatorname{dist}(t,\mathcal G_q\cap I_u)\leq\tau_q/48,
		\\
		&\supp\theta_u^{\mathrm N}\cap\supp\theta_{u'}^{\mathrm N}=\varnothing
		\quad\text{if }u\neq u',\\
		&\|\partial_t^M\theta_u^{\mathrm N}\|_0\lesssim_M\tau_q^{-M}.
	\end{aligned}
\end{equation}
For each selected $u$, set $t_u^{\mathrm N}:=\left(u+\frac12\right)\tau_q$.
Let $X_u^{\mathrm N}$ and $\Phi_u^{\mathrm N}$ be the forward and backward flows from
Subsection~\ref{subsec:forward-backward-flows}, with velocity $\bar u_{q,J_*}$ and base time $t_u^{\mathrm N}$.
We work on the interval
$
\mcI_u^{\mathrm N}
:=\left[t_u^{\mathrm N}-\frac34\tau_q,t_u^{\mathrm N}+\frac34\tau_q\right],
$
which contains $\mN_{\ell_{t,q}}(I_u)$ since $\ell_{t,q}\ll\tau_q$.
Similarly, we define
$J_u^{\mathrm N}:=\det\nabla\Phi_u^{\mathrm N}$ and $\mathbf C_u^{\mathrm N}:=\cof\nabla\Phi_u^{\mathrm N}$.

We define the spatial cutoffs using the function $\chi_0$ from
Subsection~\ref{subsec:cutoffs}. It satisfies
$
	\supp\chi_0\subset Q\left(0,\frac{9\pi}{8}\right),
$  and $
	\sum_{\kappa\in\mathbb Z^3}
	\bigl(\chi_0(z-2\pi \kappa)\bigr)^2=1.
$
We define, for $\kappa\in\mathbb Z^3$,
\begin{equation}
	\label{def:Nash-spatial-cutoffs}
	\chi_\kappa^{\mathrm N}(z)
	:=\chi_0(\lambda_q z-2\pi \kappa).
\end{equation}
The family is locally finite on $\mathbb R^3$ and satisfies
\begin{equation}
	\label{eq:Nash-spatial-cutoff-properties}
	\sum_{\kappa\in\mathbb Z^3}(\chi_\kappa^{\mathrm N})^2=1,
	\qquad
	\|\nabla^N\chi_\kappa^{\mathrm N}\|_0
	\lesssim_N\lambda_q^N.
\end{equation}
The resulting space--time cutoffs are
$
\theta_u^{\mathrm N}(t)\chi_\kappa^{\mathrm N}(\Phi_u^{\mathrm N}(t,x)).
$
We localize $S_{q,J_*}$ to these intervals by setting
$$
S_{q,J_*,u}
:=\theta_u^{\mathrm N}S_{q,J_*},
\qquad
S_{q,J_*}
=\sum_uS_{q,J_*,u}.
$$
By \eqref{def:Nash-transition-cutoff}, every positive-order derivative of $\theta_u^{\mathrm N}$ vanishes on $\mathcal G_q$.
Hence $\widetilde D_{t,q}^MS_{q,J_*,u}=\theta_u^{\mathrm N}\widetilde D_{t,q}^MS_{q,J_*}$.
For $M\leq N_{\mathrm{ind}}$ and $M+N\leq K_{J_*}^{(0)}+1$, we have
\begin{equation}\label{eq:Nash-localized-S-material-estimate}
\begin{aligned}
	\|\widetilde D_{t,q}^MS_{q,J_*,u}\|_N
	&\leq\|\widetilde D_{t,q}^MS_{q,J_*}\|_N\\
	&\lesssim A_{q,0}(\tau_q\ell_q^{-1})^\alpha
	\MM{M,N_R-3J_*,\tau_q^{-1},\ell_{t,q}^{-1}}
	\MM{N,N_{R,J_*}^{(0)}(M),\lambda_q,\ell_q^{-1}}.
\end{aligned}
\end{equation}
Here we used \eqref{eq:Nash-terminal-S-material-estimate}, whose implicit constant is independent of
$\mathfrak M,C_*$ by Proposition~\ref{prop:Nash-terminal-Newton-velocity-estimates}.
By \eqref{eq:Nash-localized-S-material-estimate}, Lemma~\ref{lem:Nash-coefficient-mixed-derivatives} applies to $S_{q,J_*,u}$ with $s=0$, $\mC_f=A_{q,0}(\tau_q\ell_q^{-1})^\alpha$, and $N_*=K_{J_*}^{(0)}+1$.
Fix a time interval $\mcI$, $k\geq1$, and $\vartheta,\sigma\in\mathbb N_0^k$.
For $N\in\mathbb N_0$ with $|\sigma|\leq N_{\mathrm{ind}}$ and $1\leq N+|\vartheta|+|\sigma|\leq K_{J_*}^{(0)}+1$, we have
\begin{equation}\label{eq:Nash-localized-S-mixed-integer-estimate}
\begin{aligned}
	\Nrm{\left(\prod_{i=1}^k\nabla^{\vartheta_i}\widetilde D_{t,q}^{\sigma_i}\right)S_{q,J_*,u}}_{N,\mcI}
	&\lesssim A_{q,0}(\tau_q\ell_q^{-1})^\alpha
	\MM{|\sigma|,N_t^{\mathrm N},\tau_q^{-1},\ell_{t,q}^{-1}}\MM{N+|\vartheta|,N_{R,J_*}^{(0)}(|\sigma|),\lambda_q,\ell_q^{-1}}.
\end{aligned}
\end{equation}
To estimate the derivatives of these cutoffs, we first record the estimates for the flow maps.
\begin{prop}
	\label{prop:Nash-flow-estimates}
	Assume the hypotheses of Proposition~\ref{prop:Nash-terminal-Newton-velocity-estimates}.
	For every selected $u$, we have
	\begin{equation}
		\begin{aligned}
			\nrm{\nabla\Phi_u^{\mathrm N}-\Id}_{0,\mcI_u^{\mathrm N}}
			+\nrm{(\nabla\Phi_u^{\mathrm N})^{-1}-\Id}_{0,\mcI_u^{\mathrm N}}
			+\nrm{\nabla X_u^{\mathrm N}-\Id}_{0,\mcI_u^{\mathrm N}}
			+\nrm{\mathbf C_u^{\mathrm N}-\Id}_{0,\mcI_u^{\mathrm N}}
			+\nrm{(\mathbf C_u^{\mathrm N})^{-1}-\Id}_{0,\mcI_u^{\mathrm N}}
			&\lesssim\lambda_{q+1}^{-3\alpha},\\
			\|J_u^{\mathrm N}-1\|_{0,\mcI_u^{\mathrm N}}+\|(J_u^{\mathrm N})^{-1}-1\|_{0,\mcI_u^{\mathrm N}}&\lesssim\lambda_{q+1}^{-3\alpha}.
		\end{aligned}
		\label{eq:Nash-flow-zero-order-estimate}
	\end{equation}
	For $M,N\in\mathbb N_0$ with either $M=0$ and $0\leq N\leq N_{\mathrm{max}}-1$, or
	$1\leq M\leq N_{\mathrm{ind}}+1$ and $M+N\leq N_{\mathrm{max}}$, the flow gradients satisfy
	\begin{equation}
		\begin{aligned}
			&\nrm{\widetilde D_{t,q}^M\nabla\Phi_u^{\mathrm N}}_{N,\mcI_u^{\mathrm N}}
			+\nrm{\widetilde D_{t,q}^M(\nabla\Phi_u^{\mathrm N})^{-1}}_{N,\mcI_u^{\mathrm N}}
			+\nrm{\partial_t^M\nabla X_u^{\mathrm N}}_{N,\mcI_u^{\mathrm N}}\\
			&\qquad\lesssim
			\MM{M,N_t^{\mathrm N},\tau_q^{-1},\ell_{t,q}^{-1}}
			\MM{N,N_{R,J_*}^{(1)}(M),\lambda_q,\ell_q^{-1}}
			\leq\MM{M,N_t^{\mathrm N},\tau_q^{-1},\ell_{t,q}^{-1}}\ell_q^{-N}.
		\end{aligned}
		\label{eq:Nash-flow-higher-order-CN-estimate}
	\end{equation}
	In the same range, the determinant, cofactor, and their inverses satisfy
	\begin{equation}
		\begin{aligned}
			&\nrm{\widetilde D_{t,q}^M J_u^{\mathrm N}}_{N,\mcI_u^{\mathrm N}}
			+\nrm{\widetilde D_{t,q}^M (J_u^{\mathrm N})^{-1}}_{N,\mcI_u^{\mathrm N}}
			+\nrm{\widetilde D_{t,q}^M\mathbf C_u^{\mathrm N}}_{N,\mcI_u^{\mathrm N}}
			+\nrm{\widetilde D_{t,q}^M(\mathbf C_u^{\mathrm N})^{-1}}_{N,\mcI_u^{\mathrm N}}\\
			&\qquad\lesssim
			\MM{M,N_t^{\mathrm N},\tau_q^{-1},\ell_{t,q}^{-1}}
			\MM{N,N_{R,J_*}^{(1)}(M),\lambda_q,\ell_q^{-1}}
			\leq\MM{M,N_t^{\mathrm N},\tau_q^{-1},\ell_{t,q}^{-1}}\ell_q^{-N}.
		\end{aligned}
		\label{eq:Nash-flow-determinant-cofactor-CN-estimate}
	\end{equation}
	The implicit constants depend only on $N_{\mathrm{max}},\alpha$ and the implicit constant in
	\eqref{eq:Nash-terminal-Newton-velocity-sharp-gradient}.
\end{prop}

\begin{proof}
	Apply Lemma~\ref{lem:flow-estimates-upsilon} with $\upsilon=\bar u_{q,J_*}$ and
	$$
	\begin{aligned}
		t_0&=t_u^{\mathrm N},\quad \tau=\tfrac34\tau_q,\quad K^*=N_R-3J_*,\quad K^\sharp=N_{\mathrm{max}},\quad N_t=N_t^{\mathrm N},\\
		\lambda_\upsilon&=\lambda_q,\quad \tilde\lambda_\upsilon=\ell_q^{-1},\quad \mu_\upsilon=\tau_q^{-1},\quad \tilde\mu_\upsilon=\ell_{t,q}^{-1},\quad \mC_\upsilon=\mC_{\Div,\upsilon}=\delta_q^{\frac12}.
	\end{aligned}
	$$
	For the gradient hypotheses, use \eqref{eq:Nash-terminal-Newton-velocity-sharp-gradient} and $N_t^{\mathrm N}\leq K^*-1$.
	The required orders satisfy $M\leq K^*\leq N_{\mathrm{ind}}$ and $M+N\leq N_{\mathrm{max}}-1$.
	Spatial interpolation at orders $N$ and $N+1$ verifies the H{\"o}lder hypotheses for $M+N\leq N_{\mathrm{max}}-2$.
	Taking the trace gives the assumptions on $\Div\bar u_{q,J_*}$.
	Choose $C_{\mathrm{mol}}\ge1$ to cover the implicit constants in \eqref{eq:Nash-terminal-Newton-velocity-sharp-gradient} throughout its stated derivative range.
	By \eqref{def:tau-q} and \eqref{eq:material-mollification-spatial-comparisons}, we may choose $a_{\mathrm{Nash}}$ so that
	$$
	\begin{aligned}
		C_{\mathrm{mol}}\tau_q\Lambda_q&\leq C_{\mathrm{mol}}\lambda_{q+1}^{-4\alpha}\leq\lambda_{q+1}^{-3\alpha}\leq c_0,\quad
		\delta_q^{\frac12}\ell_{t,q}\ell_q^{-1}=\ell_q\lambda_q\leq1.
	\end{aligned}
	$$
	Thus Lemma~\ref{lem:flow-estimates-upsilon} applies on $\mcI_u^{\mathrm N}$, with
	$\mu_{\upsilon,*}=\tau_q^{-1}$ and $\tilde\mu_{\upsilon,*}=\ell_{t,q}^{-1}$.

	By \eqref{eq:flow-gradient-close-id} and \eqref{eq:flow-gradient-CN} with $N=0$, we have
	$$
	\|\nabla\Phi_u^{\mathrm N}-\Id\|_{0,\mcI_u^{\mathrm N}}+\|(\nabla\Phi_u^{\mathrm N})^{-1}-\Id\|_{0,\mcI_u^{\mathrm N}}+\|\nabla X_u^{\mathrm N}-\Id\|_{0,\mcI_u^{\mathrm N}}
	\lesssim\tau_q\Lambda_q\leq\lambda_{q+1}^{-3\alpha}.
	$$
	For $J_u^{\mathrm N}$, use \eqref{eq:flow-jacobian-close-one} and \eqref{eq:flow-jacobian-CN} with $N=0$. The identity $\mathbf C_u^{\mathrm N}=J_u^{\mathrm N}(\nabla\Phi_u^{\mathrm N})^{-T}$ gives the remaining $C^0$ estimates.
	At $M=0$, use \eqref{eq:flow-gradient-CN} for $N\leq N_{\mathrm{max}}-1$.
	For $1\leq M\leq K^*+1$, \eqref{eq:flow-gradient-material} gives \eqref{eq:Nash-flow-higher-order-CN-estimate} after lowering the spatial threshold to $N_{R,J_*}^{(1)}(M)$.
	Here we used $\delta_q^{\frac12}\lambda_q=\Lambda_q\leq\tau_q^{-1}$ for $M\leq K^*$, and $\delta_q^{\frac12}\ell_q^{-1}\leq\ell_{t,q}^{-1}$ for $M=K^*+1>N_t^{\mathrm N}$.

	It remains to cover $K^*+1<M\leq N_{\mathrm{ind}}+1$ with $M+N\leq N_{\mathrm{max}}$. In this range, $N_{R,J_*}^{(1)}(M)=0$.
	For $0\leq M_0\leq N_{\mathrm{ind}}$ and $M_0+N_0\leq N_{\mathrm{max}}-1$, lower the material threshold in \eqref{eq:Nash-terminal-Newton-velocity-sharp-gradient} to $N_t^{\mathrm N}-1$ and use $C_{\mathrm{mol}}\tau_q\Lambda_q\leq1$ to obtain
	\begin{equation}\label{eq:Nash-terminal-velocity-gradient-Nash-scale}
	\begin{aligned}
		\|\widetilde D_{t,q}^{M_0}\nabla\bar u_{q,J_*}\|_{N_0}
		\leq\tau_q^{-1}
		\MM{M_0,N_t^{\mathrm N}-1,\tau_q^{-1},\ell_{t,q}^{-1}}
		\MM{N_0,N_{R,J_*}^{(1)}(M_0),\lambda_q,\ell_q^{-1}}.
	\end{aligned}
	\end{equation}
	By \eqref{backward flow}, $\nabla\Phi_u^{\mathrm N}$ and $(\nabla\Phi_u^{\mathrm N})^{-1}$ satisfy
	\begin{equation}\label{eq:Nash-flow-gradient-transport}
	\widetilde D_{t,q}\nabla\Phi_u^{\mathrm N}
	=-(\nabla\Phi_u^{\mathrm N})\nabla \bar u_{q,J_*},
	\qquad
	\widetilde D_{t,q}(\nabla\Phi_u^{\mathrm N})^{-1}
	=\nabla \bar u_{q,J_*}(\nabla\Phi_u^{\mathrm N})^{-1}.
	\end{equation}
	Applying $\widetilde D_{t,q}^{M-1}$ to the first equation gives
	$$
	\widetilde D_{t,q}^M\nabla\Phi_u^{\mathrm N}
	=-\sum_{r=0}^{M-1}\binom{M-1}{r}
	(\widetilde D_{t,q}^r\nabla\Phi_u^{\mathrm N})
	(\widetilde D_{t,q}^{M-1-r}\nabla\bar u_{q,J_*}).
	$$
	Using \eqref{eq:Nash-terminal-velocity-gradient-Nash-scale}, induction on $M$ gives
	$$
	\begin{aligned}
	    \|\widetilde D_{t,q}^M\nabla\Phi_u^{\mathrm N}\|_{N,\mcI_u^{\mathrm N}}
	    \lesssim\sum_{r=0}^{M-1}
	    \MM{r,N_t^{\mathrm N},\tau_q^{-1},\ell_{t,q}^{-1}}
	    \tau_q^{-1}\MM{M-1-r,N_t^{\mathrm N}-1,\tau_q^{-1},\ell_{t,q}^{-1}}\ell_q^{-N}
	    \lesssim\MM{M,N_t^{\mathrm N},\tau_q^{-1},\ell_{t,q}^{-1}}\ell_q^{-N}.
	\end{aligned}
	$$
	For the last inequality, use \eqref{eq:Newton-quadratic-weight-product}.
	The second equation in \eqref{eq:Nash-flow-gradient-transport} gives the same estimate for $(\nabla\Phi_u^{\mathrm N})^{-1}$.
	For $N=0$, \eqref{eq:pullback-flow-gradient-material} gives the $C^0$ estimate for $\partial_t^M\nabla X_u^{\mathrm N}$ directly from that for $\widetilde D_{t,q}^M(\nabla\Phi_u^{\mathrm N})^{-1}$.
	For $N\geq1$, the second inequality in \eqref{est:holder-norm-composition}, the $C^0$ estimate above,
	and \eqref{eq:flow-gradient-CN} give
	$$
	\begin{aligned}
	    \|\partial_t^M\nabla X_u^{\mathrm N}\|_{N,\mcI_u^{\mathrm N}}
	    \lesssim\|\widetilde D_{t,q}^M(\nabla\Phi_u^{\mathrm N})^{-1}\|_{N,\mcI_u^{\mathrm N}}+\|\widetilde D_{t,q}^M(\nabla\Phi_u^{\mathrm N})^{-1}\|_{1,\mcI_u^{\mathrm N}}\|\nabla X_u^{\mathrm N}\|_{N-1,\mcI_u^{\mathrm N}}
	    \lesssim\MM{M,N_t^{\mathrm N},\tau_q^{-1},\ell_{t,q}^{-1}}\ell_q^{-N}.
	\end{aligned}
	$$

	The identities $J_u^{\mathrm N}=\det\nabla\Phi_u^{\mathrm N}$ and
	$(J_u^{\mathrm N})^{-1}=\det((\nabla\Phi_u^{\mathrm N})^{-1})$
	express both coefficients as cubic polynomials in the matrix entries.
	In the stated integer ranges, \eqref{eq:Nash-flow-higher-order-CN-estimate} and \eqref{eq:holder-product-tame} with $\alpha=0$
	give, for $M_\nu,N_\nu\in\mathbb N_0$,
	$$
	\begin{aligned}
	    \|\widetilde D_{t,q}^M J_u^{\mathrm N}\|_{N,\mcI_u^{\mathrm N}}
	    +\|\widetilde D_{t,q}^M(J_u^{\mathrm N})^{-1}\|_{N,\mcI_u^{\mathrm N}}
	    &\lesssim\sum_{\substack{M_1+M_2+M_3=M\\N_1+N_2+N_3\leq N}}
	    \prod_{\nu=1}^3\Bigl(
	    \MM{M_\nu,N_t^{\mathrm N},\tau_q^{-1},\ell_{t,q}^{-1}}
	    \MM{N_\nu,N_{R,J_*}^{(1)}(M_\nu),\lambda_q,\ell_q^{-1}}\Bigr)\\
	    &\lesssim\MM{M,N_t^{\mathrm N},\tau_q^{-1},\ell_{t,q}^{-1}}
	    \MM{N,N_{R,J_*}^{(1)}(M),\lambda_q,\ell_q^{-1}}.
	\end{aligned}
	$$
	For the second inequality, we used \eqref{eq:appendix-derivative-weight-product}.
	For the cofactors, use
	$
	\mathbf C_u^{\mathrm N}=J_u^{\mathrm N}(\nabla\Phi_u^{\mathrm N})^{-T},
	$ and $
	(\mathbf C_u^{\mathrm N})^{-1}=(J_u^{\mathrm N})^{-1}(\nabla\Phi_u^{\mathrm N})^T.
	$
	The preceding estimates and \eqref{eq:Nash-flow-higher-order-CN-estimate} give
	\eqref{eq:Nash-flow-determinant-cofactor-CN-estimate}, using \eqref{eq:Nash-flow-zero-order-estimate} at order zero.
\end{proof}

Using the flow estimates above, we now estimate mixed spatial and material derivatives
of the transported cutoffs and flow coefficients.

\begin{lem}\label{lem:Nash-transported-cutoff-estimates}
Let $\theta_u^{\mathrm N}$ and $\chi_\kappa^{\mathrm N}$ be the cutoffs in
\eqref{def:Nash-transition-cutoff}--\eqref{def:Nash-spatial-cutoffs}.
Let $\Phi_u^{\mathrm N}$ be the backward flow from Proposition~\ref{prop:Nash-flow-estimates}.
For $M\in\mathbb N_0$, we have
\begin{equation}\label{eq:Nash-transported-cutoff-material-identity}
\begin{aligned}
    \widetilde D_{t,q}^M\bigl(\theta_u^{\mathrm N}\chi_\kappa^{\mathrm N}(\Phi_u^{\mathrm N})\bigr)
    &=(\partial_t^M\theta_u^{\mathrm N})\chi_\kappa^{\mathrm N}(\Phi_u^{\mathrm N}).
\end{aligned}
\end{equation}
All cutoff norms are taken over $\mcI_u^{\mathrm N}$.
For $0\leq M\leq N_{\mathrm{ind}}+1$ and $N\in\mathbb N_0$ with $M+N\leq N_{\mathrm{max}}$, we have
\begin{equation}
	\label{eq:Nash-transported-cutoff-estimate}
	\begin{aligned}
		\nrm{\widetilde D_{t,q}^M
			\bigl(\theta_u^{\mathrm N}\chi_\kappa^{\mathrm N}(\Phi_u^{\mathrm N})\bigr)}_{N,\mcI_u^{\mathrm N}}
		\lesssim\MM{M,N_t^{\mathrm N},\tau_q^{-1},\ell_{t,q}^{-1}}\MM{N,N_{R,J_*}^{(1)}(M),\lambda_q,\ell_q^{-1}}
		\leq\MM{M,N_t^{\mathrm N},\tau_q^{-1},\ell_{t,q}^{-1}}\ell_q^{-N}.
	\end{aligned}
\end{equation}
For $k\geq1$ and $\vartheta,\sigma\in\mathbb N_0^k$ with $|\sigma|\leq N_{\mathrm{ind}}+1$ and $1\leq N+|\vartheta|\leq N_{\mathrm{max}}-|\sigma|$, the ordered mixed derivatives satisfy
\begin{equation}\label{eq:Nash-transported-cutoff-gradient-estimate}
\begin{aligned}
    \Nrm{\left(\prod_{i=1}^k\nabla^{\vartheta_i}\widetilde D_{t,q}^{\sigma_i}\right)
    (\theta_u^{\mathrm N}\chi_\kappa^{\mathrm N}(\Phi_u^{\mathrm N}))}_{N,\mcI_u^{\mathrm N}}
    &\lesssim\lambda_q\MM{|\sigma|,N_t^{\mathrm N},\tau_q^{-1},\ell_{t,q}^{-1}}
    \MM{N+|\vartheta|-1,N_{R,J_*}^{(1)}(|\sigma|),\lambda_q,\ell_q^{-1}}\\
    &\leq\MM{|\sigma|,N_t^{\mathrm N},\tau_q^{-1},\ell_{t,q}^{-1}}\ell_q^{-N-|\vartheta|}.
\end{aligned}
\end{equation}
For $k\geq1$, $\vartheta,\sigma\in\mathbb N_0^k$, and $N\in\mathbb N_0$ with
$|\sigma|\leq N_{\mathrm{ind}}+1$ and $N+|\vartheta|+|\sigma|\leq N_{\mathrm{max}}-1$,
each of the fields $f=\nabla\Phi_u^{\mathrm N}$, $(\nabla\Phi_u^{\mathrm N})^{-1}$,
$J_u^{\mathrm N}$, $(J_u^{\mathrm N})^{-1}$, $\mathbf C_u^{\mathrm N}$, or $(\mathbf C_u^{\mathrm N})^{-1}$ satisfies
\begin{equation}\label{eq:Nash-flow-mixed-integer-estimate}
\begin{aligned}
    \Nrm{\left(\prod_{i=1}^k\nabla^{\vartheta_i}\widetilde D_{t,q}^{\sigma_i}\right)f}_{N,\mcI_u^{\mathrm N}}
    &\lesssim\MM{|\sigma|,N_t^{\mathrm N},\tau_q^{-1},\ell_{t,q}^{-1}}
    \MM{N+|\vartheta|,N_{R,J_*}^{(1)}(|\sigma|),\lambda_q,\ell_q^{-1}}.
\end{aligned}
\end{equation}
The implicit constants may depend on $\rho^*,\mathfrak M,C_*,N_{\mathrm{max}},\alpha$,
but are independent of $q,a,u,\kappa$.
\end{lem}

\begin{proof}
Identity~\eqref{eq:Nash-transported-cutoff-material-identity} follows from \eqref{backward flow}.
It proves \eqref{eq:Nash-transported-cutoff-estimate} when $N=0$.
For $N\geq1$, apply the second inequality in \eqref{est:holder-norm-composition}
with $F=\chi_\kappa^{\mathrm N}$ and $\Psi=\Phi_u^{\mathrm N}$.
By \eqref{eq:Nash-spatial-cutoff-properties} and \eqref{eq:Nash-flow-higher-order-CN-estimate} with material order zero, we have
$$
\begin{aligned}
    \Nrm{\nabla^N\bigl[\chi_\kappa^{\mathrm N}(\Phi_u^{\mathrm N})\bigr]}_{0,\mcI_u^{\mathrm N}}
    \lesssim\|\nabla\chi_\kappa^{\mathrm N}\|_0\|\nabla\Phi_u^{\mathrm N}\|_{N-1,\mcI_u^{\mathrm N}}
    +\|\nabla\chi_\kappa^{\mathrm N}\|_{N-1}\|\nabla\Phi_u^{\mathrm N}\|_{0,\mcI_u^{\mathrm N}}^N
    \lesssim\MM{N,N_{R,J_*}^{(1)}(0)+1,\lambda_q,\ell_q^{-1}}.
\end{aligned}
$$
Only derivatives of $\nabla\Phi_u^{\mathrm N}$ through order $N-1\leq N_{\mathrm{max}}-1$ are used.
Including the lower spatial orders and the $C^0$ composition estimate, then using
\eqref{def:Nash-transition-cutoff} and \eqref{eq:Nash-transported-cutoff-material-identity}, we obtain
\begin{equation}\label{eq:Nash-transported-cutoff-composition-estimate}
\begin{aligned}
\nrm{\widetilde D_{t,q}^M\bigl(\theta_u^{\mathrm N}\chi_\kappa^{\mathrm N}(\Phi_u^{\mathrm N})\bigr)}_{N,\mcI_u^{\mathrm N}}
\leq\|\partial_t^M\theta_u^{\mathrm N}\|_{0,\mcI_u^{\mathrm N}}
\Nrm{\chi_\kappa^{\mathrm N}(\Phi_u^{\mathrm N})}_{N,\mcI_u^{\mathrm N}}
\lesssim\tau_q^{-M}\MM{N,N_{R,J_*}^{(1)}(0)+1,\lambda_q,\ell_q^{-1}}.
\end{aligned}
\end{equation}
Since
$N_{R,J_*}^{(1)}(M)\leq N_{R,J_*}^{(1)}(0)+1$, \eqref{eq:Nash-transported-cutoff-estimate} holds.

To prove \eqref{eq:Nash-transported-cutoff-gradient-estimate}, set $M=|\sigma|$.
Apply \eqref{eq:C0-iterate} in Lemma~\ref{lem:holder-cooper-C0}
with $\upsilon=\bar u_{q,J_*}$, $f=\theta_u^{\mathrm N}\chi_\kappa^{\mathrm N}(\Phi_u^{\mathrm N})$, $M_*=M$, $\mcI=\mcI_u^{\mathrm N}$, and
$$
\begin{aligned}
	N_*&=N_{\mathrm{max}},&& N_\upsilon=N_f=N_{R,J_*}^{(1)}(M)+1,&& M_\upsilon=N_t^{\mathrm N}-1,&& M_f=N_t^{\mathrm N},&&\mC_\upsilon=\tau_q^{-1}\lambda_q^{-1},&&\mC_f=1,\\
	\lambda_\upsilon&=\lambda_f=\lambda_q,&&\tilde\lambda_\upsilon=\tilde\lambda_f=\ell_q^{-1},&&\mu_\upsilon=\mu_f=\mu_*=\tau_q^{-1}
	,&&\tilde\mu_\upsilon=\tilde\mu_f=\tilde\mu_*=\ell_{t,q}^{-1}.
\end{aligned}
$$
The velocity hypothesis follows from \eqref{eq:Nash-terminal-velocity-gradient-Nash-scale}, since $M_0<M\leq N_{\mathrm{ind}}+1$.
The function hypothesis follows from \eqref{eq:Nash-transported-cutoff-composition-estimate} for $N\geq1$.
For $N=0$, use \eqref{def:Nash-transition-cutoff} and \eqref{eq:Nash-transported-cutoff-material-identity}.
In both cases, lower the spatial threshold to the chosen value.
Apply \eqref{eq:C0-iterate} to the term in \eqref{eq:Nash-transported-cutoff-gradient-estimate},
with up to $N$ additional spatial derivatives. The zeroth-order terms are covered by $\|f\|_0\lesssim1$.
Since $N+|\vartheta|\geq1$, the spatial weight is
$
\MM{N+|\vartheta|,N_{R,J_*}^{(1)}(M)+1,\lambda_q,\ell_q^{-1}}
=\lambda_q\MM{N+|\vartheta|-1,N_{R,J_*}^{(1)}(M),\lambda_q,\ell_q^{-1}}.
$
The first inequality in \eqref{eq:Nash-transported-cutoff-gradient-estimate} follows. For the last inequality, choose $a$ sufficiently large to absorb the fixed implicit constant into $\ell_q\lambda_q$.
The same proof applies to $(\chi_\kappa^{\mathrm N})^2$ and with $\theta_u^{\mathrm N}$ replaced by $1$.

For \eqref{eq:Nash-flow-mixed-integer-estimate} at positive total orders,
apply Lemma~\ref{lem:Nash-coefficient-mixed-derivatives} with
$$
s=1,\quad M=|\sigma|,\quad\mC_f=1,\quad N_*=N_{\mathrm{max}}-1,\quad\mcI=\mcI_u^{\mathrm N}.
$$
Estimates \eqref{eq:Nash-flow-higher-order-CN-estimate}--\eqref{eq:Nash-flow-determinant-cofactor-CN-estimate}
verify its hypothesis in the stated range.
At total order zero, use \eqref{eq:Nash-flow-zero-order-estimate}.
\end{proof}

We next record the support and overlap properties of the transported cutoffs.

\begin{lem}
	\label{lem:Nash-cutoff-support-overlap}
	The transported cutoffs satisfy the following properties.
	\begin{enumerate}
		\item For every selected $u$ and $\kappa\in\mathbb Z^3$,
		\begin{equation}
			\label{eq:Nash-transported-cutoff-support}
			\supp\bigl(\theta_u^{\mathrm N}(t)\chi_\kappa^{\mathrm N}(\Phi_u^{\mathrm N}(t,x))\bigr)
			\subseteq
			\left\{
			(t,x)\in\mathbb R\times\mathbb R^3:
			t\in\supp\theta_u^{\mathrm N},\;
			\Phi_u^{\mathrm N}(t,x)\in
			Q\left(2\pi \lambda_q^{-1} \kappa,\frac{9\pi}{8}\lambda_q^{-1}\right)
			\right\}.
		\end{equation}
		The cutoff is also supported in the adapted Lagrangian cylinder from
		Subsection~\ref{subsec:adapted-cylinders}:
		\begin{equation}
			\label{eq:Nash-transported-cutoff-lagrangian-support}
			\supp\bigl(\theta_u^{\mathrm N}(t)\chi_\kappa^{\mathrm N}(\Phi_u^{\mathrm N}(t,x))\bigr)
			\subseteq
			\hat\Gamma_{\bar u_{q,J_*}}\left(
			\frac14\tau_q,2\pi\lambda_q^{-1};
			t_u^{\mathrm N},2\pi\lambda_q^{-1}\kappa
			\right).
		\end{equation}
		
		\item For every selected $u$ and $\kappa\in\mathbb Z^3$,
		\begin{equation}
			\label{eq:Nash-cutoff-cell-intersection-count}
			\#\left\{
			(u',\kappa'):\ u'\text{ selected},\ \kappa'\in\mathbb Z^3,
			\begin{aligned}
				\supp\bigl(\theta_u^{\mathrm N}\chi_\kappa^{\mathrm N}(\Phi_u^{\mathrm N})\bigr)
				\cap
				\supp\bigl(\theta_{u'}^{\mathrm N}\chi_{\kappa'}^{\mathrm N}(\Phi_{u'}^{\mathrm N})\bigr)\neq\varnothing
			\end{aligned}
			\right\}
			\leq27.
		\end{equation}
		
		\item Every point of $\mathbb R\times\mathbb R^3$ belongs to at most
		eight supports
		$
		\supp\bigl(\theta_u^{\mathrm N}\chi_\kappa^{\mathrm N}(\Phi_u^{\mathrm N})\bigr).
		$ The family of cutoff supports is locally finite.
	\end{enumerate}
\end{lem}

\begin{proof}
	Let $(t,x)$ belong to the cutoff support in \eqref{eq:Nash-transported-cutoff-support}. Set $z_\kappa=2\pi\lambda_q^{-1}\kappa$
	and $y=\Phi_u^{\mathrm N}(t,x)$. Then $|t-t_u^{\mathrm N}|<\tau_q/4$
	and $|y-z_\kappa|\leq(9\sqrt3\pi/8)\lambda_q^{-1}<2\pi\lambda_q^{-1}$.
	This proves \eqref{eq:Nash-transported-cutoff-lagrangian-support}.
	
	For \eqref{eq:Nash-cutoff-cell-intersection-count}, suppose that the
	cutoff supports indexed by $(u,\kappa)$ and $(u',\kappa')$ intersect.
	Disjointness of the temporal supports implies $u'=u$. The corresponding cubes in
	\eqref{eq:Nash-transported-cutoff-support} intersect. Since each cube has half-side length $9\pi/(8\lambda_q)$,
	we have, for each $i=1,2,3$,
	$
	\frac{2\pi}{\lambda_q}|\kappa_i-\kappa_i'|
	\leq\frac{9\pi}{4\lambda_q}.
	$
	Thus $|\kappa_i-\kappa_i'|\leq\frac98<2$. Each difference is an integer and belongs to $\{-1,0,1\}$.
	Hence there are at most $3^3=27$ possible indices $(u',\kappa')$.
	
	At a fixed $(t,x)$, at most one time index $u$ is active by
	\eqref{def:Nash-transition-cutoff}. For each cutoff support containing
	$(t,x)$, \eqref{eq:Nash-transported-cutoff-support} gives, for $i=1,2,3$,
	$
	\left|
	\kappa_i-
	\frac{\lambda_q(\Phi_u^{\mathrm N})_i(t,x)}{2\pi}
	\right|
	\leq
	\frac{9}{16}.
	$
	The interval for $\kappa_i$ has length $9/8<2$ and contains at most two
	integers. Thus at most $2^3=8$ cutoff supports contain $(t,x)$.
	There are finitely many selected $u$. For each $u$, restrict any compact
	spacetime set to $t\in\supp\theta_u^{\mathrm N}$. Its image under $\Phi_u^{\mathrm N}$
	is compact and meets only finitely many cutoff cubes, proving local finiteness.
\end{proof}

\subsection{Spatially separated Mikado flows and the momentum perturbation}
\label{subsec:Nash-separated-Mikado}

We follow the Mikado construction in \cite[Lemma~2.3]{DS17}.
The vector potentials defined below are supported in the same tubes as the flows.
Fix a nonzero radial function $\Psi=\Psi(r)\in C_c^\infty(\mathbb R^2)$
with $\supp\Psi\subset B(0,1)$, and normalize it so that
$$
\psi:=-\Delta\Psi,\qquad
\int_{\mathbb R^2}\psi^2\,\rd x=1.
$$
Then $\psi\in C_c^\infty(\mathbb R^2)$ is radial and satisfies
$
\int_{\mathbb R^2}\psi\,\rd x=0,\
\supp\psi\subset\supp\Psi.
$

Take $\Lambda$, $\varepsilon_0$, and $\gamma_\xi$ from
Lemma~\ref{Geometric lemma}, which follows from \cite[Proposition~4.1]{BMNV21}.
We have
$$
S=\sum_{\xi\in\Lambda}\gamma_\xi^2(S)\,\xi\otimes\xi
$$
for every symmetric $S\in B_{\varepsilon_0/2}(\Id)$.
The derivatives of $\gamma_\xi$ are bounded on this smaller matrix ball.

Since $\Lambda$ is finite and its directions are nonzero integer vectors,
we can choose pairwise disjoint translated closed geodesics
$\Gamma_\xi=\{\bar x_\xi+s\xi:s\in\mathbb R\}/2\pi\mathbb Z^3$.
Fix $r_0>0$ small enough that the tubes of radius $r_0$ around these geodesics
are disjoint and have orthogonal normal coordinates.
For $x\in\mathbb T^3$, define
\begin{equation}\label{def:Nash-scalar-Mikado-profiles}
\begin{aligned}
	\Psi_\xi(x):=c_\xi r_0\Psi(\operatorname{dist}(x,\Gamma_\xi)/r_0),\qquad
	\psi_\xi(x):=c_\xi r_0^{-1}\psi(\operatorname{dist}(x,\Gamma_\xi)/r_0),
	\qquad c_\xi:=\left(\frac{(2\pi)^3}{|\Gamma_\xi|}\right)^{\frac12}.
\end{aligned}
\end{equation}
Here the distance is taken on $\mathbb T^3$, and $|\Gamma_\xi|$ is the length
of the closed geodesic.
$$
\supp\Psi_\xi\cup\supp\psi_\xi
\subset\{x\in\mathbb T^3:\operatorname{dist}(x,\Gamma_\xi)<r_0\}.
$$
Define the Mikado flows and their vector potentials by
	\begin{equation}\label{def:Nash-fixed-Mikado-profiles}
	\begin{aligned}
		U_\xi(x):=\xi\psi_\xi(x),\qquad
		\widetilde U_\xi(x):=\nabla\Psi_\xi(x)\times\xi
	\end{aligned}
	\end{equation}
Finally, let $\widetilde\psi_\xi$ be the zero-mean periodic solution of
\begin{equation}\label{def:Nash-Mikado-scalar-potential}
	\Delta\widetilde\psi_\xi=\psi_\xi^2-1,\qquad
	\fint_{\mathbb T^3}\widetilde\psi_\xi\,\rd x=0.
\end{equation}

\begin{lem}
	\label{lem:Nash-spatially-separated-Mikado}
	For every $\xi\in\Lambda$, the profiles and fields defined above are smooth and periodic, with
	$\psi_\xi=-\Delta\Psi_\xi$.
	The Mikado flows and their vector potentials satisfy
	\begin{equation}
		\label{eq:Nash-Mikado-identities}
		\Div U_\xi=0,\qquad
		\Div(U_\xi\otimes U_\xi)=0,\qquad
		\curl\widetilde U_\xi=U_\xi,
	\end{equation}
	and
	\begin{equation}
		\label{eq:Nash-Mikado-averages}
		\fint_{\mathbb T^3}U_\xi\,\rd x
		=\fint_{\mathbb T^3}\widetilde U_\xi\,\rd x=0,\qquad
		\fint_{\mathbb T^3}U_\xi\otimes U_\xi\,\rd x
		=\xi\otimes\xi.
	\end{equation}
	The scalar profiles are constant along the direction $\xi$:
	$$
	\xi\cdot\nabla\Psi_\xi
	=\xi\cdot\nabla\psi_\xi
	=\xi\cdot\nabla\widetilde\psi_\xi=0.
	$$
	For distinct indices, the flows and their vector potentials have disjoint supports:
	\begin{equation}
		\label{eq:Nash-Mikado-support-separation}
		\begin{aligned}
			&(\supp U_\xi\cup\supp\widetilde U_\xi)
			\cap
			(\supp U_{\xi'}\cup\supp\widetilde U_{\xi'})
			=\varnothing
			&\text{whenever }\xi\neq\xi'.
		\end{aligned}
	\end{equation}
\end{lem}

\begin{proof}
	By \eqref{def:Nash-scalar-Mikado-profiles} and $\psi=-\Delta\Psi$,
	we have
	$\Delta\Psi_\xi(x)=-c_\xi r_0^{-1}\psi(\operatorname{dist}(x,\Gamma_\xi)/r_0)=-\psi_\xi(x)$.
	Both profiles depend only on the distance to $\Gamma_\xi$, so
	$\xi\cdot\nabla\Psi_\xi=\xi\cdot\nabla\psi_\xi=0$.
	Periodicity and \eqref{def:Nash-scalar-Mikado-profiles} give
	$$
	\begin{aligned}
		\fint_{\mathbb T^3}\psi_\xi\,\rd x
		=-\fint_{\mathbb T^3}\Delta\Psi_\xi\,\rd x=0,\quad
		\fint_{\mathbb T^3}\psi_\xi^2\,\rd x
		=\frac{c_\xi^2|\Gamma_\xi|}{(2\pi)^3}
		\int_{\mathbb R^2}\psi^2\,\rd x=1.
	\end{aligned}
	$$
	The identities in \eqref{eq:Nash-Mikado-identities} follow from
	$$
	\begin{aligned}
		\Div U_\xi=\xi\cdot\nabla\psi_\xi=0,\quad
		\Div(U_\xi\otimes U_\xi)
		=\xi(\xi\cdot\nabla\psi_\xi^2)=0,\quad
		\curl(\nabla\Psi_\xi\times\xi)
		=\nabla(\xi\cdot\nabla\Psi_\xi)
		-\xi\Delta\Psi_\xi=U_\xi.
	\end{aligned}
	$$
	For the averages in \eqref{eq:Nash-Mikado-averages}, we compute
	$$
	\begin{aligned}
		\fint_{\mathbb T^3}U_\xi\,\rd x
		=\xi\fint_{\mathbb T^3}\psi_\xi\,\rd x=0,\quad
		\fint_{\mathbb T^3}\widetilde U_\xi\,\rd x
		=\left(\fint_{\mathbb T^3}\nabla\Psi_\xi\,\rd x\right)\times\xi=0,\quad
		\fint_{\mathbb T^3}U_\xi\otimes U_\xi\,\rd x
		=(\xi\otimes\xi)\fint_{\mathbb T^3}\psi_\xi^2\,\rd x=\xi\otimes\xi.
	\end{aligned}
	$$

	Since $\fint_{\mathbb T^3}(\psi_\xi^2-1)\,\rd x=0$,
	\eqref{def:Nash-Mikado-scalar-potential} has a unique smooth zero-mean periodic solution. Its derivative along the Mikado direction satisfies
	$$
	\Delta(\xi\cdot\nabla\widetilde\psi_\xi)
	=\xi\cdot\nabla(\psi_\xi^2-1)=0,\qquad
	\fint_{\mathbb T^3}\xi\cdot\nabla\widetilde\psi_\xi\,\rd x=0.
	$$
	Thus $\xi\cdot\nabla\widetilde\psi_\xi=0$.
	Both $U_\xi$ and $\widetilde U_\xi$ are supported in the selected tube,
	which proves \eqref{eq:Nash-Mikado-support-separation}.
\end{proof}

Fix the profiles from Lemma~\ref{lem:Nash-spatially-separated-Mikado}.
Set
\begin{equation}
	\label{def:Nash-stress-level}
	e_q:=A_{q,0}
	(\tau_q\ell_q^{-1})^\alpha.
\end{equation}
On $\mcI_u^{\mathrm N}$, define
\begin{equation}
	\label{def:Nash-raw-pulled-back-matrix}
	\begin{aligned}
		T_{q,u}(t,x)
		:={}&
		(\mathbf C_u^{\mathrm N})^{-T}(t,x)
		\left(\Id-e_q^{-1}
		S_{q,J_*,u}(t,x)\right)(\mathbf C_u^{\mathrm N})^{-1}(t,x).
	\end{aligned}
\end{equation}
Using the along-flow mollification from Subsection~\ref{subsec:initial-Newton-flow-mollification}, define
\begin{equation}\label{def:Nash-flow-mollified-gluing-stress}
	\begin{aligned}
		\bar S_{q,J_*,u}(t,x)
		&:=\bigl(\eta_{\ell_{t,q}}^{\mathrm t}\ast_{X_{\bar u_{q,J_*}}}
		S_{q,J_*,u}\bigr)(t,x),
	\end{aligned}
\end{equation}
and  set
\begin{equation}
\begin{aligned}
\widetilde R_{q,u}
:={}&
(\mathbf C_u^{\mathrm N})^{-T}\left(\Id-e_q^{-1}
\bar S_{q,J_*,u}\right)
(\mathbf C_u^{\mathrm N})^{-1}.
\end{aligned}
\label{def:Nash-pulled-back-matrix}
\end{equation}
By \eqref{eq:material-mollification-hierarchy}, we may choose $a_{\mathrm{Nash}}$ so that $\ell_{t,q}\leq\tau_q/48$.
The mollification enlarges the temporal support by at most $\ell_{t,q}$, so \eqref{def:Nash-transition-cutoff} gives $\theta_u^{\mathrm N}=1$ on $\supp_t\bar S_{q,J_*,u}$.
By \eqref{def:Newton-gluing-stress} and the definition of $S_{q,J_*}$, the stress $S_{q,J_*,u}$ is symmetric, and so are $T_{q,u}$ and $\widetilde R_{q,u}$.

We now verify the smallness required by Lemma~\ref{Geometric lemma}. In this argument, all matrix $C^0$ norms use the Euclidean operator norm.
With $r=\lceil\log_2(3/\beta)\rceil$, \eqref{eq:Nash-terminal-S-material-reserve} at $M=N=0$ and $0\leq\theta_u^{\mathrm N}\leq1$ give
\begin{equation}\label{eq:Nash-normalized-gluing-stress-smallness}
\begin{aligned}
	\|e_q^{-1}\bar S_{q,J_*,u}\|_{0,\mcI_u^{\mathrm N}}
	&\leq\|\eta^{\mathrm t}\|_{L^1(\mathbb R)}\|e_q^{-1}S_{q,J_*,u}\|_0
	\lesssim_{\beta,\mathfrak M,C_*}
	(\tau_q\ell_q^{-1})^{\alpha_r-\alpha}\ell_q^{-2\alpha_r}.
\end{aligned}
\end{equation}
By \eqref{eq:Nash-flow-zero-order-estimate}, for sufficiently large $a$, we have $\|\mathbf C_u^{\mathrm N}-\Id\|_{0,\mcI_u^{\mathrm N}}<1/2$.
Then
$
\|(\mathbf C_u^{\mathrm N})^{-1}\|_{0,\mcI_u^{\mathrm N}}\leq2,\
\|(\mathbf C_u^{\mathrm N})^{-1}-\Id\|_{0,\mcI_u^{\mathrm N}}
\leq2\|\mathbf C_u^{\mathrm N}-\Id\|_{0,\mcI_u^{\mathrm N}}.
$
Definition~\eqref{def:Nash-pulled-back-matrix} gives
$$
\begin{aligned}
\widetilde R_{q,u}-\Id
=\left((\mathbf C_u^{\mathrm N})^{-T}-\Id\right)(\mathbf C_u^{\mathrm N})^{-1}
+\left((\mathbf C_u^{\mathrm N})^{-1}-\Id\right)-e_q^{-1}
(\mathbf C_u^{\mathrm N})^{-T}\bar S_{q,J_*,u}
(\mathbf C_u^{\mathrm N})^{-1}.
\end{aligned}
$$
The inverse-matrix estimates above give
$$
\begin{aligned}
\|\widetilde R_{q,u}-\Id\|_{0,\mcI_u^{\mathrm N}}
=\left\|(\mathbf C_u^{\mathrm N})^{-T}
\left(\Id-e_q^{-1}\bar S_{q,J_*,u}\right)(\mathbf C_u^{\mathrm N})^{-1}-\Id\right\|_{0,\mcI_u^{\mathrm N}}
\leq6\|\mathbf C_u^{\mathrm N}-\Id\|_{0,\mcI_u^{\mathrm N}}
+4e_q^{-1}\|\bar S_{q,J_*,u}\|_{0,\mcI_u^{\mathrm N}}<\frac{\varepsilon_0}{2}.
\end{aligned}
$$
By \eqref{eq:Nash-flow-zero-order-estimate}, \eqref{eq:Nash-normalized-gluing-stress-smallness}, and
$(\tau_q\ell_q^{-1})^{\alpha_r-\alpha}\ell_q^{-2\alpha_r}\lesssim\lambda_q^{-\frac29\alpha\beta}$, choose $a_{\mathrm{Nash}}$
so that each term on the last line is smaller than $\varepsilon_0/4$ for all $q,u$.

For each selected $u$ and $\xi\in\Lambda$, define the amplitudes by
\begin{equation}\label{def:Nash-directional-amplitude}
	\tilde a_{q,u,\xi}(t,x)
	:=\sum_{\kappa\in\mathbb Z^3}\tilde a_{q,u,\kappa,\xi}(t,x)
	:=\sum_{\kappa\in\mathbb Z^3}\theta_u^{\mathrm N}(\chi_\kappa^{\mathrm N})^2(\Phi_u^{\mathrm N})
	\left(e_q\bar\rho_{q,J_*}\right)^{\frac12}
	\gamma_\xi(\widetilde R_{q,u})
	=\theta_u^{\mathrm N}(t)\left(e_q\bar\rho_{q,J_*}(t,x)\right)^{\frac12}
	\gamma_\xi(\widetilde R_{q,u}(t,x)).
\end{equation}
Here the last equality follows from \eqref{eq:Nash-spatial-cutoff-properties}.
The geometric decomposition gives
\begin{equation}\label{eq:Nash-directional-amplitude-identity}
	\sum_{\xi\in\Lambda}\tilde a_{q,u,\xi}^2\,\xi\otimes\xi
	=e_q\bar\rho_{q,J_*}(\theta_u^{\mathrm N})^2\widetilde R_{q,u}.
\end{equation}

Define the momentum perturbation in the Nash step and its components by
\begin{equation}
	\label{def:Nash-perturbation}
	\begin{aligned}
		w_{q+1}
		&:=\sum_u\sum_{\kappa\in\mathbb Z^3}w_{q+1,u,\kappa}
		:=\lambda_{q+1}^{-1}\sum_u\sum_{\kappa\in\mathbb Z^3}\curl\Bigg(
		\sum_{\xi\in\Lambda}(\nabla\Phi_u^{\mathrm N})^T
		\tilde a_{q,u,\kappa,\xi}\widetilde U_\xi(\lambda_{q+1}\Phi_u^{\mathrm N})\Bigg).
	\end{aligned}
\end{equation}
By \eqref{eq:Nash-transported-cutoff-lagrangian-support}, we have
\begin{equation}\label{eq:Nash-localized-perturbation-cylinder-support}
\begin{aligned}
	\bigcup_{\xi\in\Lambda}\supp \tilde a_{q,u,\kappa,\xi}
	\cup\supp w_{q+1,u,\kappa}\subseteq
	\hat\Gamma_{\bar u_{q,J_*}}\left(
	\frac14\tau_q,2\pi\lambda_q^{-1};
	t_u^{\mathrm N},2\pi\lambda_q^{-1}\kappa\right).
\end{aligned}
\end{equation}
By \eqref{def:Nash-perturbation}, the momentum perturbation in the Nash step satisfies
\begin{equation}
	\label{eq:Nash-divergence-and-mean}
	\Div w_{q+1}=0,
	\qquad
	\fint_{\mathbb T^3}w_{q+1}(t,x)\,\rd x=0.
\end{equation}
We define the new state by
\begin{equation}
	\label{def:Nash-final-state}
	\rho_{q+1}:=\rho_{q,J_*},
	\qquad
	m_{q+1}:=m_{q,J_*}+w_{q+1}.
\end{equation}

For a smooth vector field $Z$, set $\Phi=\Phi_u^{\mathrm N}$ and $\mathbf C=\mathbf C_u^{\mathrm N}$. Then
$
\epsilon_{ijk}\partial_j\partial_k\Phi_s=0,\
\epsilon_{ijk}\partial_j\Phi_r\partial_k\Phi_s
=\mathbf C_{ci}\epsilon_{crs}.
$
Hence
$$
\begin{aligned}
	\bigl[\curl\bigl((\nabla\Phi)^TZ(t,\Phi)\bigr)\bigr]_i
	=\epsilon_{ijk}\partial_j\bigl(\partial_k\Phi_s Z_s(t,\Phi)\bigr)
	=\mathbf C_{ci}\epsilon_{crs}(\partial_r Z_s)(t,\Phi)
	=\bigl[\mathbf C^T(\curl Z)(t,\Phi)\bigr]_i.
\end{aligned}
$$
Using $\curl\widetilde U_\xi=U_\xi$, we define the principal part
\begin{equation}
	\label{def:Nash-principal-perturbation}
	\begin{aligned}
		w_{q+1}^{\mathrm p}
		:={}&\sum_u w_{q+1,u}^{\mathrm p}
		=\sum_u
		\sum_{\xi\in\Lambda}
		(\mathbf C_u^{\mathrm N})^T
		\tilde a_{q,u,\xi}(t,x)
		U_\xi(\lambda_{q+1}\Phi_u^{\mathrm N}(t,x)).
	\end{aligned}
\end{equation}
The pullback identity gives the divergence corrector
\begin{equation}
	\label{eq:Nash-corrector-formula}
	\begin{aligned}
		w_{q+1}^{\mathrm c}
		:={}&w_{q+1}-w_{q+1}^{\mathrm p}
		=\sum_u w_{q+1,u}^{\mathrm c}
		=\lambda_{q+1}^{-1}\sum_u
		\sum_{\xi\in\Lambda}
		\nabla \tilde a_{q,u,\xi}(t,x)
		\times\Bigl((\nabla\Phi_u^{\mathrm N})^T
		\widetilde U_\xi(\lambda_{q+1}\Phi_u^{\mathrm N}(t,x))\Bigr).
	\end{aligned}
\end{equation}
By \eqref{def:Nash-transition-cutoff} and \eqref{eq:Nash-Mikado-support-separation},
the potentials corresponding to distinct pairs $(u,\xi)$ in \eqref{def:Nash-perturbation} have disjoint supports.
The same holds for their curls and for the principal and corrector terms,
since multiplication and spatial differentiation do not enlarge support.

We first split the quadratic term as
$$
\frac{w_{q+1}^{\mathrm p}\otimes w_{q+1}^{\mathrm p}}{\rho_{q,J_*}}
=\frac{w_{q+1}^{\mathrm p}\otimes w_{q+1}^{\mathrm p}}{\bar\rho_{q,J_*}}
-\frac{\rho_{q,J_*}-\bar\rho_{q,J_*}}{\rho_{q,J_*}\bar\rho_{q,J_*}}
w_{q+1}^{\mathrm p}\otimes w_{q+1}^{\mathrm p}.
$$
We average the first term in the fast variable and include the second in $R_{q+1}^{\mathrm{mol}}$.
Fix $(t,x)$ and a selected $u$. We define the averaged principal stress by
\begin{equation}
	\label{eq:Nash-principal-quadratic-average}
	\begin{aligned}
		\bar R_{q+1,u}^{\mathrm p}
		&:=(\mathbf C_u^{\mathrm N})^T
		\left[
		\sum_{\xi\in\Lambda}
		\frac{\tilde a_{q,u,\xi}^2}{\bar\rho_{q,J_*}}
		\fint_{\mathbb T^3}
		U_\xi(y)\otimes U_\xi(y)\,\rd y
		\right]\mathbf C_u^{\mathrm N}
		=e_q(\theta_u^{\mathrm N})^2
		(\mathbf C_u^{\mathrm N})^T
		\widetilde R_{q,u}\mathbf C_u^{\mathrm N}.
	\end{aligned}
\end{equation}
For the last equality, we used \eqref{eq:Nash-Mikado-averages} and
\eqref{eq:Nash-directional-amplitude-identity}.
We compare the sum of these averages with the expression involving $T_{q,u}$.
Since $\theta_u^{\mathrm N}=1$ on $\supp S_{q,J_*,u}$,
\eqref{def:Nash-raw-pulled-back-matrix} gives
$$
\begin{aligned}
	\sum_u e_q(\theta_u^{\mathrm N})^2
	(\mathbf C_u^{\mathrm N})^TT_{q,u}\mathbf C_u^{\mathrm N}
	&=e_q\sum_u(\theta_u^{\mathrm N})^2\Id-S_{q,J_*}.
\end{aligned}
$$
The multiple of $\Id$ depends only on time and has zero divergence.
Subtracting the left-hand side from $\sum_u\bar R_{q+1,u}^{\mathrm p}$
gives the stress-mollification term in \eqref{def:Nash-flow-mollification-stress}.

We next subtract $\sum_u\bar R_{q+1,u}^{\mathrm p}$ from the quadratic term with denominator $\bar\rho_{q,J_*}$.
By \eqref{def:Nash-transition-cutoff}, \eqref{eq:Nash-Mikado-support-separation}, and
\eqref{eq:Nash-principal-quadratic-average}, we obtain
\begin{equation}\label{eq:Nash-principal-quadratic-profile-identity}
	\begin{aligned}
		\frac{w_{q+1}^{\mathrm p}\otimes w_{q+1}^{\mathrm p}}{\bar\rho_{q,J_*}}
		-\sum_u\bar R_{q+1,u}^{\mathrm p}=\sum_u\sum_{\xi\in\Lambda}
		\frac{\tilde a_{q,u,\xi}^2}{\bar\rho_{q,J_*}}
		\bigl((\mathbf C_u^{\mathrm N})^T\xi\bigr)
		\otimes\bigl((\mathbf C_u^{\mathrm N})^T\xi\bigr)
		\left[\psi_\xi^2(\lambda_{q+1}\Phi_u^{\mathrm N})-1\right].
	\end{aligned}
\end{equation}
The components of $U_\xi$ and $\widetilde U_\xi$,
their spatial derivatives, and $\psi_\xi^2-1$ are fixed
zero-mean periodic profiles.
For the profile $\psi_\xi^2-1$, we use the zero-mean periodic potential
$\widetilde\psi_\xi$ defined in \eqref{def:Nash-Mikado-scalar-potential}.

\subsection{A localized inverse-divergence operator and the new Reynolds stress}
\label{subsec:Nash-new-Reynolds-stress}

As in \cite[Proposition~7.1]{MQ25Null}, we first write the oscillatory source as the divergence of an explicit matrix plus a remainder.
We then apply a local right inverse of divergence to the remainder.
Let $X(t,\cdot)=X_\upsilon(t,\cdot;t_0)$ be the forward flow of a smooth
velocity field $\upsilon$ from time $t_0$, and let $\Phi=X^{-1}$. Set
$J=\det\nabla\Phi$ and $\mathbf C=\cof\nabla\Phi$.
Let $R^{(\upsilon)}$ be the local right inverse in
Theorem~\ref{thm:inv4divEq}, with point $(t_0,x_0)$,
time radius $\tau$, and spatial radius $\varepsilon$.
We seek a matrix $T$ satisfying
\begin{align*}
\Div T=f,\qquad
f(t,x)=J(t,x)\sum_{r=1}^L h_r(t,x)g_r(\lambda_f\Phi(t,x)).
\end{align*}
Here $g_1,\ldots,g_L$ are fixed smooth zero-mean scalar profiles on $\mathbb T^3$,
and $h_1,\ldots,h_L$ are smooth vector functions.

\begin{lem}
	\label{lem:Nash-transported-oscillatory-inverse-divergence}
	Fix $K^*,K^\sharp,N_t,K_\upsilon\in\mathbb N_0$ with $K^\sharp\geq2$, $0<\tilde\tau\leq\tau\leq1$, and $C_0>0$.
	Let $\upsilon$ be a smooth velocity field and let $g_1,\ldots,g_L$ be fixed smooth zero-mean scalar profiles on $\mathbb T^3$.
	Let $X,\Phi,J,\mathbf C$ be as above.
	Let $0<\varepsilon\leq1$, $1\leq\lambda_\upsilon\leq\tilde\lambda_\upsilon$ with $\tilde\lambda_\upsilon\geq\varepsilon^{-1}$, and
	$\lambda_f\in\mathbb N$ satisfy
	$\lambda_f\geq4\tilde\lambda_\upsilon$ and $\varepsilon\tilde\lambda_\upsilon^2\leq C_0\lambda_f$.
	Fix $\mC_\upsilon\geq0$. For $M,N\in\mathbb N_0$ with $M\leq K^*+1$ and $M+N\leq K^\sharp$, assume that
	\begin{equation}
		\label{eq:Nash-oscillatory-velocity-assumption}
		\nrm{D_{t,\upsilon}^M\nabla\upsilon}_N
		\lesssim\mC_\upsilon\lambda_\upsilon\MM{M,N_t,\tau^{-1},\tilde\tau^{-1}}\MM{N,(K_\upsilon-M)_+,\lambda_\upsilon,\tilde\lambda_\upsilon}.
	\end{equation}
	Assume that $\tau\mC_\upsilon\lambda_\upsilon\leq c_0$.
	Here $c_0>0$ is chosen sufficiently small in terms of the implicit constant in
	\eqref{eq:Nash-oscillatory-velocity-assumption} so that
	$\nrm{\nabla\Phi-\Id}_0+\nrm{\nabla X-\Id}_0\leq1/2$ for $|t-t_0|\leq\tau$.
	For $k\geq1$, $\vartheta,\sigma\in\mathbb N_0^k$, and $N\in\mathbb N_0$ with
	$|\sigma|\leq K^*$ and $N+|\vartheta|+|\sigma|\leq K^\sharp-1$, assume that each of the fields
	$F=\nabla\Phi$, $(\nabla\Phi)^{-1}$, $J$, $J^{-1}$, $\mathbf C$, or $J^{-1}\mathbf C$ satisfies
	\begin{equation}\label{eq:Nash-oscillatory-flow-mixed-assumption}
	\Nrm{\left(\prod_{i=1}^k\nabla^{\vartheta_i}D_{t,\upsilon}^{\sigma_i}\right)F}_N
	\lesssim\MM{|\sigma|,N_t,\tau^{-1},\tilde\tau^{-1}}\tilde\lambda_\upsilon^{N+|\vartheta|}.
	\end{equation}

	Let $h_1,\ldots,h_L:\mathbb R\times\mathbb R^3\to\mathbb R^3$ be smooth vector functions satisfying
	\begin{equation}\label{eq:Nash-oscillatory-source-support}
		\bigcup_{r=1}^L\supp h_r
		\subset\hat\Gamma_\upsilon(\tau,\varepsilon;t_0,x_0).
	\end{equation}
	Fix $\mC_f\geq0$. For $M,N\in\mathbb N_0$ with $M\leq K^*$ and $M+N\leq K^\sharp$, assume that
	\begin{equation}
		\label{eq:Nash-oscillatory-amplitude-assumption}
		\sum_{r=1}^L\nrm{D_{t,\upsilon}^M h_r}_{N}
		\lesssim \mC_f\MM{M,N_t,\tau^{-1},\tilde\tau^{-1}}\tilde\lambda_\upsilon^{N}.
	\end{equation}
	Define
	$
	f(t,x)=J(t,x)\sum_{r=1}^L h_r(t,x)g_r(\lambda_f\Phi(t,x)).
	$
	Suppose that $\int_{\mathbb R^3}f(t,x)\,\rd x=0$ for $|t-t_0|\leq\tau$.
	Then, there is a matrix $\mathcal R_{\lambda_f}^{(\upsilon)}[f]$,
	linear in the amplitudes $h_r$, such that
	$$
	\Div\mathcal R_{\lambda_f}^{(\upsilon)}[f]=f,\qquad
	\supp_t\mathcal R_{\lambda_f}^{(\upsilon)}[f]\subset\bigcup_{r=1}^L\supp_t h_r,\qquad \supp\mathcal R_{\lambda_f}^{(\upsilon)}[f]\subset\hat\Gamma_\upsilon(\tau,\varepsilon;t_0,x_0)
	$$
	 For $M,N\in\mathbb N_0$ with $M\leq K^*$ and $M+N\leq K^\sharp-2$, we have
	\begin{equation}
		\label{eq:Nash-transported-oscillatory-inverse-divergence}
		\nrm{D_{t,\upsilon}^{M}\mathcal R_{\lambda_f}^{(\upsilon)}[f]}_{N}
		\lesssim
		\mC_f\MM{M,N_t,\tau^{-1},\tilde\tau^{-1}}\lambda_f^{N-1}.
	\end{equation}
	The implicit constant may depend on $K^*,K^\sharp,C_0,c_0$ and the constants in
	\eqref{eq:Nash-oscillatory-velocity-assumption}--%
	\eqref{eq:Nash-oscillatory-flow-mixed-assumption} and
	\eqref{eq:Nash-oscillatory-amplitude-assumption}.
\end{lem}

\begin{proof}
	For each $r$, let $G_r$ be the zero-mean periodic solution of $\Delta G_r=g_r$. Define
	$$
	\begin{aligned}
		Q_{il}&:=\mathbf C_{sl}\sum_{r=1}^L\left[
		\lambda_f^{-1}(h_r)_i(\partial_sG_r)(\lambda_f\Phi)
		-\lambda_f^{-2}J^{-1}\mathbf C_{sj}\partial_{j}(h_r)_i G_r(\lambda_f\Phi)\right],\\
		U&:=\sum_{r=1}^L\mathbf C_{sj}\partial_{j}\left(J^{-1}\mathbf C_{sk}\partial_{k}h_r\right)G_r(\lambda_f\Phi).
	\end{aligned}
	$$
	Extend $Q$ and $U$ by zero outside $[t_0-\tau,t_0+\tau]$.
	The Piola identity and $D_{t,\upsilon}\Phi=0$ give
	\begin{equation}\label{eq:Nash-Piola-material-identities}
	\partial_l\mathbf C_{sl}=0,\qquad
	J^{-1}\mathbf C_{sl}\partial_l\Phi_j=\delta_{sj},\qquad
	[D_{t,\upsilon},J^{-1}\mathbf C_{sl}\partial_l]=0.
	\end{equation}
	For the last identity, differentiate $J^{-1}\mathbf C=(\nabla\Phi)^{-T}$ and use $D_{t,\upsilon}\Phi=0$ to obtain
	$$
	D_{t,\upsilon}(J^{-1}\mathbf C_{sl})
	=J^{-1}\mathbf C_{sj}\partial_j\upsilon_l.
	$$
	By \eqref{eq:first-order-material-spatial-commutator}, we compute
	$$
	\begin{aligned}
	[D_{t,\upsilon},J^{-1}\mathbf C_{sl}\partial_l]
	=D_{t,\upsilon}(J^{-1}\mathbf C_{sl})\partial_l
	-J^{-1}\mathbf C_{sl}(\partial_l\upsilon_j)\partial_j
	=J^{-1}\mathbf C_{sj}(\partial_j\upsilon_l)\partial_l
	-J^{-1}\mathbf C_{sl}(\partial_l\upsilon_j)\partial_j=0.
	\end{aligned}
	$$
	The first two identities in \eqref{eq:Nash-Piola-material-identities} and $\Delta G_r=g_r$ give
	\begin{equation}\label{eq:Nash-twofold-profile-divergence}
		\begin{aligned}
			&\partial_l\left(\mathbf C_{sl}\left[
			\lambda_f^{-1}(h_r)_i(\partial_sG_r)(\lambda_f\Phi)
			-\lambda_f^{-2}J^{-1}\mathbf C_{sj}\partial_{j}(h_r)_i G_r(\lambda_f\Phi)\right]\right)\\
			&\qquad=J(h_r)_i g_r(\lambda_f\Phi)
			-\lambda_f^{-2}\mathbf C_{sj}\partial_{j}\left(J^{-1}\mathbf C_{sk}\partial_{k}(h_r)_i\right)G_r(\lambda_f\Phi).
		\end{aligned}
	\end{equation}
	Summing over $r$ gives $\Div Q=f-\lambda_f^{-2}U$.
	Since $\supp Q\cup\supp U\subset\bigcup_{r=1}^L\supp h_r$, integrating this identity gives
	$$
	\int_{\mathbb R^3}U(t,x)\,\rd x
	=\lambda_f^2\int_{\mathbb R^3}f(t,x)\,\rd x=0.
	$$

	We first estimate $U$ and $Q$ for $M,N\in\mathbb N_0$ with $M\leq K^*$ and $M+N\leq K^\sharp-2$.
	By \eqref{eq:Nash-Piola-material-identities}, $D_{t,\upsilon}$ commutes with $J^{-1}\mathbf C_{sl}\partial_l$.
	Since $D_{t,\upsilon}G_r(\lambda_f\Phi)=0$, we have
	$$
	D_{t,\upsilon}^M(J^{-1}U)
	=\sum_{r=1}^LJ^{-1}\mathbf C_{sj}\partial_j\left(
	J^{-1}\mathbf C_{sk}\partial_kD_{t,\upsilon}^Mh_r\right)G_r(\lambda_f\Phi).
	$$
	For $0\leq k\leq N$, apply the second inequality in \eqref{est:holder-norm-composition} when $k\geq1$ and the $C^0$ composition estimate when $k=0$.
	Using \eqref{eq:Nash-oscillatory-flow-mixed-assumption} with $F=\nabla\Phi$ and $|\sigma|=|\vartheta|=0$, and $\lambda_f\geq4\tilde\lambda_\upsilon$, we obtain
	$$
	\|G_r(\lambda_f\Phi)\|_k+\|(\partial_sG_r)(\lambda_f\Phi)\|_k\lesssim\lambda_f^k.
	$$
	Expand the inner spatial derivative and use $U=J(J^{-1}U)$. Then
	\eqref{eq:holder-product-tame} with $\alpha=0$, \eqref{eq:Nash-oscillatory-flow-mixed-assumption},
	and \eqref{eq:Nash-oscillatory-amplitude-assumption} give
	\begin{equation}\label{eq:Nash-oscillatory-remainder-source-estimate}
	\begin{aligned}
	\|D_{t,\upsilon}^MU\|_N
	&\lesssim\sum_{\substack{M_1+M_2=M\\N_1+\cdots+N_5=N}}\sum_{r=1}^L
	\|D_{t,\upsilon}^{M_1}J\|_{N_1}\|J^{-1}\mathbf C\|_{N_2}
	\|G_r(\lambda_f\Phi)\|_{N_5}\Bigl(
	\|J^{-1}\mathbf C\|_{N_3+1}\|D_{t,\upsilon}^{M_2}h_r\|_{N_4+1}
	+\|J^{-1}\mathbf C\|_{N_3}\|D_{t,\upsilon}^{M_2}h_r\|_{N_4+2}\Bigr)\\
	&\lesssim\mC_f\MM{M,N_t,\tau^{-1},\tilde\tau^{-1}}
	\sum_{k=0}^N\tilde\lambda_\upsilon^{k+2}\lambda_f^{N-k}\\
	&\lesssim\mC_f\tilde\lambda_\upsilon^2\MM{M,N_t,\tau^{-1},\tilde\tau^{-1}}\lambda_f^N.
	\end{aligned}
	\end{equation}
	The derivatives of $h_r$ have total order at most $M+N+2\leq K^\sharp$. The spatial order on $J^{-1}\mathbf C$ is at most $N+1\leq K^\sharp-1$.

	For $Q$, \eqref{eq:Nash-Piola-material-identities} and $D_{t,\upsilon}\Phi=0$ give
	$$
	\begin{aligned}
	D_{t,\upsilon}^MQ_{il}
	=\sum_{M_1+M_2=M}\binom{M}{M_1}\sum_{r=1}^L
	(D_{t,\upsilon}^{M_1}\mathbf C_{sl})\Bigl[
	\lambda_f^{-1}(D_{t,\upsilon}^{M_2}h_r)_i(\partial_sG_r)(\lambda_f\Phi)
	-\lambda_f^{-2}J^{-1}\mathbf C_{sj}\partial_j(D_{t,\upsilon}^{M_2}h_r)_iG_r(\lambda_f\Phi)\Bigr].
	\end{aligned}
	$$
	Applying \eqref{eq:holder-product-tame} with $\alpha=0$ to this expansion gives
	\begin{equation}\label{eq:Nash-oscillatory-principal-inverse-estimate}
	\begin{aligned}
	\|D_{t,\upsilon}^MQ\|_N
	&\lesssim\lambda_f^{-1}
	\sum_{\substack{M_1+M_2=M\\N_1+N_2+N_3=N}}\sum_{r=1}^L
	\|D_{t,\upsilon}^{M_1}\mathbf C\|_{N_1}\|D_{t,\upsilon}^{M_2}h_r\|_{N_2}
	\|(\nabla G_r)(\lambda_f\Phi)\|_{N_3}\\
	&\quad+\lambda_f^{-2}
	\sum_{\substack{M_1+M_2=M\\N_1+\cdots+N_4=N}}\sum_{r=1}^L
	\|D_{t,\upsilon}^{M_1}\mathbf C\|_{N_1}\|J^{-1}\mathbf C\|_{N_2}\|D_{t,\upsilon}^{M_2}h_r\|_{N_3+1}\|G_r(\lambda_f\Phi)\|_{N_4}\\
	&\lesssim\mC_f\MM{M,N_t,\tau^{-1},\tilde\tau^{-1}}
	\sum_{k=0}^N\left(\lambda_f^{-1}\tilde\lambda_\upsilon^k
	+\lambda_f^{-2}\tilde\lambda_\upsilon^{k+1}\right)\lambda_f^{N-k}\\
	&\lesssim\mC_f\MM{M,N_t,\tau^{-1},\tilde\tau^{-1}}\lambda_f^{N-1}\left(1+\frac{\tilde\lambda_\upsilon}{\lambda_f}\right).
	\end{aligned}
	\end{equation}

	The source $U$ is smooth and has zero spatial mean.
	The support inclusion $\supp U\subset\hat\Gamma_\upsilon(\tau,\varepsilon;t_0,x_0)$ follows from \eqref{eq:Nash-oscillatory-source-support}.
	Use the operator $R^{(\upsilon)}[U]$ defined by \eqref{eq:Rjl0}--\eqref{eq:Rjl2}, with $\bar\tau=\tau$, $\bar\varepsilon=\varepsilon$, and the fixed averaging kernel.
	The divergence computation in the proof of Proposition~\ref{prop:Rjl} gives $\Div R^{(\upsilon)}[U]=U$.
	Define
	\begin{equation}\label{def:Nash-transported-local-inverse}
		\mathcal R_{\lambda_f}^{(\upsilon)}[f]:=Q+\lambda_f^{-2}R^{(\upsilon)}[U].
	\end{equation}
	Then $\Div\mathcal R_{\lambda_f}^{(\upsilon)}[f]=f$.
	Follow the argument for \eqref{eq:inv4divEq:CN-from-integer-source-zero}--\eqref{eq:inv4divEq:CN-from-integer-source}
	in Theorem~\ref{thm:inv4divEq} with
	$$
	\begin{aligned}
		K^*&=\max\{1,M\},&& K^\sharp=\max\{3,M+2,M+N+1\},&& K_U=0,&& M_*=M,\\
		\lambda_U&=\tilde\lambda_U=\lambda_f,&&\mu_U=\tau^{-1},&&\tilde\mu_U=\tilde\tau^{-1},&&\mC_U=C\mC_f\tilde\lambda_\upsilon^2.
	\end{aligned}
	$$
	Here $C\geq1$ absorbs the implicit constant in \eqref{eq:Nash-oscillatory-remainder-source-estimate}.
	The flow estimates needed in this argument follow from \eqref{eq:Nash-oscillatory-flow-mixed-assumption}.
	The argument uses \eqref{eq:inv4divEq:integer-source-assumption}
	only for $0\leq j\leq M$ and $j+k\leq M+N$,
	so \eqref{eq:Nash-oscillatory-remainder-source-estimate} gives
	\begin{equation}\label{eq:Nash-oscillatory-remainder-inverse-estimate}
	\lambda_f^{-2}\|D_{t,\upsilon}^MR^{(\upsilon)}[U]\|_N
	\lesssim\mC_f\MM{M,N_t,\tau^{-1},\tilde\tau^{-1}}
	\varepsilon\tilde\lambda_\upsilon^2\lambda_f^{N-2}.
	\end{equation}
	Combining \eqref{eq:Nash-oscillatory-principal-inverse-estimate} and \eqref{eq:Nash-oscillatory-remainder-inverse-estimate} gives
	\begin{equation}\label{eq:Nash-local-profile-integer-estimate}
		\begin{aligned}
			\|D_{t,\upsilon}^M\mathcal R_{\lambda_f}^{(\upsilon)}[f]\|_N
			&\lesssim\mC_f\MM{M,N_t,\tau^{-1},\tilde\tau^{-1}}\lambda_f^{N-1}\left(1+\frac{\tilde\lambda_\upsilon}{\lambda_f}
			+\frac{\varepsilon\tilde\lambda_\upsilon^2}{\lambda_f}\right)
			\lesssim\mC_f\MM{M,N_t,\tau^{-1},\tilde\tau^{-1}}\lambda_f^{N-1}.
		\end{aligned}
	\end{equation}
	This proves \eqref{eq:Nash-transported-oscillatory-inverse-divergence}.
	The support claims follow from
	the support inclusion for $Q,U$ above,
	\eqref{eq:supp4Rjl}, and \eqref{eq:time-supp-R}.
\end{proof}

Next, we use Lemma~\ref{lem:Nash-transported-oscillatory-inverse-divergence} to construct the transport and oscillation errors
$R_{q+1}^{\mathrm{trans}}$ and $R_{q+1}^{\mathrm{osc}}$.
We then define the new Reynolds stress $R_{q+1}$.
For each $(u,\kappa)$, take
$$
\begin{aligned}
	&\upsilon=\bar u_{q,J_*},\quad \Phi=\Phi_u^{\mathrm N},\quad X=X_u^{\mathrm N},\quad t_0=t_u^{\mathrm N},\quad x_0=2\pi\lambda_q^{-1}\kappa,\quad	\tau=\tfrac34\tau_q,\quad \tilde\tau=\ell_{t,q},\quad \varepsilon=2\pi\lambda_q^{-1},\quad C_0=2\pi,\\
	&K_\upsilon=N_R-3J_*-1,\quad \lambda_\upsilon=\lambda_q,\quad \tilde\lambda_\upsilon=\ell_q^{-1},\quad \lambda_f=\lambda_{q+1},\quad \mC_\upsilon=\delta_q^{\frac12}.
\end{aligned}
$$
Use the construction in \eqref{def:Nash-transported-local-inverse}
with these parameters.
For the profile representations of the sources below, write
\begin{equation}\label{def:Nash-cell-oscillatory-inverse}
	\mathcal R_{q,u,\kappa}^{\mathrm N}[F]
	:=\mathcal R_{\lambda_{q+1}}^{(\bar u_{q,J_*})}[F].
\end{equation}
For sufficiently large $a_{\mathrm{Nash}}$, the scale and short-time conditions
in Lemma~\ref{lem:Nash-transported-oscillatory-inverse-divergence} follow from
\eqref{def:lambda-delta}, \eqref{def:parameter-mollification}, and
\eqref{def:tau-q}--\eqref{eq:material-mollification-hierarchy}. 

Using \eqref{def:Nash-perturbation}, define the transport source by
\begin{equation}\label{def:Nash-transport-source}
\begin{aligned}
	F_{q+1}^{\mathrm{trans}}
	:=\sum_{u,\kappa}F_{q+1,u,\kappa}^{\mathrm{trans}}
	:=\sum_{u,\kappa}\Bigl[\partial_t w_{q+1,u,\kappa}+\Div\bigl(\bar u_{q,J_*}\otimes w_{q+1,u,\kappa}+w_{q+1,u,\kappa}\otimes \bar u_{q,J_*}\bigr)\Bigr].
\end{aligned}
\end{equation}
Since $w_{q+1,u,\kappa}$ is the curl of a compactly supported potential,
\eqref{def:Nash-transport-source} gives
$\int_{\mathbb R^3}F_{q+1,u,\kappa}^{\mathrm{trans}}(t,x)\,\rd x=0.$
Using $\Div w_{q+1,u,\kappa}=0$ and $\widetilde D_{t,q}\Phi_u^{\mathrm N}=0$, we expand \eqref{def:Nash-perturbation} in \eqref{def:Nash-transport-source}.
This expresses $(J_u^{\mathrm N})^{-1}F_{q+1,u,\kappa}^{\mathrm{trans}}$ as a sum of vector amplitudes
times scalar profiles at $\lambda_{q+1}\Phi_u^{\mathrm N}$.
The profiles are $\psi_\xi$ and the components of $\widetilde U_\xi$.
Since $\psi_\xi=-\Delta\Psi_\xi$, periodicity and \eqref{eq:Nash-Mikado-averages} give zero mean for all these profiles.
Applying the localized inverse to these sources, we define the transport error by
\begin{equation}
	\label{def:Nash-transport-stress}
	R_{q+1}^{\mathrm{trans}}
	:=\sum_{u,\kappa}R_{q+1,u,\kappa}^{\mathrm{trans}}:=\sum_{u,\kappa}\mathcal R_{q,u,\kappa}^{\mathrm N}[F_{q+1,u,\kappa}^{\mathrm{trans}}].
\end{equation}

We next define the oscillation source by
\begin{equation}
	\label{def:Nash-oscillation-source}
\begin{aligned}
	F_{q+1}^{\mathrm{osc}}
	:=\sum_{u,\kappa}F_{q+1,u,\kappa}^{\mathrm{osc}}
	:=\sum_{u,\kappa}\Div\Bigg((\chi_\kappa^{\mathrm N})^2(\Phi_u^{\mathrm N})\left(\frac{w_{q+1,u}^{\mathrm p}\otimes w_{q+1,u}^{\mathrm p}}{\bar\rho_{q,J_*}}-\bar R_{q+1,u}^{\mathrm p}\right)\Bigg).
\end{aligned}
\end{equation}
The tensors in \eqref{def:Nash-oscillation-source} are compactly supported by \eqref{eq:Nash-transported-cutoff-lagrangian-support}, so
$$\int_{\mathbb R^3}F_{q+1,u,\kappa}^{\mathrm{osc}}(t,x)\,\rd x=0.$$
Using \eqref{eq:Nash-Piola-material-identities} and $\xi\cdot\nabla\psi_\xi=0$, we obtain
\begin{equation}\label{eq:Nash-oscillation-profile-source}
	\begin{aligned}
		\left(F_{q+1,u,\kappa}^{\mathrm{osc}}\right)_i
		&=\sum_{\xi\in\Lambda}
		(\mathbf C_u^{\mathrm N})_{sl}\partial_l\left(
		\frac{(\chi_\kappa^{\mathrm N})^2(\Phi_u^{\mathrm N})\tilde a_{q,u,\xi}^2}{\bar\rho_{q,J_*}}
		(\mathbf C_u^{\mathrm N})_{ri}
		\xi_r\xi_s\right)
		\left[\psi_\xi^2(\lambda_{q+1}\Phi_u^{\mathrm N})-1\right]\\
		&=\sum_{\xi\in\Lambda}
		(\mathbf C_u^{\mathrm N})_{sl}\partial_l\left(
		\frac{(\chi_\kappa^{\mathrm N})^2(\Phi_u^{\mathrm N})\tilde a_{q,u,\xi}^2}{\bar\rho_{q,J_*}}
		(\mathbf C_u^{\mathrm N})_{ri}
		\xi_r\xi_s\right)
		g_\xi(\lambda_{q+1}\Phi_u^{\mathrm N}).
	\end{aligned}
\end{equation}
For Lemma~\ref{lem:Nash-transported-oscillatory-inverse-divergence}, take
$
g_\xi:=\psi_\xi^2-1,\
G_\xi:=\widetilde\psi_\xi,
$
so that $\Delta G_\xi=g_\xi$ by \eqref{def:Nash-Mikado-scalar-potential}.
With these profiles and potentials, we define the oscillation error by
\begin{equation}
	\label{def:Nash-oscillation-stress}
	R_{q+1}^{\mathrm{osc}}
	:=\sum_{u,\kappa}R_{q+1,u,\kappa}^{\mathrm{osc}}:=\sum_{u,\kappa}\mathcal R_{q,u,\kappa}^{\mathrm N}[F_{q+1,u,\kappa}^{\mathrm{osc}}].
\end{equation}

The support inclusions \eqref{eq:Nash-transported-cutoff-support}--%
\eqref{eq:Nash-transported-cutoff-lagrangian-support} verify
\eqref{eq:Nash-oscillatory-source-support}.
By Lemma~\ref{lem:Nash-transported-oscillatory-inverse-divergence}, both localized errors are supported in
$\hat\Gamma_{\bar u_{q,J_*}}(\tfrac34\tau_q,2\pi\lambda_q^{-1};
t_u^{\mathrm N},2\pi\lambda_q^{-1}\kappa)
\cap(\supp\theta_u^{\mathrm N}\times\mathbb R^3)$.
Fix $(t,x)\in\mathbb R\times\mathbb R^3$.
By \eqref{def:Nash-transition-cutoff}, at most one $u$ satisfies $t\in\supp\theta_u^{\mathrm N}$.
The cylinder condition gives
$|\lambda_q\Phi_u^{\mathrm N}(t,x)/(2\pi)-\kappa|<1$,
so each coordinate of $\kappa$ has at most two possible integer values. Thus
\begin{equation}\label{eq:Nash-local-inverse-cylinder-overlap}
	\#\left\{(u,\kappa):(t,x)\in
	\supp R_{q+1,u,\kappa}^{\mathrm{trans}}
	\cup\supp R_{q+1,u,\kappa}^{\mathrm{osc}}\right\}\leq8.
\end{equation}
In particular, the sums in \eqref{def:Nash-transport-stress} and \eqref{def:Nash-oscillation-stress}
each contain at most eight nonzero terms at $(t,x)$.
The compact-set argument in the proof of Lemma~\ref{lem support of R under cutoffs},
with spatial spacing $2\pi\lambda_q^{-1}$, also proves local finiteness of this family of supports.
The zero extensions are smooth since $\theta_u^{\mathrm N}\in C_c^\infty(I_u)$.
Since $\lambda_q,\lambda_{q+1}\in\mathbb N$, the periodicity argument in the proof of Lemma~\ref{lem support of R under cutoffs} applies.
Thus $R_{q+1}^{\mathrm{trans}}$ and $R_{q+1}^{\mathrm{osc}}$ are smooth periodic matrix fields.

We define the mollification error by
\begin{equation}
	\label{def:Nash-flow-mollification-stress}
	\begin{aligned}
		R_{q+1}^{\mathrm{mol}}
		:=&S_{q,J_*}
		+\sum_u\bar R_{q+1,u}^{\mathrm p}
		-\sum_u e_q(\theta_u^{\mathrm N})^2\Id\\
		&+(u_{q,J_*}-\bar u_{q,J_*})\otimes w_{q+1}
		+w_{q+1}\otimes(u_{q,J_*}-\bar u_{q,J_*})-\frac{\rho_{q,J_*}-\bar\rho_{q,J_*}}{\rho_{q,J_*}\bar\rho_{q,J_*}}
		w_{q+1}^{\mathrm p}\otimes w_{q+1}^{\mathrm p}.
	\end{aligned}
\end{equation}
The quadratic terms involving $w_{q+1}^{\mathrm c}$ are collected in the divergence corrector error
\begin{equation}
	\label{def:Nash-corrector-stress}
	R_{q+1}^{\mathrm{corr}}
	:=\frac{1}{\rho_{q,J_*}}
	\left(w_{q+1}^{\mathrm p}\otimes w_{q+1}^{\mathrm c}
	+w_{q+1}^{\mathrm c}\otimes w_{q+1}^{\mathrm p}
	+w_{q+1}^{\mathrm c}\otimes w_{q+1}^{\mathrm c}\right).
\end{equation}
Adding these four errors to $R_{q,J_*}$ and $P_{q+1,J_*}$, we define the new Reynolds stress by
\begin{equation}
	\label{def:Nash-final-Reynolds-stress}
	\begin{aligned}
		R_{q+1}:={}&R_{q,J_*}+P_{q+1,J_*}
		+R_{q+1}^{\mathrm{trans}}
		+R_{q+1}^{\mathrm{osc}}
		+R_{q+1}^{\mathrm{mol}}
		+R_{q+1}^{\mathrm{corr}}.
	\end{aligned}
\end{equation}

Using Lemma~\ref{lem:Nash-transported-oscillatory-inverse-divergence} in
\eqref{def:Nash-transport-stress} and \eqref{def:Nash-oscillation-stress},
we combine the definitions above with \eqref{eq:Nash-principal-quadratic-profile-identity} to obtain
$$
\begin{aligned}
	\Div\left(R_{q+1}^{\mathrm{trans}}+R_{q+1}^{\mathrm{osc}}
	+R_{q+1}^{\mathrm{mol}}+R_{q+1}^{\mathrm{corr}}\right)
	=\partial_t w_{q+1}
	+\Div\Bigg(u_{q,J_*}\otimes w_{q+1}
	+w_{q+1}\otimes u_{q,J_*}+\frac{w_{q+1}\otimes w_{q+1}}{\rho_{q,J_*}}+S_{q,J_*}\Bigg).
\end{aligned}
$$
	By this identity, \eqref{eq:Nash-divergence-and-mean}, and \eqref{eq:Newton-substep-Euler-Reynolds-system} with $j=J_*$,
	the state in \eqref{def:Nash-final-state} and the stress \eqref{def:Nash-final-Reynolds-stress} satisfy
	\begin{equation}
		\label{eq:Nash-output-system}
		\left\{
		\begin{aligned}
			&\partial_t\rho_{q+1}+\Div m_{q+1}=0,\\
			&\partial_t m_{q+1}
			+\Div\left(\frac{m_{q+1}\otimes m_{q+1}}{\rho_{q+1}}\right)
			+\nabla p(\rho_{q+1})
			=\Div R_{q+1}.
		\end{aligned}
		\right.
	\end{equation}

\subsection{Estimates for the momentum perturbation and Reynolds stress}
\label{subsec:Nash-estimates}

\begin{prop}
	\label{prop:Nash-perturbation-estimates}
	Assume the hypotheses of this section.
	For each selected $u$, $\kappa\in\mathbb Z^3$, $\xi\in\Lambda$, and $M,N\in\mathbb N_0$ with
	$M\leq N_{\mathrm{ind}}$ and $M+N\leq K_{J_*}^{(0)}+1$, the amplitudes satisfy
	\begin{equation}\label{eq:Nash-amplitude-estimate}
	\begin{aligned}
		\|\widetilde D_{t,q}^M\tilde a_{q,u,\xi}\|_{N,\mcI_u^{\mathrm N}}
		+\|\widetilde D_{t,q}^M\tilde a_{q,u,\kappa,\xi}\|_{N,\mcI_u^{\mathrm N}}
		&\lesssim A_{q,0}^{\frac12}(\tau_q\ell_q^{-1})^{\frac\alpha2}
		\MM{M,N_t^{\mathrm N},\tau_q^{-1},\ell_{t,q}^{-1}}
		\MM{N,N_{R,J_*}^{(1)}(M),\lambda_q,\ell_q^{-1}}.
	\end{aligned}
	\end{equation}
	For $M,N\in\mathbb N_0$ with $M\leq N_{\mathrm{ind}}$ and $M+N\leq K_{J_*}^{(0)}$,
	the principal perturbation, the corrector, and the total perturbation satisfy
	\begin{align}
		\nrm{\widetilde D_{t,q}^Mw_{q+1}^{\mathrm p}}_{N}
		&\lesssim
		A_{q,0}^{\frac12}
		(\tau_q\ell_q^{-1})^{\frac\alpha2}\MM{M,N_t^{\mathrm N},\tau_q^{-1},\ell_{t,q}^{-1}}
		\lambda_{q+1}^{N},
		\label{eq:Nash-principal-estimate}\\
		\nrm{\widetilde D_{t,q}^Mw_{q+1}^{\mathrm c}}_{N}
		&\lesssim
		A_{q,0}^{\frac12}
		(\tau_q\ell_q^{-1})^{\frac\alpha2}\MM{M,N_t^{\mathrm N},\tau_q^{-1},\ell_{t,q}^{-1}}
		\MM{1,N_{R,J_*}^{(1)}(M),\lambda_q,\ell_q^{-1}}
		\lambda_{q+1}^{N-1},
		\label{eq:Nash-corrector-estimate}\\
		\nrm{\widetilde D_{t,q}^Mw_{q+1}}_{N}
		&\lesssim
		A_{q,0}^{\frac12}
		(\tau_q\ell_q^{-1})^{\frac\alpha2}\MM{M,N_t^{\mathrm N},\tau_q^{-1},\ell_{t,q}^{-1}}
		\lambda_{q+1}^{N}.
		\label{eq:Nash-total-perturbation-estimate}
	\end{align}
	The perturbation also satisfies
	\begin{equation}\label{eq:Nash-perturbation-temporal-support}
	\supp_t w_{q+1}\subset\bigcup_u\supp\theta_u^{\mathrm N}.
	\end{equation}
	The implicit constants may depend on $\rho^*,\mathfrak M,C_*,N_{\mathrm{max}},\alpha$.
	They are independent of $q,a,u,\kappa,\xi$.
\end{prop}

\begin{proof}[Proof of Proposition~\ref{prop:Nash-perturbation-estimates}]
	We first prove \eqref{eq:Nash-amplitude-estimate}, then estimate the principal perturbation and the corrector.
	For the coefficient estimates below, we use the spatial threshold $N_{R,J_*}^{(1)}(M)$ and the material threshold
	$N_t^{\mathrm N}$ from \eqref{def:Nash-terminal-material-threshold}.
	As in the proof of Lemma~\ref{lem:initial-Newton-flow-mollification}, estimate the spatial compositions in
	\eqref{def:Nash-flow-mollified-gluing-stress} by the second inequality in \eqref{est:holder-norm-composition} when $N\geq1$.
	For $N=0$, use the $C^0$ composition estimate.
	By \eqref{eq:Nash-localized-S-material-estimate} and \eqref{eq:Nash-flow-higher-order-CN-estimate}, we obtain, for
	$M\leq N_{\mathrm{ind}}$ and $M+N\leq K_{J_*}^{(0)}+1$,
	\begin{equation}\label{eq:Nash-time-mollified-gluing-estimate}
	\begin{aligned}
	\|\widetilde D_{t,q}^M\bar S_{q,J_*,u}\|_{N,\mcI_u^{\mathrm N}}
	&\lesssim e_q\MM{M,N_t^{\mathrm N},\tau_q^{-1},\ell_{t,q}^{-1}}
	\MM{N,N_{R,J_*}^{(1)}(M),\lambda_q,\ell_q^{-1}}.
	\end{aligned}
	\end{equation}
	Apply \eqref{eq:holder-product-tame} with $\alpha=0$ to \eqref{def:Nash-pulled-back-matrix}.
	For $M\leq N_{\mathrm{ind}}$ and $M+N\leq K_{J_*}^{(0)}+1$, the resulting estimate is
	\begin{equation}\label{eq:Nash-pulled-back-matrix-integer-estimate}
	\begin{aligned}
	\|\widetilde D_{t,q}^M\widetilde R_{q,u}\|_{N,\mcI_u^{\mathrm N}}
	&\lesssim\sum_{\substack{M_1+M_2+M_3=M\\N_1+N_2+N_3=N}}
	\Nrm{\widetilde D_{t,q}^{M_1}[(\mathbf C_u^{\mathrm N})^{-T}]}_{N_1,\mcI_u^{\mathrm N}}\Nrm{\widetilde D_{t,q}^{M_2}\bigl(\Id-e_q^{-1}\bar S_{q,J_*,u}\bigr)}_{N_2,\mcI_u^{\mathrm N}}
	\Nrm{\widetilde D_{t,q}^{M_3}[(\mathbf C_u^{\mathrm N})^{-1}]}_{N_3,\mcI_u^{\mathrm N}}\\
	&\lesssim\MM{M,N_t^{\mathrm N},\tau_q^{-1},\ell_{t,q}^{-1}}
	\MM{N,N_{R,J_*}^{(1)}(M),\lambda_q,\ell_q^{-1}}.
	\end{aligned}
	\end{equation}
	Here the last inequality uses \eqref{eq:Nash-flow-determinant-cofactor-CN-estimate} and
	\eqref{eq:Nash-time-mollified-gluing-estimate}.
	The flow estimates apply since $K_{J_*}^{(0)}+1\leq N_{\mathrm{max}}-1$.

	By \eqref{eq:Nash-flow-zero-order-estimate} and \eqref{eq:Nash-normalized-gluing-stress-smallness},
	we have $\widetilde R_{q,u}\in B_{\varepsilon_0/2}(\Id)$.
	The derivatives of $\gamma_\xi$ from Lemma~\ref{Geometric lemma} are bounded on this ball.
	We estimate the resulting spatial compositions by \eqref{eq:faa-di-bruno} and \eqref{est:holder-norm-composition}. At spatial order zero, we use the bounded derivatives of $\gamma_\xi$.
	For $1\leq M+N\leq K_{J_*}^{(0)}+1$ and $M\leq N_{\mathrm{ind}}$, the differentiated terms satisfy
	$$
	\begin{aligned}
		\Nrm{\nabla^N\widetilde D_{t,q}^M\bigl(\gamma_\xi(\widetilde R_{q,u})\bigr)}_{0,\mcI_u^{\mathrm N}}\lesssim\sum_{k=1}^{M+N}
		\sum_{\substack{M_1+\cdots+M_k=M\\N_1+\cdots+N_k=N\\M_i+N_i\geq1}}
		\prod_{i=1}^k\Nrm{\widetilde D_{t,q}^{M_i}\widetilde R_{q,u}}_{N_i,\mcI_u^{\mathrm N}}.
	\end{aligned}
	$$
	Estimates~\eqref{eq:Nash-pulled-back-matrix-integer-estimate} and
	\eqref{eq:appendix-derivative-weight-product}, together with
	$\|\gamma_\xi(\widetilde R_{q,u})\|_{0,\mcI_u^{\mathrm N}}\lesssim1$, give, for
	$M\leq N_{\mathrm{ind}}$ and $M+N\leq K_{J_*}^{(0)}+1$,
	\begin{equation}\label{eq:Nash-geometric-coefficient-estimate}
	\begin{aligned}
		\Nrm{\widetilde D_{t,q}^M\bigl(\gamma_\xi(\widetilde R_{q,u})\bigr)}_{N,\mcI_u^{\mathrm N}}
		&\lesssim\MM{M,N_t^{\mathrm N},\tau_q^{-1},\ell_{t,q}^{-1}}
		\MM{N,N_{R,J_*}^{(1)}(M),\lambda_q,\ell_q^{-1}}.
	\end{aligned}
	\end{equation}

	Differentiating \eqref{def:Nash-directional-amplitude} gives, in the same derivative range,
	$$
	\begin{aligned}
	\|\widetilde D_{t,q}^M\tilde a_{q,u,\xi}\|_{N,\mcI_u^{\mathrm N}}
	&\lesssim e_q^{\frac12}\sum_{\substack{r+M_1+M_2=M\\N_1+N_2=N}}
	\|\partial_t^r\theta_u^{\mathrm N}\|_{0,\mcI_u^{\mathrm N}}\Nrm{\widetilde D_{t,q}^{M_1}[(\bar\rho_{q,J_*})^{\frac12}]}_{N_1,\mcI_u^{\mathrm N}}
	\Nrm{\widetilde D_{t,q}^{M_2}\bigl(\gamma_\xi(\widetilde R_{q,u})\bigr)}_{N_2,\mcI_u^{\mathrm N}}\\
	&\lesssim A_{q,0}^{\frac12}(\tau_q\ell_q^{-1})^{\frac\alpha2}
	\MM{M,N_t^{\mathrm N},\tau_q^{-1},\ell_{t,q}^{-1}}
	\MM{N,N_{R,J_*}^{(1)}(M),\lambda_q,\ell_q^{-1}}.
	\end{aligned}
	$$
	In the last inequality, we used \eqref{eq:Nash-mollified-density-coefficient-estimate} with $s=\tfrac12$,
	\eqref{def:Nash-transition-cutoff}, and
	\eqref{eq:Nash-geometric-coefficient-estimate}, together with
	$e_q^{\frac12}=A_{q,0}^{\frac12}(\tau_q\ell_q^{-1})^{\frac{\alpha}{2}}$.
	For the localized amplitude, \eqref{def:Nash-directional-amplitude} and $\widetilde D_{t,q}\Phi_u^{\mathrm N}=0$ give
	$$
	\widetilde D_{t,q}^M\tilde a_{q,u,\kappa,\xi}
	=(\chi_\kappa^{\mathrm N})^2(\Phi_u^{\mathrm N})\widetilde D_{t,q}^M\tilde a_{q,u,\xi}.
	$$
	Use \eqref{eq:Nash-transported-cutoff-estimate} with $\theta_u^{\mathrm N}$ replaced by $1$ and $\chi_\kappa^{\mathrm N}$ by $(\chi_\kappa^{\mathrm N})^2$.
	Then \eqref{eq:holder-product-tame} with $\alpha=0$ gives
	$$
	\begin{aligned}
	\|\widetilde D_{t,q}^M\tilde a_{q,u,\kappa,\xi}\|_{N,\mcI_u^{\mathrm N}}
	&\lesssim\sum_{N_0=0}^N
	\Nrm{(\chi_\kappa^{\mathrm N})^2(\Phi_u^{\mathrm N})}_{N-N_0,\mcI_u^{\mathrm N}}
	\|\widetilde D_{t,q}^M\tilde a_{q,u,\xi}\|_{N_0,\mcI_u^{\mathrm N}}\\
	&\lesssim A_{q,0}^{\frac12}(\tau_q\ell_q^{-1})^{\frac\alpha2}
	\MM{M,N_t^{\mathrm N},\tau_q^{-1},\ell_{t,q}^{-1}}
	\MM{N,N_{R,J_*}^{(1)}(M),\lambda_q,\ell_q^{-1}}.
	\end{aligned}
	$$
	Estimate~\eqref{eq:Nash-amplitude-estimate} follows.

	To prove \eqref{eq:Nash-principal-estimate}, take $M\leq N_{\mathrm{ind}}$ and $M+N\leq K_{J_*}^{(0)}$.
	Since $\widetilde D_{t,q}\Phi_u^{\mathrm N}=0$, material derivatives in \eqref{def:Nash-principal-perturbation}
	act only on $\tilde a_{q,u,\xi}$ and $\mathbf C_u^{\mathrm N}$.
	For $1\leq N\leq K_{J_*}^{(0)}$, the second inequality in \eqref{est:holder-norm-composition} gives
	$$
	\begin{aligned}
		\Nrm{\nabla^N\bigl[U_\xi(\lambda_{q+1}\Phi_u^{\mathrm N})\bigr]}_{0,\mcI_u^{\mathrm N}}
		+\Nrm{\nabla^N\bigl[\widetilde U_\xi(\lambda_{q+1}\Phi_u^{\mathrm N})\bigr]}_{0,\mcI_u^{\mathrm N}}\lesssim\lambda_{q+1}\|\nabla\Phi_u^{\mathrm N}\|_{N-1,\mcI_u^{\mathrm N}}
		+\lambda_{q+1}^N\|\nabla\Phi_u^{\mathrm N}\|_{0,\mcI_u^{\mathrm N}}^N
		\lesssim\lambda_{q+1}^N.
	\end{aligned}
	$$
	For the last inequality, we used \eqref{eq:Nash-flow-higher-order-CN-estimate} at material order zero.
	Including the lower spatial orders and the $C^0$ composition estimate, we obtain, for $0\leq N\leq K_{J_*}^{(0)}$,
	\begin{equation}\label{eq:Nash-transported-profile-spatial-estimate}
	\Nrm{U_\xi(\lambda_{q+1}\Phi_u^{\mathrm N})}_{N,\mcI_u^{\mathrm N}}
	+\Nrm{\widetilde U_\xi(\lambda_{q+1}\Phi_u^{\mathrm N})}_{N,\mcI_u^{\mathrm N}}
	\lesssim\lambda_{q+1}^N.
	\end{equation}
	Estimates~\eqref{eq:Nash-flow-determinant-cofactor-CN-estimate}, \eqref{eq:Nash-amplitude-estimate}, and
	\eqref{eq:holder-product-tame} with $\alpha=0$ give
	\begin{equation}\label{eq:Nash-principal-integer-estimate}
	\begin{aligned}
		\|\widetilde D_{t,q}^Mw_{q+1,u}^{\mathrm p}\|_N
		&\lesssim\sum_{\xi\in\Lambda}\sum_{\substack{M_1+M_2=M\\N_1+N_2+N_3=N}}
		\Nrm{\widetilde D_{t,q}^{M_1}\bigl((\mathbf C_u^{\mathrm N})^T\bigr)}_{N_1,\mcI_u^{\mathrm N}}
		\nrm{\widetilde D_{t,q}^{M_2}\tilde a_{q,u,\xi}}_{N_2,\mcI_u^{\mathrm N}}
		\Nrm{U_\xi(\lambda_{q+1}\Phi_u^{\mathrm N})}_{N_3,\mcI_u^{\mathrm N}}\\
		&\lesssim e_q^{\frac12}
		\sum_{\substack{M_1+M_2=M\\N_1+N_2+N_3=N}}
		\ell_q^{-N_1-N_2}\lambda_{q+1}^{N_3}
		\prod_{i=1}^2\MM{M_i,N_t^{\mathrm N},\tau_q^{-1},\ell_{t,q}^{-1}}\\
		&\lesssim A_{q,0}^{\frac12}(\tau_q\ell_q^{-1})^{\frac\alpha2}
		\MM{M,N_t^{\mathrm N},\tau_q^{-1},\ell_{t,q}^{-1}}\lambda_{q+1}^N.
	\end{aligned}
	\end{equation}
	Here the last inequality uses \eqref{eq:appendix-derivative-weight-product} and $\ell_q^{-1}\leq\lambda_{q+1}$.
	Summing over $u$ and using the disjoint temporal supports gives
	\eqref{eq:Nash-principal-estimate}.

	For \eqref{eq:Nash-corrector-estimate}, take $M\leq N_{\mathrm{ind}}$ and $M+N\leq K_{J_*}^{(0)}$.
	For a smooth scalar amplitude $\tilde a$ and $M\in\mathbb N_0$, we have
	\begin{equation}\label{eq:Nash-amplitude-gradient-material-identity}
		\widetilde D_{t,q}^M\bigl((\nabla\Phi_u^{\mathrm N})^{-T}\nabla \tilde a\bigr)
		=(\nabla\Phi_u^{\mathrm N})^{-T}\nabla\widetilde D_{t,q}^M\tilde a.
	\end{equation}
	This follows from the commutation identity in \eqref{eq:Nash-Piola-material-identities} with $\Phi=\Phi_u^{\mathrm N}$.
	The cross product satisfies
	\begin{equation}\label{eq:Nash-cofactor-cross-product-identity}
	\nabla \tilde a\times((\nabla\Phi_u^{\mathrm N})^TZ)
	=(\mathbf C_u^{\mathrm N})^T\bigl(((\nabla\Phi_u^{\mathrm N})^{-T}\nabla \tilde a)\times Z\bigr).
	\end{equation}
	Thus the $C^N$ estimate for $w_{q+1}^{\mathrm c}$ requires total order $M+N+1\leq K_{J_*}^{(0)}+1$, with no additional material derivative.
	For $0\leq M_0\leq M$ and $N\in\N_0$, we have
	\begin{equation}\label{eq:Nash-slow-spatial-derivative-factor}
	\MM{N+1,N_{R,J_*}^{(1)}(M_0),\lambda_q,\ell_q^{-1}}
	\leq\MM{1,N_{R,J_*}^{(1)}(M),\lambda_q,\ell_q^{-1}}\ell_q^{-N}.
	\end{equation}
	For $\tilde a=\tilde a_{q,u,\xi}$ or $\tilde a=\tilde a_{q,u,\kappa,\xi}$,
	\eqref{eq:Nash-amplitude-estimate}, \eqref{eq:Nash-amplitude-gradient-material-identity}, and
	\eqref{eq:Nash-slow-spatial-derivative-factor} give
	\begin{equation}\label{eq:Nash-amplitude-transported-gradient-estimate}
	\begin{aligned}
		\Nrm{\widetilde D_{t,q}^M\bigl((\nabla\Phi_u^{\mathrm N})^{-T}\nabla \tilde a\bigr)}_{N,\mcI_u^{\mathrm N}}
		&\lesssim\sum_{N_1+N_2=N}
		\Nrm{(\nabla\Phi_u^{\mathrm N})^{-T}}_{N_1,\mcI_u^{\mathrm N}}
		\nrm{\widetilde D_{t,q}^M\tilde a}_{N_2+1,\mcI_u^{\mathrm N}}\\
		&\lesssim A_{q,0}^{\frac12}(\tau_q\ell_q^{-1})^{\frac\alpha2}
		\MM{M,N_t^{\mathrm N},\tau_q^{-1},\ell_{t,q}^{-1}}
		\MM{1,N_{R,J_*}^{(1)}(M),\lambda_q,\ell_q^{-1}}\ell_q^{-N}.
	\end{aligned}
	\end{equation}
	Here the flow factor is estimated by \eqref{eq:Nash-flow-higher-order-CN-estimate} with material order zero.
	Use \eqref{eq:Nash-cofactor-cross-product-identity} and \eqref{eq:Nash-amplitude-transported-gradient-estimate} in
	\eqref{eq:Nash-corrector-formula}, together with \eqref{eq:Nash-transported-profile-spatial-estimate}.
	These estimates give, uniformly in $u$,
	$$
	\begin{aligned}
		\|\widetilde D_{t,q}^Mw_{q+1,u}^{\mathrm c}\|_N
		&\lesssim\lambda_{q+1}^{-1}\sum_{\xi\in\Lambda}\sum_{\substack{M_1+M_2=M\\N_1+N_2+N_3=N}}
		\Nrm{\widetilde D_{t,q}^{M_1}\bigl((\mathbf C_u^{\mathrm N})^T\bigr)}_{N_1,\mcI_u^{\mathrm N}}\Nrm{\widetilde D_{t,q}^{M_2}\bigl((\nabla\Phi_u^{\mathrm N})^{-T}\nabla \tilde a_{q,u,\xi}\bigr)}_{N_2,\mcI_u^{\mathrm N}}
		\Nrm{\widetilde U_\xi(\lambda_{q+1}\Phi_u^{\mathrm N})}_{N_3,\mcI_u^{\mathrm N}}\\
		&\lesssim A_{q,0}^{\frac12}(\tau_q\ell_q^{-1})^{\frac\alpha2}
		\MM{M,N_t^{\mathrm N},\tau_q^{-1},\ell_{t,q}^{-1}}
		\MM{1,N_{R,J_*}^{(1)}(M),\lambda_q,\ell_q^{-1}}\lambda_{q+1}^{N-1}.
	\end{aligned}
	$$
	Summing over $u$ and using the disjoint temporal supports proves \eqref{eq:Nash-corrector-estimate}.
	Adding \eqref{eq:Nash-principal-estimate}, we obtain \eqref{eq:Nash-total-perturbation-estimate}.
	The factor $\theta_u^{\mathrm N}$ in each summand gives \eqref{eq:Nash-perturbation-temporal-support}.
\end{proof}

\begin{prop}
	\label{prop:Nash-stress-estimates}
	Under the assumptions of Proposition~\ref{prop:Nash-perturbation-estimates},
	let $M,N\in\mathbb N_0$.
	The transport and oscillation estimates below hold for $M\leq N_{\mathrm{ind}}-1$
	and, respectively, $M+N\leq K_{J_*}^{(0)}-3$ and $M+N\leq K_{J_*}^{(0)}-2$.
	Estimate~\eqref{eq:Nash-corrector-stress-estimate} holds for all $M+N\leq N_{\mathrm{ind}}$.
	\begin{align}
		\nrm{\widetilde D_{t,q}^MR_{q+1}^{\mathrm{trans}}}_{N}
		&\lesssim
		A_{q,0}^{\frac12}
		(\tau_q\ell_q^{-1})^{\frac\alpha2}\tau_q^{-1}\MM{M,N_t^{\mathrm N}-1,\tau_q^{-1},\ell_{t,q}^{-1}}
		\lambda_{q+1}^{N-1},
		\label{eq:Nash-transport-stress-estimate}\\
		\nrm{\widetilde D_{t,q}^MR_{q+1}^{\mathrm{osc}}}_{N}
		&\lesssim
		A_{q,0}(\tau_q\ell_q^{-1})^\alpha\MM{M,N_t^{\mathrm N},\tau_q^{-1},\ell_{t,q}^{-1}}
		\MM{1,N_{R,J_*}^{(1)}(M),\lambda_q,\ell_q^{-1}}
		\lambda_{q+1}^{N-1},
		\label{eq:Nash-oscillation-stress-estimate}\\
		\nrm{\widetilde D_{t,q}^MR_{q+1}^{\mathrm{corr}}}_{N}
		&\lesssim
		A_{q,0}(\tau_q\ell_q^{-1})^\alpha\MM{M,N_t^{\mathrm N},\tau_q^{-1},\ell_{t,q}^{-1}}
		\MM{1,N_{R,J_*}^{(1)}(M),\lambda_q,\ell_q^{-1}}
		\lambda_{q+1}^{N-1}.
		\label{eq:Nash-corrector-stress-estimate}
	\end{align}
	For $M+N\leq N_{\mathrm{ind}}-4$, the mollification estimates for $M=0$ and $M\geq1$ are, respectively,
	\begin{align}
		\nrm{R_{q+1}^{\mathrm{mol}}}_{N}
		&\lesssim
		A_{q,0}(\tau_q\ell_q^{-1})^\alpha\ell_q^4\lambda_q^4
		\lambda_{q+1}^{N}+A_{q,0}^{\frac12}(\tau_q\ell_q^{-1})^{\frac\alpha2}
		\tau_q^{-1}\lambda_{q+1}^{N-1},
		\label{eq:Nash-flow-stress-zero-material-estimate}\\
		\nrm{\widetilde D_{t,q}^MR_{q+1}^{\mathrm{mol}}}_{N}
		&\lesssim
		\MM{M,N_t^{\mathrm N}-1,\tau_q^{-1},\ell_{t,q}^{-1}}
		\Bigl(A_{q,0}(\tau_q\ell_q^{-1})^\alpha\ell_q^4\lambda_q^4\lambda_{q+1}^{N}
		+A_{q,0}^{\frac12}(\tau_q\ell_q^{-1})^{\frac\alpha2}
		\tau_q^{-1}\lambda_{q+1}^{N-1}\Bigr).
		\label{eq:Nash-flow-stress-positive-material-estimate}
	\end{align}

	For $M+N\leq N_R$, the full stress $R_{q+1}$ satisfies
	\begin{equation}
		\begin{aligned}
			\nrm{\widetilde D_{t,q}^MR_{q+1}}_N
			\lesssim\lambda_q^{-\alpha}\lambda_{q+1}^{-8\alpha}\delta_{q+2}
			\ell_{t,q}^{-M}\ell_q^{-N}
			+\bigl(A_{q,0}^{\frac12}(\tau_q\ell_q^{-1})^{\frac\alpha2}\tau_q^{-1}
			+A_{q,0}(\tau_q\ell_q^{-1})^\alpha\lambda_q\bigr)
			\MM{M,N_t^{\mathrm N}-1,\tau_q^{-1},\ell_{t,q}^{-1}}\lambda_{q+1}^{N-1}.
		\end{aligned}
		\label{eq:closure-Reynolds-background-material}
	\end{equation}
	The temporal support of $R_{q+1}$ satisfies
	\begin{equation}\label{eq:closure-stress-support-enlargement}
	\begin{aligned}
		\supp_tR_{q+1}
		\subset\supp_tR_{q,J_*}\cup\supp_tP_{q+1,J_*}
		\cup\bigcup_u\supp\theta_u^{\mathrm N}
		\subset\mN_{(2J_*+2)\tau_q}(\supp_tR_q).
	\end{aligned}
	\end{equation}
	The implicit constants may depend on
	$\beta,b,\alpha,\rho^*,\mathfrak M,C_*,N_R,N_{\mathrm{ind}},N_{\mathrm{max}}$.
	They are independent of $q$, $a$, $u$, and $\kappa$, and uniform in the stated derivative ranges.
\end{prop}

\begin{proof}[Proof of Proposition~\ref{prop:Nash-stress-estimates}]
	We first estimate the coefficients in the transport and oscillation sources.
	Fix $(u,\kappa)$.
	For $\xi\in\Lambda$ and $1\leq j\leq3$, define the vector amplitudes
	$$
	\begin{aligned}
		h_\xi^{\mathrm p}
		&:=\widetilde D_{t,q}\bigl[(\nabla\Phi_u^{\mathrm N})^{-1}\tilde a_{q,u,\kappa,\xi}\xi\bigr]+\tilde a_{q,u,\kappa,\xi}\bigl((\nabla\Phi_u^{\mathrm N})^{-1}\xi\cdot\nabla\bigr)\bar u_{q,J_*},\\
		h_{\xi,j}^{\mathrm c}
		&:=\lambda_{q+1}^{-1}\widetilde D_{t,q}\Bigl[
		(\nabla\Phi_u^{\mathrm N})^{-1}
		\bigl(((\nabla\Phi_u^{\mathrm N})^{-T}\nabla\tilde a_{q,u,\kappa,\xi})\times e_j\bigr)\Bigr]+\lambda_{q+1}^{-1}\Bigl[
		(\nabla\Phi_u^{\mathrm N})^{-1}
		\bigl(((\nabla\Phi_u^{\mathrm N})^{-T}\nabla\tilde a_{q,u,\kappa,\xi})\times e_j\bigr)\cdot\nabla\Bigr]\bar u_{q,J_*},
	\end{aligned}
	$$
	where $e_j$ are the standard coordinate vectors.
	By \eqref{Dt Jm} with $\sigma=1$, we have
	$\widetilde D_{t,q}J_u^{\mathrm N}=-(\Div\bar u_{q,J_*})J_u^{\mathrm N}$.
	Using $U_\xi=\xi\psi_\xi$, $\Div w_{q+1,u,\kappa}=0$, and $\widetilde D_{t,q}\Phi_u^{\mathrm N}=0$ in
	\eqref{def:Nash-transport-source}, together with \eqref{eq:Nash-cofactor-cross-product-identity}, we obtain
	\begin{equation}\label{eq:Nash-transport-source-decomposition}
	\begin{aligned}
		(J_u^{\mathrm N})^{-1}F_{q+1,u,\kappa}^{\mathrm{trans}}
		&=\widetilde D_{t,q}\bigl[(J_u^{\mathrm N})^{-1}w_{q+1,u,\kappa}\bigr]+\bigl((J_u^{\mathrm N})^{-1}w_{q+1,u,\kappa}\cdot\nabla\bigr)\bar u_{q,J_*}\\
		&=\sum_{\xi\in\Lambda}h_\xi^{\mathrm p}\psi_\xi(\lambda_{q+1}\Phi_u^{\mathrm N})+\sum_{\xi\in\Lambda}\sum_{j=1}^3h_{\xi,j}^{\mathrm c}(\widetilde U_\xi)_j(\lambda_{q+1}\Phi_u^{\mathrm N}).
	\end{aligned}
	\end{equation}

	Take $M,N\in\mathbb N_0$ with $M\leq N_{\mathrm{ind}}-1$ and $M+N\leq K_{J_*}^{(0)}-1$. By
	\eqref{eq:Nash-flow-higher-order-CN-estimate}, \eqref{eq:Nash-terminal-velocity-gradient-Nash-scale},
	\eqref{eq:Nash-amplitude-estimate}, and \eqref{eq:Nash-amplitude-transported-gradient-estimate}, we have the following estimates.
	All norms in the next display are taken over $\mcI_u^{\mathrm N}$.
	$$
	\begin{aligned}
		\sum_{\xi\in\Lambda}\|\widetilde D_{t,q}^{M}h_\xi^{\mathrm p}\|_{N}
		&\lesssim\sum_{\xi\in\Lambda}
		\sum_{\substack{M_1+M_2=M+1\\N_1+N_2=N}}
		\Nrm{\widetilde D_{t,q}^{M_1}(\nabla\Phi_u^{\mathrm N})^{-1}}_{N_1}
		\|\widetilde D_{t,q}^{M_2}\tilde a_{q,u,\kappa,\xi}\|_{N_2}\\
		&\quad+\sum_{\xi\in\Lambda}
		\sum_{\substack{M_1+M_2+M_3=M\\N_1+N_2+N_3=N}}
		\Nrm{\widetilde D_{t,q}^{M_1}(\nabla\Phi_u^{\mathrm N})^{-1}}_{N_1}
		\|\widetilde D_{t,q}^{M_2}\tilde a_{q,u,\kappa,\xi}\|_{N_2}
		\|\widetilde D_{t,q}^{M_3}\nabla\bar u_{q,J_*}\|_{N_3}\\
		&\lesssim e_q^{\frac12}\tau_q^{-1}
		\MM{M,N_t^{\mathrm N}-1,\tau_q^{-1},\ell_{t,q}^{-1}}\ell_q^{-N},\\
		\sum_{\xi\in\Lambda}\sum_{j=1}^3\|\widetilde D_{t,q}^{M}h_{\xi,j}^{\mathrm c}\|_{N}
		&\lesssim\lambda_{q+1}^{-1}\sum_{\xi\in\Lambda}
		\sum_{\substack{M_1+M_2=M+1\\N_1+N_2=N}}
		\Nrm{\widetilde D_{t,q}^{M_1}(\nabla\Phi_u^{\mathrm N})^{-1}}_{N_1}
		\Nrm{\widetilde D_{t,q}^{M_2}\bigl[
		(\nabla\Phi_u^{\mathrm N})^{-T}\nabla\tilde a_{q,u,\kappa,\xi}\bigr]}_{N_2}\\
		&\quad+\lambda_{q+1}^{-1}\sum_{\xi\in\Lambda}
		\sum_{\substack{M_1+M_2+M_3=M\\N_1+N_2+N_3=N}}
		\Nrm{\widetilde D_{t,q}^{M_1}(\nabla\Phi_u^{\mathrm N})^{-1}}_{N_1}
		\Nrm{\widetilde D_{t,q}^{M_2}\bigl[
		(\nabla\Phi_u^{\mathrm N})^{-T}\nabla\tilde a_{q,u,\kappa,\xi}\bigr]}_{N_2}
		\|\widetilde D_{t,q}^{M_3}\nabla\bar u_{q,J_*}\|_{N_3}\\
		&\lesssim e_q^{\frac12}\tau_q^{-1}
		\MM{M,N_t^{\mathrm N}-1,\tau_q^{-1},\ell_{t,q}^{-1}}
		\lambda_{q+1}^{-1}\ell_q^{-N-1}.
	\end{aligned}
	$$
	Combining these estimates gives
	$$
	\begin{aligned}
		\sum_{\xi\in\Lambda}\|\widetilde D_{t,q}^{M}h_\xi^{\mathrm p}\|_{N,\mcI_u^{\mathrm N}}
		+\sum_{\xi\in\Lambda}\sum_{j=1}^3\|\widetilde D_{t,q}^{M}h_{\xi,j}^{\mathrm c}\|_{N,\mcI_u^{\mathrm N}}
		&\lesssim e_q^{\frac12}\tau_q^{-1}
		\MM{M,N_t^{\mathrm N}-1,\tau_q^{-1},\ell_{t,q}^{-1}}
		\ell_q^{-N}\left(1+\lambda_{q+1}^{-1}\ell_q^{-1}\right)\\
		&\lesssim A_{q,0}^{\frac12}(\tau_q\ell_q^{-1})^{\frac\alpha2}\tau_q^{-1}
		\MM{M,N_t^{\mathrm N}-1,\tau_q^{-1},\ell_{t,q}^{-1}}\ell_q^{-N}.
	\end{aligned}
	$$
	Here we used, for $0\leq M_1\leq M$,
	$$
	\MM{M_1,N_t^{\mathrm N},\tau_q^{-1},\ell_{t,q}^{-1}}
	\MM{M-M_1+1,N_t^{\mathrm N},\tau_q^{-1},\ell_{t,q}^{-1}}
	\leq\tau_q^{-1}\MM{M,N_t^{\mathrm N}-1,\tau_q^{-1},\ell_{t,q}^{-1}}.
	$$
	The differentiated amplitudes require material order at most $M+1\leq N_{\mathrm{ind}}$ and
	total order at most $M+N+2\leq K_{J_*}^{(0)}+1$.
	For the transport source, choose
	$N_t=N_t^{\mathrm N}-1$ and $\mC_f=A_{q,0}^{\frac12}(\tau_q\ell_q^{-1})^{\frac\alpha2}\tau_q^{-1}$.

	For the oscillation source, \eqref{def:Nash-directional-amplitude} gives
	$
	\bar\rho_{q,J_*}^{-1}\tilde a_{q,u,\xi}^2
	=e_q(\theta_u^{\mathrm N})^2\gamma_\xi^2(\widetilde R_{q,u}).
	$
	Using this identity, we write \eqref{eq:Nash-oscillation-profile-source} as
	\begin{equation}\label{eq:Nash-oscillation-source-decomposition}
	\begin{aligned}
	    (J_u^{\mathrm N})^{-1}F_{q+1,u,\kappa}^{\mathrm{osc}}
	    &=\sum_{\xi\in\Lambda}h_\xi^{\mathrm o} g_\xi(\lambda_{q+1}\Phi_u^{\mathrm N}),\\
		(h_\xi^{\mathrm o})_i&:=e_q(J_u^{\mathrm N})^{-1}(\mathbf C_u^{\mathrm N})_{sl}\partial_l\left(
		(\theta_u^{\mathrm N})^2(\chi_\kappa^{\mathrm N})^2(\Phi_u^{\mathrm N})
		\gamma_\xi^2(\widetilde R_{q,u})(\mathbf C_u^{\mathrm N})_{ri}\xi_r\xi_s\right).
	\end{aligned}
	\end{equation}
	By \eqref{eq:Nash-Piola-material-identities},
	$\widetilde D_{t,q}$ commutes with $(J_u^{\mathrm N})^{-1}(\mathbf C_u^{\mathrm N})_{sl}\partial_l$.
	Using \eqref{eq:Nash-transported-cutoff-estimate}, \eqref{eq:Nash-flow-mixed-integer-estimate},
	\eqref{eq:Nash-geometric-coefficient-estimate}, and \eqref{eq:Nash-slow-spatial-derivative-factor},
	we obtain, for $M,N\in\mathbb N_0$ with $M\leq N_{\mathrm{ind}}-1$ and $M+N\leq K_{J_*}^{(0)}$,
	$$
	\begin{aligned}
		\sum_{\xi\in\Lambda}\|\widetilde D_{t,q}^{M}h_\xi^{\mathrm o}\|_{N,\mcI_u^{\mathrm N}}
		&\lesssim e_q\sum_{\xi\in\Lambda}\sum_{N_1+N_2=N}
		\|(\nabla\Phi_u^{\mathrm N})^{-T}\|_{N_1,\mcI_u^{\mathrm N}}\Nrm{\widetilde D_{t,q}^{M}\left[
		(\theta_u^{\mathrm N})^2(\chi_\kappa^{\mathrm N})^2(\Phi_u^{\mathrm N})
		\gamma_\xi^2(\widetilde R_{q,u})(\mathbf C_u^{\mathrm N})^T\right]}_{N_2+1,\mcI_u^{\mathrm N}}\\
		&\lesssim e_q\MM{M,N_t^{\mathrm N},\tau_q^{-1},\ell_{t,q}^{-1}}
		\sum_{N_1+N_2=N}\|(\nabla\Phi_u^{\mathrm N})^{-T}\|_{N_1,\mcI_u^{\mathrm N}}
		\MM{N_2+1,N_{R,J_*}^{(1)}(M),\lambda_q,\ell_q^{-1}}\\
		&\lesssim e_q\MM{M,N_t^{\mathrm N},\tau_q^{-1},\ell_{t,q}^{-1}}
		\MM{1,N_{R,J_*}^{(1)}(M),\lambda_q,\ell_q^{-1}}\ell_q^{-N}.
	\end{aligned}
	$$
	The coefficient estimates require total order at most $M+N+1\leq K_{J_*}^{(0)}+1$.
	For this source, choose
	$
		N_t=N_t^{\mathrm N},\ \mC_f=A_{q,0}(\tau_q\ell_q^{-1})^\alpha \MM{1,N_{R,J_*}^{(1)}(M),\lambda_q,\ell_q^{-1}}.
	$
	
	We now apply Lemma~\ref{lem:Nash-transported-oscillatory-inverse-divergence} to the two sources.
	Take $M\leq N_{\mathrm{ind}}-1$, with $M+N\leq K_{J_*}^{(0)}-3$ for the transport source and $M+N\leq K_{J_*}^{(0)}-2$ for the oscillation source.
	In both cases, use $K^*=M$, $K^\sharp=M+N+2$, the parameters fixed before \eqref{def:Nash-cell-oscillatory-inverse},
	and the choices of $N_t$ and $\mC_f$ made above for each source.
	The amplitude estimates cover total order $K^\sharp$.
	For the oscillation source, the choice of $\mC_f$ covers all lower material orders by \eqref{eq:Nash-slow-spatial-derivative-factor}.
	The velocity and flow hypotheses follow from \eqref{eq:Nash-terminal-Newton-velocity-sharp-gradient}
	and \eqref{eq:Nash-flow-mixed-integer-estimate}.
	The scale, support, and mean conditions were verified in Subsection~\ref{subsec:Nash-new-Reynolds-stress}.
	Summing and using \eqref{eq:Nash-local-inverse-cylinder-overlap} proves
	\eqref{eq:Nash-transport-stress-estimate}--\eqref{eq:Nash-oscillation-stress-estimate}.

	For $M+N\leq N_{\mathrm{ind}}$, estimate $\rho_{q,J_*}^{-1}$ by
	\eqref{eq:Newton-state-zero-order-bounds} at order zero and
	\eqref{eq:Nash-density-coefficient-material-estimate} with $s=-1$ at positive orders.
	Applying \eqref{eq:holder-product-tame} with $\alpha=0$ to \eqref{def:Nash-corrector-stress}
	and using \eqref{eq:Nash-principal-estimate}--\eqref{eq:Nash-corrector-estimate}, we obtain
	$$
	\begin{aligned}
	    \nrm{\widetilde D_{t,q}^MR_{q+1}^{\mathrm{corr}}}_N
	    &\lesssim\sum_{\substack{M_0+M_1+M_2=M\\N_0+N_1+N_2=N}}
	    \|\widetilde D_{t,q}^{M_0}\rho_{q,J_*}^{-1}\|_{N_0}
	    \|\widetilde D_{t,q}^{M_1}w_{q+1}^{\mathrm c}\|_{N_1}\left(\|\widetilde D_{t,q}^{M_2}w_{q+1}^{\mathrm p}\|_{N_2}
	    +\|\widetilde D_{t,q}^{M_2}w_{q+1}^{\mathrm c}\|_{N_2}\right)\\
	    &\lesssim A_{q,0}(\tau_q\ell_q^{-1})^\alpha
	    \MM{1,N_{R,J_*}^{(1)}(M),\lambda_q,\ell_q^{-1}}\lambda_{q+1}^{N-1}\left(1+\ell_q^{-1}\lambda_{q+1}^{-1}\right)
	    \sum_{M_0+M_1+M_2=M}
	    \prod_{i=0}^{2}\MM{M_i,N_t^{\mathrm N},\tau_q^{-1},\ell_{t,q}^{-1}}\\
	    &\lesssim A_{q,0}(\tau_q\ell_q^{-1})^\alpha
	    \MM{M,N_t^{\mathrm N},\tau_q^{-1},\ell_{t,q}^{-1}}
	    \MM{1,N_{R,J_*}^{(1)}(M),\lambda_q,\ell_q^{-1}}\lambda_{q+1}^{N-1}.
	\end{aligned}
	$$
	Thus \eqref{eq:Nash-corrector-stress-estimate} holds.

	For the mollification error, \eqref{def:Nash-pulled-back-matrix}, \eqref{eq:Nash-principal-quadratic-average}, and
	\eqref{def:Nash-flow-mollification-stress} give
	\begin{equation}\label{eq:Nash-flow-stress-pulled-back-identity}
	\begin{aligned}
	R_{q+1}^{\mathrm{mol}}
	={}&\sum_u(\theta_u^{\mathrm N})^2
	\left(S_{q,J_*,u}-\bar S_{q,J_*,u}\right)\\
	&+(u_{q,J_*}-\bar u_{q,J_*})\otimes w_{q+1}
	+w_{q+1}\otimes(u_{q,J_*}-\bar u_{q,J_*})
	-\frac{\rho_{q,J_*}-\bar\rho_{q,J_*}}{\rho_{q,J_*}\bar\rho_{q,J_*}}
	w_{q+1}^{\mathrm p}\otimes w_{q+1}^{\mathrm p}.
	\end{aligned}
	\end{equation}
	Here $\theta_u^{\mathrm N}=1$ on the supports of $S_{q,J_*,u}$ and $\bar S_{q,J_*,u}$.

	For the velocity difference, the following comparison holds for every $M\in\mathbb N_0$,
\begin{equation*}
	\begin{aligned}
		\frac{\ell_q^4\delta_q^{\frac12}\MM{M,N_t^{\mathrm N}+4,\tau_q^{-1},\ell_{t,q}^{-1}}\MM{4,N_{R,J_*}^{(1)}(M),\lambda_q,\ell_q^{-1}}}
		{\lambda_{q+1}^{-1}\tau_q^{-1}\MM{M,N_t^{\mathrm N}-1,\tau_q^{-1},\ell_{t,q}^{-1}}}
		\leq\frac{\ell_q^4\lambda_q^4\delta_q^{\frac12}}{\lambda_{q+1}^{-1}\tau_q^{-1}}
		=\tau_q\Lambda_q\leq1.
	\end{aligned}
\end{equation*}
	For the first inequality, lower the material threshold to $N_R-3J_*-1$ and apply
	\eqref{eq:Newton-one-threshold-weight-shift} four times.
	The equality uses $(\ell_q\lambda_q)^4=\lambda_q/\lambda_{q+1}$ from \eqref{def:parameter-mollification}.
	Taking $L=4$ in \eqref{eq:Nash-drift-mollification-difference} gives, for sufficiently large $a$ and $M+N\leq N_{\mathrm{ind}}-4$,
	\begin{equation}\label{eq:Nash-drift-mollification-transport-bound}
	\begin{aligned}
		\|\widetilde D_{t,q}^M(u_{q,J_*}-\bar u_{q,J_*})\|_N
		&\lesssim\ell_q^4\delta_q^{\frac12}\MM{M,N_t^{\mathrm N}+4,\tau_q^{-1},\ell_{t,q}^{-1}}
		\MM{N+4,N_{R,J_*}^{(1)}(M),\lambda_q,\ell_q^{-1}}\\
		&\leq\lambda_{q+1}^{-1}\tau_q^{-1}\MM{M,N_t^{\mathrm N}-1,\tau_q^{-1},\ell_{t,q}^{-1}}\ell_q^{-N}.
	\end{aligned}
	\end{equation}
	For $M+N\leq N_{\mathrm{ind}}-4$, take $L=4$ in
	\eqref{eq:Nash-density-mollification-difference} and use the threshold comparison above.
	For the stress difference, repeat the argument in the proof of
	Lemma~\ref{lem:initial-Newton-flow-mollification} with $L=4$ and $\widetilde D_{t,q}^M$, using
	\eqref{eq:Nash-localized-S-material-estimate} and \eqref{eq:Nash-flow-higher-order-CN-estimate}
	for the source and flow estimates, respectively.
	The required orders satisfy $M+4\leq N_{\mathrm{ind}}$ and
	$M+N+4\leq N_{\mathrm{ind}}\leq K_{J_*}^{(0)}+1$.
	Thus
	\begin{equation}\label{eq:Nash-separate-mollification-remainders}
	\begin{aligned}
	\|\widetilde D_{t,q}^M(\rho_{q,J_*}-\bar\rho_{q,J_*})\|_{N,\mcI_u^{\mathrm N}}
	&\lesssim(\ell_q\lambda_q)^4\delta_q
	\MM{M,N_t^{\mathrm N}-1,\tau_q^{-1},\ell_{t,q}^{-1}}\ell_q^{-N},\\
	e_q^{-1}\|\widetilde D_{t,q}^M(\bar S_{q,J_*,u}
	-S_{q,J_*,u})\|_{N,\mcI_u^{\mathrm N}}
	&\lesssim(\ell_{t,q}/\tau_q)^4
	\MM{M,N_t^{\mathrm N},\tau_q^{-1},\ell_{t,q}^{-1}}\ell_q^{-N}\\
	&\leq(\ell_q\lambda_q)^4
	\MM{M,N_t^{\mathrm N}-1,\tau_q^{-1},\ell_{t,q}^{-1}}\ell_q^{-N}.
	\end{aligned}
	\end{equation}
	For sufficiently large $a$, the last inequality follows from \eqref{eq:material-mollification-spatial-comparisons} and the definitions of the scales.

	The reciprocal factors are controlled by \eqref{eq:Newton-state-zero-order-bounds},
	\eqref{eq:Nash-mollified-terminal-zero-order-bounds},
	\eqref{eq:Nash-density-coefficient-material-estimate}, and
	\eqref{eq:Nash-mollified-density-material-estimate}, with $s=-1$ in the last two estimates.
	For $M+N\leq N_{\mathrm{ind}}-4$, apply \eqref{eq:holder-product-tame} with $\alpha=0$ to obtain
	\begin{equation}\label{eq:Nash-density-mollification-quotient-estimate}
	\begin{aligned}
	    \Nrm{\widetilde D_{t,q}^M\left(
	    \frac{\rho_{q,J_*}-\bar\rho_{q,J_*}}{\rho_{q,J_*}\bar\rho_{q,J_*}}\right)}_N
	    &\lesssim\sum_{\substack{M_0+M_1+M_2=M\\N_0+N_1+N_2=N}}
	    \Nrm{\widetilde D_{t,q}^{M_0}(\rho_{q,J_*}-\bar\rho_{q,J_*})}_{N_0}
	    \Nrm{\widetilde D_{t,q}^{M_1}(\rho_{q,J_*}^{-1})}_{N_1}
	    \Nrm{\widetilde D_{t,q}^{M_2}(\bar\rho_{q,J_*}^{-1})}_{N_2}\\
	    &\lesssim(\ell_q\lambda_q)^4\delta_q
	    \MM{M,N_t^{\mathrm N}-1,\tau_q^{-1},\ell_{t,q}^{-1}}\ell_q^{-N}.
	\end{aligned}
	\end{equation}
	In the same range, we obtain
	$$
	\begin{aligned}
	    &\left\|\widetilde D_{t,q}^M\left(
	    (u_{q,J_*}-\bar u_{q,J_*})\otimes w_{q+1}
	    +w_{q+1}\otimes(u_{q,J_*}-\bar u_{q,J_*})\right)\right\|_N\\
	    &\quad\lesssim\sum_{\substack{M_0+M_1=M\\N_0+N_1=N}}
	    \Nrm{\widetilde D_{t,q}^{M_0}(u_{q,J_*}-\bar u_{q,J_*})}_{N_0}
	    \nrm{\widetilde D_{t,q}^{M_1}w_{q+1}}_{N_1}\lesssim e_q^{\frac12}\tau_q^{-1}
	    \MM{M,N_t^{\mathrm N}-1,\tau_q^{-1},\ell_{t,q}^{-1}}\lambda_{q+1}^{N-1},\\
	    &\left\|\widetilde D_{t,q}^M\left(
	    \frac{\rho_{q,J_*}-\bar\rho_{q,J_*}}{\rho_{q,J_*}\bar\rho_{q,J_*}}
	    w_{q+1}^{\mathrm p}\otimes w_{q+1}^{\mathrm p}\right)\right\|_N\\
	    &\quad\lesssim\sum_{\substack{M_0+M_1+M_2=M\\N_0+N_1+N_2=N}}
	    \Nrm{\widetilde D_{t,q}^{M_0}\left(
	    \frac{\rho_{q,J_*}-\bar\rho_{q,J_*}}{\rho_{q,J_*}\bar\rho_{q,J_*}}\right)}_{N_0}
	    \nrm{\widetilde D_{t,q}^{M_1}w_{q+1}^{\mathrm p}}_{N_1}
	    \nrm{\widetilde D_{t,q}^{M_2}w_{q+1}^{\mathrm p}}_{N_2}
	    \lesssim(\ell_q\lambda_q)^4\delta_q e_q
	    \MM{M,N_t^{\mathrm N}-1,\tau_q^{-1},\ell_{t,q}^{-1}}\lambda_{q+1}^N,
	\end{aligned}
	$$
	where we used \eqref{eq:Newton-quadratic-weight-product}, \eqref{eq:Nash-principal-estimate},
	\eqref{eq:Nash-total-perturbation-estimate}, \eqref{eq:Nash-drift-mollification-transport-bound},
	and \eqref{eq:Nash-density-mollification-quotient-estimate}.
	Adding the stress difference in \eqref{eq:Nash-separate-mollification-remainders}, we conclude
	$$
	\begin{aligned}
		\|\widetilde D_{t,q}^MR_{q+1}^{\mathrm{mol}}\|_N
		&\lesssim\Bigl(e_q(\ell_q\lambda_q)^4\lambda_{q+1}^N
		+e_q^{\frac12}\tau_q^{-1}\lambda_{q+1}^{N-1}\Bigr)
		\MM{M,N_t^{\mathrm N}-1,\tau_q^{-1},\ell_{t,q}^{-1}}.
	\end{aligned}
	$$
	This proves \eqref{eq:Nash-flow-stress-zero-material-estimate}--\eqref{eq:Nash-flow-stress-positive-material-estimate}.

	The derivative conditions in Subsection~\ref{subsec:Newton-scheme-parameters} give
	$
	N_R+9\leq N_{\mathrm{ind}}-1.
	$
	These conditions ensure that all the error estimates above apply for $M+N\leq N_R$.
	For the oscillation and corrector terms, \eqref{eq:material-mollification-spatial-comparisons} gives, for sufficiently large $a$,
	\begin{equation}\label{eq:closure-material-spatial-compensation}
	\begin{aligned}
		&\frac{\MM{1,N_{R,J_*}^{(1)}(M),\lambda_q,\ell_q^{-1}}}{\lambda_q}
		\frac{\MM{M,N_t^{\mathrm N},\tau_q^{-1},\ell_{t,q}^{-1}}}
		{\MM{M,N_t^{\mathrm N}-1,\tau_q^{-1},\ell_{t,q}^{-1}}}\leq
		\begin{cases}
			1,&M\leq N_t^{\mathrm N}+2,\\
			\dfrac{\ell_{t,q}\ell_q^{-1}}{\tau_q\lambda_q},&M\geq N_t^{\mathrm N}+3
		\end{cases}
		\leq1.
	\end{aligned}
	\end{equation}
	Using $(\ell_q\lambda_q)^4=\lambda_q/\lambda_{q+1}$, we obtain
	$$
	\begin{aligned}
		\left\|\widetilde D_{t,q}^M\left(
		R_{q+1}^{\mathrm{trans}}+R_{q+1}^{\mathrm{osc}}
		+R_{q+1}^{\mathrm{corr}}+R_{q+1}^{\mathrm{mol}}\right)\right\|_N\lesssim\left(e_q^{\frac12}\tau_q^{-1}+e_q\lambda_q\right)
		\MM{M,N_t^{\mathrm N}-1,\tau_q^{-1},\ell_{t,q}^{-1}}\lambda_{q+1}^{N-1}.
	\end{aligned}
	$$
	Adding \eqref{eq:Nash-retained-Newton-stresses-material} in \eqref{def:Nash-final-Reynolds-stress}
	and using \eqref{def:Nash-stress-level} proves \eqref{eq:closure-Reynolds-background-material}.

	The factors $\theta_u^{\mathrm N}$ and Lemma~\ref{lem:Nash-transported-oscillatory-inverse-divergence}
	show that the four errors have temporal support in $\bigcup_u\supp\theta_u^{\mathrm N}$.
	Each selected interval $I_u$ has length $\tau_q/2$ and meets $\mathcal G_q$.
	By \eqref{def:Nash-transition-cutoff} and \eqref{eq:Newton-terminal-temporal-support}, we have
$$
\bigcup_u\supp\theta_u^{\mathrm N}
\subset\bigcup_u I_u
\subset\mN_{\tau_q/2}(\mathcal G_q)
\subset\mN_{(2J_*+2)\tau_q}(\supp_tR_q).
$$
	Combining this inclusion with \eqref{eq:Newton-terminal-temporal-support} in
	\eqref{def:Nash-final-Reynolds-stress} proves \eqref{eq:closure-stress-support-enlargement}.
\end{proof}

\subsection{Estimates at the next stage}
\label{subsec:closure-density-momentum}
We now estimate $\rho_{q+1}-\rho_q$, $m_{q+1}-m_q$, $\rho_{q+1}$, $m_{q+1}$, and $R_{q+1}$ with $D_{t,q+1}$.
The intermediate constants may depend on
$\beta,b,\alpha,\rho^*,\mathfrak M,C_*,N_R,N_{\mathrm{ind}},N_{\mathrm{max}}$,
but not on $q$ or $a$.
For fixed $C,\varepsilon,\eta>0$ and sufficiently large $a$, \eqref{def:lambda-delta} gives
\begin{equation}\label{eq:closure-fixed-constant-absorption}
	C\lambda_q^{-\varepsilon}\leq Ca^{-\varepsilon}\leq\eta,\qquad q\in\mathbb N_0.
\end{equation}
By \eqref{def:lambda-delta} and \eqref{def:admissible-alpha-losses}, we also have
\begin{equation}\label{eq:closure-state-smallness-factors}
\begin{aligned}
    \lambda_q^{-(b-1)(1-3\beta)+4b\alpha}
    \leq\lambda_q^{-20\alpha},\quad
    \lambda_q^{-(b-1)(1-2\beta)+4b\alpha}
    \leq\lambda_q^{-4\alpha},\quad
    \left(\frac{\lambda_q}{\lambda_{q+1}}\right)^{1-\beta}
    \lesssim\lambda_q^{-(b-1)(1-\beta)}
    \leq\lambda_q^{-16\alpha}.
\end{aligned}
\end{equation}

\subsubsection{Estimates with the new material derivative}
\label{subsec:closure-material-derivative-conversion}
\label{subsec:closure-four-Nash-stresses}

The following lemma transfers estimates from $\widetilde D_{t,q}$ to $D_{t,q+1}$.

\begin{lem}
\label{lem:Nash-material-derivative-change}
For $a>a_{\mathrm{Nash}}$ and $M,N\in\mathbb N_0$ in the indicated ranges, we have
\begin{equation}\label{eq:Nash-material-change-normalized-coefficients}
\begin{aligned}
	&\nrm{\widetilde D_{t,q}^{M}\left(u_{q,J_*}-\bar u_{q,J_*}+w_{q+1}/\rho_{q+1}\right)}_{N}
	\lesssim_{N_{\mathrm{max}}}\delta_{q+1}^{\frac12}\Lambda_{q+1}^{M}\lambda_{q+1}^{N},&& M+N\leq N_{\mathrm{ind}},\\
	&\nrm{\widetilde D_{t,q}^{M}\nabla\bar u_{q,J_*}}_{N}
	\lesssim_{N_{\mathrm{max}}}\Lambda_{q+1}^{M+1}\lambda_{q+1}^{N},&& M+N\leq N_{\mathrm{ind}}-1.
\end{aligned}
\end{equation}
Let $K\in\mathbb N$ with $1\leq K\leq N_{\mathrm{ind}}+1$, $A>0$, and $C\geq1$.
Suppose that a smooth scalar, vector, or tensor field $F$ satisfies
\begin{equation}\label{ass:Nash-material-derivative-change}
	\nrm{\widetilde D_{t,q}^MF}_N
	\leq CA\Lambda_{q+1}^M\lambda_{q+1}^N,
	\qquad 1\leq M+N\leq K.
\end{equation}
Then, 
\begin{equation}
	\nrm{D_{t,q+1}^MF}_N
	\lesssim CA\Lambda_{q+1}^M\lambda_{q+1}^N,\qquad1\leq M+N\leq K
	\label{eq:closure-material-derivative-comparison}
\end{equation}
where the implicit constant depends only on $N_{\mathrm{max}}$.
\end{lem}

\begin{proof}
We first prove \eqref{eq:Nash-material-change-normalized-coefficients}.
For $1\leq M+N\leq N_{\mathrm{ind}}$, \eqref{eq:Nash-mollified-terminal-velocity-integer-estimate} and \eqref{eq:Nash-terminal-velocity-material-estimate} give
\begin{equation}\label{eq:closure-background-velocity-at-next-scale}
\begin{aligned}
	\|\widetilde D_{t,q}^M\bar u_{q,J_*}\|_N+\|\widetilde D_{t,q}^Mu_{q,J_*}\|_N
	&\lesssim\delta_q^{\frac12}
	\MM{M,N_R-3J_*,\tau_q^{-1},\ell_{t,q}^{-1}}
	\MM{N,N_{R,J_*}^{(0)}(M),\lambda_q,\ell_q^{-1}}\\
	&\lesssim\lambda_q^{-(b-1)(1-2\beta)+4b\alpha}
	\delta_{q+1}^{\frac12}\Lambda_{q+1}^M\lambda_{q+1}^N\\
	&\leq\lambda_q^{-4\alpha}\delta_{q+1}^{\frac12}\Lambda_{q+1}^M\lambda_{q+1}^N.
\end{aligned}
\end{equation}
Here we used \eqref{def:lambda-delta}, \eqref{def:parameter-mollification}, \eqref{def:tau-q}, and \eqref{eq:material-mollification-hierarchy}.
Estimates \eqref{eq:closure-fixed-constant-absorption} and \eqref{eq:closure-state-smallness-factors} make the coefficient in \eqref{eq:closure-background-velocity-at-next-scale} as small as needed.

By \eqref{def:tau-q} and the definition of $A_{q,0}$, we have
\begin{align}
	A_{q,0}^{\frac12}(\tau_q\ell_q^{-1})^{\frac{\alpha}{2}}
	\lesssim\lambda_q^{-4\alpha
	+\frac\alpha2\left(\frac{b-1}{4}+\beta-4b\alpha\right)}\delta_{q+1}^{\frac12}\leq\lambda_q^{-\frac{15\alpha}{4}}\delta_{q+1}^{\frac12}.
	\label{eq:closure-Nash-increment-ratio}
\end{align}
For the second inequality, we used $b-1<2\beta$ and $\beta<1/3$.
For $w_{q+1}$, \eqref{eq:material-mollification-hierarchy},
\eqref{eq:Nash-total-perturbation-estimate}, and \eqref{eq:closure-Nash-increment-ratio} give
\begin{equation}\label{eq:closure-Nash-perturbation-at-new-scale}
\begin{aligned}
	\|\widetilde D_{t,q}^Mw_{q+1}\|_N
	&\lesssim A_{q,0}^{\frac12}(\tau_q\ell_q^{-1})^{\frac\alpha2}
	\MM{M,N_t^{\mathrm N},\tau_q^{-1},\ell_{t,q}^{-1}}\lambda_{q+1}^N\\
	&\lesssim A_{q,0}^{\frac12}(\tau_q\ell_q^{-1})^{\frac\alpha2}
	\Lambda_{q+1}^M\lambda_{q+1}^N,
	&& M+N\leq N_{\mathrm{ind}}.
\end{aligned}
\end{equation}
Using the zero-order estimates in \eqref{eq:Newton-state-zero-order-bounds}
and \eqref{eq:Nash-density-coefficient-material-estimate} with $s=-1$,
we obtain, for $M+N\leq N_{\mathrm{ind}}$,
\begin{equation}\label{eq:closure-Nash-density-weighted-perturbation}
	\begin{aligned}
	\Nrm{\widetilde D_{t,q}^M(w_{q+1}\rho_{q+1}^{-1})}_N
	&\lesssim\sum_{\substack{M_1+M_2=M\\N_1+N_2=N}}
	\|\widetilde D_{t,q}^{M_1}(\rho_{q+1}^{-1})\|_{N_1}
	\|\widetilde D_{t,q}^{M_2}w_{q+1}\|_{N_2}\\
	&\lesssim A_{q,0}^{\frac12}(\tau_q\ell_q^{-1})^{\frac\alpha2}
	\MM{M,N_t^{\mathrm N},\tau_q^{-1},\ell_{t,q}^{-1}}\lambda_{q+1}^N,
	\end{aligned}
\end{equation}
where we used \eqref{eq:Newton-quadratic-weight-product} with $K=N_t^{\mathrm N}$.
At $M=N=0$, use \eqref{eq:Nash-drift-mollification-difference} with $L=4$ and \eqref{eq:closure-Nash-density-weighted-perturbation} to obtain
\begin{equation}\label{eq:closure-zero-order-velocity-difference}
\begin{aligned}
	\|u_{q,J_*}-\bar u_{q,J_*}+w_{q+1}/\rho_{q+1}\|_0
	&\lesssim(\ell_q\lambda_q)^4\delta_q^{\frac12}
	+A_{q,0}^{\frac12}(\tau_q\ell_q^{-1})^{\frac\alpha2}\leq\lambda_q^{-3\alpha}\delta_{q+1}^{\frac12}.
\end{aligned}
\end{equation}
Here $(\ell_q\lambda_q)^4\delta_q^{\frac12}=(\lambda_q/\lambda_{q+1})^{1-\beta}\delta_{q+1}^{\frac12}$ by \eqref{def:lambda-delta} and \eqref{def:parameter-mollification}. The last inequality follows from \eqref{eq:closure-fixed-constant-absorption}, \eqref{eq:closure-state-smallness-factors}, and \eqref{eq:closure-Nash-increment-ratio}, after increasing $a_{\mathrm{Nash}}$.
For $1\leq M+N\leq N_{\mathrm{ind}}$, use \eqref{eq:closure-fixed-constant-absorption} and \eqref{eq:closure-Nash-increment-ratio} in \eqref{eq:closure-background-velocity-at-next-scale} and \eqref{eq:closure-Nash-density-weighted-perturbation}. Together with \eqref{eq:closure-zero-order-velocity-difference}, we obtain the first inequality in \eqref{eq:Nash-material-change-normalized-coefficients} after increasing $a_{\mathrm{Nash}}$.

For $M+N\leq N_{\mathrm{ind}}-1$, \eqref{eq:Nash-terminal-Newton-velocity-sharp-gradient} gives
$$
\begin{aligned}
	\|\widetilde D_{t,q}^M\nabla\bar u_{q,J_*}\|_N
	\lesssim\Lambda_q
	\MM{M,N_R-3J_*-1,\tau_q^{-1},\ell_{t,q}^{-1}}
	\MM{N,N_{R,J_*}^{(1)}(M),\lambda_q,\ell_q^{-1}}
	\lesssim\left(\frac{\lambda_q}{\lambda_{q+1}}\right)^{1-\beta}
	\Lambda_{q+1}^{M+1}\lambda_{q+1}^N.
\end{aligned}
$$
Here the second inequality uses \eqref{def:lambda-delta}, \eqref{def:parameter-mollification}, and \eqref{eq:material-mollification-hierarchy}.
By \eqref{eq:closure-fixed-constant-absorption}, the second inequality in \eqref{eq:Nash-material-change-normalized-coefficients} holds after increasing $a_{\mathrm{Nash}}$.

It remains to prove \eqref{eq:closure-material-derivative-comparison}.
The material derivatives satisfy
$$
D_{t,q+1}=\widetilde D_{t,q}
+\left(u_{q,J_*}-\bar u_{q,J_*}+\frac{w_{q+1}}{\rho_{q+1}}\right)\cdot\nabla.
$$
Choose $C_0\geq1$ to cover both implicit constants in \eqref{eq:Nash-material-change-normalized-coefficients}.
Apply Lemma~\ref{lem:source-material-derivative-change}
with $\upsilon=\bar u_{q,J_*}$, $\omega=u_{q,J_*}-\bar u_{q,J_*}+w_{q+1}/\rho_{q+1}$, and
$$
	J=K,\quad N_t=0,\quad \widetilde{\mC}_0=C_0,\quad
	\mC_\omega=\delta_{q+1}^{\frac12},\quad \mC_F=CA,\quad
	\lambda=\tilde\lambda=\lambda_{q+1},\quad \mu=\tilde\mu=\Lambda_{q+1}.
$$
The scale conditions follow from $\delta_{q+1}^{\frac12}\lambda_{q+1}=\Lambda_{q+1}$.
Estimate~\eqref{eq:Nash-material-change-normalized-coefficients} verifies the coefficient hypotheses,
which require total orders at most $K-1\leq N_{\mathrm{ind}}$ for $\omega$
and $K-2\leq N_{\mathrm{ind}}-1$ for $\nabla\bar u_{q,J_*}$.
The gradient hypothesis is void when $K=1$.
The source hypothesis is \eqref{ass:Nash-material-derivative-change}.
The constant in Lemma~\ref{lem:source-material-derivative-change} depends only on $K$ and $C_0$.
Since $K\leq N_{\mathrm{max}}$ and $C_0$ depends only on $N_{\mathrm{max}}$,
this proves \eqref{eq:closure-material-derivative-comparison} with the stated dependence.
\end{proof}

\begin{prop}
\label{prop:closure-state-material}
For $a>a_{\mathrm{Nash}}$ and $M,N\in\mathbb N_0$, the density and momentum perturbations satisfy
\begin{equation}\label{eq:closure-final-density-perturbation-material}
\begin{aligned}
    \nrm{D_{t,q+1}^M(\rho_{q+1}-\rho_q)}_N
    &\leq
    \frac12\lambda_q^{-6\alpha}\delta_{q+1}\Lambda_{q+1}^M\lambda_{q+1}^N,
    && M+N\leq N_{\mathrm{ind}}.
\end{aligned}
\end{equation}
\begin{equation}\label{eq:closure-final-perturbation-material}
\begin{aligned}
	\nrm{D_{t,q+1}^M(m_{q+1}-m_q)}_N
	&\lesssim_{\beta,\alpha,N_{\mathrm{max}}}
	\lambda_q^{-3\alpha}\delta_{q+1}^{\frac12}\Lambda_{q+1}^M\lambda_{q+1}^N,
	&& M+N\leq N_{\mathrm{ind}}.
\end{aligned}
\end{equation}
The updated density and momentum satisfy
\begin{align}
	\nrm{D_{t,q+1}^M\rho_{q+1}}_N
	&\leq
	\delta_{q+1}
	\lambda_{q+1}^N\Lambda_{q+1}^M,
	&&1\leq M+N\leq N_{\mathrm{ind}},
	\label{eq:closure-final-density-material}\\
	\nrm{D_{t,q+1}^Mm_{q+1}}_N
	&\lesssim_{\beta,\alpha,N_{\mathrm{max}}}
	\lambda_q^{-3\alpha}\delta_{q+1}^{\frac12}
	\lambda_{q+1}^N\Lambda_{q+1}^M,
	&&1\leq M+N\leq N_{\mathrm{ind}}.
	\label{eq:closure-final-momentum-material}
\end{align}
The Reynolds stress satisfies
\begin{equation}
	\nrm{D_{t,q+1}^MR_{q+1}}_N
	\leq
	\delta_{q+2}
	\lambda_{q+1}^{N-8\alpha}\Lambda_{q+1}^M,
	\qquad 1\leq M+N\leq N_R.
	\label{eq:closure-final-Reynolds-material}
\end{equation}
At order zero, we have
\begin{equation}
	\nrm{R_{q+1}}_0
	\leq
	\lambda_{q+1}^{-8\alpha}\delta_{q+2}.
	\label{eq:closure-Reynolds-zero}
\end{equation}
\end{prop}

\begin{proof}
For the perturbations and the updated state, we take $K=N_{\mathrm{ind}}$ and $C=1$
in Lemma~\ref{lem:Nash-material-derivative-change}.
The implicit constant in \eqref{eq:closure-material-derivative-comparison} depends only on $N_{\mathrm{max}}$.
For $M+N\leq N_{\mathrm{ind}}$, \eqref{def:parameter-mollification}, \eqref{eq:material-mollification-hierarchy},
\eqref{eq:Nash-terminal-Newton-momentum-increment-material}, \eqref{def:Nash-final-state}, and
\eqref{eq:closure-Nash-perturbation-at-new-scale} give
\begin{equation}\label{eq:closure-momentum-increment-background-material}
\begin{aligned}
	\nrm{\widetilde D_{t,q}^M(m_{q+1}-m_q)}_N
	&\leq\nrm{\widetilde D_{t,q}^M(m_{q,J_*}-m_q)}_N
	+\nrm{\widetilde D_{t,q}^Mw_{q+1}}_N\\
	&\lesssim_{\beta,b,\alpha,N_{\mathrm{max}},\rho^*,\mathfrak M,C_*}
	\left[A_{q,0}\delta_q^{-\frac12}\lambda_{q+1}^{-3\alpha}
	+A_{q,0}^{\frac12}(\tau_q\ell_q^{-1})^{\frac\alpha2}\right]
	\MM{M,N_t^{\mathrm N},\tau_q^{-1},\ell_{t,q}^{-1}}\lambda_{q+1}^N\\
	&\lesssim_{\beta,b,\alpha,N_{\mathrm{max}},\rho^*,\mathfrak M,C_*}
	\left(\lambda_q^{-11\alpha}+\lambda_q^{-\frac{15\alpha}{4}}\right)
	\delta_{q+1}^{\frac12}\MM{M,N_t^{\mathrm N},\tau_q^{-1},\ell_{t,q}^{-1}}\lambda_{q+1}^N\\
	&\leq\lambda_q^{-3\alpha}\delta_{q+1}^{\frac12}\Lambda_{q+1}^M\lambda_{q+1}^N.
\end{aligned}
\end{equation}
The third and fourth inequalities use \eqref{eq:Newton-terminal-velocity-difference-amplitudes},
\eqref{eq:closure-fixed-constant-absorption}, and \eqref{eq:closure-Nash-increment-ratio}, after increasing $a_{\mathrm{Nash}}$.
For $F=m_{q+1}-m_q$, take $A=\lambda_q^{-3\alpha}\delta_{q+1}^{\frac12}$ in Lemma~\ref{lem:Nash-material-derivative-change}.
This proves \eqref{eq:closure-final-perturbation-material} at positive total orders.
At $M=N=0$, use \eqref{eq:closure-momentum-increment-background-material} directly.

For $M+N\leq N_{\mathrm{ind}}$, \eqref{eq:Nash-terminal-Newton-density-increment-material},
\eqref{def:parameter-mollification}, \eqref{eq:material-mollification-hierarchy},
and $A_{q,0}=C_*\lambda_q^{-8\alpha}\delta_{q+1}$ give
\begin{equation}\label{eq:closure-density-increment-mollified-material}
\begin{aligned}
    \nrm{\widetilde D_{t,q}^M(\rho_{q+1}-\rho_q)}_N
    &\lesssim_{\beta,b,\alpha,N_{\mathrm{max}},\rho^*,\mathfrak M,C_*}
    A_{q,0}\ell_q^{-\alpha/2}
    \MM{M,N_R-3J_*,\tau_q^{-1},\ell_{t,q}^{-1}}
    \MM{N,N_{R,J_*}^{(0)}(M),\lambda_q,\ell_q^{-1}}\\
    &\leq\lambda_q^{-7\alpha}\delta_{q+1}\Lambda_{q+1}^M\lambda_{q+1}^N.
\end{aligned}
\end{equation}
Here $A_{q,0}\ell_q^{-\alpha/2}\lesssim_{b,\alpha,C_*}
\lambda_q^{-8\alpha+(3+b)\alpha/8}\delta_{q+1}$ and $b<1+2\beta<5/3$.
The last inequality follows from \eqref{eq:closure-fixed-constant-absorption}, after increasing $a_{\mathrm{Nash}}$.
Apply Lemma~\ref{lem:Nash-material-derivative-change} with $F=\rho_{q+1}-\rho_q$ and $A=\lambda_q^{-7\alpha}\delta_{q+1}$.
Absorbing the factor $\lambda_q^{-\alpha}$ by \eqref{eq:closure-fixed-constant-absorption}
proves \eqref{eq:closure-final-density-perturbation-material} with coefficient $\tfrac12$ at positive total orders, after enlarging $a_{\mathrm{Nash}}$.
At $M=N=0$, use \eqref{eq:closure-density-increment-mollified-material} and
$\lambda_q^{-\alpha}\leq\tfrac12$, also ensured by this choice of $a_{\mathrm{Nash}}$.

For the updated state, use $\rho_{q+1}=\rho_{q,J_*}$ and $m_{q+1}=m_{q,J_*}+w_{q+1}$ in \eqref{def:Nash-final-state}.
For $1\leq M+N\leq N_{\mathrm{ind}}$, \eqref{def:parameter-mollification}, \eqref{def:tau-q},
\eqref{eq:material-mollification-hierarchy}, and $N_R-3J_*\geq1$ give
\begin{equation}\label{eq:closure-positive-order-weight-comparison}
\begin{aligned}
	\MM{M,N_R-3J_*,\tau_q^{-1},\ell_{t,q}^{-1}}
	\MM{N,N_{R,J_*}^{(0)}(M),\lambda_q,\ell_q^{-1}}\leq\tau_q^{-1}\Lambda_{q+1}^{M-1}\lambda_{q+1}^N.
\end{aligned}
\end{equation}
For $M=0<N$, \eqref{eq:closure-positive-order-weight-comparison} follows from
$\MM{N,N_R-3J_*,\lambda_q,\ell_q^{-1}}\leq\lambda_q\lambda_{q+1}^{N-1}$
and $\lambda_q\lambda_{q+1}^{-1}\leq\tau_q^{-1}\Lambda_{q+1}^{-1}$.
Estimates
\eqref{eq:Nash-density-coefficient-material-estimate} with $s=1$,
\eqref{eq:Nash-terminal-momentum-material-estimate}, and
\eqref{eq:closure-Nash-perturbation-at-new-scale} give
\begin{align}
	\nrm{\widetilde D_{t,q}^M\rho_{q+1}}_N
	&\lesssim_{\beta,b,\alpha,N_{\mathrm{max}},\rho^*,\mathfrak M,C_*}\delta_q
	\MM{M,N_R-3J_*,\tau_q^{-1},\ell_{t,q}^{-1}}
	\MM{N,N_{R,J_*}^{(0)}(M),\lambda_q,\ell_q^{-1}}\notag\\
	&\lesssim_{\beta,b,\alpha,N_{\mathrm{max}},\rho^*,\mathfrak M,C_*}\delta_q\tau_q^{-1}\Lambda_{q+1}^{M-1}\lambda_{q+1}^N\notag\\
	&\leq\lambda_q^{-20\alpha}\delta_{q+1}
	\Lambda_{q+1}^M\lambda_{q+1}^N,
	\label{eq:closure-density-input-at-new-scale}\\
	\nrm{\widetilde D_{t,q}^Mm_{q+1}}_N
	&\leq\nrm{\widetilde D_{t,q}^Mm_{q,J_*}}_N
	+\nrm{\widetilde D_{t,q}^Mw_{q+1}}_N\notag\\*
	&\lesssim_{\beta,b,\alpha,N_{\mathrm{max}},\rho^*,\mathfrak M,C_*}\left[\tau_q^{-1}\Lambda_{q+1}^{-1}\delta_q^{\frac12}
	+A_{q,0}^{\frac12}(\tau_q\ell_q^{-1})^{\frac\alpha2}\right]
	\Lambda_{q+1}^M\lambda_{q+1}^N\notag\\
	&\lesssim_{\beta,b,\alpha,N_{\mathrm{max}},\rho^*,\mathfrak M,C_*}\left(\lambda_q^{-4\alpha}+\lambda_q^{-\frac{15\alpha}{4}}\right)
	\delta_{q+1}^{\frac12}\Lambda_{q+1}^M\lambda_{q+1}^N\notag\\
	&\leq\lambda_q^{-3\alpha}\delta_{q+1}^{\frac12}\Lambda_{q+1}^M\lambda_{q+1}^N.
	\label{eq:closure-momentum-input-at-new-scale}
\end{align}
Here we used \eqref{def:lambda-delta}, \eqref{def:tau-q},
\eqref{eq:closure-state-smallness-factors}, \eqref{eq:closure-Nash-increment-ratio},
and \eqref{eq:closure-positive-order-weight-comparison}.
The last inequalities follow from \eqref{eq:closure-fixed-constant-absorption} after enlarging $a_{\mathrm{Nash}}$.
For the density estimate, \eqref{def:admissible-alpha-losses} gives
$(b-1)(1-3\beta)-4b\alpha>20\alpha$.

For $F=\rho_{q+1}$, apply Lemma~\ref{lem:Nash-material-derivative-change} with $A=\lambda_q^{-20\alpha}\delta_{q+1}$,
using \eqref{eq:closure-density-input-at-new-scale}.
Absorbing the factor $\lambda_q^{-20\alpha}$ by \eqref{eq:closure-fixed-constant-absorption}
proves \eqref{eq:closure-final-density-material}, after increasing $a_{\mathrm{Nash}}$.
For $F=m_{q+1}$, take $A=\lambda_q^{-3\alpha}\delta_{q+1}^{\frac12}$ in the same lemma.
Estimate~\eqref{eq:closure-momentum-input-at-new-scale} gives \eqref{eq:closure-final-momentum-material}.
Its implicit constant comes from \eqref{eq:closure-material-derivative-comparison};
the other constants have been absorbed into $a_{\mathrm{Nash}}$.

We next estimate the coefficients in \eqref{eq:closure-Reynolds-background-material}.
Definitions~\eqref{def:lambda-delta}, \eqref{def:parameter-mollification}, \eqref{def:tau-q},
and the definition of $A_{q,0}$ give
	\begin{align}
		\frac{A_{q,0}^{\frac12}
			(\tau_q\ell_q^{-1})^{\frac\alpha2}
			\tau_q^{-1}\lambda_{q+1}^{-1}}
		{\lambda_{q+1}^{-8\alpha}\delta_{q+2}}
		&\lesssim_{\beta,b,\alpha,C_*}
		\lambda_q^{-(b-1)(1-\beta(2b+1))+(12b-4)\alpha+\frac\alpha2\left(\frac{b-1}{4}+\beta-4b\alpha\right)}\leq\lambda_q^{-15\alpha},
		\label{eq:closure-transport-stress-ratio}\\
		\frac{A_{q,0}(\tau_q\ell_q^{-1})^\alpha
			\lambda_q\lambda_{q+1}^{-1}}
		{\lambda_{q+1}^{-8\alpha}\delta_{q+2}}
		&\lesssim_{\beta,b,\alpha,C_*}
		\lambda_q^{-(b-1)(1-2\beta b-8\alpha)+\alpha\left(\frac{b-1}{4}+\beta-4b\alpha\right)}\leq\lambda_q^{-15\alpha}.
		\label{eq:closure-oscillation-corrector-ratio}
	\end{align}
Here \eqref{def:admissible-alpha-losses} and $b-1<2\beta$ give
$0<\frac{b-1}{4}+\beta-4b\alpha<2\beta<\frac23$.
Using $$1-\beta(2b+1)=3\left(\frac13-\beta(2b-1)\right)+4\beta(b-1),$$
the two exponents are less than $-(12b+\frac{11}{3})\alpha$ and $-(8b+\frac{22}{3})\alpha$, respectively.
Both are less than $-15\alpha$ since $b>1$.
For $M+N\leq N_R$, \eqref{def:parameter-mollification}, \eqref{eq:material-mollification-hierarchy},
\eqref{eq:closure-transport-stress-ratio}, and \eqref{eq:closure-oscillation-corrector-ratio} give
\begin{equation}\label{eq:closure-Reynolds-normalized-input}
\begin{aligned}
	\|\widetilde D_{t,q}^MR_{q+1}\|_N
	&\lesssim_{\beta,b,\alpha,N_{\mathrm{max}},\rho^*,\mathfrak M,C_*}\lambda_{q+1}^{-8\alpha}\delta_{q+2}
	\Bigl[\lambda_q^{-\alpha}\ell_{t,q}^{-M}\ell_q^{-N}
	+\lambda_q^{-15\alpha}
	\MM{M,N_t^{\mathrm N}-1,\tau_q^{-1},\ell_{t,q}^{-1}}\lambda_{q+1}^N\Bigr]\\
	&\lesssim_{\beta,b,\alpha,N_{\mathrm{max}},\rho^*,\mathfrak M,C_*}\bigl(\lambda_q^{-\alpha}+\lambda_q^{-15\alpha}\bigr)
	\delta_{q+2}\lambda_{q+1}^{N-8\alpha}\Lambda_{q+1}^M.
\end{aligned}
\end{equation}
For $1\leq M+N\leq N_R$, apply Lemma~\ref{lem:Nash-material-derivative-change} to
\eqref{eq:closure-Reynolds-normalized-input} with $F=R_{q+1}$ and $K=N_R$ to obtain
$$
\begin{aligned}
	\|D_{t,q+1}^MR_{q+1}\|_N
	\lesssim_{\beta,b,\alpha,N_{\mathrm{max}},\rho^*,\mathfrak M,C_*}\bigl(\lambda_q^{-\alpha}+\lambda_q^{-15\alpha}\bigr)
	\delta_{q+2}\lambda_{q+1}^{N-8\alpha}\Lambda_{q+1}^M
	\leq\delta_{q+2}\lambda_{q+1}^{N-8\alpha}\Lambda_{q+1}^M.
\end{aligned}
$$
The last inequality follows from \eqref{eq:closure-fixed-constant-absorption} after increasing $a_{\mathrm{Nash}}$.
This proves \eqref{eq:closure-final-Reynolds-material}.
At $M=N=0$, the same absorption in \eqref{eq:closure-Reynolds-normalized-input} gives \eqref{eq:closure-Reynolds-zero}.

\end{proof}
\section{Proof of the main iterative proposition}
\label{sec:main-iteration-proof}
\begin{proof}[Proof of Proposition~\ref{prop:main-iteration}]

	Let $(\rho_q,m_q,R_q)$ satisfy the hypotheses of Proposition~\ref{prop:main-iteration}.
	The triple $(\rho_{q+1},m_{q+1},R_{q+1})$ defined by
	\eqref{def:Nash-final-state} and \eqref{def:Nash-final-Reynolds-stress} is a smooth solution of
	\eqref{eq:compressible-Euler-Reynolds} by \eqref{eq:Nash-output-system}.
	Let $C\geq1$ be a common constant for the estimates
	\eqref{eq:closure-final-density-perturbation-material}--\eqref{eq:closure-final-momentum-material}
	over their stated ranges.
	By the proof of Proposition~\ref{prop:closure-state-material}, $C$ depends only on
	$\beta,\alpha,N_{\mathrm{max}}$.
	Fix $\mathfrak M:=2C$, choose the construction constants, and then take
	$a_0\geq a_{\mathrm{Nash}}$ large enough for the estimates below.

	The Reynolds stress $R_{q+1}$ satisfies \eqref{ind:R-zero}--\eqref{ind:R-material}
	by Proposition~\ref{prop:closure-state-material}. We first verify the inductive estimates \eqref{ind:rho-zero}--\eqref{ind:m-material} and the density condition \eqref{ind:rho-gap} at stage $q+1$.
	We then prove the estimates \eqref{inc:rho}--\eqref{inc:m} and check the support conditions \eqref{ind:R-support} and \eqref{inc:support}.

Estimate~\eqref{eq:closure-final-perturbation-material} gives
\begin{equation}\label{eq:closure-momentum-increment-with-iteration-coefficient}
\begin{aligned}
	\nrm{D_{t,q+1}^M(m_{q+1}-m_q)}_N
	&\leq\frac{\mathfrak M}{2}\lambda_q^{-3\alpha}\delta_{q+1}^{\frac12}
	\Lambda_{q+1}^M\lambda_{q+1}^N,
	&& M+N\leq N_{\mathrm{ind}}.
\end{aligned}
\end{equation}
Using \eqref{ind:m-zero} and \eqref{eq:closure-momentum-increment-with-iteration-coefficient} at $M=N=0$, we obtain
$$
\begin{aligned}
	\nrm{m_{q+1}}_0
	\leq\nrm{m_q}_0+\nrm{m_{q+1}-m_q}_0
	\leq\mathfrak M(1-\delta_q^{\frac12})+\frac{\mathfrak M}{2}\delta_{q+1}^{\frac12}
	\leq\mathfrak M(1-\delta_{q+1}^{\frac12}).
\end{aligned}
$$
The last inequality follows from \eqref{def:lambda-delta} for sufficiently large $a$.
Since $\rho_{q+1}=\rho_{q,J_*}$, \eqref{eq:Newton-state-zero-order-bounds} gives
\eqref{ind:rho-zero} and \eqref{ind:rho-gap}.
The density estimate \eqref{eq:closure-final-density-material} gives \eqref{ind:rho-material}.
Estimate~\eqref{eq:closure-final-momentum-material}, with $C=\mathfrak M/2$ and
$\lambda_q^{-3\alpha}\leq1$, proves \eqref{ind:m-material} at stage $q+1$.

We now verify \eqref{inc:rho}--\eqref{inc:m}.
For the time derivatives, write $\partial_t=D_{t,q+1}-\rho_{q+1}^{-1}m_{q+1}\cdot\nabla$.
Using \eqref{ind:m-zero} and \eqref{ind:rho-gap} at stage $q+1$,
we obtain from \eqref{eq:closure-final-density-perturbation-material} and
\eqref{eq:closure-momentum-increment-with-iteration-coefficient}
\begin{equation}\label{eq:closure-increment-time-estimates}
\begin{aligned}
    \nrm{\partial_t(\rho_{q+1}-\rho_q)}_0
    &\leq\nrm{D_{t,q+1}(\rho_{q+1}-\rho_q)}_0
    +3\nrm{m_{q+1}}_0\nrm{\rho_{q+1}^{-1}}_0\nrm{\rho_{q+1}-\rho_q}_1\\
    &\leq\frac12\left(1+\frac{6\mathfrak M}{\rho^*}\right)
    \lambda_q^{-6\alpha}\delta_{q+1}\lambda_{q+1},\\
    \nrm{\partial_t(m_{q+1}-m_q)}_0
    &\leq\nrm{D_{t,q+1}(m_{q+1}-m_q)}_0
    +3\nrm{m_{q+1}}_0\nrm{\rho_{q+1}^{-1}}_0\nrm{m_{q+1}-m_q}_1\\
    &\leq\frac{\mathfrak M}{2}\left(1+\frac{6\mathfrak M}{\rho^*}\right)
    \lambda_q^{-3\alpha}\delta_{q+1}^{\frac12}\lambda_{q+1}.
\end{aligned}
\end{equation}
Here $\Lambda_{q+1}\leq\lambda_{q+1}$ by \eqref{def:lambda-delta}.
Combining \eqref{eq:closure-final-density-perturbation-material},
\eqref{eq:closure-momentum-increment-with-iteration-coefficient}, and
\eqref{eq:closure-increment-time-estimates}, we obtain
\begin{equation}\label{eq:closure-normalized-increments}
\begin{aligned}
	&\nrm{\rho_{q+1}-\rho_q}_0+\lambda_{q+1}^{-1}\nrm{\rho_{q+1}-\rho_q}_1
	+\lambda_{q+1}^{-1}\nrm{\partial_t(\rho_{q+1}-\rho_q)}_0
	\leq\frac12\left(3+\frac{6\mathfrak M}{\rho^*}\right)\lambda_q^{-6\alpha}\delta_{q+1}
	\leq\mathfrak M\delta_{q+1},\\
	&\nrm{m_{q+1}-m_q}_0+\lambda_{q+1}^{-1}\nrm{m_{q+1}-m_q}_1
	+\lambda_{q+1}^{-1}\nrm{\partial_t(m_{q+1}-m_q)}_0
	\leq\frac{\mathfrak M}{2}\left(3+\frac{6\mathfrak M}{\rho^*}\right)
	\lambda_q^{-3\alpha}\delta_{q+1}^{\frac12}
	\leq\mathfrak M\delta_{q+1}^{\frac12}.
\end{aligned}
\end{equation}
The last inequalities follow from \eqref{eq:closure-fixed-constant-absorption} for sufficiently large $a$.
Equation~\eqref{eq:closure-normalized-increments} proves \eqref{inc:rho}--\eqref{inc:m}.

It remains to check \eqref{ind:R-support} and \eqref{inc:support}.
For sufficiently large $a$, definitions~\eqref{def:lambda-delta} and \eqref{def:tau-q} give
$$
\begin{aligned}
(2J_*+2)\tau_q
\leq(2J_*+2)\lambda_{q+1}^{-3\alpha}\Lambda_q^{-1}
<(1-\Lambda_q\Lambda_{q+1}^{-1})\Lambda_q^{-1}
=\Lambda_q^{-1}-\Lambda_{q+1}^{-1}.
\end{aligned}
$$
By \eqref{ind:R-support} and \eqref{eq:closure-stress-support-enlargement}, we obtain
$$
\begin{aligned}
	\supp_tR_{q+1}
	&\subset\mN_{(2J_*+2)\tau_q}(\supp_tR_q)\\
	&\subset
	\bigl[-2+\Lambda_{q+1}^{-1},-1-\Lambda_{q+1}^{-1}\bigr]
	\cup
	\bigl[1+\Lambda_{q+1}^{-1},2-\Lambda_{q+1}^{-1}\bigr].
\end{aligned}
$$
For the perturbations, \eqref{eq:Newton-terminal-temporal-support}, \eqref{def:Nash-final-state},
\eqref{eq:Nash-perturbation-temporal-support}, and \eqref{eq:closure-stress-support-enlargement} give
\begin{equation}
	\begin{aligned}
		\supp_t(\rho_{q+1}-\rho_q,m_{q+1}-m_q)
		&\subset
		\supp_t(\rho_{q,J_*}-\rho_q,m_{q,J_*}-m_q)
		\cup\supp_t w_{q+1}\\
		&\subset
		\mN_{(2J_*+2)\tau_q}(\supp_tR_q)\\
		&\subset
		(-2,-1)\cup(1,2).
	\end{aligned}
\end{equation}
These inclusions prove \eqref{ind:R-support} at stage $q+1$ and \eqref{inc:support}.
\end{proof}

\appendix
\addtocontents{toc}{\protect\setcounter{tocdepth}{1}}

	\section{H{\"o}lder and Besov spaces}
	\label{app:holder-besov-spaces}

	\subsection{H{\"o}lder spaces}
	\label{sec:holder-spaces}

	We define the H{\"o}lder norms and record the interpolation and
	product estimates used in the paper. For the Littlewood--Paley
	characterization and related estimates, see \cite{BCD11}.
	
	Let $d\in\mathbb N$ and $N\in\mathbb N_0$. For a function
	$f:\mathbb T^d\to\mathbb R$, we define
	$$
	\|f\|_0:=\sup_{x\in\mathbb T^d}|f(x)|.
	$$
	For a multi-index $\theta$, we write $\partial^\theta f$ for the
	corresponding spatial derivative. We set
	$
	[f]_0:=\|f\|_0,
	$
	and, for $N\geq1$,
	$$
	[f]_N
	:=
	\sup_{|\theta|=N}\|\partial^\theta f\|_0.
	$$
	For $0<\alpha<1$, we define the H{\"o}lder seminorm
	$$
	[f]_{N+\alpha}
	:=
	\sup_{|\theta|=N}
	\sup_{x\neq y}
	\frac{
		|\partial^\theta f(x)-\partial^\theta f(y)|
	}{
		|x-y|^\alpha
	}.
	$$
	Here $|x-y|$ denotes the distance on $\mathbb T^d$.
	The corresponding H{\"o}lder norms are
	$$
	\|f\|_N
	:=
	\sum_{j=0}^N [f]_j,
	\qquad
	\|f\|_{N+\alpha}
	:=
	\|f\|_N+[f]_{N+\alpha}.
	$$
	We use the same notation for vector and tensor fields, taking norms
componentwise. For flow matrices, we define the $C^0$ norm using the equivalent
Euclidean operator norm.
For functions on $\mathbb R^d$, we use the same definitions with suprema
taken over $\mathbb R^d$.
In composition estimates involving whole-space norms of an outer function
$F:\Omega\to\mathbb R$, a fixed smooth extension to $\mathbb R^m$
is assumed, where $\Omega\subset\mathbb R^m$. The norms of $F$ in
these estimates refer to this extension.
Unless stated otherwise, time-dependent
H{\"o}lder norms are spatial norms taken uniformly in time.
For a time interval
	$\mcI\subset\mathbb R$, we write
	$$
	\|f\|_{N,\mcI}
	:=
	\sup_{t\in \mcI}\|f(t,\cdot)\|_N,
	\qquad
	\|f\|_{N+\alpha,\mcI}
	:=
	\sup_{t\in \mcI}\|f(t,\cdot)\|_{N+\alpha}.
	$$
	When the time interval is clear from the context, we write
	$$
	\|f\|_N:=\sup_t\|f(t,\cdot)\|_N,
	\qquad
	\|f\|_{N+\alpha}:=\sup_t\|f(t,\cdot)\|_{N+\alpha}.
	$$
	We also use $\nabla^N f$ for the collection of all
	spatial derivatives of order $N$. All estimates involving $\nabla^N f$
	are understood componentwise.
	
	For $M\in\mathbb N_0$, we use the convention
	$\|f\|_{M+0}:=\|f\|_M$. The interpolation inequality takes the following form. If
	$N,N_1,N_2\in\mathbb N_0$, $0\leq\alpha,\alpha_1,\alpha_2<1$,
	$0\leq\vartheta\leq1$, and
	$
	N+\alpha
	=
	\vartheta(N_1+\alpha_1)
	+(1-\vartheta)(N_2+\alpha_2),
	$
	then
	$$
	\|f\|_{N+\alpha}
	\leq
	C
	\|f\|_{N_1+\alpha_1}^{\vartheta}
	\|f\|_{N_2+\alpha_2}^{1-\vartheta},
	$$
	where $C$ depends only on the dimension and the interpolation orders and exponents.
	The standard product estimate is
	\begin{equation}\label{eq:holder-product-tame}
	\|fg\|_{N+\alpha}
	\leq
	C\Big(
	\|f\|_{N+\alpha}\|g\|_0
	+
	\|f\|_0\|g\|_{N+\alpha}
	\Big)\leq
	C\Big(
	\|f\|_{N+\alpha}\|g\|_\alpha
	+
	\|f\|_\alpha\|g\|_{N+\alpha}
	\Big).
	\end{equation}
	For scalar functions $g_1,\ldots,g_k$, the telescoping identity gives, for $k\geq1$ and $0<\alpha<1$,
	\begin{equation}
		\left\|\prod_{i=1}^k g_i\right\|_{\alpha}
		\leq
		\sum_{i=1}^k
		\|g_i\|_{\alpha}
		\prod_{j\neq i}\|g_j\|_0 .
		\label{eq:holder-product-finite}
	\end{equation}
	For every $N\in\mathbb N_0$, we have
	\begin{equation}
		\left\|\prod_{i=1}^k g_i\right\|_{N+\alpha}
		\lesssim_{k,N}
		\sum_{\substack{\gamma_1,\ldots,\gamma_k\in\mathbb N_0^d\\
				|\gamma_1|+\cdots+|\gamma_k|\leq N}}
		\sum_{i=1}^k
		\|\partial^{\gamma_i}g_i\|_{\alpha}
		\prod_{j\neq i}
		\|\partial^{\gamma_j}g_j\|_0 .
		\label{eq:holder-product-finite-Cnalpha}
	\end{equation}

	For $k,N\in\mathbb N_0$ and $0<\eta\leq\alpha<1$, the embedding $C^{N+\eta}\hookrightarrow C^N$
	and \eqref{eq:holder-product-finite-Cnalpha} give
	\begin{equation}
		\begin{aligned}
			\left\|F\prod_{\nu=1}^{k}g_\nu\right\|_{N+\alpha}
			\lesssim_{k,N}\sum_{N_0+\cdots+N_k=N}
			\Bigg(&\|F\|_{N_0+\alpha}\prod_{\nu=1}^{k}\|g_\nu\|_{N_\nu}
			+\|F\|_{N_0+\eta}\sum_{i=1}^{k}\|g_i\|_{N_i+\alpha}
			\prod_{\substack{\nu=1\\\nu\neq i}}^{k}\|g_\nu\|_{N_\nu}\Bigg).
		\end{aligned}
		\label{eq:holder-product-lower-exponent}
	\end{equation}
	The sum is over nonnegative integer spatial orders. The same estimates
	hold for tensor contractions and for norms taken uniformly in time.
	\subsection{Littlewood--Paley decomposition and Besov spaces}\label{app:littlewood-paley-besov}

We define the nonhomogeneous Littlewood--Paley decomposition and
Besov norms on the torus, following \cite{BCD11,Zhao25}. Vector- and
tensor-valued norms are taken componentwise.

\begin{df}
	Let $d\geq2$ and let $\mathbb T^d=(\mathbb R/2\pi\mathbb Z)^d$. For $f\in\mathcal D'(\mathbb T^d)$ and $k\in\mathbb Z^d$, define
	$$
	\widehat f(k)=\frac{1}{(2\pi)^d}\left\langle f,e^{-i k\cdot x}\right\rangle.
	$$
	Let $\psi:\mathbb R^d\to\mathbb R$ be a smooth radial function such that $\supp\psi\subset B_{3/2}(0)$ and $\psi(\xi)=1$ for $\xi\in B_1(0)$. For $j\geq0$, set
	$$
	\chi_j(\xi)=\psi(2^{-j}\xi)-\psi(2^{-j+1}\xi).
	$$
	For $j\geq0$, define
	\begin{equation}\label{eq:app-lp-projections}
		\Delta_j f(x)=\sum_{k\in\mathbb Z^d}\chi_j(k)\widehat f(k)e^{i k\cdot x}.
	\end{equation}
	The zero-frequency projection is defined by
	$$
	\Delta_{-1}f=\widehat f(0).
	$$
	We set $\Delta_j f=0$ for $j<-1$. For $j\geq-1$, the low-frequency cutoff operators are defined by
	\begin{equation}\label{eq:app-low-frequency-cutoff}
		S_j f=\sum_{i\leq j}\Delta_i f.
	\end{equation}
	Equivalently, for $j\geq0$,
	$$
	S_j f(x)=\sum_{k\in\mathbb Z^d}\psi(2^{-j}k)\widehat f(k)e^{i k\cdot x}.
	$$
	Then, in the sense of distributions,
	$
	f=\sum_{j\geq-1}\Delta_j f.
	$
\end{df}
For $j,\ell\geq0$,
$$
\supp\widehat{\Delta_j f}\subset \mathbb Z^d\cap\{2^{j-1}\leq |k|\leq2^{j+1}\},\qquad \Delta_j\Delta_\ell f=0\quad\text{if }|j-\ell|>1.
$$

\begin{df}\label{def:app-besov-spaces}
	Let $s\in\mathbb R$ and $1\leq p,q\leq\infty$. The nonhomogeneous Besov space $B^s_{p,q}(\mathbb T^d)$ consists of all $f\in\mathcal D'(\mathbb T^d)$ such that
	\begin{equation}\label{eq:app-besov-norm}
		\|f\|_{B^s_{p,q}}
		=
		\left\|
		\left(
		2^{js}\|\Delta_j f\|_{L^p(\mathbb T^d)}
		\right)_{j\geq-1}
		\right\|_{\ell^q}
		<\infty.
	\end{equation}
	In particular, when $q=\infty$,
	\begin{equation}\label{eq:app-besov-infty-norm}
		\|f\|_{B^s_{p,\infty}}
		=
		\sup_{j\geq-1}2^{js}\|\Delta_j f\|_{L^p(\mathbb T^d)}.
	\end{equation}
	For zero-mean distributions, $\Delta_{-1}f=0$, and the above norm reduces to
	$$
	\|f\|_{B^s_{p,q}}
	=
	\left\|
	\left(
	2^{js}\|\Delta_j f\|_{L^p(\mathbb T^d)}
	\right)_{j\geq0}
	\right\|_{\ell^q}.
	$$
	We apply the periodic inverse Laplacian $\Delta^{-1}$ to the zero-mean part of each field and choose the solution to have zero mean.
\end{df}

\begin{rem}\label{rem:app-holder-besov-equivalence}
	For $0<\alpha<1$, one has $\|f\|_{B^\alpha_{\infty,\infty}}\sim \|f\|_{C^\alpha}$. For $N\in\mathbb N$ and $0<\alpha<1$, one also has $\|f\|_{B^{N+\alpha}_{\infty,\infty}}\sim \|f\|_{C^{N+\alpha}}$. 
\end{rem}

\begin{lem}\label{lem:app-basic-besov-estimates}
	Let $s\in\mathbb R$ and $1\leq p,q\leq\infty$.
	
	\emph{(i)} Let $m\in\mathbb R$ and let $T$ be the periodic Fourier multiplier defined on the nonzero modes by $\widehat{Tf}(k)=\sigma(k)\widehat f(k)$. The symbol $\sigma\in C^\infty(\mathbb R^d\setminus\{0\})$ satisfies $|\partial_\xi^\gamma \sigma(\xi)|\leq C_\gamma|\xi|^{m-|\gamma|}$ for every multi-index $\gamma$ and $|\xi|\geq1/2$. Let $f\in B^s_{p,q}(\mathbb T^d)$. If $\sigma$ extends smoothly to the origin, set $\widehat{Tf}(0)=\sigma(0)\widehat f(0)$. Otherwise, assume that $f$ has zero mean and set $\widehat{Tf}(0)=0$. Then
	\begin{equation}\label{eq:app-multiplier-besov-estimate}
		\|Tf\|_{B^{s-m}_{p,q}}\leq C\|f\|_{B^s_{p,q}}.
	\end{equation}
	
	\emph{(ii)} If $s>0$ and $f,g\in B^s_{\infty,\infty}(\mathbb T^d)$, then
	\begin{equation}\label{eq:app-besov-product-estimate}
		\|fg\|_{B^s_{\infty,\infty}}
		\leq C_s\left(
		\|f\|_{L^\infty}\|g\|_{B^s_{\infty,\infty}}
		+\|g\|_{L^\infty}\|f\|_{B^s_{\infty,\infty}}
		\right).
	\end{equation}
	
	\emph{(iii)} In particular, for $N\in\mathbb N$, $0<\alpha<1$, and $f,g\in B^{N+\alpha}_{\infty,\infty}(\mathbb T^d)$,
	\begin{equation}\label{eq:app-besov-tame-product-estimate}
		\|fg\|_{B^{N+\alpha}_{\infty,\infty}}
		\leq C_{N,\alpha}\left(
		\|f\|_{B^{N+\alpha}_{\infty,\infty}}\|g\|_{B^\alpha_{\infty,\infty}}
		+\|f\|_{B^\alpha_{\infty,\infty}}\|g\|_{B^{N+\alpha}_{\infty,\infty}}
		\right).
	\end{equation}
\end{lem}

	The following lemma estimates the spatial norm of a vector field using its divergence,
curl, and spatial mean. We apply it recursively to material derivatives.

	\begin{lem}
		\label{lem:periodic-Hodge-estimates}
		Fix $0<\alpha<1$ and let $Y$ be a smooth vector field on
		$\mathbb T^3$. Then, for every integer $L\geq1$,
		$$
		\|Y\|_{L+\alpha}
		\leq
		C_{L,\alpha}\left(
		\|\Div Y\|_{L-1+\alpha}
		+\|\curl Y\|_{L-1+\alpha}
		+\left|\fint_{\mathbb T^3}Y\,\rd x\right|
		\right).
		$$
		At the lowest spatial order, one has
		$$
		\|Y\|_{\alpha}
		\leq
		C_\alpha\left(
		\|\Div Y\|_{B^{-1+\alpha}_{\infty,\infty}}
		+\|\curl Y\|_{B^{-1+\alpha}_{\infty,\infty}}
		+\left|\fint_{\mathbb T^3}Y\,\rd x\right|
		\right).
		$$
	\end{lem}
	
	\begin{proof}
		Set
		$
		\overline Y:=\fint_{\mathbb T^3}Y\,\rd x.
		$
		The vector identity
		$$
		\curl\curl Y
		=
		\nabla\Div Y-\Delta Y
		$$
		and the fact that $Y-\overline Y$ has zero spatial mean give
		$$
		Y-\overline Y
		=
		(-\Delta)^{-1}\curl(\curl Y)
		-
		(-\Delta)^{-1}\nabla(\Div Y).
		$$
		The operators
		$
		(-\Delta)^{-1}\curl
		$
		and
		$
		(-\Delta)^{-1}\nabla
		$
		are periodic Fourier multipliers of order minus one. Hence Remark~\ref{rem:app-holder-besov-equivalence} and the periodic
		multiplier estimate \eqref{eq:app-multiplier-besov-estimate} imply, for every integer $L\geq1$,
		$$
		\|Y-\overline Y\|_{L+\alpha}
		\leq
		C_{L,\alpha}\left(
		\|\Div Y\|_{L-1+\alpha}
		+\|\curl Y\|_{L-1+\alpha}
		\right).
		$$
		Since
		$
		\|\overline Y\|_{L+\alpha}=|\overline Y|,
		$
		this proves the first estimate.
		
		The same multiplier estimates in $B^s_{\infty,\infty}$ give
		$$
		\|Y-\overline Y\|_{B^\alpha_{\infty,\infty}}
		\leq
		C_\alpha\left(
		\|\Div Y\|_{B^{-1+\alpha}_{\infty,\infty}}
		+\|\curl Y\|_{B^{-1+\alpha}_{\infty,\infty}}
		\right).
		$$
		Using
		$
		B^\alpha_{\infty,\infty}\approx C^\alpha
		$
		completes the proof.
	\end{proof}
	\section{Tools for convex integration}
\subsection{Estimates for transport equations}
We state the following transport estimates in dimension $d$ with initial time $t_0$, see \cite[Proposition~D.1]{BDIS15} and
\cite[Proposition~B.1]{BDSV19}.
\begin{lem}\label{estimate of the transport equation}
	Let $\mcI\subset\R$ be a time interval, let $t_0\in\mcI$, and let $f$
	be a smooth solution of
	\begin{equation}
		\left\{
		\begin{aligned}
			& \partial_tf+\upsilon\cdot\nabla f=g,\\
			&f|_{t=t_0}=f_0,
		\end{aligned}\right.
	\end{equation}
	on $\mcI\times\T^d$, where $\upsilon$ and $g$ are smooth. Assume that
	$t\in\mcI$, $t\geq t_0$, and $|t-t_0|[\upsilon]_1\leq1$. Then
	\begin{align}
		\|f(t)\|_0\leq&\|f_0\|_0+\int_{t_0}^{t}\|g(\tau,\cdot)\|_0\rd \tau,\label{transport equation 0}\\
		{}[f(t)]_{\alpha}\lesssim&\left([f_0]_{\alpha}+\int_{t_0}^{t}\|g(\tau,\cdot)\|_{\alpha}\rd \tau\right),\label{transport equation 1}
	\end{align}
	for $0\leq\alpha<1$. Moreover, for every integer $N\geq1$ and
	$0\leq\alpha<1$,
	\begin{equation}
		\begin{aligned}
			\relax[f(t)]_{N+\alpha}
			\lesssim{}&
			[f_0]_{N+\alpha}
			+
			|t-t_0|[\upsilon]_{N+\alpha}[f_0]_1+
			\int_{t_0}^{t}
			\left(
			[g(\tau,\cdot)]_{N+\alpha}
			+
			(t-\tau)[\upsilon]_{N+\alpha}
			[g(\tau,\cdot)]_1
			\right)\rd\tau.
		\end{aligned}
		\label{transport equation N}
	\end{equation}
	Let $\Phi(t,\cdot)$ denote the backward flow map, inverse to the forward flow $X$ generated by
	$\upsilon$, normalized by $X(t_0,x)=x$. Under the same assumption,
	\begin{equation}
		\|\nabla \Phi(t,\cdot)-\Id\|_0\lesssim |t-t_0|[\upsilon]_1,\quad 	[ \Phi(t,\cdot)]_N\lesssim |t-t_0|[\upsilon]_N,\quad N\geq 2.\label{PhiN}
	\end{equation}
	The estimates for $t\leq t_0$ follow by reversing time.
\end{lem}

\subsection{Geometric lemma}

We use the geometric decomposition in
\cite[Proposition~4.1]{BMNV21}, with the directions rescaled to lie in
$\mathbb Z^3$ and the coefficients rescaled accordingly. 

\begin{lem}
	\label{Geometric lemma}
	Let $B_\varepsilon(\Id)$ denote the ball of radius $\varepsilon>0$
	centered at the identity in $\mathbb S^{3\times3}$ with respect to the operator norm.
	There exist a constant $\varepsilon_0>0$, a finite set
	$\Lambda\subset\mathbb{Z}^3\setminus\{0\}$, and smooth positive
	functions
	$
	\gamma_\xi:
	B_{\varepsilon_0}(\Id)\to\mathbb{R},
	$
	for each $\xi\in\Lambda$, such that
	$$
	R
	=
	\sum_{\xi\in\Lambda}
	\gamma_\xi^2(R)\,\xi\otimes\xi,
	$$
	for every $R\in B_{\varepsilon_0}(\Id)$.
\end{lem}

\subsection{Singular integral operators}

We use the following class of periodic Calder\'on--Zygmund operators.
Let
$
K:\mathbb{R}^d\setminus\{0\}\to\mathbb{R}
$
be a kernel of the form
$
K(z)
=
\frac{\Omega(z/|z|)}{|z|^d},
$
where $\Omega\in C^\infty(\mathbb{S}^{d-1})$ satisfies
$
\int_{\mathbb{S}^{d-1}}\Omega(\theta)\,\rd S(\theta)=0.
$
Up to an additive constant, its periodicization is
$$
K_{\mathbb{T}^d}(z)
:=
K(z)
+
\sum_{n\in\mathbb{Z}^d\setminus\{0\}}
\bigl(
K(z+2\pi n)-K(2\pi n)
\bigr).
$$
Choose the additive constant so that $K_{\mathbb{T}^d}$ has zero mean in the principal value sense.
For a smooth $\mathbb{T}^d$-periodic function $f$, set
$$
T_Kf(x)
:=
\mathrm{p.v.}
\int_{\mathbb{T}^d}
K_{\mathbb{T}^d}(x-y)f(y)\,\rd y.
$$
Then $T_K1=0$, and $T_K$ acts componentwise on vector- and tensor-valued functions.
By \cite[Proposition~C.1]{BDSV19}, $T_K$ is bounded on $C^\alpha$ for $0<\alpha<1$.
The operator is also bounded on $C^{N+\alpha}$ for every $N\in\mathbb N_0$, since it commutes with spatial derivatives.
The proofs of this estimate and the following proposition apply in dimension $d$.

\begin{prop}[{cf.\ \cite[Proposition~D.1]{BDSV19}}]
	\label{prop:CZ commutator}
	Let $\alpha\in(0,1)$ and $N\in\mathbb N_0$, and let $T_K$ be the
	periodic Calder\'on--Zygmund operator defined above, and let
	$Y\in C^{N+1+\alpha}(\mathbb{T}^d;\mathbb{R}^d)$ be a vector field.
	Then, for every $f\in C^{N+\alpha}(\mathbb{T}^d)$,
	$$
	\|[T_K,Y\cdot\nabla]f\|_{N+\alpha}
	\lesssim
	\|Y\|_{1+\alpha}\|f\|_{N+\alpha}
	+
	\|Y\|_{N+1+\alpha}\|f\|_{\alpha}.
	$$
	The commutator is understood in the sense of distributions.
	The estimate is uniform in $f$ and $Y$.
\end{prop}

\subsection{Fa\`a di Bruno and composition estimates}
\label{sec:faa-di-bruno-composition}

We recall the multivariable Fa\`a di Bruno formula
\cite{BMNV21}. Let $m,N\in\mathbb N$ and let $g(x)=f(h(x))$, where $f:\mathbb R^m\to\mathbb R$ and $h:\mathbb R^d\to\mathbb R^m$ are $C^N$ functions. Let $\theta\in\mathbb N_0^d$ satisfy $|\theta|=N$, and let $\beta\in\mathbb N_0^m$ satisfy $1\leq|\beta|\leq N$. We define
\begin{align*}
	\mathfrak P(\theta,\beta)
	&= \Bigg\{ (k_1,\ldots,k_N; \ell_1,\ldots,\ell_N) \in (\N_0^m)^N \times (\N_0^d)^N \colon \exists s \mbox{ with }1\leq s \leq N \mbox{ s.t. } \\
	&\qquad   \qquad  |k_j|, |\ell_j| > 0 \mbox{ for }1 \leq j \leq s, \, k_j=\ell_j=0 \mbox{ for }s<j\leq N, \, 0 \prec \ell_{1} \prec \ldots \prec \ell_s,  \sum_{j=1}^s k_j =  \beta, \sum_{j=1}^s |k_j| \ell_j = \theta \Bigg\} \, .
\end{align*}
Here $\prec$ denotes a fixed total ordering of the multi-indices, with
$0\prec\ell$ for every nonzero $\ell$. We use the notation
$$
(\partial^\ell h)^k := \prod_{i=1}^m (\partial^\ell h_i)^{k_i}, \qquad k!:=\prod_{i=1}^m k_i!, \qquad \ell!:=\prod_{i=1}^d\ell_i!.
$$
Then
\begin{equation}
	\partial^\theta g(x) = \theta! \sum_{1\leq|\beta|\leq N} (\partial^\beta f)(h(x)) \sum_{\mathfrak P(\theta,\beta)} \prod_{j=1}^s \frac{(\partial^{\ell_j}h(x))^{k_j}} {k_j!(\ell_j!)^{|k_j|}} .
	\label{eq:faa-di-bruno}
\end{equation}

We use the composition estimates in \cite[Proposition~A.1]{GR24}
(see also \cite[Proposition~A.1]{BDSV19}).

\begin{lem}
	\label{lem:BDSV-composition}
	Let $F:\Omega\rightarrow\mathbb R$ and
	$\Psi:\mathbb R^n\rightarrow\Omega$ be smooth functions, with
	$\Omega\subset\mathbb R^m$ open. Then, for every
	$N\in\mathbb N\setminus\{0\}$, the following estimates hold with
	implicit constants depending only on $N,m,n$,
	\begin{equation}
		\begin{aligned}
			\|\nabla^N(F\circ\Psi)\|_0
			&\lesssim
			\|\nabla F\|_0\|\nabla\Psi\|_{N-1}
			+\|\nabla F\|_{N-1}\|\Psi\|_0^{N-1}\|\Psi\|_N,\\
			\|\nabla^N(F\circ\Psi)\|_0
			&\lesssim
			\|\nabla F\|_0\|\nabla\Psi\|_{N-1}
			+\|\nabla F\|_{N-1}\|\nabla\Psi\|_0^{N}.
		\end{aligned}
		\label{est:holder-norm-composition}
	\end{equation}
	Each estimate holds whenever the norms on its right-hand side are finite.

\end{lem}

\begin{lem}\label{lem:holder-composition-general}
	Fix $m,n\in\mathbb N$, $0<\alpha<1$, and an open set
	$\Omega\subset\mathbb R^m$. Suppose that $F:\Omega\to\mathbb R$
	and $\Psi:\mathbb R^n\to\Omega$ are smooth and that the norms
	appearing below are finite. Then, for each integer $N\geq1$,
	\begin{equation}
		\begin{aligned}
			{}[\nabla^N(F\circ\Psi)]_{\alpha}
			\lesssim{}&
			\sum_{k=1}^N
			[\nabla^k F]_{\alpha}\|\nabla\Psi\|_0^\alpha
			\sum_{\substack{M_1+\cdots+M_k=N\\ M_j\geq1}}
			\prod_{j=1}^k\|\nabla^{M_j}\Psi\|_0
			\\
			&+
			\sum_{k=1}^N
			\|\nabla^k F\|_0
			\sum_{\substack{M_1+\cdots+M_k=N\\ M_j\geq1}}
			\sum_{j=1}^k
			[\nabla^{M_j}\Psi]_{\alpha}
			\prod_{\ell\neq j}\|\nabla^{M_\ell}\Psi\|_0.
		\end{aligned}
		\label{est:holder-composition-general}
	\end{equation}
	For $N=0$, one has
	\begin{equation}
		[F\circ\Psi]_{\alpha}
		\lesssim
		[F]_{\alpha}\|\nabla\Psi\|_0^\alpha.
		\label{est:holder-composition-zero}
	\end{equation}
	We use the same extension and norm conventions as in
	Lemma~\ref{lem:BDSV-composition}.

\end{lem}

\begin{proof}
	By \eqref{eq:faa-di-bruno}, every term in $\nabla^N(F\circ\Psi)$ is a
	constant multiple of a tensor contraction of
	$
	((\nabla^kF)\circ\Psi)\prod_{j=1}^k\nabla^{M_j}\Psi,
	$
	where $M_1+\cdots+M_k=N$ and $M_j\geq1$. Apply the H{\"o}lder seminorm product inequality
	to this contraction
	and use
	$
	[(\nabla^kF)\circ\Psi]_\alpha
	\leq[\nabla^kF]_\alpha\|\nabla\Psi\|_0^\alpha,
	$
	together with
	$
	\|(\nabla^kF)\circ\Psi\|_0\leq\|\nabla^kF\|_0
	$
	to obtain the two terms in \eqref{est:holder-composition-general}.
	When $N=0$, \eqref{est:holder-composition-zero} follows from
	$
	|\Psi(x)-\Psi(y)|\leq\|\nabla\Psi\|_0|x-y|.
	$
\end{proof}

\begin{lem}
	\label{lem:composition-scale}
	Fix $0<\alpha<1$, $N_*\in\N_0$, and a time interval $\mcI\subset\R$.
	Let $\Phi:\mcI\times\T^d\to\T^d$ and
	$f:\mcI\times\T^d\to\R$ be smooth. We define the composition pointwise
	in time by
	$
	(f\circ\Phi)(t,x):=f(t,\Phi(t,x)).
	$
	Assume that for some $\lambda_\Phi\geq1$,
	\begin{equation}
		\|\nabla\Phi\|_{r,\mcI}
		\lesssim
		\lambda_\Phi^r,
		\qquad
		\|\nabla\Phi\|_{r+\alpha,\mcI}
		\lesssim
		\lambda_\Phi^{r+\alpha},
		\label{eq:Phi-scale-holder}
	\end{equation}
	for every integer $0\leq r\leq N_*$. Then, for every integer
	$0\leq N\leq N_*$,
	\begin{align}
		\|f\circ\Phi\|_{N,\mcI}
		\lesssim&
		\sum_{M=0}^N
		\lambda_\Phi^{N-M}\|f\|_{M,\mcI},
		\label{eq:composition-scale-CN}\\
		\|f\circ\Phi\|_{N+\alpha,\mcI}
		\lesssim&
		\sum_{M=0}^N
		\lambda_\Phi^{N-M}\|f\|_{M+\alpha,\mcI}
		+
		\sum_{M=0}^N
		\lambda_\Phi^{N+\alpha-M}\|f\|_{M,\mcI}.
		\label{eq:composition-scale-holder}
	\end{align}
	
	Moreover, fix an integer $0\leq N\leq N_*$ and suppose that, for some
	$\mC_f\geq0$ and $1\leq\lambda_f\leq\tilde\lambda_f$, one has
	\begin{equation}
		\|f\|_{M,\mcI}
		\lesssim
		\mC_f
		\MM{M,N_*,\lambda_f,\tilde\lambda_f},
		\qquad
		\|f\|_{M+\alpha,\mcI}
		\lesssim
		\mC_f
		\MM{M,N_*,\lambda_f,\tilde\lambda_f}
		\tilde\lambda_f^\alpha
		\label{eq:f-scale-assumption}
	\end{equation}
	for every integer $0\leq M\leq N$. Set
	$
	\lambda_*:=\max\{\lambda_\Phi,\lambda_f\},\
	\tilde\lambda_*:=\max\{\lambda_\Phi,\tilde\lambda_f\}.
	$
	Then
	\begin{equation}
		\|f\circ\Phi\|_{N+\alpha,\mcI}
		\lesssim
		\mC_f
		\MM{N,N_*,\lambda_*,\tilde\lambda_*}
		\tilde\lambda_*^\alpha.
		\label{eq:composition-scale-geometric}
	\end{equation}
	The constants depend on the fixed constants in \eqref{eq:Phi-scale-holder};
	in \eqref{eq:composition-scale-geometric}, they may also depend on the
	constant in \eqref{eq:f-scale-assumption}. With these constants fixed,
	they are independent of $f$, $\mC_f$, and the scale parameters.

	The same conclusions hold
	componentwise when $f$ is vector- or tensor-valued.
\end{lem}

\begin{proof}
	Let $N\geq1$. By \eqref{eq:faa-di-bruno} and
	\eqref{eq:Phi-scale-holder}, each product of $M$ derivatives of $\Phi$
	with total order $N$ is bounded by a constant times $\lambda_\Phi^{N-M}$. Thus
	$$
	\|\nabla^N(f\circ\Phi)\|_{0,\mcI}
	\lesssim
	\sum_{M=1}^N
	\lambda_\Phi^{N-M}\|f\|_{M,\mcI}.
	$$
	Applying the same argument with $N$ replaced by each
	$1\leq r\leq N$, and using $\lambda_\Phi\geq1$, controls all spatial
	derivatives of orders at most $N$. Together with
	$
	\|f\circ\Phi\|_{0,\mcI}\leq\|f\|_{0,\mcI},
	$
	we obtain \eqref{eq:composition-scale-CN}.
	
	We next prove the H{\"o}lder estimate. For $N=0$,
	\eqref{eq:composition-scale-holder} follows from
	\eqref{est:holder-composition-zero}--\eqref{eq:Phi-scale-holder} and
	the corresponding $C^0$ estimate. Let $N\geq1$. At a fixed time, apply
	Lemma~\ref{lem:holder-composition-general} to lifts on $\mathbb R^d$.
	The norms without a time subscript below are taken at this fixed time.
	Using
	$
	\|\nabla^r\Phi\|_0\lesssim\lambda_\Phi^{r-1},\
	[\nabla^r\Phi]_\alpha\lesssim\lambda_\Phi^{r-1+\alpha},
	$
	we have
	$$
	[\nabla^N(f\circ\Phi)]_\alpha
	\lesssim
	\sum_{M=1}^N
	\lambda_\Phi^{N-M}\|f\|_{M+\alpha,\mcI}
	+
	\sum_{M=1}^N
	\lambda_\Phi^{N+\alpha-M}\|f\|_{M,\mcI}.
	$$
	Taking the supremum in time and adding the $C^N$ bound
	\eqref{eq:composition-scale-CN} proves
	\eqref{eq:composition-scale-holder}, since $\lambda_\Phi^\alpha\geq1$.
	
	Substitute \eqref{eq:f-scale-assumption} into
	\eqref{eq:composition-scale-holder}. Since $M\leq N\leq N_*$, the
	definition of $\mathcal M$ gives
	$
	\MM{M,N_*,\lambda_f,\tilde\lambda_f}=\lambda_f^M
	$
	and
	$
	\MM{N,N_*,\lambda_*,\tilde\lambda_*}=\lambda_*^N.
	$
	With the definitions of $\lambda_*$ and $\tilde\lambda_*$, we obtain
	\eqref{eq:composition-scale-geometric} after summing over $M$.
\end{proof}

\subsection{Iterated material derivatives and commutator estimates}
\label{sec:holder-material-commutators}
In this subsection,
we establish $C^N$ and $C^{N+\alpha}$ estimates for ordered products,
commutators, and mixed spatial and material derivatives.
We begin with products of differential operators and powers of their sum,
following \cite[Lemma~A.10]{BMNV21} (see also \cite[Lemma~A.6]{GKN24}).
We then estimate commutators with material derivatives and use them to
control ordered mixed derivatives, allowing different derivative thresholds
for the velocity and the function.
The proofs adapt \cite[Lemmas~A.12--A.14]{BMNV21}
(see also \cite[Lemma~A.7]{GKN24}) to the total-order hypotheses below.
Finally, we adapt \cite[Lemma~A.2]{GKN26Strong} to estimate material derivatives
after a change of velocity under total-order hypotheses.
In a product of
operators, we use the notation
$$\left(\prod_{i=1}^kT_i\right)F=T_1(T_2(\cdots T_kF)).$$

\begin{lem}
	\label{lem:operator-product-estimates}
	Fix $\widetilde N,\widetilde M\in\mathbb N_0$, $N_*\in\mathbb N$, and a time interval
	$\mcI\subset \mathbb R$. Let $\upsilon:\mcI\times\mathbb T^d\to \mathbb R^d$
	and $f:\mcI\times\mathbb T^d\to\mathbb R$ be smooth,
	and let $B$ be a derivation on smooth functions satisfying
	$
	B(fg)=(Bf)g+f(Bg).
	$
	Set
	$
	A:=\upsilon\cdot\nabla .
	$
	For every integer $k\geq1$ and $\theta,\beta\in\mathbb N_0^k$, we write
	$
	|\theta|:=\sum_{i=1}^k\theta_i,
	|\beta|:=\sum_{i=1}^k\beta_i .
	$
	Fix $\mC_{\upsilon}\geq0$,
	$1\leq \lambda_{\upsilon}\leq \tilde\lambda_{\upsilon}$, and
	$1\leq \mu_{\upsilon}\leq \tilde\mu_{\upsilon}$.
	Assume that, whenever $|\theta|+|\beta|\leq N_*$, $\upsilon$ satisfies
	\begin{equation}
		\left\|
		\left(
		\prod_{i=1}^k \nabla^{\theta_i}B^{\beta_i}
		\right)\upsilon
		\right\|_{0,\mcI}
		\lesssim
		\mC_{\upsilon}\,
		\MM{|\theta|,\widetilde N,\lambda_{\upsilon},\tilde\lambda_{\upsilon}}
		\MM{|\beta|,\widetilde M,\mu_{\upsilon},\tilde\mu_{\upsilon}}.
		\label{eq:operator-product-upsilon}
	\end{equation}
	Fix $\mC_f\geq0$,
	$1\leq \lambda_f\leq \tilde\lambda_f$, and
	$1\leq \mu_f\leq \tilde\mu_f$.
	Assume also that, whenever $|\theta|+|\beta|\leq N_*$, $f$ satisfies
	\begin{equation}
		\left\|
		\left(
		\prod_{i=1}^k \nabla^{\theta_i}B^{\beta_i}
		\right)f
		\right\|_{0,\mcI}
		\lesssim
		\mC_f\,
		\MM{|\theta|,\widetilde N,\lambda_f,\tilde\lambda_f}
		\MM{|\beta|,\widetilde M,\mu_f,\tilde\mu_f}.
		\label{eq:operator-product-f}
	\end{equation}
	Define
	$$
	\lambda:=\max\{\lambda_f,\lambda_{\upsilon}\},
	\quad
	\tilde\lambda:=\max\{\tilde\lambda_f,\tilde\lambda_{\upsilon}\},\quad \mu:=\max\{\mu_f,\mu_{\upsilon}\},
	\quad
	\tilde\mu:=\max\{\tilde\mu_f,\tilde\mu_{\upsilon}\}.
	$$
	Then, for every integer $N\geq0$ such that
	$N+|\theta|+|\beta|\leq N_*$, one has
	\begin{align}
		\left\|
		\left(
		\prod_{i=1}^k A^{\theta_i}B^{\beta_i}
		\right)f
		\right\|_{N,\mcI}
		&\lesssim
		\mC_f\mC_{\upsilon}^{|\theta|}
		\MM{N+|\theta|,\widetilde N,\lambda,\tilde\lambda}
		\MM{|\beta|,\widetilde M,\mu,\tilde\mu}
		\label{eq:operator-product-main}
		\\
		&\lesssim
		\mC_f
		\MM{N,\widetilde N,\lambda,\tilde\lambda}
		(\mC_{\upsilon}\tilde\lambda)^{|\theta|}
		\MM{|\beta|,\widetilde M,\mu,\tilde\mu}
		\notag
		\\
		&\lesssim
		\mC_f
		\MM{N,\widetilde N,\lambda,\tilde\lambda}
		\MM{|\theta|+|\beta|,\widetilde M,\max\{\mu,\mC_{\upsilon}\tilde\lambda\},\max\{\tilde\mu,\mC_{\upsilon}\tilde\lambda\}}.
		\label{eq:operator-product-absorbed}
	\end{align}
	Consequently, for all integers $N,M\geq0$ with $N+M\leq N_*$,
	\begin{equation}
		\left\|
		(A+B)^M f
		\right\|_{N,\mcI}
		\lesssim
		\mC_f
		\MM{N,\widetilde N,\lambda,\tilde\lambda}
		\MM{M ,\widetilde M,\ \max\{\mu,\mC_{\upsilon}\tilde\lambda\},\max\{\tilde\mu,\mC_{\upsilon}\tilde\lambda\}}.
		\label{eq:operator-product-A-plus-B}
	\end{equation}

\end{lem}

\begin{proof}
	For every integer $r\geq1$,
	\begin{equation}
		A^r
		=
		(\upsilon\cdot\nabla)^r
		=
		\sum_{j=1}^r
		F_{j,r}\nabla^j,
		\qquad
		F_{j,r}
		=
		\sum_{\substack{
				\zeta\in\mathbb N_0^r\\
				|\zeta|=r-j
		}}
		c_{r,j,\zeta}
		\prod_{\ell=1}^r
		\left(
		\nabla^{\zeta_\ell}\upsilon
		\right),
		\label{eq:operator-product-expansion}
	\end{equation}
	where the $c_{r,j,\zeta}$ are combinatorial coefficient tensors.
	We use the conventions
	$
	F_{0,0}:=1,
	A^0=F_{0,0}\nabla^0
	$.
	Pairs $(\theta_i,\beta_i)=(0,0)$ may be omitted. If all pairs are trivial, \eqref{eq:operator-product-main} follows
	from \eqref{eq:operator-product-f}. Hence, we may assume that
	no pair is trivial, so that
	$
	k\leq|\theta|+|\beta|\leq N_*.
	$
	For $0\leq r\leq N$, \eqref{eq:operator-product-expansion} gives
	\begin{align}
		&\nabla^r
		\prod_{i=1}^k
		A^{\theta_i}B^{\beta_i}
		=
		\sum_{\substack{
				\gamma\in\mathbb N_0^k\\
				\gamma_i=0\ \text{if }\theta_i=0\\
				1\leq\gamma_i\leq\theta_i\ \text{if }\theta_i\geq1
		}}
		\nabla^r
		\prod_{i=1}^k
		\left(
		F_{\gamma_i,\theta_i}
		\nabla^{\gamma_i}B^{\beta_i}
		\right)
		\notag\\
		=&
		\sum_{\substack{
				\gamma\in\mathbb N_0^k\\
				\gamma_i=0\ \text{if }\theta_i=0\\
				1\leq\gamma_i\leq\theta_i\ \text{if }\theta_i\geq1
		}}
		\sum_{\substack{
				0\leq N'\leq r+|\gamma|\\
				0\leq M'\leq|\beta|
		}}
		\sum_{\substack{
				\delta,\kappa\in\mathbb N_0^k\\
				|\delta|=r+|\gamma|-N'\\
				|\kappa|=|\beta|-M'
		}}
		\left(
		\prod_{i=1}^k
		\left(
		\sum_{\substack{
				\delta_i',\kappa_i'\in\mathbb N_0^k\\
				|\delta_i'|=\delta_i\\
				|\kappa_i'|=\kappa_i
		}}
		\widetilde c_{(\ldots)}
		\left(
		\prod_{\ell_i=1}^k
		\nabla^{\delta_{i,\ell_i}'}
		B^{\kappa_{i,\ell_i}'}
		\right)
		F_{\gamma_i,\theta_i}
		\right)
		\right)
		\left(
		\sum_{\substack{
				\mathbf h,\mathbf N\in\mathbb N_0^k\\
				|\mathbf h|=N'\\
				|\mathbf N|=M'
		}}
		\overline c_{(\ldots)}
		\prod_{s=1}^k
		\nabla^{h_s}B^{N_s}
		\right),
		\label{eq:operator-product-operator-expansion}
	\end{align}
	where the coefficients
	$\widetilde c_{(\ldots)}$ and $\overline c_{(\ldots)}$ are combinatorial coefficient tensors.
	Applying the operator expansion to $f$ in
	\eqref{eq:operator-product-expansion}--%
	\eqref{eq:operator-product-operator-expansion}, we have
	\begin{align}
		\nabla^r
		\left(
		\prod_{i=1}^k
		A^{\theta_i}B^{\beta_i}
		\right)f
		=&
		\sum_{\substack{
				\gamma\in\mathbb N_0^k\\
				\gamma_i=0\ \text{if }\theta_i=0\\
				1\leq\gamma_i\leq\theta_i\ \text{if }\theta_i\geq1
		}}
		\sum_{\substack{
				0\leq N'\leq r+|\gamma|\\
				0\leq M'\leq|\beta|
		}}
		\sum_{\substack{
				\delta,\kappa\in\mathbb N_0^k\\
				|\delta|=r+|\gamma|-N'\\
				|\kappa|=|\beta|-M'
		}}
		\left(
		\sum_{\substack{
				\mathbf h,\mathbf N\in\mathbb N_0^k\\
				|\mathbf h|=N'\\
				|\mathbf N|=M'
		}}
		\overline c_{(\ldots)}
		\left(
		\prod_{s=1}^k
		\nabla^{h_s}B^{N_s}
		\right)f
		\right)
		\notag\\
		&\hspace{50pt}\cdot
		\left(
		\prod_{i=1}^k
		\left(
		\sum_{\substack{
				\delta_i',\kappa_i'\in\mathbb N_0^k\\
				|\delta_i'|=\delta_i\\
				|\kappa_i'|=\kappa_i
		}}
		\widetilde c_{(\ldots)}
		\left(
		\prod_{\ell_i=1}^k
		\nabla^{\delta_{i,\ell_i}'}
		B^{\kappa_{i,\ell_i}'}
		\right)
		\left(
		\sum_{\substack{
				\zeta_i\in\mathbb N_0^{\theta_i}\\
				|\zeta_i|=\theta_i-\gamma_i
		}}
		c_{(\ldots)}
		\prod_{r_i=1}^{\theta_i}
		\left(
		\nabla^{\zeta_{i,r_i}}\upsilon
		\right)
		\right)
		\right)
		\right).
		\label{eq:operator-product-full-expansion}
	\end{align}
	After expanding the derivatives of the coefficient products, each term
	in \eqref{eq:operator-product-full-expansion} contains $|\theta|$ factors
	involving $\upsilon$ and one factor involving $f$. Their total spatial
	and $B$-orders are $r+|\theta|$ and $|\beta|$, respectively. Since
	$r+|\theta|+|\beta|\leq N_*$, \eqref{eq:operator-product-upsilon}--\eqref{eq:operator-product-f} apply to every factor.
	For $k\geq1$, $N_i,K_i\in\mathbb N_0$, and $1\leq\lambda_i\leq\tilde\lambda_i$, we use
	\begin{equation}\label{eq:appendix-derivative-weight-product}
	\prod_{i=1}^k\MM{N_i,K_i,\lambda_i,\tilde\lambda_i}
	\leq\MM{\sum_{i=1}^k N_i,\min_iK_i,\max_i\lambda_i,\max_i\tilde\lambda_i}.
	\end{equation}
	By monotonicity in the scales and thresholds, the left-hand side is at most $\prod_i\MM{N_i,K,\max_i\lambda_i,\max_i\tilde\lambda_i}$, where $K:=\min_iK_i$.
	The inequality then follows from $\sum_i(N_i-K)_+\leq(\sum_i N_i-K)_+$.
	In particular, \eqref{eq:appendix-derivative-weight-product} gives
	\begin{align*}
		\prod_{i=1}^k
		\MM{
			\theta_i-\gamma_i+\delta_i,
			\widetilde N,
			\lambda,
			\tilde\lambda
		}
		&\leq
		\MM{
			|\theta|-|\gamma|+|\delta|,
			\widetilde N,
			\lambda,
			\tilde\lambda
		}=
		\MM{
			r+|\theta|-N',
			\widetilde N,
			\lambda,
			\tilde\lambda
		},\\
		\prod_{i=1}^k
		\MM{
			\kappa_i,
			\widetilde M,
			\mu,
			\tilde\mu
		}
		&\leq
		\MM{
			|\kappa|,
			\widetilde M,
			\mu,
			\tilde\mu
		}=
		\MM{
			|\beta|-M',
			\widetilde M,
			\mu,
			\tilde\mu
		}.
	\end{align*}
	We apply \eqref{eq:operator-product-upsilon}--\eqref{eq:operator-product-f}
	to each term of \eqref{eq:operator-product-full-expansion} and obtain
	\begin{align*}
		\left\|
		\nabla^r
		\left(
		\prod_{i=1}^k
		A^{\theta_i}B^{\beta_i}
		\right)f
		\right\|_{0,\mcI}\lesssim
		\mC_f\mC_{\upsilon}^{|\theta|}
		\sum_{\substack{
				0\leq N'\leq r+|\theta|\\
				0\leq M'\leq|\beta|
		}}
		\MM{N',\widetilde N,\lambda,\tilde\lambda}
		\MM{
			r+|\theta|-N',
			\widetilde N,
			\lambda,
			\tilde\lambda
		}
		\MM{M',\widetilde M,\mu,\tilde\mu}
		\MM{
			|\beta|-M',
			\widetilde M,
			\mu,
			\tilde\mu
		}.
	\end{align*}
	Estimate~\eqref{eq:appendix-derivative-weight-product} and the definition of $\mathcal M$ give \eqref{eq:operator-product-main} after summing
	over all components and $0\leq r\leq N$. The number of terms is bounded
	in terms of $d$ and $N_*$.
	We next absorb the extra spatial factors into the material scales. Estimate \eqref{eq:operator-product-absorbed} follows from
	\begin{align*}
		\MM{
			N+|\theta|,
			\widetilde N,
			\lambda,
			\tilde\lambda
		}
		&\leq
		\MM{
			N,
			\widetilde N,
			\lambda,
			\tilde\lambda
		}
		\tilde\lambda^{|\theta|},\\
		(\mC_{\upsilon}\tilde\lambda)^{|\theta|}
		\MM{
			|\beta|,
			\widetilde M,
			\mu,
			\tilde\mu
		}
		&\leq
		\MM{
			|\theta|+|\beta|,
			\widetilde M,
			\max\{\mu,\mC_{\upsilon}\tilde\lambda\},
			\max\{\tilde\mu,\mC_{\upsilon}\tilde\lambda\}
		}.
	\end{align*}

	Finally, for $M\geq1$, we expand
	$$
	(A+B)^M
	=
	\sum_{\sigma\in\{A,B\}^M}
	\sigma_1\cdots\sigma_M.
	$$
	Each ordered product can be written as $\prod_{i=1}^k A^{\theta_i}B^{\beta_i}$
	with $|\theta|+|\beta|=M$. Apply \eqref{eq:operator-product-absorbed}
	to each ordered product and sum the at most $2^M$ estimates to obtain
	\eqref{eq:operator-product-A-plus-B}. For $M=0$, use \eqref{eq:operator-product-f}.
\end{proof}

\begin{lem}
	\label{lem:operator-product-estimates-Calpha-loss}
	Fix $0<\alpha<1$, $\widetilde N,\widetilde M\in\mathbb N_0$, $N_*\in\mathbb N$,
	and a time interval $\mcI\subset\mathbb R$.
	Let $\upsilon:\mcI\times\mathbb T^d\to\mathbb R^d$
	and $f:\mcI\times\mathbb T^d\to\mathbb R$ be smooth,
	and let $B$ be a derivation on smooth functions satisfying
	$
	B(fg)=(Bf)g+f(Bg).
	$
	Set
	$
	A:=\upsilon\cdot\nabla.
	$
	For every integer $k\geq1$ and $\theta,\beta\in\mathbb N_0^k$, we write
	$
	|\theta|:=\sum_{i=1}^k\theta_i,
	|\beta|:=\sum_{i=1}^k\beta_i.
	$
	Fix $\mC_{\upsilon}\geq0$,
	$1\leq\lambda_{\upsilon}\leq\tilde\lambda_{\upsilon}$, and
	$1\leq\mu_{\upsilon}\leq\tilde\mu_{\upsilon}$.
	Assume that, whenever $|\theta|+|\beta|\leq N_*$, $\upsilon$ satisfies
	\begin{align}
		\left\|
		\left(
		\prod_{i=1}^k \nabla^{\theta_i}B^{\beta_i}
		\right)\upsilon
		\right\|_{0,\mcI}
		&\lesssim
		\mC_{\upsilon}
		\MM{|\theta|,\widetilde N,\lambda_{\upsilon},\tilde\lambda_{\upsilon}}
		\MM{|\beta|,\widetilde M,\mu_{\upsilon},\tilde\mu_{\upsilon}},
		\label{eq:operator-product-Calpha-loss-upsilon-C0}
		\\
		\left\|
		\left(
		\prod_{i=1}^k \nabla^{\theta_i}B^{\beta_i}
		\right)\upsilon
		\right\|_{\alpha,\mcI}
		&\lesssim
		\mC_{\upsilon}
		\La{|\theta|,\widetilde N,\lambda_{\upsilon},\tilde\lambda_{\upsilon}}
		\MM{|\theta|,\widetilde N,\lambda_{\upsilon},\tilde\lambda_{\upsilon}}
		\MM{|\beta|,\widetilde M,\mu_{\upsilon},\tilde\mu_{\upsilon}}.
		\label{eq:operator-product-Calpha-loss-upsilon}
	\end{align}
	
	Fix $\mC_f\geq0$,
	$1\leq\lambda_f\leq\tilde\lambda_f$, and
	$1\leq\mu_f\leq\tilde\mu_f$.
	Assume that, whenever $|\theta|+|\beta|\leq N_*$, $f$ satisfies
	\begin{equation}
		\left\|
		\left(
		\prod_{i=1}^k \nabla^{\theta_i}B^{\beta_i}
		\right)f
		\right\|_{\alpha,\mcI}
		\lesssim
		\mC_f
		\La{|\theta|,\widetilde N,\lambda_f,\tilde\lambda_f}
		\MM{|\theta|,\widetilde N,\lambda_f,\tilde\lambda_f}
		\MM{|\beta|,\widetilde M,\mu_f,\tilde\mu_f}.
		\label{eq:operator-product-Calpha-loss-f}
	\end{equation}
	
	Define
	$$
	\lambda:=\max\{\lambda_f,\lambda_{\upsilon}\},
	\qquad
	\tilde\lambda:=\max\{\tilde\lambda_f,\tilde\lambda_{\upsilon}\},
	\qquad
	\mu:=\max\{\mu_f,\mu_{\upsilon}\},
	\qquad
	\tilde\mu:=\max\{\tilde\mu_f,\tilde\mu_{\upsilon}\},
	$$
	and set
	$$
	\mu_*:=\max\{\mu,\mC_{\upsilon}\tilde\lambda\},
	\qquad
	\tilde\mu_*:=\max\{\tilde\mu,\mC_{\upsilon}\tilde\lambda\}.
	$$
	Then, for every integer $N\geq0$ such that
	$N+|\theta|+|\beta|\leq N_*$, one has
	\begin{align}
		\left\|
		\left(
		\prod_{i=1}^k A^{\theta_i}B^{\beta_i}
		\right)f
		\right\|_{N+\alpha,\mcI}
		&\lesssim
		\mC_f\mC_{\upsilon}^{|\theta|}
		\Lat{N+|\theta|,\widetilde N,\lambda,\tilde\lambda}
		\MM{N+|\theta|,\widetilde N,\lambda,\tilde\lambda}
		\MM{|\beta|,\widetilde M,\mu,\tilde\mu}
		\label{eq:operator-product-Calpha-loss-main}
		\\
		&\lesssim
		\mC_f
		\Lat{N+|\theta|,\widetilde N,\lambda,\tilde\lambda}
		\MM{N,\widetilde N,\lambda,\tilde\lambda}
		(\mC_{\upsilon}\tilde\lambda)^{|\theta|}
		\MM{|\beta|,\widetilde M,\mu,\tilde\mu}
		\notag
		\\
		&\lesssim
		\mC_f
		\Lat{N+|\theta|,\widetilde N,\lambda,\tilde\lambda}
		\MM{N,\widetilde N,\lambda,\tilde\lambda}
		\MM{|\theta|+|\beta|,\widetilde M,\mu_*,\tilde\mu_*}.
		\label{eq:operator-product-Calpha-loss-absorbed}
	\end{align}
	
	Consequently, for all integers $N,M\geq0$ with $N+M\leq N_*$,
	\begin{equation}
		\left\|
		(A+B)^M f
		\right\|_{N+\alpha,\mcI}
		\lesssim
		\mC_f
		\Lat{N+M,\widetilde N,\lambda,\tilde\lambda}
		\MM{N,\widetilde N,\lambda,\tilde\lambda}
		\MM{M,\widetilde M,\mu_*,\tilde\mu_*}.
		\label{eq:operator-product-Calpha-loss-A-plus-B}
	\end{equation}

\end{lem}

\begin{proof}
	We use \eqref{eq:operator-product-full-expansion} with the notation from the proof of
	Lemma~\ref{lem:operator-product-estimates}. It suffices to
	estimate
	$
	\left\|
	\nabla^r
	\left(
	\prod_{i=1}^k A^{\theta_i}B^{\beta_i}
	\right)f
	\right\|_{\alpha,\mcI}
	$
	for $0\leq r\leq N$.
	We first estimate the coefficient factors in $C^0$ and $C^\alpha$.
	For the $i$-th coefficient factor in
	\eqref{eq:operator-product-full-expansion}, assumptions
	\eqref{eq:operator-product-Calpha-loss-upsilon-C0}--\eqref{eq:operator-product-Calpha-loss-upsilon} give
	$$
	\begin{aligned}
		\left\|
		\left(
		\prod_{\ell=1}^k
		\nabla^{\delta_{i,\ell}'}
		B^{\kappa_{i,\ell}'}
		\right)
		F_{\gamma_i,\theta_i}
		\right\|_{0,\mcI}
		&\lesssim
		\mC_{\upsilon}^{\theta_i}
		\MM{\theta_i-\gamma_i+\delta_i,\widetilde N,\lambda_{\upsilon},\tilde\lambda_{\upsilon}}
		\MM{\kappa_i,\widetilde M,\mu_{\upsilon},\tilde\mu_{\upsilon}},
		\\
		\left\|
		\left(
		\prod_{\ell=1}^k
		\nabla^{\delta_{i,\ell}'}
		B^{\kappa_{i,\ell}'}
		\right)
		F_{\gamma_i,\theta_i}
		\right\|_{\alpha,\mcI}
		&\lesssim
		\mC_{\upsilon}^{\theta_i}
		\La{\theta_i-\gamma_i+\delta_i,\widetilde N,\lambda_{\upsilon},\tilde\lambda_{\upsilon}}
		\MM{\theta_i-\gamma_i+\delta_i,\widetilde N,\lambda_{\upsilon},\tilde\lambda_{\upsilon}}
		\MM{\kappa_i,\widetilde M,\mu_{\upsilon},\tilde\mu_{\upsilon}}.
	\end{aligned}
	$$
	These estimates follow from \eqref{eq:holder-product-finite} applied componentwise to the differentiated
	$\upsilon$-factors. 
	
	We next estimate each term $F\prod_i g_i$ in \eqref{eq:operator-product-full-expansion}, where $F$ denotes
	the differentiated $f$-factor and $g_i$ are the coefficient factors.
	Apply \eqref{eq:holder-product-lower-exponent} with $N=0$ and $\eta=\alpha$:
	$$
	\left\|F\prod_i g_i\right\|_\alpha
	\lesssim\|F\|_\alpha\prod_i\|g_i\|_0
	+\|F\|_\alpha\sum_{i_0}\|g_{i_0}\|_\alpha\prod_{i\neq i_0}\|g_i\|_0.
	$$
	For the first term, \eqref{eq:appendix-derivative-weight-product},
	combined with assumptions \eqref{eq:operator-product-Calpha-loss-upsilon-C0} and \eqref{eq:operator-product-Calpha-loss-f}, gives
	$$
	\begin{aligned}
		\|F\|_\alpha\prod_i\|g_i\|_0
		\lesssim&\mC_f\mC_{\upsilon}^{|\theta|}
		\La{N',\widetilde N,\lambda,\tilde\lambda}
		\MM{N',\widetilde N,\lambda,\tilde\lambda}
		\MM{r+|\theta|-N',\widetilde N,\lambda,\tilde\lambda}
		\MM{M',\widetilde M,\mu,\tilde\mu}
		\MM{|\beta|-M',\widetilde M,\mu,\tilde\mu}\\
		\leq&
		\mC_f\mC_{\upsilon}^{|\theta|}
		\La{r+|\theta|,\widetilde N,\lambda,\tilde\lambda}
		\MM{r+|\theta|,\widetilde N,\lambda,\tilde\lambda}
		\MM{|\beta|,\widetilde M,\mu,\tilde\mu}.
	\end{aligned}
	$$
	
	For each $i_0$ in the second term, \eqref{eq:appendix-derivative-weight-product} and
	assumptions \eqref{eq:operator-product-Calpha-loss-upsilon-C0}--\eqref{eq:operator-product-Calpha-loss-f} give
	$$
	\begin{aligned}
		\|F\|_\alpha\|g_{i_0}\|_\alpha\prod_{i\neq i_0}\|g_i\|_0
		\lesssim&\mC_f\mC_{\upsilon}^{|\theta|}
		\La{\theta_{i_0}-\gamma_{i_0}+\delta_{i_0},\widetilde N,\lambda,\tilde\lambda}
		\La{N',\widetilde N,\lambda,\tilde\lambda}
		\MM{r+|\theta|,\widetilde N,\lambda,\tilde\lambda}
		\MM{|\beta|,\widetilde M,\mu,\tilde\mu}\\
		\leq&\mC_f\mC_{\upsilon}^{|\theta|}
		\Lat{r+|\theta|,\widetilde N,\lambda,\tilde\lambda}
		\MM{r+|\theta|,\widetilde N,\lambda,\tilde\lambda}
		\MM{|\beta|,\widetilde M,\mu,\tilde\mu}.
	\end{aligned}
	$$
	Here the H{\"o}lder factors satisfy
	$
		\La{\theta_{i_0}-\gamma_{i_0}+\delta_{i_0},\widetilde N,\lambda,\tilde\lambda}
		\La{N',\widetilde N,\lambda,\tilde\lambda}
		\leq
		\Lat{r+|\theta|,\widetilde N,\lambda,\tilde\lambda}.
	$
	Summing over the terms in \eqref{eq:operator-product-full-expansion}, the components, and
	$0\leq r\leq N$ proves \eqref{eq:operator-product-Calpha-loss-main}.
	The scale comparisons used for \eqref{eq:operator-product-absorbed}
	also give
	\eqref{eq:operator-product-Calpha-loss-absorbed}.

	Finally, for $M\geq1$, expand $(A+B)^M$ as in the proof of
	\eqref{eq:operator-product-A-plus-B} and apply
	\eqref{eq:operator-product-Calpha-loss-absorbed} to each ordered product.
	This proves \eqref{eq:operator-product-Calpha-loss-A-plus-B}. For $M=0$, use
	\eqref{eq:operator-product-Calpha-loss-f}.
\end{proof}

For the $C^N$ estimates of commutators and ordered mixed derivatives,
we adapt the proof of \cite[Lemma~A.14]{BMNV21} (see also \cite[Lemma~A.7]{GKN24}) with $p=\infty$
to the total-order hypotheses below.

\begin{lem}
	\label{lem:holder-cooper-C0}
	Fix $N_\upsilon,M_\upsilon,N_f,M_f\in\mathbb N_0$, $N_*\in\mathbb N$,
	an integer $0\leq M_*\leq N_*$, and a time interval
	$\mcI\subset\mathbb R$. Let $\upsilon:\mcI\times\mathbb T^d\to\mathbb R^d$
	and $f:\mcI\times\mathbb T^d\to\mathbb R$ be smooth. Set
	$D_{t,\upsilon}:=\partial_t+\upsilon\cdot\nabla$.
	Fix $\mC_{\upsilon}\geq0$, $1\leq\lambda_{\upsilon}\leq\tilde\lambda_{\upsilon}$,
	and $1\leq\mu_{\upsilon}\leq\tilde\mu_{\upsilon}$. Assume that, for every
	$M,N\in\mathbb N_0$ satisfying $M<M_*$ and $M+N\leq N_*-1$,
	\begin{equation}
		\|D_{t,\upsilon}^M\nabla\upsilon\|_{N,\mcI}\lesssim\mC_{\upsilon}\MM{N+1,N_\upsilon,\lambda_{\upsilon},\tilde\lambda_{\upsilon}}\MM{M,M_\upsilon,\mu_{\upsilon},\tilde\mu_{\upsilon}}.
		\label{eq:C0-upsilon}
	\end{equation}
	Fix $\mC_f\geq0$, $1\leq\lambda_f\leq\tilde\lambda_f$, and
	$1\leq\mu_f\leq\tilde\mu_f$. Assume that, for every $M,N\in\mathbb N_0$ satisfying $M\leq M_*$ and $1\leq M+N\leq N_*$,
	\begin{equation}
		\|D_{t,\upsilon}^M f\|_{N,\mcI}\lesssim\mC_f\MM{N,N_f,\lambda_f,\tilde\lambda_f}\MM{M,M_f,\mu_f,\tilde\mu_f}.
		\label{eq:C0-f}
	\end{equation}
	
	Set $\widetilde N:=\min\{N_\upsilon,N_f\}$ and $\widetilde M:=\min\{M_\upsilon+1,M_f\}$.
	Define
	$$
	\lambda:=\max\{\lambda_f,\lambda_{\upsilon}\},\qquad \tilde\lambda:=\max\{\tilde\lambda_f,\tilde\lambda_{\upsilon}\},\qquad \mu:=\max\{\mu_f,\mu_{\upsilon}\},\qquad \tilde\mu:=\max\{\tilde\mu_f,\tilde\mu_{\upsilon}\},
	$$
	and
	$$
	\mu_*:=\max\left\{\mu,\mC_{\upsilon}\MM{1,N_\upsilon,\lambda_{\upsilon},\tilde\lambda_{\upsilon}}\right\},\qquad \tilde\mu_*:=\max\left\{\tilde\mu,\mC_{\upsilon}\MM{1,N_\upsilon,\lambda_{\upsilon},\tilde\lambda_{\upsilon}}\right\}.
	$$
	Then, for integers $1\leq M\leq M_*$ and $r,N\geq0$ satisfying $M+N+r\leq N_*$, one has
	\begin{align}
		\|[D_{t,\upsilon}^M,\nabla^r]f\|_{N,\mcI}&\lesssim\mC_f\mC_{\upsilon}\MM{1,N_\upsilon,\lambda_{\upsilon},\tilde\lambda_{\upsilon}}\MM{N+r,\widetilde N,\lambda,\tilde\lambda}\MM{M-1,\min\{M_\upsilon,M_f\},\mu_*,\tilde\mu_*}
		\label{eq:C0-commutator-strong}\\
		&\lesssim\mC_f\MM{N+r,\widetilde N,\lambda,\tilde\lambda}\MM{M,\widetilde M,\mu_*,\tilde\mu_*}.
		\label{eq:C0-commutator}
	\end{align}
	Moreover, for every integer $k\geq1$ and $\theta,\beta\in\mathbb N_0^k$ satisfying $|\beta|\leq M_*$ and $1\leq|\theta|+|\beta|\leq N_*$, one has
	\begin{equation}
		\left\|\left(\prod_{i=1}^k\nabla^{\theta_i}D_{t,\upsilon}^{\beta_i}\right)f\right\|_{0,\mcI}\lesssim\mC_f\MM{|\theta|,\widetilde N,\lambda,\tilde\lambda}\MM{|\beta|,\widetilde M,\mu_*,\tilde\mu_*}.
		\label{eq:C0-iterate}
	\end{equation}
	If $\|f\|_{0,\mcI}\lesssim\mC_f$, \eqref{eq:C0-iterate} also holds when $|\theta|+|\beta|=0$.
	The lemma also holds on $\mathbb R^d$ when the stated norms are finite.
\end{lem}

The proof of Lemma~\ref{lem:holder-cooper-C0} uses the following two algebraic identities, with
$D_{t,\upsilon}$ and $\upsilon$ as above.

\begin{lem}
	\label{lem:Komatsu-material-spatial}
	Fix $M,r\in\mathbb N_0$ and coordinate derivatives
	$\partial_{i_1},\ldots,\partial_{i_r}$. Set
	$\partial^\gamma:=\prod_{j=1}^r\partial_{i_j}$, where $|\gamma|=r$.
	For every coordinate derivative $\partial_i$, define
	$$
	\operatorname{ad}_{D_{t,\upsilon}}^0(\partial_i):=\partial_i,\qquad \operatorname{ad}_{D_{t,\upsilon}}^{s+1}(\partial_i):=\left[D_{t,\upsilon},\operatorname{ad}_{D_{t,\upsilon}}^{s}(\partial_i)\right].
	$$
	Then
	\begin{equation}
		[D_{t,\upsilon}^M,\partial^\gamma]=\sum_{\substack{\nu\in\mathbb N_0^r\\1\leq|\nu|\leq M}}\frac{M!}{\nu!(M-|\nu|)!}\left(\prod_{j=1}^r\operatorname{ad}_{D_{t,\upsilon}}^{\nu_j}(\partial_{i_j})\right)D_{t,\upsilon}^{M-|\nu|}.
		\label{eq:Komatsu-material-spatial}
	\end{equation}
	Here $|\nu|:=\sum_{j=1}^r\nu_j$ and $\nu!:=\prod_{j=1}^r\nu_j!$.
	The product in \eqref{eq:Komatsu-material-spatial} follows the order
	defining $\partial^\gamma$.
	When $r=0$, all empty products are understood to be the identity.
\end{lem}

\begin{proof}
	This follows from Komatsu's formula \cite[Lemma~A.12]{BMNV21}
	by reversing the order of the coordinate derivatives and reindexing the sum.
\end{proof}

\begin{lem}
	\label{lem:iterated-material-spatial}
	For every integer $s\geq1$ and every coordinate derivative $\partial_i$, one has
	\begin{equation}
		\operatorname{ad}_{D_{t,\upsilon}}^{s}(\partial_i)=\sum_{k=1}^sF_{s,k}^{(i)}\cdot\nabla,
		\label{eq:iterated-ad-material-spatial}
	\end{equation}
	where each $F_{s,k}^{(i)}:\mcI\times\mathbb T^d\to\mathbb R^d$ is a column vector field of the form
	\begin{equation}
		F_{s,k}^{(i)}=\sum_{\substack{\zeta\in\mathbb N_0^k\\|\zeta|=s-k}}c_{s,k,\zeta}\left(\prod_{j=1}^{k-1}D_{t,\upsilon}^{\zeta_j}\nabla\upsilon\right)D_{t,\upsilon}^{\zeta_k}\partial_i\upsilon.
		\label{eq:iterated-ad-coefficient}
	\end{equation}
	Here $(\nabla\upsilon)_{jl}=\partial_l\upsilon_j$, and $c_{s,k,\zeta}$ are combinatorial constants. In \eqref{eq:iterated-ad-coefficient}, the matrices are multiplied in the displayed order and act on the column vector to their right.
\end{lem}

\begin{proof}
	We argue by induction on $s$, following \cite[Lemma~A.13]{BMNV21}.
	For every smooth vector field $G$, one has
	\begin{equation}
		[D_{t,\upsilon},G\cdot\nabla]=\left(D_{t,\upsilon}G-(\nabla\upsilon)G\right)\cdot\nabla.
		\label{eq:commutator-vector-field}
	\end{equation}
	This follows from
	$[D_{t,\upsilon},\partial_j]=-(\partial_j\upsilon)\cdot\nabla$.
	
	For $s=1$, take $F_{1,1}^{(i)}=-\partial_i\upsilon$.
	Assume \eqref{eq:iterated-ad-material-spatial}--\eqref{eq:iterated-ad-coefficient}
	hold at order $s$, and set $F_{s,k}^{(i)}:=0$ for $k<1$ or $k>s$.
	By \eqref{eq:commutator-vector-field}, define
	\begin{equation}
		F_{s+1,k}^{(i)}:=D_{t,\upsilon}F_{s,k}^{(i)}-(\nabla\upsilon)F_{s,k-1}^{(i)},\qquad 1\leq k\leq s+1.
		\label{eq:iterated-ad-coefficient-recurrence}
	\end{equation}
	For $1\leq k\leq s$, differentiating a term of
	\eqref{eq:iterated-ad-coefficient} preserves its $k$ factors and increases
	the total material order from $s-k$ to $s+1-k$.
	The first term in \eqref{eq:iterated-ad-coefficient-recurrence} vanishes
	for $k=s+1$. The second vanishes for $k=1$. Otherwise, left
	multiplication by $\nabla\upsilon$ adds a factor of material order zero,
	giving $k$ factors of total material order $s+1-k$.
	Thus the recurrence preserves \eqref{eq:iterated-ad-coefficient} and
	yields \eqref{eq:iterated-ad-material-spatial} at order $s+1$.
\end{proof}

\begin{proof}[Proof of Lemma~\ref{lem:holder-cooper-C0}]
	
	We first prove \eqref{eq:C0-commutator-strong}--\eqref{eq:C0-commutator}. We then apply \eqref{eq:C0-commutator} to prove \eqref{eq:C0-iterate} by induction on $k$.
	The commutator vanishes for $r=0$, so assume $r\geq1$. Fix a spatial multi-index
	$\gamma\in\mathbb N_0^d$ with $|\gamma|=r$, and choose coordinate
	derivatives $\partial_{i_1},\ldots,\partial_{i_r}$ such that
	$\partial^\gamma=\prod_{j=1}^r\partial_{i_j}$.
	We estimate $F_{s,k}^{(i)}$ in
	\eqref{eq:iterated-ad-material-spatial}. By
	\eqref{eq:appendix-derivative-weight-product}, \eqref{eq:C0-upsilon}, and \eqref{eq:iterated-ad-coefficient}, for every
	$1\leq k\leq s\leq M_*$ and $N\geq0$
	satisfying $N+s\leq N_*$ one has
	\begin{equation}
		\|F_{s,k}^{(i)}\|_{N,\mcI}\lesssim\mC_{\upsilon}^k\MM{1,N_\upsilon,\lambda_{\upsilon},\tilde\lambda_{\upsilon}}^k\MM{N,(N_\upsilon-1)_+,\lambda_{\upsilon},\tilde\lambda_{\upsilon}}\MM{s-k,M_\upsilon,\mu_{\upsilon},\tilde\mu_{\upsilon}}.
		\label{eq:iterated-ad-coefficient-C0}
	\end{equation}
	Each term in \eqref{eq:iterated-ad-coefficient} contains $k$ factors involving $\nabla\upsilon$, with total material order $s-k$.
	For each factor involving $\nabla\upsilon$, the spatial weights satisfy
	$$\MM{n+1,N_\upsilon,\lambda_{\upsilon},\tilde\lambda_{\upsilon}}=\MM{1,N_\upsilon,\lambda_{\upsilon},\tilde\lambda_{\upsilon}}\MM{n,(N_\upsilon-1)_+,\lambda_{\upsilon},\tilde\lambda_{\upsilon}}.$$
	For such a factor with material and additional spatial orders $M,n$, we have $M<M_*$ and $M+n\leq N+s-1\leq N_*-1$, as required by \eqref{eq:C0-upsilon}.
	
	We next expand the product of iterated commutators. Fix $\nu\in\mathbb N_0^r$ with $|\nu|\geq1$. Expanding each factor by
	\eqref{eq:iterated-ad-material-spatial} gives
	\begin{equation}
		\prod_{p=1}^r\operatorname{ad}_{D_{t,\upsilon}}^{\nu_p}(\partial_{i_p})=\sum_{j=1}^rG_{j,\nu}\nabla^j,
		\label{eq:product-of-iterated-ad}
	\end{equation}
	where $G_{j,\nu}$ denotes a finite collection of rank-$j$ tensor coefficients, with the contraction against $\nabla^j$ understood componentwise.
	Set
	$
	I_\nu:=\left\{p\in\{1,\ldots,r\}:\nu_p\geq1\right\}.
	$
	By \eqref{eq:iterated-ad-material-spatial}, each component of $G_{j,\nu}$ is a sum of tensor contractions of
	$$
	\prod_{p\in I_\nu}\nabla^{\delta_p}F_{\nu_p,k_p}^{(i_p)},\quad 1\leq k_p\leq\nu_p,\quad \sum_{p\in I_\nu}\delta_p=r-j.
	$$
	The factors with $\nu_p=0$ are the coordinate derivatives $\partial_{i_p}$ and introduce no factor involving $\nabla\upsilon$.
	Set $K:=\sum_{p\in I_\nu}k_p$. Then
	$
	1\leq K\leq|\nu|,\ \sum_{p\in I_\nu}(\nu_p-k_p)=|\nu|-K.
	$
	Estimates~\eqref{eq:appendix-derivative-weight-product} and \eqref{eq:iterated-ad-coefficient-C0} give, for
	$|\nu|\leq M_*$ and $N+r-j+|\nu|\leq N_*$,
	\begin{equation}
		\|G_{j,\nu}\|_{N,\mcI}\lesssim\sum_{K=1}^{|\nu|}\mC_{\upsilon}^K\MM{1,N_\upsilon,\lambda_{\upsilon},\tilde\lambda_{\upsilon}}^K\MM{N+r-j,(N_\upsilon-1)_+,\lambda_{\upsilon},\tilde\lambda_{\upsilon}}\MM{|\nu|-K,M_\upsilon,\mu_{\upsilon},\tilde\mu_{\upsilon}}.
		\label{eq:product-of-iterated-ad-C0}
	\end{equation}
	
	We now estimate $[D_{t,\upsilon}^M,\partial^\gamma]f$ using \eqref{eq:product-of-iterated-ad-C0}.
	Fix $1\leq M\leq M_*$ and $N\geq0$ with $M+N+r\leq N_*$. By \eqref{eq:Komatsu-material-spatial} and \eqref{eq:product-of-iterated-ad}, we have
	$$
	[D_{t,\upsilon}^M,\partial^\gamma]f=\sum_{\substack{\nu\in\mathbb N_0^r\\1\leq|\nu|\leq M}}\frac{M!}{\nu!(M-|\nu|)!}\sum_{j=1}^rG_{j,\nu}\nabla^jD_{t,\upsilon}^{M-|\nu|}f.
	$$
	Estimates~\eqref{eq:C0-f} and \eqref{eq:product-of-iterated-ad-C0} give
	$$
	\begin{aligned}
		\left\|[D_{t,\upsilon}^M,\partial^\gamma]f\right\|_{N,\mcI}
		&\lesssim\sum_{\substack{\nu\in\mathbb N_0^r\\1\leq|\nu|\leq M}}\sum_{j=1}^r\sum_{N_1+N_2\leq N}
		\|G_{j,\nu}\|_{N_1,\mcI}\|D_{t,\upsilon}^{M-|\nu|}f\|_{N_2+j,\mcI}\\
		&\lesssim\mC_f\sum_{\substack{\nu\in\mathbb N_0^r\\1\leq|\nu|\leq M}}\sum_{j=1}^r\sum_{N_1+N_2\leq N}\sum_{K=1}^{|\nu|}
		\mC_{\upsilon}^K\MM{1,N_\upsilon,\lambda_{\upsilon},\tilde\lambda_{\upsilon}}^K\MM{N_1+r-j,(N_\upsilon-1)_+,\lambda_{\upsilon},\tilde\lambda_{\upsilon}}\\
		&\quad\cdot\MM{N_2+j,N_f,\lambda_f,\tilde\lambda_f}\MM{|\nu|-K,M_\upsilon,\mu_{\upsilon},\tilde\mu_{\upsilon}}\MM{M-|\nu|,M_f,\mu_f,\tilde\mu_f}.
	\end{aligned}
	$$
	Since $j\geq1$, the derivatives of $f$ in the sum have positive spatial order $N_2+j$.
	Their total order satisfies $(M-|\nu|)+(N_2+j)\leq M+N+r-1\leq N_*-1$, so \eqref{eq:C0-f} applies.
	For $\widetilde N\geq1$, $\MM{N_2+j,N_f,\lambda_f,\tilde\lambda_f}=\lambda_f\MM{N_2+j-1,N_f-1,\lambda_f,\tilde\lambda_f}$.
	By monotonicity, both spatial thresholds may be replaced by zero when $\widetilde N=0$.
	Estimate~\eqref{eq:appendix-derivative-weight-product} gives
	$$
	\begin{aligned}
		&\MM{1,N_\upsilon,\lambda_{\upsilon},\tilde\lambda_{\upsilon}}^K\MM{N_1+r-j,(N_\upsilon-1)_+,\lambda_{\upsilon},\tilde\lambda_{\upsilon}}\MM{N_2+j,N_f,\lambda_f,\tilde\lambda_f}\\
		&\quad\cdot\MM{|\nu|-K,M_\upsilon,\mu_{\upsilon},\tilde\mu_{\upsilon}}\MM{M-|\nu|,M_f,\mu_f,\tilde\mu_f}\\
		\leq&\MM{1,N_\upsilon,\lambda_{\upsilon},\tilde\lambda_{\upsilon}}^K\MM{N+r,\widetilde N,\lambda,\tilde\lambda}\MM{M-K,\min\{M_\upsilon,M_f\},\mu,\tilde\mu}.
	\end{aligned}
	$$
	Summing over $\nu$, $j$, and $N_1,N_2$ gives
	\begin{equation}
		\left\|[D_{t,\upsilon}^M,\partial^\gamma]f\right\|_{N,\mcI}\lesssim\mC_f\MM{N+r,\widetilde N,\lambda,\tilde\lambda}\sum_{K=1}^M\left(\mC_{\upsilon}\MM{1,N_\upsilon,\lambda_{\upsilon},\tilde\lambda_{\upsilon}}\right)^K\MM{M-K,\min\{M_\upsilon,M_f\},\mu,\tilde\mu}.
		\label{eq:C0-commutator-before-absorption}
	\end{equation}
	
	To estimate the sum over $K$, use the definitions of $\mu_*$ and $\tilde\mu_*$. For every $1\leq K\leq M$,
	$$
	\begin{aligned}
		\left(\mC_{\upsilon}\MM{1,N_\upsilon,\lambda_{\upsilon},\tilde\lambda_{\upsilon}}\right)^K\MM{M-K,\min\{M_\upsilon,M_f\},\mu,\tilde\mu}\leq\mC_{\upsilon}\MM{1,N_\upsilon,\lambda_{\upsilon},\tilde\lambda_{\upsilon}}\MM{M-1,\min\{M_\upsilon,M_f\},\mu_*,\tilde\mu_*}.
	\end{aligned}
	$$
	After summing over $K$ in the range $1\leq K\leq M\leq N_*$, we obtain
	\begin{equation}
		\begin{aligned}
		&\sum_{K=1}^M\left(\mC_{\upsilon}\MM{1,N_\upsilon,\lambda_{\upsilon},\tilde\lambda_{\upsilon}}\right)^K\MM{M-K,\min\{M_\upsilon,M_f\},\mu,\tilde\mu}\\
		&\qquad\lesssim\mC_{\upsilon}\MM{1,N_\upsilon,\lambda_{\upsilon},\tilde\lambda_{\upsilon}}\MM{M-1,\min\{M_\upsilon,M_f\},\mu_*,\tilde\mu_*}.
		\end{aligned}
		\label{eq:C0-material-sum-absorption}
	\end{equation}
	Estimates~\eqref{eq:C0-commutator-before-absorption}--%
	\eqref{eq:C0-material-sum-absorption} prove \eqref{eq:C0-commutator-strong}, uniformly over $|\gamma|=r$.
	Since $\min\{M_\upsilon,M_f\}\geq\widetilde M-1$, \eqref{eq:C0-commutator} follows from
	$$
	\mC_{\upsilon}\MM{1,N_\upsilon,\lambda_{\upsilon},\tilde\lambda_{\upsilon}}\MM{M-1,\min\{M_\upsilon,M_f\},\mu_*,\tilde\mu_*}\lesssim\MM{M,\widetilde M,\mu_*,\tilde\mu_*}.
	$$
	
	We now apply \eqref{eq:C0-commutator} to prove \eqref{eq:C0-iterate} by induction on $k$, as in \cite[Lemma~A.14]{BMNV21}.
	Fix spatial multi-indices $\gamma_i$ with $|\gamma_i|=\theta_i$ and omit the identity operators with $(\theta_i,\beta_i)=(0,0)$, so that $k\leq|\theta|+|\beta|\leq N_*$.
	For $k=1$, use \eqref{eq:C0-f} with $M=\beta_1$ and $N=\theta_1$.
	For $k\geq2$, assume \eqref{eq:C0-iterate} for products of at most $k-1$ mixed derivative operators. Then
	\begin{align}
		\left(\prod_{j=1}^k\partial^{\gamma_j}D_{t,\upsilon}^{\beta_j}\right)f
		=\partial^{\gamma_1+\gamma_2}D_{t,\upsilon}^{\beta_1+\beta_2}\left(\prod_{j=3}^k\partial^{\gamma_j}D_{t,\upsilon}^{\beta_j}\right)f
		+\partial^{\gamma_1}[D_{t,\upsilon}^{\beta_1},\partial^{\gamma_2}]D_{t,\upsilon}^{\beta_2}\left(\prod_{j=3}^k\partial^{\gamma_j}D_{t,\upsilon}^{\beta_j}\right)f.
		\label{eq:C0-iterate-BMNV-decomposition}
	\end{align}
	The first term in \eqref{eq:C0-iterate-BMNV-decomposition} involves at most $k-1$ mixed derivative operators, with total spatial and material orders $|\theta|$ and $|\beta|$.
	It therefore satisfies \eqref{eq:C0-iterate} by induction.
	
	For the second term in \eqref{eq:C0-iterate-BMNV-decomposition}, set
	$
	h:=D_{t,\upsilon}^{\beta_2}\left(\prod_{j=3}^k\partial^{\gamma_j}D_{t,\upsilon}^{\beta_j}\right)f.
	$
	Write $S:=\sum_{j=3}^k\theta_j$ and $T:=\beta_2+\sum_{j=3}^k\beta_j$.
	Each spatial derivative of $D_{t,\upsilon}^Mh$ involves at most $k-1$ mixed derivative operators acting on $f$.
	For $M+T\leq M_*$ and $1\leq M+N+S+T\leq N_*$, use \eqref{eq:C0-iterate} by induction when $M+S+T\geq1$.
	When $M=S=T=0$, we have $h=f$ and use \eqref{eq:C0-f}. In both cases, we obtain
	\begin{equation}
		\|D_{t,\upsilon}^Mh\|_{N,\mcI}\lesssim\mC_f\MM{N+S,\widetilde N,\lambda,\tilde\lambda}\MM{M+T,\widetilde M,\mu_*,\tilde\mu_*}.
		\label{eq:C0-iterate-tail-bound}
	\end{equation}
	If $\beta_1=0$ or $\theta_2=0$, the commutator term vanishes. Otherwise, use the identity
	$$
	\MM{N+s,K,\lambda,\tilde\lambda}=\MM{s,K,\lambda,\tilde\lambda}\MM{N,(K-s)_+,\lambda,\tilde\lambda},\qquad N,s,K\in\mathbb N_0.
	$$
	Use the preceding identity to write \eqref{eq:C0-iterate-tail-bound} as \eqref{eq:C0-f} for $h$, with $\mC_h:=\mC_f\MM{S,\widetilde N,\lambda,\tilde\lambda}\MM{T,\widetilde M,\mu_*,\tilde\mu_*}$.
	The spatial and material thresholds are $(\widetilde N-S)_+$ and $(\widetilde M-T)_+$, respectively.
	The scale parameters for $h$ are $\lambda,\tilde\lambda,\mu_*,\tilde\mu_*$.
	The total and material order limits are $N_*-S-T$ and $\min\{M_*-T,N_*-S-T\}$, respectively.
	The velocity estimate \eqref{eq:C0-upsilon} holds on this reduced range. Since the new thresholds do not exceed $N_\upsilon$ and $M_\upsilon+1$, they are also the common thresholds in \eqref{eq:C0-commutator}.
	The orders required in \eqref{eq:C0-commutator} satisfy
	$
	\beta_1+T=|\beta|\leq M_*,\ \beta_1+\theta_1+\theta_2+S+T=|\theta|+|\beta|\leq N_*.
	$
	Apply \eqref{eq:C0-commutator} to $h$ with $M=\beta_1$, $r=\theta_2$, and $N=\theta_1$ to obtain
	$$
	\begin{aligned}
	\left\|\partial^{\gamma_1}[D_{t,\upsilon}^{\beta_1},\partial^{\gamma_2}]h\right\|_{0,\mcI}
	\lesssim\mC_h\MM{\theta_1+\theta_2,(\widetilde N-S)_+,\lambda,\tilde\lambda}\MM{\beta_1,(\widetilde M-T)_+,\mu_*,\tilde\mu_*}=\mC_f\MM{|\theta|,\widetilde N,\lambda,\tilde\lambda}\MM{|\beta|,\widetilde M,\mu_*,\tilde\mu_*}.
	\end{aligned}
	$$
	This completes the induction for \eqref{eq:C0-iterate}.
\end{proof}

\begin{lem}
	\label{lem:holder-cooper-Calpha}
	Fix $0<\alpha<1$, $N_\upsilon,M_\upsilon,N_f,M_f\in\mathbb N_0$, $N_*\in\mathbb N$,
	an integer $0\leq M_*\leq N_*$, and a time interval $\mcI\subset\mathbb R$.
	Let $\upsilon:\mcI\times\mathbb T^d\to\mathbb R^d$
	and $f:\mcI\times\mathbb T^d\to\mathbb R$ be smooth. Set
	$D_{t,\upsilon}:=\partial_t+\upsilon\cdot\nabla$.
	The constants $\mC_\upsilon,\mC_f$, scale parameters, and $\widetilde N,\widetilde M$ are as in
	Lemma~\ref{lem:holder-cooper-C0}.
	Assume that, for every $M,N\in\mathbb N_0$ satisfying $M<M_*$ and $M+N\leq N_*-1$,
	\begin{align}
		\|D_{t,\upsilon}^M\nabla\upsilon\|_{N,\mcI}
		&\lesssim\mC_{\upsilon}\MM{N+1,N_\upsilon,\lambda_{\upsilon},\tilde\lambda_{\upsilon}}\MM{M,M_\upsilon,\mu_{\upsilon},\tilde\mu_{\upsilon}},
		\label{eq:Calpha-upsilon-CN}
		\\
		\|D_{t,\upsilon}^M\nabla\upsilon\|_{N+\alpha,\mcI}
		&\lesssim\mC_{\upsilon}\La{N+1,N_\upsilon,\lambda_{\upsilon},\tilde\lambda_{\upsilon}}\MM{N+1,N_\upsilon,\lambda_{\upsilon},\tilde\lambda_{\upsilon}}\MM{M,M_\upsilon,\mu_{\upsilon},\tilde\mu_{\upsilon}}.
		\label{eq:Calpha-upsilon}
	\end{align}
	Assume, moreover, that
	\begin{equation}
		\|D_{t,\upsilon}^M f\|_{N+\alpha,\mcI}
		\lesssim\mC_f\La{N,N_f,\lambda_f,\tilde\lambda_f}\MM{N,N_f,\lambda_f,\tilde\lambda_f}\MM{M,M_f,\mu_f,\tilde\mu_f}
		\label{eq:Calpha-f}
	\end{equation}
	for every $M,N\in\mathbb N_0$ satisfying $M\leq M_*$ and $1\leq M+N\leq N_*$.
	
	Set
	$$
	\lambda:=\max\{\lambda_f,\lambda_{\upsilon}\},
	\qquad
	\tilde\lambda:=\max\{\tilde\lambda_f,\tilde\lambda_{\upsilon}\},
	\qquad
	\mu:=\max\{\mu_f,\mu_{\upsilon}\},
	\qquad
	\tilde\mu:=\max\{\tilde\mu_f,\tilde\mu_{\upsilon}\},
	$$
	and
	$$
	\mu_*:=\max\left\{\mu,\mC_{\upsilon}\MM{1,N_\upsilon,\lambda_{\upsilon},\tilde\lambda_{\upsilon}}\right\},
	\qquad
	\tilde\mu_*:=\max\left\{\tilde\mu,\mC_{\upsilon}\MM{1,N_\upsilon,\lambda_{\upsilon},\tilde\lambda_{\upsilon}}\right\}.
	$$
	Then, for integers $1\leq M\leq M_*$ and $r,N\geq0$ satisfying $M+N+r\leq N_*$, one has
	\begin{align}
		\|[D_{t,\upsilon}^M,\nabla^r]f\|_{N+\alpha,\mcI}
		&\lesssim\mC_f\mC_{\upsilon}\MM{1,N_\upsilon,\lambda_{\upsilon},\tilde\lambda_{\upsilon}}\Lat{N+r,\widetilde N,\lambda,\tilde\lambda}\MM{N+r,\widetilde N,\lambda,\tilde\lambda}\nonumber\\
		&\quad\cdot\MM{M-1,\min\{M_\upsilon,M_f\},\mu_*,\tilde\mu_*}
		\label{eq:Calpha-commutator-strong}
		\\
		&\lesssim\mC_f\Lat{N+r,\widetilde N,\lambda,\tilde\lambda}\MM{N+r,\widetilde N,\lambda,\tilde\lambda}\MM{M,\widetilde M,\mu_*,\tilde\mu_*}.
		\label{eq:Calpha-commutator}
	\end{align}
	For $M=0$, the commutator vanishes and \eqref{eq:Calpha-commutator} remains valid.
	
	Moreover, for every integer $k\geq1$ and $\theta,\beta\in\mathbb N_0^k$ satisfying $|\beta|\leq M_*$ and $1\leq|\theta|+|\beta|\leq N_*$, one has
	\begin{equation}
		\left\|\left(\prod_{i=1}^k\nabla^{\theta_i}D_{t,\upsilon}^{\beta_i}\right)f\right\|_{\alpha,\mcI}
		\lesssim\mC_f\Lat{|\theta|,\widetilde N,\lambda,\tilde\lambda}\MM{|\theta|,\widetilde N,\lambda,\tilde\lambda}\MM{|\beta|,\widetilde M,\mu_*,\tilde\mu_*}.
		\label{eq:Calpha-iterate}
	\end{equation}
	
	Suppose, in addition, that
	\begin{equation}
		\|D_{t,\upsilon}^M f\|_{N,\mcI}
		\lesssim\mC_f\MM{N,N_f,\lambda_f,\tilde\lambda_f}\MM{M,M_f,\mu_f,\tilde\mu_f}
		\label{eq:Calpha-f-CN-extra}
	\end{equation}
	for $M,N\in\mathbb N_0$ with $M\leq M_*$ and $1\leq M+N\leq N_*$. Then, for integers
	$1\leq M\leq M_*$ and $r,N\geq0$ with $M+N+r\leq N_*$, the commutator estimates improve to
	\begin{align}
		\|[D_{t,\upsilon}^M,\nabla^r]f\|_{N+\alpha,\mcI}
		&\lesssim\mC_f\mC_{\upsilon}\MM{1,N_\upsilon,\lambda_{\upsilon},\tilde\lambda_{\upsilon}}\La{N+r,\widetilde N,\lambda,\tilde\lambda}\MM{N+r,\widetilde N,\lambda,\tilde\lambda}\nonumber\\
		&\quad\cdot\MM{M-1,\min\{M_\upsilon,M_f\},\mu_*,\tilde\mu_*}
		\label{eq:Calpha-commutator-strong-oneLambda}
		\\
		&\lesssim\mC_f\La{N+r,\widetilde N,\lambda,\tilde\lambda}\MM{N+r,\widetilde N,\lambda,\tilde\lambda}\MM{M,\widetilde M,\mu_*,\tilde\mu_*},
		\label{eq:Calpha-commutator-oneLambda}
	\end{align}
	and, for integers $k\geq1$ and $\theta,\beta\in\mathbb N_0^k$ with
	$|\beta|\leq M_*$ and $1\leq|\theta|+|\beta|\leq N_*$, one also has
	\begin{equation}
		\left\|\left(\prod_{i=1}^k\nabla^{\theta_i}D_{t,\upsilon}^{\beta_i}\right)f\right\|_{\alpha,\mcI}
		\lesssim\mC_f\La{|\theta|,\widetilde N,\lambda,\tilde\lambda}\MM{|\theta|,\widetilde N,\lambda,\tilde\lambda}\MM{|\beta|,\widetilde M,\mu_*,\tilde\mu_*}.
		\label{eq:Calpha-iterate-oneLambda}
	\end{equation}
	If \eqref{eq:Calpha-f} also holds at $M=N=0$, the ordered mixed derivative estimates include $|\theta|+|\beta|=0$.
	The lemma also holds on $\mathbb R^d$ when the stated norms are finite.
\end{lem}

\begin{proof}
	We first assume \eqref{eq:Calpha-f-CN-extra} and prove the commutator and ordered mixed derivative estimates with the factor $\Xi_\alpha$.
	We then remove this assumption by increasing $\mC_f$, which gives the estimates with $\Xi_\alpha^2$.
	We begin with \eqref{eq:Calpha-commutator-strong-oneLambda}--\eqref{eq:Calpha-commutator-oneLambda}.
	The commutator vanishes for $r=0$, so assume $r\geq1$.
	For $N\geq0$, $1\leq j\leq r$, $1\leq|\nu|\leq M_*$, and $N+r-j+|\nu|\leq N_*$, \eqref{eq:holder-product-finite-Cnalpha} and
	\eqref{eq:Calpha-upsilon-CN}--\eqref{eq:Calpha-upsilon} give the
	following $C^{N+\alpha}$ estimate for $G_{j,\nu}$ in \eqref{eq:product-of-iterated-ad}.
	\begin{equation}
		\begin{aligned}
			\|G_{j,\nu}\|_{N+\alpha,\mcI}
			&\lesssim\sum_{K=1}^{|\nu|}\mC_{\upsilon}^K\MM{1,N_\upsilon,\lambda_{\upsilon},\tilde\lambda_{\upsilon}}^K\La{N+r-j+1,N_\upsilon,\lambda_{\upsilon},\tilde\lambda_{\upsilon}}\\
			&\quad\cdot
			\MM{N+r-j,(N_\upsilon-1)_+,\lambda_{\upsilon},\tilde\lambda_{\upsilon}}
			\MM{|\nu|-K,M_\upsilon,\mu_{\upsilon},\tilde\mu_{\upsilon}}.
		\end{aligned}
		\label{eq:product-of-iterated-ad-Calpha}
	\end{equation}
	In the products defining $G_{j,\nu}$, one factor involving $\nabla\upsilon$ is estimated in $C^\alpha$ and the others in $C^0$.
	
	We next estimate the commutator. Fix $1\leq M\leq M_*$ and $N\geq0$ with $M+N+r\leq N_*$. For $|\gamma|=r$, \eqref{eq:holder-product-finite-Cnalpha},
	\eqref{eq:Komatsu-material-spatial}, and \eqref{eq:product-of-iterated-ad} give
	$$
	\begin{aligned}
		\left\|[D_{t,\upsilon}^M,\partial^\gamma]f\right\|_{N+\alpha,\mcI}
		\lesssim
		\sum_{\substack{\nu\in\mathbb N_0^r\\1\leq|\nu|\leq M}}
		\sum_{j=1}^r
		\sum_{N_1+N_2\leq N}
		\Big(
		\|G_{j,\nu}\|_{N_1+\alpha,\mcI}\|D_{t,\upsilon}^{M-|\nu|}f\|_{N_2+j,\mcI}
		+\|G_{j,\nu}\|_{N_1,\mcI}\|D_{t,\upsilon}^{M-|\nu|}f\|_{N_2+j+\alpha,\mcI}
		\Big).
	\end{aligned}
	$$
	For the first term, use \eqref{eq:Calpha-f-CN-extra} and \eqref{eq:product-of-iterated-ad-Calpha}.
	For the second, use \eqref{eq:product-of-iterated-ad-C0} and \eqref{eq:Calpha-f}.
	The H{\"o}lder factors in the first and second terms satisfy, respectively,
	$$
	\begin{aligned}
	\La{N_1+r-j+1,N_\upsilon,\lambda_{\upsilon},\tilde\lambda_{\upsilon}}
	&\leq\La{N+r,\widetilde N,\lambda,\tilde\lambda},\\
	\La{N_2+j,N_f,\lambda_f,\tilde\lambda_f}
	&\leq\La{N+r,\widetilde N,\lambda,\tilde\lambda}.
	\end{aligned}
	$$
	Use the argument for \eqref{eq:C0-commutator-before-absorption}.
	With \eqref{eq:C0-material-sum-absorption}, we obtain \eqref{eq:Calpha-commutator-strong-oneLambda}.
	Comparing the material derivative weights as in the proof of \eqref{eq:C0-commutator} gives \eqref{eq:Calpha-commutator-oneLambda}.

	For \eqref{eq:Calpha-iterate-oneLambda}, we use \eqref{eq:C0-iterate-BMNV-decomposition} and induction on $k$, with $|\gamma_i|=\theta_i$ and identity operators omitted.
	For $k=1$, apply \eqref{eq:Calpha-f} with $M=\beta_1$ and $N=\theta_1$.
	For $k\geq2$, assume \eqref{eq:Calpha-iterate-oneLambda} for products of at most $k-1$ mixed derivative operators.
	The first term in \eqref{eq:C0-iterate-BMNV-decomposition} satisfies \eqref{eq:Calpha-iterate-oneLambda} by induction.

	If $\beta_1=0$ or $\theta_2=0$, the commutator term vanishes.
	Otherwise, use $h$, $S$, and $T$ from the proof of \eqref{eq:C0-iterate-tail-bound}.
	For $M+T\leq M_*$ and $1\leq M+N+S+T\leq N_*$, \eqref{eq:C0-iterate-tail-bound} gives the $C^N$ estimate for $D_{t,\upsilon}^Mh$.
	The H{\"o}lder estimate follows from \eqref{eq:Calpha-iterate-oneLambda} by induction when $M+S+T\geq1$, and from \eqref{eq:Calpha-f} when $M=S=T=0$. Thus
	\begin{equation}\label{eq:ordered-product-commutator-input-estimates}
	\begin{aligned}
		\|D_{t,\upsilon}^Mh\|_{N,\mcI}
		&\lesssim\mC_f\MM{N+S,\widetilde N,\lambda,\tilde\lambda}\MM{M+T,\widetilde M,\mu_*,\tilde\mu_*},
		\\
		\|D_{t,\upsilon}^Mh\|_{N+\alpha,\mcI}
		&\lesssim\mC_f\La{N+S,\widetilde N,\lambda,\tilde\lambda}\MM{N+S,\widetilde N,\lambda,\tilde\lambda}\MM{M+T,\widetilde M,\mu_*,\tilde\mu_*}.
	\end{aligned}
	\end{equation}

	Use the amplitude, thresholds, and scales specified after
	\eqref{eq:C0-iterate-tail-bound} for $h$. The required derivative orders satisfy
	$$
		\beta_1+T=|\beta|\leq M_*,\quad
		\beta_1+\theta_1+\theta_2+S+T=|\theta|+|\beta|\leq N_*.
	$$
	The H{\"o}lder factor for $h$ has spatial threshold $(\widetilde N-S)_+$, since $\La{N+S,\widetilde N,\lambda,\tilde\lambda}=\La{N,(\widetilde N-S)_+,\lambda,\tilde\lambda}$.
	The velocity estimates \eqref{eq:Calpha-upsilon-CN}--\eqref{eq:Calpha-upsilon} hold on this reduced range.
	Using \eqref{eq:ordered-product-commutator-input-estimates}, apply \eqref{eq:Calpha-commutator-oneLambda} to $h$ with $M=\beta_1$, $r=\theta_2$, and $N=\theta_1$ to obtain
	$$
	\begin{aligned}
		\left\|\partial^{\gamma_1}[D_{t,\upsilon}^{\beta_1},\partial^{\gamma_2}]h\right\|_{\alpha,\mcI}
		&\lesssim\mC_f\La{|\theta|,\widetilde N,\lambda,\tilde\lambda}
		\MM{|\theta|,\widetilde N,\lambda,\tilde\lambda}
		\MM{|\beta|,\widetilde M,\mu_*,\tilde\mu_*}.
	\end{aligned}
	$$
	This completes the induction for \eqref{eq:Calpha-iterate-oneLambda}.

	We now prove \eqref{eq:Calpha-commutator-strong}--\eqref{eq:Calpha-iterate} without assuming \eqref{eq:Calpha-f-CN-extra}.
	The scales and implicit constants in \eqref{eq:Calpha-commutator-strong-oneLambda}--%
	\eqref{eq:Calpha-iterate-oneLambda} are independent of $\mC_f$.
	For \eqref{eq:Calpha-commutator-strong}--\eqref{eq:Calpha-commutator}, fix $1\leq M\leq M_*$ and $N,r\geq0$ with $M+N+r\leq N_*$. Set $\mC_f':=\mC_f\La{N+r,\widetilde N,\lambda,\tilde\lambda}$.
	The derivatives of $f$ in the commutator sum have spatial orders $N_2+j\leq N+r$.
	Let $m,n\in\mathbb N_0$ satisfy $m\leq M_*$, $1\leq m+n\leq N_*$, and $n\leq N+r$.
	By \eqref{eq:Calpha-f}, we have
	$$
	\begin{aligned}
	\|D_{t,\upsilon}^m f\|_{n,\mcI}
	&\lesssim\mC_f\La{n,N_f,\lambda_f,\tilde\lambda_f}\MM{n,N_f,\lambda_f,\tilde\lambda_f}\MM{m,M_f,\mu_f,\tilde\mu_f}\\
	&\leq\mC_f'\MM{n,N_f,\lambda_f,\tilde\lambda_f}\MM{m,M_f,\mu_f,\tilde\mu_f}.
	\end{aligned}
	$$
	Repeat the proof of \eqref{eq:Calpha-commutator-strong-oneLambda}--\eqref{eq:Calpha-commutator-oneLambda}
	with amplitude $\mC_f'$, using this $C^n$ estimate and \eqref{eq:Calpha-f}.
	Substituting the definition of $\mC_f'$ gives \eqref{eq:Calpha-commutator-strong}--\eqref{eq:Calpha-commutator}.

	For \eqref{eq:Calpha-iterate} with $|\theta|\geq1$, set $\mC_f':=\mC_f\La{|\theta|,\widetilde N,\lambda,\tilde\lambda}$.
	For the commutator term in \eqref{eq:C0-iterate-BMNV-decomposition},
	\eqref{eq:Calpha-commutator-oneLambda} requires $C^N$ estimates for the material derivatives of $h$.
	The corresponding spatial orders of $f$ are at most $N+S$, and $0\leq N\leq\theta_1+\theta_2,
	\ N+S\leq\theta_1+\theta_2+S\leq|\theta|$. Thus
	$$
		\mC_f\La{N+S,N_f,\lambda_f,\tilde\lambda_f}
		\leq\mC_f\La{|\theta|,\widetilde N,\lambda,\tilde\lambda}
		=\mC_f'.
	$$
	The required integer-order estimates for $f$ follow from \eqref{eq:Calpha-f} with amplitude $\mC_f'$.
	We repeat the inductions for \eqref{eq:C0-iterate} and
	\eqref{eq:Calpha-iterate-oneLambda} for products of total spatial order
	at most $|\theta|$, using this fixed amplitude. We obtain
	$$
	\begin{aligned}
		\left\|\left(\prod_{i=1}^k\nabla^{\theta_i}D_{t,\upsilon}^{\beta_i}\right)f\right\|_{\alpha,\mcI}
		&\lesssim\mC_f'\La{|\theta|,\widetilde N,\lambda,\tilde\lambda}\MM{|\theta|,\widetilde N,\lambda,\tilde\lambda}\MM{|\beta|,\widetilde M,\mu_*,\tilde\mu_*}\\
		&=\mC_f\Lat{|\theta|,\widetilde N,\lambda,\tilde\lambda}\MM{|\theta|,\widetilde N,\lambda,\tilde\lambda}\MM{|\beta|,\widetilde M,\mu_*,\tilde\mu_*}.
	\end{aligned}
	$$
	This proves \eqref{eq:Calpha-iterate}.
	At $|\theta|=0$, use \eqref{eq:Calpha-f} directly.
\end{proof}

The following lemma uses \eqref{eq:C0-iterate} and \eqref{eq:Calpha-iterate-oneLambda} to estimate iterates of the material derivative after a change of velocity.
It is a variant of \cite[Lemma~A.2]{GKN26Strong} with total-order hypotheses.

\begin{lem}
	\label{lem:total-order-upgrade-material}
	Fix $\widetilde N,\widetilde M\in\mathbb N_0$, $N_0\in\mathbb N$, and use the notation and parameter restrictions of
	Lemma~\ref{lem:holder-cooper-C0} with $M_*=N_*=N_0$, $N_\upsilon=N_f=\widetilde N$, and $M_\upsilon=M_f=\widetilde M$. Assume that $\upsilon$ and $f$
	satisfy \eqref{eq:C0-upsilon} and \eqref{eq:C0-f} for all $M,N\in\mathbb N_0$
	with $M+N\leq N_0$. Let $\omega:\mcI\times\mathbb T^d\to\mathbb R^d$
	be a smooth vector field, and write
	$$
	D_{t,\upsilon+\omega}=D_{t,\upsilon}+\omega\cdot\nabla
	=\partial_t+(\upsilon+\omega)\cdot\nabla.
	$$
	Assume that
	\begin{equation}
		\left\|\left(\prod_{i=1}^k\nabla^{\theta_i}D_{t,\upsilon}^{\beta_i}\right)\omega\right\|_{0,\mcI}
		\lesssim\mC_\omega\MM{|\theta|,\widetilde N,\lambda_\omega,\tilde\lambda_\omega}\MM{|\beta|,\widetilde M,\mu_\omega,\tilde\mu_\omega}
		\label{eq:C0-w-total-order}
	\end{equation}
	for every integer $k\geq1$ and $\theta,\beta\in\mathbb N_0^k$ satisfying $|\theta|+|\beta|\leq N_0$, where $\mC_\omega\geq0$, $1\leq\lambda_\omega\leq\tilde\lambda_\omega$, and $1\leq\mu_\omega\leq\tilde\mu_\omega$.
	Set
	$$
	\lambda:=\max\{\lambda_f,\lambda_{\upsilon},\lambda_\omega\},\
	\tilde\lambda:=\max\{\tilde\lambda_f,\tilde\lambda_{\upsilon},\tilde\lambda_\omega\},\ \mu:=\max\{\mu_f,\mu_{\upsilon},\mu_\omega,\mC_{\upsilon}\tilde\lambda_{\upsilon},\mC_\omega\tilde\lambda\},\
	\tilde\mu:=\max\{\tilde\mu_f,\tilde\mu_{\upsilon},\tilde\mu_\omega,\mC_{\upsilon}\tilde\lambda_{\upsilon},\mC_\omega\tilde\lambda\}.
	$$
	Then, for all $N,M\in\mathbb N_0$ satisfying $N+M\leq N_0$,
	\begin{equation}
		\|D_{t,\upsilon+\omega}^M f\|_{N,\mcI}
		\lesssim\mC_f\MM{N,\widetilde N,\lambda,\tilde\lambda}\MM{M,\widetilde M,\mu,\tilde\mu}.
		\label{eq:C0-upgrade-material-total-order}
	\end{equation}

	For an integer $1\leq J\leq N_0$, \eqref{eq:C0-upgrade-material-total-order} also holds when
	$M\leq J$ and $1\leq M+N\leq N_0$ under the following assumptions on $\nabla\upsilon$, $f$, and $\omega$.
	It is enough to assume \eqref{eq:C0-upsilon} for $M<J$ and
	$M+N\leq N_0-2$, and \eqref{eq:C0-f} for $M\leq J$ and
	$1\leq M+N\leq N_0$. In place of \eqref{eq:C0-w-total-order}, it is enough to assume
	$$
	\nrm{D_{t,\upsilon}^M\omega}_{N,\mcI}\lesssim\mC_\omega\MM{N,\widetilde N,\lambda_\omega,\tilde\lambda_\omega}\MM{M,\widetilde M,\mu_\omega,\tilde\mu_\omega}
	$$
	for $M<J$ and $M+N\leq N_0-1$.
	The condition on $\nabla\upsilon$ is imposed only for $N_0\geq2$.
	
	For the H{\"o}lder estimate, assume \eqref{eq:C0-upsilon}--\eqref{eq:C0-f} for $M+N\leq N_0$. Fix $0<\alpha<1$. Suppose that
	\eqref{eq:Calpha-upsilon-CN}--%
	\eqref{eq:Calpha-f} and \eqref{eq:Calpha-f-CN-extra} hold whenever
	$N+M\leq N_0$. Assume, in addition to
	\eqref{eq:C0-w-total-order}, that
	\begin{equation}
		\left\|\left(\prod_{i=1}^k\nabla^{\theta_i}D_{t,\upsilon}^{\beta_i}\right)\omega\right\|_{\alpha,\mcI}
		\lesssim\mC_\omega\La{|\theta|,\widetilde N,\lambda_\omega,\tilde\lambda_\omega}\MM{|\theta|,\widetilde N,\lambda_\omega,\tilde\lambda_\omega}\MM{|\beta|,\widetilde M,\mu_\omega,\tilde\mu_\omega}
		\label{eq:Calpha-w-total-order}
	\end{equation}
	whenever $|\theta|+|\beta|\leq N_0$. Then, for all $N,M\in\mathbb N_0$ satisfying $N+M\leq N_0$,
	\begin{equation}
		\|D_{t,\upsilon+\omega}^M f\|_{N+\alpha,\mcI}
		\lesssim\mC_f\Lat{N+M,\widetilde N,\lambda,\tilde\lambda}\MM{N,\widetilde N,\lambda,\tilde\lambda}\MM{M,\widetilde M,\mu,\tilde\mu}.
		\label{eq:Calpha-upgrade-material-total-order}
	\end{equation}
	\end{lem}

\begin{proof}
	We first prove \eqref{eq:C0-upgrade-material-total-order} assuming \eqref{eq:C0-w-total-order}.
	Set $A:=\omega\cdot\nabla$ and $B:=D_{t,\upsilon}$. With \eqref{eq:C0-upsilon}--\eqref{eq:C0-f} for $M+N\leq N_0$,
	Lemma~\ref{lem:holder-cooper-C0} gives \eqref{eq:C0-iterate} for $f$, including order zero.
	Together with \eqref{eq:C0-w-total-order}, this verifies the hypotheses of
	Lemma~\ref{lem:operator-product-estimates}. Applying
	\eqref{eq:operator-product-A-plus-B} proves \eqref{eq:C0-upgrade-material-total-order}.

	We next prove \eqref{eq:C0-upgrade-material-total-order} for $M\leq J$ and $1\leq M+N\leq N_0$, using the stated estimates for $\nabla\upsilon$, $f$, and $\omega$. The case $M=0$ follows from \eqref{eq:C0-f}.
	For $M\geq1$, we first use the proof of \eqref{eq:C0-iterate} to estimate the mixed derivatives of $f$
	with $|\beta|\leq J$ and $1\leq|\theta|+|\beta|\leq N_0$.
	Use \eqref{eq:Komatsu-material-spatial} and \eqref{eq:iterated-ad-coefficient} to expand the commutator terms.
	Each nonzero term has total material order at most $|\beta|-1$ and contains a spatial derivative of $f$.
	For each factor involving $\nabla\upsilon$, its material and additional spatial orders $m,n$ satisfy
	$
		m\leq|\beta|-1\leq J-1,\
		m+n\leq|\theta|+|\beta|-2\leq N_0-2.
	$
	We next estimate the mixed derivatives of $\omega$. When $N_0\geq2$, apply Lemma~\ref{lem:holder-cooper-C0} with $f$ replaced by $\omega$ and
	$
	N_*=N_0-1,\ M_*=J-1.
	$
	We obtain \eqref{eq:C0-iterate} for $\omega$ when $|\beta|\leq J-1$ and $|\theta|+|\beta|\leq N_0-1$.
	When $N_0=1$, only the assumed $C^0$ estimate for $\omega$ is needed.

	We now estimate an ordered product with $|\theta|+|\beta|=M$. Each factor involving $\omega$ in
	\eqref{eq:operator-product-full-expansion} has material and spatial orders $m,n$ satisfying
	$
		m\leq M-1\leq J-1,\
		m+n\leq M+N-1\leq N_0-1.
	$
	The factor involving $f$ has total order at most $M+N$.
	It has positive spatial order when $|\theta|\geq1$. When $|\theta|=0$, the product is $B^Mf$.
	The proof of \eqref{eq:operator-product-absorbed} therefore gives
	$$
	\begin{aligned}
		\left\|\left(\prod_{i=1}^kA^{\theta_i}B^{\beta_i}\right)f\right\|_{N,\mcI}
		\lesssim\mC_f\MM{N,\widetilde N,\lambda,\tilde\lambda}
		(\mC_\omega\tilde\lambda)^{|\theta|}\MM{|\beta|,\widetilde M,\mu,\tilde\mu}
		\leq\mC_f\MM{N,\widetilde N,\lambda,\tilde\lambda}\MM{M,\widetilde M,\mu,\tilde\mu}.
	\end{aligned}
	$$
	Summing over the ordered products in $(A+B)^M$ proves \eqref{eq:C0-upgrade-material-total-order}.

	To prove \eqref{eq:Calpha-upgrade-material-total-order}, apply Lemma~\ref{lem:holder-cooper-Calpha} to obtain \eqref{eq:Calpha-iterate-oneLambda}
	for $f$, including order zero. Together with \eqref{eq:C0-w-total-order} and \eqref{eq:Calpha-w-total-order},
	this verifies the hypotheses of Lemma~\ref{lem:operator-product-estimates-Calpha-loss}.
	Its conclusion \eqref{eq:operator-product-Calpha-loss-A-plus-B} gives \eqref{eq:Calpha-upgrade-material-total-order}.
\end{proof}

\begin{lem}
	\label{lem:velocity-to-commuted-gradient}
	Fix $\widetilde N,\widetilde M\in\mathbb N_0$, $N_*\in\mathbb N$, and a time interval
	$\mcI\subset\mathbb R$. Let $\upsilon:\mcI\times\mathbb T^d\to\mathbb R^d$
	be smooth and set $D_{t,\upsilon}:=\partial_t+\upsilon\cdot\nabla$.
	Fix $\mC_\upsilon\geq0$, $1\leq\lambda_\upsilon\leq\tilde\lambda_\upsilon$,
	and $1\leq\mu_\upsilon\leq\tilde\mu_\upsilon$. Suppose that
	\begin{equation}
		\|D_{t,\upsilon}^M\upsilon\|_{N,\mcI}
		\lesssim\mC_\upsilon\MM{N,\widetilde N,\lambda_\upsilon,\tilde\lambda_\upsilon}\MM{M,\widetilde M,\mu_\upsilon,\tilde\mu_\upsilon}
		\label{eq:velocity-before-commuted-gradient}
	\end{equation}
	for every $M,N\in\mathbb N_0$ satisfying $M+N\leq N_*+1$. Set
	$$
	\mu_{\upsilon,*}:=\max\{\mu_\upsilon,\mC_\upsilon\tilde\lambda_\upsilon\},
	\qquad
	\tilde\mu_{\upsilon,*}:=\max\{\tilde\mu_\upsilon,\mC_\upsilon\tilde\lambda_\upsilon\}.
	$$
	Then, for every $M,N\in\mathbb N_0$ satisfying $M+N\leq N_*$,
	\begin{equation}
		\|D_{t,\upsilon}^M\nabla\upsilon\|_{N,\mcI}
		+
		\|\nabla D_{t,\upsilon}^M\upsilon\|_{N,\mcI}
		\lesssim\mC_\upsilon\MM{N+1,\widetilde N,\lambda_\upsilon,\tilde\lambda_\upsilon}\MM{M,\widetilde M,\mu_{\upsilon,*},\tilde\mu_{\upsilon,*}}.
		\label{eq:velocity-derived-commuted-gradient}
	\end{equation}
	Moreover, for every integer $k\geq1$ and $\theta,\beta\in\mathbb N_0^k$ satisfying $|\theta|+|\beta|\leq N_*$,
	\begin{equation}
		\left\|\left(\prod_{i=1}^k\nabla^{\theta_i}D_{t,\upsilon}^{\beta_i}\right)\upsilon\right\|_{0,\mcI}
		\lesssim\mC_\upsilon\MM{|\theta|,\widetilde N,\lambda_\upsilon,\tilde\lambda_\upsilon}\MM{|\beta|,\widetilde M,\mu_{\upsilon,*},\tilde\mu_{\upsilon,*}}.
		\label{eq:velocity-derived-ordered-mixed}
	\end{equation}

\end{lem}

\begin{proof}
	We first estimate $\nabla D_{t,\upsilon}^M\upsilon$ directly from
	\eqref{eq:velocity-before-commuted-gradient} at spatial order $N+1$.
	We next estimate $D_{t,\upsilon}^M\nabla\upsilon$ by induction on $M$.
	The case $M=0$ follows from \eqref{eq:velocity-before-commuted-gradient} at spatial order $N+1$.

	Fix $M\geq1$ and $N\geq0$ with $M+N\leq N_*$, and assume \eqref{eq:velocity-derived-commuted-gradient}
	at smaller material orders. Apply Lemma~\ref{lem:holder-cooper-C0}
	componentwise with $f=\upsilon$, $r=1$, and total order limit $N_*+1$.
	Choose
	$$
	M_*=M,\quad N_\upsilon=N_f=\widetilde N,\quad M_\upsilon=M_f=\widetilde M,\quad \mC_f=\mC_\upsilon,
	$$
	and use scales $\lambda_\upsilon,\tilde\lambda_\upsilon,\mu_{\upsilon,*},\tilde\mu_{\upsilon,*}$.
	The induction hypothesis gives \eqref{eq:C0-upsilon} for $m<M$ and
	$m+n\leq N_*$. Estimate~\eqref{eq:velocity-before-commuted-gradient}
	implies \eqref{eq:C0-f}.
	Since
	$
	\mC_\upsilon\MM{1,\widetilde N,\lambda_\upsilon,\tilde\lambda_\upsilon}
	\leq\mC_\upsilon\tilde\lambda_\upsilon
	\leq\mu_{\upsilon,*},
	$
	the material scales in \eqref{eq:C0-commutator} are $\mu_{\upsilon,*}$ and $\tilde\mu_{\upsilon,*}$.
	Combining \eqref{eq:C0-commutator} with
	$$D_{t,\upsilon}^M\nabla\upsilon
	=\nabla D_{t,\upsilon}^M\upsilon+[D_{t,\upsilon}^M,\nabla]\upsilon$$
	and the estimate for $\nabla D_{t,\upsilon}^M\upsilon$ proves
	\eqref{eq:velocity-derived-commuted-gradient}.

	Finally, for \eqref{eq:velocity-derived-ordered-mixed}, apply
	Lemma~\ref{lem:holder-cooper-C0} again with the same thresholds and
	scales, taking both the total and material order limits to be $N_*$.
	The hypotheses follow from \eqref{eq:velocity-before-commuted-gradient}--%
	\eqref{eq:velocity-derived-commuted-gradient}.
	Estimate~\eqref{eq:C0-iterate} gives the conclusion at positive total
	order. The case $|\theta|+|\beta|=0$ follows from
	\eqref{eq:velocity-before-commuted-gradient} at $M=N=0$.
\end{proof}

	\section{Mollification estimates}
\label{app:mollification-estimates}

We use spatial mollifiers with finitely many vanishing moments, as in
\cite{BMNV21,Isett17,NV23}. The material derivative estimates use the expansions in
Subsection~\ref{sec:holder-material-commutators}. 

Fix the spatial dimension $d\in\N$ and $N_{\mathrm{mom}}\in\N$. Choose $\eta^{(1)}\in C_c^\infty(\R;\R)$, depending only on $d$ and $N_{\mathrm{mom}}$ and not necessarily nonnegative, such that
$$
\supp\eta^{(1)}\subset\left[-\frac{1}{2\sqrt d},\frac{1}{2\sqrt d}\right],\qquad
\int_{\R}\eta^{(1)}(s)\rd s=1,\qquad
\int_{\R}\eta^{(1)}(s)s^k\rd s=0,\qquad 1\leq k\leq N_{\mathrm{mom}},
$$
and define
$
\eta(x):=\prod_{j=1}^{d}\eta^{(1)}(x_j).
$
Then $\eta\in C_c^\infty(\R^d;\R)$, $\supp\eta\subset B_1(0)$, and
\begin{equation}
	\int_{\R^d}\eta(x)\rd x=1,\qquad
	\int_{\R^d}\eta(x)x^\sgm\rd x=0,\qquad
	1\leq|\sgm|\leq N_{\mathrm{mom}}.
	\label{eq:appendix-mollifier-moments}
\end{equation}
For each $\ell>0$, set
$$
\eta_\ell(x):=\ell^{-d}\eta(x/\ell),\qquad
P_{\leq\ell^{-1}}f:=\eta_\ell*f,\qquad
P_{>\ell^{-1}}f:=f-P_{\leq\ell^{-1}}f.
$$
We identify functions on $\T^d$ with their periodic extensions to
$\R^d$.

\begin{prop}
	\label{prop:appendix-spatial-mollification}
	Fix a time interval $\mcI\subset\R$, $\ell>0$, and $\alpha\in(0,1)$.
	For every smooth scalar, vector, or tensor field $f$ on $\mcI\times\T^d$
	and all integers $N,L\geq0$,
	\begin{equation}
		\begin{aligned}
			\nrm{P_{\leq\ell^{-1}}f}_{N,\mcI}
			&\leq\|\eta\|_{L^1(\mathbb R^d)}\nrm{f}_{N,\mcI},
			&&N\geq0,\\
			\nrm{P_{\leq\ell^{-1}}f}_{N+\alpha,\mcI}
			&\leq\|\eta\|_{L^1(\mathbb R^d)}\nrm{f}_{N+\alpha,\mcI},
			&&N\geq0,\\
			\nrm{P_{>\ell^{-1}}f}_{N,\mcI}
			&\lesssim\ell^L\nrm{f}_{N+L,\mcI},
			&&N\geq0,\quad 0\leq L\leq N_{\mathrm{mom}},\\
			\nrm{P_{>\ell^{-1}}f}_{N+\alpha,\mcI}
			&\lesssim\ell^L\nrm{f}_{N+L+\alpha,\mcI},
			&&N\geq0,\quad 0\leq L\leq N_{\mathrm{mom}}.
		\end{aligned}
		\label{est:appendix-spatial-mollification}
	\end{equation}
	
	For the estimates with two spatial scales, fix $K^*,K^\sharp\in\N$ with
	$K^\sharp>K^*$, $M\in\N_0$ with $M\leq K^*$, and
	$\lambda_F\geq1$ satisfying
	$
	\ell\lambda_F\leq1.
	$
	Fix $\mC_{F,M}>0$, and suppose that $F$ is a smooth scalar, vector, or tensor field on $\mcI\times\T^d$ satisfying
	\begin{equation}
		\nrm{F}_{N,\mcI}
		\lesssim\mC_{F,M}\lambda_F^N,
		\qquad
		0\leq N\leq K^*-M.
		\label{ass:appendix-spatial-mollification-F}
	\end{equation}
	Then, for every $0\leq N\leq K^\sharp$,
	\begin{align}
		\nrm{P_{\leq\ell^{-1}}F}_{N,\mcI}
		&\lesssim\mC_{F,M}\MM{N,K^*-M,\lambda_F,\ell^{-1}},
		\label{est:appendix-spatial-mollification-F-integer}\\
		\nrm{P_{\leq\ell^{-1}}F}_{N+\alpha,\mcI}
		&\lesssim\mC_{F,M}\La{N,K^*-M,\lambda_F,\ell^{-1}}\MM{N,K^*-M,\lambda_F,\ell^{-1}}.
		\label{est:appendix-spatial-mollification-F-holder}
	\end{align}
	If, in addition, $K^*-M\leq N_{\mathrm{mom}}$, then
	\begin{align}
		\nrm{P_{>\ell^{-1}}F}_{N,\mcI}
		&\lesssim\mC_{F,M}\ell^{K^*-M-N}\lambda_F^{K^*-M},
		&&0\leq N\leq K^*-M,
		\label{est:appendix-spatial-mollification-F-high-integer}\\
		\nrm{P_{>\ell^{-1}}F}_{N+\alpha,\mcI}
		&\lesssim\mC_{F,M}\ell^{K^*-M-N-\alpha}\lambda_F^{K^*-M},
		&&0\leq N\leq K^*-M-1.
		\label{est:appendix-spatial-mollification-F-high-holder}
	\end{align}
	The constants in \eqref{est:appendix-spatial-mollification} may depend on $N_{\mathrm{mom}},\alpha$. Those in \eqref{est:appendix-spatial-mollification-F-integer}--\eqref{est:appendix-spatial-mollification-F-high-holder} may also depend on $K^\sharp$ and are independent of
	$F,M,\lambda_F,\ell$, and $\mC_{F,M}$, provided the constant in \eqref{ass:appendix-spatial-mollification-F} is uniform in these quantities.
\end{prop}

\begin{proof}
	Young's inequality applied to spatial derivatives and their
	differences gives the first two estimates in
	\eqref{est:appendix-spatial-mollification}, since $\|\eta_\ell\|_{L^1(\mathbb R^d)}=\|\eta\|_{L^1(\mathbb R^d)}$.
	For $L=0$, the last two estimates follow from the first two and
	$P_{>\ell^{-1}}=I-P_{\leq\ell^{-1}}$.

	Assume $1\leq L\leq N_{\mathrm{mom}}$. Let $N\geq0$ and fix a multi-index $\beta$ with $|\beta|\leq N$. Since spatial derivatives commute with $P_{>\ell^{-1}}$,
	$
	\partial^\beta P_{>\ell^{-1}}f=P_{>\ell^{-1}}\partial^\beta f.
	$
	Using Taylor's theorem and \eqref{eq:appendix-mollifier-moments}, we obtain
	$$
	\partial^\beta P_{>\ell^{-1}}f(x)
	=
	-\sum_{|\sgm|=L}\frac{L}{\sgm!}\int_{\R^d}\eta_\ell(h)(-h)^\sgm\int_0^1(1-\theta)^{L-1}\partial^{\beta+\sgm}f(x-\theta h)\rd\theta\rd h.
	$$
	Since $|h|\lesssim\ell$ on $\supp\eta_\ell$, it follows that
	$$
	\nrm{P_{>\ell^{-1}}f}_{N,\mcI}
	\lesssim
	\ell^L\nrm{f}_{N+L,\mcI}.
	$$
	
	For the H{\"o}lder estimate, let $|\beta|=N$. Taking the difference of
	the integral remainder at $x,y\in\T^d$ gives
	$$
	\frac{\left|\partial^\beta P_{>\ell^{-1}}f(x)-\partial^\beta P_{>\ell^{-1}}f(y)\right|}{|x-y|^\alpha}
	\lesssim
	\ell^L\nrm{f}_{N+L+\alpha,\mcI}.
	$$
	Combining this with the preceding $C^N$ estimate proves the last estimate in \eqref{est:appendix-spatial-mollification}.
	
	We now prove \eqref{est:appendix-spatial-mollification-F-integer}.
	Fix $0\leq N\leq K^\sharp$. Let $\beta$ be a multi-index with $|\beta|\leq N$, and choose a multi-index $\gamma\leq\beta$ such that
	$
	|\gamma|=\min\{|\beta|,K^*-M\}.
	$
	Integration by parts places the remaining spatial derivatives on the kernel and gives
	$
	\partial^\beta P_{\leq\ell^{-1}}F
	=
	\partial^{\beta-\gamma}\eta_\ell*\partial^\gamma F.
	$
	By \eqref{ass:appendix-spatial-mollification-F}, we have
	$$
	\begin{aligned}
		\nrm{\partial^\beta P_{\leq\ell^{-1}}F}_{0,\mcI}
		\lesssim
		\mC_{F,M}\ell^{-(|\beta|-K^*+M)_+}\lambda_F^{\min\{|\beta|,K^*-M\}}=
		\mC_{F,M}\MM{|\beta|,K^*-M,\lambda_F,\ell^{-1}}.
	\end{aligned}
	$$
	We obtain \eqref{est:appendix-spatial-mollification-F-integer}.
	
	To obtain the H{\"o}lder estimate, we first extend \eqref{est:appendix-spatial-mollification-F-integer}
	to spatial order $K^\sharp+1$ by the same argument.
	The definition of $\La{\cdot}$ gives
	$$
	\La{N,K^*-M,\lambda_F,\ell^{-1}}\MM{N,K^*-M,\lambda_F,\ell^{-1}}
	=
	\begin{cases}
		\lambda_F^{N+\alpha}, &0\leq N<K^*-M,\\
		\lambda_F^{K^*-M}\ell^{-(N-K^*+M)-\alpha}, &K^*-M\leq N\leq K^\sharp.
	\end{cases}
	$$
	Interpolating the spatial derivative estimates at orders $N$ and $N+1$ therefore proves
	\eqref{est:appendix-spatial-mollification-F-holder} for
	$0\leq N\leq K^\sharp$.
	
	We finish with the estimates for the mollification difference. Assume that $K^*-M\leq N_{\mathrm{mom}}$. For $0\leq N\leq K^*-M$, apply the third estimate in \eqref{est:appendix-spatial-mollification} with $L=K^*-M-N$. Then
	$$
	\nrm{P_{>\ell^{-1}}F}_{N,\mcI}
	\lesssim
	\ell^{K^*-M-N}\nrm{F}_{K^*-M,\mcI}
	\lesssim
	\mC_{F,M}\ell^{K^*-M-N}\lambda_F^{K^*-M}.
	$$
	Finally, for $0\leq N\leq K^*-M-1$, interpolation between \eqref{est:appendix-spatial-mollification-F-high-integer} at orders $N$ and $N+1$ yields \eqref{est:appendix-spatial-mollification-F-high-holder}.
	\begin{equation*}
		\nrm{P_{>\ell^{-1}}F}_{N+\alpha,\mcI}
		\lesssim
		\nrm{P_{>\ell^{-1}}F}_{N,\mcI}^{1-\alpha}
		\nrm{P_{>\ell^{-1}}F}_{N+1,\mcI}^{\alpha}
		\lesssim
		\mC_{F,M}\ell^{K^*-M-N-\alpha}\lambda_F^{K^*-M}.\qedhere
	\end{equation*}
\end{proof}

We next use the commutator argument in Lemma~\ref{lem:holder-cooper-C0}
to estimate material derivatives of gradients with a spatial threshold
depending on the material order.

\begin{lem}
	\label{lem:appendix-scale-dependent-commuted-gradient}
	Fix a time interval $\mcI\subset\R$ and $K_\upsilon,K_G,K^\sharp,M,N_t\in\N_0$
	with $K_G\leq K_\upsilon$, $M\leq K_\upsilon$, and $M\leq K^\sharp-1$.
	Let $\upsilon:\mcI\times\T^d\to\R^d$ be smooth, and let $G$
	be a smooth scalar, vector, or tensor field on $\mcI\times\T^d$.
	Set $D_{t,\upsilon}:=\partial_t+\upsilon\cdot\nabla$.
	Fix $\mC_\upsilon,\mC_G>0$, $1\leq\lambda\leq\widetilde\lambda$,
	and $1\leq\mu\leq\tilde\mu$. Assume that
	$\mC_\upsilon\lambda\leq\mu$. Suppose that, for $m,N\in\N_0$,
	\begin{equation}
		\nrm{D_{t,\upsilon}^mG}_{N,\mcI}
		\lesssim
		\mC_G\MM{m,N_t,\mu,\tilde\mu}\MM{N,(K_G-m)_+,\lambda,\widetilde\lambda}
		\label{ass:appendix-scale-dependent-uncommuted-field}
	\end{equation}
	whenever $0\leq m\leq M$, $N\geq1$, and $m+N\leq K^\sharp$. Assume further
	that
	\begin{equation}
		\nrm{D_{t,\upsilon}^m\nabla\upsilon}_{N,\mcI}
		\lesssim
		\mC_\upsilon\MM{m,N_t,\mu,\tilde\mu}\MM{N+1,K_\upsilon-m,\lambda,\widetilde\lambda}
		\label{ass:appendix-scale-dependent-lower-gradient}
	\end{equation}
	whenever $0\leq m<M$ and $m+N\leq K^\sharp-1$. Then, for every
	$N\in\N_0$ satisfying $M+N\leq K^\sharp-1$,
	\begin{equation}
		\nrm{D_{t,\upsilon}^M\nabla G}_{N,\mcI}
		+
		\nrm{\nabla D_{t,\upsilon}^MG}_{N,\mcI}
		\lesssim
		\mC_G\MM{M,N_t,\mu,\tilde\mu}\MM{N+1,(K_G-M)_+,\lambda,\widetilde\lambda}.
		\label{est:appendix-scale-dependent-commuted-gradient}
	\end{equation}
	If $G=\upsilon$, $\mC_G=\mC_\upsilon$, and $K_G=K_\upsilon$, assumption
	\eqref{ass:appendix-scale-dependent-lower-gradient} may be omitted. In this
	case, \eqref{est:appendix-scale-dependent-commuted-gradient} holds at every
	admissible material order by induction on $M$.

	The implicit constants may depend on $K^\sharp$ and the constants in \eqref{ass:appendix-scale-dependent-uncommuted-field}--\eqref{ass:appendix-scale-dependent-lower-gradient}.
\end{lem}

\begin{proof}
	The estimate for $\nabla D_{t,\upsilon}^MG$ follows from
	\eqref{ass:appendix-scale-dependent-uncommuted-field} at spatial order
	$N+1$. This also proves the case $M=0$.

	For $M\geq1$, we next estimate the commutator using the proof of \eqref{eq:C0-commutator}
	componentwise with $f=G$ and $r=1$. Take
	$$
	N_*=K^\sharp,\quad M_*=M,\quad N_\upsilon=K_\upsilon-M+1,\quad N_f=(K_G-M)_+,\quad M_\upsilon=M_f=N_t,
	$$
	with amplitudes $\mC_\upsilon,\mC_G$ and scales
	$\lambda,\widetilde\lambda,\mu,\tilde\mu$.
	That proof uses the source estimate only at positive spatial orders.
	The estimates for $D_{t,\upsilon}^mG$ and
	$D_{t,\upsilon}^m\nabla\upsilon$ follow from
	\eqref{ass:appendix-scale-dependent-uncommuted-field} and
	\eqref{ass:appendix-scale-dependent-lower-gradient}, respectively.
	Indeed, $K_\upsilon-m\geq N_\upsilon$ for $m<M$ and
	$(K_G-m)_+\geq N_f$ for $m\leq M$.
	Since $K_G\leq K_\upsilon$, the common spatial threshold is $N_f$.
	The common material threshold is $N_t$.
	Moreover, $N_\upsilon\geq1$ and $\mC_\upsilon\lambda\leq\mu$,
	so the material scales remain $\mu,\tilde\mu$.
	Using $M+N+1\leq K^\sharp$, we obtain
	$$
	\nrm{[D_{t,\upsilon}^M,\nabla]G}_{N,\mcI}
	\lesssim
	\mC_G\MM{M,N_t,\mu,\tilde\mu}
	\MM{N+1,(K_G-M)_+,\lambda,\widetilde\lambda}.
	$$
	The identity
	$D_{t,\upsilon}^M\nabla G
	=\nabla D_{t,\upsilon}^MG+[D_{t,\upsilon}^M,\nabla]G$
	and the estimate for $\nabla D_{t,\upsilon}^MG$ prove
	\eqref{est:appendix-scale-dependent-commuted-gradient}.

	If $G=\upsilon$, $\mC_G=\mC_\upsilon$, and $K_G=K_\upsilon$,
	we proceed by induction on $M$. The case $M=0$ was proved above.
	For $m<M\leq K_\upsilon$, the induction hypothesis gives
	\eqref{ass:appendix-scale-dependent-lower-gradient}, since
	$(K_\upsilon-m)_+=K_\upsilon-m$.
\end{proof}

\begin{prop}
	\label{prop:appendix-material-mollification}
	Fix a time interval $\mcI\subset\R$, $\ell>0$, $K^*\in\N$ with $K^*\geq2$,
	$K_\upsilon,K_f\in\N_0$ with $K_\upsilon\geq\max\{1,K_f\}$, and $\alpha\in(0,1)$. Let
	$\upsilon:\mcI\times\T^d\to\R^d$ be a smooth vector field and $f$
	a smooth scalar, vector, or tensor field. Define
	$
	D_{t,\upsilon}:=\partial_t+\upsilon\cdot\nabla.
	$
	Fix $\mC_\upsilon,\mC_f>0$, $\lambda_\upsilon,\lambda_f\geq1$,
	$1\leq\mu_\upsilon\leq\tilde\mu_\upsilon$, and $1\leq\mu_f\leq\tilde\mu_f$. Set
	$$
		\lambda_*:=\max\{\lambda_\upsilon,\lambda_f\},\ \
		\mu_{\upsilon,*}:=\max\{\mu_\upsilon,\mC_\upsilon\lambda_\upsilon\},\ \
		\tilde\mu_{\upsilon,*}:=\max\{\tilde\mu_\upsilon,\mu_{\upsilon,*}\},\ \
		\mu_*:=\max\{\mu_f,\mu_{\upsilon,*}\},\ \
		\tilde\mu_*:=\max\{\tilde\mu_f,\tilde\mu_{\upsilon,*}\}.
	$$
	Assume that $\ell\lambda_*\leq1$.
	Suppose that, for all $M,N\in\N_0$ in the indicated ranges,
	\begin{equation}
		\begin{aligned}
			\nrm{D_{t,\upsilon}^M\nabla\upsilon}_{N,\mcI}
			+\nrm{\nabla D_{t,\upsilon}^M\upsilon}_{N,\mcI}
			&\lesssim\mC_\upsilon\lambda_\upsilon^{N+1}\MM{M,K_\upsilon-1,\mu_\upsilon,\tilde\mu_\upsilon},
			&&M+N\leq K^*-1,\\
			\nrm{D_{t,\upsilon}^Mf}_{N,\mcI}
			&\lesssim\mC_f\lambda_f^N\MM{M,K_f,\mu_f,\tilde\mu_f},
			&&M+N\leq K^*.
		\end{aligned}
		\label{ass:appendix-material-mollification}
	\end{equation}
	Then, for $M,N\in\N_0$ and every integer $0\leq L\leq N_{\mathrm{mom}}$,
	the following estimates hold in the indicated ranges.
	\begin{equation}
		\begin{aligned}
			\nrm{D_{t,\upsilon}^MP_{\leq\ell^{-1}}f}_{N,\mcI}
			&\lesssim\mC_f\lambda_*^N\MM{M,K_f,\mu_*,\tilde\mu_*},
			&&M+N\leq K^*,\\
			\nrm{D_{t,\upsilon}^MP_{\leq\ell^{-1}}f}_{N+\alpha,\mcI}
			&\lesssim\mC_f\lambda_*^{N+\alpha}\MM{M,K_f,\mu_*,\tilde\mu_*},
			&&M+N\leq K^*-1,\\
			\nrm{D_{t,\upsilon}^MP_{>\ell^{-1}}f}_{N,\mcI}
			&\lesssim\mC_f\ell^L\lambda_*^{N+L}\MM{M,K_f,\mu_*,\tilde\mu_*},
			&&M+N+L\leq K^*,\\
			\nrm{D_{t,\upsilon}^MP_{>\ell^{-1}}f}_{N+\alpha,\mcI}
			&\lesssim\mC_f\ell^L\lambda_*^{N+L+\alpha}\MM{M,K_f,\mu_*,\tilde\mu_*},
			&&M+N+L\leq K^*-1.
		\end{aligned}
		\label{est:appendix-material-mollification}
	\end{equation}
	
	For larger total orders, let $K^\sharp\in\N$ satisfy $K^\sharp>K^*$, and suppose that the velocity field satisfies
	\begin{equation}
		\nrm{D_{t,\upsilon}^M\nabla\upsilon}_{N,\mcI}
		+\nrm{\nabla D_{t,\upsilon}^M\upsilon}_{N,\mcI}
		\lesssim\mC_\upsilon\lambda_\upsilon\MM{N,K^*-1-M,\lambda_\upsilon,\ell^{-1}}\MM{M,K_\upsilon-1,\mu_{\upsilon,*},\tilde\mu_{\upsilon,*}}.
		\label{ass:appendix-extended-velocity}
	\end{equation}
	This estimate is assumed for $0\leq M\leq K^*-1$ and $M+N\leq K^\sharp-1$. Then
	\begin{equation}
		\begin{aligned}
			\nrm{D_{t,\upsilon}^MP_{\leq\ell^{-1}}f}_{N,\mcI}
			&\lesssim\mC_f\MM{N,K^*-M,\lambda_*,\ell^{-1}}\MM{M,K_f,\mu_*,\tilde\mu_*},
			&&0\leq M\leq K^*,\quad M+N\leq K^\sharp,\\
			\nrm{D_{t,\upsilon}^MP_{\leq\ell^{-1}}f}_{N+\alpha,\mcI}
			&\lesssim\mC_f\La{N,K^*-M,\lambda_*,\ell^{-1}}\\
			&\quad\cdot\MM{N,K^*-M,\lambda_*,\ell^{-1}}
			\MM{M,K_f,\mu_*,\tilde\mu_*},
			&&0\leq M\leq K^*,\quad M+N\leq K^\sharp-1.
		\end{aligned}
		\label{est:appendix-extended-spatial-mollification}
	\end{equation}

	The two-scale estimates below use the same fields, amplitudes, and order parameters. They require the following scale and derivative assumptions instead.
	Fix $1\leq\lambda\leq\ell^{-1}$ and
	$1\leq\mu\leq\tilde\mu$ such that
	$\mC_\upsilon\lambda\leq\mu$.
	Suppose that the following estimates hold for $M,N\in\N_0$ in the ranges specified below.
	\begin{equation}\label{ass:appendix-two-scale-material-mollification}
	\begin{aligned}
		\nrm{D_{t,\upsilon}^M\nabla\upsilon}_{N,\mcI}
		+\nrm{\nabla D_{t,\upsilon}^M\upsilon}_{N,\mcI}
		&\lesssim\mC_\upsilon\lambda\MM{M,K_\upsilon-1,\mu,\tilde\mu}
		\MM{N,(K_\upsilon-1-M)_+,\lambda,\ell^{-1}},\\
		\nrm{D_{t,\upsilon}^Mf}_{N,\mcI}
		&\lesssim\mC_f\MM{M,K_f,\mu,\tilde\mu}
		\MM{N,(K_f-M)_+,\lambda,\ell^{-1}}.
	\end{aligned}
	\end{equation}
	The velocity estimates are assumed for $M+N\leq K^*-1$, and the estimate for $f$ for $1\leq M+N\leq K^*$.
	Then, for $1\leq L\leq N_{\mathrm{mom}}$ and $M+N+L\leq K^*$,
	\begin{equation}\label{est:appendix-two-scale-material-mollification}
		\nrm{D_{t,\upsilon}^MP_{>\ell^{-1}}f}_{N,\mcI}
		\lesssim\mC_f\ell^L\MM{M,K_f,\mu,\tilde\mu}
		\MM{N+L,(K_f-M)_+,\lambda,\ell^{-1}}.
	\end{equation}
	For $M+N+L\leq K^*-1$, one also has
	\begin{equation}\label{est:appendix-two-scale-material-mollification-holder}
	\begin{aligned}
		\nrm{D_{t,\upsilon}^MP_{>\ell^{-1}}f}_{N+\alpha,\mcI}
		&\lesssim\mC_f\ell^L\MM{M,K_f,\mu,\tilde\mu}\La{N+L,(K_f-M)_+,\lambda,\ell^{-1}}
		\MM{N+L,(K_f-M)_+,\lambda,\ell^{-1}}.
	\end{aligned}
	\end{equation}
	In this two-scale setting, let $K^\sharp\in\N$ satisfy $K^\sharp>K^*$. Suppose also that $K_f\leq K^*$ and that
	the velocity estimates in \eqref{ass:appendix-two-scale-material-mollification}
	hold for $0\leq M\leq K^*-1$ and $M+N\leq K^\sharp-1$.
	Then, for $0\leq M\leq K^*$,
	\begin{equation}\label{est:appendix-two-scale-extended-spatial-mollification}
		\begin{aligned}
			\nrm{D_{t,\upsilon}^MP_{\leq\ell^{-1}}f}_{N,\mcI}
			&\lesssim\mC_f\MM{N,(K_f-M)_+,\lambda,\ell^{-1}}
			\MM{M,K_f,\mu,\tilde\mu},\\
			\nrm{D_{t,\upsilon}^MP_{\leq\ell^{-1}}f}_{N+\alpha,\mcI}
			&\lesssim\mC_f\La{N,(K_f-M)_+,\lambda,\ell^{-1}}\MM{N,(K_f-M)_+,\lambda,\ell^{-1}}
			\MM{M,K_f,\mu,\tilde\mu},
		\end{aligned}
	\end{equation}
	where $1\leq M+N\leq K^\sharp$ for the $C^N$ estimate and
	$1\leq M+N\leq K^\sharp-1$ for the H{\"o}lder estimate.

	The implicit constants may depend on $d,K^*,K^\sharp,N_{\mathrm{mom}},\alpha$
	and the constants in \eqref{ass:appendix-material-mollification}, \eqref{ass:appendix-extended-velocity}, or \eqref{ass:appendix-two-scale-material-mollification}, as required for each estimate.
\end{prop}

\begin{rem}
	\label{rem:appendix-alternative-velocity-assumptions}
	For \eqref{est:appendix-material-mollification}, the first estimate in
	\eqref{ass:appendix-material-mollification} may be replaced by
	$$
	\nrm{D_{t,\upsilon}^M\upsilon}_{N,\mcI}
	\lesssim
	\mC_\upsilon\lambda_\upsilon^N\MM{M,K_\upsilon-1,\mu_{\upsilon,*},\tilde\mu_{\upsilon,*}},
	\qquad
	N\geq1,\quad M+N\leq K^*.
	$$
	Apply Lemma~\ref{lem:appendix-scale-dependent-commuted-gradient}
	with $G=\upsilon$ and $\mC_G=\mC_\upsilon$.
	Set both spatial thresholds and the total-order limit in the lemma to $K^*$.
	Take $N_t=K_\upsilon-1$, $\lambda=\widetilde\lambda=\lambda_\upsilon$,
	$\mu=\mu_{\upsilon,*}$, and $\tilde\mu=\tilde\mu_{\upsilon,*}$.
	The lemma gives
	$$
	\nrm{D_{t,\upsilon}^M\nabla\upsilon}_{N,\mcI}
	+
	\nrm{\nabla D_{t,\upsilon}^M\upsilon}_{N,\mcI}
	\lesssim
	\mC_\upsilon\lambda_\upsilon^{N+1}\MM{M,K_\upsilon-1,\mu_{\upsilon,*},\tilde\mu_{\upsilon,*}},
	\qquad
	M+N\leq K^*-1.
	$$
	
	Likewise, \eqref{ass:appendix-extended-velocity} may be replaced by
	$$
	\nrm{D_{t,\upsilon}^M\upsilon}_{N,\mcI}
	\lesssim
	\mC_\upsilon\MM{M,K_\upsilon-1,\mu_{\upsilon,*},\tilde\mu_{\upsilon,*}}
	\MM{N,K^*-M,\lambda_\upsilon,\ell^{-1}},
	\qquad
	N\geq1,\quad 0\leq M\leq K^*-1,\quad M+N\leq K^\sharp.
	$$
	Use the same spatial thresholds and material parameters, with
	$\lambda=\lambda_\upsilon$, $\widetilde\lambda=\ell^{-1}$,
	and total-order limit $K^\sharp$.
	The lemma then gives \eqref{ass:appendix-extended-velocity}.
\end{rem}

The assumptions \eqref{ass:appendix-material-mollification}
and spatial interpolation give
\begin{equation}
	\begin{aligned}
		\nrm{D_{t,\upsilon}^M\nabla\upsilon}_{N+\alpha,\mcI}
		+\nrm{\nabla D_{t,\upsilon}^M\upsilon}_{N+\alpha,\mcI}
		&\lesssim_\alpha\mC_\upsilon\lambda_\upsilon^{N+1+\alpha}\MM{M,K_\upsilon-1,\mu_\upsilon,\tilde\mu_\upsilon},
		&&M+N\leq K^*-2,\\
		\nrm{D_{t,\upsilon}^Mf}_{N+\alpha,\mcI}
		&\lesssim_\alpha\mC_f\lambda_f^{N+\alpha}\MM{M,K_f,\mu_f,\tilde\mu_f},
		&&M+N\leq K^*-1.
	\end{aligned}
	\label{est:appendix-material-mollification-interpolated-holder}
\end{equation}
Under \eqref{ass:appendix-two-scale-material-mollification}, spatial interpolation gives
\begin{equation}\label{est:appendix-two-scale-material-mollification-interpolated-holder}
	\begin{aligned}
		\nrm{D_{t,\upsilon}^M\nabla\upsilon}_{N+\alpha,\mcI}
		+\nrm{\nabla D_{t,\upsilon}^M\upsilon}_{N+\alpha,\mcI}
		&\lesssim\mC_\upsilon\lambda
		\MM{M,K_\upsilon-1,\mu,\tilde\mu}\La{N,(K_\upsilon-1-M)_+,\lambda,\ell^{-1}}\\
		&\quad\cdot
		\MM{N,(K_\upsilon-1-M)_+,\lambda,\ell^{-1}},\\
		\nrm{D_{t,\upsilon}^Mf}_{N+\alpha,\mcI}
		&\lesssim\mC_f\MM{M,K_f,\mu,\tilde\mu}\La{N,(K_f-M)_+,\lambda,\ell^{-1}}\\
		&\quad\cdot\MM{N,(K_f-M)_+,\lambda,\ell^{-1}}.
	\end{aligned}
\end{equation}
The velocity estimates hold for $M+N\leq K^*-2$, and the estimate for $f$ for
$1\leq M+N\leq K^*-1$.

We use the following commutator estimates to prove
Proposition~\ref{prop:appendix-material-mollification}.
The single-commutator estimates follow the kernel argument of
\cite[Lemma~A.4]{GK22}. The proof below treats iterated commutators
with compactly supported kernels, using vanishing moments for the higher-order
estimates and spatial differences for the H{\"o}lder estimates.

\begin{lem}
	\label{lem:multiple-mollifier-commutator}
	Fix a time interval $\mcI\subset\R$, $0<\ell\leq1$, and $r\in\N_0$.
	Let $\phi_1,\ldots,\phi_r$ be smooth scalar functions on
	$\mcI\times\T^d$, and let $g$ be a smooth scalar, vector, or tensor
	field on the same domain. Define
	$$
	\mathcal C_0(P_{\leq\ell^{-1}})g:=P_{\leq\ell^{-1}}g,\qquad
	\mathcal C_r(\phi_1,\ldots,\phi_r;P_{\leq\ell^{-1}})g:=[\phi_r,\mathcal C_{r-1}(\phi_1,\ldots,\phi_{r-1};P_{\leq\ell^{-1}})]g
	$$
	for every $r\geq1$. Then
	$$
	\mathcal C_r(\phi_1,\ldots,\phi_r;P_{\leq\ell^{-1}})g(x)=\int_{\R^d}\eta_\ell(h)\prod_{\nu=1}^{r}\bigl(\phi_\nu(x)-\phi_\nu(x-h)\bigr)g(x-h)\rd h.
	$$
	
	For every $N\in\N_0$, every multi-index $\gamma\in\N_0^d$, and every $\alpha\in(0,1)$,
	\begin{equation}
		\begin{aligned}
			\nrm{\mathcal C_r(\phi_1,\ldots,\phi_r;P_{\leq\ell^{-1}})\partial^\gamma g}_{N,\mcI}
			&\lesssim\ell^r\sum_{N_0+\cdots+N_r=N}\nrm{g}_{N_0+|\gamma|,\mcI}\prod_{\nu=1}^{r}\nrm{\nabla\phi_\nu}_{N_\nu,\mcI},\\
			\nrm{\mathcal C_r(\phi_1,\ldots,\phi_r;P_{\leq\ell^{-1}})\partial^\gamma g}_{N+\alpha,\mcI}
			&\lesssim\ell^r\sum_{N_0+\cdots+N_r=N}\nrm{g}_{N_0+|\gamma|+\alpha,\mcI}\prod_{\nu=1}^{r}\nrm{\nabla\phi_\nu}_{N_\nu,\mcI}\\
			&\quad+\ell^r\sum_{s=1}^{r}\sum_{N_0+\cdots+N_r=N}\nrm{g}_{N_0+|\gamma|,\mcI}\nrm{\nabla\phi_s}_{N_s+\alpha,\mcI}\prod_{\nu\neq s}\nrm{\nabla\phi_\nu}_{N_\nu,\mcI}.
		\end{aligned}
		\label{est:multiple-mollifier-commutator}
	\end{equation}

	If, in addition, $L\in\N$ satisfies $1\leq r\leq L\leq N_{\mathrm{mom}}$, then the vanishing moments in \eqref{eq:appendix-mollifier-moments} give the stronger estimates
	\begin{equation}
		\begin{aligned}
			\nrm{\mathcal C_r(\phi_1,\ldots,\phi_r;P_{\leq\ell^{-1}})\partial^\gamma g}_{N,\mcI}
			&\lesssim\ell^L\sum_{N_0+\cdots+N_r=N+L-r}\nrm{g}_{N_0+|\gamma|,\mcI}\prod_{\nu=1}^{r}\nrm{\nabla\phi_\nu}_{N_\nu,\mcI},\\
			\nrm{\mathcal C_r(\phi_1,\ldots,\phi_r;P_{\leq\ell^{-1}})\partial^\gamma g}_{N+\alpha,\mcI}
			&\lesssim\ell^L\sum_{N_0+\cdots+N_r=N+L-r}\nrm{g}_{N_0+|\gamma|+\alpha,\mcI}\prod_{\nu=1}^{r}\nrm{\nabla\phi_\nu}_{N_\nu,\mcI}\\
			&\quad+\ell^L\sum_{s=1}^{r}\sum_{N_0+\cdots+N_r=N+L-r}\nrm{g}_{N_0+|\gamma|,\mcI}\nrm{\nabla\phi_s}_{N_s+\alpha,\mcI}\prod_{\nu\neq s}\nrm{\nabla\phi_\nu}_{N_\nu,\mcI}.
		\end{aligned}
		\label{est:higher-order-multiple-mollifier-commutator}
	\end{equation}

	More generally, fix $N_g\in\{0,\ldots,N\}$. Then
	\begin{equation}
		\begin{aligned}
			\nrm{\mathcal C_r(\phi_1,\ldots,\phi_r;P_{\leq\ell^{-1}})\partial^\gamma g}_{N,\mcI}
			&\lesssim_{N,r}\sum_{\substack{N_0+\cdots+N_r=N_g\\ M_0+\cdots+M_r=N-N_g}}
			\ell^{r-M_0-\#\{\nu\in\{1,\ldots,r\}:M_\nu\geq1\}}
			\nrm{g}_{N_0+|\gamma|,\mcI}
			\prod_{\nu=1}^{r}\nrm{\nabla\phi_\nu}_{N_\nu+(M_\nu-1)_+,\mcI},\\
			\nrm{\mathcal C_r(\phi_1,\ldots,\phi_r;P_{\leq\ell^{-1}})\partial^\gamma g}_{N+\alpha,\mcI}
			&\lesssim_{N,r,\alpha}\sum_{\substack{N_0+\cdots+N_r=N_g\\ M_0+\cdots+M_r=N-N_g}}
			\ell^{r-M_0-\#\{\nu\in\{1,\ldots,r\}:M_\nu\geq1\}}\\
			&\hspace{56pt}\cdot\left(
			\nrm{g}_{N_0+|\gamma|+\alpha,\mcI}
			\prod_{\nu=1}^{r}\nrm{\nabla\phi_\nu}_{N_\nu+(M_\nu-1)_+,\mcI}\right.\\
			&\hspace{60pt}\left.
			+\nrm{g}_{N_0+|\gamma|,\mcI}
			\sum_{s=1}^{r}\nrm{\nabla\phi_s}_{N_s+(M_s-1)_++\alpha,\mcI}
			\prod_{\nu\neq s}\nrm{\nabla\phi_\nu}_{N_\nu+(M_\nu-1)_+,\mcI}
			\right).
		\end{aligned}
		\label{est:general-multiple-mollifier-commutator}
	\end{equation}

	For a smooth vector field $\upsilon$ on $\mcI\times\T^d$, one also has
	\begin{equation}
		\begin{aligned}
			D_{t,\upsilon}\bigl(\mathcal C_r(\phi_1,\ldots,\phi_r;P_{\leq\ell^{-1}})g\bigr)
			=&\mathcal C_r(\phi_1,\ldots,\phi_r;P_{\leq\ell^{-1}})D_{t,\upsilon}g
			+\sum_{s=1}^{r}\mathcal C_r(\phi_1,\ldots,D_{t,\upsilon}\phi_s,\ldots,\phi_r;P_{\leq\ell^{-1}})g\\
			&+\sum_{j=1}^{d}\mathcal C_{r+1}(\phi_1,\ldots,\phi_r,\upsilon_j;P_{\leq\ell^{-1}})\partial_jg\\
			&-\sum_{s=1}^{r}\sum_{j=1}^{d}\mathcal C_r(\phi_1,\ldots,\widehat{\phi_s},\ldots,\phi_r,\upsilon_j;P_{\leq\ell^{-1}})\bigl((\partial_j\phi_s)g\bigr).
		\end{aligned}
		\label{eq:Dt-multiple-mollifier-commutator}
	\end{equation}
	Here $\widehat{\phi_s}$ means that the factor $\phi_s$ is omitted.
\end{lem}

\begin{proof}
	We prove the spatial estimates from the kernel representation, then compute the material derivative.
	The kernel representation follows by induction on $r$ from the
	definition of $\mathcal C_r$.
	For $h\in\R^d$, write $\Delta_h\phi(x):=\phi(x)-\phi(x-h)$. We first prove \eqref{est:general-multiple-mollifier-commutator}. Fix $|\beta|=N$ and choose $\kappa\leq\beta$ with $|\kappa|=N_g$. Applying $\partial^\kappa$ with $h$ fixed gives
	$$
	\partial^\kappa\bigl(\mathcal C_r(\phi_1,\ldots,\phi_r;P_{\leq\ell^{-1}})\partial^\gamma g\bigr)(x)
	=
	\sum_{\kappa_0+\cdots+\kappa_r=\kappa}c_{\kappa_0,\ldots,\kappa_r}\int_{\R^d}\eta_\ell(h)\prod_{\nu=1}^{r}\Delta_h\partial^{\kappa_\nu}\phi_\nu(x)\,\partial^{\kappa_0+\gamma}g(x-h)\rd h.
	$$
	After the change of variables $y=x-h$, apply $\partial^{\beta-\kappa}$ and then return to $h=x-y$. Every resulting term has the form
	$$
	\int_{\R^d}\partial^{\sgm_0}\eta_\ell(h)
	\prod_{\substack{1\leq\nu\leq r\\ \sgm_\nu=0}}\Delta_h\partial^{\kappa_\nu}\phi_\nu(x)
	\prod_{\substack{1\leq\nu\leq r\\ \sgm_\nu\neq0}}\partial^{\kappa_\nu+\sgm_\nu}\phi_\nu(x)\,
	\partial^{\kappa_0+\gamma}g(x-h)\rd h,
	$$
	where
	$
	\kappa_0+\cdots+\kappa_r=\kappa,\ \sgm_0+\cdots+\sgm_r=\beta-\kappa.
	$
	Set $N_\nu:=|\kappa_\nu|$ and $M_\nu:=|\sgm_\nu|$. Then
	$$
	N_0+\cdots+N_r=N_g,\quad M_0+\cdots+M_r=N-N_g.
	$$
	If $M_\nu=0$, then
	$$
	\left|\Delta_h\partial^{\kappa_\nu}\phi_\nu(x)\right|\lesssim|h|\nrm{\nabla\phi_\nu}_{N_\nu,\mcI}.
	$$
	If $M_\nu\geq1$, then
	$$
	\left|\partial^{\kappa_\nu+\sgm_\nu}\phi_\nu(x)\right|\lesssim\nrm{\nabla\phi_\nu}_{N_\nu+M_\nu-1,\mcI}.
	$$
	The scaling of the kernel gives
	$$
	\int_{\R^d}\left|\partial^{\sgm_0}\eta_\ell(h)\right||h|^{r-\#\{\nu\in\{1,\ldots,r\}:M_\nu\geq1\}}\rd h
	\lesssim
	\ell^{r-M_0-\#\{\nu\in\{1,\ldots,r\}:M_\nu\geq1\}}.
	$$
	This proves the $C^0$ estimate for derivatives of order exactly $N$.
	
	For an arbitrary multi-index $\widetilde\beta$ with $|\widetilde\beta|=n<N$, set $N_g':=\min\{N_g,n\}$ and repeat the preceding argument with $\widetilde\beta$ and $N_g'$ in place of $\beta$ and $N_g$, respectively. To estimate this expression by the right-hand side of \eqref{est:general-multiple-mollifier-commutator}, set
	$$
	\begin{aligned}
		N_0':=N_0+N_g-N_g',\quad
		M_0':=M_0+(N-N_g)-(n-N_g'),\quad
		N_0'+\sum_{\nu=1}^rN_\nu=N_g,\quad
		M_0'+\sum_{\nu=1}^rM_\nu=N-N_g.
	\end{aligned}
	$$
	Here $N_0'\geq N_0$ and $M_0'\geq M_0$. Increasing $N_0$ enlarges the norm of $g$, and increasing $M_0$ enlarges the kernel factor because $0<\ell\leq1$.
	Each term is therefore bounded by the right-hand side of \eqref{est:general-multiple-mollifier-commutator}.
	Summing over the derivative orders up to $N$ gives the $C^N$ estimate.
	
	For the H{\"o}lder estimate, use \eqref{eq:holder-product-finite}. For $x,z\in\T^d$, we have
	$$
	\left|\partial^{\kappa_0+\gamma}g(x-h)-\partial^{\kappa_0+\gamma}g(z-h)\right|
	\lesssim
	|x-z|^\alpha\nrm{g}_{N_0+|\gamma|+\alpha,\mcI}.
	$$
	For the difference of the $s$th coefficient, when $M_s=0$, use
	$$
	\left|\Delta_h\partial^{\kappa_s}\phi_s(x)-\Delta_h\partial^{\kappa_s}\phi_s(z)\right|
	\lesssim
	|h||x-z|^\alpha\nrm{\nabla\phi_s}_{N_s+\alpha,\mcI}.
	$$
	If $M_s\geq1$, use
	$$
	\left|\partial^{\kappa_s+\sgm_s}\phi_s(x)-\partial^{\kappa_s+\sgm_s}\phi_s(z)\right|
	\lesssim
	|x-z|^\alpha\nrm{\nabla\phi_s}_{N_s+M_s-1+\alpha,\mcI}.
	$$
	Estimating the remaining factors as in the $C^N$ estimate and integrating in $h$ proves the second estimate in \eqref{est:general-multiple-mollifier-commutator}.
	Taking $N_g=N$ gives $M_0=\cdots=M_r=0$, so \eqref{est:general-multiple-mollifier-commutator} reduces to \eqref{est:multiple-mollifier-commutator}.
	
	We next prove \eqref{est:higher-order-multiple-mollifier-commutator}. The case $L=r$ is contained in \eqref{est:multiple-mollifier-commutator}, so assume that $L>r$. Fix $|\beta|=N$. After applying $\partial^\beta$ with $h$ fixed, consider a term corresponding to $\boldsymbol\kappa=(\kappa_0,\ldots,\kappa_r)$ with
	$
	\kappa_0+\cdots+\kappa_r=\beta.
	$
	For each $1\leq\nu\leq r$, use
	$$
	\Delta_h\partial^{\kappa_\nu}\phi_\nu(x)
	=
	\sum_{j_\nu=1}^{d}h_{j_\nu}\int_0^1\partial_{j_\nu}\partial^{\kappa_\nu}\phi_\nu(x-\theta_\nu h)\rd\theta_\nu.
	$$
	For each $j_1,\ldots,j_r\in\{1,\ldots,d\}$, set
	$$
	B_{j_1,\ldots,j_r}^{\boldsymbol\kappa}(x,h)
	:=
	\left(\prod_{\nu=1}^{r}\int_0^1\partial_{j_\nu}\partial^{\kappa_\nu}\phi_\nu(x-\theta_\nu h)\rd\theta_\nu\right)
	\partial^{\kappa_0+\gamma}g(x-h).
	$$
	The corresponding term is
	$$
	c_{\kappa_0,\ldots,\kappa_r}\sum_{j_1,\ldots,j_r=1}^{d}\int_{\R^d}\eta_\ell(h)h_{j_1}\cdots h_{j_r}B_{j_1,\ldots,j_r}^{\boldsymbol\kappa}(x,h)\rd h.
	$$
	Taylor's theorem in the second variable at $h=0$ gives
	$$
	\begin{aligned}
		B_{j_1,\ldots,j_r}^{\boldsymbol\kappa}(x,h)
		=\sum_{|\zeta|\leq L-r-1}\frac{h^\zeta}{\zeta!}(\partial_h^\zeta B_{j_1,\ldots,j_r}^{\boldsymbol\kappa})(x,0)
		+(L-r)\sum_{|\zeta|=L-r}\frac{h^\zeta}{\zeta!}\int_0^1(1-\vartheta)^{L-r-1}(\partial_h^\zeta B_{j_1,\ldots,j_r}^{\boldsymbol\kappa})(x,\vartheta h)\rd\vartheta.
	\end{aligned}
	$$
	For every $|\zeta|\leq L-r-1$, the monomial $h_{j_1}\cdots h_{j_r}h^\zeta$ has degree between $r$ and $L-1$. Hence, by scaling and \eqref{eq:appendix-mollifier-moments},
	$$
	\int_{\R^d}\eta_\ell(h)h_{j_1}\cdots h_{j_r}h^\zeta\rd h=0.
	$$
	Thus only the integral remainder contributes. For every $|\zeta|=L-r$ and $0\leq\vartheta\leq1$, we have
	$$
	\left|(\partial_h^\zeta B_{j_1,\ldots,j_r}^{\boldsymbol\kappa})(x,\vartheta h)\right|
	\lesssim
	\sum_{L_0+\cdots+L_r=L-r}
	\nrm{g}_{N_0+L_0+|\gamma|,\mcI}
	\prod_{\nu=1}^{r}\nrm{\nabla\phi_\nu}_{N_\nu+L_\nu,\mcI},
	$$
	where $N_\nu:=|\kappa_\nu|$ for $0\leq\nu\leq r$. The factor $h^\zeta$, together with $h_{j_1}\cdots h_{j_r}$, is bounded by $|h|^L$. Summing over the finitely many terms and using
	$
	\int_{\R^d}|\eta_\ell(h)||h|^L\rd h\lesssim\ell^L
	$
	gives the $C^0$ estimate in \eqref{est:higher-order-multiple-mollifier-commutator} for derivatives of order $N$. For $|\beta|=n<N$, the total index is $n+L-r$. Increasing $N_0$ by $N-n$ and using the definition of the $C^k$ norms proves the stated estimate.
	
	For the H{\"o}lder estimate, \eqref{eq:holder-product-finite} gives, for $x,z\in\T^d$,
	$$
	\begin{aligned}
		&\left|(\partial_h^\zeta B_{j_1,\ldots,j_r}^{\boldsymbol\kappa})(x,\vartheta h)
		-(\partial_h^\zeta B_{j_1,\ldots,j_r}^{\boldsymbol\kappa})(z,\vartheta h)\right|\\
		&\qquad\lesssim
		|x-z|^\alpha
		\sum_{L_0+\cdots+L_r=L-r}
		\Bigg(
		\nrm{g}_{N_0+L_0+|\gamma|+\alpha,\mcI}
		\prod_{\nu=1}^{r}\nrm{\nabla\phi_\nu}_{N_\nu+L_\nu,\mcI}
		+\nrm{g}_{N_0+L_0+|\gamma|,\mcI}
		\sum_{s=1}^{r}
		\nrm{\nabla\phi_s}_{N_s+L_s+\alpha,\mcI}
		\prod_{\nu\neq s}\nrm{\nabla\phi_\nu}_{N_\nu+L_\nu,\mcI}
		\Bigg).
	\end{aligned}
	$$
	Integrating against $|\eta_\ell(h)||h|^L$ gives the H{\"o}lder seminorm estimate.
	Together with the first estimate in \eqref{est:higher-order-multiple-mollifier-commutator}, we obtain the second estimate there.
	
	It remains to prove \eqref{eq:Dt-multiple-mollifier-commutator}. With $h$ fixed,
	\begin{align*}
		\partial_t\bigl(\mathcal C_r(\phi_1,\ldots,\phi_r;P_{\leq\ell^{-1}})g\bigr)
		={}&\mathcal C_r(\phi_1,\ldots,\phi_r;P_{\leq\ell^{-1}})\partial_tg
		+\sum_{s=1}^{r}\mathcal C_r(\phi_1,\ldots,\partial_t\phi_s,\ldots,\phi_r;P_{\leq\ell^{-1}})g,\\
		\partial_j\bigl(\mathcal C_r(\phi_1,\ldots,\phi_r;P_{\leq\ell^{-1}})g\bigr)
		={}&\mathcal C_r(\phi_1,\ldots,\phi_r;P_{\leq\ell^{-1}})\partial_jg
		+\sum_{s=1}^{r}\mathcal C_r(\phi_1,\ldots,\partial_j\phi_s,\ldots,\phi_r;P_{\leq\ell^{-1}})g.
	\end{align*}
	We use the identity
	\begin{equation}
		\chi\mathcal C_r(\phi_1,\ldots,\phi_r;P_{\leq\ell^{-1}})g-\mathcal C_r(\phi_1,\ldots,\phi_r;P_{\leq\ell^{-1}})(\chi g)
		=\mathcal C_{r+1}(\phi_1,\ldots,\phi_r,\chi;P_{\leq\ell^{-1}})g,
		\label{eq:multiple-commutator-input-multiplication}
	\end{equation}
	which follows from the definition of $\mathcal C_{r+1}$.
	For $1\leq s\leq r$, we further use
	\begin{equation}
		\begin{aligned}
			&\chi\mathcal C_r(\phi_1,\ldots,\phi_{s-1},\psi,\phi_{s+1},\ldots,\phi_r;P_{\leq\ell^{-1}})g
			-\mathcal C_r(\phi_1,\ldots,\phi_{s-1},\chi\psi,\phi_{s+1},\ldots,\phi_r;P_{\leq\ell^{-1}})g\\
			=&-\mathcal C_r(\phi_1,\ldots,\widehat{\phi_s},\ldots,\phi_r,\chi;P_{\leq\ell^{-1}})(\psi g),
		\end{aligned}
		\label{eq:multiple-commutator-slot-multiplication}
	\end{equation}
	which follows from the kernel representation of $\mathcal C_r$ and
	the identity
	$$
	\chi(x)\bigl(\psi(x)-\psi(x-h)\bigr)
	-\bigl(\chi(x)\psi(x)-\chi(x-h)\psi(x-h)\bigr)
	=
	-\bigl(\chi(x)-\chi(x-h)\bigr)\psi(x-h).
	$$
	Apply \eqref{eq:multiple-commutator-input-multiplication} with
	$\chi=\upsilon_j$ to the term containing $\partial_jg$, and
	\eqref{eq:multiple-commutator-slot-multiplication} with
	$\chi=\upsilon_j$ and $\psi=\partial_j\phi_s$ to the coefficient terms.
	We obtain
	\eqref{eq:Dt-multiple-mollifier-commutator}.
\end{proof}

\begin{proof}[Proof of Proposition~\ref{prop:appendix-material-mollification}]
	We first expand the iterated commutators and use them to prove \eqref{est:appendix-material-mollification}.
	We then extend the spatial derivative range and treat the estimates with two spatial scales.
	For every $s\in\N_0$, define
	$$
	\operatorname{ad}_{D_{t,\upsilon}}^0(P_{\leq\ell^{-1}}):=P_{\leq\ell^{-1}},\qquad
	\operatorname{ad}_{D_{t,\upsilon}}^{s+1}(P_{\leq\ell^{-1}}):=[D_{t,\upsilon},\operatorname{ad}_{D_{t,\upsilon}}^s(P_{\leq\ell^{-1}})].
	$$
	Induction on $M$, using the definition of the iterated commutator
	and Pascal's identity, gives
	\begin{equation}
		D_{t,\upsilon}^MP_{\leq\ell^{-1}}f
		=
		\sum_{s=0}^{M}\binom{M}{s}\operatorname{ad}_{D_{t,\upsilon}}^s(P_{\leq\ell^{-1}})D_{t,\upsilon}^{M-s}f.
		\label{eq:material-mollifier-expansion}
	\end{equation}
	
	Let $s\geq1$. We place the spatial derivatives outside the material
	derivatives of $\upsilon$ using the identity
	$$
	[D_{t,\upsilon},\partial^\kappa]h
	=
	-\sum_{0<\sgm\leq\kappa}\binom{\kappa}{\sgm}(\partial^\sgm\upsilon)\cdot\nabla\partial^{\kappa-\sgm}h.
	$$
	In particular,
	$$
	D_{t,\upsilon}\bigl(\partial^\kappa D_{t,\upsilon}^{\sigma}\upsilon_j\bigr)
	=
	\partial^\kappa D_{t,\upsilon}^{\sigma+1}\upsilon_j
	-\sum_{0<\sgm\leq\kappa}\binom{\kappa}{\sgm}(\partial^\sgm\upsilon)\cdot\nabla\partial^{\kappa-\sgm}D_{t,\upsilon}^{\sigma}\upsilon_j.
	$$
	Using these identities, we prove by induction on $s$ that $\operatorname{ad}_{D_{t,\upsilon}}^s(P_{\leq\ell^{-1}})g$ is a finite linear combination of expressions of the following form. The coefficients and number of terms depend only on $d$ and $s$.
	$$
	\mathcal C_r\bigl(\partial^{\kappa_1}D_{t,\upsilon}^{\sigma_1}\upsilon_{j_1},\ldots,\partial^{\kappa_r}D_{t,\upsilon}^{\sigma_r}\upsilon_{j_r};P_{\leq\ell^{-1}}\bigr)
	\left(\left(\prod_{k=1}^{r_0}\partial^{\beta_k}D_{t,\upsilon}^{m_k}\upsilon_{i_k}\right)\partial^\gamma g\right),
	$$
	where
	\begin{align*}
		&r_0\in\N_0,\qquad 1\leq r\leq s,\qquad
		m_k,\sigma_\nu\in\N_0,\qquad
		i_k,j_\nu\in\{1,\ldots,d\},\\
		&\beta_k,\kappa_\nu,\gamma\in\N_0^d,\qquad
		|\beta_k|\geq1,\qquad |\gamma|\geq1,\\
		&\mathfrak m:=\sum_{k=1}^{r_0}m_k+\sum_{\nu=1}^{r}\sigma_\nu+r_0+r=s,\qquad
		q:=|\gamma|+\sum_{k=1}^{r_0}(|\beta_k|-1)+\sum_{\nu=1}^{r}|\kappa_\nu|=r.
	\end{align*}
	
	For $s=1$, one has
	$$
	[D_{t,\upsilon},P_{\leq\ell^{-1}}]g
	=
	\sum_{j=1}^{d}\mathcal C_1(\upsilon_j;P_{\leq\ell^{-1}})\partial_jg.
	$$
	Thus the representation holds with $r_0=0$, $r=1$, $\sigma_1=0$, $\kappa_1=0$, and $\gamma=e_j$. In particular, $\mathfrak m=q=r=1$.
	
	Suppose that the representation holds for $s$, and consider one term in this representation. Set
	$$
	A:=\prod_{k=1}^{r_0}\partial^{\beta_k}D_{t,\upsilon}^{m_k}\upsilon_{i_k},\qquad
	\phi_\nu:=\partial^{\kappa_\nu}D_{t,\upsilon}^{\sigma_\nu}\upsilon_{j_\nu},
	$$
	and define
	$$
	Tg:=\mathcal C_r(\phi_1,\ldots,\phi_r;P_{\leq\ell^{-1}})(A\partial^\gamma g).
	$$
	Since
	$
	\operatorname{ad}_{D_{t,\upsilon}}^{s+1}(P_{\leq\ell^{-1}})
	=
	[D_{t,\upsilon},\operatorname{ad}_{D_{t,\upsilon}}^s(P_{\leq\ell^{-1}})],
	$
	it suffices to compute the commutator $[D_{t,\upsilon},T]g$. By \eqref{eq:Dt-multiple-mollifier-commutator},
	$$
	\begin{aligned}
		\bigl[D_{t,\upsilon},T\bigr]g
		={}&\mathcal C_r(\phi_1,\ldots,\phi_r;P_{\leq\ell^{-1}})
		\bigl(D_{t,\upsilon}(A\partial^\gamma g)-A\partial^\gamma D_{t,\upsilon}g\bigr)
		+\sum_{\nu=1}^{r}\mathcal C_r(\phi_1,\ldots,D_{t,\upsilon}\phi_\nu,\ldots,\phi_r;P_{\leq\ell^{-1}})(A\partial^\gamma g)\\
		&+\sum_{j=1}^{d}\mathcal C_{r+1}(\phi_1,\ldots,\phi_r,\upsilon_j;P_{\leq\ell^{-1}})\partial_j(A\partial^\gamma g)
		-\sum_{\nu=1}^{r}\sum_{j=1}^{d}
		\mathcal C_r(\phi_1,\ldots,\widehat{\phi_\nu},\ldots,\phi_r,\upsilon_j;P_{\leq\ell^{-1}})
		\bigl((\partial_j\phi_\nu)A\partial^\gamma g\bigr).
	\end{aligned}
	$$
	For the four types of terms on the right, we compute the changes in
	$\mathfrak m$, $q$, and $r$.
	First,
	$$
	D_{t,\upsilon}(A\partial^\gamma g)-A\partial^\gamma D_{t,\upsilon}g
	=
	(D_{t,\upsilon}A)\partial^\gamma g+A[D_{t,\upsilon},\partial^\gamma]g.
	$$
	Commuting $D_{t,\upsilon}$ past $\partial^{\beta_k}$ in a factor of
	$A$ gives a term with $m_k$ replaced by $m_k+1$, together with
	commutator terms. In the first term, $\mathfrak m$ increases by one,
	while $r_0$, $r$, and $q$ remain unchanged. Each commutator term replaces one differentiated velocity factor by two. Hence $r_0$ and $\mathfrak m$ increase by one, and
	$
	(|\sgm|-1)+\bigl(|\beta_k-\sgm|+1-1\bigr)=|\beta_k|-1
	$
	shows that $q$ remains unchanged.
	Each term in $A[D_{t,\upsilon},\partial^\gamma]g$ introduces one new differentiated velocity factor in the input. Thus $r_0$ and $\mathfrak m$ increase by one, while
	$
	(|\sgm|-1)+|\gamma-\sgm|+1=|\gamma|
	$
	shows that $q$ remains unchanged.
	
	Commuting $D_{t,\upsilon}$ past $\partial^{\kappa_\nu}$ in $\phi_\nu$
	gives a term with $\sigma_\nu$ replaced by $\sigma_\nu+1$, together
	with commutator terms. In the first term, $\mathfrak m$ increases by
	one, while $r_0$, $r$, and $q$ remain unchanged.
	The commutator terms contain products of two velocity factors.
	For a product $B\psi$, we use
	$$
	\begin{aligned}
		&\mathcal C_r(\phi_1,\ldots,\phi_{\nu-1},B\psi,\phi_{\nu+1},\ldots,\phi_r;P_{\leq\ell^{-1}})h\\
		=&\mathcal C_r(\phi_1,\ldots,\phi_{\nu-1},\psi,\phi_{\nu+1},\ldots,\phi_r;P_{\leq\ell^{-1}})(Bh)
		+\mathcal C_{r+1}(\phi_1,\ldots,\phi_{\nu-1},\psi,\phi_{\nu+1},\ldots,\phi_r,B;P_{\leq\ell^{-1}})h\\
		&+\mathcal C_r(\phi_1,\ldots,\phi_{\nu-1},\phi_{\nu+1},\ldots,\phi_r,B;P_{\leq\ell^{-1}})(\psi h),
	\end{aligned}
	$$
	which follows from \eqref{eq:multiple-commutator-input-multiplication}--%
	\eqref{eq:multiple-commutator-slot-multiplication}. For the product terms in $D_{t,\upsilon}\phi_\nu$, take
	$
	B=\partial^\sgm\upsilon_j
	$
	and
	$
	\psi=\partial_j\partial^{\kappa_\nu-\sgm}D_{t,\upsilon}^{\sigma_\nu}\upsilon_{j_\nu}.
	$
	In the first term, the input is $Bh$. Hence $r_0$ and $\mathfrak m$ increase by one, and $r$ remains unchanged. Also,
	$
	(|\sgm|-1)+|\kappa_\nu-\sgm|+1=|\kappa_\nu|,
	$
	so $q$ remains unchanged. In the second term, $B$ becomes an additional commutator coefficient. Thus $r$ and $\mathfrak m$ increase by one, and
	$
	|\sgm|+|\kappa_\nu-\sgm|+1=|\kappa_\nu|+1
	$
	shows that $q$ also increases by one. In the third term, $B$ is a commutator coefficient and the input is $\psi h$. The material order $\sigma_\nu$ is now counted in $\sum m_k$ rather than $\sum \sigma_\nu$. The indices $r_0$ and $\mathfrak m$ increase by one, and $r$ remains unchanged. Since
	$
	|\sgm|+\bigl(|\kappa_\nu-\sgm|+1-1\bigr)=|\kappa_\nu|,
	$
	the value of $q$ remains unchanged.
	
	We now consider
	$
	\mathcal C_{r+1}(\phi_1,\ldots,\phi_r,\upsilon_j;P_{\leq\ell^{-1}})\partial_j(A\partial^\gamma g).
	$
	In all these terms, $r$ and $\mathfrak m$ increase by one. If the derivative acts on $g$, then $|\gamma|$ increases by one. If it acts on one of the factors in $A$, then $|\beta_k|-1$ increases by one. Thus $q$ also increases by one.
	
	Finally, consider
	$
	\mathcal C_r(\phi_1,\ldots,\widehat{\phi_\nu},\ldots,\phi_r,\upsilon_j;P_{\leq\ell^{-1}})
	\bigl((\partial_j\phi_\nu)A\partial^\gamma g\bigr).
	$
	Here $\upsilon_j$ replaces $\phi_\nu$ as a commutator coefficient, and
	$
	\partial_j\phi_\nu=\partial^{\kappa_\nu+e_j}D_{t,\upsilon}^{\sigma_\nu}\upsilon_{j_\nu}
	$
	multiplies $A\partial^\gamma g$ in the input. Hence $r$ remains unchanged, while $r_0$ and $\mathfrak m$ increase by one. The material order $\sigma_\nu$ is now counted in $\sum m_k$ rather than $\sum \sigma_\nu$, and
	$
	|\kappa_\nu+e_j|-1=|\kappa_\nu|
	$
	shows that $q$ remains unchanged.
	
	In each term, $\mathfrak m$ increases by one, while $q$ and $r$
	either remain unchanged or both increase by one. Thus
	$
	\mathfrak m=s+1,\ q=r.
	$
	The index ranges and the requirement that at least one spatial derivative act on $g$ are preserved,
	which proves the representation by induction.
	
	We now estimate the terms in this expansion. Let $s\geq1$. Since $|\gamma|\geq1$ in every term of the expansion, only spatial derivatives of $g$ of positive order are needed. For the $C^N$ estimate, assume that $N+s\leq K^*$ and
	$$
	\nrm{g}_{n,\mcI}\lesssim\mC_g\lambda_f^n,\qquad 1\leq n\leq N+s.
	$$
	For the H{\"o}lder estimate, assume in addition that $N+s\leq K^*-1$ and
	$$
	\nrm{g}_{n+\alpha,\mcI}\lesssim\mC_g\lambda_f^{n+\alpha},\qquad 1\leq n\leq N+s.
	$$
	The index relations imply
	$$
	\begin{aligned}
		\sum_{k=1}^{r_0}m_k+\sum_{\nu=1}^{r}\sigma_\nu
		+\sum_{k=1}^{r_0}(|\beta_k|-1)+\sum_{\nu=1}^{r}|\kappa_\nu|
		=
		s-r_0-r+r-|\gamma|
		=
		s-r_0-|\gamma|
		\leq s-1.
	\end{aligned}
	$$
	Apply \eqref{est:multiple-mollifier-commutator} and distribute the $N$ spatial derivatives among the input and the coefficient factors. The total number of material and additional spatial derivatives applied to the velocity gradients is at most $N+s-1$. Since $|\gamma|\leq q=r\leq s$, the order required of $g$ is at most $N+|\gamma|\leq N+s$. If $N+s\leq K^*$, all velocity norms of integer order lie in the range of \eqref{ass:appendix-material-mollification}. If $N+s\leq K^*-1$, all H{\"o}lder velocity norms are controlled by \eqref{est:appendix-material-mollification-interpolated-holder}.
	
	Set $k=r_0+r$. Since $\mathfrak m=s$, the material weights satisfy
	$$
	\begin{aligned}
	(\mC_\upsilon\lambda_\upsilon)^k
	\prod_{i=1}^{r_0}\MM{m_i,K_\upsilon-1,\mu_\upsilon,\tilde\mu_\upsilon}\prod_{\nu=1}^{r}\MM{\sigma_\nu,K_\upsilon-1,\mu_\upsilon,\tilde\mu_\upsilon}\leq\mu_*^k\MM{s-k,K_\upsilon-1,\mu_*,\tilde\mu_*}
	\leq\MM{s,K_f,\mu_*,\tilde\mu_*}.
	\end{aligned}
	$$
	Here the last inequality uses $K_\upsilon\geq\max\{1,K_f\}$.
	For a fixed term in the expansion, \eqref{est:multiple-mollifier-commutator}
	gives
	$$
	\begin{aligned}
	&\left\|\mathcal C_r\left(\partial^{\kappa_1}D_{t,\upsilon}^{\sigma_1}\upsilon_{j_1},\ldots,
	\partial^{\kappa_r}D_{t,\upsilon}^{\sigma_r}\upsilon_{j_r};P_{\leq\ell^{-1}}\right)
	\left(\left(\prod_{k=1}^{r_0}\partial^{\beta_k}D_{t,\upsilon}^{m_k}\upsilon_{i_k}\right)\partial^\gamma g\right)\right\|_{N,\mcI}\\
	&\quad\lesssim\ell^r
	\sum_{N_0+\sum_{k=1}^{r_0}N_k'+\sum_{\nu=1}^rN_\nu''=N}
	\nrm{g}_{N_0+|\gamma|,\mcI}
	\prod_{k=1}^{r_0}\nrm{\nabla D_{t,\upsilon}^{m_k}\upsilon_{i_k}}_{N_k'+|\beta_k|-1,\mcI}
	\prod_{\nu=1}^r\nrm{\nabla D_{t,\upsilon}^{\sigma_\nu}\upsilon_{j_\nu}}_{N_\nu''+|\kappa_\nu|,\mcI}.
	\end{aligned}
	$$
	By \eqref{ass:appendix-material-mollification} and the assumed
	$C^n$ estimates for $g$, the preceding norm is at most a constant times
	$$
	\begin{aligned}
		\mC_g\ell^r\lambda_*^N\lambda_f^{|\gamma|}
		\lambda_\upsilon^{q-|\gamma|}
		(\mC_\upsilon\lambda_\upsilon)^{r_0+r}
		\MM{s-r_0-r,K_\upsilon-1,\mu_\upsilon,\tilde\mu_\upsilon}.
	\end{aligned}
	$$
	The $N$ spatial derivatives contribute $\lambda_*^N$, and
	$$
	\sum_{k=1}^{r_0}(|\beta_k|-1)+\sum_{\nu=1}^{r}|\kappa_\nu|
	=
	q-|\gamma|.
	$$
	By the second inequality in \eqref{est:multiple-mollifier-commutator},
	the $C^{N+\alpha}$ norm is at most a constant times
	$$
	\begin{aligned}
	&\mC_g\ell^r\lambda_*^N\lambda_f^{|\gamma|}
	\lambda_\upsilon^{q-|\gamma|}
	(\mC_\upsilon\lambda_\upsilon)^{r_0+r}
	\MM{s-r_0-r,K_\upsilon-1,\mu_\upsilon,\tilde\mu_\upsilon}
	\bigl(\lambda_f^\alpha+\lambda_\upsilon^\alpha\bigr).
	\end{aligned}
	$$
	Using $q=r$ and the estimate for the material weights above,
	we estimate each term by a constant multiple of
	$
	\mC_g\ell^r\lambda_*^{N+r}\MM{s,K_f,\mu_*,\tilde\mu_*}
	$
	in $C^N$, and the same estimate with one additional factor $\lambda_*^\alpha$ in $C^{N+\alpha}$. Since $r\geq1$ and $\ell\lambda_*\leq1$,
	$
	\ell^r\lambda_*^r=(\ell\lambda_*)^r\leq\ell\lambda_*.
	$
	Consequently,
	\begin{equation}
		\begin{aligned}
			\nrm{\operatorname{ad}_{D_{t,\upsilon}}^s(P_{\leq\ell^{-1}})g}_{N,\mcI}
			&\lesssim\mC_g\ell\lambda_*^{N+1}\MM{s,K_f,\mu_*,\tilde\mu_*},
			&&N+s\leq K^*,\\
			\nrm{\operatorname{ad}_{D_{t,\upsilon}}^s(P_{\leq\ell^{-1}})g}_{N+\alpha,\mcI}
			&\lesssim\mC_g\ell\lambda_*^{N+1+\alpha}\MM{s,K_f,\mu_*,\tilde\mu_*},
			&&N+s\leq K^*-1.
		\end{aligned}
		\label{est:iterated-material-mollifier}
	\end{equation}
	
	We next use the vanishing moments to obtain the factor $\ell^L$. Fix $1\leq L\leq N_{\mathrm{mom}}$. For every $s\geq1$, we have
	\begin{equation}
		\begin{aligned}
			\nrm{\operatorname{ad}_{D_{t,\upsilon}}^s(P_{\leq\ell^{-1}})g}_{N,\mcI}
			&\lesssim\mC_g\ell^L\lambda_*^{N+L}\MM{s,K_f,\mu_*,\tilde\mu_*},
			&&N+s+L-1\leq K^*,\\
			\nrm{\operatorname{ad}_{D_{t,\upsilon}}^s(P_{\leq\ell^{-1}})g}_{N+\alpha,\mcI}
			&\lesssim\mC_g\ell^L\lambda_*^{N+L+\alpha}\MM{s,K_f,\mu_*,\tilde\mu_*},
			&&N+s+L-1\leq K^*-1.
		\end{aligned}
		\label{est:higher-order-iterated-material-mollifier}
	\end{equation}
	Here both estimates assume
	$
	\nrm{g}_{n,\mcI}\lesssim\mC_g\lambda_f^n,
	$
	and the second additionally assumes
	$
	\nrm{g}_{n+\alpha,\mcI}\lesssim\mC_g\lambda_f^{n+\alpha},
	$
	for $1\leq n\leq N+\max\{s,L\}$.
	
	Consider a fixed term with commutator order $r$. If $r\leq L$, apply \eqref{est:higher-order-multiple-mollifier-commutator} with $\gamma=0$ in that estimate to the input
	$
	\left(\prod_{k=1}^{r_0}\partial^{\beta_k}D_{t,\upsilon}^{m_k}\upsilon_{i_k}\right)\partial^\gamma g.
	$
	Distribute $N+L-r$ spatial derivatives among $g$ and the differentiated velocity factors. These derivatives contribute a spatial factor of at most $\lambda_*^{N+L-r}$. With $\mathfrak m=s$ and $q=r$, the term has norm at most a constant times
	$
	\mC_g\ell^L\lambda_*^{N+L}\MM{s,K_f,\mu_*,\tilde\mu_*}
	$
	in $C^N$, and the same estimate with one additional factor $\lambda_*^\alpha$ in $C^{N+\alpha}$.
	The order required of $g$ is at most
	$
	N+L-r+|\gamma|\leq N+L,
	$
	since $|\gamma|\leq q=r$. The index relation
	$$
	\sum_{k=1}^{r_0}m_k+\sum_{\nu=1}^{r}\sigma_\nu
	+\sum_{k=1}^{r_0}(|\beta_k|-1)+\sum_{\nu=1}^{r}|\kappa_\nu|
	\leq s-1
	$$
	shows that the total number of material and additional spatial derivatives applied to the velocity gradients is at most
	$
	(s-1)+(N+L-r)\leq N+s+L-2.
	$
	These orders lie in the assumed range for $g$ and in \eqref{ass:appendix-material-mollification} for the velocity gradients. For the H{\"o}lder estimate, use \eqref{est:appendix-material-mollification-interpolated-holder}.
	
	If $r>L$, apply \eqref{est:multiple-mollifier-commutator}. The $C^N$ norm is then at most a constant times
	$
	\mC_g\ell^r\lambda_*^{N+r}\MM{s,K_f,\mu_*,\tilde\mu_*}.
	$
	For the H{\"o}lder norm, include the factor $\lambda_*^\alpha$. The order required of $g$ is at most $N+|\gamma|\leq N+s$. The total order on the velocity gradients is at most $N+s-1$. Moreover,
	$
	\ell^r\lambda_*^{N+r}
	=
	\ell^L\lambda_*^{N+L}(\ell\lambda_*)^{r-L}
	\leq
	\ell^L\lambda_*^{N+L}.
	$
	This proves \eqref{est:higher-order-iterated-material-mollifier}.
	
	In \eqref{eq:material-mollifier-expansion}, estimate the term with $s=0$ using Proposition~\ref{prop:appendix-spatial-mollification}, \eqref{ass:appendix-material-mollification}, and \eqref{est:appendix-material-mollification-interpolated-holder}. For $1\leq s\leq M$, set
	$
	g=D_{t,\upsilon}^{M-s}f.
	$
	The estimates \eqref{ass:appendix-material-mollification} and
	\eqref{est:appendix-material-mollification-interpolated-holder} give
	the $C^n$ and $C^{n+\alpha}$ estimates for $g$ required for
	\eqref{est:iterated-material-mollifier}, with
	$\mC_g:=\mC_f\MM{M-s,K_f,\mu_f,\tilde\mu_f}$, for
	$M+N\leq K^*$ and $M+N\leq K^*-1$, respectively. Indeed,
	$n+M-s\leq N+M$ for every $0\leq n\leq N+s$.
	We use \eqref{eq:appendix-derivative-weight-product} in the form
	$
	\MM{s,K_f,\mu_*,\tilde\mu_*}\MM{M-s,K_f,\mu_f,\tilde\mu_f}
	\leq
	\MM{M,K_f,\mu_*,\tilde\mu_*}.
	$
	Apply \eqref{est:iterated-material-mollifier} for $1\leq s\leq M$,
	use $\ell\lambda_*\leq1$, and include the term $s=0$ to prove
	the first two estimates in \eqref{est:appendix-material-mollification}.
	
	For $L=0$, the last two estimates follow from
	\eqref{ass:appendix-material-mollification},
	\eqref{est:appendix-material-mollification-interpolated-holder},
	the first two estimates, and $P_{>\ell^{-1}}=I-P_{\leq\ell^{-1}}$.
	It remains to estimate the mollification difference for $1\leq L\leq N_{\mathrm{mom}}$. Since
	$
	P_{>\ell^{-1}}=I-P_{\leq\ell^{-1}},
	$
	we obtain from \eqref{eq:material-mollifier-expansion}
	$$
	D_{t,\upsilon}^MP_{>\ell^{-1}}f
	=
	P_{>\ell^{-1}}D_{t,\upsilon}^Mf
	-\sum_{s=1}^{M}\binom{M}{s}
	\operatorname{ad}_{D_{t,\upsilon}}^s(P_{\leq\ell^{-1}})D_{t,\upsilon}^{M-s}f.
	$$
	For $M+N+L\leq K^*$, \eqref{est:appendix-spatial-mollification} and \eqref{ass:appendix-material-mollification} give
	\begin{equation}\label{eq:material-mollification-spatial-tail-integer-estimate}
	\nrm{P_{>\ell^{-1}}D_{t,\upsilon}^Mf}_{N,\mcI}
	\lesssim
	\mC_f\ell^L\lambda_f^{N+L}\MM{M,K_f,\mu_f,\tilde\mu_f}.
	\end{equation}
	If $M+N+L\leq K^*-1$, interpolate \eqref{eq:material-mollification-spatial-tail-integer-estimate} at spatial orders $N$ and $N+1$ to obtain
	\begin{equation}\label{eq:material-mollification-spatial-tail-holder-estimate}
	\nrm{P_{>\ell^{-1}}D_{t,\upsilon}^Mf}_{N+\alpha,\mcI}
	\lesssim
	\mC_f\ell^L\lambda_f^{N+L+\alpha}\MM{M,K_f,\mu_f,\tilde\mu_f}.
	\end{equation}
	
	For the remaining terms, use the same $g$ and $\mC_g$. The $C^n$
	and $C^{n+\alpha}$ estimates now hold for $0\leq n\leq N+\max\{s,L\}$ under
	$M+N+L\leq K^*$ and $M+N+L\leq K^*-1$, respectively, since
	$
	n+M-s\leq N+M+\max\{s,L\}-s\leq N+M+L-1\leq K^*.
	$
	The remaining conditions in \eqref{est:higher-order-iterated-material-mollifier} hold because
	$
	N+s+L-1\leq N+M+L-1.
	$
	We may therefore apply \eqref{est:higher-order-iterated-material-mollifier}. Estimate~\eqref{eq:appendix-derivative-weight-product} gives
	$$
	\MM{s,K_f,\mu_*,\tilde\mu_*}\MM{M-s,K_f,\mu_f,\tilde\mu_f}
	\leq
	\MM{M,K_f,\mu_*,\tilde\mu_*}
	$$
	Combining this inequality with \eqref{est:higher-order-iterated-material-mollifier} and \eqref{eq:material-mollification-spatial-tail-integer-estimate}--\eqref{eq:material-mollification-spatial-tail-holder-estimate}, and summing over $s$, proves the last two estimates in \eqref{est:appendix-material-mollification}.
	
	We next prove \eqref{est:appendix-extended-spatial-mollification} and \eqref{est:appendix-two-scale-extended-spatial-mollification}, starting with the $C^N$ estimates. Fix $M,N\in\N_0$ such that
	$
	0\leq M\leq K^*
	$
	and
	$
	M+N\leq K^\sharp,
	$
	and use the scales specified for each case.
	For the two-scale case, set $(\lambda_*,\mu_*,\tilde\mu_*)=(\lambda,\mu,\tilde\mu)$ and require $M+N\geq1$. Replace $K^*-M$ by $(K_f-M)_+$ throughout the following argument.
	We first estimate the term $s=0$ by the kernel argument for \eqref{est:appendix-spatial-mollification-F-integer},
	transferring spatial derivatives above $K^*-M$ to the kernel.
	For the two-scale case with $M=0$, leave at least one spatial derivative on $f$
	at positive spatial orders. Use $\|P_{\leq\ell^{-1}}f\|_{0,\mcI}\lesssim\|f\|_{1,\mcI}$ for the $C^0$ part.
	Then \eqref{ass:appendix-material-mollification} or \eqref{ass:appendix-two-scale-material-mollification} gives
	\begin{equation}\label{eq:extended-material-mollification-uncommuted-estimate}
	\nrm{P_{\leq\ell^{-1}}D_{t,\upsilon}^Mf}_{N,\mcI}
	\lesssim
	\mC_f\MM{N,K^*-M,\lambda_*,\ell^{-1}}\MM{M,K_f,\mu_*,\tilde\mu_*}.
	\end{equation}
	
	For the commutator terms, let $1\leq s\leq M$, set
	$
	g:=D_{t,\upsilon}^{M-s}f,\
	\mC_g:=\mC_f\MM{M-s,K_f,\mu_*,\tilde\mu_*},
	$
	and consider a term of the form
	$$
	\mathcal C_r\bigl(\partial^{\kappa_1}D_{t,\upsilon}^{\sigma_1}\upsilon_{j_1},\ldots,\partial^{\kappa_r}D_{t,\upsilon}^{\sigma_r}\upsilon_{j_r};P_{\leq\ell^{-1}}\bigr)
	\left(\left(\prod_{k=1}^{r_0}\partial^{\beta_k}D_{t,\upsilon}^{m_k}\upsilon_{i_k}\right)\partial^\gamma g\right),
	$$
	where $\mathfrak m=s$ and $q=r$. Set
	$
	N_g:=\min\{N,K^*-M\}.
	$
	We apply \eqref{est:general-multiple-mollifier-commutator} to the entire input
	$
	\left(\prod_{k=1}^{r_0}\partial^{\beta_k}D_{t,\upsilon}^{m_k}\upsilon_{i_k}\right)\partial^\gamma g
	$
	with $\gamma=0$ in that estimate. Distribute the $N_0$ input derivatives and write
	$
	N_{0,g}+\sum_{k=1}^{r_0}N_{0,k}=N_0.
	$
	Since $1\leq|\gamma|\leq q=r\leq s$,
	$
	N_{0,g}+|\gamma|+M-s\leq N_g+M\leq K^*,
	$
	so all the required norms of $g$ are controlled by \eqref{ass:appendix-material-mollification} or \eqref{ass:appendix-two-scale-material-mollification}. The index relations also give
	$$
	\sum_{k=1}^{r_0}m_k+\sum_{\nu=1}^{r}\sigma_\nu
	+\sum_{k=1}^{r_0}(|\beta_k|-1)+\sum_{\nu=1}^{r}|\kappa_\nu|
	=
	s-r_0-|\gamma|
	\leq s-1.
	$$
	For each resulting term, the total number of material and additional spatial derivatives applied to the velocity gradients is therefore at most
	$$
	\sum_{k=1}^{r_0}\bigl(m_k+N_{0,k}+|\beta_k|-1\bigr)
	+
	\sum_{\nu=1}^{r}\bigl(\sigma_\nu+N_\nu+(M_\nu-1)_++|\kappa_\nu|\bigr)
	\leq
	N+s-1.
	$$
	Here we used
	$
	N_{0,g}+\sum_{k=1}^{r_0}N_{0,k}=N_0,\
	N_0+\sum_{\nu=1}^{r}N_\nu=N_g,\
	\sum_{\nu=1}^{r}(M_\nu-1)_+\leq N-N_g.
	$
	We have $s\leq M\leq K^*$ and $N+s-1\leq K^\sharp-1$. Use \eqref{ass:appendix-extended-velocity} or the extended range of \eqref{ass:appendix-two-scale-material-mollification}. Each material order $c\in\{m_k,\sigma_\nu\}$ satisfies $c\leq s-1$, so
	$
	K^*-1-c\geq K^*-M,\ (K_\upsilon-1-c)_+\geq(K_f-M)_+.
	$
	Lower the spatial thresholds of $g$ and the velocity gradients to $K^*-M$ in the first case and to $(K_f-M)_+$ in the two-scale case.
	
	For a fixed summand in \eqref{est:general-multiple-mollifier-commutator}, define
	$$
	\mathfrak d
	:=
	M_0+\#\{\nu\in\{1,\ldots,r\}:M_\nu\geq1\}.
	$$
	Since $M_0+\cdots+M_r=N-N_g$,
	$
	0\leq\mathfrak d\leq N-N_g=(N-(K^*-M))_+.
	$
	The commutator estimate contributes $\ell^{r-\mathfrak d}$, while the total spatial index entering the $\mathcal M$-factors is
	$$
		N_0+\sum_{\nu=1}^{r}\bigl(N_\nu+(M_\nu-1)_+\bigr)+q
		=
		N_g+\sum_{\nu=1}^{r}(M_\nu-1)_++r
		=
		N-M_0-\#\{\nu\in\{1,\ldots,r\}:M_\nu\geq1\}+r
		=
		N-\mathfrak d+r.
	$$
	Estimate~\eqref{eq:appendix-derivative-weight-product} gives
	$
	\prod_{j=1}^{J}\MM{N_j,K^*-M,\lambda_*,\ell^{-1}}
	\leq
	\MM{N_1+\cdots+N_J,K^*-M,\lambda_*,\ell^{-1}}
	$
	for all nonnegative integers $N_1,\ldots,N_J$. Absorb the velocity prefactors into $\mu_*$ and use $\mathfrak m=s$. The $C^N$ norm of each summand is at most a constant times
	$
	\mC_g\MM{s,K_f,\mu_*,\tilde\mu_*}
	\ell^{r-\mathfrak d}
	\MM{N-\mathfrak d+r,K^*-M,\lambda_*,\ell^{-1}}.
	$
	
	The required scale inequality is
	\begin{equation}
		\ell^{r-\mathfrak d}
		\MM{N-\mathfrak d+r,K^*-M,\lambda_*,\ell^{-1}}
		\leq
		\MM{N,K^*-M,\lambda_*,\ell^{-1}}.
		\label{est:appendix-extended-MM-count}
	\end{equation}
	If $N\geq K^*-M$, then $\mathfrak d\leq N-(K^*-M)$, so $N-\mathfrak d\geq K^*-M$ and
	$$
	\begin{aligned}
		\ell^{r-\mathfrak d}
		\MM{N-\mathfrak d+r,K^*-M,\lambda_*,\ell^{-1}}
		=
		\ell^{r-\mathfrak d}
		\lambda_*^{K^*-M}
		\ell^{-(N-\mathfrak d+r-(K^*-M))}
		=
		\lambda_*^{K^*-M}\ell^{-(N-(K^*-M))}
		=
		\MM{N,K^*-M,\lambda_*,\ell^{-1}}.
	\end{aligned}
	$$
	If $N<K^*-M$, then $\mathfrak d=0$ and
	$$
	\ell^r\MM{N+r,K^*-M,\lambda_*,\ell^{-1}}
	=
	(\ell\lambda_*)^{\min\{r,K^*-M-N\}}
	\MM{N,K^*-M,\lambda_*,\ell^{-1}}
	\leq
	\MM{N,K^*-M,\lambda_*,\ell^{-1}}.
	$$
	Summing over the expansion and using \eqref{est:appendix-extended-MM-count}, we obtain
	$$
	\nrm{
		\operatorname{ad}_{D_{t,\upsilon}}^s
		(P_{\leq\ell^{-1}})g
	}_{N,\mcI}
	\lesssim
	\mC_f
	\MM{N,K^*-M,\lambda_*,\ell^{-1}}
	\MM{M,K_f,\mu_*,\tilde\mu_*}.
	$$
	Together with \eqref{eq:extended-material-mollification-uncommuted-estimate}, \eqref{eq:material-mollifier-expansion} proves the $C^N$ estimates in \eqref{est:appendix-extended-spatial-mollification} and \eqref{est:appendix-two-scale-extended-spatial-mollification}.

	For the H{\"o}lder estimate, assume in addition that
	$
	M+N\leq K^\sharp-1.
	$
	Interpolate the $C^N$ and $C^{N+1}$ estimates to obtain
	$$
	\begin{aligned}
		\nrm{D_{t,\upsilon}^MP_{\leq\ell^{-1}}f}_{N+\alpha,\mcI}
		&\lesssim_\alpha
		\nrm{D_{t,\upsilon}^MP_{\leq\ell^{-1}}f}_{N,\mcI}^{1-\alpha}
		\nrm{D_{t,\upsilon}^MP_{\leq\ell^{-1}}f}_{N+1,\mcI}^{\alpha}\\
		&\lesssim
		\mC_f
		\MM{N,K^*-M,\lambda_*,\ell^{-1}}^{1-\alpha}
		\MM{N+1,K^*-M,\lambda_*,\ell^{-1}}^{\alpha}
		\MM{M,K_f,\mu_*,\tilde\mu_*}\\
		&=
		\mC_f
		\La{N,K^*-M,\lambda_*,\ell^{-1}}
		\MM{N,K^*-M,\lambda_*,\ell^{-1}}
		\MM{M,K_f,\mu_*,\tilde\mu_*}.
	\end{aligned}
	$$
	The H{\"o}lder estimates in \eqref{est:appendix-extended-spatial-mollification} and \eqref{est:appendix-two-scale-extended-spatial-mollification} follow.
	It remains to prove \eqref{est:appendix-two-scale-material-mollification}.
	For $\sum_i M_i=M$ and $\sum_i N_i=N$, \eqref{eq:appendix-derivative-weight-product} gives
	\begin{equation}\label{eq:appendix-material-mollification-weight-product}
	\prod_i\left(\MM{M_i,K_f,\mu,\tilde\mu}\MM{N_i,(K_f-M_i)_+,\lambda,\ell^{-1}}\right)
	\leq\MM{M,K_f,\mu,\tilde\mu}\MM{N,(K_f-M)_+,\lambda,\ell^{-1}}.
	\end{equation}
	Indeed, the material weights satisfy $\prod_i\MM{M_i,K_f,\mu,\tilde\mu}
	\leq\MM{M,K_f,\mu,\tilde\mu}$, and $(K_f-M_i)_+\geq(K_f-M)_+$ gives
	$$
	\prod_i\MM{N_i,(K_f-M_i)_+,\lambda,\ell^{-1}}
	\leq\prod_i\MM{N_i,(K_f-M)_+,\lambda,\ell^{-1}}
	\leq\MM{N,(K_f-M)_+,\lambda,\ell^{-1}}.
	$$
	For each velocity-gradient factor, use $\mC_\upsilon\lambda\leq\mu$ and $K_\upsilon\geq\max\{1,K_f\}$ to obtain
	$$
	\mu\MM{m,K_\upsilon-1,\mu,\tilde\mu}\MM{n,(K_\upsilon-1-m)_+,\lambda,\ell^{-1}}
	\leq\MM{m+1,K_f,\mu,\tilde\mu}\MM{n,(K_f-m-1)_+,\lambda,\ell^{-1}}.
	$$
	Use \eqref{eq:material-mollifier-expansion} with $P_{>\ell^{-1}}=I-P_{\leq\ell^{-1}}$.
	By \eqref{est:appendix-spatial-mollification}, the term with $s=0$ satisfies
	$$
	\begin{aligned}
	\nrm{P_{>\ell^{-1}}D_{t,\upsilon}^Mf}_{N,\mcI}
	\lesssim\ell^L\nrm{D_{t,\upsilon}^Mf}_{N+L,\mcI}
	\lesssim\mC_f\ell^L\MM{M,K_f,\mu,\tilde\mu}
	\MM{N+L,(K_f-M)_+,\lambda,\ell^{-1}}.
	\end{aligned}
	$$
	For $1\leq s\leq M$, set $g=D_{t,\upsilon}^{M-s}f$.
	For $1\leq n\leq K^*-M+s$, the hypothesis on $f$ in
	\eqref{ass:appendix-two-scale-material-mollification} gives
	$$
	\nrm{g}_{n,\mcI}
	\lesssim\mC_f\MM{M-s,K_f,\mu,\tilde\mu}
	\MM{n,(K_f-M+s)_+,\lambda,\ell^{-1}}.
	$$
	Consider a fixed term in the expansion with commutator order $r$.
	If $r\leq L$, apply \eqref{est:higher-order-multiple-mollifier-commutator}, as in the proof of
	\eqref{est:higher-order-iterated-material-mollifier}.
	Distribute the $N$ spatial derivatives and the additional $L-r$ derivatives among $g$ and the differentiated velocity factors.
	The estimate for each velocity-gradient factor and $\mathfrak m=s$ give total material index $M$. The spatial indices sum to $N+L$ by $q=r$.
	By \eqref{eq:appendix-material-mollification-weight-product}, this term satisfies the same estimate as for $s=0$.
	If $r>L$, apply \eqref{est:multiple-mollifier-commutator} and use
	$$
	\ell^r\MM{N+r,(K_f-M)_+,\lambda,\ell^{-1}}
	\leq\ell^L\MM{N+L,(K_f-M)_+,\lambda,\ell^{-1}}.
	$$
	The required total orders of $f$ and the velocity gradients are at most
	$M-s+N+\max\{s,L\}\leq M+N+L-1$ and $N+s+L-2\leq M+N+L-2$, respectively.
	Thus \eqref{ass:appendix-two-scale-material-mollification} covers every factor, and
	$$
	\begin{aligned}
	\nrm{\operatorname{ad}_{D_{t,\upsilon}}^s(P_{\leq\ell^{-1}})D_{t,\upsilon}^{M-s}f}_{N,\mcI}
	&\lesssim\mC_f\ell^L\MM{M,K_f,\mu,\tilde\mu}
	\MM{N+L,(K_f-M)_+,\lambda,\ell^{-1}}.
	\end{aligned}
	$$
	Summing over $s$ proves \eqref{est:appendix-two-scale-material-mollification}.
	For $M+N+L\leq K^*-1$, use this estimate at spatial orders $N$ and $N+1$.
	Spatial interpolation gives
	$$
	\begin{aligned}
		\nrm{D_{t,\upsilon}^MP_{>\ell^{-1}}f}_{N+\alpha,\mcI}
		&\lesssim_\alpha
		\nrm{D_{t,\upsilon}^MP_{>\ell^{-1}}f}_{N,\mcI}^{1-\alpha}
		\nrm{D_{t,\upsilon}^MP_{>\ell^{-1}}f}_{N+1,\mcI}^{\alpha}\\
		&\lesssim\mC_f\ell^L\MM{M,K_f,\mu,\tilde\mu}
		\MM{N+L,(K_f-M)_+,\lambda,\ell^{-1}}^{1-\alpha}
		\MM{N+L+1,(K_f-M)_+,\lambda,\ell^{-1}}^{\alpha}\\
		&=\mC_f\ell^L\MM{M,K_f,\mu,\tilde\mu}
		\La{N+L,(K_f-M)_+,\lambda,\ell^{-1}}
		\MM{N+L,(K_f-M)_+,\lambda,\ell^{-1}}.
	\end{aligned}
	$$
	We obtain \eqref{est:appendix-two-scale-material-mollification-holder}.
\end{proof}

\begin{lem} 
	\label{lem:appendix-change-material-derivative} 
	Fix a time interval $\mcI\subset\R$ and $K,J\in\N$ with
	$1\leq J\leq K$. Let $u,\overline u:\mcI\times\T^d\to\R^d$ be
	smooth. Set 
	$ 
	D_{t,u}:=\partial_t+u\cdot\nabla, 
	D_{t,\overline u}:=\partial_t+\overline u\cdot\nabla. 
	$ 
	Fix $0<\mC_u\leq1$ and $\lambda,\mu\geq1$ with
	$\mC_u\lambda\leq\mu$. 
	Assume that, for every $j,N\in\N_0$ with $0\leq j<J$ and 
	$j+N\leq K-1$, 
	$$ 
	\begin{aligned} 
		\nrm{D_{t,\overline u}^j(\overline u-u)}_{N,\mcI} 
		&\lesssim 
		\mC_u\mu^j\lambda^N,\\ 
		\nrm{D_{t,\overline u}^j\nabla\overline u}_{N,\mcI} 
		+ 
		\nrm{\nabla D_{t,\overline u}^j\overline u}_{N,\mcI} 
		&\lesssim 
		\mC_u\mu^j\lambda^{N+1}. 
	\end{aligned} 
	$$ 
	Let $H$ be a smooth scalar, vector, or tensor field on $\mcI\times\T^d$.
	Fix $\mC_H>0$ and assume that
	$$ 
	\nrm{D_{t,u}^{M}H}_{N,\mcI} 
	\lesssim 
	\mC_H\mu^{M}\lambda^N, 
	\qquad 
	1\leq M+N\leq K. 
	$$ 
	Then 
	\begin{align} 
		\nrm{D_{t,\overline u}^{M_1}D_{t,u}^{M_2}H}_{N,\mcI} 
		&\lesssim 
		\mC_H\mu^{M_1+M_2}\lambda^N, 
		&& 
		M_1\leq J,\quad 
		1\leq M_1+M_2+N\leq K,
		\label{est:appendix-change-material-derivative}\\
		\nrm{D_{t,\overline u}^{M_1}\nabla D_{t,u}^{M_2}H}_{N,\mcI} 
		+ 
		\nrm{\nabla D_{t,\overline u}^{M_1}D_{t,u}^{M_2}H}_{N,\mcI} 
		&\lesssim 
		\mC_H\mu^{M_1+M_2}\lambda^{N+1}, 
		&& 
		M_1\leq J,\quad 
		M_1+M_2+N\leq K-1.
		\label{est:appendix-change-material-derivative-gradient} 
	\end{align} 

\end{lem}

\begin{proof}
	We prove both estimates by induction on $M_1$, uniformly in $M_2$. The
	case $M_1=0$ follows from the assumptions on $H$. Let $1\leq M_1\leq J$
	and suppose that both conclusions hold at every smaller value of $M_1$.
	Since
	$D_{t,\overline u}=D_{t,u}+(\overline u-u)\cdot\nabla$,
	$$
	\begin{aligned}
		D_{t,\overline u}^{M_1}D_{t,u}^{M_2}H
		=
		D_{t,\overline u}^{M_1-1}D_{t,u}^{M_2+1}H
		+
		\sum_{j=0}^{M_1-1}\binom{M_1-1}{j}
		D_{t,\overline u}^j(\overline u-u)\cdot
		D_{t,\overline u}^{M_1-1-j}\nabla D_{t,u}^{M_2}H.
	\end{aligned}
	$$
	The first term satisfies \eqref{est:appendix-change-material-derivative}
	by induction. The remaining terms satisfy
	$$
	\begin{aligned}
		\nrm{
			D_{t,\overline u}^j(\overline u-u)\cdot
			D_{t,\overline u}^{M_1-1-j}\nabla D_{t,u}^{M_2}H
		}_{N,\mcI}\lesssim
		\sum_{N_1+N_2=N}
		\mC_u\mu^j\lambda^{N_1}
		\mC_H\mu^{M_1-1-j+M_2}\lambda^{N_2+1}
		\lesssim
		\mC_H\mu^{M_1+M_2}\lambda^N.
	\end{aligned}
	$$
	Here $j+N_1\leq K-1$ and
	$M_1-1-j+M_2+N_2\leq K-1$, while
	$\mC_u\lambda\leq\mu$. Hence
	\eqref{est:appendix-change-material-derivative} holds at order $M_1$.
	
	For the gradient estimate, apply
	Lemma~\ref{lem:appendix-scale-dependent-commuted-gradient} with
	$\upsilon=\overline u$, $\mC_\upsilon=\mC_u$,
	$G=D_{t,u}^{M_2}H$,
	$\mC_G=\mC_H\mu^{M_2}$, and $\widetilde\lambda=\lambda$, taking
	$K_\upsilon=K_G=K$, $K^\sharp=K-M_2$, $N_t=0$, and $\tilde\mu=\mu$ in that lemma. Estimate~\eqref{est:appendix-change-material-derivative},
	already proved for every material order not exceeding $M_1$, verifies
	\eqref{ass:appendix-scale-dependent-uncommuted-field}. The assumed estimates
	for $D_{t,\overline u}^j\nabla\overline u$, $j<M_1$, verify
	\eqref{ass:appendix-scale-dependent-lower-gradient}. We obtain
	\eqref{est:appendix-change-material-derivative-gradient}, which completes
	the induction.
\end{proof}

The next lemma requires estimates for material derivatives only along the original velocity field.

\begin{lem}
	\label{lem:source-material-derivative-change}
	Fix a time interval $\mcI\subset\R$ and $K,J\in\N$ with
	$1\leq J\leq K$. Let $\upsilon,\omega:\mcI\times\T^d\to\R^d$
	be smooth. Set
	$D_{t,\upsilon}:=\partial_t+\upsilon\cdot\nabla$ and
	$D_{t,\upsilon+\omega}:=D_{t,\upsilon}+\omega\cdot\nabla$.
	Fix $N_t\in\N_0$, $\widetilde{\mC}_0\geq1$, and $\mC_\omega,\mC_F>0$.
	Let $F$ be a smooth scalar, vector, or tensor field on $\mcI\times\T^d$.
	Fix $1\leq\lambda\leq\tilde\lambda$ and $1\leq\mu\leq\tilde\mu$ such that
	$
	\mC_\omega\lambda\leq\mu,\ \mC_\omega\tilde\lambda\leq\tilde\mu.
	$
	For $M,N\in\N_0$ in the ranges specified below, assume that
	\begin{equation}\label{ass:source-two-scale-material-change}
	\begin{aligned}
	\nrm{D_{t,\upsilon}^M\omega}_{N,\mcI}
	&\leq\widetilde{\mC}_0\mC_\omega\MM{M,N_t,\mu,\tilde\mu}\MM{N,(N_t-M)_+,\lambda,\tilde\lambda},\\
	\nrm{D_{t,\upsilon}^M\nabla\upsilon}_{N,\mcI}
	&\leq\widetilde{\mC}_0\mu\MM{M,N_t,\mu,\tilde\mu}\MM{N,(N_t-M)_+,\lambda,\tilde\lambda},\\
	\nrm{D_{t,\upsilon}^MF}_{N,\mcI}
	&\leq\mC_F\MM{M,N_t,\mu,\tilde\mu}\MM{N,(N_t-M)_+,\lambda,\tilde\lambda}.
	\end{aligned}
	\end{equation}
	The first two estimates are assumed for $M<J$, with $M+N\leq K-1$ and $M+N\leq K-2$, respectively.
	The gradient condition is imposed only for $K\geq2$.
	The third is assumed for $M\leq J$ and $1\leq M+N\leq K$.
	Then, for $M\leq J$ and $1\leq M+N\leq K$,
	\begin{equation}\label{eq:source-material-derivative-change}
	\nrm{D_{t,\upsilon+\omega}^MF}_{N,\mcI}
	\lesssim\mC_F\MM{M,N_t,\mu,\tilde\mu}\MM{N,(N_t-M)_+,\lambda,\tilde\lambda}.
	\end{equation}
	The implicit constant may depend on $K,\widetilde{\mC}_0$ but not on the amplitudes or scales.
	If the hypothesis on $F$ also holds at $M=N=0$, so does the conclusion.
\end{lem}

\begin{proof}
	The case $M=0$ follows from the hypothesis on $F$. For $1\leq M\leq J$, expand $D_{t,\upsilon+\omega}^MF$ as in the proof of
	Lemma~\ref{lem:total-order-upgrade-material}, with $f=F$ and $N_0=K$.
	A term with $k$ factors involving $\omega$ and $h$ factors involving
	$\nabla\upsilon$ has total material order $M-k-h$ and total spatial order at most $N+k$,
	counting only additional spatial derivatives on $\nabla\upsilon$.
	By \eqref{eq:appendix-material-mollification-weight-product}
	with $K_f=N_t$ and $\ell=\tilde\lambda^{-1}$, its norm is at most a fixed constant times
	$$
	\begin{aligned}
	&\mC_F\mC_\omega^k\mu^h\MM{M-k-h,N_t,\mu,\tilde\mu}
	\MM{N+k,(N_t-M+k+h)_+,\lambda,\tilde\lambda}\\
	&\quad\leq\mC_F\mC_\omega^k\MM{M-k,N_t,\mu,\tilde\mu}
	\MM{N+k,(N_t-M+k)_+,\lambda,\tilde\lambda}\\
	&\quad\leq\mC_F\MM{M,N_t,\mu,\tilde\mu}\MM{N,(N_t-M)_+,\lambda,\tilde\lambda}.
	\end{aligned}
	$$
	For the last inequality, apply $k$ times the identity
	$$
	\frac{\MM{r-1,N_t,\mu,\tilde\mu}\MM{n+1,(N_t-r+1)_+,\lambda,\tilde\lambda}}
	{\MM{r,N_t,\mu,\tilde\mu}\MM{n,(N_t-r)_+,\lambda,\tilde\lambda}}
	=\begin{cases}\lambda/\mu,&1\leq r\leq N_t,\\
	\tilde\lambda/\tilde\mu,&r>N_t.
	\end{cases}
	$$
	The factor $\widetilde{\mC}_0^{k+h}$ is absorbed in the implicit constant since $k+h\leq M\leq K$.
	The expansion requires total orders at most $M+N$, $M+N-1$, and $M+N-2$
	on $F$, $\omega$, and $\nabla\upsilon$, respectively. Their material orders
	are at most $M$, $M-1$, and $M-1$. Hence all factors lie in the stated ranges.
	Summation proves \eqref{eq:source-material-derivative-change}.
\end{proof}

We now estimate the velocity obtained by dividing the mollified momentum
by the mollified density.

\begin{lem}
	\label{lem:appendix-mollified-quotient-velocity}
	Fix a time interval $\mcI\subset\R$, $\ell>0$, $K^*,K^\sharp\in\N$ with
	$K^\sharp>K^*\geq2$ and $K^*\leq N_{\mathrm{mom}}$,
	and $\varrho_*>0$. Let
	$\varrho:\mcI\times\T^d\to\R$ and $w:\mcI\times\T^d\to\R^d$ be
	smooth, and set
	$$
	\begin{gathered}
		\upsilon:=\frac{w}{\varrho},
		\qquad
		\overline\varrho:=P_{\leq\ell^{-1}}\varrho,
		\qquad
		\overline w:=P_{\leq\ell^{-1}}w,
		\qquad
		\overline\upsilon:=\frac{\overline w}{\overline\varrho},\qquad
		D_{t,\upsilon}:=\partial_t+\upsilon\cdot\nabla,
		\qquad
		D_{t,\overline\upsilon}:=\partial_t+\overline\upsilon\cdot\nabla.
	\end{gathered}
	$$
	Let $0<\mC_\varrho,\mC_w\leq1$ and $\lambda,\mu\geq1$, set
	$\mC_\upsilon:=\max\{\mC_\varrho,\mC_w\}$, and assume that
	$\mC_\upsilon\lambda\leq\mu$ and $\ell\lambda\leq1$. Suppose that there
	exist constants $0<c_0\leq1$ and $\widetilde{\mC}_0\geq1$ such that
	$$
	c_0\leq\varrho,\overline\varrho\leq c_0^{-1}
	\qquad\text{on }\mcI\times\T^d,
	$$
	and
	$$
	\nrm{\varrho-\varrho_*}_{0,\mcI}+\nrm{w}_{0,\mcI}\leq\widetilde{\mC}_0.
	$$
	Assume further that, for every $M,N\in\N_0$ with 
	$1\leq M+N\leq K^*$, 
	$$ 
	\begin{aligned} 
		\nrm{D_{t,\upsilon}^M(\varrho-\varrho_*)}_{N,\mcI} 
		\leq\mC_\varrho\mu^M\lambda^N,\quad
		\nrm{D_{t,\upsilon}^Mw}_{N,\mcI} 
		\leq\mC_w\mu^M\lambda^N. 
	\end{aligned} 
	$$

	Then
	$$
	\nrm{\overline\varrho-\varrho_*}_{0,\mcI}
	+
	\nrm{\overline w}_{0,\mcI}
	+
	\nrm{\overline\upsilon}_{0,\mcI}
	\lesssim\widetilde{\mC}_0.
	$$
	For every $M,N\in\N_0$ satisfying $M\leq K^*$ and
	$1\leq M+N\leq K^\sharp$, one has
	\begin{equation}
		\begin{aligned}
			\nrm{D_{t,\overline\upsilon}^{\,M}(\overline\varrho-\varrho_*)}_{N,\mcI}
			&\lesssim
			\mC_\varrho\mu^M
			\MM{N,K^*-M,\lambda,\ell^{-1}},\\
			\nrm{D_{t,\overline\upsilon}^{\,M}\overline w}_{N,\mcI}
			&\lesssim
			\mC_w\mu^M
			\MM{N,K^*-M,\lambda,\ell^{-1}},\\
			\nrm{D_{t,\overline\upsilon}^{\,M}\overline\upsilon}_{N,\mcI}
			&\lesssim
			\mC_\upsilon\mu^M
			\MM{N,K^*-M,\lambda,\ell^{-1}}.
		\end{aligned}
		\label{est:appendix-mollified-quotient-main}
	\end{equation}
	For every $M,N\in\N_0$ satisfying $M\leq K^*$ and
	$M+N\leq K^\sharp-1$, one also has
	$$
	\nrm{D_{t,\overline\upsilon}^{\,M}\nabla\overline\upsilon}_{N,\mcI}
	+
	\nrm{\nabla D_{t,\overline\upsilon}^{\,M}\overline\upsilon}_{N,\mcI}
	\lesssim
	\mC_\upsilon\mu^M\MM{N+1,K^*-M,\lambda,\ell^{-1}}.
	$$

	For every $N\in\N_0$ with $N\leq K^*$, the spatial mollification differences satisfy
	\begin{equation}
		\begin{aligned}
			\nrm{\varrho-\overline\varrho}_{N,\mcI}
			&\lesssim
			\mC_\varrho\ell^{K^*-N}\lambda^{K^*},\\
			\nrm{w-\overline w}_{N,\mcI}
			&\lesssim
			\mC_w\ell^{K^*-N}\lambda^{K^*},\\
			\nrm{\upsilon-\overline\upsilon}_{N,\mcI}
			&\lesssim
			\mC_\upsilon\ell^{K^*-N}\lambda^{K^*}.
		\end{aligned}
		\label{est:appendix-mollified-quotient-high-frequency}
	\end{equation}
	For every $M,N\in\N_0$ satisfying $1\leq M\leq K^*$ and
	$M+N\leq K^*$, the material derivatives of the mollification differences satisfy
	\begin{equation}
		\begin{aligned}
			\nrm{D_{t,\overline\upsilon}^{\,M}(\varrho-\overline\varrho)}_{N,\mcI}
			&\lesssim
			\mC_\varrho\mu^M\ell^{K^*-M-N}\lambda^{K^*-M},\\
			\nrm{D_{t,\overline\upsilon}^{\,M}(w-\overline w)}_{N,\mcI}
			&\lesssim
			\mC_w\mu^M\ell^{K^*-M-N}\lambda^{K^*-M},\\
			\nrm{D_{t,\overline\upsilon}^{\,M}(\upsilon-\overline\upsilon)}_{N,\mcI}
			&\lesssim
			\mC_\upsilon\mu^M\ell^{K^*-M-N}\lambda^{K^*-M}.
		\end{aligned}
		\label{est:appendix-mollified-quotient-material-high-frequency}
	\end{equation}

	If, in addition,
	$$
	\partial_t\varrho+\Div w=0,
	$$
	then, for every $M,N\in\N_0$ satisfying $M\leq K^*-1$ and
	$M+N\leq K^\sharp-1$, one has
	$$
	\nrm{D_{t,\overline\upsilon}^{\,M}\Div\overline\upsilon}_{N,\mcI}
	\lesssim
	\mC_\varrho\mu^{M+1}
	\MM{N,K^*-M-1,\lambda,\ell^{-1}}.
	$$

	For every $\alpha\in(0,1)$, the estimates corresponding to
	$M=N=0$ are
	$$
	\begin{aligned}
		\nrm{\overline\varrho-\varrho_*}_{\alpha,\mcI}
		\lesssim
		\mC_\varrho^\alpha\lambda^\alpha,\quad
		\nrm{\overline w}_{\alpha,\mcI}
		\lesssim
		\mC_w^\alpha\lambda^\alpha,\quad
		\nrm{\overline\upsilon}_{\alpha,\mcI}
		\lesssim
		\mC_\upsilon^\alpha\lambda^\alpha.
	\end{aligned}
	$$
	For every $M,N\in\N_0$ satisfying $M\leq K^*$ and
	$1\leq M+N\leq K^\sharp-1$, the corresponding H{\"o}lder estimates are
	\begin{equation}
		\begin{aligned}
			\nrm{D_{t,\overline\upsilon}^{\,M}(\overline\varrho-\varrho_*)}_{N+\alpha,\mcI}
			&\lesssim
			\mC_\varrho\mu^M
			\La{N,K^*-M,\lambda,\ell^{-1}}
			\MM{N,K^*-M,\lambda,\ell^{-1}},\\
			\nrm{D_{t,\overline\upsilon}^{\,M}\overline w}_{N+\alpha,\mcI}
			&\lesssim
			\mC_w\mu^M
			\La{N,K^*-M,\lambda,\ell^{-1}}
			\MM{N,K^*-M,\lambda,\ell^{-1}},\\
			\nrm{D_{t,\overline\upsilon}^{\,M}\overline\upsilon}_{N+\alpha,\mcI}
			&\lesssim
			\mC_\upsilon\mu^M
			\La{N,K^*-M,\lambda,\ell^{-1}}
			\MM{N,K^*-M,\lambda,\ell^{-1}}.
		\end{aligned}
		\label{est:appendix-mollified-quotient-holder}
	\end{equation}
	For every $M,N\in\N_0$ satisfying $M\leq K^*$ and
	$M+N\leq K^\sharp-2$, one also has
	$$
	\begin{aligned}
		\nrm{D_{t,\overline\upsilon}^{\,M}\nabla\overline\upsilon}_{N+\alpha,\mcI}
		+
		\nrm{\nabla D_{t,\overline\upsilon}^{\,M}\overline\upsilon}_{N+\alpha,\mcI}
		\lesssim
		\mC_\upsilon\mu^M
		\La{N+1,K^*-M,\lambda,\ell^{-1}}
		\MM{N+1,K^*-M,\lambda,\ell^{-1}}.
	\end{aligned}
	$$

	For every $N\in\N_0$ with $N\leq K^*-1$, the spatial H{\"o}lder estimates for the mollification differences are
	\begin{equation}
		\begin{aligned}
			\nrm{\varrho-\overline\varrho}_{N+\alpha,\mcI}
			&\lesssim
			\mC_\varrho\ell^{K^*-N-\alpha}\lambda^{K^*},\\
			\nrm{w-\overline w}_{N+\alpha,\mcI}
			&\lesssim
			\mC_w\ell^{K^*-N-\alpha}\lambda^{K^*},\\
			\nrm{\upsilon-\overline\upsilon}_{N+\alpha,\mcI}
			&\lesssim
			\mC_\upsilon\ell^{K^*-N-\alpha}\lambda^{K^*}.
		\end{aligned}
		\label{est:appendix-mollified-quotient-high-frequency-holder}
	\end{equation}
	For every $M,N\in\N_0$ satisfying $1\leq M\leq K^*-1$ and
	$M+N\leq K^*-1$, the material derivatives of the mollification differences satisfy
	\begin{equation}
		\begin{aligned}
			\nrm{D_{t,\overline\upsilon}^{\,M}(\varrho-\overline\varrho)}_{N+\alpha,\mcI}
			&\lesssim
			\mC_\varrho\mu^M\ell^{K^*-M-N-\alpha}\lambda^{K^*-M},\\
			\nrm{D_{t,\overline\upsilon}^{\,M}(w-\overline w)}_{N+\alpha,\mcI}
			&\lesssim
			\mC_w\mu^M\ell^{K^*-M-N-\alpha}\lambda^{K^*-M},\\
			\nrm{D_{t,\overline\upsilon}^{\,M}(\upsilon-\overline\upsilon)}_{N+\alpha,\mcI}
			&\lesssim
			\mC_\upsilon\mu^M\ell^{K^*-M-N-\alpha}\lambda^{K^*-M}.
		\end{aligned}
		\label{est:appendix-mollified-quotient-material-high-frequency-holder}
	\end{equation}

	Under the continuity equation, for every $M,N\in\N_0$ satisfying
	$M\leq K^*-1$ and $M+N\leq K^\sharp-2$, one also has
	$$
	\nrm{D_{t,\overline\upsilon}^{\,M}\Div\overline\upsilon}_{N+\alpha,\mcI}
	\lesssim
	\mC_\varrho\mu^{M+1}
	\La{N,K^*-M-1,\lambda,\ell^{-1}}
	\MM{N,K^*-M-1,\lambda,\ell^{-1}}.
	$$
	The implicit constants may depend on $K^*,K^\sharp,c_0,\widetilde{\mC}_0,\alpha$, but not on
	$\mC_\varrho,\mC_w,\mC_\upsilon,$
	$\lambda,\ell$, or $\mu$.

\end{lem}

\begin{proof}
	We first estimate the original velocity $\upsilon$ and prove the spatial $C^N$ estimates.
	We then prove the $C^N$ estimates along $\overline\upsilon$ by induction on $M$.
	We finish with the endpoint $M+N=K^*$ and the divergence and H{\"o}lder estimates.
	The lower bound
	$\varrho\geq c_0$ gives
	$$
	\nrm{\upsilon}_{0,\mcI} = \Nrm{\frac{w}{\varrho}}_{0,\mcI} \lesssim \nrm{w}_{0,\mcI} \lesssim \widetilde{\mC}_0.
	$$
	For every $1\leq N\leq K^*$, apply the Fa\`a di Bruno formula \eqref{eq:faa-di-bruno}
	to the reciprocal function on a neighborhood of $[c_0,c_0^{-1}]$
	to obtain
	$$
	\begin{aligned}
		\nrm{\nabla^N(\varrho^{-1})}_{0,\mcI} \lesssim \sum_{k=1}^N \sum_{\substack{N_1+\cdots+N_k=N\\N_j\geq1}} \prod_{j=1}^k\nrm{\nabla^{N_j}(\varrho-\varrho_*)}_{0,\mcI}
		\lesssim \sum_{k=1}^N\mC_\varrho^k\lambda^N \lesssim \mC_\varrho\lambda^N.
	\end{aligned}
	$$
	In particular, for every $0\leq N\leq K^*$,
	$$
	\nrm{\varrho^{-1}}_{N,\mcI} \lesssim \lambda^N.
	$$
	For $1\leq N\leq K^*$, differentiating $\upsilon=w\varrho^{-1}$ yields
	\begin{equation}\label{eq:quotient-velocity-spatial-estimate}
	\begin{aligned}
		\nrm{\upsilon}_{N,\mcI} \lesssim \nrm{\upsilon}_{0,\mcI} + \sum_{N_1+N_2=N} \nrm{\nabla^{N_1}w}_{0,\mcI} \nrm{\nabla^{N_2}(\varrho^{-1})}_{0,\mcI}
		\lesssim \mC_w\lambda + \mC_w\lambda^N + \widetilde{\mC}_0\mC_\varrho\lambda^N + \mC_w\mC_\varrho\lambda^N \lesssim \mC_\upsilon\lambda^N.
	\end{aligned}
	\end{equation}
	
	We next estimate the material derivatives of the original velocity.
	For $M\geq1$, applying $D_{t,\upsilon}^M$ to
	$\varrho\upsilon=w$ gives
	\begin{equation}\label{eq:quotient-material-product-identity}
	\varrho D_{t,\upsilon}^M\upsilon = D_{t,\upsilon}^Mw - \upsilon D_{t,\upsilon}^M\varrho - \sum_{j=1}^{M-1}\binom{M}{j} D_{t,\upsilon}^j\varrho\,D_{t,\upsilon}^{M-j}\upsilon.
	\end{equation}
	We proceed by induction on $M$, estimating $D_{t,\upsilon}^M\upsilon$ for
	$M+N\leq K^*$ and both $D_{t,\upsilon}^M\nabla\upsilon$ and
	$\nabla D_{t,\upsilon}^M\upsilon$ for $M+N\leq K^*-1$.
	By \eqref{eq:quotient-material-product-identity}, the assumptions, and the induction hypothesis, we have
	\begin{equation}\label{eq:quotient-velocity-material-estimate}
	\begin{aligned}
		\nrm{D_{t,\upsilon}^M\upsilon}_{N,\mcI} &\lesssim \sum_{N_0+N_1=N} \nrm{\varrho^{-1}}_{N_0,\mcI} \Bigg( \nrm{D_{t,\upsilon}^Mw}_{N_1,\mcI}+ \sum_{N_2+N_3=N_1} \nrm{\upsilon}_{N_2,\mcI} \nrm{D_{t,\upsilon}^M(\varrho-\varrho_*)}_{N_3,\mcI}\\
		&\hspace{90pt} + \sum_{j=1}^{M-1} \sum_{N_2+N_3=N_1} \nrm{D_{t,\upsilon}^j(\varrho-\varrho_*)}_{N_2,\mcI} \nrm{D_{t,\upsilon}^{M-j}\upsilon}_{N_3,\mcI} \Bigg)\\
		&\lesssim \sum_{N_0+N_1=N}\lambda^{N_0} \Bigg( \mC_w\mu^M\lambda^{N_1} + \mC_\varrho\mu^M\lambda^{N_1} + \sum_{j=1}^{M-1} \mC_\varrho\mC_\upsilon\mu^M\lambda^{N_1} \Bigg)\lesssim \mC_\upsilon\mu^M\lambda^N.
	\end{aligned}
	\end{equation}
	Apply Lemma~\ref{lem:appendix-scale-dependent-commuted-gradient} with
	$G=\upsilon$, $K_\upsilon=K_G=K^\sharp=K^*$, $N_t=0$,
	$\widetilde\lambda=\lambda$, and $\tilde\mu=\mu$.
	Estimates~\eqref{est:appendix-scale-dependent-commuted-gradient}, \eqref{eq:quotient-velocity-spatial-estimate}, and \eqref{eq:quotient-velocity-material-estimate} give, for $M+N\leq K^*-1$,
	$$
	\nrm{D_{t,\upsilon}^M\nabla\upsilon}_{N,\mcI} + \nrm{\nabla D_{t,\upsilon}^M\upsilon}_{N,\mcI} \lesssim \mC_\upsilon\mu^M\lambda^{N+1}.
	$$
	
	We now turn to the mollified fields and prove the conclusions at $M=0$. By
	\eqref{est:appendix-spatial-mollification} with $N=0$ and the lower bound for
	$\overline\varrho$, we have
	$$
	\begin{aligned}
		\nrm{\overline\varrho-\varrho_*}_{0,\mcI}+\nrm{\overline w}_{0,\mcI}+\nrm{\overline\upsilon}_{0,\mcI}
		\lesssim \nrm{\varrho-\varrho_*}_{0,\mcI}+\nrm{w}_{0,\mcI}+\nrm{\overline w}_{0,\mcI}\nrm{\overline\varrho^{-1}}_{0,\mcI}
		\lesssim \widetilde{\mC}_0.
	\end{aligned}
	$$
	Let $H=\varrho-\varrho_*$ or $H=w$, and set
	$\mC_H=\mC_\varrho$ or $\mC_H=\mC_w$, respectively. For
	$1\leq N\leq K^\sharp$,
	the assumption at spatial order one gives
	$\nrm{H}_{0,\mcI}\leq\nrm{H}_{1,\mcI}\lesssim\mC_H\lambda$.
	Using \eqref{est:appendix-spatial-mollification-F-integer} for $\nabla H$, we obtain
	$$
	\begin{aligned}
		\nrm{P_{\leq\ell^{-1}}H}_{N,\mcI} &\lesssim \nrm{P_{\leq\ell^{-1}}H}_{0,\mcI} + \nrm{P_{\leq\ell^{-1}}\nabla H}_{N-1,\mcI}
		\lesssim \mC_H\lambda + \mC_H\lambda\MM{N-1,K^*-1,\lambda,\ell^{-1}}
		\lesssim \mC_H\MM{N,K^*,\lambda,\ell^{-1}}.
	\end{aligned}
	$$
	The Fa\`a di Bruno formula \eqref{eq:faa-di-bruno} yields, for
	$1\leq N\leq K^\sharp$,
	$$
	\begin{aligned}
		\nrm{\nabla^N(\overline\varrho^{-1})}_{0,\mcI}
		\lesssim \sum_{k=1}^N\sum_{\substack{N_1+\cdots+N_k=N\\N_j\geq1}}\prod_{j=1}^k\nrm{\nabla^{N_j}(\overline\varrho-\varrho_*)}_{0,\mcI}
		\lesssim \sum_{k=1}^N\mC_\varrho^k\MM{N,K^*,\lambda,\ell^{-1}}
		\lesssim \mC_\varrho\MM{N,K^*,\lambda,\ell^{-1}}.
	\end{aligned}
	$$
	The density lower bound and $\nrm{w}_{1,\mcI}\lesssim\mC_w\lambda$ also give
	$\nrm{\overline\upsilon}_{0,\mcI}\lesssim\mC_w\lambda$. Thus, for every
	$1\leq N\leq K^\sharp$,
	$$
	\begin{aligned}
		\nrm{\overline\upsilon}_{N,\mcI}
		\lesssim \sum_{N_1+N_2\leq N}\nrm{\nabla^{N_1}\overline w}_{0,\mcI}\nrm{\nabla^{N_2}(\overline\varrho^{-1})}_{0,\mcI}
		\lesssim\mC_w\MM{N,K^*,\lambda,\ell^{-1}}+\widetilde{\mC}_0\mC_\varrho\MM{N,K^*,\lambda,\ell^{-1}}
		\lesssim \mC_\upsilon\MM{N,K^*,\lambda,\ell^{-1}}.
	\end{aligned}
	$$
	These estimates prove \eqref{est:appendix-mollified-quotient-main}
	at $M=0$.
	The gradient estimate at $M=0$ follows by taking one additional spatial derivative.
	
	We next estimate the spatial mollification differences. For $0\leq N\leq K^*$, \eqref{est:appendix-spatial-mollification} with $L=K^*-N$ gives
	$$
	\begin{aligned}
		\nrm{\varrho-\overline\varrho}_{N,\mcI} \lesssim \ell^{K^*-N}\nrm{\varrho-\varrho_*}_{K^*,\mcI} \lesssim \mC_\varrho\ell^{K^*-N}\lambda^{K^*},\quad
		\nrm{w-\overline w}_{N,\mcI} \lesssim \ell^{K^*-N}\nrm{w}_{K^*,\mcI} \lesssim \mC_w\ell^{K^*-N}\lambda^{K^*}.
	\end{aligned}
	$$
	The velocity difference satisfies
	\begin{equation}
		\overline\upsilon-\upsilon
		=
		\frac{\overline w-w}{\overline\varrho}
		+
		w\frac{\varrho-\overline\varrho}{\varrho\overline\varrho}.
		\label{eq:appendix-mollified-quotient-velocity-difference}
	\end{equation}
	The preceding estimates for $\varrho^{-1}$ and $\overline\varrho^{-1}$ give, for every $0\leq N\leq K^*$,
	$$
	\nrm{\overline\varrho^{-1}}_{N,\mcI} + \nrm{(\varrho\overline\varrho)^{-1}}_{N,\mcI} \lesssim \lambda^N.
	$$
	Therefore, for every $0\leq N\leq K^*$,
	$$
	\begin{aligned}
		\nrm{\overline\upsilon-\upsilon}_{N,\mcI} &\lesssim \sum_{N_1+N_2=N} \nrm{\overline w-w}_{N_1,\mcI} \nrm{\overline\varrho^{-1}}_{N_2,\mcI}+ \sum_{N_1+N_2+N_3=N} \nrm{w}_{N_1,\mcI} \nrm{\varrho-\overline\varrho}_{N_2,\mcI} \nrm{(\varrho\overline\varrho)^{-1}}_{N_3,\mcI}\\
		&\lesssim \mC_\upsilon \sum_{j=0}^N\ell^{K^*-j}\lambda^{K^*}\lambda^{N-j}\lesssim \mC_\upsilon\ell^{K^*-N}\lambda^{K^*}.
	\end{aligned}
	$$
	Thus \eqref{est:appendix-mollified-quotient-high-frequency} holds.
	
	We now prove the estimates along $\overline\upsilon$ by induction on $M$.
	Fix $1\leq M\leq K^*$,
	and assume that the main estimates, the gradient estimate, and
	the material derivative estimates for the mollification differences in the range
	$m+N\leq K^*-1$ hold at every material order $0\leq m<M$. In particular,
	whenever
	$0\leq m<M$ and $m+N\leq K^*-1$,
	$$
	\begin{aligned}
		\nrm{D_{t,\overline\upsilon}^m(\varrho-\overline\varrho)}_{N,\mcI} &\lesssim \mC_\varrho\mu^m\ell^{K^*-m-N}\lambda^{K^*-m},\\
		\nrm{D_{t,\overline\upsilon}^m(w-\overline w)}_{N,\mcI} &\lesssim \mC_w\mu^m\ell^{K^*-m-N}\lambda^{K^*-m},\\
		\nrm{D_{t,\overline\upsilon}^m(\upsilon-\overline\upsilon)}_{N,\mcI} &\lesssim \mC_\upsilon\mu^m\ell^{K^*-m-N}\lambda^{K^*-m}.
	\end{aligned}
	$$
	For $m=0$, these estimates follow from
	\eqref{est:appendix-mollified-quotient-high-frequency}. For $m\geq1$,
	they are part of the induction hypothesis. Whenever
	$0\leq m<M$ and $m+N\leq K^\sharp-1$,
	$$
	\nrm{D_{t,\overline\upsilon}^m\nabla\overline\upsilon}_{N,\mcI}
	+
	\nrm{\nabla D_{t,\overline\upsilon}^m\overline\upsilon}_{N,\mcI}
	\lesssim
	\mC_\upsilon\mu^m\MM{N+1,K^*-m,\lambda,\ell^{-1}}.
	$$
	Since $\ell\lambda\leq1$, whenever $0\leq m<M$ and
	$m+N\leq K^*-1$,
	\begin{equation}\label{eq:appendix-quotient-material-change-coefficient-estimates}
	\begin{aligned}
		\nrm{D_{t,\overline\upsilon}^m(\upsilon-\overline\upsilon)}_{N,\mcI} &\lesssim \mC_\upsilon\mu^m\lambda^N (\ell\lambda)^{K^*-m-N} \lesssim \mC_\upsilon\mu^m\lambda^N,\\
		\nrm{D_{t,\overline\upsilon}^m\nabla\overline\upsilon}_{N,\mcI} + \nrm{\nabla D_{t,\overline\upsilon}^m\overline\upsilon}_{N,\mcI} &\lesssim \mC_\upsilon\mu^m\lambda^{N+1}.
	\end{aligned}
	\end{equation}
	
	We first change the material derivative on the unmollified fields.
	Take $H=\varrho-\varrho_*$ or $H=w$, with
	$\mC_H=\mC_\varrho$ or $\mC_H=\mC_w$, respectively. Estimates~\eqref{eq:appendix-quotient-material-change-coefficient-estimates} verify the velocity
	hypotheses of Lemma~\ref{lem:appendix-change-material-derivative} with
	$K=K^*$, $J=M$, $u=\upsilon$, $\overline u=\overline\upsilon$, and
	$\mC_u=\mC_\upsilon$. For the two choices of $H$, the lemma yields,
	whenever $M+N\leq K^*$,
	$$
	\begin{aligned}
		\nrm{D_{t,\overline\upsilon}^M(\varrho-\varrho_*)}_{N,\mcI} \lesssim \mC_\varrho\mu^M\lambda^N,\qquad
		\nrm{D_{t,\overline\upsilon}^Mw}_{N,\mcI} \lesssim \mC_w\mu^M\lambda^N.
	\end{aligned}
	$$
	For the same two choices of $H$,
	\eqref{est:appendix-change-material-derivative-gradient} gives, whenever
	$M+N\leq K^*-1$,
	$$
	\begin{aligned}
		\nrm{D_{t,\overline\upsilon}^M\nabla(\varrho-\varrho_*)}_{N,\mcI}
		+\nrm{\nabla D_{t,\overline\upsilon}^M(\varrho-\varrho_*)}_{N,\mcI}
		&\lesssim \mC_\varrho\mu^M\lambda^{N+1},\\
		\nrm{D_{t,\overline\upsilon}^M\nabla w}_{N,\mcI}
		+\nrm{\nabla D_{t,\overline\upsilon}^Mw}_{N,\mcI}
		&\lesssim \mC_w\mu^M\lambda^{N+1}.
	\end{aligned}
	$$
	
	For the mollified density and momentum, fix $N\in\N_0$ with $M+N\leq K^\sharp$. Applying
	\eqref{eq:material-mollifier-expansion} to $H$, we obtain
	\begin{equation}
		D_{t,\overline\upsilon}^MP_{\leq\ell^{-1}}H
		=
		\sum_{s=0}^M\binom{M}{s}
		\operatorname{ad}_{D_{t,\overline\upsilon}}^s(P_{\leq\ell^{-1}})
		D_{t,\overline\upsilon}^{M-s}H.
		\label{eq:appendix-new-material-mollifier-expansion}
	\end{equation}
	Lemma~\ref{lem:appendix-change-material-derivative} gives, for every $0\leq n\leq K^*-M$,
	$$
	\nrm{D_{t,\overline\upsilon}^MH}_{n,\mcI} \lesssim \mC_H\mu^M\lambda^n.
	$$
	By \eqref{est:appendix-spatial-mollification-F-integer}, the mollified term satisfies
	$$
	\nrm{P_{\leq\ell^{-1}}D_{t,\overline\upsilon}^MH}_{N,\mcI} \lesssim \mC_H\mu^M\MM{N,K^*-M,\lambda,\ell^{-1}}.
	$$
	For $1\leq s\leq M$, set
	$g:=D_{t,\overline\upsilon}^{M-s}H$. Then, for every
	$1\leq n\leq K^*-M+s$,
	$$
	\nrm{g}_{n,\mcI} \lesssim \mC_H\mu^{M-s}\lambda^n.
	$$
	Only spatial derivatives of $g$ of positive order occur in the commutator expansion.
	Every velocity factor in the commutator has material order at most
	$s-1<M$, and
	$
	N+s-1 \leq N+M-1 \leq K^\sharp-1.
	$
	Estimates~\eqref{est:general-multiple-mollifier-commutator} and
	\eqref{est:appendix-extended-MM-count} with the common spatial threshold
	$K^*-M$ give
	$$
	\begin{aligned}
		\nrm{ \operatorname{ad}_{D_{t,\overline\upsilon}}^s(P_{\leq\ell^{-1}}) D_{t,\overline\upsilon}^{M-s}H }_{N,\mcI}
		\lesssim \mC_H\mu^{M-s}\mu^s \MM{N,K^*-M,\lambda,\ell^{-1}}
		\lesssim \mC_H\mu^M\MM{N,K^*-M,\lambda,\ell^{-1}}.
	\end{aligned}
	$$
	Consequently,
	$$
	\begin{aligned}
		\nrm{D_{t,\overline\upsilon}^MP_{\leq\ell^{-1}}H}_{N,\mcI} &\lesssim \nrm{P_{\leq\ell^{-1}}D_{t,\overline\upsilon}^MH}_{N,\mcI}+ \sum_{s=1}^M \nrm{ \operatorname{ad}_{D_{t,\overline\upsilon}}^s(P_{\leq\ell^{-1}}) D_{t,\overline\upsilon}^{M-s}H }_{N,\mcI}
		\lesssim \mC_H\mu^M\MM{N,K^*-M,\lambda,\ell^{-1}}.
	\end{aligned}
	$$
	Taking $H=\varrho-\varrho_*$ and $H=w$ proves the first two estimates
	in \eqref{est:appendix-mollified-quotient-main} at material order $M$
	for $M+N\leq K^\sharp$.
	
	For $\overline\upsilon$, differentiating
	$\overline w=\overline\varrho\,\overline\upsilon$ gives
	\begin{equation}\label{eq:mollified-quotient-material-product-identity}
	\begin{aligned}
		\overline\varrho\,D_{t,\overline\upsilon}^M\overline\upsilon
		=D_{t,\overline\upsilon}^M\overline w-\overline\upsilon\,D_{t,\overline\upsilon}^M\overline\varrho
		-\sum_{j=1}^{M-1}\binom{M}{j}D_{t,\overline\upsilon}^j\overline\varrho\,D_{t,\overline\upsilon}^{M-j}\overline\upsilon.
	\end{aligned}
	\end{equation}
	If $K_j\geq K^*-M$ and $N_1+\cdots+N_J=N$, then
	$
	\prod_{j=1}^J\MM{N_j,K_j,\lambda,\ell^{-1}} \leq \MM{N,K^*-M,\lambda,\ell^{-1}}.
	$
	The density lower bound and the preceding spatial estimates for $\overline\varrho^{-1}$ give, for every $0\leq N\leq K^\sharp$,
	$$
	\nrm{\overline\varrho^{-1}}_{N,\mcI} \lesssim \MM{N,K^*,\lambda,\ell^{-1}}.
	$$
	Using these estimates in \eqref{eq:mollified-quotient-material-product-identity}, we obtain
	$$
	\begin{aligned}
		\nrm{D_{t,\overline\upsilon}^M\overline\upsilon}_{N,\mcI}
		&\lesssim
		\sum_{N_0+N_1=N}
		\nrm{\overline\varrho^{-1}}_{N_0,\mcI}
		\Bigg(
		\nrm{D_{t,\overline\upsilon}^M\overline w}_{N_1,\mcI}+
		\sum_{N_2+N_3=N_1}
		\nrm{\overline\upsilon}_{N_2,\mcI}
		\nrm{D_{t,\overline\upsilon}^M(\overline\varrho-\varrho_*)}_{N_3,\mcI}\\
		&\hspace{90pt}+
		\sum_{j=1}^{M-1}
		\sum_{N_2+N_3=N_1}
		\nrm{D_{t,\overline\upsilon}^j(\overline\varrho-\varrho_*)}_{N_2,\mcI}
		\nrm{D_{t,\overline\upsilon}^{M-j}\overline\upsilon}_{N_3,\mcI}
		\Bigg)\\
		&\lesssim
		\mC_w\mu^M\MM{N,K^*-M,\lambda,\ell^{-1}}+
		\mC_\varrho\mu^M\MM{N,K^*-M,\lambda,\ell^{-1}}+
		\sum_{j=1}^{M-1}
		\mC_\varrho\mC_\upsilon\mu^M
		\MM{N,K^*-M,\lambda,\ell^{-1}}\\
		&\lesssim
		\mC_\upsilon\mu^M\MM{N,K^*-M,\lambda,\ell^{-1}}.
	\end{aligned}
	$$
	Thus \eqref{est:appendix-mollified-quotient-main} holds at material
	order $M$.

	For the velocity gradients, apply Lemma~\ref{lem:appendix-scale-dependent-commuted-gradient} with
	$G=\overline\upsilon$, $K_\upsilon=K_G=K^*$, $N_t=0$,
	$\widetilde\lambda=\ell^{-1}$, and $\tilde\mu=\mu$.
	Use the velocity estimate in \eqref{est:appendix-mollified-quotient-main},
	already proved through material order $M$. The lemma then gives,
	for $M+N\leq K^\sharp-1$,
	$$
	\nrm{D_{t,\overline\upsilon}^M\nabla\overline\upsilon}_{N,\mcI}
	+
	\nrm{\nabla D_{t,\overline\upsilon}^M\overline\upsilon}_{N,\mcI}
	\lesssim
	\mC_\upsilon\mu^M\MM{N+1,K^*-M,\lambda,\ell^{-1}}.
	$$
	
	If $M=K^*$, the induction step is complete. To estimate the mollification differences, suppose
	$1\leq M\leq K^*-1$, fix $N\in\N_0$ with $M+N\leq K^*-1$, and set
	$L:=K^*-M-N$. Then $1\leq L\leq K^*\leq N_{\mathrm{mom}}$. From
	$P_{>\ell^{-1}}=I-P_{\leq\ell^{-1}}$,
	$$
	D_{t,\overline\upsilon}^MP_{>\ell^{-1}}H=P_{>\ell^{-1}}D_{t,\overline\upsilon}^MH-\sum_{s=1}^M\binom{M}{s}\operatorname{ad}_{D_{t,\overline\upsilon}}^s(P_{\leq\ell^{-1}})D_{t,\overline\upsilon}^{M-s}H.
	$$
	Estimate \eqref{est:appendix-spatial-mollification} gives
	$$
	\begin{aligned}
		\nrm{P_{>\ell^{-1}}D_{t,\overline\upsilon}^MH}_{N,\mcI} \lesssim \ell^L\nrm{D_{t,\overline\upsilon}^MH}_{N+L,\mcI}
		\lesssim \mC_H\mu^M\ell^L\lambda^{N+L}
		= \mC_H\mu^M\ell^{K^*-M-N}\lambda^{K^*-M}.
	\end{aligned}
	$$
	For $1\leq s\leq M$,
	$
	N+\max\{s,L\}\leq K^*-M+s, \ N+s+L-2\leq K^*-2.
	$
	Every velocity factor has material order at most $s-1<M$.
	Apply \eqref{est:higher-order-iterated-material-mollifier} with
	$g=D_{t,\overline\upsilon}^{M-s}H$, $\mC_g=\mC_H\mu^{M-s}$,
	$K_\upsilon=1$, $K_f=0$, $\lambda_\upsilon=\lambda_f=\lambda$, and
	$\mu_\upsilon=\tilde\mu_\upsilon=\mu_f=\tilde\mu_f=\mu$. We obtain
	$$
	\begin{aligned}
		\nrm{ \operatorname{ad}_{D_{t,\overline\upsilon}}^s(P_{\leq\ell^{-1}}) D_{t,\overline\upsilon}^{M-s}H }_{N,\mcI}
		\lesssim \mC_H\mu^{M-s}\mu^s\ell^L\lambda^{N+L}
		\lesssim \mC_H\mu^M\ell^{K^*-M-N}\lambda^{K^*-M}.
	\end{aligned}
	$$
	Summing these estimates and taking $H=\varrho-\varrho_*$ and $H=w$
	proves the first two estimates in
	\eqref{est:appendix-mollified-quotient-material-high-frequency}
	at material order $M$ for $M+N\leq K^*-1$.
	
	For the velocity difference, we first estimate the reciprocal densities. For $\widetilde\varrho=\varrho$ or $\widetilde\varrho=\overline\varrho$
	and $j\geq1$, repeated differentiation gives
	\begin{equation}\label{eq:quotient-reciprocal-material-identity}
	D_{t,\overline\upsilon}^j(\widetilde\varrho^{-1})
	=
	\sum_{k=1}^j\widetilde\varrho^{-k-1}
	\sum_{\substack{\sigma_1+\cdots+\sigma_k=j\\\sigma_i\geq1}}
	c_{\sigma_1,\ldots,\sigma_k}
	\prod_{i=1}^kD_{t,\overline\upsilon}^{\sigma_i}\widetilde\varrho.
	\end{equation}
	For $1\leq j\leq M$ and $j+N\leq K^\sharp$,
	\begin{equation}\label{eq:mollified-density-reciprocal-material-estimate}
	\begin{aligned}
		\nrm{D_{t,\overline\upsilon}^j(\overline\varrho^{-1})}_{N,\mcI}
		&\lesssim
		\sum_{k=1}^j
		\sum_{\substack{\sigma_1+\cdots+\sigma_k=j\\\sigma_i\geq1}}
		\sum_{N_0+\cdots+N_k=N}
		\nrm{\overline\varrho^{-k-1}}_{N_0,\mcI}
		\prod_{i=1}^k
		\nrm{D_{t,\overline\upsilon}^{\sigma_i}\overline\varrho}_{N_i,\mcI}\\
		&\lesssim
		\sum_{k=1}^j
		\mu^j
		\sum_{N_0+\cdots+N_k=N}
		\MM{N_0,K^*,\lambda,\ell^{-1}}
		\prod_{i=1}^k\MM{N_i,K^*-j,\lambda,\ell^{-1}}\\
		&\lesssim
		\mu^j\MM{N,K^*-j,\lambda,\ell^{-1}}.
	\end{aligned}
	\end{equation}

	For $1\leq j\leq M$ and $j+N\leq K^*$,
	\begin{equation}\label{eq:density-reciprocal-changed-material-estimate}
	\begin{aligned}
		\nrm{D_{t,\overline\upsilon}^j(\varrho^{-1})}_{N,\mcI}
		\lesssim
		\sum_{k=1}^j
		\sum_{\substack{\sigma_1+\cdots+\sigma_k=j\\\sigma_i\geq1}}
		\sum_{N_0+\cdots+N_k=N}
		\nrm{\varrho^{-k-1}}_{N_0,\mcI}
		\prod_{i=1}^k
		\nrm{D_{t,\overline\upsilon}^{\sigma_i}\varrho}_{N_i,\mcI}
		\lesssim
		\sum_{k=1}^j
		\sum_{N_0+\cdots+N_k=N}
		\lambda^{N_0}
		\mC_\varrho^k\mu^j
		\prod_{i=1}^k\lambda^{N_i}
		\lesssim
		\mu^j\lambda^N.
	\end{aligned}
	\end{equation}
	The estimates at $j=0$ follow from the spatial derivative estimates for the reciprocal densities. For
	$0\leq j\leq M$ and $j+N\leq K^*$, multiplying these estimates gives
	\begin{equation}\label{eq:density-reciprocal-product-material-estimate}
	\begin{aligned}
		\nrm{D_{t,\overline\upsilon}^j((\varrho\overline\varrho)^{-1})}_{N,\mcI}
		\lesssim \sum_{\substack{j_1+j_2=j\\N_1+N_2=N}}\nrm{D_{t,\overline\upsilon}^{j_1}(\varrho^{-1})}_{N_1,\mcI}\nrm{D_{t,\overline\upsilon}^{j_2}(\overline\varrho^{-1})}_{N_2,\mcI}
		\lesssim \sum_{\substack{j_1+j_2=j\\N_1+N_2=N}}\mu^{j_1}\lambda^{N_1}\mu^{j_2}\lambda^{N_2}\lesssim \mu^j\lambda^N.
	\end{aligned}
	\end{equation}
	If $1\leq j\leq M$ and $j+N\leq K^*$, every term contains a derivative
	of $\varrho-\varrho_*$ or $\overline\varrho-\varrho_*$ of positive order.
	Each such derivative has the amplitude factor $\mC_\varrho$, so $\nrm{D_{t,\overline\upsilon}^j(\varrho^{-1})}_{N,\mcI} + \nrm{D_{t,\overline\upsilon}^j((\varrho\overline\varrho)^{-1})}_{N,\mcI} \lesssim \mC_\varrho\mu^j\lambda^N$.
	
	Applying $D_{t,\overline\upsilon}^M$ to
	\eqref{eq:appendix-mollified-quotient-velocity-difference} and taking
	the $C^N$ norm gives
	$$ 
	\begin{aligned} 
		\nrm{D_{t,\overline\upsilon}^M(\overline\upsilon-\upsilon)}_{N,\mcI} 
		&\lesssim 
		\sum_{M_1+M_2=M}\sum_{N_1+N_2=N} 
		\nrm{D_{t,\overline\upsilon}^{M_1}(\overline w-w)}_{N_1,\mcI} 
		\nrm{D_{t,\overline\upsilon}^{M_2}(\overline\varrho^{-1})}_{N_2,\mcI}\\ 
		&\quad+ 
		\sum_{M_1+M_2+M_3=M}\sum_{N_1+N_2+N_3=N} 
		\nrm{D_{t,\overline\upsilon}^{M_1}w}_{N_1,\mcI} 
		\nrm{D_{t,\overline\upsilon}^{M_2}(\varrho-\overline\varrho)}_{N_2,\mcI} 
		\nrm{D_{t,\overline\upsilon}^{M_3}((\varrho\overline\varrho)^{-1})}_{N_3,\mcI}. 
	\end{aligned} 
	$$ 
	For either $H=\varrho-\varrho_*$ or $H=w$, let
	$\mC_H=\mC_\varrho$ or $\mC_H=\mC_w$, respectively. Combine the
	first two estimates in \eqref{est:appendix-mollified-quotient-high-frequency} with those in
	\eqref{est:appendix-mollified-quotient-material-high-frequency}, already proved for the density and momentum through order $M$. For $0\leq j\leq M$ and $j+n\leq M+N$, we obtain
	$$
	\begin{aligned}
		\nrm{D_{t,\overline\upsilon}^j(H-P_{\leq\ell^{-1}}H)}_{n,\mcI}
		\lesssim
		\mC_H\mu^j\ell^{K^*-j-n}\lambda^{K^*-j}
		=
		\mC_H\mu^j\ell^L\lambda^{n+L}(\ell\lambda)^{M+N-j-n}
		\lesssim
		\mC_H\mu^j\ell^L\lambda^{n+L}.
	\end{aligned}
	$$
	Using $L=K^*-M-N$, $\ell\lambda\leq1$, and $M+N\leq K^*-1$,
	we estimate the remaining factors in their stated ranges. If $M_1+M_2=M$ and
	$N_1+N_2=N$, then
	$$ 
	\begin{aligned} 
		\nrm{D_{t,\overline\upsilon}^{M_1}(\overline w-w)}_{N_1,\mcI} 
		\nrm{D_{t,\overline\upsilon}^{M_2}(\overline\varrho^{-1})}_{N_2,\mcI}\lesssim 
		\mC_w\mu^{M_1}\ell^L\lambda^{N_1+L}
		\mu^{M_2}\lambda^{N_2}
		=
		\mC_w\mu^M\ell^L\lambda^{N+L}. 
	\end{aligned} 
	$$ 
	Similarly, if $M_1+M_2+M_3=M$ and $N_1+N_2+N_3=N$, then
	$$ 
	\begin{aligned} 
		\nrm{D_{t,\overline\upsilon}^{M_1}w}_{N_1,\mcI} 
		\nrm{D_{t,\overline\upsilon}^{M_2}(\varrho-\overline\varrho)}_{N_2,\mcI} 
		\nrm{D_{t,\overline\upsilon}^{M_3}((\varrho\overline\varrho)^{-1})}_{N_3,\mcI}\lesssim 
		\mu^{M_1}\lambda^{N_1}
		\mC_\varrho\mu^{M_2}\ell^L\lambda^{N_2+L}
		\mu^{M_3}\lambda^{N_3}
		=
		\mC_\varrho\mu^M\ell^L\lambda^{N+L}. 
	\end{aligned} 
	$$ 
	Summing and using $\mC_w+\mC_\varrho\lesssim\mC_\upsilon$ proves
	the last estimate in
	\eqref{est:appendix-mollified-quotient-material-high-frequency}
	for $M+N\leq K^*-1$ and closes the induction on $M$.
	
	With the main estimates proved through material order $K^*$,
	identity~\eqref{eq:quotient-reciprocal-material-identity} proves \eqref{eq:mollified-density-reciprocal-material-estimate}--\eqref{eq:density-reciprocal-product-material-estimate} for every
	$0\leq j\leq K^*$ in their stated total-order ranges.
	It remains to treat the endpoint $M+N=K^*$ with $M\geq1$. Let
	$H=\varrho-\varrho_*$ or $H=w$, set
	$\mC_H=\mC_\varrho$ or $\mC_H=\mC_w$, respectively, and write
	$\overline H=P_{\leq\ell^{-1}}H$. Estimates~\eqref{est:appendix-change-material-derivative}
	and \eqref{est:appendix-mollified-quotient-main} give
	$$
	\begin{aligned}
		\nrm{D_{t,\overline\upsilon}^M(H-\overline H)}_{N,\mcI}
		&\lesssim
		\nrm{D_{t,\overline\upsilon}^MH}_{N,\mcI}
		+
		\nrm{D_{t,\overline\upsilon}^M\overline H}_{N,\mcI}
		\lesssim
		\mC_H\mu^M\lambda^N,
	\end{aligned}
	$$
	where $\MM{N,K^*-M,\lambda,\ell^{-1}}=\lambda^N$.
	Since $M\geq1$, every decomposition $M_1+M_2=M$ has
	$M_1\geq1$ or $M_2\geq1$. The corresponding factor contributes
	$\mC_w$ or $\mC_\varrho$, respectively. Therefore
	$$
	\begin{aligned}
		\nrm{D_{t,\overline\upsilon}^M\upsilon}_{N,\mcI}
		&\lesssim
		\sum_{M_1+M_2=M}\sum_{N_1+N_2=N}
		\nrm{D_{t,\overline\upsilon}^{M_1}w}_{N_1,\mcI}
		\nrm{D_{t,\overline\upsilon}^{M_2}(\varrho^{-1})}_{N_2,\mcI}
		\lesssim
		\mC_\upsilon\mu^M\lambda^N.
	\end{aligned}
	$$
	Together with
	\eqref{est:appendix-mollified-quotient-main}, we have
	$$
	\nrm{D_{t,\overline\upsilon}^M(\upsilon-\overline\upsilon)}_{N,\mcI}
	\lesssim
	\mC_\upsilon\mu^M\lambda^N.
	$$
	Since
	$\ell^{K^*-M-N}\lambda^{K^*-M}=\lambda^N$ when $M+N=K^*$, these are
	the endpoint cases of
	\eqref{est:appendix-mollified-quotient-material-high-frequency}.
	
	For the divergence estimate, suppose that $\partial_t\varrho+\Div w=0$. Then
	$
	\partial_t\overline\varrho+\Div\overline w=0, \ \Div\overline\upsilon = -\overline\varrho^{-1}D_{t,\overline\upsilon}\overline\varrho.
	$
	Hence
	$$
	D_{t,\overline\upsilon}^M\Div\overline\upsilon = - \sum_{j=0}^M\binom{M}{j} D_{t,\overline\upsilon}^j(\overline\varrho^{-1}) D_{t,\overline\upsilon}^{M+1-j}\overline\varrho.
	$$
	For $M\leq K^*-1$ and $M+N\leq K^\sharp-1$,
	$$
	\begin{aligned}
		\nrm{D_{t,\overline\upsilon}^M\Div\overline\upsilon}_{N,\mcI}
		\lesssim
		\sum_{j=0}^M
		\sum_{N_1+N_2=N}
		\nrm{D_{t,\overline\upsilon}^j(\overline\varrho^{-1})}_{N_1,\mcI}
		\nrm{D_{t,\overline\upsilon}^{M+1-j}\overline\varrho}_{N_2,\mcI}
		\lesssim
		\mC_\varrho\mu^{M+1}
		\MM{N,K^*-M-1,\lambda,\ell^{-1}}.
	\end{aligned}
	$$
	
	Finally, we use spatial interpolation for the H{\"o}lder estimates. Let $\alpha\in(0,1)$. Take $H=\overline\varrho-\varrho_*$,
	$\overline w$, or $\overline\upsilon$.
	Set $\mC_H=\mC_\varrho$, $\mC_w$, or $\mC_\upsilon$, respectively.
	Interpolation gives, for
	$M\leq K^*$ and $1\leq M+N\leq K^\sharp-1$,
	$$
	\begin{aligned}
		\nrm{D_{t,\overline\upsilon}^MH}_{N+\alpha,\mcI}
		\lesssim
		\nrm{D_{t,\overline\upsilon}^MH}_{N,\mcI}^{1-\alpha}
		\nrm{D_{t,\overline\upsilon}^MH}_{N+1,\mcI}^{\alpha}
		\lesssim
		\mC_H\mu^M
		\La{N,K^*-M,\lambda,\ell^{-1}}
		\MM{N,K^*-M,\lambda,\ell^{-1}}.
	\end{aligned}
	$$
	For $M\leq K^*$ and $M+N\leq K^\sharp-2$, the gradient estimates give
	$$
	\begin{aligned}
		\nrm{D_{t,\overline\upsilon}^M\nabla\overline\upsilon}_{N+\alpha,\mcI}
		+
		\nrm{\nabla D_{t,\overline\upsilon}^M\overline\upsilon}_{N+\alpha,\mcI}
		&\lesssim
		\Bigl(
		\nrm{D_{t,\overline\upsilon}^M\nabla\overline\upsilon}_{N,\mcI}
		+
		\nrm{\nabla D_{t,\overline\upsilon}^M\overline\upsilon}_{N,\mcI}
		\Bigr)^{1-\alpha}
		\Bigl(
		\nrm{D_{t,\overline\upsilon}^M\nabla\overline\upsilon}_{N+1,\mcI}
		+
		\nrm{\nabla D_{t,\overline\upsilon}^M\overline\upsilon}_{N+1,\mcI}
		\Bigr)^\alpha\\
		&\lesssim
		\mC_\upsilon\mu^M
		\La{N+1,K^*-M,\lambda,\ell^{-1}}
		\MM{N+1,K^*-M,\lambda,\ell^{-1}}.
	\end{aligned}
	$$
	Under the continuity equation, the same interpolation proves the
	divergence estimate for
	$M\leq K^*-1$ and $M+N\leq K^\sharp-2$.
	
	Take $H=\varrho$, $w$, or $\upsilon$, with $\overline H=\overline\varrho$,
	$\overline w$, or $\overline\upsilon$, respectively.
	Use the corresponding amplitude $\mC_H=\mC_\varrho$, $\mC_w$, or $\mC_\upsilon$.
	Interpolation between the $C^N$ and $C^{N+1}$ estimates gives, for
	$0\leq M\leq K^*-1$ and $M+N\leq K^*-1$,
	$$
	\begin{aligned}
		\nrm{D_{t,\overline\upsilon}^M(H-\overline H)}_{N+\alpha,\mcI}
		\lesssim
		\nrm{D_{t,\overline\upsilon}^M(H-\overline H)}_{N,\mcI}^{1-\alpha}
		\nrm{D_{t,\overline\upsilon}^M(H-\overline H)}_{N+1,\mcI}^{\alpha}
		\lesssim
		\mC_H\mu^M
		\ell^{K^*-M-N-\alpha}\lambda^{K^*-M}.
	\end{aligned}
	$$

	For the $C^\alpha$ estimate at $M=0$, take
	$G=\overline\varrho-\varrho_*$, $\overline w$, or $\overline\upsilon$, with
	$\mC_G=\mC_\varrho$, $\mC_w$, or $\mC_\upsilon$, respectively. Then
	\begin{align*}
		\nrm{G}_{\alpha,\mcI}
		&\lesssim
		\nrm{G}_{0,\mcI}^{1-\alpha}
		\nrm{G}_{1,\mcI}^{\alpha}
		\lesssim
		\widetilde{\mC}_0^{1-\alpha}\mC_G^\alpha\lambda^\alpha
		\lesssim
		\mC_G^\alpha\lambda^\alpha.\qedhere
	\end{align*}
	
\end{proof}
	\section{Well-posedness and estimates for the linearized Euler equations}

We prove well-posedness and estimates for the linearized equation
used in the Newton step. We adapt the Picard iteration in the proof of
\cite[Proposition~E.1]{GR24} to a background field that may have nonzero divergence.
In three dimensions, we use a zero-mean vector potential in place of the
two-dimensional stream function estimated in \cite[Lemma~3.7]{GR24}.
Let $I\subset\mathbb R$ be a nondegenerate compact interval and let $t_0\in I$. For
$t\in I$, we denote by
$
I(t,t_0):=[\min\{t,t_0\},\max\{t,t_0\}]
$
the interval with endpoints $t$ and $t_0$.

On $I\times\mathbb T^3$, consider
\begin{equation}\label{eq:linearized-non-divfree}
	\begin{cases}
		\partial_t \omega+\Div (U\otimes \omega+\omega\otimes U)+\nabla \pi=F,\\
		\Div\omega=0,\\
		\omega|_{t=t_0}=\omega_0 .
	\end{cases}
\end{equation}
Here $U$ and $F$ are prescribed vector fields, with $\Div F=0$.
The unknowns are the divergence-free vector field $\omega$ and the scalar
pressure $\pi$.
Equivalently, since $\Div \omega=0$, the first equation in
\eqref{eq:linearized-non-divfree} can be written as
\begin{equation}\label{eq:linearized-expanded-intro}
	\partial_t\omega
	+U\cdot\nabla\omega
	+\omega\cdot\nabla U
	+\omega\,\Div U
	+\nabla\pi
	=F .
\end{equation}
Taking the divergence of the first equation in
\eqref{eq:linearized-non-divfree} and using $\Div \omega=\Div F=0$,
we obtain the elliptic equation
\begin{equation}\label{eq:pressure-linearized}
	-\Delta \pi
	=
	\Div \Div (U\otimes \omega+\omega\otimes U).
\end{equation}
With the normalization
$
\fint_{\mathbb T^3}\pi(t,x)\,\rd x=0
$
for each $t\in I$, the pressure is uniquely determined by
\eqref{eq:pressure-linearized}.

We write $\mathbb P_H$ for the periodic Leray projection onto divergence-free,
zero-mean vector fields. For a zero-mean vector field $h$ on $\mathbb T^3$,
define the periodic vector-potential operator
\begin{equation}\label{eq:vector-potential-operator}
	\mathcal B h:=-\Delta^{-1}\curl h.
\end{equation}
Then
$
\Div\mathcal B h=0
$
and
$
\fint_{\mathbb T^3}\mathcal B h\,\rd x=0.
$
Moreover,
$
\curl\mathcal B h=\mathbb P_H h.
$
In particular, if $\Div h=0$, then
$
\curl\mathcal B h=h.
$
For such $h$, $\mathcal B h$ is the unique zero-mean divergence-free
vector potential. Indeed, the difference of two such potentials has zero
divergence, curl, and mean, and hence vanishes.

\begin{prop}\label{prop:linearized-equation}
	Fix a nondegenerate compact interval $I\subset\mathbb R$, $t_0\in I$,
	an integer $N\geq1$, and $0<\alpha<1$. Let
	$
	U\in C(I;C^{N+1+\alpha}(\mathbb T^3))
	$
	and
	$
	F\in C(I;C^{N+\alpha}(\mathbb T^3)).
	$
	Assume that $F$ is divergence-free and has zero spatial mean. Let
	$
	\omega_0\in C^{N+\alpha}(\mathbb T^3)
	$
	be divergence-free and have zero spatial mean.
	Assume also that $U$, $F$, and $\omega_0$ admit smooth spatial approximations
	converging in the stated norms, uniformly in time for $U$ and $F$.
	Then there exists a solution $(\omega,\pi)$ of
	\begin{equation}\label{eq:linearized-equation}
		\left\{
		\begin{aligned}
			&\partial_t\omega+\Div(U\otimes\omega+\omega\otimes U)+\nabla\pi=F,\\
			&\Div\omega=0,\\
			&\omega|_{t=t_0}=\omega_0.
		\end{aligned}
		\right.
	\end{equation}
	The solution satisfies
	$
	\omega\in C(I;C^{N+\alpha}(\mathbb T^3))\cap C^1(I;C^{N-1+\alpha}(\mathbb T^3)).
	$
	The field $\omega$ is unique in this class. With the normalization
	$
	\fint_{\mathbb T^3}\pi(t,x)\,\rd x=0,
	$
	the pressure is unique and
	$
	\pi\in C(I;C^{N+\alpha}(\mathbb T^3)).
	$
	
	For $t\geq t_0$, one has the low-order H{\"o}lder estimate
	\begin{equation}\label{eq:holder-low-estimate}
		\|\omega(t,\cdot)\|_{\alpha}
		\leq
		C_{\alpha}\exp\left(C_\alpha\int_{t_0}^{t}\|\nabla U(\tau,\cdot)\|_\alpha\,\rd\tau\right)
		\left(
		\|\omega_0\|_{\alpha}
		+\int_{t_0}^{t}\|F(\tau,\cdot)\|_{\alpha}\,\rd\tau
		\right).
	\end{equation}
	For the same times, the high-order H{\"o}lder estimate is
	\begin{equation}\label{eq:holder-estimate}
		\begin{aligned}
			\|\omega(t,\cdot)\|_{N+\alpha}
			\leq{}&
			C_{N,\alpha}\exp\left(C_{N,\alpha}\int_{t_0}^{t}\|\nabla U(\tau,\cdot)\|_\alpha\,\rd\tau\right)\\
			&\cdot
			\left(
			\|\omega_0\|_{N+\alpha}
			+\int_{t_0}^{t}\|F(\tau,\cdot)\|_{N+\alpha}\,\rd\tau
			+\int_{t_0}^{t}
			\|\nabla U(\tau,\cdot)\|_{N+\alpha}
			\left(
			\|\omega_0\|_{\alpha}
			+\int_{t_0}^{\tau}\|F(s,\cdot)\|_{\alpha}\,\rd s
			\right)\rd\tau
			\right).
		\end{aligned}
	\end{equation}
	For $t\leq t_0$, replace $\int_{t_0}^{t}$ by $\int_t^{t_0}$ and
	$\int_{t_0}^{\tau}$ by $\int_\tau^{t_0}$ in these estimates.
	
	Set $\Psi:=\mathcal B\omega$, so that $\curl\Psi=\omega$. Then
	\begin{equation}\label{eq:vector-potential-holder-estimate}
		\|\Psi(t,\cdot)\|_{N+1+\alpha}
		\leq C_{N,\alpha}\|\omega(t,\cdot)\|_{N+\alpha}.
	\end{equation}
	In addition, for every $t\in I$, the vector potential satisfies
	\begin{equation}\label{eq:besov-estimate}
		\begin{aligned}
			\|\Psi(t,\cdot)\|_{\alpha}
			\leq{}&
			C_{\alpha}\exp\left(
			C_{\alpha}\int_{I(t,t_0)}
			\|\nabla U(\tau,\cdot)\|_{\alpha}\,\rd\tau
			\right)
			\left(
			\|\mathcal B\omega_0\|_{\alpha}
			+\int_{I(t,t_0)}
			\|\mathcal B F(\tau,\cdot)\|_{\alpha}\,\rd\tau
			\right).
		\end{aligned}
	\end{equation}
	
	Finally, the normalized pressure satisfies
	\begin{align}
		\|\pi(t,\cdot)\|_{N+\alpha}
		\leq& C_{N,\alpha}\left(
		\|U(t,\cdot)\|_{N+\alpha}\|\omega(t,\cdot)\|_{\alpha}
		+\|U(t,\cdot)\|_{\alpha}\|\omega(t,\cdot)\|_{N+\alpha}
		\right),\label{eq:pressure-holder-estimate}\\
		\|\pi(t,\cdot)\|_{N+\alpha}
		\leq& C_{N,\alpha}\big(
		\|\nabla U(t,\cdot)\|_{\alpha}
		\|\Psi(t,\cdot)\|_{N+\alpha}+\|\nabla U(t,\cdot)\|_{N+\alpha}
		\|\Psi(t,\cdot)\|_{\alpha}
		\big).\label{eq:pressure-besov-estimate}
	\end{align}
	The constants $C_\alpha$ depend only on $\alpha$, and $C_{N,\alpha}$ only on $N$ and $\alpha$.
	They are independent of $I$, $t_0$, $U$, $F$, and $\omega_0$.
\end{prop}

\begin{proof}
	By time reversal, it suffices to consider $t\geq t_0$.
	We first construct $\omega$ by Picard iteration and extend it to $I$.
	We then estimate $\omega$, $\mathcal B\omega$, and $\pi$.
	Integrating \eqref{eq:linearized-equation} in space shows that any solution
	with the prescribed initial data has zero spatial mean. Applying $\mathbb P_H$
	and using $\mathbb P_H F=F$, we obtain
	\begin{equation}\label{eq:projected-linearized}
		\partial_t\omega+\mathbb P_H\Div(U\otimes\omega+\omega\otimes U)=F .
	\end{equation}
	Since $\Div\omega=0$, this is equivalent to
	\begin{equation}\label{eq:projected-transport-form}
		\partial_t\omega+U\cdot\nabla\omega
		=F-\mathbb P_H(\omega\cdot\nabla U+\omega\Div U)-[\mathbb P_H,U\cdot\nabla]\omega .
	\end{equation}

	Constant-coefficient directional derivatives commute with $\mathbb P_H$.
	For the mean projection, integration by parts gives
	$\fint_{\mathbb T^3}U\cdot\nabla h\,\rd x=-\fint_{\mathbb T^3}(\Div U)h\,\rd x$.
	For every vector field $h\in C^{N+\alpha}(\mathbb T^3)$, apply
	Proposition~\ref{prop:CZ commutator} to $U$ minus its spatial mean.
	Using Poincar\'e's inequality, \eqref{eq:holder-product-tame}, and the estimate for
	the mean-projection commutator, we obtain
	\begin{align}
		\|\mathbb P_H(h\cdot\nabla U)\|_{\alpha}
		+\|\mathbb P_H(h\Div U)+[\mathbb P_H,U\cdot\nabla]h\|_{\alpha}
		&\leq C_\alpha\|\nabla U\|_{\alpha}\|h\|_{\alpha},
		\label{eq:LA-low-estimate}\\
		\|\mathbb P_H(h\cdot\nabla U)\|_{N+\alpha}
		+\|\mathbb P_H(h\Div U)+[\mathbb P_H,U\cdot\nabla]h\|_{N+\alpha}&\leq C_{N,\alpha}\|\nabla U\|_{\alpha}\|h\|_{N+\alpha}
		+C_{N,\alpha}\|\nabla U\|_{N+\alpha}\|h\|_{\alpha}.
		\label{eq:LA-high-estimate}
	\end{align}

	We solve on short intervals $J=[t_-,t_+]\subset I$.
	Let $\omega_-$ denote the initial value on $J$, starting with
	$\omega_0$ at $t_0$.
	Set $\omega^{(0)}(t,\cdot)=\omega_-$ and, for $m\geq0$, define
	$\omega^{(m+1)}$ by
	\begin{equation}\label{eq:omega-picard-iteration}
		\left\{
		\begin{aligned}
			&\partial_t\omega^{(m+1)}+U\cdot\nabla\omega^{(m+1)}
			+\mathbb P_H(\omega^{(m+1)}\Div U)
			+[\mathbb P_H,U\cdot\nabla]\omega^{(m+1)}
			=F-\mathbb P_H(\omega^{(m)}\cdot\nabla U),\\
			&\omega^{(m+1)}|_{t=t_-}=\omega_- .
		\end{aligned}
		\right.
	\end{equation}
	By the approximation assumption, the transport operators are strongly
	continuous on the closure of smooth fields in $C^{N+\alpha}$.
	The other operators in \eqref{eq:omega-picard-iteration} preserve this
	space by smooth approximation and are bounded there by \eqref{eq:LA-high-estimate}.
	For each $m$, solve the integral equation along the flow of $U$ by
	contraction for sufficiently small $|J|$. We obtain a unique solution
	$
	\omega^{(m+1)}\in C(J;C^{N+\alpha})\cap C^1(J;C^{N-1+\alpha})
	$.
	
	We verify that the divergence-free and zero-mean conditions are preserved.
	Applying $I-\mathbb P_H$ to \eqref{eq:omega-picard-iteration} and combining
	the transport and commutator terms gives
	$$
	\begin{aligned}
		\partial_t(I-\mathbb P_H)\omega^{(m+1)}
		+U\cdot\nabla(I-\mathbb P_H)\omega^{(m+1)}
		-[\mathbb P_H,U\cdot\nabla](I-\mathbb P_H)\omega^{(m+1)}
		=0,
	\end{aligned}
	$$
	with $(I-\mathbb P_H)\omega^{(m+1)}|_{t=t_-}=0$.
	By the commutator estimate used to prove \eqref{eq:LA-high-estimate},
	the same short-time argument shows that
	$
	(I-\mathbb P_H)\omega^{(m+1)}=0.
	$
	Thus $\omega^{(m+1)}$ is divergence-free and has zero spatial mean. Equation
	\eqref{eq:omega-picard-iteration} is equivalent to
	$$
	\partial_t\omega^{(m+1)}
	+\mathbb P_H\Div(\omega^{(m+1)}\otimes U+U\otimes\omega^{(m)})=F .
	$$
	
	Let
	$
	\delta\omega^{(m+1)}=\omega^{(m+1)}-\omega^{(m)}.
	$
	Then, for $m\geq1$,
	$$
	\left\{
	\begin{aligned}
		&\partial_t\delta\omega^{(m+1)}+U\cdot\nabla\delta\omega^{(m+1)}
		+\mathbb P_H(\delta\omega^{(m+1)}\Div U)
		+[\mathbb P_H,U\cdot\nabla]\delta\omega^{(m+1)}
		=-\mathbb P_H(\delta\omega^{(m)}\cdot\nabla U),\\
		&\delta\omega^{(m+1)}|_{t=t_-}=0 .
	\end{aligned}
	\right.
	$$
	Take $|J|\|U\|_{C(J;C^{N+1+\alpha})}\leq1$ and apply
	Lemma~\ref{estimate of the transport equation} componentwise by smooth
	spatial approximation. Using \eqref{eq:LA-low-estimate}--%
	\eqref{eq:LA-high-estimate}, we obtain
	$$
	\begin{aligned}
		\|\delta\omega^{(m+1)}\|_{C(J;C^{N+\alpha})}
		\leq{}&
		C_{N,\alpha}|J|
		\exp\left(C_{N,\alpha}\int_J\|\nabla U(\tau,\cdot)\|_\alpha\,\rd\tau\right)
		\|U\|_{C(J;C^{N+1+\alpha})}
		\left(
		\|\delta\omega^{(m)}\|_{C(J;C^{N+\alpha})}
		+\|\delta\omega^{(m+1)}\|_{C(J;C^{N+\alpha})}
		\right).
	\end{aligned}
	$$
	Choose $J$ so that the coefficient of the sum on the right is at most
	$1/4$. We absorb the last term on the right and conclude that
	$$
	\|\delta\omega^{(m+1)}\|_{C(J;C^{N+\alpha})}
	\leq
	\frac12
	\|\delta\omega^{(m)}\|_{C(J;C^{N+\alpha})}.
	$$
	Hence $\{\omega^{(m)}\}_{m\geq0}$ is Cauchy in $C(J;C^{N+\alpha})$.
	Let $\omega$ be its limit. Pass to the limit in
	\eqref{eq:omega-picard-iteration} to obtain \eqref{eq:projected-transport-form}.
	Equation~\eqref{eq:omega-picard-iteration} also gives convergence of
	$\partial_t\omega^{(m+1)}$ in $C(J;C^{N-1+\alpha})$. Thus
	$
	\omega\in C(J;C^{N+\alpha})\cap C^1(J;C^{N-1+\alpha}).
	$
	The limit is divergence-free, has zero spatial mean, and belongs to the
	closure of smooth fields in $C^{N+\alpha}$.
	
	For uniqueness, let $\omega_1$ and $\omega_2$ be two solutions on $J$
	with the same initial data. Their difference
	$\delta\omega=\omega_1-\omega_2$ satisfies
	\eqref{eq:projected-transport-form} with zero forcing and zero initial data.
	Lemma~\ref{estimate of the transport equation}, together with \eqref{eq:LA-low-estimate}, gives, after possibly shortening $J$,
	$$
	\|\delta\omega\|_{C(J;C^\alpha)}
	\leq
	\frac12\|\delta\omega\|_{C(J;C^\alpha)}.
	$$
	Hence $\delta\omega=0$ on $J$.
	
	The length of $J$ depends only on $N$, $\alpha$, and
	$\|U\|_{C(I;C^{N+1+\alpha})}$.
	Repeating the construction on consecutive intervals and reversing time
	proves existence and uniqueness on all of $I$.
	
	Let $\mathcal F$ denote the right-hand side of
	\eqref{eq:projected-transport-form}.
	By \eqref{eq:LA-low-estimate},
	\begin{equation}\label{eq:linearized-transport-source-low-estimate}
	\|\mathcal F(\tau,\cdot)\|_{\alpha}
	\leq C_{\alpha}\|F(\tau,\cdot)\|_{\alpha}
	+C_{\alpha}\|\nabla U(\tau,\cdot)\|_{\alpha}\|\omega(\tau,\cdot)\|_{\alpha}.
	\end{equation}
	The characteristic estimate gives
	$$
	\begin{aligned}
		\|\omega(t,\cdot)\|_\alpha
		\leq
		C_\alpha\exp\left(C_\alpha\int_{t_0}^t\|\nabla U(s,\cdot)\|_\alpha\,\rd s\right)
		\|\omega_0\|_\alpha
		+C_\alpha\int_{t_0}^t
		\exp\left(C_\alpha\int_\tau^t\|\nabla U(s,\cdot)\|_\alpha\,\rd s\right)
		\|\mathcal F(\tau,\cdot)\|_\alpha\,\rd\tau.
	\end{aligned}
	$$
	Substituting \eqref{eq:linearized-transport-source-low-estimate} and applying
	Gronwall's inequality proves \eqref{eq:holder-low-estimate}.
	
	For a spatial multi-index $\theta$ with $|\theta|=N$, the commutator satisfies
	$$
	\begin{aligned}
	\|[\partial^\theta,U\cdot\nabla]\omega\|_\alpha
	\leq\sum_{0<\beta\leq\theta}\binom{\theta}{\beta}
	\|(\partial^\beta U)\cdot\nabla\partial^{\theta-\beta}\omega\|_\alpha
	\lesssim_{N,\alpha}\|\nabla U\|_\alpha\|\omega\|_{N+\alpha}
	+\|\nabla U\|_{N+\alpha}\|\omega\|_\alpha.
	\end{aligned}
	$$
	Here the last inequality follows from \eqref{eq:holder-product-tame} and interpolation.
	By \eqref{eq:LA-high-estimate},
	\begin{equation}\label{eq:linearized-transport-source-high-estimate}
	\|\mathcal F(\tau,\cdot)\|_{N+\alpha}
	\leq C_{N,\alpha}\|F(\tau,\cdot)\|_{N+\alpha}
	+C_{N,\alpha}\|\nabla U(\tau,\cdot)\|_{\alpha}\|\omega(\tau,\cdot)\|_{N+\alpha}
	+C_{N,\alpha}\|\nabla U(\tau,\cdot)\|_{N+\alpha}\|\omega(\tau,\cdot)\|_{\alpha}.
	\end{equation}
	Applying the characteristic estimate to \eqref{eq:projected-transport-form}
	and its spatial derivatives gives
	$$
	\begin{aligned}
		\|\omega(t,\cdot)\|_{N+\alpha}
		\leq{}&
		C_{N,\alpha}\exp\left(C_{N,\alpha}\int_{t_0}^t\|\nabla U(s,\cdot)\|_\alpha\,\rd s\right)
		\|\omega_0\|_{N+\alpha}\\
		&+C_{N,\alpha}\int_{t_0}^t
		\exp\left(C_{N,\alpha}\int_\tau^t\|\nabla U(s,\cdot)\|_\alpha\,\rd s\right)\Bigl(
		\|\mathcal F(\tau,\cdot)\|_{N+\alpha}
		+\|\nabla U(\tau,\cdot)\|_{N+\alpha}\|\omega(\tau,\cdot)\|_\alpha
		\Bigr)\,\rd\tau.
	\end{aligned}
	$$
	Substitute \eqref{eq:linearized-transport-source-high-estimate} and apply
	Gronwall's inequality. Using \eqref{eq:holder-low-estimate}, we obtain
	\eqref{eq:holder-estimate} with $C_{N,\alpha}\geq C_\alpha$ sufficiently large.
	
	We now construct the pressure. Since $\omega$ is divergence-free and has zero
	spatial mean, \eqref{eq:projected-linearized} shows that the vector field
	$
	\partial_t\omega+\Div(U\otimes\omega+\omega\otimes U)-F
	$
	has zero spatial mean and vanishing Leray projection. It is therefore a
	gradient field. This determines the unique zero-mean pressure $\pi$ for which
	$(\omega,\pi)$ solves \eqref{eq:linearized-equation}.
	
	Set $\Psi:=\mathcal B\omega$. Remark~\ref{rem:app-holder-besov-equivalence}, together with the multiplier estimate \eqref{eq:app-multiplier-besov-estimate} with $m=-1$, gives
	\eqref{eq:vector-potential-holder-estimate}.
	Define the source vector potential by
	\begin{equation}\label{eq:source-vector-potential}
		G:=\mathcal B F=-\Delta^{-1}\curl F
	\end{equation}
	so that
	$
	\curl G=F.
	$
	For a divergence-free, zero-mean vector field $g\in C^{N+\alpha}(\mathbb T^3)$,
	define
	\begin{equation}\label{eq:vector-potential-lower-order-def}
		\mathcal B_Ug
		:=
		\mathcal B\Div\left(
		U\otimes\curl g+\curl g\otimes U
		\right)
		-U\cdot\nabla g.
	\end{equation}
	Applying $\mathcal B$ to \eqref{eq:projected-linearized}, and using
	$
	\mathcal B\mathbb P_H=\mathcal B
	$
	and
	$
	\omega=\curl\Psi,
	$
	we obtain
	\begin{equation}\label{eq:transport-vector-potential-equation}
		\left\{
		\begin{aligned}
			&\partial_t\Psi+U\cdot\nabla\Psi+\mathcal B_U\Psi=G,\\
			&\Psi|_{t=t_0}=\Psi_0:=\mathcal B\omega_0.
		\end{aligned}
		\right.
	\end{equation}
	
	To cancel the first-order terms in
	\eqref{eq:vector-potential-lower-order-def}, set
	$
	w:=\curl g.
	$
	Since $\Div w=0$, one has
	$$
	\Div(U\otimes w+w\otimes U)
	=
	\sum_{j=1}^3\partial_j(U_jw)+(w\cdot\nabla)U.
	$$
	Define
	$
	T_j:=\mathcal B\partial_j\curl,\  j=1,2,3.
	$
	Each $T_j$ is a constant-coefficient Fourier multiplier of order one. Since $g$ is divergence-free and has zero spatial mean,
	$
	T_jg
	=
	\partial_j\mathcal B\curl g
	=
	\partial_jg.
	$
	Using
	$$
	\curl(U_jg)=U_j\curl g+\nabla U_j\times g,
	$$
	we obtain
	$$
	\mathcal B\partial_j(U_jw)
	=
	T_j(U_jg)
	-
	\mathcal B\partial_j(\nabla U_j\times g).
	$$
	
	Define the tensor field $\cC_Ug$ by
	$
	(\cC_Ug)_{ik}:=(g\times\nabla U_i)_k.
	$
	The terms containing second derivatives of $U_i$ vanish by the
	symmetry of its Hessian and the antisymmetry of the cross product. Hence
	$$
	\Div(\cC_Ug)
	=
	(\curl g\cdot\nabla)U
	=
	(w\cdot\nabla)U.
	$$
	Since
	$
	T_j(U_jg)
	=
	U_jT_jg+[T_j,U_j]g
	=
	U_j\partial_jg+[T_j,U_j]g,
	$
	the term $U\cdot\nabla g$ cancels in
	\eqref{eq:vector-potential-lower-order-def}, giving
	\begin{equation}\label{eq:vector-potential-lower-order-decomposition}
		\begin{aligned}
			\mathcal B_Ug
			=
			\sum_{j=1}^3[T_j,U_j]g-
			\sum_{j=1}^3\mathcal B\partial_j(\nabla U_j\times g)
			+
			\mathcal B\Div(\cC_Ug).
		\end{aligned}
	\end{equation}
	
	Since $T_j=\mathbb P_H\partial_j$, we have
	$[T_j,U_j]g=[\mathbb P_H,U_j\partial_j]g
	+\mathbb P_H((\partial_jU_j)g)$.
	Subtracting the spatial mean of $U_j$ leaves the first term unchanged.
	Apply Proposition~\ref{prop:CZ commutator}, the mean-projection commutator
	estimate, and Poincar\'e's inequality to this term. Using
	\eqref{eq:holder-product-tame} for the second term, we obtain
	$$
	\|[T_j,U_j]g\|_{\alpha}
	\lesssim_\alpha
	\|\nabla U\|_{\alpha}
	\|g\|_{\alpha}.
	$$
	Since $C^\alpha$ is an algebra and
	$\mathcal B\partial_j$ and $\mathcal B\Div$ are of order zero,
	the decomposition \eqref{eq:vector-potential-lower-order-decomposition} gives
	\begin{equation}\label{eq:lower-order-besov-bound-low}
		\|\mathcal B_Ug\|_{\alpha}
		\leq
		C_\alpha
		\|\nabla U\|_{\alpha}
		\|g\|_{\alpha}.
	\end{equation}
	The same argument with \eqref{eq:holder-product-tame} and Proposition~\ref{prop:CZ commutator} gives
	\begin{equation}\label{eq:vector-potential-lower-order-bound-high}
		\|\mathcal B_Ug\|_{N+\alpha}
		\leq
		C_{N,\alpha}\|\nabla U\|_{\alpha}\|g\|_{N+\alpha}
		+
		C_{N,\alpha}\|\nabla U\|_{N+\alpha}\|g\|_{\alpha}.
	\end{equation}
	
	To prove \eqref{eq:besov-estimate}, use the characteristic estimate for
	\eqref{eq:transport-vector-potential-equation} with source
	$G-\mathcal B_U\Psi$. By \eqref{eq:lower-order-besov-bound-low}, we have
	$$
	\begin{aligned}
		\|\Psi(t,\cdot)\|_{\alpha}
		\leq{}&
		C_\alpha\exp\left(C_\alpha\int_{t_0}^t\|\nabla U(s,\cdot)\|_\alpha\,\rd s\right)
		\|\Psi_0\|_\alpha\\
		&+C_\alpha\int_{t_0}^t
		\exp\left(C_\alpha\int_\tau^t\|\nabla U(s,\cdot)\|_\alpha\,\rd s\right)
		\left(
		\|G(\tau,\cdot)\|_\alpha
		+C_\alpha\|\nabla U(\tau,\cdot)\|_\alpha\|\Psi(\tau,\cdot)\|_\alpha
		\right)\rd\tau.
	\end{aligned}
	$$
	Gronwall's inequality proves \eqref{eq:besov-estimate}, using
	$G=\mathcal B F$ and $\Psi_0=\mathcal B\omega_0$.
	
	We finish with the pressure estimates, using \eqref{eq:pressure-linearized}.
	Apply \eqref{eq:app-multiplier-besov-estimate} with $m=0$ to the zero-mean
	part of $U\otimes\omega+\omega\otimes U$.
	Using \eqref{eq:holder-product-tame} and
	Remark~\ref{rem:app-holder-besov-equivalence}, we obtain
	\eqref{eq:pressure-holder-estimate} and $\pi\in C(I;C^{N+\alpha}(\mathbb T^3))$.
	For the estimate in terms of $\Psi$, $\Div\omega=0$ and commutation of spatial
	derivatives give
	$$
	\Div\Div(U\otimes\omega+\omega\otimes U)=2\Div((\omega\cdot\nabla)U).
	$$
	 Since
	$
	\Div(\cC_U\Psi)=(\curl\Psi\cdot\nabla)U=(\omega\cdot\nabla)U,
	$
	the zero-mean normalization of $\pi$ gives
	\begin{equation}\label{eq:pressure-vector-potential-identity}
		\pi=-2\Delta^{-1}\Div\Div(\cC_U\Psi).
	\end{equation}
	Apply \eqref{eq:app-multiplier-besov-estimate} with $m=0$ to the zero-mean
	part of $\cC_U\Psi$. Using \eqref{eq:holder-product-tame} and
	Remark~\ref{rem:app-holder-besov-equivalence}, we have
	$$
	\begin{aligned}
		\|\pi(t,\cdot)\|_{N+\alpha}
		\leq C_{N,\alpha}\big(
		&\|\nabla U(t,\cdot)\|_{\alpha}\|\Psi(t,\cdot)\|_{N+\alpha}+\|\nabla U(t,\cdot)\|_{N+\alpha}\|\Psi(t,\cdot)\|_{\alpha}
		\big).
	\end{aligned}
	$$
	This proves \eqref{eq:pressure-besov-estimate}.
\end{proof}

\begin{rem}\label{rem:smooth-data-smooth-solution}
	For smooth data satisfying the hypotheses of Proposition~\ref{prop:linearized-equation},
	the solution with zero-mean pressure satisfies
	$
	\omega\in C^\infty(I\times\mathbb T^3),
	\
	\pi\in C^\infty(I\times\mathbb T^3).
	$
	Applying the proposition
	for every $N$ and using uniqueness gives
	$\omega\in C(I;C^\infty(\mathbb T^3))$. Repeated differentiation of
	\eqref{eq:projected-linearized}
	proves smoothness in time. Equation \eqref{eq:pressure-linearized} and
	the zero-mean normalization imply that $\pi$ is smooth.

\end{rem}

\noindent{\bf Acknowledgment.}
P.~Qu acknowledges support from the NSFC  Grants No.~62588101 and No.~12431007. S.~Mao acknowledges support from the Fundamental Research Funds for the Central Universities  (No. 2232026D46).\\

\noindent\textbf{Competing interests.}
The authors have no relevant financial or non-financial interests to disclose.\\

\noindent\textbf{Data availability.}
No datasets were generated or analysed during the current study.\\

\noindent\textbf{Declaration of AI use.}
ChatGPT (OpenAI) was used to assist with language editing and checks of proof steps and cross-references. The authors take full responsibility for the content.

\end{document}